\documentclass[a4paper,twoside,openany,index=totoc]{scrbook}
\usepackage[ngerman,USenglish]{babel}
\usepackage{tikz}
\usepackage{graphicx}
\usepackage{xcolor}
\usepackage[utf8]{inputenc}

\usepackage{amsmath}
\usepackage{amssymb}
\usepackage{amsfonts}
\usepackage{amsthm}
\usepackage{bm}
\usepackage{mathrsfs}

\usepackage{url}

\usepackage{makeidx}
\makeindex

\usepackage[in,headings]{fullpage}
\usepackage{setspace}
\usepackage[citestyle=numeric,style=numeric,backend=biber]{biblatex}
\usepackage{csquotes}     

\usepackage[singlespacing=true, headsepline=0.4pt]{scrlayer-scrpage}
\clearpairofpagestyles
\rohead{\pagemark}    
\lehead{\pagemark}    

\setuptoc{toc}{totoc}

\usepackage{hyperref}

\newcommand{\BTB}{\,\raisebox{.68em}{\resizebox{.67em}{!}{\ensuremath{\bm{\sim}}}}\hspace{-.83em}B} 
\newcommand{\SBTB}{\raisebox{.5em}{\resizebox{.5em}{!}{\ensuremath{\bm{\sim}}}}\hspace{-.63em}B} 

\makeatletter
\newcommand\Biggg{\bBigg@{4}}
\newcommand\biggg{\bBigg@{3}}
\makeatother

\newcommand\smallO{
  \mathchoice
    {{\scriptstyle\mathcal{O}}}
    {{\scriptstyle\mathcal{O}}}
    {{\scriptscriptstyle\mathcal{O}}}
    {\scalebox{.7}{$\scriptscriptstyle\mathcal{O}$}}
  }

\begin{document}

\newcounter{DSBcount}[section]
\renewcommand{\theDSBcount}{\arabic{chapter}.\arabic{section}.\arabic{DSBcount}}
\newcounter{Einlcount}[section]
\renewcommand{\theEinlcount}{(\Roman{Einlcount})}
\newcounter{Figcount}[section]
\renewcommand{\theFigcount}{\arabic{chapter}.\arabic{section}.\arabic{Figcount}}
\newcounter{Tabcount}[section]
\renewcommand{\theTabcount}{\arabic{chapter}.\arabic{section}.\arabic{Tabcount}}
\newcounter{titlenr}        
\newcounter{Anhangcount}[section]
\renewcommand{\theAnhangcount}{\arabic{Anhangcount}}
\newcounter{AnhNrcount}[section]
\renewcommand{\theAnhNrcount}{(A.\arabic{AnhNrcount})}
\newcounter{AnhSatzcount}[section]
\renewcommand{\theAnhSatzcount}{A.\arabic{AnhSatzcount}}
\newcounter{Dummycount}[section]
\renewcommand{\theDummycount}{\arabic{section}.\arabic{Dummycount}}

\definecolor{skyblue}{rgb}{0.53, 0.81, 0.92}
\definecolor{lightblue}{rgb}{0.68, 0.85, 0.9}
\definecolor{lightskyblue}{rgb}{0.53, 0.81, 0.98}
\definecolor{lightcyan}{rgb}{0.88, 1.0, 1.0}
\definecolor{pastelred}{rgb}{1.0, 0.41, 0.38}
\definecolor{brilliantlavender}{rgb}{0.96, 0.73, 1.0}
\definecolor{lavenderrose}{rgb}{0.98, 0.63, 0.89}
\definecolor{lightcoral}{rgb}{0.94, 0.5, 0.5}

\setcounter{section}{0}
\setcounter{titlenr}{0}
\setcounter{DSBcount}{0}
\setcounter{Einlcount}{0}
\setcounter{Figcount}{0}
\setcounter{Tabcount}{0}

\renewcommand{\theenumi}{\alph{enumi}}     
\renewcommand{\labelenumi}{\theenumi)}
\renewcommand{\theenumii}{\roman{enumii}}  
\renewcommand{\labelenumii}{\theenumii)}

\thispagestyle{empty}
\vspace*{-1.04cm}  

\begin{center}

\LARGE
\textbf{EDGEWORTH EXPANSIONS FOR}\\[0.3cm] 
\textbf{LINEAR RANK STATISTICS}\\[0.3cm]
\textbf{USING STEIN'S METHOD}

\vspace{1.3cm}

\large
by\\[0.2cm]
\Large
WALTER SCHNELLER

\vspace{12.2cm}
\Large
\textcopyright\ Copyright by Walter Schneller 2025\\[0.2cm]
All Rights Reserved

\end{center}
\setcounter{page}{0}
                       
\thispagestyle{empty}
\vspace*{-1.04cm} 
\begin{flushleft}
\begin{tabular}{@{}l@{\hspace*{2ex}}l}
Walter Schneller\\
\href{https://www.thws.de/en/}{Technische Hochschule W{\"u}rzburg-Schweinfurt}\\
\href{https://www.thws.de/en/about-thws/faculties/applied-natural-sciences-and-humanities/}
{Fakult{\"a}t Angewandte Natur- und Geisteswissenschaften}\\
Postal address: M{\"u}nzstrasse 12, 97070 W{\"u}rzburg, Germany\\
e-mail: \href{mailto:walter.schneller@thws.de}{\nolinkurl{walter.schneller@thws.de}}
\end{tabular}
\end{flushleft}
\vspace*{11.3cm}
This work is an English translation of an updated version of the doctoral thesis {''Edgeworth-Entwicklungen
f{\"u}r lineare Rangstatistiken''}, published at the Technische Universit{\"a}t Berlin in the year 1987.
\vspace{0.7cm}
\begin{flushleft}
Promotionsausschuss (Doctoral Committee):\\[3ex]
\begin{tabular}{@{}l@{\hspace*{2ex}}l}
Vorsitzender:&Prof. Dr. K. Lehmann\\[0.5ex]
Berichter:&Prof. Dr. E. Bolthausen (Principal Adviser)\\[0.5ex]
&Prof. Dr. R. Scha{\ss}berger
\end{tabular}\\[3ex]
\begin{tabular}{@{}l@{\hspace*{2ex}}l}
Tag der Einreichung (Date of Submission):&17.11.1986\\[0.5ex]
Tag der m{\"u}ndlichen Pr{\"u}fung (Date of Oral Examination):&04.03.1987
\end{tabular}
\end{flushleft}
\setcounter{page}{0}
              
\thispagestyle{empty}
\chapter*{Abstract}

\thispagestyle{empty}    

Let $\mathscr{F}_{\!A}$ be the distribution function of the standardized linear rank statistic for the matrix $A$, 
where it is assumed that the rank vector is equally distributed on the permutations. In this work, 
we analyze the conditions to be imposed on $A$ so that $\mathscr{F}_{\!A}$ 
has first and second order Edgeworth expansions with asymptotically {''sufficiently''} small remainder terms. 
The methods used are the Stein method combined with an extension of techniques 
that go back to Bolt\-hau\-sen (1984).\\[2.8ex]
The conditions obtained are very similar to the necessary and sufficient conditions 
found by Bickel and Robinson (1982). In this paper they investigated
when the distribution function of a sum of independent and identically distributed random variables has Edgeworth expansions with asymptotically small enough residues. However, these conditions are often difficult to prove directly.\\[2.8ex]
For simple linear rank statistics, however, it is possible to use a result from van Zwet (1982) to verify these assumptions. Thus we obtain conditions on $A$ for the validity of Edgeworth expansions, which on the one hand are very easy to prove and on the other hand are much more general than all previously known conditions.\\[2.8ex]
Finally, this result is applied to the special case of approximating and exact scores.
\cleardoublepage
\setcounter{page}{1}
\rehead{Contents}    
\tableofcontents
\rehead{Introduction}    

\chapter*{Introduction, summary and acknowledgements of the original 1987 version}

\addcontentsline{toc}{chapter}{\text{Introduction, summary and acknowledgements of the original 1987 version}}

Let $A$ be an $n{\times}n$\,-\,matrix of real numbers, let $\pi$ be uniformly distributed on the 
set $\mathscr{P}_{n}$ 
\index{permutation sets!$\mathscr{P}_{n}$}of permutations of $\{1,\ldots,n\}$ and
\[T_{A} = \sum_{i=1}^n a_{i\pi(i)}.\]
In this work, we investigate when the distribution function 
$\mathscr{F}_{\!A}$\index{rank statistic}\index{rank statistic!distribution function}
\index{rank statistic!distribution function!$F_{A}$, $\mathscr{F}_{A}$}of
\index{rank statistic!$T_{A}$, $\mathscr{T}_{A}$}
\[\mathscr{T}_{A} = \dfrac{T_{A} - E(T_{A})}{\sqrt{\text{Var}(T_{A})}}\]
possesses \textbf{Edgeworth expansions of first and second order}\index{Edgeworth expansion} such that the error made, when
approximating
$\mathscr{F}_{\!A}$
by these expansions, is asymptotically {''sufficiently''} small. Before we go into this in more detail, however, we would like to make a few brief comments on the statistical relevance of $T_{A}$.\\[2.8ex]
These random variables are of particular interest in nonparametric statistics. They arise in the following
situation:\\[2.8ex]
Let $X_{1},\ldots,X_{n}$ be independent random variables with continuous distribution functions
$F_{1},\ldots,$ $F_{n}$. Furthermore, if
$X_{1 : n} < X_{2 : n} < \ldots < 
X_{n : n}$\index{order statistic}\index{order statistic!$X_{1 : n}$, $X_{2 : n}$, $X_{n : n}$} 
denotes the order of the sequence $X_{1}, X_{2},\ldots,X_{n}$ in increasing size, 
then the rank $R_{j}$ of $X_{j}$ is defined by $X_{j} = X_{R_{j} : n}$. The mapping
\[S_{A} = \sum_{i=1}^n a_{iR_{i}}\]
is now called a \textbf{linear rank statistic}\index{rank statistic!linear}. If the matrix $A$ 
has coefficients of the form $a_{ij} = e_{i} \, d_{j}$, 
then $S_{A}$ is called a \textbf{simple linear rank statistic}\index{rank statistic!linear!simple}. 
This statistic can be used to test the null hypothesis
$H_{0} : F_{1} = F_{2} = \ldots = F_{n}$\index{null hypothesis $H_{0}$} 
against different classes of alternatives. The structure of the alternatives is reflected
in the choice of the regression constants $e_{i}$\index{regression constants} and the scores $d_{j}$\index{scores}. 
A well-known example is the \textbf{two-sample statistic}\index{regression constants!two-sample statistic}, 
for which $e_{1} = \ldots = e_{m} = 1$, $e_{m+1} = \ldots = e_{n} = 0$, and which is used to test $H_{0}$ against alternatives of the form $F_{1} = \ldots = F_{m}$, $F_{m+1} = \ldots = F_{n}$.
Further information can be found, for example, in the monograph by H{\'a}jek, \v{S}id{\'a}k and Sen \cite{sidak1999theory}. 
Under the null hypothesis $H_{0}$\index{null hypothesis $H_{0}$}, 
the random vector $(R_{1},\ldots,R_{n})$ is always uniformly distributed on the set $\mathscr{P}_{n}$ 
\index{permutation sets!$\mathscr{P}_{n}$}of permutations 
of $\{1,\ldots,n\}$, so that $S_{A}$ is distributed as $T_{A}$.\\[2.8ex]
The problem of determining critical values for these tests as simply as possible then leads, among other things, to the question of the conditions under which the distribution of $\mathscr{T}_{A}$ is approximated 
by the standard normal distribution for large $n$. 
This question and the related question of the convergence rate of this approximation have been investigated by a large number of authors. See the detailed overview of Does \cite{does1982higher}, section 1.2, 
and the paper of Bolthausen \cite{Bolthausen1984}.\\[2.8ex]
The desire for more accurate approximations than those provided by the normal approximation then takes us 
to the Edgeworth expansions of $\mathscr{F}_{\!A}$.
To date, the papers of Bickel \& Van Zwet \cite{10.1214/aos/1176344305} and Robinson \cite{10.1214/aos/1176344306} 
exist on this topic, but they only consider two-sample statistics\index{regression constants!two-sample statistic}.
Furthermore, there is a paper of Does \cite{does1982higher} (resp. \cite{10.1214/aos/1176346166}) 
for simple linear rank statistics.
Does, however, only considers the case of approximating scores\index{scores!approximating}, 
i.e. $d_{j} = J\bigl(\dfrac{j}{n+1}\bigr)$, where
$J : (0, 1) \rightarrow \mathbb{R}$ is a function that must fulfil certain regularity assumptions.\\[2.8ex]
In this work, we mostly consider completely general $a_{ij}$. The results we obtain are
presented in the following overview of the content. We also indicate the techniques used.\\[2.8ex]
Stein's method plays a central role in this work. This method was introduced by C. Stein \cite{stein1972bound} 
to prove results of the Berry-Ess\'een type. We will use this method to 
establish Edgeworth expansions for $\mathscr{F}_{\!A}$ (i.e. determine these expansions and estimate the error of the approximation by these expansions).\\[2.8ex]
To show how we proceed \textbf{principally}, we introduce in chapter \ref{EWD_Kap1} Edgeworth expansions
for sums of independent and identically distributed random variables $X_{n}$, $n \in \mathbb{N}$. 
We obtain a result (cf. Theorem \ref{EWD_1_1_03}, page \pageref{EWD_1_1_03}), which essentially 
originates from Bickel and Robinson \cite{10.1214/aop/1176993873}. 
However, their proof is based on Fourier analytical methods.\\[2.8ex]
What is remarkable about this result is that it does not assume the validity of Cram\'er's condition (i.e. $\limsup\limits_{|t| \rightarrow \infty} E(e^{itX_{1}}) < 1$)\index{Cram\'er's condition}, as comparable results do. Instead, the following more difficult condition is required:\\[2.8ex]
There exists a constant $\mathcal{C} > 0$ such that\\[2.3ex]
\refstepcounter{Einlcount}
\label{EWD_Einl_1}
\text{\hspace*{-0.8ex}\theEinlcount}
\begin{tabular}{@{\hspace*{4ex}}r@{\hspace*{0.8ex}}c@{\hspace*{0.8ex}}l}
$\sup\limits_{z \in \mathbb{R}}\, \Big|\bigl(F_{n}(z + y) - F_{n}(z)\bigr) - \bigl(F_{n}(z) - F_{n}(z - y)\bigr)\Big|$
                   &$\leq$&$\mathcal{C}\Bigl(\dfrac{1}{n} + y^2 \Bigr)$\\[2ex]
&&for all $0 \leq y \leq \dfrac{1}{\sqrt{n}}$, $n \in \mathbb{N}$.
\end{tabular}\\[2.8ex]
Here, $F_{n}$ denotes the distribution function of the standardized sum of $X_{1},\ldots,X_{n}$. In
chapter \ref{EWD_Kap4}, section \ref{EWD_Kap4_Sec2}, we show that \ref{EWD_Einl_1} really follows from
Cram\'er's condition.\\[2.8ex]
When proving the Theorem \ref{EWD_1_1_03} using Stein's method, some purely analytical results are required. 
We do not prove them in chapter \ref{EWD_Kap1} for reasons of {''clarity''} and since similar results are also used in chapter \ref{EWD_Kap3}. Instead, we list these analytical auxiliary {''considerations''} in chapter \ref{EWD_Kap2} 
and prove them there.\\[2.8ex]
In the central chapter \ref{EWD_Kap3} we derive Edgeworth expansions for $\mathscr{F}_{\!A}$. 
In order to present the results obtained in more detail, we define
\[\bm\breve{a}_{ij} = a_{ij} - \dfrac{1}{n} \sum_{l=1}^n a_{il} - \dfrac{1}{n} \sum_{k=1}^n a_{kj} + \dfrac{1}{n^2}
\sum_{k,l=1}^n a_{kl}\]
and\index{matrix!standardized!$\bm\hat{a}_{ij}$}\index{matrix!$\sigma_{\hspace*{-0.3ex}A}$}
\[\sigma_{\!A}^{2} = \dfrac{1}{n-1}\sum_{i,j=1}^n \bm\breve{a}_{ij}^{2}\ ,\ 
\hspace*{3.5ex}\bm\hat{a}_{ij} = \dfrac{1}{\sigma_{\!A}} \bm\breve{a}_{ij}\ .\]
The first and second order expansions used then have the following form:\\[2.8ex]
$\begin{array}{@{\hspace*{4ex}}l@{\hspace*{0.8ex}}c@{\hspace*{0.8ex}}l} 
e_{1,A}(x)&=&\Phi(x) - \psi(x)\, \dfrac{\lambda_{1,A}}{6}\, (x^2-1)\ ,\\[2ex]
e_{2,A}(x)&=&\Phi(x) - \psi(x)\, \biggl\{ \dfrac{\lambda_{1,A}}{6}\, (x^2-1) +   \dfrac{\lambda_{2,A}}{24}\, (x^3-3 x)
+  \dfrac{\lambda_{1,A}^{2}}{72}\, (x^5-10x^3 + 15x)\biggr\}\ .
\end{array}$\index{Edgeworth expansion}\index{Edgeworth expansion!$e_{1,A}$, $e_{2,A}$}\\[2.8ex]
Here $\Phi$ denotes the distribution function and $\psi$ the density of the standard normal distribution.
Furthermore are defined\\[2.3ex]
$\begin{array}{@{\hspace*{4ex}}l@{\hspace*{0.8ex}}c@{\hspace*{0.8ex}}l}
\lambda_{1,A}&=&\displaystyle{\dfrac{1}{n}\sum\limits_{i,j=1}^n \bm\hat{a}_{ij}^{3}}\ \ \text{and}
\end{array}$\\[3ex]
$\begin{array}{@{\hspace*{4ex}}l@{\hspace*{0.8ex}}c@{\hspace*{0.8ex}}l}
\lambda_{2,A}&=&\displaystyle{\dfrac{1}{n}\sum\limits_{i,j=1}^n \bm\hat{a}_{ij}^{4} + \dfrac{3}{n} - 
\dfrac{3}{n^2} \sum\limits_{i,j,k=1}^n \bigl(\bm\hat{a}_{ij}^{2} \bm\hat{a}_{ik}^{2} + 
\bm\hat{a}_{ij}^{2} \bm\hat{a}_{kj}^{2}\bigr)}\ .
\end{array}$\index{Edgeworth expansion!$\lambda_{1,A}$, $\lambda_{2,A}$}\\[2.8ex]
If $\mathscr{F}_{\!A}$ now fulfils a condition that has great analogies to \ref{EWD_Einl_1}, 
we can show the existence of a constant $K_{1} > 0$ that depends only on 
constants from this condition, such that\\[2.3ex]
\refstepcounter{Einlcount}
\label{EWD_Einl_2}
\text{\hspace*{-0.8ex}\theEinlcount}
\hspace*{4ex}$\displaystyle{\sup\limits_{z \in \mathbb{R}} \Big| \mathscr{F}_{\!A}(z) - e_{1,A}(z)\Big| 
                   \leq K_{1}\, \dfrac{1}{n}\sum\limits_{i,j=1}^n |\bm\hat{a}_{ij}|^{4}}$\hfill
(cf. Theorem \ref{EWD_3_1_10}, page \pageref{EWD_3_1_10}).\\[2.8ex]
Under somewhat more complicated conditions, the following also holds\\[2.3ex]
\refstepcounter{Einlcount}
\label{EWD_Einl_3}
\text{\hspace*{-0.8ex}\theEinlcount}
\hspace*{4ex}$\displaystyle{\sup\limits_{z \in \mathbb{R}} \Big| \mathscr{F}_{\!A}(z) - e_{2,A}(z)\Big| 
                   \leq K_{2}\, \dfrac{1}{n}\sum\limits_{i,j=1}^n |\bm\hat{a}_{ij}|^{5}}$\hfill
(cf. Theorem \ref{EWD_3_1_13}, page \pageref{EWD_3_1_13}).\\[2ex]
It should be noted that in many typical application cases 
$\displaystyle{\dfrac{1}{n}\sum\limits_{i,j=1}^n |\bm\hat{a}_{ij}|^{4}}$ 
is of the order of magnitude of $\dfrac{1}{n}$ and
$\displaystyle{\dfrac{1}{n}\sum\limits_{i,j=1}^n |\bm\hat{a}_{ij}|^{5}}$ 
is of the order of magnitude of $\dfrac{1}{n^{3/2}}$.\\[2.8ex]
The proof of \ref{EWD_Einl_2} and \ref{EWD_Einl_3} takes up the whole chapter \ref{EWD_Kap3}. 
In section \ref{EWD_Kap3_Sec1} we formulate the Theorems \ref{EWD_3_1_10} and \ref{EWD_3_1_13} and make a few basic remarks about them. The sections \ref{EWD_Kap3_Sec2}, \ref{EWD_Kap3_Sec3} and \ref{EWD_Kap3_Sec4} 
contain some essential preliminary considerations for the proof of the two theorems, 
such as the truncation technique in section \ref{EWD_Kap3_Sec4}. 
The section \ref{EWD_Kap3_Sec5} is one of the most important in this work.
In this section we extend some techniques going back to Bolthausen \cite{Bolthausen1984}, 
which ensure that the structure of the proof we learned in chapter \ref{EWD_Kap1} is also portable to linear rank statistics (under the null hypothesis $H_{0}$)\index{null hypothesis $H_{0}$}.
In section \ref{EWD_Kap3_Sec6} we use these techniques for the first time 
to prove another ''auxiliary result''. This is:\\[2.8ex]
If $k \in \mathbb{N}$, there exists a constant $C(k) > 0$ that depends only on $k$, such that\\[2.3ex]
\refstepcounter{Einlcount}
\label{EWD_Einl_4}
\text{\hspace*{-0.8ex}\theEinlcount}
\hfill$\displaystyle{\sup\limits_{z \in \mathbb{R}}\, 
\bigg| E\Bigl(|\,\mathscr{T}_{A}|^k\, 1_{(- \infty, z]}(\,\mathscr{T}_{A}) \Bigr) 
- \int_{- \infty}^{z} |x|^k \psi(x) dx \bigg| \leq
C(k)\,\biggl( \dfrac{1}{n}\sum\limits_{i,j=1}^n |\bm\hat{a}_{ij}|^{3} +
\dfrac{1}{n}\sum\limits_{i,j=1}^n |\bm\hat{a}_{ij}|^{k+4} \biggr)}$\\[2ex]
\hspace*{1ex}\hfill (cf. Theorem \ref{EWD_3_6_01}, page \pageref{EWD_3_6_01}).\\[2ex]
This result for $k = 1,2$ is then used together with the techniques from section
\ref{EWD_Kap3_Sec5} and Stein's method to derive the Theorem \ref{EWD_3_1_10} in section \ref{EWD_Kap3_Sec7} 
and the Theorem \ref{EWD_3_1_13} in section \ref{EWD_Kap3_Sec8}.\\[2.8ex]
However, the assumptions of the Theorems \ref{EWD_3_1_10} and \ref{EWD_3_1_13} have the disadvantage that they are often difficult to prove directly. Therefore, in chapter \ref{EWD_Kap4} 
we give easily verifiable conditions for simple linear rank statistics
under which \ref{EWD_Einl_2} or \ref{EWD_Einl_3} hold (except for small correction factors of
size $(\log n)^2$ or $n^{\epsilon}$, $\epsilon > 0$).
These conditions\index{conditions of van Zwet} originate from van Zwet \cite{vanZwet1982}, 
who used them to derive an estimate of the characteristic function of $\mathscr{F}_{\!A}$ 
(see Theorem \ref{EWD_4_3_01}, page \pageref{EWD_4_3_01}). 
From this estimate, we obtain the conditions of the two theorems 
(modulo the above correction factors) in section \ref{EWD_Kap4_Sec1} 
and section \ref{EWD_Kap4_Sec3}. It follows that\\[2.3ex]
\refstepcounter{Einlcount}
\label{EWD_Einl_5}
\text{\hspace*{-0.8ex}\theEinlcount}
\hspace*{4ex}$\displaystyle{\sup\limits_{z \in \mathbb{R}} \Big| \mathscr{F}_{\!A}(z) - e_{1,A}(z)\Big| 
                   \leq K_{3}\, (\log n)^2\, \dfrac{1}{n}\sum\limits_{i,j=1}^n |\bm\hat{a}_{ij}|^{4}}$\\[2.8ex]
and for any $\epsilon > 0$\\[2.3ex]
\refstepcounter{Einlcount}
\label{EWD_Einl_6}
\text{\hspace*{-0.8ex}\theEinlcount}
\hspace*{4ex}$\displaystyle{\sup\limits_{z \in \mathbb{R}} \Big| \mathscr{F}_{\!A}(z) - e_{2,A}(z)\Big| 
                   \leq K_{4}\, n^{\epsilon}\, \dfrac{1}{n}\sum\limits_{i,j=1}^n |\bm\hat{a}_{ij}|^{5}}$
\hfill (cf. Theorem \ref{EWD_4_3_13}, page \pageref{EWD_4_3_13}).\\[2.3ex]
Here, $K_{3}$ and $K_{4}$ depend only on the constants from van Zwet's conditions, 
and $K_{4}$ additionally depends on $\epsilon$.
If these conditions are strengthened slightly, it is even the case that
$\displaystyle{\dfrac{1}{n}\sum\limits_{i,j=1}^n |\bm\hat{a}_{ij}|^{4}}$ 
is of the order of magnitude of $\dfrac{1}{n}$ 
and $\displaystyle{\dfrac{1}{n}\sum\limits_{i,j=1}^n |\bm\hat{a}_{ij}|^{5}}$ 
is of the order of magnitude of $\dfrac{1}{n^{3/2}}$
(cf. Lemma \mbox{\rule[0ex]{0ex}{3.4ex}\ref{EWD_3_1_18}, \ref{EWD_3_1_18_BWa}) in combination with}
Lemma \ref{EWD_4_3_06}, \ref{EWD_4_3_06_BWc})).\\[2.8ex]
These results are then applied to the special case of approximating and exact scores 
in section \ref{EWD_Kap4_Sec4} of chapter \ref{EWD_Kap4} (see Theorem \ref{EWD_4_4_21}, page \pageref{EWD_4_4_21} 
and Theorem \ref{EWD_4_4_35}, page \pageref{EWD_4_4_35}).
We obtain more far-reaching results than Does (cf. Remark \ref{EWD_4_4_28}, \ref{EWD_4_4_28_BWc}), 
page \pageref{EWD_4_4_28_BWc}).\\[2.8ex]
I received the suggestion for this 
doctoral thesis from Prof. Dr. E. Bolthausen. 
I would like to thank him for this and for his valuable expert advice. 
I would also like to thank Mrs F. Siwak and Mrs G. Lindner-Rapha\"{e}l 
for the careful preparation of the original manuscript (prepared with typewriters).  

\newpage

\thispagestyle{empty}

\refstepcounter{Einlcount}
\label{EWD_Einl_7}

\vspace*{7ex}

\usekomafont{disposition}\usekomafont{chapter}{Comments on the updated version of 2025}

\vspace*{3ex}

\addcontentsline{toc}{chapter}{\text{Comments on the updated version of 2025}}

\normalsize

\normalfont

The entire work was written with \LaTeX\ and then translated to English with the help
of an artificial intelligence machine translation tool.
In addition, it was partly extended for better comprehensibility 
and (printing) errors were eliminated. Most of the changes were made in the sections 
\ref{EWD_Kap3_Sec5} and \ref{EWD_Kap4_Sec4}.\\[2.8ex] 
In section \ref{EWD_Kap3_Sec5} the too simple definition of the distribution of 
$\underline{I} = (\,I_{1},\ldots,I_{16}\,)$ was corrected (cf. (\ref{EWD_3_5_07}) - (\ref{EWD_3_5_10})). 
The distribution of $(I_{1}, I_{2}, I_{3}, I_{4})$ now corresponds exactly with the distribution 
of $(I_{1}, I_{2}, I_{3}, I_{4})$ in Bolthausen's paper \cite{Bolthausen1984} 
(derivable from (\ref{EWD_3_5_11})).\\[2.8ex]
In section \ref{EWD_Kap4_Sec4} the term {''condition $V_{\alpha}$''} was introduced (cf. Lemma \ref{EWD_4_4_08}). 
This term allows a better and shorter formulation of the results of this section.\\[2.8ex]
In addition, more precise estimates were made in the sections \ref{EWD_Kap1_Sec4} and \ref{EWD_Kap1_Sec5},
so that the constant $\mathcal{K}$ in Theorem \ref{EWD_1_1_03} 
has become smaller and slightly better, but certainly still not optimal.\\[2.8ex]
Furthermore, it was added in section \ref{EWD_Kap3_Sec1} that the condition (\ref{EWD_3_1_15}) for $B = \bm\hat{A}$ 
can also be derived from (\ref{EWD_3_1_16}) with another constant 
(cf. Proposition \ref{EWD_3_1_20}, \ref{EWD_3_1_20_BWc}) and \ref{EWD_3_1_20_BWd})).\\[2.8ex]
Last but not least, the proof of (\ref{EWD_3_8_32}) in subsection \ref{EWD_Kap3_Sec8_3} 
was improved and shortened.
We now use that the function $r'_{z}$ is a 
linear combination of functions $\in \mathcal{H}$ (cf. (\ref{EWD_0_1_01})).\index{function!$\mathcal{H}$}\\[2.8ex]
Essential ideas from this work have been published in \cite{10.1214/aos/1176347258} for first-order 
Edgeworth expansions and with the restriction to standardized matrices $\bm\hat{A}$
whose elements are bounded in absolute value by one.\\[2.8ex]
Moreover, the techniques of Bolthausen \cite{Bolthausen1984} and here can also be used to give a 
short proof for the normal approximation of $\mathscr{T}_{A}$ 
for large $n$ (see Schneller \cite{Schneller1988}).\\[2.8ex]
Finally, it should be noted that up-to-date bibliographies on the topic of convergence rates for 
combinatorial central limit theorems can be found, for example, in 
Chen and Fang \cite{10.3150/13-BEJ569} and Roos \cite{doi:10.1137/S0040585X97T990794}.
                        
\rehead{Notations}    
\setcounter{chapter}{-1}
\chapter[Notations and preliminary remarks]{Notations and preliminary remarks}\label{EWD_Kap0}

\setcounter{section}{1}

\textbf{a) Numbers, sequences, quantities}\\[1ex]
For $x \in \mathbb{R}$ we define:\\[2ex]
\hspace*{12.1ex}
\begin{tabular}{@{}l@{\hspace*{20ex}}l@{}}
$x^{+} = \text{max}\{x,\,0\}$&(positive part);\\[1.5ex]
$\lfloor x \rfloor = \text{max}\{k \in \mathbb{Z} : k \leq x\}$&(floor function).
\end{tabular}\\[2.5ex]
Furthermore, if $(a_{n})_{n \in \mathbb{N}}$ and $(b_{n})_{n \in \mathbb{N}}$ are sequences of real numbers
(with $b_{n} \not= 0$), we use the notations\\[2.5ex]
\hspace*{12.1ex}
\begin{tabular}{@{}l@{\hspace*{2ex}}c@{\hspace*{2ex}}l}
$a_{n} \sim b_{n}$&$\Leftrightarrow$&$\lim\limits_{n \rightarrow \infty} \dfrac{a_{n}}{b_{n}} = 1$;\\[2.5ex]
$a_{n} = \smallO(b_{n})$&$\Leftrightarrow$&$\lim\limits_{n \rightarrow \infty} \dfrac{a_{n}}{b_{n}} = 0$;\\[2.5ex]
$a_{n} = \mathcal{O}(b_{n})$&$\Leftrightarrow$&the sequence $\Bigl(\dfrac{a_{n}}{b_{n}}\Bigr)_{n \in \mathbb{N}}$
is bounded.
\end{tabular}\\[2.5ex]
In the following, we also use the symbol $\mathbb{N}_{0}$ 
for the set of natural numbers including zero.\\[3.3ex]
\textbf{b) Functions}\index{function}\\[1ex]
In this work, $\psi(x)$ denotes the density and $\Phi(x)$ the distribution function of the standard
normal distribution.\index{standard normal distribution}\index{standard normal distribution!density}
\index{standard normal distribution!density!$\psi(x)$}\index{standard normal distribution!distribution function}\index{standard normal distribution!distribution function!$\Phi(x)$}In addition, 
we write $\Phi(h)$ or $\Phi(h(x))$ for the expected value of the function $h$ with respect 
to the standard normal distribution.\index{standard normal distribution!expected value of $h$}\index{standard normal distribution!expected value of $h$!$\Phi(h)$, $\Phi(h(x))$}\\[2ex]
Furthermore, we define\\[2ex]
\refstepcounter{DSBcount}
\label{EWD_0_1_01}
\text{\hspace*{-0.8ex}(\theDSBcount)}
\hspace*{4ex}
\begin{tabular}{@{}l}
$\mathcal{H} = \Bigl\{\,h : \mathbb{R} \rightarrow [\,0,\,1\,]\,:\,\text{$h$ is monotonically decreasing}$\\
\hspace*{50.5ex}$\text{and continuous from the left}\,\Bigr\}$.
\end{tabular}\index{function!monotonically decreasing}\index{function!continuous from the left}\index{function!$\mathcal{H}$}\\[2ex]
If $F : \mathbb{R} \rightarrow \mathbb{R}$ is a function, we set\\[2ex]
\hspace*{12.1ex}
\begin{tabular}{@{}l}
$||F||$\hspace*{2ex}or\hspace*{2ex}$||F(z)|| = \sup\limits_{z \in \mathbb{R}} |F(z)|$.
\end{tabular}\index{function!${"|}{"|}F{"|}{"|}$, ${"|}{"|}F(z){"|}{"|}$}\\[2ex]
Next, for $y \in \mathbb{R}$ we denote by $\Delta_{y}F$ the function 
$\mathbb{R} \rightarrow \mathbb{R}$, which is defined by
\index{difference}\index{difference!$\Delta_{y}F$}\\[2ex]
\refstepcounter{DSBcount}
\label{EWD_0_1_02}
\text{\hspace*{-0.8ex}(\theDSBcount)}
\hspace*{3.3ex}
\begin{tabular}{@{}l}
$\Delta_{y}F(z) = F(z + y) - F(z)$\hspace*{2ex}for $z \in \mathbb{R}$
\end{tabular}\\[2ex]
$\Delta_{y}F$ is the \textbf{first difference}\index{difference!first} of $F$ with respect to $y$. The
\textbf{\textit{k}th difference}\index{difference!\textit{k}th} of $F$ with respect to $y$ 
is then obtained by applying the operator $\Delta_{y}$ $k$ times, i.e.\\[2ex]
\refstepcounter{DSBcount}
\label{EWD_0_1_03}
\text{\hspace*{-0.8ex}(\theDSBcount)}
\hspace*{4ex}
\begin{tabular}{@{}l@{\hspace*{0.8ex}}c@{\hspace*{0.8ex}}l}
$\Delta_{y}^{k}F(z)$&$=$&$\displaystyle{\Delta_{y}\bigl(\Delta_{y}\ldots(\Delta_{y}F)\bigr)(z)}$\\[2ex]
&$=$&$\displaystyle{\sum\limits_{j = 0}^{k}\,(-1)^{k-j}\, \dbinom{k}{j}\, 
F(z + jy)}$\hspace*{2ex}for $z \in \mathbb{R}$, $k \in \mathbb{N}$.
\end{tabular}\index{difference!$\Delta_{y}^{k}F$}\\[2.5ex]
With this notation, the \textbf{interpolating polynomial}\index{} 
of $F$ at the points $z$, $z + y$, \ldots , $z + ky$ of degree at most $k \in \mathbb{N}_{0}$ 
can be written as\\[2ex] 
\refstepcounter{DSBcount}
\label{EWD_0_1_04}
\text{\hspace*{-0.8ex}(\theDSBcount)}
\hspace*{4ex}
\begin{tabular}{@{}l@{\hspace*{0.8ex}}c@{\hspace*{0.8ex}}l}
$P_{y}^{k}(x;z,F)$&$=$&$\displaystyle{F(z)\, +\, \sum\limits_{s = 1}^{k}\, \Delta_{y}^{s}F(z)\, \prod\limits_{i = 1}^{s}\,
\dfrac{x - z - (i - 1)y}{iy}}$\\[3ex]
&&(where $\displaystyle{\sum \ldots = 0}$\hspace*{2ex}for $y = 0$).
\end{tabular}\index{polynomial!interpolating}\index{polynomial!interpolating!$P_{y}^{k}(x;z,F)$}\\[3.3ex]
\textbf{c) Constants}\\[1ex]
Throughout the paper, Latin letters (e.g. $c$, $c_{1}$, $C_{1}$, $e$, $E$, $n_{0}$, \ldots) and Greek letters (e.g. $\epsilon_{0}$, $\delta$, \ldots) are used to denote those constants that
depend only on the theorem and proof in which they appear. 
If there are also dependencies on quantities from the theorem or proof, this is indicated in a subordinate clause 
or by the use of round brackets (e.g. $C(k)$).\\[2ex]
In contrast, the letters $\mathcal{C}$, $\mathcal{K}$, $\mathcal{E}$, $\mathfrak{K}$, $\mathfrak{L}$ (including 
indices and embellishments) denote constants that are considered fixed for the entire work. If these constants depend on other mathematical quantities, this is also stated in a subordinate clause or indicated by round brackets.\\[3.3ex]
\textbf{d) Further notations}\\[1ex]
In H{\"o}lder's inequality, which is frequently used in this work, we denote the exponents with $p$ and $q$, 
where $p > 1$, $q > 1$ and \mbox{\rule[0ex]{0ex}{4.4ex}$\dfrac{1}{p} + \dfrac{1}{q} = 1$}. 
For example, an application of this inequality with the exponents $p$ and 
\mbox{\rule[0ex]{0ex}{2.8ex}$q = \dfrac{p}{p-1}$} to finite sequences $x_{1},\ldots,x_{\nu}$
\index{H{\"o}lder's inequality!for finite sequences using length $\nu$}
\index{H{\"o}lder's inequality!$p$, $q$, $\nu$}gives\\[2ex]
\refstepcounter{DSBcount}
\label{EWD_0_1_05}
\text{\hspace*{-0.8ex}(\theDSBcount)}
\hspace*{4ex}
\begin{tabular}{@{}l}
$\displaystyle{\Bigl(\,\sum\limits_{i = 1}^{\nu} |x_{i}|\,\Bigr)^p 
= \Bigl(\,\sum\limits_{i = 1}^{\nu} 1\,|x_{i}|\,\Bigr)^p
\leq \nu^{p-1}\,\sum\limits_{i = 1}^{\nu} |x_{i}|^p}$.
\end{tabular}\\[4ex]
Further important notations in connection with $n{\times}n$\,-\,matrices $A$ 
can be found at the beginning of the section \ref{EWD_Kap3_Sec1} of chapter \ref{EWD_Kap3}.
 
\rehead{Edgeworth expansions for iid random variables}    
\chapter[Edgeworth expansions of first order for sums of iid random variables]{Edgeworth expansions of first order for sums of independent and identically distributed random variables}\label{EWD_Kap1}

\section{Introduction and results}\label{EWD_Kap1_Sec1}

In the central part of this work, Edgeworth expansions for linear rank statistics are established 
using Stein's method (see chapter \ref{EWD_Kap3}). Since this requires somewhat extensive considerations, 
we will demonstrate in this chapter how to obtain Edgeworth expansions for independent and identically distributed (iid) random variables using Stein's method. However, for the sake of clarity, we will postpone 
some important analytical estimates required for this to the next chapter \ref{EWD_Kap2}.\\[2.8ex]
In the following, let $X_{i}$, $i \in \mathbb{N}$, be a sequence of independent and identically distributed
random variables\index{random variable} with a common distribution function $F$ and\\[2.5ex]
\refstepcounter{DSBcount}
\label{EWD_1_1_01}
\text{\hspace*{-0.8ex}(\theDSBcount)}
\hspace*{4ex}
$E(X_{i}) = 0,\ E(X_{i}^2) = 1\ \ \text{and}\ \ \beta_{4} = E(X_{i}^4) < \infty$.\\[2.5ex]
Then, of course, $\mu_{3} = E(X_{i}^3)$ and 
$\beta_{3} = E(|X_{i}|^3)$ also exist.\index{random variable!moment}
\index{random variable!moment!$\mu_{3}$, $\beta_{3}$, $\beta_{4}$}
Furthermore, let\\ 
$S_{n} = \displaystyle{\dfrac{1}{\sqrt{n}} \sum\limits_{i = 1}^n X_{i}}$\index{random variable!$X_{i}$, $S_{n}$} 
for $n \in \mathbb{N}$,
$F_{n}$ the distribution function of $S_{n}$ and\\[2ex]
\refstepcounter{DSBcount}
\label{EWD_1_1_02}
\text{\hspace*{-0.8ex}(\theDSBcount)}
\hspace*{4ex}\index{Edgeworth expansion}\index{Edgeworth expansion!$e_{n}$}
$e_{n}(x) = \Phi(x) - \dfrac{\mu_{3}}{6\,\sqrt{n}} (x^2 -1) \psi(x)\ \ \ \text{for}\ x \in \mathbb{R}$.\\[2.5ex]
With these notations we show:\\[4ex]  
\refstepcounter{DSBcount}
\label{EWD_1_1_03}
\textbf{\hspace*{-0.8ex}\theDSBcount\ Theorem}\index{Theorem!for the iid-case}\\[0.8ex]
Suppose that there exists a constant $\mathcal{C} > 0$ such that\\[2.5ex]
\refstepcounter{DSBcount}
\label{EWD_1_1_04}
\text{\hspace*{-0.8ex}(\theDSBcount)}
\hspace*{4ex}
$\displaystyle{||\Delta_{y}^{2}\, F_{n}|| \leq \mathcal{C}\Bigl(\dfrac{1}{n} + y^2 \Bigr)}$
\hspace*{2ex}for all $0 \leq y \leq \dfrac{1}{\sqrt{n}}$, $n \in \mathbb{N}$.\\[2ex]
Then there exists a constant $\mathcal{K} > 0$ depending only on $\beta_{3}$, $\beta_{4}$ and $\mathcal{C}$ 
such that\\[2ex]
\refstepcounter{DSBcount}
\label{EWD_1_1_05}
\text{\hspace*{-0.8ex}(\theDSBcount)}
\hspace*{4ex}
$\displaystyle{||F_{n} - e_{n}|| \leq \dfrac{\mathcal{K}}{n}}$
\hspace*{2ex}for all $n \in \mathbb{N}$.\\[4ex]
\refstepcounter{DSBcount}
\label{EWD_1_1_06}
\textbf{\hspace*{-0.8ex}\theDSBcount\ Remarks}
\begin{enumerate}
\item\label{EWD_1_1_06_BWa}
The proof of the theorem shows that\\[2ex]
\refstepcounter{DSBcount}
\label{EWD_1_1_07}
\text{\hspace*{-0.8ex}(\theDSBcount)}
\hspace*{4ex}
$\displaystyle{\mathcal{K} = \bigl(\,2 + \beta_{4}\,\bigr)\,\mathcal{C} 
+ \bigl(\,3 + 11\,\beta_{3} + 13\,\beta_{4} + 9\,\beta_{3}\,\beta_{4}\,\bigr)}$\\[2ex]
can be selected. See (\ref{EWD_1_5_02}) and the considerations before it.\\[0.8ex]
In this proof, however, the primary aim was to present our application of Stein's method and 
only secondarily to minimize $\mathcal{K}$. This constant is therefore 
probably still far from being an optimal constant. 
\item\label{EWD_1_1_06_BWb}
In chapter \ref{EWD_Kap4}, section \ref{EWD_Kap4_Sec2} we show, that  
from Cram\'er's condition\index{Cram\'er's condition} 
(i.e. $\limsup\limits_{|t| \rightarrow \infty} E(e^{itX_{1}}) < 1$)
we get the existence of a constant $\mathcal{C} > 0$ 
that fulfills the condition (\ref{EWD_1_1_04}).\\[0.8ex]
In particular, such a $\mathcal{C} > 0$ exists if the distribution of $X_{1}$ has a density.
See also Lemma \ref{EWD_2_4_02}, \ref{EWD_2_4_02_BWb}) for a different approach
in the case of an existing density.
\item\label{EWD_1_1_06_BWc}
If, on the other hand, $P(X_{1} = -1) = \dfrac{1}{2} = P(X_{1} = 1)$, neither (\ref{EWD_1_1_04}) nor (\ref{EWD_1_1_05}) apply, since we obtain with the help of Stirling's formula\index{Stirling's formula}:\\[2ex]
\begin{tabular}{@{\hspace*{6.3ex}}l@{\hspace*{0.8ex }}c@{\hspace*{0.8ex}}l@{\hspace*{3ex}}l@{}}
$\displaystyle{||F_{2n} - e_{2n}||}$&$\geq$&
$\displaystyle{\text{max}\Bigl\{|F_{2n}(0) - e_{2n}(0)|,\, |F_{2n}(0-) - e_{2n}(0)|\Bigr\}}$\\[2ex] 
&$\geq$&$\displaystyle{\dfrac{1}{2}\, |F_{2n}(0) - F_{2n}(0-)|}$\\[2ex]
&$=$&$\displaystyle{\dfrac{1}{2}\, P(S_{2n} = 0)} = 
\dfrac{1}{2}\, \dfrac{\dbinom{2n}{n}}{2^{2n}} \sim 
\dfrac{1}{2}\, \dfrac{1}{\sqrt{\pi n}}$\hspace*{2ex}for $n \rightarrow \infty$.
\end{tabular}\vspace*{1ex}
\item\label{EWD_1_1_06_BWd}
From (\ref{EWD_1_1_05}) we can conclude the other way round:\\[2ex]
\begin{tabular}{@{\hspace*{6.3ex}}l@{\hspace*{0.8ex }}c@{\hspace*{0.8ex}}l@{\hspace*{3.8ex}}l@{}}
$\displaystyle{||\Delta_{y}^{2}\, F_{n}||}$&$\leq$&
$\displaystyle{4\,\dfrac{\mathcal{K}}{n} + ||\Delta_{y}^{2}\, e_{n}||}$\\[2ex]
&$\leq$&$\displaystyle{4\,\dfrac{\mathcal{K}}{n} + 
\Bigl( \dfrac{1}{4} + \dfrac{1}{5}\,\dfrac{\beta_{3}}{\sqrt{n}} \Bigr)\, y^2}$
&\hspace*{-5.7ex}(see part \ref{EWD_1_2_02_BWb)}) of the proof of Proposition \ref{EWD_1_2_02})\\[2.5ex]
&$\leq$&$\displaystyle{\Bigl(4 \mathcal{K} + 
\dfrac{1}{4} + \dfrac{1}{5}\,\beta_{3} \Bigr)\, \Bigl(\dfrac{1}{n} + y^2 \Bigr)}$
&for all $y \in \mathbb{R}$, $n \in \mathbb{N}$.
\end{tabular}\vspace*{1ex}
\item\label{EWD_1_1_06_BWe}
Bickel and Robinson \cite{10.1214/aop/1176993873}, cf. Theorem and Note (ii), have shown an analogous result 
with a shorter Fourier-analytic proof. In their result, however, $\mathcal{K}$ depends on 
$\mathcal{C}$ and $F$ (instead of $\beta_{3}$, $\beta_{4}$).
\end{enumerate}
In the following proof of Theorem \ref{EWD_1_1_03}, let $n$ be fixed. Furthermore, we can assume $n \geq 6$ 
without loss of generality. For $n \leq 5$ and $\mathcal{K}$ from (\ref{EWD_1_1_07})
we trivially obtain (cf. Lemma \ref{EWD_2_2_08}, \ref{EWD_2_2_08_BWe}) and $\beta_{3} \geq 1$):\\[2ex]
\refstepcounter{DSBcount}
\label{EWD_1_1_08}
\text{\hspace*{-0.8ex}(\theDSBcount)}
\hspace*{3.8ex}
$\displaystyle{||F_{n} - e_{n}|| \leq  
||F_{n} - \Phi|| + \dfrac{\beta_{3}}{6\,\sqrt{n}}\,||\,(x^2 - 1)\,\psi\,||
\leq 1 + \dfrac{1}{15}\,\dfrac{\beta_{3}}{\sqrt{n}}
\leq \dfrac{1}{n}\,6\,\beta_{3}
\leq \dfrac{1}{n}\,\mathcal{K}}$.\\[2.5ex]
By using Stein's method, we now show that from\\[2.5ex]
\refstepcounter{DSBcount}
\label{EWD_1_1_09}
\text{\hspace*{-0.8ex}(\theDSBcount)}
\hspace*{4ex}
$\displaystyle{||\Delta_{y}^{2}\, F_{n}|| \leq \mathcal{C}\Bigl(\dfrac{1}{n} + y^2 \Bigr)}$
\hspace*{2ex}for $0 \leq y \leq \dfrac{1}{\sqrt{n}}$\hspace*{2ex}and\\[2ex]
\refstepcounter{DSBcount}
\label{EWD_1_1_10}
\text{\hspace*{-0.8ex}(\theDSBcount)}
\hspace*{2.8ex}
$\displaystyle{||\Delta_{\frac{1}{\sqrt{n-1}}}^{2}\, F_{n-1}|| \leq \mathcal{C}\, \dfrac{2}{n-1}}$\\[2.5ex]
we get the inequality\\[2ex]
\hspace*{12.1ex}$\displaystyle{||F_{n} - e_{n}|| \leq \dfrac{\mathcal{K}}{n}}$
\hspace*{2ex}with $\mathcal{K}$ according to (\ref{EWD_1_1_07})\\[2.5ex]
The corresponding proof is divided into four parts. The central considerations take place in the third section, 
while the next section contains a preliminary consideration.

\section{Transition to smooth functions}\label{EWD_Kap1_Sec2}

The starting point of this section is the following trivial equation\\[2ex]
\hspace*{12.1ex}$\displaystyle{F_{n}(z) - e_{n}(z) = \int\limits_{\mathbb{R}} 1_{(- \infty,\, z\,]}\,dF_{n} -  
\int\limits_{\mathbb{R}} 1_{(- \infty,\, z\,]}\,e'_{n}\,dx}$.\\[2ex]
In order to apply Stein's method, we must first replace the discontinuous functions
$1_{(- \infty,\, z\,]}$, $z \in \mathbb{R}$, by certain differentiable functions $q_{z}$, $z \in \mathbb{R}$.
The following piecewise defined functions prove to be suitable. For $z \in \mathbb{R}$ let 
(cf. (\ref{EWD_2_3_02}) with $\lambda = \dfrac{1}{\sqrt{n}}$)\\[2ex]
\refstepcounter{DSBcount}
\label{EWD_1_2_01}
\text{\hspace*{-0.8ex}(\theDSBcount)}
\hspace*{4ex}
$\displaystyle{q_{z}(x) =
\left\{
\begin{array}{cl}
1& \hspace*{4ex}
\text{for}\ x \leq z,\\[1ex]
1 - \dfrac{n}{2}\,\Bigl(\,x - z\,\Bigr)^2& \hspace*{4ex}
\text{for}\ z \leq x \leq z + \dfrac{1}{\sqrt{n}},\\[2ex]
\dfrac{n}{2}\,\Bigl(\,z + \dfrac{2}{\sqrt{n}} - x\,\Bigr)^2& \hspace*{4ex}
\text{for}\ z + \dfrac{1}{\sqrt{n}} \leq x \leq z + \dfrac{2}{\sqrt{n}},\\[2ex]
0 & \hspace*{4ex}
\text{for}\ z + \dfrac{2}{\sqrt{n}} \leq x\,.
\end{array}  \right.}$
\index{function!smooth!$q_{z}$}\\[3ex]
We get the following result:\\[4ex]  
\refstepcounter{DSBcount}
\label{EWD_1_2_02}
\textbf{\hspace*{-0.8ex}\theDSBcount\ Proposition}\\[2ex]
$\displaystyle{||F_{n} - e_{n}|| \leq\,
\sup\limits_{z \in \mathbb{R}}\,\Big|\int\limits_{\mathbb{R}} q_{z}\,dF_{n} - \int\limits_{\mathbb{R}} q_{z}\,e'_{n}\,dx\Big|
+ \dfrac{1}{n}\,\Bigl(\dfrac{7}{12}\,\mathcal{C} + \dfrac{1}{48} + 
\dfrac{1}{60}\,\dfrac{\beta_{3}}{\sqrt{n}}\Bigr)\,}$.\\[3ex]
\textbf{Proof:}
\begin{enumerate}
\item
Integration by parts\index{integration!by parts for Lebesgue-Stieltjes integrals} 
for Lebesgue-Stieltjes integrals (cf. e.g. \cite{athreya2006measure}, Theorem 5.2.3) gives\\[2ex]
\begin{tabular}{@{}l@{\hspace*{0.8ex}}l@{}}
&$\displaystyle{\Big|F_{n}(z) - 
\int\limits_{\mathbb{R}} q_{\textstyle{\bigl(z - \frac{1}{\sqrt{n}}\bigr)}}\,dF_{n}\Big|}$\\[3.5ex]
$=$&$\displaystyle{\Big|\int\limits_{\mathbb{R}} 1_{(- \infty,\, z\,]} - 
q_{\textstyle{\bigl(z - \frac{1}{\sqrt{n}}\bigr)}}\,dF_{n}\Big|}$
\end{tabular}\\[4.5ex]
\begin{tikzpicture}[domain=0:14]
\draw (0,3) -- (6, 3);
\draw (0,6) -- (6, 6);
\draw (3,2.8) -- (3, 6);
\draw[dotted] (1.5,3.2) -- (1.5,6);
\draw[dotted] (0,4.5) -- (3,4.5);
\draw (1.5,2.8) -- (1.5,3.2);
\draw (4.5,2.8) -- (4.5,3.2);
\draw (8,3) -- (14, 3);
\draw (11,1.5) -- (11,4.5); 
\draw[dotted] (11,4.5) -- (13.75,4.5);
\draw[dotted] (11,1.5) -- (13.5,1.5);
\draw[red,thick] (0,6) -- (1.5,6); 
\draw[red,thick] (1.5,6) .. controls (2.25,6) and (2.25,6) .. (3,4.5);
\draw[red,thick] (3,4.5) .. controls (3.75,3) and (3.75,3) .. (4.5,3);
\draw[red,thick] (4.5,3) -- (6,3);
\draw[blue,thick] (8,3) -- (9.5,3);
\draw[blue,thick] (9.5,3) .. controls (10.25,3) and (10.25,3) .. (11,4.5);
\draw[blue,thick] (11,1.5) .. controls (11.75,3) and (11.75,3) .. (12.5,3);
\draw[blue,thick] (12.5,3) -- (14,3);
\draw (9.5,2.8) -- (9.5,3.2);
\draw (12.5,2.8) -- (12.5,3.2);
\node[] at (-0.2,3) {$0$};
\node[] at (-0.2,4.5) {$\dfrac{1}{2}$};
\node[] at (-0.2,6) {$1$};
\node[] at (3,2.6) {$z$};
\node[] at (1.5,2.5) {$z - \frac{1}{\sqrt{n}}$};
\node[] at (4.5,2.5) {$z + \frac{1}{\sqrt{n}}$};
\node[] at (10.8,2.7) {$z$};
\node[] at (9.5,2.5) {$z - \frac{1}{\sqrt{n}}$};
\node[] at (12.5,2.5) {$z + \frac{1}{\sqrt{n}}$};
\node[] at (13.9,4.5) {$\dfrac{1}{2}$};
\node[] at (13.8,1.5) {$- \dfrac{1}{2}$};
\node[] at (7,4.5) {, hence};
\node[red] at (4.1,4.5) {$q_{\textstyle{\bigl(z - \frac{1}{\sqrt{n}}\bigr)}}$};
\node[blue] at (10.1,5) {$1_{(- \infty,\, z\,]} - q_{\textstyle{\bigl(z - \frac{1}{\sqrt{n}}\bigr)}}$};
\end{tikzpicture}\\[2ex]
\refstepcounter{Figcount}
\label{EWDFig_1_2_01}
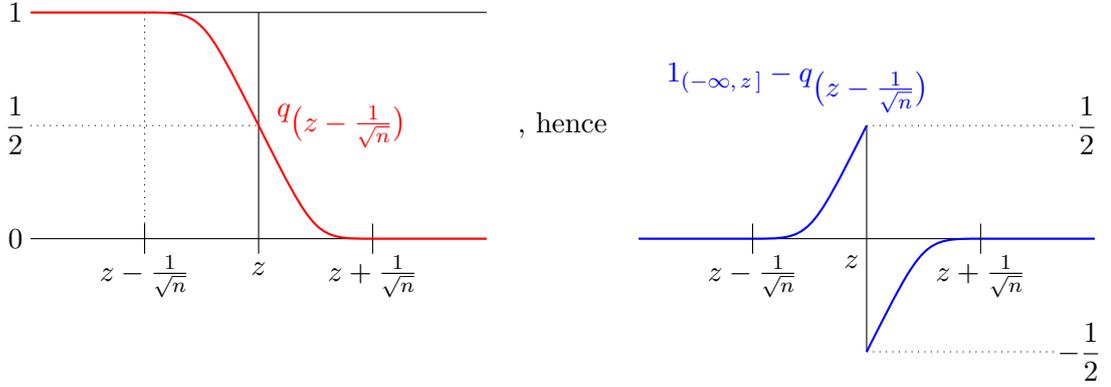
\textbf{Fig. \theFigcount:} Sketch of the functions used\\[5ex]
\begin{tabular}{@{}l@{\hspace*{0.8ex}}l@{}}
$=$&$\displaystyle{\bigg| \int\limits_{\textstyle{\bigl(z - \frac{1}{\sqrt{n}},\,z\,\bigr]}}
\dfrac{n}{2}\,\Bigl( x - \bigl(z - \dfrac{1}{\sqrt{n}}\bigr) \Bigr)^2\,dF_{n}(x)\, - 
\int\limits_{\textstyle{\bigl(z,\,z + \frac{1}{\sqrt{n}}\,\bigr]}}
\dfrac{n}{2}\,\Bigl( \bigl(z + \dfrac{1}{\sqrt{n}}\bigr) - x \Bigr)^2\,dF_{n}(x)\, \bigg|}$\\[7.5ex]
$=$&$\displaystyle{\bigg|\ \dfrac{1}{2} \, F_{n}(z)\, - \int\limits_{\textstyle{\bigl(z - \frac{1}{\sqrt{n}},\,z\,\bigr]}}
n\,\Bigl( x - \bigl(z - \dfrac{1}{\sqrt{n}}\bigr) \Bigr)\,F_{n}(x)\,dx\,\, +}$\\[6.5ex]
&\hspace*{34.5ex}$\displaystyle{\dfrac{1}{2} \, F_{n}(z)\, - 
\int\limits_{\textstyle{\bigl(z,\,z + \frac{1}{\sqrt{n}}\,\bigr]}}
n\,\Bigl( \bigl(z + \dfrac{1}{\sqrt{n}}\bigr) - x \Bigr)\,F_{n}(x)\,dx\,\bigg|}$\\[7.5ex]
$=$&$\displaystyle{\bigg|\ F_{n}(z)\, + \int\limits_{z}^{z - 1/\sqrt{n}}
n\,\Bigl( \dfrac{1}{\sqrt{n}} - \bigl(z - x\bigr) \Bigr)\,F_{n}(x)\,dx -  
\int\limits_{z}^{z + 1/\sqrt{n}}
n\,\Bigl( \dfrac{1}{\sqrt{n}} - \bigl(x - z\bigr) \Bigr)\,F_{n}(x)\,dx\,\bigg|}$.
\end{tabular}\\[3ex]
We now substitute\index{integration!by substitution} in the first integral
$y = z - x$ and in the second integral $y = x - z$. Because of
$\displaystyle{\int\limits_{0}^{1/\sqrt{n}} n\,\Bigl( \dfrac{1}{\sqrt{n}} - y \Bigr)\, dy = \dfrac{1}{2}}$ 
we then obtain\\[2ex]
\begin{tabular}{@{}l@{\hspace*{0.8ex}}l@{}}
$=$&$\displaystyle{\bigg| \int\limits_{0}^{1/\sqrt{n}}n\,\Bigl( \dfrac{1}{\sqrt{n}} - y \Bigr)\,
\Bigl\{ - F_{n}(z + y) + 2\,F_{n}(z) - F_{n}(z - y)  \Bigr\} \,dy\,\bigg|}$\\[4ex]
$\leq$&$\displaystyle{\int\limits_{0}^{1/\sqrt{n}}n\,\Bigl( \dfrac{1}{\sqrt{n}} - y \Bigr)\,
\mathcal{C}\Bigl(\dfrac{1}{n} + y^2 \Bigr) \,dy}$\hspace*{36.5ex}(because of (\ref{EWD_1_1_09}))\\[4.2ex]
$=$&$\dfrac{7}{12}\ \dfrac{\mathcal{C}}{n}$.
\end{tabular}\vspace*{2ex}
\item\label{EWD_1_2_02_BWb)}
In addition, we have for all $y \in \mathbb{R}$\\[2ex]
\begin{tabular}{@{}l@{\hspace*{0.8ex}}c@{\hspace*{0.8ex}}l@{\hspace*{2.6ex}}r@{}}
$||\Delta_{y}^{2}\, e_{n}||$&$\leq$&$y^2\,||e''_{n}||$
&(cf. Lemma \ref{EWD_2_4_02}, \ref{EWD_2_4_02_BWb}))\\[2ex]
&$\leq$&$y^2\,\Bigl(\,||x\,\psi(x)|| + \dfrac{\beta_{3}}{6\,\sqrt{n}}\,||(x^4 - 6x^2+3)\,\psi(x)||\,\Bigr)$\\[3ex]
&$\leq$&$y^2\,\Bigl(\,\dfrac{1}{4} + \dfrac{1}{5}\,\dfrac{\beta_{3}}{\sqrt{n}}\,\Bigr)$
&(cf. Lemma \ref{EWD_2_2_08}, \ref{EWD_2_2_08_BWb}) and \ref{EWD_2_2_08_BWh})).
\end{tabular}\\[2.5ex]
By proceeding in the same way as in part a) with $e_{n}(z)$ instead of $F_{n}(z)$, we get\\[2ex]
\begin{tabular}{@{}l@{\hspace*{0.8ex}}l@{}}
&$\displaystyle{\Big|e_{n}(z) - 
\int\limits_{\mathbb{R}} q_{\textstyle{\bigl(z - \frac{1}{\sqrt{n}}\bigr)}}\,e'_{n}\,dx\Big|}$\\[3.7ex]
$\leq$&$\displaystyle{\int\limits_{0}^{1/\sqrt{n}}n\,\Bigl( \dfrac{1}{\sqrt{n}} - y \Bigr)\,
y^2\,\Bigl(\,\dfrac{1}{4} + \dfrac{1}{5}\,\dfrac{\beta_{3}}{\sqrt{n}}\,\Bigr)\,dy}$\\[4.2ex]
$=$&$\dfrac{1}{n}\,\Bigl(\,\dfrac{1}{48} + \dfrac{1}{60}\,\dfrac{\beta_{3}}{\sqrt{n}}\,\Bigr)$.
\end{tabular}\\[-2.4ex]
\hspace*{1ex}\hfill$\Box$\vspace*{1ex}
\end{enumerate}

\section{The central step}\label{EWD_Kap1_Sec3}

According to Proposition \ref{EWD_1_2_02}, the term\\[2ex]
\hspace*{12.1ex}$\displaystyle{\sup\limits_{z \in \mathbb{R}}\,\Big|\int\limits_{\mathbb{R}} q_{z}\,dF_{n} - 
\int\limits_{\mathbb{R}} q_{z}\,e'_{n}\,dx\Big|}$\\[2.5ex]
still needs to be estimated. To do this, we use Stein's method to derive an equation of the type\\[2ex]
\refstepcounter{DSBcount}
\label{EWD_1_3_01}
\text{\hspace*{-0.8ex}(\theDSBcount)}
\hspace*{4ex}
$\displaystyle{\int\limits_{\mathbb{R}} q_{z}\,dF_{n} - 
\int\limits_{\mathbb{R}} q_{z}\,\psi\,dx = T(q_{z}) + R(q_{z})}$\\[2.5ex]
for each $z \in \mathbb{R}$. For the two right-hand summands, we can then show that\\[2.5ex]
\hspace*{12.1ex}$\displaystyle{\Big|\,T(q_{z}) - 
\dfrac{\mu_{3}}{6\,\sqrt{n}} \int\limits_{\mathbb{R}} q_{z}(x)\,(3\,x - x^3)\,\psi(x)\,dx\Big| 
\leq \dfrac{c_{1}}{n}}$\hspace*{5ex}(see section \ref{EWD_Kap1_Sec4})\\[1.5ex] 
and\\[1.5ex]
\hspace*{12.1ex}$\displaystyle{|R(q_{z})| \leq \dfrac{c_{2}\,\mathcal{C} + c_{3}}{n}}$\hspace*{5ex}(see 
section \ref{EWD_Kap1_Sec5})\\[2ex] 
holds.\\[2.8ex]
In order to be able to formulate the summands $T(q_{z})$ and $R(q_{z})$ precisely, we define\\[2ex]
\refstepcounter{DSBcount}
\label{EWD_1_3_02}
\text{\hspace*{-0.8ex}(\theDSBcount)}
\hspace*{4ex}
$\displaystyle{f_{z}(x) = \psi(x)^{-1}\,\int\limits_{- \infty}^{x}\,\bigl(q_{z}(y) - 
\Phi(q_{z})\,\bigr)\,\psi(y)\,dy}$.\index{function!smooth!$f_{z}$}\\[2ex]
This function is a solution of 
the differential equation (or, in this context, Stein's equation\index{Stein's equation})\\[2ex]
\refstepcounter{DSBcount}
\label{EWD_1_3_03}
\text{\hspace*{-0.8ex}(\theDSBcount)}
\hspace*{4ex}
$f'_{z}(x) - x\,f_{z}(x) = q_{z}(x) - \Phi(q_{z})$\\[2ex]
and also has some very useful analytical properties (cf. chapter \ref{EWD_Kap2}).\\[2.8ex]
The precise form of (\ref{EWD_1_3_01}) is now:\\[4ex]
\refstepcounter{DSBcount}
\label{EWD_1_3_04}
\textbf{\hspace*{-0.8ex}\theDSBcount\ Proposition}\\[0.8ex]
For every $z \in \mathbb{R}$ is valid:\\[2ex]
\hspace*{12.1ex}$\displaystyle{\int\limits_{\mathbb{R}} q_{z}\,dF_{n} - 
\int\limits_{\mathbb{R}} q_{z}\,\psi\,dx = T(q_{z}) + R(q_{z})}$,\\[2ex]
where\\[2ex]
\hspace*{1ex}\hfill\begin{tabular}{@{}l@{\hspace*{0.8ex}}c@{\hspace*{0.8ex}}l@{}}
$\displaystyle{T(q_{z})}$&$=$&
$\displaystyle{\dfrac{\mu_{3}}{2\,\sqrt{n}}\,\Bigl(\,- E\bigl(S_{n}\,f'_{z}(S_{n})\bigr)\,\Bigr)}$,\\[3ex]
$\displaystyle{R(q_{z})}$&$=$&
$\displaystyle{\dfrac{1}{\sqrt{n}}\,E\biggl(\,\int\limits_{0}^{1} \Bigl[\,X_{n} + 
\dfrac{\mu_{3}}{2}\,X_{n}^2 - (1 - t)\,X_{n}^{3}\,\Bigr]\,
\Bigl(f''_{z}\bigl(\,S_{n-1}^{n} + t\,\dfrac{X_{n}}{\sqrt{n}}\,\bigr) - f''_{z}\bigl(\,S_{n-1}^{n}\,\bigr)\Bigl)\,
dt\biggl)}$.
\end{tabular}\\[1.7ex]
Here we write $\displaystyle{S_{n-1}^{n} = \dfrac{1}{\sqrt{n}}\,\sum_{i = 1}^{n-1}\,X_{i}}$.
\index{random variable!$S_{n-1}^{n}$}\\[3.8ex]
\textbf{Proof:}\\[0.8ex]
Throughout the proof, let $z \in \mathbb{R}$ be fixed. 
We can therefore set $q = q_{z}$ and $f = f_{z}$.\\[2.8ex]
Using (\ref{EWD_1_3_03}), we then obtain\\[2ex]
\hspace*{12.1ex}\begin{tabular}{@{}l@{\hspace*{0.8ex}}c@{\hspace*{0.8ex}}l@{}}
$\displaystyle{\int\limits_{\mathbb{R}} q\,dF_{n} - 
\int\limits_{\mathbb{R}} q\,\psi\,dx}$&$=$&
$E\Bigl(q(S_{n}) - \Phi(q)\Bigr)$\\[2ex]
&$=$&$E\Bigl(f'(S_{n}) - S_{n}\,f(S_{n})\Bigr)$\\[2.5ex]
&$=$&$E\Bigl(f'(S_{n})\Bigr) - E\Bigl(\sqrt{n}\,X_{n}\,f(S_{n})\Bigr)$.
\end{tabular}\\[2.5ex]
For the last equation, it was used that the $X_{i}$ are independent and identically distributed.
It should also be noted that the above chain of equations forms the core of Stein's method.\\[2.8ex]
We now consider the terms $E\Bigl(f'(S_{n})\Bigr)$ and $E\Bigl(\sqrt{n}\,X_{n}\,f(S_{n})\Bigr)$ separately.
A Taylor expansion\index{Taylor's theorem} of $f'$ and $f$ 
about $S_{n-1}^{n}$ then gives, if we additionally use $E(X_{n}) = 0$,
$E(X_{n}^2) = 1$ and the independence of $X_{n}$ and $S_{n-1}^{n}$:\\[2.5ex]
\hspace*{3ex}\begin{tabular}{@{}r@{\hspace*{0.8ex}}c@{\hspace*{0.8ex}}l@{}}
$E\Bigl(f'(S_{n})\Bigr)$&$=$&
$\displaystyle{E\biggl(\,f'(S_{n-1}^{n}) + \dfrac{X_{n}}{\sqrt{n}}\,f''(S_{n-1}^{n})
+ \dfrac{X_{n}}{\sqrt{n}}\, \int\limits_{0}^{1}\,\Bigl(\,f''\bigl(S_{n-1}^{n} + t\,\dfrac{X_{n}}{\sqrt{n}}\bigr)
- f''(S_{n-1}^{n})\,\Bigr)\,dt\,\biggr)}$\\[3ex]
&$=$&
$\displaystyle{E\Bigl(\,f'(S_{n-1}^{n})\,\Bigr) + 
E\biggl(\,\dfrac{X_{n}}{\sqrt{n}}\, \int\limits_{0}^{1}\,\Bigl(\,f''\bigl(S_{n-1}^{n} + t\,\dfrac{X_{n}}{\sqrt{n}}\bigr)
- f''(S_{n-1}^{n})\,\Bigr)\,dt\,\biggr)}$,
\end{tabular}\\[3ex]
\begin{tabular}{@{}r@{\hspace*{0.8ex}}c@{\hspace*{0.8ex}}l@{}}
$E\Bigl(\sqrt{n}\,X_{n}\,f(S_{n})\Bigr)$&$=$&
$\displaystyle{E\biggl(\,\sqrt{n}\,X_{n}\,\biggl[\,
f(S_{n-1}^{n}) + \dfrac{X_{n}}{\sqrt{n}}\,f'(S_{n-1}^{n}) + \dfrac{X_{n}^{2}}{2n}\,f''(S_{n-1}^{n})}$\\[3ex] 
&&\hspace{17.1ex}$\displaystyle{
+\ \dfrac{X_{n}^{2}}{n}\, \int\limits_{0}^{1}\,(1 - t)\,\Bigl(\,f''\bigl(S_{n-1}^{n} + t\,\dfrac{X_{n}}{\sqrt{n}}\bigr)
- f''(S_{n-1}^{n})\,\Bigr)\,dt\,\biggr]\,\biggr)}$\\[4.5ex]
&$=$&$\displaystyle{E\Bigl(\,f'(S_{n-1}^{n})\,\Bigr) + 
\dfrac{\mu_{3}}{2\sqrt{n}}\,E\Bigl(\,f''(S_{n-1}^{n})\,\Bigr)}$\\[3ex]
&&\hspace{14.5ex}$\displaystyle{+\ E\biggl(\,
\dfrac{X_{n}^{3}}{\sqrt{n}}\, \int\limits_{0}^{1}\,(1 - t)\,\Bigl(\,f''\bigl(S_{n-1}^{n} 
+ t\,\dfrac{X_{n}}{\sqrt{n}}\bigr)
- f''(S_{n-1}^{n})\,\Bigr)\,dt\,\biggr)}$.
\end{tabular}\\[3ex]
Combining these calculations yields\\[3ex]
\refstepcounter{DSBcount}
\label{EWD_1_3_05}
\text{\hspace*{-0.8ex}(\theDSBcount)}
\hspace*{4ex}
\begin{tabular}[t]{@{}l@{}}
$\displaystyle{\int\limits_{\mathbb{R}} q\,dF_{n} - 
\int\limits_{\mathbb{R}} q\,\psi\,dx =
\dfrac{\mu_{3}}{2\sqrt{n}}\,E\Bigl(\,- f''(S_{n-1}^{n})\,\Bigr)}$
\end{tabular}\\[3.5ex]
\hspace*{4ex}
\begin{tabular}[t]{@{}l@{}}
\hspace*{16ex}$\displaystyle{+\ \dfrac{1}{\sqrt{n}}\,E\biggl(\,
\int\limits_{0}^{1}\,\Bigl[\,X_{n} - (1 - t)\,X_{n}^{3}\,\Bigr]\,
\Bigl(\,f''\bigl(S_{n-1}^{n} 
+ t\,\dfrac{X_{n}}{\sqrt{n}}\bigr)
- f''(S_{n-1}^{n})\,\Bigr)\,dt\,\biggr)}$.
\end{tabular}\\[2.5ex]
The degree of the derivative of $f$ in the term $\dfrac{\mu_{3}}{2\sqrt{n}}\,E\Bigl(\,- f''(S_{n-1}^{n})\,\Bigr)$ is now being reduced again. The following calculation is used for this:\\[2.5ex]
\hspace*{4ex}\begin{tabular}{@{}c@{\hspace*{0.8ex}}l@{}}
&$E\Bigl(S_{n}\,f'(S_{n})\Bigr)$\\[3ex]
$=$&$E\Bigl(\sqrt{n}\,X_{n}\,f'(S_{n})\Bigr)$\\[2ex]
$=$&
$\displaystyle{E\biggl(\,\sqrt{n}\,X_{n}\,\biggl[\,
f'(S_{n-1}^{n}) + \dfrac{X_{n}}{\sqrt{n}}\,f''(S_{n-1}^{n}) +
\dfrac{X_{n}}{\sqrt{n}}\, \int\limits_{0}^{1}\,\Bigl(\,f''\bigl(S_{n-1}^{n} + t\,\dfrac{X_{n}}{\sqrt{n}}\bigr)
- f''(S_{n-1}^{n})\,\Bigr)\,dt\,\biggr]\,\biggr)}$\\[3.5ex]
$=$&
$\displaystyle{E\Bigl(\,f''(S_{n-1}^{n})\,\Bigr) + 
E\biggl(\,X_{n}^{2}\, \int\limits_{0}^{1}\,\Bigl(\,f''\bigl(S_{n-1}^{n} + t\,\dfrac{X_{n}}{\sqrt{n}}\bigr)
- f''(S_{n-1}^{n})\,\Bigr)\,dt\,\biggr)}$.
\end{tabular}\\[2.2ex]
If we now insert the corresponding expression for $E\Bigl(\,f''(S_{n-1}^{n})\,\Bigr)$ into
(\ref{EWD_1_3_05}), we obtain the assertion of the proposition.\hfill$\Box$

\section{Deducing the correct expansion}\label{EWD_Kap1_Sec4}

This section deals with the term
$T(q_{z}) = \dfrac{\mu_{3}}{2\,\sqrt{n}}\,\Bigl(\,- E\bigl(S_{n}\,f'_{z}(S_{n})\bigr)\,\Bigr)$.
We show the following result.\\[4ex]
\refstepcounter{DSBcount}
\label{EWD_1_4_01}
\textbf{\hspace*{-0.8ex}\theDSBcount\ Proposition}\\[2ex]
$\sup\limits_{z \in \mathbb{R}}\,\bigg|\,E\Bigl(\,S_{n}\,f'_{z}(S_{n})\,\Bigr) - \Phi\Bigl(\,x\,f'_{z}(x)\,\Bigr)\,\bigg|
\leq 21,2\,\dfrac{\beta_{3}}{\sqrt{n}}$.\\[4.2ex]
Since\\[1ex]
\hspace*{12.1ex}$\displaystyle{\Phi\bigl(x\,f'_{z}(x)\bigr) = 
\dfrac{1}{3}\,\int\limits_{\mathbb{R}} q_{z}(x)\,(3\,x - x^3)\,\psi(x)\,dx}$
\hspace*{3ex}for all $z \in \mathbb{R}$\\[2.5ex]
(cf. Lemma \ref{EWD_2_2_05}, \ref{EWD_2_2_05_BWc})), we get from 
Proposition \ref{EWD_1_4_01} the estimate:\\[2.5ex]
\refstepcounter{DSBcount}
\label{EWD_1_4_02}
\text{\hspace*{-0.8ex}(\theDSBcount)}
\hspace*{4ex}
$\displaystyle{\sup\limits_{z \in \mathbb{R}}\,\bigg|\,T(q_{z}) + 
\dfrac{\mu_{3}}{6\sqrt{n}}\,\int\limits_{\mathbb{R}} q_{z}(x)\,(3\,x - x^3)\,\psi(x)\,dx\,\bigg|
\leq 10,6\,\dfrac{\beta_{3}^{2}}{n} \leq 10,6\,\dfrac{\beta_{4}}{n}}$.\\[3ex]
For the second inequality, the following application of H{\"o}lder's inequality\index{H{\"o}lder's inequality!for
random variables} is used:\\[2.3ex] 
\hspace*{12.1ex}$\beta_{3} = E\Bigl( |X_{i}|\,|X_{i}|^2\Bigr) \leq 
\Bigl(E(X_{i}^{2})\Bigr)^{1/2}\,\beta_{4}^{1/2} = \beta_{4}^{1/2}$.\\[2.8ex]
We take the most important step in proving Proposition \ref{EWD_1_4_01} by proving the following theorem,
which is interesting in itself. This theorem is also useful in the next section when estimating $R(q_{z})$.\\[4ex]
\refstepcounter{DSBcount}
\label{EWD_1_4_03}
\textbf{\hspace*{-0.8ex}\theDSBcount\ Theorem}\index{Theorem!for the iid-case}\\[2ex]
$\displaystyle{\sup\limits_{z \in \mathbb{R}}\,\bigg|\,
E\Bigl(\,|S_{n}|\, 1_{(- \infty,\, z\,]}(S_{n})\,\Bigr) - \Phi\Bigl(\,|x|\, 1_{(- \infty,\, z\,]}(x)\,\Bigr)\,\bigg|
\leq 10,6\,\dfrac{\beta_{3}}{\sqrt{n}}}$.\\[3.8ex]
\textbf{Proof:}\\[0.8ex]
We proceed as in the proof of Bolthausen's paper (cf. \cite{Bolthausen1984}, pp. 380, 381) using the method of Stein
and a recursion argument. Furthermore, in the following we assume $n \geq  75$ without loss of generality. 
This is possible because we have $E\bigl(|S_{n}|\bigr) \leq 1$, $\Phi\bigl(|x|\bigr) \leq 1$ 
and $\beta_{3} \geq 1$.\\[2.8ex]
For $x, z \in \mathbb{R}$ and $\lambda > 0$ we now define the functions\\[2ex]
\hspace*{12.1ex}\begin{tabular}{@{}l@{\hspace*{0.8ex}}c@{\hspace*{0.8ex}}l@{}}
$p_{z}^{\lambda}(x)$&$=$&
$\displaystyle{|x|\,\Bigl\{\,1_{(- \infty,\, z\,]}(x) + 
\dfrac{z + \lambda - x}{\lambda}\,1_{(\,z,\, z + \lambda\,]}(x)\,\Bigr\} }$,\\[2ex]
$p_{z}^{0}(x)$&$=$&
$\displaystyle{|x|\,1_{(- \infty,\, z\,]}(x)}$.  
\end{tabular}\\[2.5ex]
Furthermore, let\\[2ex]
\hspace*{12.1ex}\begin{tabular}{@{}l@{\hspace*{0.8ex}}c@{\hspace*{0.8ex}}l@{}}
$\mathscr{L}(n, \beta_{3})$&$=$&
$\Bigl\{\,\underline{X} = (X_{1},\ldots,X_{n})$ : $X_{1},\ldots,X_{n}$ are independent and identically\\[0.7ex]
&&\hspace*{1.8ex}distributed with $E(X_{i}) = 0$, $E(X_{i}^{2}) = 1$ and $E(|X_{i}|^{3}) = \beta_{3}\,\Bigr\}$,\\[2ex]
$\delta(\lambda, \beta_{3}, n)$&$=$&
$\sup\,\Bigl\{\,\Big|\,E\bigl(p_{z}^{\lambda}(S_{n})\bigr) - 
\Phi\bigl(p_{z}^{\lambda}\bigr)\,\Big|\,:\,z \in \mathbb{R},\,
\underline{X} \in \mathscr{L}(n, \beta_{3})\,\Bigr\}$,\\[2ex]
$\delta(\beta_{3}, n)$&$=$&$\delta(0, \beta_{3}, n)$.
\end{tabular}\\[2.5ex]
By using $p_{z - \lambda}^{\lambda} \leq p_{z}^{0} \leq p_{z}^{\lambda}$ and
$\sup\limits_{x \in \mathbb{R}}\,|x|\,\psi(x) = \dfrac{1}{\sqrt{2 \,\pi \, e}}$ 
(cf. Lemma \ref{EWD_2_2_08}, \ref{EWD_2_2_08_BWb})) we obtain\\[2.5ex]
\hspace*{12.1ex}\begin{tabular}{@{}r@{\hspace*{0.8ex}}c@{\hspace*{0.8ex}}l@{}}
$E\bigl(p_{z}^{0}(S_{n})\bigr) - 
\Phi\bigl(p_{z}^{0}\bigr)$&$\leq$&$E\bigl(p_{z}^{\lambda}(S_{n})\bigr) - 
\Phi\bigl(p_{z}^{\lambda}\bigr) + \Phi\bigl(p_{z}^{\lambda}\bigr) - \Phi\bigl(p_{z}^{0}\bigr)$\\[2ex]
&$\leq$&$\displaystyle{\delta(\lambda, \beta_{3}, n) + \dfrac{1}{\sqrt{2 \,\pi \, e}}\,
\int\limits_{z}^{z + \lambda}\dfrac{z + \lambda - x}{\lambda}\,dx}$\hspace*{7ex}and\\[4ex]
$- \Bigl(\,E\bigl(p_{z}^{0}(S_{n})\bigr) - 
\Phi\bigl(p_{z}^{0}\bigr)\,\Bigr)$&$\leq$&$\Phi\bigl(p_{z - \lambda}^{\lambda}\bigr) 
- E\bigl(p_{z - \lambda}^{\lambda}(S_{n})\bigr) 
+ \Phi\bigl(p_{z}^{0}\bigr) - \Phi\bigl(p_{z - \lambda}^{\lambda}\bigr)$
\end{tabular}\\[2ex]
\hspace*{12.1ex}\begin{tabular}{@{}r@{\hspace*{0.8ex}}c@{\hspace*{0.8ex}}l@{}}
\hspace*{24.8ex}&$\leq$&$\displaystyle{\delta(\lambda, \beta_{3}, n) + \dfrac{1}{\sqrt{2 \,\pi \, e}}\,
\int\limits_{z - \lambda}^{z} 1 - \dfrac{z - x}{\lambda}\,dx}$.
\end{tabular}\\[1.5ex]
Together this leads to\\[2ex]
\refstepcounter{DSBcount}
\label{EWD_1_4_04}
\text{\hspace*{-0.8ex}(\theDSBcount)}
\hspace*{4ex}
$\delta(\beta_{3}, n)\, \leq\, \delta(\lambda, \beta_{3}, n) + \dfrac{\lambda}{2}\,\dfrac{1}{\sqrt{2 \,\pi \, e}}$.\\[2.5ex]
In order to estimate $\delta(\beta_{3}, n)$, we will first estimate $\delta(\lambda, \beta_{3}, n)$. 
To do this, we consider\\[2ex]
\hspace*{12.1ex}$\displaystyle{\Big|\,E\bigl(p_{z}^{\lambda}(S_{n})\bigr) - 
\Phi\bigl(p_{z}^{\lambda}\bigr)\,\Big|}$\hspace*{4ex}for \textbf{fixed} $\lambda > 0$
\hspace*{1.4ex}and\hspace*{1.4ex}
$z \in \mathbb{R}$.\\[2.5ex]
To simplify the notations we will write $p = p_{z}^{\lambda}$ in the following and define\\[2ex]
\hspace*{12.1ex}$\displaystyle{d(x) = 
\psi(x)^{-1}\,\int\limits_{- \infty}^{x}\,\bigl(p(y) - \Phi(p)\,\bigr)\,\psi(y)\,dy}$.\\[2.2ex]
Then Stein's equation\index{Stein's equation} is\\[2ex] 
\hspace*{12.1ex}$d'(x) - x\,d(x) = p(x) - \Phi(p)$.\\[2.2ex]
This leads to\\[2ex]
\hspace*{12.1ex}\begin{tabular}{@{}c@{\hspace*{0.8ex}}l@{}}
&$\displaystyle{\Big|\,E\bigl(p(S_{n})\bigr) - 
\Phi\bigl(p\bigr)\,\Big|}$\\[2.3ex]
$=$&$\displaystyle{\Big|\,E\bigl(d'(S_{n})\bigr) - E\bigl(S_{n}\,d(S_{n})\bigr)\,\Big|}$\\[2.3ex]
$=$&$\displaystyle{\Big|\,E\bigl(d'(S_{n})\bigr) - \sqrt{n}\,E\bigl(X_{n}\,d(S_{n})\bigr)\,\Big|}$\\[1.7ex]
$\leq$&
$\displaystyle{E\Bigl(\,\Big|\,d'(S_{n}) - d'(S_{n-1}^{n})\,\Big|\,\Bigr) 
+ E\Bigl(\,X_{n}^{2}\,\int\limits_{0}^{1}\,\Big|\,d'\bigl(\,S_{n-1}^{n} + t\,\dfrac{X_{n}}{\sqrt{n}}\,\bigr)
- d'\bigl(\,S_{n-1}^{n}\,\bigr)\,\Big|\,dt\,\Bigr)}$\\[2.5ex]
$=$&$A_{1} + A_{2}$.
\end{tabular}\\[2.5ex]
To obtain the inequality in the above calculations, we used the independence of $X_{n}$ and 
\mbox{\rule[0ex]{0ex}{4.8ex}$\displaystyle{S_{n-1}^{n} = \dfrac{1}{\sqrt{n}}\,\sum_{i = 1}^{n-1}\,X_{i}}$} and
$E(X_{n}) = 0$, $E(X_{n}^{2}) = 1$, as in the last section.\\[2.8ex]
For the estimate of $A_{1}$ we now use\\[2ex]
\hspace*{12.1ex}\begin{tabular}{@{}c@{\hspace*{0.8ex}}l@{}}
&$|\,d'(x+y) - d'(x)\,|$\\[1.7ex]
$\leq$&
$\displaystyle{|y|\,\biggl\{\,2 + \sqrt{\dfrac{2}{\pi}}\,|x| + x^2 + |x\,y| + 
\dfrac{|x|}{\lambda}\,\int\limits_{0}^{1}\,1_{(z,\, z+\lambda\,]}(x + s\,y)\,ds\biggr\}}$
\hspace*{4ex}for $x, y \in \mathbb{R}$
\end{tabular}\\[2.5ex]
(cf. Lemma \ref{EWD_2_3_05}, \ref{EWD_2_3_05_BWc}), second inequality), so that for $A_{1}$ follows\\[2.5ex]
\hspace*{8.6ex}\begin{tabular}{@{}l@{\hspace*{0.8ex}}c@{\hspace*{0.8ex}}l@{}}
$A_{1}$&$\leq$&
$\displaystyle{E\biggl(\,\dfrac{|X_{n}|}{\sqrt{n}}\,\biggl\{\,
2 + \sqrt{\dfrac{2}{\pi}}\,\big|S_{n-1}^{n}\big| + \bigl( S_{n-1}^{n}\bigr)^{2}
+ \big|S_{n-1}^{n}\big|\,\dfrac{|X_{n}|}{\sqrt{n}}}$\\[2.5ex]
&&\hspace*{37.1ex}
$\displaystyle{+\ \dfrac{\big|S_{n-1}^{n}\big|}{\lambda}\,
\int\limits_{0}^{1}\,1_{(z,\, z+\lambda\,]}\bigl(S_{n-1}^{n} + 
s\,\dfrac{X_{n}}{\sqrt{n}}\bigr)\,ds\,\biggr\}\,\biggr)}$
\end{tabular}\\ 
\hspace*{11.3ex}\begin{tabular}{@{}l@{\hspace*{0.8ex}}c@{\hspace*{0.8ex}}l@{}}
&$\leq$&
$\displaystyle{3,9134\,\dfrac{\beta_{3}}{\sqrt{n}} + 
\dfrac{1}{\lambda\,\sqrt{n}}\,\int\limits_{0}^{1}\ \int\limits_{\mathbb{R}}\, |x_{n}|\,
E\Bigl(\,\big|S_{n-1}^{n}\big| \,1_{(z,\, z+\lambda\,]}\bigl(S_{n-1}^{n} + s\,\dfrac{x_{n}}{\sqrt{n}}\bigr)
\,\Bigr)}$\\[-1ex]
&&\hspace*{64.2ex}
$\,PX_{n}^{-1}(dx_{n})\,ds$.
\end{tabular}\\[2.5ex]
Here again, the independence of $X_{n}$ and $S_{n-1}^{n}$, and $E\bigl((S_{n-1}^{n})^2\bigr) \leq 1$,
$\beta_{3} \geq 1$ and $n \geq 75$ are used.\\[2.8ex]
Let us denote for a short moment with $a = \sqrt{\dfrac{n}{n-1}}\,\bigl(z - s\,\dfrac{x_{n}}{\sqrt{n}}\bigr)$ 
and\\
$b = \sqrt{\dfrac{n}{n-1}} \, \lambda \, \leq \, \sqrt{\dfrac{75}{74}} \, \lambda$, 
then\\[2ex]
\refstepcounter{DSBcount}
\label{EWD_1_4_05}
\text{\hspace*{-0.8ex}(\theDSBcount)}
\hspace*{4ex}
\begin{tabular}{@{}c@{\hspace*{0.8ex}}l@{}}
&$\displaystyle{E\Bigl(\,\big|S_{n-1}^{n}\big| \,1_{(z,\, z+\lambda\,]}\bigl(S_{n-1}^{n} + 
s\,\dfrac{x_{n}}{\sqrt{n}}\bigr)\,\Bigr)}$\\[2ex]
$\leq$&$\displaystyle{E\Bigl(\,\big|S_{n-1}\big| \,1_{(a,\, a+b\,]}\bigl(S_{n-1}\bigr)\,\Bigr)}$\\[2.3ex]
$\leq$&$\displaystyle{\Big|\,E\bigl(p_{a+b}^{0}(S_{n-1})\bigr) - 
\Phi\bigl(p_{a+b}^{0}\bigr)\,\Big| + \Big|\,E\bigl(p_{a}^{0}(S_{n-1})\bigr) - 
\Phi\bigl(p_{a}^{0}\bigr)\,\Big|
+ b\,\sup\limits_{z \in \mathbb{R}}\,|x|\,\psi(x)}$\\[2.3ex]
$\leq$&$\displaystyle{2\,\delta(\beta_{3}, n - 1) + 0,2437\,\lambda}$.
\end{tabular}\\[2.5ex]
This leads to\\[2ex]
\hspace*{8.4ex}$\displaystyle{A_{1}
\,\leq\,4,1571\,\dfrac{\beta_{3}}{\sqrt{n}} + 
\dfrac{2}{\lambda}\,\dfrac{\beta_{3}}{\sqrt{n}}\,\delta(\beta_{3}, n - 1)}$.\\[2.5ex]
For the estimate of $A_{2}$, on the other hand, we use\\[2ex]
\hspace*{12.1ex}\begin{tabular}{@{}c@{\hspace*{0.8ex}}l@{}}
&$|\,d'(x+y) - d'(x)\,|$\\[1.7ex]
$\leq$&
$\displaystyle{|y|\,\biggl\{\,2 +  
\dfrac{|x|}{\lambda}\,\int\limits_{0}^{1}\,1_{(z,\, z+\lambda\,]}(x + s\,y)\,ds\biggr\}
+ |x|\,|\,d(x+y) - d(x)\,|}$
\hspace*{4ex}for $x, y \in \mathbb{R}$
\end{tabular}\\[2ex]
(cf. Lemma \ref{EWD_2_3_05}, \ref{EWD_2_3_05_BWc}), first inequality), from which follows\\[2ex]
\hspace*{8.6ex}\begin{tabular}{@{}l@{\hspace*{0.8ex}}c@{\hspace*{0.8ex}}l@{}}
$A_{2}$&$\leq$&
$\displaystyle{E\biggl(\,X_{n}^{2}\,\int\limits_{0}^{1}\,t\,\dfrac{|X_{n}|}{\sqrt{n}}\,\biggl\{\,
2 + \dfrac{\big|S_{n-1}^{n}\big|}{\lambda}\,
\int\limits_{0}^{1}\,1_{(z,\, z+\lambda\,]}\bigl(S_{n-1}^{n} + 
s\,t\,\dfrac{X_{n}}{\sqrt{n}}\bigr)\,ds\,\biggr\}\,dt\,\biggr)}$
\end{tabular}\\[3ex]
\hspace*{8.6ex}\begin{tabular}{@{}l@{\hspace*{0.8ex}}c@{\hspace*{0.8ex}}l@{}}
&&\hspace*{2ex}$\displaystyle{+\
E\biggl(\,1_{\,\{|X_{n}|\,\leq\,1,4477\,\sqrt{n}\,\}}\,X_{n}^{2}\,\big|S_{n-1}^{n}\big|\,
\int\limits_{0}^{1}\,\Big|\,d\bigl(\,S_{n-1}^{n} + t\,\dfrac{X_{n}}{\sqrt{n}}\,\bigr)
- d\bigl(\,S_{n-1}^{n}\,\bigr)\,\Big|\,dt\,\biggr)}$\\[3ex]
&&\hspace*{2ex}$\displaystyle{+\
E\biggl(\,1_{\,\{|X_{n}|\,>\,1,4477\,\sqrt{n}\,\}}\,X_{n}^{2}\,\big|S_{n-1}^{n}\big|\,
\int\limits_{0}^{1}\,\Big|\,d\bigl(\,S_{n-1}^{n} + t\,\dfrac{X_{n}}{\sqrt{n}}\,\bigr)
- d\bigl(\,S_{n-1}^{n}\,\bigr)\,\Big|\,dt\,\biggr)}$\\[3ex]
\hspace*{2.7ex}&$=$&$A_{3} + A_{4} + A_{5}$.
\end{tabular}\\[2.5ex]
Using Fubini's theorem\index{integration!using Fubini's theorem} and 
a calculation analogous to (\ref{EWD_1_4_05}), we further obtain\\[2ex]
\hspace*{8.6ex}\begin{tabular}{@{}l@{\hspace*{0.8ex}}c@{\hspace*{0.8ex}}l@{}}
$A_{3}$&
$=$&$\displaystyle{\int\limits_{0}^{1}\,t\,\int\limits_{0}^{1}\,\int\limits_{\mathbb{R}}
\,\dfrac{|x_{n}|^{3}}{\sqrt{n}}\,E\biggl(\,2 + \dfrac{\big|S_{n-1}^{n}\big|}{\lambda}\,
1_{(z,\, z+\lambda\,]}\bigl(S_{n-1}^{n} + 
s\,t\,\dfrac{x_{n}}{\sqrt{n}}\bigr)\,\biggr)\,PX_{n}^{-1}(dx_{n})\,ds\,dt}$\\[4ex]
&$\leq$&
$\displaystyle{\int\limits_{0}^{1}\,t\,\int\limits_{0}^{1}\,\int\limits_{\mathbb{R}}
\,\dfrac{|x_{n}|^{3}}{\sqrt{n}}\,\biggl(\,2 + \dfrac{1}{\lambda}\,
\Bigl(\,2\,\delta(\beta_{3}, n - 1) + 0,2437\,\lambda\,\Bigr)\,\biggr)\,PX_{n}^{-1}(dx_{n})\,ds\,dt}$\\[4ex]
&$=$&$\displaystyle{\int\limits_{0}^{1}\,t\,dt\,
\biggl(\,2,2437\,\dfrac{\beta_{3}}{\sqrt{n}} + 
\dfrac{2}{\lambda}\,\dfrac{\beta_{3}}{\sqrt{n}}\,\delta(\beta_{3}, n - 1)\,\biggr)}$
\end{tabular}\\
\hspace*{11.3ex}\begin{tabular}{@{}l@{\hspace*{0.8ex}}c@{\hspace*{0.8ex}}l@{}}
&$\leq$&$1,1219\,\dfrac{\beta_{3}}{\sqrt{n}} + 
\dfrac{1}{\lambda}\,\dfrac{\beta_{3}}{\sqrt{n}}\,\delta(\beta_{3}, n - 1)$.
\end{tabular}\\[2.5ex]
For the estimate of $A_{4}$ we use\\[2ex]
\hspace*{12.1ex}$\displaystyle{|\,d(x+y) - d(x)\,| \,\leq\,
|y|\,\sup\limits_{0 < \theta < 1}\, |d'(x + \theta\,y)| \,\leq\,
|y|\,\Bigl\{\,|x| + |y| + \sqrt{\dfrac{2}{\pi}}\,\Bigr\}}$\\[2.5ex]
(cf. the estimate of $f_{1}'$ in Corollary \ref{EWD_2_1_12}, \ref{EWD_2_1_12_BWb})). We then get\\[2ex]
\hspace*{8.6ex}\begin{tabular}{@{}l@{\hspace*{0.8ex}}c@{\hspace*{0.8ex}}l@{}}
$A_{4}$&$\leq$&
$\displaystyle{E\Bigl(\,1_{\,\{|X_{n}|\,\leq\,1,4477\,\sqrt{n}\,\}}\,\dfrac{|X_{n}|^{3}}{\sqrt{n}}
\,\big|S_{n-1}^{n}\big|\,
\int\limits_{0}^{1}\,t\,\Bigl\{\,\big|S_{n-1}^{n}\big| + t\,\dfrac{|X_{n}|}{\sqrt{n}} + 
\sqrt{\dfrac{2}{\pi}}\,\Bigr\}\,dt\,\Bigr)}$\\[3.5ex]
&$\leq$&$\displaystyle{E\Bigl(\,1_{\,\{|X_{n}|\,\leq\,1,4477\,\sqrt{n}\,\}}\,\dfrac{|X_{n}|^{3}}{\sqrt{n}}
\,\big|S_{n-1}^{n}\big|\,\Bigl\{\,\dfrac{1}{2}\,\big|S_{n-1}^{n}\big| + \dfrac{1}{3}\,1,4477 + 
\dfrac{1}{2}\,\sqrt{\dfrac{2}{\pi}}\,\Bigr\}\,\Bigr)}$\\[3.5ex]
&$\leq$&$\displaystyle{E\Bigl(\,1_{\,\{|X_{n}|\,\leq\,1,4477\,\sqrt{n}\,\}}\,\dfrac{|X_{n}|^{3}}{\sqrt{n}}\,\Bigr)\,
E\Bigl(\,\big|S_{n-1}^{n}\big|\,\Bigl\{\,\dfrac{1}{2}\,\big|S_{n-1}^{n}\big| + 0,8816\,\Bigr\}\,\Bigr)}$\\[3.5ex]
&$\leq$&$\displaystyle{1,3816\,E\Bigl(\,1_{\,\{|X_{n}|\,\leq\,1,4477\,\sqrt{n}\,\}}\,
\dfrac{|X_{n}|^{3}}{\sqrt{n}}\,\Bigr)}$
\hspace*{24ex}$\Bigl(\,E\bigl((S_{n-1}^{n})^2\bigr) \leq 1\,\Bigr)$.
\end{tabular}\\[2.8ex]
For the estimate of $A_{5}$, on the other hand, we use\\[2ex]
\hspace*{12.1ex}$\displaystyle{|\,d(x+y) - d(x)\,| \,\leq\, |\,d(x+y)\,| + |\,d(x)\,|\, \leq\, 2}$\\[2.5ex]
(cf. the estimate of $f_{1}$ in Corollary \ref{EWD_2_1_12}, \ref{EWD_2_1_12_BWb})). The result is\\[2ex]
\hspace*{8.6ex}\begin{tabular}{@{}l@{\hspace*{0.8ex}}c@{\hspace*{0.8ex}}l@{}}
$A_{5}$&$\leq$&$\displaystyle{2\,E\Bigl(\,1_{\,\{|X_{n}|\,>\,1,4477\,\sqrt{n}\,\}}\,X_{n}^{2}\,\Bigr)\,
E\Bigl(\,\big|S_{n-1}^{n}\big|\,\Bigr)}$\\[3.5ex]
&$\leq$&$\displaystyle{1,3816\,E\Bigl(\,1_{\,\{|X_{n}|\,>\,1,4477\,\sqrt{n}\,\}}\,\dfrac{|X_{n}|^{3}}{\sqrt{n}}\,\Bigr)}$
\hspace*{22.7ex}$\Bigl(\,\dfrac{2}{1,4477} \leq 1,3816\,\Bigr)$.
\end{tabular}\\[2.5ex]
If we summarize the results for $A_{3}$, $A_{4}$ and $A_{5}$, we obtain\\[2ex]
\hspace*{8.6ex}\begin{tabular}{@{}l@{\hspace*{0.8ex}}c@{\hspace*{0.8ex}}l@{}}
$A_{2}$&$\leq$&$\displaystyle{1,1219\,\dfrac{\beta_{3}}{\sqrt{n}} + 
\dfrac{1}{\lambda}\,\dfrac{\beta_{3}}{\sqrt{n}}\,\delta(\beta_{3}, n - 1)
+ 1,3816\,\dfrac{\beta_{3}}{\sqrt{n}}}$\\[3.5ex]
&$=$&$\displaystyle{2,5035\,\dfrac{\beta_{3}}{\sqrt{n}} + 
\dfrac{1}{\lambda}\,\dfrac{\beta_{3}}{\sqrt{n}}\,\delta(\beta_{3}, n - 1)}$.
\end{tabular}\\[3ex]
Together with the estimate of $A_{1}$ and using (\ref{EWD_1_4_04}), we get for all $n \geq 75$\\[2.5ex]
\hspace*{12.1ex}$\displaystyle{
\delta(\beta_{3}, n) \,\leq\, 6,6606\,\dfrac{\beta_{3}}{\sqrt{n}} + 
\dfrac{3}{\lambda}\,\dfrac{\beta_{3}}{\sqrt{n}}\,\delta(\beta_{3}, n - 1) + 
\dfrac{\lambda}{2}\,\dfrac{1}{\sqrt{2 \,\pi \, e}}}$.\\[3ex]
If we now choose $\lambda = 16,26\,\dfrac{\beta_{3}}{\sqrt{n}}$, we obtain for all $n \geq 75$\\[2.5ex]
\hspace*{12.1ex}$\displaystyle{
\delta(\beta_{3}, n) \,\leq\, 8,6279\,\dfrac{\beta_{3}}{\sqrt{n}} + 0,1846\,\delta(\beta_{3}, n - 1)}$\\[2.5ex]
and thus\\[2ex]
\hspace*{12.1ex}$\displaystyle{
\delta(\beta_{3}, n) \,\leq\, 8,6279\,\beta_{3}\,\sum\limits_{i = 0}^{n - 75}\, \dfrac{0,1846^{i}}{\sqrt{n-i}}
+ 0,1846^{n - 74}\,\delta(\beta_{3}, 74)}$.\\[3ex]
For further considerations, we note that for all $n \geq 75$ we have\\[2.5ex]
\hspace*{12.1ex}$\displaystyle{
\dfrac{n}{n-i} = \dfrac{n}{n-1} \,\, \dfrac{n-1}{n-2} \,\,\ldots\,\, \dfrac{n-i+1}{n-i} \leq 
\Bigl(\,\dfrac{75}{74}\,\Bigr)^{i}}$
\hspace*{4ex} for $0 \leq i \leq n - 74$.\\[3ex]
Moreover, due to $\delta(\beta_{3}, 74) \leq 1$, $\beta_{3} \geq 1$ and $\sqrt{74} \leq 8,6279$, 
we also have for all $n \geq 75$\\[2.5ex]
\hspace*{12.1ex}$\displaystyle{
0,1846^{n - 74}\,\delta(\beta_{3}, 74) \leq 0,1846^{n - 74}\,1 \leq
0,1846^{n - 74}\,8,6279\,\dfrac{\beta_{3}}{\sqrt{n - (n - 74)}}}$.\\[3ex]
So all in all, using the summation formula for the geometric series\index{geometric series} leads to\\[2.5ex]
\hspace*{12.1ex}\begin{tabular}{@{}l@{\hspace*{0.8ex}}c@{\hspace*{0.8ex}}l@{}}
$\displaystyle{\delta(\beta_{3}, n)}$&$\leq$&$\displaystyle{8,6279\,\dfrac{\beta_{3}}{\sqrt{n}}
\,\sum\limits_{i = 0}^{n - 74}\, \biggl(\,0,1846\,\sqrt{\dfrac{75}{74}}\,\biggr)^{i}}$
\end{tabular}\\[3.5ex] 
\hspace*{12.1ex}\begin{tabular}{@{}l@{\hspace*{0.8ex}}c@{\hspace*{0.8ex}}l@{}}
\hspace*{7.7ex}&$\leq$&$\displaystyle{8,6279\,\dfrac{\beta_{3}}{\sqrt{n}}
\,\,\sum\limits_{i = 0}^{\infty}\, 0,1859^{i}
\,\leq\, 10,5981\,\dfrac{\beta_{3}}{\sqrt{n}}}$.
\end{tabular}\\[3ex]
This proves the theorem.\hfill$\Box$\\[4ex]
\refstepcounter{DSBcount}
\label{EWD_1_4_06}
\textbf{\hspace*{-0.8ex}\theDSBcount\ Remarks}
\begin{enumerate}
\item
The above proof shows that for Theorem \ref{EWD_1_4_03} and thus all results of this section
only $\beta_{3} < \infty$ instead of $\beta_{4} < \infty$ has to be assumed.
\item
The estimate\\[2ex]
\hspace*{6.3ex}$\displaystyle{\sup\limits_{z \in \mathbb{R}}\,\bigg|\,
E\Bigl(\,|S_{n}|\, 1_{(- \infty,\, z\,]}(S_{n})\,\Bigr) - \Phi\Bigl(\,|x|\, 1_{(- \infty,\, z\,]}(x)\,\Bigr)\,\bigg|
\leq C\,\dfrac{\beta_{3}}{\sqrt{n}}}$,\\[2.5ex]
where $C > 0$ is an absolute constant, can also be deduced from a more general result of 
Sweeting \cite{10.1214/aop/1176995888} (see Theorem 1, page 39 and the following example).
In contrast to here, however, Sweeting uses the operator method, which goes back to Trotter \cite{Trotter1959}.
\end{enumerate}
\vspace*{3ex}
The Theorem \ref{EWD_1_4_03} also holds 
without the vertical bars for the absolute values of $S_{n}$ and $x$.\\[4ex]
\refstepcounter{DSBcount}
\label{EWD_1_4_07}
\textbf{\hspace*{-0.8ex}\theDSBcount\ Corollary}\\[2ex]
$\displaystyle{\sup\limits_{z \in \mathbb{R}}\,\bigg|\,
E\Bigl(\,S_{n}\, 1_{(- \infty,\, z\,]}(S_{n})\,\Bigr) - \Phi\Bigl(\,x\, 1_{(- \infty,\, z\,]}(x)\,\Bigr)\,\bigg|
\leq 10,6\,\dfrac{\beta_{3}}{\sqrt{n}}}$.\\[3.8ex]
\textbf{Proof:}\\[0.8ex]
For $z \leq 0$ the above estimate follows directly from Theorem \ref{EWD_1_4_03}.
For $z \geq 0$, on the other hand, we obtain\\[2ex]
\hspace*{12.1ex}\begin{tabular}{@{}c@{\hspace*{0.8ex}}l@{\hspace*{-2.8ex}}r@{}}
&$\displaystyle{\bigg|\,
E\Bigl(\,S_{n}\, 1_{(- \infty,\, z\,]}(S_{n})\,\Bigr) - \Phi\Bigl(\,x\, 1_{(- \infty,\, z\,]}(x)\,\Bigr)\,\bigg|}$\\[2ex]
$=$&$\displaystyle{\bigg|\,
E\Bigl(\,S_{n}\, 1_{(\,z,\,\infty\,)}(S_{n})\,\Bigr) - \Phi\Bigl(\,x\, 1_{(\,z,\,\infty\,)}(x)\,\Bigr)\,\bigg|}$
&$\Bigr($\,since $E(S_{n}) = 0 = \Phi(id_{\mathbb{R}})$\,$\Bigl)$\\[2ex]
$=$&$\displaystyle{\bigg|\,
E\Bigl(\,|S_{n}|\, 1_{(\,z,\,\infty\,)}(S_{n})\,\Bigr) - \Phi\Bigl(\,|x|\, 1_{(\,z,\,\infty\,)}(x)\,\Bigr)\,\bigg|}$
&$\Bigr($\,since $z \geq 0$\,$\Bigl)$\\[2ex]
$=$&$\displaystyle{\bigg|\,
E\Bigl(\,| - S_{n}|\, 1_{(\,- \infty,\,- z\,)}( - S_{n})\,\Bigr) - 
\Phi\Bigl(\,|x|\, 1_{(\,- \infty,\,- z\,)}(x)\,\Bigr)\,\bigg|}$\\[2ex]
$\leq$&$\displaystyle{10,6\,\dfrac{\beta_{3}}{\sqrt{n}}}$.
\end{tabular}\\[2.5ex] 
The inequality is valid, because with $X_{1},\ldots,X_{n}$ also $-X_{1},\ldots,-X_{n}$
fulfil the conditions of Theorem \ref{EWD_1_4_03}.\hfill$\Box$\\[2.8ex]
For the proof of Proposition \ref{EWD_1_4_01} we need a slightly more general version of the last corollary.
This reads as follows:\\[4ex]
\refstepcounter{DSBcount}
\label{EWD_1_4_08}
\textbf{\hspace*{-0.8ex}\theDSBcount\ Corollary}\\[0.8ex]
Let $\mathcal{H}$ be defined 
according to (\ref{EWD_0_1_01})\index{function!$\mathcal{H}$}. Then\\[2ex]
$\displaystyle{\sup\limits_{h \in \mathcal{H}}\,\,\bigg|\,
E\Bigl(\,S_{n}\, h(S_{n})\,\Bigr) - \Phi\Bigl(\,x\, h(x)\,\Bigr)\,\bigg|
\leq 10,6\,\dfrac{\beta_{3}}{\sqrt{n}}}$.\\[3.8ex]
\textbf{Proof:}\\[0.8ex]
According to Corollary \ref{EWD_1_4_07},\\[2ex]
\refstepcounter{DSBcount}
\label{EWD_1_4_09}
\text{\hspace*{-0.8ex}(\theDSBcount)}
\hspace*{4ex}
$\displaystyle{\bigg|\,E\Bigl(\,S_{n}\, h(S_{n})\,\Bigr) - \Phi\Bigl(\,x\, h(x)\,\Bigr)\,\bigg|
\leq 10,6\,\dfrac{\beta_{3}}{\sqrt{n}}}$\\[2.5ex]
is valid for the functions $h = 1_{(- \infty,\, z\,]}$, $z \in \mathbb{R}$. 
From this we can derive (\ref{EWD_1_4_09}) for convex combinations\index{function!convex combinations} 
of these functions and then, using 
majorized convergence\index{integration!using majorized convergence}, 
for all $h \in \mathcal{H}$.\\[1.5ex] 
To apply the theorem of majorized convergence, we note that for all $h \in \mathcal{H}$ the inequalities 
$\big|\,S_{n}\,h(S_{n})\,\big| \leq \big|\,S_{n}\,\big|$
and $\big|\,x\,h(x)\,\psi(x)\,\big| \leq \big|\,x\,\psi(x)\,\big|$ hold
due to $0 \leq h \leq 1$.
\hspace*{1ex}\hfill$\Box$\\[4ex]
We now obtain Proposition \ref{EWD_1_4_01}, since $f'_{z}(x) = h_{1}(x) - h_{2}(x)$ with\\[2ex]
\hspace*{12.1ex}$h_{1}(x) = q_{z}(x) \in \mathcal{H}$
\hspace*{2ex}and\hspace*{2ex}
$h_{2}(x) = -\,\,\bigl(\,x\,f_{z}(x) - \Phi(q_{z})\,\bigr) \in \mathcal{H}$\\[2.5ex]
(cf. Lemma \ref{EWD_2_1_13}, \ref{EWD_2_1_13_BWa}) and (\ref{EWD_1_3_03})).

\section[The estimate of $R(q_{z})$]{The estimate of {\boldmath $R(q_{z})$}}\label{EWD_Kap1_Sec5}

Finally, in this section, we show the following estimate:\\[4ex]
\refstepcounter{DSBcount}
\label{EWD_1_5_01}
\textbf{\hspace*{-0.8ex}\theDSBcount\ Proposition}\\[2ex]
$\displaystyle{\sup\limits_{z \in \mathbb{R}}\,\Big|\,R(q_{z})\,\Big|\,\leq\,
\dfrac{1}{n}\,\Bigl(\,\dfrac{1}{2} + \dfrac{5}{12}\,\beta_{4}\,\Bigr)\,
\Bigl(\,\dfrac{2\,n}{n-1}\,\mathcal{C} + 5,2 + 21,2\,\beta_{3}\,\Bigr)}$.\\[3.8ex]
\textbf{Proof:}\\[0.8ex]
For $x, y, z \in \mathbb{R}$ we use\\[2ex]
\hspace*{12.1ex}\begin{tabular}{@{}c@{\hspace*{0.8ex}}l@{}}
&$\displaystyle{\Big|\,\bigl(\,f_{z}'' - q_{z}'\,\bigr)(x+y) - \bigl(\,f_{z}'' - q_{z}'\,\bigr)(x)\,\Big|}$\\[2ex]
$\leq$&$\displaystyle{|y|\,\biggl\{
3 + \dfrac{\sqrt{2\pi}}{4}\,|x| + x^2 + \sqrt{n}\,|x|\,\int\limits_{0}^{1}\,1_{\bigl(\,z,\,z\, 
+ \textstyle{\frac{2}{\sqrt{n}}}\,\bigr]}(x + s\,y)\,ds\,\biggr\}}$
\end{tabular}\\[2.2ex]
(cf. Lemma \ref{EWD_2_3_06}, \ref{EWD_2_3_06_BWd}) with $\lambda = \dfrac{1}{\sqrt{n}}\,$).
Thus we obtain for all $z \in \mathbb{R}$\\[2ex]
\begin{tabular}{@{\hspace*{2ex}}l@{\hspace*{0.8ex}}c@{\hspace*{0.8ex}}l@{}}
$\displaystyle{\Big|\,R(q_{z})\,\Big|}$&$=$&
$\displaystyle{\dfrac{1}{\sqrt{n}}\,\bigg|\,E\biggl(\,\int\limits_{0}^{1} \Bigl[\,X_{n} + 
\dfrac{\mu_{3}}{2}\,X_{n}^2 - (1 - t)\,X_{n}^{3}\,\Bigr]\,
\Bigl(f''_{z}\bigl(\,S_{n-1}^{n} + t\,\dfrac{X_{n}}{\sqrt{n}}\,\bigr) - f''_{z}\bigl(\,S_{n-1}^{n}\,\bigr)\Bigl)\,
dt\biggl)\,\bigg|}$\\[4.5ex]
&$\leq$&$\dfrac{1}{n}\,\bigl(\,A_{1} + A_{2} + A_{3}\,\bigr)$,
\end{tabular}\\[2ex]
where\\[2ex]
\begin{tabular}{@{\hspace*{2ex}}l@{\hspace*{0.8ex}}c@{\hspace*{0.8ex}}l@{}}
$A_{1}$&$=$&$\displaystyle{E\biggl(\,\int\limits_{0}^{1}\,
\Bigl[\,t\,X_{n}^{2} + 
t\,\dfrac{\beta_{3}}{2}\,|X_{n}|^3 + (t-t^2)\,X_{n}^{4}\,\Bigr]\,dt\,
\Bigl\{\,3 + \dfrac{\sqrt{2\pi}}{4}\,|S_{n-1}^{n}| + (S_{n-1}^{n})^2\,\Bigr\} 
\,\biggr)}$,\\[4ex]
$A_{2}$&$=$&$\displaystyle{\int\limits_{0}^{1}\,\int\limits_{0}^{1}\,
E\biggl(\,\Bigl[\,t\,X_{n}^{2} + 
t\,\dfrac{\beta_{3}}{2}\,|X_{n}|^3 + (t-t^2)\,X_{n}^{4}\,\Bigr]}$\\[2ex]
&&\hspace*{16.4ex}$\displaystyle{\sqrt{n}\,|S_{n-1}^{n}|\,1_{\bigl(\,z,\,z\, 
+ \textstyle{\frac{2}{\sqrt{n}}}\,\bigr]}\Bigl(\,S_{n-1}^{n} + 
s\,t\,\dfrac{X_{n}}{\sqrt{n}}\,\Bigr)\,\biggr)\,ds\,dt}$,\\[3ex]
$A_{3}$&$=$&$\displaystyle{\int\limits_{0}^{1}\,\int\limits_{0}^{1}\,
\bigg|\,E\biggl(\,\Bigl[\,t\,X_{n}^{2} + 
t\,\dfrac{\mu_{3}}{2}\,X_{n}^{3} - (t - t^2)\,X_{n}^{4}\,\Bigr]}$\\[2ex]
&&\hspace*{16.4ex}$\displaystyle{n\,\Bigl(\,1_{\bigl(\,z\, 
+ \textstyle{\frac{1}{\sqrt{n}}},\,\scriptstyle{z}\, 
+ \textstyle{\frac{2}{\sqrt{n}}}\,\bigr]}\,-\,1_{\bigl(\,z,\,z\, 
+ \textstyle{\frac{1}{\sqrt{n}}}\,\bigr]}\,\Bigr)\,\Bigl(\,S_{n-1}^{n} + 
s\,t\,\dfrac{X_{n}}{\sqrt{n}}\,\Bigr)\,\biggr)\,\bigg|\,ds\,dt}$
\end{tabular}\\[3ex]
(for $A_{3}$ see also Lemma \ref{EWD_2_3_06}, \ref{EWD_2_3_06_BWb})).\\[2.8ex]
Due to the independence of $X_{n}$ and $S_{n-1}^{n}$, and due to $E(X_{n}^{2}) = 1$, 
$E\bigl( (S_{n-1}^{n})^2\bigr) \leq 1$ and $\beta_{3}^{2} \leq \beta_{4}$ then follows\\[2ex]
\begin{tabular}{@{\hspace*{2ex}}l@{\hspace*{0.8ex}}c@{\hspace*{0.8ex}}l@{}}
$A_{1}$&$=$&$\displaystyle{
E\biggl(\,\dfrac{1}{2}\,X_{n}^{2} + \dfrac{1}{2}\,\dfrac{\beta_{3}}{2}\,|X_{n}|^3 + 
\dfrac{1}{6}\,X_{n}^{4}\,\biggr)\,E\biggl(\,3 + \dfrac{\sqrt{2\pi}}{4}\,|S_{n-1}^{n}| + (S_{n-1}^{n})^2\,\biggr)}$\\[4ex]
&$\leq$&$\displaystyle{\Bigl(\,4 + \dfrac{\sqrt{2\pi}}{4}\,\Bigr)\,
\Bigl(\,\dfrac{1}{2} + \dfrac{5}{12}\,\beta_{4}\,\Bigr)}$.
\end{tabular}\\[3ex]
Furthermore, using Fubini's theorem\index{integration!using Fubini's theorem} and
Theorem \ref{EWD_1_4_03} and Lemma \ref{EWD_2_2_08}, \ref{EWD_2_2_08_BWb}), we obtain\\[2ex]
\begin{tabular}{@{\hspace*{2ex}}l@{\hspace*{0.8ex}}c@{\hspace*{0.8ex}}l@{}}
$A_{2}$&$=$&$\displaystyle{\int\limits_{0}^{1}\,\int\limits_{0}^{1}\,\int\limits_{\mathbb{R}}
\,\Bigl[\,t\,x_{n}^{2} + 
t\,\dfrac{\beta_{3}}{2}\,|x_{n}|^3 + (t-t^2)\,x_{n}^{4}\,\Bigr]}$\\[2.5ex]
&&\hspace*{17.8ex}$\displaystyle{\sqrt{n}\,E\biggl(\,|S_{n-1}^{n}|\,1_{\bigl(\,z,\,z\, 
+ \textstyle{\frac{2}{\sqrt{n}}}\,\bigr]}\Bigl(\,S_{n-1}^{n} + 
s\,t\,\dfrac{x_{n}}{\sqrt{n}}\,\Bigr)\,\biggr)\,PX_{n}^{-1}(dx_{n})\,ds\,dt}$\\[4ex]
&$\leq$&$\Bigr(\,21,2\,\beta_{3} + \dfrac{2}{\sqrt{2\,\pi\,e}}\,\Bigl)\,
\Bigl(\,\dfrac{1}{2} + \dfrac{5}{12}\,\beta_{4}\,\Bigr)$.
\end{tabular}\\[3ex]
Finally, by applying again Fubini's theorem\index{integration!using Fubini's theorem} and then
condition (\ref{EWD_1_1_10}), we get\\[2ex]
\begin{tabular}{@{\hspace*{2ex}}l@{\hspace*{0.8ex}}c@{\hspace*{0.8ex}}l@{}}
$A_{3}$&$\leq$&$\displaystyle{\int\limits_{0}^{1}\,\int\limits_{0}^{1}\,\int\limits_{\mathbb{R}}
\,\Bigl[\,t\,x_{n}^{2} + 
t\,\dfrac{\beta_{3}}{2}\,|x_{n}|^3 + (t-t^2)\,x_{n}^{4}\,\Bigr]}$\\[5ex]
&&\hspace*{1.1ex}$\displaystyle{n\,\bigg|\,E\biggl(\,\Bigl(\,1_{\bigl(\,z\, 
+ \textstyle{\frac{1}{\sqrt{n}}},\,\scriptstyle{z}\, 
+ \textstyle{\frac{2}{\sqrt{n}}}\,\bigr]}\,-\,1_{\bigl(\,z,\,z\, 
+ \textstyle{\frac{1}{\sqrt{n}}}\,\bigr]}\,\Bigr)\,\Bigl(\,S_{n-1}^{n} + 
s\,t\,\dfrac{x_{n}}{\sqrt{n}}\,\Bigr)\,\bigg|\,PX_{n}^{-1}(dx_{n})\,ds\,dt}$
\end{tabular}\\[4ex]
\begin{tabular}{@{\hspace*{2ex}}l@{\hspace*{0.8ex}}c@{\hspace*{0.8ex}}l@{}}
\hspace*{2.8ex}&$\leq$&$\dfrac{2\,n}{n-1}\,\mathcal{C}\,
\Bigl(\,\dfrac{1}{2} + \dfrac{5}{12}\,\beta_{4}\,\Bigr)$.
\end{tabular}\\[-2.4ex]
\hspace*{1ex}\hfill$\Box$\\[4ex]
Overall, summarizing the results of this chapter 
(cf. Proposition \ref{EWD_1_2_02}, Proposition \ref{EWD_1_3_04}, (\ref{EWD_1_4_02}) 
and Proposition \ref{EWD_1_5_01}) gives\\[2ex]
\hspace*{12.1ex}\begin{tabular}{@{}l@{\hspace*{0.8ex}}c@{\hspace*{0.8ex}}l@{}}
$||F_{n} - e_{n}||$&$\leq$&$\displaystyle{\dfrac{1}{n}\,\biggl\{\Bigl(\,\dfrac{7}{12}\,\mathcal{C} + \dfrac{1}{48} + 
\dfrac{1}{60}\,\dfrac{\beta_{3}}{\sqrt{n}}\,\Bigr)\,+\,10,6\,\beta_{4}}$\\[2.5ex]
&&\hspace*{16.6ex}$\displaystyle{+\ \Bigl(\,\dfrac{1}{2} + \dfrac{5}{12}\,\beta_{4}\,\Bigr)\,
\Bigl(\,\dfrac{2\,n}{n-1}\,\mathcal{C}
\,+\,5,2\,+\,21,2\,\beta_{3}\,\Bigr)\,\biggr\}}$.
\end{tabular}\\[2.5ex]
For $n \geq 6$ $\Bigl(\,\Leftrightarrow \dfrac{2\,n}{n-1} \leq \dfrac{12}{5}\,\Bigr)$ this leads to\\[2ex]
\refstepcounter{DSBcount}
\label{EWD_1_5_02}
\text{\hspace*{-0.8ex}(\theDSBcount)}
\hspace*{4ex}
$\displaystyle{||F_{n} - e_{n}|| \, \leq \, \dfrac{1}{n}\,\Bigl\{\,\bigl(\,2 + \beta_{4}\,\bigr)\,\mathcal{C} 
+ \bigl(\,3 + 11\,\beta_{3} + 13\,\beta_{4} + 9\,\beta_{3}\,\beta_{4}\,\bigr)\,\Bigr\}}$.\\[2.8ex]
Because of (\ref{EWD_1_1_08}) this completes the proof of Theorem \ref{EWD_1_1_03}.                                
\rehead{Some analytical tools}    
\chapter[Some analytical tools]{Some analytical tools}\label{EWD_Kap2}

\section[Elementary properties of the solution $f_{k}$ of Stein's equation for the
function $|x|^{k}\,q(x)$]
{Elementary properties of the solution {\boldmath $f_{k}$} of Stein's equation for the
function {\boldmath $\displaystyle{|x|^{k}\,q(x)}$}}\label{EWD_Kap2_Sec1}

As we have seen in chapter \ref{EWD_Kap1}, some properties of the solution $f_{k}$ of Stein's equation for the function $|x|^{k}\,q(x)$ are required when applying Stein's method. These are provided in this and the
following sections.\\[2.8ex]
First of all, let $q\,:\,\mathbb{R} \rightarrow \mathbb{R}$ be measurable and bounded. For 
$k \in \mathbb{N}_{0}$ and $x \in \mathbb{R}$ we \textbf{define} \\[2ex]
\hspace*{12.1ex}$\displaystyle{f_{k}(x) = \psi(x)^{-1}\,\int\limits_{-\,\infty}^{x}
\Bigl(\,|y|^{k}\,q(y) - \Phi\bigl(\,|v|^{k}\,q(v)\,\bigr)\,\Bigr)\, \psi(y)\,dy}$,\\[2ex]
\refstepcounter{DSBcount}
\label{EWD_2_1_01}
\text{\hspace*{-0.8ex}(\theDSBcount)}
\hspace*{4ex}
$f_{k}'(x) = x\,f_{k}(x) + |x|^{k}\,q(x) - \Phi\bigl(\,|v|^{k}\,q(v)\,\bigr)$
\hfill(Stein's equation\index{Stein's equation}),\\[3.8ex]
\hspace*{12.1ex}$f = f_{0}$, $f' = f_{0}'$.\\[2.8ex]
In the case of a continuous function $q$ (see e.g. $p_{z,0}^{\lambda}$, $q_{z}^{\lambda,0}$ 
and $r_{z}^{\lambda}$ in section \ref{EWD_Kap2_Sec3}), $f_{k}$ is therefore differentiable 
and the derivative is just the above $f_{k}'$. Now the following holds:\\[4ex]
\refstepcounter{DSBcount}
\label{EWD_2_1_02}
\textbf{\hspace*{-0.8ex}\theDSBcount\ Lemma}\\[0.8ex]
Let $q \in \mathcal{H}$ (cf. (\ref{EWD_0_1_01}) for the definition of $\mathcal{H}$),
$k \in \mathbb{N}_{0}$ and $f_{k}$, $f_{k}'$, $f$, $f'$ 
as above.\index{function!$\mathcal{H}$}\\[2.8ex]
Furthermore, for $a \in \mathbb{R}$ let\\[2.5ex] 
\hspace*{12.1ex}$\displaystyle{\eta_{k}(a) = 
\sum_{i = 1}^{\lfloor k/2 \rfloor}\,a^{k + 1 - 2\,i}\,\prod_{j = 1}^{i-1}\,(k + 1 - 2\,j)
+ \prod_{j = 1}^{\lfloor k/2 \rfloor}\,(k + 1 - 2\,j)}$.\\[2.5ex]
As usual, $\displaystyle{\sum\limits_{\emptyset}\ldots = 0}$ and
$\displaystyle{\prod\limits_{\emptyset}\ldots = 1}$.\\[2.8ex]
Then, for all $k \in \mathbb{N}_{0}$ and $x \in \mathbb{R}$,
\begin{enumerate}
\item\label{EWD_2_1_02_BWa}
$|f_{k}(x)| \leq \eta_{k}(|x|)$,
\item\label{EWD_2_1_02_BWb}
$|f_{k}'(x)| \leq |x|\,\eta_{k}(|x|) + \Phi(|y|^{k})$.
\end{enumerate}
\vspace*{2.5ex}
\textbf{Proof:}\\[0.8ex]
First of all, we note that we may assume $q = 1_{(- \infty,\, z\,]}$, $z \in \mathbb{R}$, without loss of generality. 
From this we obtain the assertion for $q \in \mathcal{H}$ by first going to the 
convex combinations\index{function!convex combinations} of the functions 
$1_{(- \infty,\, z\,]}$, $z \in \mathbb{R}$, and then to their limits 
via e.g. majorized convergence\index{integration!using majorized convergence}.\\[2.8ex] 
In the following, let $z \in \mathbb{R}$ be fixed and $q = 1_{(- \infty,\, z\,]}$. Then\\[2ex]
\begin{tabular}{@{\hspace*{6ex}}l@{\hspace*{0.8ex}}c@{\hspace*{0.8ex}}l@{}}
$f_{k}(x)$&$=$&$\displaystyle{\psi(x)^{-1}\,\int\limits_{-\,\infty}^{x}
\Bigl(\,|y|^{k}\,1_{(- \infty,\, z\,]}(y) - 
\Phi\bigl(\,|v|^{k}\,1_{(- \infty,\, z\,]}(v)\,\bigr)\,\Bigr)\, \psi(y)\,dy}$,\\[4ex]
&$=$&$\displaystyle{\psi(x)^{-1} \cdot
\left\{
\begin{array}{cl}
\Phi\Bigl(\,|y|^{k}\,1_{(- \infty,\, x\,]}(y)\,\Bigr)
- \Phi\Bigl(\,|y|^{k}\,1_{(- \infty,\, z\,]}(y)\,\Bigr)\,\Phi(x)& \hspace*{4ex}
\text{for}\ x \leq z\, ,\\[3ex]
\Phi\Bigl(\,|y|^{k}\,1_{(- \infty,\, z\,]}(y)\,\Bigr)\,\Bigl(\,1 - \Phi(x)\,\Bigr)& \hspace*{4ex}
\text{for}\ x \geq z\, .
\end{array}  \right.}$
\end{tabular}\\[3ex]
With the definition\\[2ex]
\hspace*{12.1ex}$S_{k}(x) = \dfrac{1 - \Phi(x)}{\psi(x)}\,\Phi\Bigl(\,|y|^{k}\,1_{(- \infty,\, x\,]}(y)\,\Bigr)$,\\[2.5ex]
we can show the estimates\\[2ex]
\refstepcounter{DSBcount}
\label{EWD_2_1_03}
\text{\hspace*{-0.8ex}(\theDSBcount)}
\hspace*{4ex}
$- S_{k}(-x)\, \leq\, f_{k}(x)\, \leq\, S_{k}(x)$\hspace*{2.4ex}for all $x \in \mathbb{R}$, 
$k \in \mathbb{N}_{0}$.\\[2.5ex]
We obtain the upper estimate if we replace $z$ with $x$ in the representations of $f_{k}(x)$ 
both in the case of $x \leq z$ and in the case of $x \geq z$.\\[2.8ex]
Furthermore, since $f_{k}(x) \geq 0$ for $x \geq z$,
it is sufficient for the lower estimate to consider the case $x \leq z$:\\[2ex]
\hspace*{12.1ex}\begin{tabular}{@{}l@{\hspace*{0.8ex}}c@{\hspace*{0.8ex}}l@{}}
$f_{k}(x)$&$=$&$\displaystyle{\psi(x)^{-1}\,\biggl[\,\Phi\Bigl(\,|y|^{k}\,1_{(- \infty,\, x\,]}(y)\,\Bigr)
- \Phi\Bigl(\,|y|^{k}\,1_{(- \infty,\, z\,]}(y)\,\Bigr)\,\Phi(x)\,\biggr]}$\\[3ex]
&$\geq$&$\displaystyle{\psi(x)^{-1}\,\Phi(x)\,\biggl[\,\Phi\Bigl(\,|y|^{k}\,1_{(- \infty,\, x\,]}(y)\,\Bigr)
- \Phi\Bigl(\,|y|^{k}\,\Bigr)\,\biggr]}$\\[3ex]
&$=$&$\displaystyle{\psi(x)^{-1}\,\Phi(x)\,
\biggl[\,-\, \Phi\Bigl(\,|y|^{k}\,1_{(\,x,\, \infty\,)}(y)\,\Bigr)\,\biggr]}$
\end{tabular}\\[2.8ex]
\hspace*{12.1ex}\begin{tabular}{@{}l@{\hspace*{0.8ex}}c@{\hspace*{0.8ex}}l@{}}
\hspace*{5.4ex}&$=$&
$\displaystyle{-\, \psi(-x)^{-1}\,\Bigl(\,1 - \Phi(-x)\,\Bigr)\,\Phi\Bigl(\,|y|^{k}\,1_{(- \infty,\, - x\,]}(y)\,\Bigr) 
= - S_{k}(-x)}$.
\end{tabular}\\[3.3ex]
This shows (\ref{EWD_2_1_03}). From (\ref{EWD_2_1_03}) we obtain\\[2.4ex]
\refstepcounter{DSBcount}
\label{EWD_2_1_04}
\text{\hspace*{-0.8ex}(\theDSBcount)}
\hspace*{4ex}
$|f_{k}(x)| = \text{max}\,\bigl\{\,- f_{k}(x),\,f_{k}(x)\,\bigr\}
\leq \text{max}\,\bigl\{\,S_{k}(-x),\,S_{k}(x)\,\bigr\}$\hspace*{2.4ex}for all $x \in \mathbb{R}$, 
$k \in \mathbb{N}_{0}$.\\[2.5ex]
Part \ref{EWD_2_1_02_BWa}) is thus verified by the proof of the following estimate:\\[2.4ex]
\refstepcounter{DSBcount}
\label{EWD_2_1_05}
\text{\hspace*{-0.8ex}(\theDSBcount)}
\hspace*{4ex}
$0\, \leq\, S_{k}(x)\, \leq\, \eta_{k}(|x|)$\hspace*{2.4ex}for all $x \in \mathbb{R}$, 
$k \in \mathbb{N}_{0}$.\\[2.5ex]
To prove (\ref{EWD_2_1_05}), we first use that integration by parts\index{integration!by parts} 
(together with $\psi'(x) = - x\,\psi(x)$) 
gives the following recursion formulas:\\[2.4ex]
\hspace*{12.1ex}$\displaystyle{\Phi\Bigl(\,y^{k\, +\, 2}\,1_{(- \infty,\, x\,]}(y)\,\Bigr) =
- x^{k\, +\, 1}\,\psi(x)\, + (\,k\, +\, 1\,)\,\Phi\Bigl(\,y^{k}\,1_{(- \infty,\, x\,]}(y)\,\Bigr)}$
\hspace*{4ex}and\\[2ex]
\hspace*{12.1ex}$\displaystyle{
\Phi\Bigl(\,|y|^{k\,+\,2}\,\Bigr) = (\,k\, +\, 1\,)\,\Phi\Bigl(\,|y|^{k}\,\Bigr)}$
\hspace*{2.4ex}for all $x \in \mathbb{R}$, $k \in \mathbb{N}_{0}$.\\[3ex]
Furthermore, the following formulas follow from the definition of $\eta_{k}(x)$:\\[2.5ex]
\hspace*{12.1ex}$\displaystyle{
\eta_{k\, +\, 2}(x) = x^{k\, +\, 1}} + (k + 1)\,\eta_{k}(x)$\hspace*{4ex}and\\[1.5ex]
\hspace*{12.1ex}$\displaystyle{
\eta_{k}(0) = \prod_{j = 1}^{\lfloor k/2 \rfloor}\,(k + 1 - 2\,j)}$
\hspace*{2.4ex}for all $x \in \mathbb{R}$, $k \in \mathbb{N}_{0}$.\\[2.5ex]
Thus, by induction we get\\[2.8ex]
\refstepcounter{DSBcount}
\label{EWD_2_1_06}
\text{\hspace*{-0.8ex}(\theDSBcount)}
\hspace*{4ex}
$\displaystyle{
\Phi\Bigl(\,y^{k}\,1_{(- \infty,\, x\,]}(y)\,\Bigr) =
\left\{
\begin{array}{ll@{}}
-\psi(x)\,\eta_{k}(x)& \hspace*{2.4ex}
\text{for}\ k = 1, 3, 5, \ldots\\[2ex]
- \psi(x)\,\eta_{k}(x) + \Bigl(\,\psi(x) + \Phi(x)\,\Bigr)\,\eta_{k}(0)& \hspace*{2.4ex}
\text{for}\ k = 0, 2, 4, \ldots\,
\end{array}  \right.}$\\[3ex]
and\\[2ex]
\refstepcounter{DSBcount}
\label{EWD_2_1_07}
\text{\hspace*{-0.8ex}(\theDSBcount)}
\hspace*{4ex}
$\displaystyle{
\Phi\Bigl(\,|y|^{k}\,\Bigr) =
\left\{
\begin{array}{ll@{}}
\sqrt{\dfrac{2}{\pi}}\,\eta_{k}(0)& \hspace*{2.4ex}
\text{for}\ k = 1, 3, 5, \ldots\\[2ex]
\eta_{k}(0)& \hspace*{2.4ex}
\text{for}\ k = 0, 2, 4, \ldots\,\,.
\end{array}  \right.}$\\[3.5ex]
In addition, we observe the following properties of $\eta_{k}(x)$:\\[2ex]
\refstepcounter{DSBcount}
\label{EWD_2_1_08}
\text{\hspace*{-0.8ex}(\theDSBcount)}
\hspace*{4ex}
$\eta_{k}(x) = \eta_{k}(|x|)$\hspace*{4ex}for \textbf{odd} $k$ and $x \in \mathbb{R}$,\\[2.5ex]
\refstepcounter{DSBcount}
\label{EWD_2_1_09}
\text{\hspace*{-0.8ex}(\theDSBcount)}
\hspace*{4ex}
$-\,\bigl(\,\eta_{k}(x) - \eta_{k}(0)\,\bigr) = \eta_{k}(|x|) - \eta_{k}(0)$
\hspace*{4ex}for \textbf{even} $k$ and $x \leq 0$.\\[2.5ex]
Now (\ref{EWD_2_1_05}) follows directly from (\ref{EWD_2_1_06}) and (\ref{EWD_2_1_08}) 
for \textbf{odd} $k$ and $x \leq 0$.\\[2.8ex] 
Moreover, from the well-known inequality
(cf. e.g. C. Stein \cite{stein1972bound}, inequality (2.67), p. 594)\\[2ex]
\refstepcounter{DSBcount}
\label{EWD_2_1_10}
\text{\hspace*{-0.8ex}(\theDSBcount)}
\hspace*{2.8ex}
$\displaystyle{\dfrac{1 - \Phi(x)}{\psi(x)}\, \leq\, \dfrac{1}{x}}$
\hspace*{2.4ex}for all $x > 0$,\\[2.5ex]
it also follows that the function\\[2ex]
\hspace*{12.1ex}$\displaystyle{b(y) = \dfrac{1 - \Phi(y)}{\psi(y)}}$
\hspace*{2.4ex}is monotonically decreasing all over 
$\mathbb{R}$\index{function!monotonically decreasing}.\\[2.5ex]
Using this and then (\ref{EWD_2_1_07}), we get for \textbf{odd} $k$ and $x \geq 0$ the estimate\\[2ex]
\hspace*{12.1ex}$\displaystyle{S_{k}(x)\, \leq\, \dfrac{1}{2\,\psi(0)}\,\Phi\bigl(\,|y|^{k}\,\bigr)\,
=\, \eta_{k}(0)\, \leq\, \eta_{k}(|x|)}$.\\[2.5ex]
Furthermore, to prove (\ref{EWD_2_1_05}) for \textbf{even} $k$ we still need\\[2ex]
\refstepcounter{DSBcount}
\label{EWD_2_1_11}
\text{\hspace*{-0.8ex}(\theDSBcount)}
\hspace*{2.8ex}
$\displaystyle{\dfrac{\Phi(x)\,\bigl(\,1 - \Phi(x)\,\bigr)}{\psi(x)}\,
\leq\, \text{min}\,\Bigl\{\,\dfrac{1}{4\,\psi(x)},\,\dfrac{\Phi(|x|)}{|x|}\,\Bigr\}\, \leq\, 1}$
\hspace*{2.4ex}for all $x \in \mathbb{R}$.\\[2.5ex]
The first inequality of (\ref{EWD_2_1_11}) is valid because the function 
$g(x) = \Phi(x)\,\bigl(\,1-\Phi(x)\,\bigr)$ has its maximum at the point $x = 0$ with $g(0) = \dfrac{1}{4}$.
In addition, (\ref{EWD_2_1_10}) is needed again.\\[1.5ex]
For the second inequality of (\ref{EWD_2_1_11}), we use $\psi(x) \geq \dfrac{1}{4}$ for $|x| \leq 0,96$ and 
$\Phi(|x|) \leq |x|$ for $|x| \geq 0,79$.\\[2.8ex] 
Because of (\ref{EWD_2_1_06}) and (\ref{EWD_2_1_11}) it now follows for \textbf{even} $k$ and $x \geq 0$\\[2.5ex]
\hspace*{12.1ex}\begin{tabular}{@{}l@{\hspace*{0.8ex}}c@{\hspace*{0.8ex}}l@{}}
$S_{k}(x)$&$=$&$\displaystyle{\dfrac{1 - \Phi(x)}{\psi(x)}\,
\Bigl[\,- \psi(x)\,\eta_{k}(|x|) + \Bigl(\,\psi(x) + \Phi(x)\,\Bigr)\,\eta_{k}(0)\,\Bigr]\,}$\\[2.5ex]
&$\leq$&$\displaystyle{\dfrac{1 - \Phi(x)}{\psi(x)}\,
\Bigl[\,- \psi(x)\,\eta_{k}(|x|) + \Bigl(\,\psi(x) + \Phi(x)\,\Bigr)\,\eta_{k}(|x|)\,\Bigr]\,}$\\[2.5ex]
&$=$&$\displaystyle{\dfrac{1 - \Phi(x)}{\psi(x)}\,\Phi(x)\,\eta_{k}(|x|)\, \leq \eta_{k}(|x|)\,}$.
\end{tabular}\\[3ex]
For \textbf{even} $k$ and $x \leq 0$, on the other hand, (\ref{EWD_2_1_06}) and (\ref{EWD_2_1_11}) and additionally (\ref{EWD_2_1_09}) provide\\[2.5ex]
\hspace*{12.1ex}\begin{tabular}{@{}l@{\hspace*{0.8ex}}c@{\hspace*{0.8ex}}l@{}}
$S_{k}(x)$&$=$&$\displaystyle{\dfrac{1 - \Phi(x)}{\psi(x)}\,
\Bigl[\,- \psi(x)\,\Bigl(\,\eta_{k}(x) - \eta_{k}(0)\,\Bigr)\, + \Phi(x)\,\eta_{k}(0)\,\Bigr]\,}$\\[2.5ex]
&$=$&$\displaystyle{\dfrac{1 - \Phi(x)}{\psi(x)}\,
\Bigl[\,\psi(x)\,\Bigl(\,\eta_{k}(|x|) - \eta_{k}(0)\,\Bigr)\, + \Phi(x)\,\eta_{k}(0)\,\Bigr]\,}$\\[2.5ex]
&$=$&$\displaystyle{\Bigl(\,1 - \Phi(x)\,\Bigr)\,\Bigl(\,\eta_{k}(|x|) - \eta_{k}(0)\,\Bigr)\, +
\dfrac{\bigl(\,1 - \Phi(x)\,\bigr)\,\Phi(x)}{\psi(x)}\,\eta_{k}(0)}$
\end{tabular}\\[2.3ex]
\hspace*{12.1ex}\begin{tabular}{@{}l@{\hspace*{0.8ex}}c@{\hspace*{0.8ex}}l@{}}
\hspace*{5.6ex}&$\leq$&$\displaystyle{\eta_{k}(|x|) - \eta_{k}(0) + \eta_{k}(0) = \eta_{k}(|x|)}$.
\end{tabular}\\[2.5ex]
This concludes the proof of (\ref{EWD_2_1_05}) for all possible cases.\\[2.8ex] 
Multiplication of (\ref{EWD_2_1_03}) with $x > 0$ or $x < 0$ also provides\\[2ex]
\hspace*{12.1ex}\begin{tabular}{@{}r@{\hspace*{0.8ex}}c@{\hspace*{0.8ex}}l@{\hspace*{0.8ex}}c@{\hspace*{0.8ex}}l
@{\hspace*{4ex}}l@{}}
$- x\,S_{k}(-x)$&$\leq$&$x\,f_{k}(x)$&$\leq$&$x\,S_{k}(x)$&for $x > 0$\hspace*{4ex}and\\[2ex]
$x\,S_{k}(x)$&$\leq$&$x\,f_{k}(x)$&$\leq$&$- x\,S_{k}(- x)$&for $x < 0$,
\end{tabular}\\[2.5ex]
from which, in summary, for $x \in \mathbb{R}$ follows\\[2ex]
\hspace*{12.1ex}$- |x|\,S_{k}(- |x|)\, \leq\, x\,f_{k}(x)\, \leq\, |x|\,S_{k}(|x|)$.\\[2.5ex] 
Thus, by applying the definition (\ref{EWD_2_1_01}) and the inequalities
(\ref{EWD_2_1_10}) and (\ref{EWD_2_1_05}) we get \\[2ex]
\hspace*{12.1ex}\begin{tabular}{@{}l@{\hspace*{0.8ex}}c@{\hspace*{0.8ex}}l@{}}
$f_{k}'(x)$&$=$&$x\,f_{k}(x) + |x|^{k}\,q(x) - \Phi\bigl(\,|v|^{k}\,q(v)\,\bigr)$\\[2ex]
&$\leq$&$|x|\,S_{k}(|x|) + |x|^{k}$\\[2.2ex]
&$\leq$&$\displaystyle{\left\{
\begin{array}{l@{\hspace*{0.8ex}}l@{\hspace*{0.8ex}}ll@{}}
|x|\,\dfrac{1}{|x|}\,\Phi(|y|^k) + |x|^{k}& \leq& \Phi(|y|^k) + |x|\,\eta_{k}(|x|)
& \hspace*{3.4ex}
\text{for}\ x \not= 0\ \text{and}\ k \not= 0,\\[3ex]
|x|\,\eta_{k}(|x|) + |x|^{k}& \leq& |x|\,\eta_{k}(|x|) + \Phi(|y|^k)
& \hspace*{3.4ex}
\text{for}\ x = 0\ \text{or}\ k = 0.
\end{array}  \right.}$
\end{tabular}\\[2.5ex]
On the other hand, using (\ref{EWD_2_1_01}) and (\ref{EWD_2_1_05}) we obtain\\[2ex]
\hspace*{12.1ex}\begin{tabular}{@{}l@{\hspace*{0.8ex}}c@{\hspace*{0.8ex}}l@{}}
$f_{k}'(x)$&$\geq$&$- |x|\,S_{k}(- |x|) - \Phi(|y|^k)$\\[2ex]
&$\geq$&$- |x|\,\eta_{k}(\big| - |x|\big|) - \Phi(|y|^k)$\\[2ex]
&$=$&$- |x|\,\eta_{k}(|x|) - \Phi(|y|^k)$.
\end{tabular}\\[2.5ex]
All in all, this proves the part \ref{EWD_2_1_02_BWb}).\hfill$\Box$\\[4ex]
For $k = 0$ we find better estimates of $f_{k}$ and $f_{k}'$ in the literature.
These and simpler, but not as good 
estimates of $f_{k}$ and $f_{k}'$ for $k \geq 2$ are given below.\\[4ex]
\refstepcounter{DSBcount}
\label{EWD_2_1_12}
\textbf{\hspace*{-0.8ex}\theDSBcount\ Corollary}\\[0.8ex]
Suppose that the conditions and notations of Lemma \ref{EWD_2_1_02} apply.
Then, for all $x \in \mathbb{R}$,
\begin{enumerate}
\item\label{EWD_2_1_12_BWa}
$|f(x)| \leq \dfrac{\sqrt{2\pi}}{4} = 0,6267$
\hspace*{3ex}and\hspace*{3ex}
$|f'(x)| \leq 1$,
\item\label{EWD_2_1_12_BWb}
$|f_{1}(x)| \leq 1$
\hspace*{3ex}and\hspace*{3ex}
$|f_{1}'(x)| \leq |x| + \sqrt{\dfrac{2}{\pi}}$,
\item\label{EWD_2_1_12_BWc}
$|f_{k}(x)| \leq B_{k}\,\bigl(\,1 + |x|^{k-1}\,\bigr)$
\hspace*{3ex}and\hspace*{3ex}
$|f_{k}'(x)| \leq 2\,B_{k}\,\bigl(\,1 + |x|^{k}\,\bigr)$\hspace*{3ex}for $k \geq 2$,\\[2ex]
where $\displaystyle{B_{k} = 
\Bigl\lfloor\,\dfrac{k}{2}\,\Bigr\rfloor\,\prod_{j = 1}^{\lfloor k/2 \rfloor}\,(k + 1 - 2\,j)}$.
\end{enumerate}
\vspace*{2.5ex}
\textbf{Proof:}\\[0.8ex]
Part a): The calculations for $|f'(x)| \leq 1$ can be found in C. Stein \cite{stein1972bound}, 
proof of the inequality (2.66) from Lemma 2.5, p. 594, 595.
Alternatively, this inequality can also be found in
Chen, L. H. Y., Goldstein, L. and Shao, Q. \cite{Chen2011}, Lemma 2.3, p. 16 
(proof in the Appendix to Chapter 2, pp. 37, 38). 
The proof of $|f(x)| \leq \sqrt{2\pi}/4$ can also be found there.\\[2.8ex]
Part b): Due to $\eta_{1} \equiv 1$ and (\ref{EWD_2_1_07}), this follows from Lemma \ref{EWD_2_1_02}.\\[2.8ex]
Part c): At first, we can use $|x|^{l} \leq 1 + |x|^{r}$ for $0 < l < r$ to obtain the estimate\\[2ex]
\hspace*{12.1ex}$\displaystyle{\sum_{i = 1}^{\lfloor k/2 \rfloor}\,|x|^{k + 1 - 2\,i} \leq |x|^{k-1} + 
\sum_{i = 2}^{\lfloor k/2 \rfloor}\,\bigl(\,1 + |x|^{k-1}\,\bigr)
= \Bigl\lfloor\,\dfrac{k}{2}\,\Bigr\rfloor\,|x|^{k-1} 
+ \Bigl\lfloor\,\dfrac{k}{2}\,\Bigr\rfloor - 1}$\hspace*{3.2ex}for $k \geq 2$.\\[2.5ex]
From this and Lemma \ref{EWD_2_1_02} we get\\[2ex]
\hspace*{9.4ex}\begin{tabular}{@{}l@{\hspace*{0.8ex}}c@{\hspace*{0.8ex}}l@{}}
$|f_{k}(x)| \leq \eta_{k}(|x|)$&$\leq$&$\displaystyle{B_{k}\,|x|^{k-1} + 
\Bigl(\,\Bigl\lfloor\,\dfrac{k}{2}\,\Bigr\rfloor - 1\,\Bigr)\,\prod_{j = 1}^{\lfloor k/2 \rfloor}\,(k + 1 - 2\,j) + 
\prod_{j = 1}^{\lfloor k/2 \rfloor}\,(k + 1 - 2\,j)}$\\[3ex]
&$=$&$B_{k}\,\bigl(\,|x|^{k-1} + 1\,\bigr)$\hspace*{4ex}for $k \geq 2$.
\end{tabular}\\[2.5ex]
Finally, we additionally use $\Phi(|y|^k) \leq B_{k}$ for $k \geq 2$ (cf. (\ref{EWD_2_1_07})) 
to estimate $|f_{k}'(x)|$.\hfill$\Box$\\[4ex]
Next, we prove monotonicity properties of some of the functions considered here.\\[4ex]
\refstepcounter{DSBcount}
\label{EWD_2_1_13}
\textbf{\hspace*{-0.8ex}\theDSBcount\ Lemma}\\[0.8ex]
Let $q \in \mathcal{H}_{c} = \Bigl\{\, h \in \mathcal{H}\hspace*{0.1ex}:\hspace*{0.1ex}
\text{there exists a $K = K(h) \in \mathbb{R}$ 
so that $h$ is continuous in 
$\mathbb{R} \setminus [-K,\,K]$}\, \Bigr\}$.\index{function!$\mathcal{H}_{c}$}\vspace*{0.5ex}
Furthermore, let $q(+ \infty) = \lim\limits_{x \rightarrow + \infty} q(x)$
and $q(- \infty) = \lim\limits_{x \rightarrow - \infty} q(x)$, and\\[2ex]
\hspace*{12.1ex}$\displaystyle{f^{1}(x) = \psi(x)^{-1}\,\int\limits_{-\,\infty}^{x}
\Bigl(\,y\,q(y) - \Phi\bigl(\,v\,q(v)\,\bigr)\,\Bigr)\, \psi(y)\,dy}$\\[2.5ex]
and $f$ as above (see (\ref{EWD_2_1_01})). Then
\begin{enumerate}
\item\label{EWD_2_1_13_BWa}
The mapping $x \rightarrow x\,f(x)$ is \textbf{monotonically increasing}\index{function!monotonically increasing} 
with\\[2ex]
\hspace*{4ex}$\displaystyle{\lim\limits_{x \rightarrow - \infty} x\,f(x) = \Phi(q) - q(- \infty)}$
\hspace*{4ex}and\hspace*{4ex}
$\displaystyle{\lim\limits_{x \rightarrow + \infty} x\,f(x) = \Phi(q) - q(+ \infty)}$,\vspace*{0.5ex}
\item\label{EWD_2_1_13_BWb}
The mapping $x \rightarrow f^{1}(x)$ is \textbf{monotonically increasing}\index{function!monotonically increasing} 
with\\[2ex]
\hspace*{4ex}$\displaystyle{\lim\limits_{x \rightarrow - \infty} f^{1}(x) = - q(- \infty)}$
\hspace*{4ex}and\hspace*{4ex}
$\displaystyle{\lim\limits_{x \rightarrow + \infty} f^{1}(x) = - q(+ \infty)}$.
\end{enumerate}
\vspace*{2.5ex}
\textbf{Proof:}\\[0.8ex]
An application of l'Hospital's rule\index{l'Hospital's rule} yields\\[1.5ex]
\hspace*{0.5ex}\begin{tabular}{@{}l@{\hspace*{0.8ex}}c@{\hspace*{0.8ex}}l@{}}
$\displaystyle{\lim\limits_{x \rightarrow \pm \infty} x\,f(x)}$&$=$&
$\displaystyle{\lim\limits_{x \rightarrow \pm \infty}\,\dfrac{\displaystyle{\int\limits_{-\,\infty}^{x}
\Bigl(\,q(y) - \Phi(q)\,\Bigr)\, \psi(y)\,dy}}{\dfrac{1}{x}\,\psi(x)} =
\lim\limits_{x \rightarrow \pm \infty}\,\dfrac{\Bigl(\,q(x) - \Phi(q)\,\Bigr)\, \psi(x)}
{- \dfrac{1}{x^2}\,\psi(x) - \psi(x)}}$\\[4ex]
&$=$&$\Phi(q) - q(\pm \infty)$
\end{tabular}\\[2ex]
and\\[2ex]
\hspace*{0.5ex}\begin{tabular}{@{}l@{\hspace*{0.8ex}}c@{\hspace*{0.8ex}}l@{}}
$\displaystyle{\lim\limits_{x \rightarrow \pm \infty} f^{1}(x)}$&$=$&
$\displaystyle{\lim\limits_{x \rightarrow \pm \infty}\,\dfrac{\displaystyle{\int\limits_{-\,\infty}^{x}
\Bigl(\,y\,q(y) - \Phi\bigl(\,v\,q(v)\,\bigr)\,\Bigr)\, \psi(y)\,dy}}{\psi(x)} =
\lim\limits_{x \rightarrow \pm \infty}\,\dfrac{\Bigl(\,x\,q(x) - \Phi\bigl(\,v\,q(v)\,\bigr)\,\Bigr)\, \psi(x)}
{- x\,\psi(x)}}$\\[4ex]
&$=$&$- q(\pm \infty)$.
\end{tabular}\\[2.5ex]
To prove that the functions $x\,f(x)$ and $f^{1}(x)$ are monotonically increasing, we can assume 
as in the proof of Lemma \ref{EWD_2_1_02} without loss of generality, that 
$q = 1_{(- \infty,\, z\,]}$ with a fixed $z \in \mathbb{R}$.\\[2.8ex]
We then get\\[2ex]
\hspace*{12.1ex}$\displaystyle{x\,f(x) = x\,\psi(x)^{-1} \cdot
\left\{
\begin{array}{ll@{}}
\Phi(x)\,\bigl(\,1 - \Phi(z)\,\bigr)& \hspace*{2.4ex}
\text{for all}\ x \leq z,\\[2ex]
\Phi(z)\,\bigl(\,1 - \Phi(x)\,\bigr)& \hspace*{2.4ex}
\text{for all}\ x \geq z.
\end{array}  \right.}$\\[3ex]
Differentiating this function at the point $x \not= z$ gives\\[2.5ex]
\hspace*{12.1ex}$\displaystyle{\dfrac{d}{dx}\,x\,f(x) = 
\left\{
\begin{array}{ll@{}}
\bigl(\,1 - \Phi(z)\,\bigr)\,T(-x)& \hspace*{2.4ex}
\text{for all}\ x < z,\\[2ex]
\Phi(z)\,T(x)& \hspace*{2.4ex}
\text{for all}\ x > z,
\end{array}  \right.}$\\[2.5ex]
where\\[2ex]
\hspace*{12.1ex}$\displaystyle{T(x) = (1 + x^2)\,\dfrac{1 - \Phi(x)}{\psi(x)} - x}$.\\[2.5ex]
Since the function $x \rightarrow x\,f(x)$ is continuous in $z$, we obtain the desired monotonicity 
if we show $T(x) \geq 0$ for $x \in \mathbb{R}$.\\[2.8ex]
For $x \leq 0$, however, this is clear, and for $x > 0$ this follows from the inequality\\[2.5ex]
\refstepcounter{DSBcount}
\label{EWD_2_1_14}
\text{\hspace*{-0.8ex}(\theDSBcount)}
\hspace*{2.8ex}
$\displaystyle{\dfrac{1 - \Phi(x)}{\psi(x)} \geq \dfrac{x}{1 + x^2}}$
\hspace*{2.4ex}for all $x > 0$
\hspace*{5.8ex}(cf. e.g. McKean \cite{mckean1969stochastic}, page 4, Problem 1).\\[2.5ex]
In addition, we have\\[2ex]
\hspace*{12.1ex}\begin{tabular}{@{}l@{\hspace*{0.8ex}}c@{\hspace*{0.8ex}}l@{}}
$f^{1}(x)$&$=$&$\displaystyle{\psi(x)^{-1} \cdot
\left\{
\begin{array}{ll@{}}
\Phi\bigl(\,y\,1_{(- \infty,\, x\,]}(y)\,\bigr) - \Phi\bigl(\,y\,1_{(- \infty,\, z\,]}(y)\,\bigr)\,\Phi(x)
& \hspace*{2.4ex}
\text{for all}\ x \leq z,\\[2.5ex]
\Phi\bigl(\,y\,1_{(- \infty,\, z\,]}(y)\,\bigr)\,\bigl(\,1 - \Phi(x)\,\bigr)& \hspace*{2.4ex}
\text{for all}\ x \geq z,
\end{array}  \right.}$\\[6.5ex]
&$=$&$\displaystyle{
\left\{
\begin{array}{ll@{}}
\psi(z)\,\dfrac{\Phi(x)}{\psi(x)} - 1
& \hspace*{2.4ex}
\text{for all}\ x \leq z,\\[2.5ex]
- \psi(z)\,\dfrac{1 - \Phi(x)}{\psi(x)}& \hspace*{2.4ex}
\text{for all}\ x \geq z.
\end{array}  \right.}$
\end{tabular}\\[3ex]
Differentiating this function at the point $x \not= z$ gives\\[2.5ex]
\refstepcounter{DSBcount}
\label{EWD_2_1_15}
\text{\hspace*{-0.8ex}(\theDSBcount)}
\hspace*{2.8ex}
$\displaystyle{\dfrac{d}{dx}\,f^{1}(x) = 
\left\{
\begin{array}{ll@{}}
\psi(z)\,\bigl(\,1 - V(-x)\,\bigr)& \hspace*{2.4ex}
\text{for all}\ x < z,\\[2ex]
\psi(z)\,\bigl(\,1 - V(x)\,\bigr)& \hspace*{2.4ex}
\text{for all}\ x > z,
\end{array}  \right.}$\\[2.5ex]
where\\[2ex]
\refstepcounter{DSBcount}
\label{EWD_2_1_16}
\text{\hspace*{-0.8ex}(\theDSBcount)}
\hspace*{2.8ex}
$\displaystyle{V(x) = \dfrac{x\,\bigl(\,1 - \Phi(x)\,\bigr)}{\psi(x)}}$.\\[2.5ex]
Since the function $x \rightarrow f^{1}(x)$ is continuous in $z$, we obtain the desired monotonicity 
analogously to above if we show $V(x) \leq 1$ for $x \in \mathbb{R}$.\\[2.8ex]
This is clear for $x \leq 0$, and for $x > 0$ this follows from the inequality (\ref{EWD_2_1_10}), 
which has already been used several times.\hfill$\Box$\\[4ex]
From the lemma just proved we obtain further important estimates.\\[3.5ex]
\refstepcounter{DSBcount}
\label{EWD_2_1_17}
\textbf{\hspace*{-0.8ex}\theDSBcount\ Corollary}\\[0.8ex]
Suppose that the conditions and notations of Lemma \ref{EWD_2_1_13} apply.
Furthermore, we \textbf{define}\\[2ex]
\hspace*{12.1ex}$(f^{1})'(x) = x\,f^{1}(x) + x\,q(x) - \Phi\bigl(\,v\,q(v)\,\bigr)$
\hfill(Stein's equation\index{Stein's equation})\\[2.5ex]
analogously to (\ref{EWD_2_1_01}). Then, for all $x \in \mathbb{R}$,
\begin{enumerate}
\item\label{EWD_2_1_17_BWa}
$|x\,f(x)| \leq \text{max}\,\bigl\{\,q(- \infty) - \Phi(q),\,\Phi(q) - q(+ \infty)\,\bigr\}\leq 1$,
\item\label{EWD_2_1_17_BWb}
$|f^{1}(x)| \leq \text{max}\,\bigl\{\,q(- \infty),\,q(+ \infty)\,\bigr\}\leq 1$,
\item\label{EWD_2_1_17_BWc}
$|(f^{1})'(x)| \leq |x| + \dfrac{1}{\sqrt{2\,\pi}}$.
\end{enumerate}
\vspace*{2.5ex}
\textbf{Proof:}\\[0.8ex]
Parts a) and b) follow directly from Lemma \ref{EWD_2_1_13}, so that only part c) remains to be proved.\\[2.8ex]
As in the proof of Lemma \ref{EWD_2_1_02} and Lemma \ref{EWD_2_1_13}, we can assume without loss of generality, that 
$q = 1_{(- \infty,\, z\,]}$ with a fixed $z \in \mathbb{R}$.\\[2.8ex]
Then $(f^{1})'(x)$ has the representation (\ref{EWD_2_1_15}) with the only difference that $x \leq z$ instead 
of $x < z$ holds. Moreover, it follows from the (proof of) Lemma \ref{EWD_2_1_13} that $(f^{1})'(x) \geq 0$ 
for all $x \in \mathbb{R}$.\\[2.8ex]
To prove part c), we consider the following four cases. We use the function $V(x)$ 
according to (\ref{EWD_2_1_16}).\\[2.8ex]
\textit{Case 1:} $x \leq 0 \leq z$ or $x \leq z \leq 0$.\\[0.8ex]
In this case, we have $(-x) \geq 0$ and therefore $V(-x) \geq 0$, from which 
$(f^{1})'(x) \leq \psi(z) \leq \psi(0) = \dfrac{1}{\sqrt{2\,\pi}}$ follows.\\[2.8ex]
\textit{Case 2:} $0 \leq z < x$ or $z < 0 \leq x$.\\[0.8ex]
Now $V(x) \geq 0$ is valid, from which again $(f^{1})'(x) \leq \psi(z) \leq \psi(0) 
= \dfrac{1}{\sqrt{2\,\pi}}$ follows.\\[2.8ex]
\textit{Case 3:} $0 \leq x \leq z$.\\[0.8ex]
In this case, we get\\[2ex]
\refstepcounter{DSBcount}
\label{EWD_2_1_18}
\text{\hspace*{-0.8ex}(\theDSBcount)}
\hspace*{2.8ex}
$\psi(z) \leq \psi(x) \leq \psi(0) = \dfrac{1}{\sqrt{2\,\pi}}$\\[2ex]
and\\[1.5ex]
\hspace*{12.1ex}$V(-x) = \dfrac{(-x)\,\bigl(\,1-\Phi(-x)\,\bigr)}{\psi(-x)}
= \dfrac{(-x)\,\Phi\bigl(\,-(-x)\,\bigr)}{\psi(x)} = \dfrac{(-x)\,\Phi(x)}{\psi(x)}$.\\[2.5ex]
It follows that\\[2ex]
\hspace*{12.1ex}\begin{tabular}{@{}l@{\hspace*{0.8ex}}c@{\hspace*{0.8ex}}l@{\hspace*{22.6ex}}r@{}}
$(f^{1})'(x)$&$=$&$\psi(z)\,\Bigl(\,1 + \dfrac{x\,\Phi(x)}{\psi(x)}\,\Bigr)$\\[3ex]
&$\leq$&$\psi(z) + |x|\, \dfrac{\psi(z)}{\psi(x)}$&($\Phi(x) \leq 1$
\hspace*{1.4ex}and\hspace*{1.4ex} 
$x = |x|$)\\[3ex]
&$\leq$&$\psi(0) + |x| = |x| + \dfrac{1}{\sqrt{2\,\pi}}$&($\dfrac{\psi(z)}{\psi(x)} \leq 1$).
\end{tabular}\\[2.8ex]
\textit{Case 4:} $z < x \leq 0$.\\[0.8ex]
In case 4, the inequality chain (\ref{EWD_2_1_18}) holds analogously to case 3.
It follows that\\[2.5ex]
\hspace*{12.1ex}\begin{tabular}{@{}l@{\hspace*{0.8ex}}c@{\hspace*{0.8ex}}l@{\hspace*{11.4ex}}r@{}}
$(f^{1})'(x)$&$=$&$\psi(z)\,\Bigl(\,1 - x\,\dfrac{\bigl(\,1-\Phi(x)\,\bigr)}{\psi(x)}\,\Bigr)$\\[3ex]
&$\leq$&$\psi(z) + |x|\, \dfrac{\psi(z)}{\psi(x)}$&($\bigl(\,1 - \Phi(x)\,\bigr) \leq 1$
\hspace*{1.4ex}and\hspace*{1.4ex} 
$-x = |x|$)\\[3ex]
&$\leq$&$\psi(0) + |x| = |x| + \dfrac{1}{\sqrt{2\,\pi}}$&($\dfrac{\psi(z)}{\psi(x)} \leq 1$). $\Box$
\end{tabular}\vspace*{2ex}

\section[A relation between the Hermite polynomials and the functions $q(x)$ and $f(x)$]
{A relation between the Hermite polynomials and the functions 
{\boldmath $q(x)$} and {\boldmath $f(x)$}}\label{EWD_Kap2_Sec2}

Let $n \in \mathbb{N}_{0}$. We denote by\\[2ex]
\hspace*{12.1ex}$\displaystyle{H_{n}(x) = (-1)^{n}\,e^{x^2/2}\,\dfrac{d^{n}}{dx^{n}}\,e^{-\,x^2/2}}$\\[2.5ex]
the $n$th-order \textbf{Hermite polynomial}\index{polynomial!Hermite polynomial}. 
We obtain this Hermite polynomial if we differentiate the function $\psi(x)$ $n$ times, i.e.
trivially holds \\[2ex]
\refstepcounter{DSBcount}
\label{EWD_2_2_01}
\text{\hspace*{-0.8ex}(\theDSBcount)}
\hspace*{4ex}
$\displaystyle{\dfrac{d^{n}}{dx^{n}}\,\psi(x) = (-1)^{n}\,\psi(x)\,H_{n}(x)}$.\\[2.5ex]
This immediately leads to the following useful formulas\\[2.5ex] 
\refstepcounter{DSBcount}
\label{EWD_2_2_02}
\text{\hspace*{-0.8ex}(\theDSBcount)}
\hspace*{4ex}
$H_{n+1}(x)\,\psi(x) = (-1)^{n+1}\,\Bigl(\,\dfrac{d^{n}}{dx^{n}}\,\psi(x)\,\Bigr)'
= \Bigl(\,-\, H_{n}(x)\,\psi(x)\,\Bigr)'$\\[2.5ex]
and\\[2ex]
\refstepcounter{DSBcount}
\label{EWD_2_2_03}
\text{\hspace*{-0.8ex}(\theDSBcount)}
\hspace*{4ex}
$\displaystyle{\int\limits_{\mathbb{R}}H_{n+1}(y)\,\psi(y)\,dy = 
\Bigl[\,-\, H_{n}(x)\,\psi(x)\,\Bigr]_{- \infty}^{+ \infty} = 0}$
\hspace*{4ex}for all $n \in \mathbb{N}_{0}$.\\[2.5ex]
The first nine Hermite polynomials are given as examples:\\[2ex]
\hspace*{12.1ex}$H_{0}(x) = 1$,\hspace*{3ex}$H_{1}(x) = x$,\hspace*{3ex}$H_{2}(x) = x^{2} - 1$,
\hspace*{3ex}$H_{3}(x) = x^{3} - 3\,x$,\\[2ex]
\hspace*{12.1ex}$H_{4}(x) = x^{4} - 6\,x^{2} + 3$,\hspace*{3ex}$H_{5}(x) = x^{5} - 10\,x^{3} 
+ 15\,x$,\\[2ex]
\hspace*{12.1ex}$H_{6}(x) = x^{6} - 15\,x^{4} + 45\,x^{2} - 15$,
\hspace*{3ex}$H_{7}(x) = x^{7} - 21\,x^{5} + 105\,x^{3} - 105\,x$,\\[2ex]
\hspace*{12.1ex}$H_{8}(x) = 
x^{8} - 28\,x^{6} + 210\,x^{4} - 420\,x^{2} + 105$.\index{polynomial!Hermite polynomial!$H_{0},\ldots,H_{8}$}\\[2.5ex]
Further Hermite polynomials can be specified by using the following recursion formula, which can be proved by simple calculations (using the product rule for $n$-fold derivatives):\\[2ex]
\refstepcounter{DSBcount}
\label{EWD_2_2_04}
\text{\hspace*{-0.8ex}(\theDSBcount)}
\hspace*{4ex}
$H_{n+1}(x) = (- n)\,H_{n-1}(x) + x\,H_{n}(x)$
\hspace*{4ex}for $n \in \mathbb{N}$, $x \in \mathbb{R}$.\\[2.5ex]
Hermite polynomials are of interest for Edgeworth expansions because they are part of these expansions,
cf. (\ref{EWD_1_1_02}) in chapter \ref{EWD_Kap1} 
$\bigl($and (\ref{EWD_3_1_06}), (\ref{EWD_3_1_07}) in chapter \ref{EWD_Kap3}$\bigr)$:\\[2.5ex] 
\hspace*{12.1ex}$e_{n}(x) = \Phi(x) - \dfrac{\mu_{3}}{6\,\sqrt{n}}\,H_{2}(x)\,\psi(x) 
= \Phi(x) - \dfrac{\mu_{3}}{6\,\sqrt{n}}\,\psi''(x)$.\\[2.8ex]
The next lemma is therefore an important tool that we need to determine the correct
Edgeworth expansions when applying Stein's method.\\[4ex]
\refstepcounter{DSBcount}
\label{EWD_2_2_05}
\textbf{\hspace*{-0.8ex}\theDSBcount\ Lemma}\\[0.8ex]
\begin{tabular}{@{}r@{\hspace*{0.8ex}}c@{\hspace*{0.8ex}}l@{}}
Let $\mathcal{H}_{d}$&$=$&$\Bigl\{\,h \in \mathcal{H}\,:\,h\ \text{is differentiable with a bounded 
derivative}\ h'\, \Bigr\}$.\index{function!$\mathcal{H}_{d}$, $\mathcal{H}_{e}$, $\mathcal{H}_{e1}$, $\mathcal{H}_{e2}$}
Furthermore, let\\[1ex]
$q \in \mathcal{H}_{e}$&$=$&$\mathcal{H}_{e1} \cup \mathcal{H}_{e2} =
\Bigl\{\, h\,:\,h \in \mathcal{H}_{d}\, \Bigr\} \cup
\Bigl\{\, id_{\mathbb{R}} \cdot h\,:\,h \in \mathcal{H}_{d}\, \Bigr\}$ and
\end{tabular}\\[2ex]
\hspace*{12.1ex}$\displaystyle{f(x) = \psi(x)^{-1}\,\int\limits_{-\,\infty}^{x}
\bigl(\,q(y) - \Phi(q)\,\bigr)\, \psi(y)\,dy}$.\\[2.5ex]
Then
\begin{enumerate}
\item\label{EWD_2_2_05_BWa}
$\displaystyle{\int\limits_{\mathbb{R}}f(y)\,H_{n}(y)\,\psi(y)\,dy = 
- \dfrac{1}{n+1}\,\int\limits_{\mathbb{R}}q(y)\,H_{n+1}(y)\,\psi(y)\,dy}$,
\item\label{EWD_2_2_05_BWb}
$\displaystyle{\int\limits_{\mathbb{R}}f(y)\,H_{n}(y)\,\psi(y)\,dy = 
\int\limits_{\mathbb{R}}f'(y)\,H_{n-1}(y)\,\psi(y)\,dy =
\int\limits_{\mathbb{R}}f''(y)\,H_{n-2}(y)\,\psi(y)\,dy}$.
\end{enumerate} 
In particular, we obtain the special cases:\vspace*{0.5ex}
\begin{enumerate}
\addtocounter{enumi}{+2}
\item\label{EWD_2_2_05_BWc}
$\displaystyle{\int\limits_{\mathbb{R}}f'(y)\,y\,\psi(y)\,dy = 
\dfrac{1}{3}\,\int\limits_{\mathbb{R}}q(y)\,\bigl(\,3\,y - y^3\,\bigr)\,\psi(y)\,dy}$,
\item\label{EWD_2_2_05_BWd}
$\displaystyle{\int\limits_{\mathbb{R}}f''(y)\,y\,\psi(y)\,dy = 
\dfrac{1}{4}\,\int\limits_{\mathbb{R}}q(y)\,\bigl(\,- y^4 + 6\,y^2 - 3\,\bigr)\,\psi(y)\,dy}$,
\item\label{EWD_2_2_05_BWe}
$\displaystyle{\dfrac{1}{3}\,\int\limits_{\mathbb{R}}f'(y)\,y\,\bigl(\,3\,y - y^3\,\bigr)\,\psi(y)\,dy
+ \int\limits_{\mathbb{R}}f''(y)\,y\,\psi(y)\,dy}$\\[2ex] 
\hspace*{19ex}$\displaystyle{=\ 
\dfrac{1}{18}\,\int\limits_{\mathbb{R}}q(y)\,\bigl(\,y^6 - 15\,y^4 + 45\,y^2 - 15\,\bigr)\,\psi(y)\,dy}$.
\end{enumerate}
\vspace*{2.5ex}
\textbf{Proof:}\\[0.8ex]
First of all, we note that the integrals \\[2ex]
\hspace*{7.1ex}$\displaystyle{\int\limits_{\mathbb{R}}f(y)\,p(y)\,\psi(y)\,dy}$,
\hspace*{3ex}$\displaystyle{\int\limits_{\mathbb{R}}f'(y)\,p(y)\,\psi(y)\,dy}$
\hspace*{3ex}and\hspace*{3ex}
$\displaystyle{\int\limits_{\mathbb{R}}f''(y)\,p(y)\,\psi(y)\,dy}$\\[2.5ex]
exist for all polynomials $p$. This follows for the first two integrals in the case of $q \in \mathcal{H}_{e1}$
because of Corollary \ref{EWD_2_1_12}, \ref{EWD_2_1_12_BWa}) and in the case of $q \in \mathcal{H}_{e2}$
because of Corollary \ref{EWD_2_1_17}, \ref{EWD_2_1_17_BWb}) and \ref{EWD_2_1_17_BWc}).\\[2.8ex]
Differentiating Stein's equation $f'(x) = x\,f(x) + q(x) - \Phi(q)$ (see (\ref{EWD_2_1_01})) gives\\[2ex]
\hspace*{12.1ex}$f''(x) = f(x) + x\,f'(x) + q'(x)$,\\[2ex]
from which, together with the above references, we can conclude that the third integral also exists.\\[2.5ex]
The above references also show that\\[2ex]
\refstepcounter{DSBcount}
\label{EWD_2_2_06}
\text{\hspace*{-0.8ex}(\theDSBcount)}
\hspace*{4ex}
$\displaystyle{\lim\limits_{x \rightarrow \pm \infty}\,f(x)\,H_{n}(x)\,\psi(x) = 0}$\hspace*{4ex}and\\[2.5ex]
\refstepcounter{DSBcount}
\label{EWD_2_2_07}
\text{\hspace*{-0.8ex}(\theDSBcount)}
\hspace*{4ex}
$\displaystyle{\lim\limits_{x \rightarrow \pm \infty}\,f'(x)\,H_{n}(x)\,\psi(x) = 0}$\hspace*{4ex}for
all $n \in \mathbb{N}_{0}$.\\[2.5ex]
We now obtain the assertions of the lemma as follows:\\[3ex]
\begin{tabular}{@{}l@{\hspace*{4ex}}c@{\hspace*{0.8ex}}l@{\hspace*{-51ex}}r@{}}
a)&&$\displaystyle{\int\limits_{\mathbb{R}}f(y)\,H_{n}(y)\,\psi(y)\,dy}$\\[4ex]
&$=$&$\displaystyle{- \dfrac{1}{n+1}\,\biggl(\,\int\limits_{\mathbb{R}}f(y)\,H_{n+2}(y)\,\psi(y)\,dy
- \int\limits_{\mathbb{R}}y\,f(y)\,H_{n+1}(y)\,\psi(y)\,dy\,\biggr)}$&(cf. (\ref{EWD_2_2_04}))\\[4ex]
&$=$&$\displaystyle{- \dfrac{1}{n+1}\,\biggl(\,\int\limits_{\mathbb{R}}f'(y)\,H_{n+1}(y)\,\psi(y)\,dy
- \int\limits_{\mathbb{R}}y\,f(y)\,H_{n+1}(y)\,\psi(y)\,dy\,\biggr)}$\\[4.5ex]
&&&(integration by parts\index{integration!by parts} using 
(\ref{EWD_2_2_02}) (with $n + 1$ instead of $n$) and (\ref{EWD_2_2_06}))
\end{tabular}

\newpage

\hspace*{2.1ex}\begin{tabular}{@{}l@{\hspace*{3ex}}c@{\hspace*{0.8ex}}l@{\hspace*{14.2ex}}r@{}}
&$=$&$\displaystyle{- \dfrac{1}{n+1}\,\int\limits_{\mathbb{R}}
\Bigl(f'(y) - y\,f(y) + \Phi(q)\,\Bigr)\,H_{n+1}(y)\,\psi(y)\,dy}$&(cf. (\ref{EWD_2_2_03}))\\[4.5ex]
&$=$&$\displaystyle{- \dfrac{1}{n+1}\,\int\limits_{\mathbb{R}}
q(y)\,H_{n+1}(y)\,\psi(y)\,dy}$&(cf. Stein's equation\index{Stein's equation}).
\end{tabular}\\[2ex]
\begin{tabular}{@{}l@{\hspace*{4ex}}l@{}}
b)&Double integration by parts\index{integration!by parts} as in \ref{EWD_2_2_05_BWa}) using 
(\ref{EWD_2_2_02}) and (\ref{EWD_2_2_06}), (\ref{EWD_2_2_07}).\\[2ex]
c)&Straightforward conclusion from \ref{EWD_2_2_05_BWa}) and \ref{EWD_2_2_05_BWb}) with $n = 2$.\\[2ex]
d)&Straightforward conclusion from \ref{EWD_2_2_05_BWa}) and \ref{EWD_2_2_05_BWb}) with $n = 3$.
\end{tabular}\\[2ex]
\begin{tabular}{@{}l@{\hspace*{4ex}}c@{\hspace*{0.8ex}}l@{\hspace*{-5.5ex}}r@{}}
e)&&$\displaystyle{\dfrac{1}{3}\,\int\limits_{\mathbb{R}}f'(y)\,y\,\bigl(\,3\,y - y^3\,\bigr)\,\psi(y)\,dy
+ \int\limits_{\mathbb{R}}f''(y)\,y\,\psi(y)\,dy}$\\[4.5ex]
&$=$&$\displaystyle{\int\limits_{\mathbb{R}}f'(y)\,\Bigl[\,H_{2}(y) - \dfrac{1}{3}\,y\,H_{3}(y)\,\Bigr]\,\psi(y)\,dy}$
&(cf. \ref{EWD_2_2_05_BWb}) with $n = 3$)\\[4.5ex]
&$=$&$- \dfrac{1}{3}\,\displaystyle{\int\limits_{\mathbb{R}}f'(y)\,H_{4}(y)\,\psi(y)\,dy}$
&(cf. (\ref{EWD_2_2_04}) with $n = 3$)\\[4.5ex]
&$=$&$\dfrac{1}{18}\,\displaystyle{\int\limits_{\mathbb{R}}q(y)\,H_{6}(y)\,\psi(y)\,dy}$
&(cf. \ref{EWD_2_2_05_BWb}) with $n = 5$ and \ref{EWD_2_2_05_BWa}) with $n =5$). $\Box$
\end{tabular}\\[4ex]
We conclude this section with some easy-to-calculate 
estimates of terms such as $||\,H_{n}\,\psi\,||$, etc.\\[4ex]
\refstepcounter{DSBcount}
\label{EWD_2_2_08}
\textbf{\hspace*{-0.8ex}\theDSBcount\ Lemma}
\vspace*{1ex}
\begin{enumerate}
\item\label{EWD_2_2_08_BWa}
$||\,H_{0}\,\psi\,|| = ||\,\psi\,|| = \psi(0) = \dfrac{1}{\sqrt{2\,\pi}} = 0,39894 \leq \dfrac{2}{5}$,\vspace*{0.8ex}
\item\label{EWD_2_2_08_BWb}
$||\,H_{1}\,\psi\,|| = ||\,x\,\psi(x)\,|| = \psi(1) = \dfrac{1}{\sqrt{2\,\pi\,e}} = 
0,24197 \leq \dfrac{1}{4}$,\vspace*{0.8ex}
\item\label{EWD_2_2_08_BWc}
$||\,x\,\bigl(H_{1}\,\psi\bigr)(x)\,|| = ||\,x^2\,\psi(x)\,|| = 2\,\psi(\sqrt{2}) = \dfrac{\sqrt{2}}{e\,\sqrt{\pi}} 
= 0,29353 \leq \dfrac{3}{10}$,\vspace*{0.8ex}
\item\label{EWD_2_2_08_BWd}
$||\,x^2\,\bigl(H_{1}\,\psi\bigr)(x)\,|| = ||\,x^3\,\psi(x)\,|| = 3^{3/2}\,\psi(\sqrt{3}) 
= 0,46254 \leq \dfrac{7}{15}$,\vspace*{0.8ex}
\item\label{EWD_2_2_08_BWe}
$||\,H_{2}\,\psi\,|| = \psi(0) = \dfrac{1}{\sqrt{2\,\pi}} = 0,39894 \leq \dfrac{2}{5}$,\vspace*{0.8ex}
\item\label{EWD_2_2_08_BWf}
$||\,H_{3}\,\psi\,|| = \Big|\,\bigl(H_{3}\,\psi\bigr)\Bigl(\sqrt{3 - \sqrt{6}}\Bigr)\,\Big| 
= \Big|\,\bigl(H_{3}\,\psi\bigr)\bigl(0,74196\bigr)\,\Big|
= 0,55059 \leq \dfrac{5}{9}$,\vspace*{0.8ex}
\item\label{EWD_2_2_08_BWg}
$||\,x\,\bigl(H_{3}\,\psi\bigr)(x)\,|| = 
2\,\psi(1) = \dfrac{\sqrt{2}}{\sqrt{\pi\,e}} = 
0,48394 \leq \dfrac{1}{2}$,\vspace*{0.8ex}
\item\label{EWD_2_2_08_BWh}
$||\,H_{4}\,\psi\,|| = 3\,\psi(0) = \dfrac{3}{\sqrt{2\,\pi}} = 1,19683 \leq \dfrac{6}{5}$,\vspace*{0.8ex}
\item\label{EWD_2_2_08_BWi}
$||\,x\,\bigl(H_{4}\,\psi\bigr)(x)\,|| = 
\Big|\,1,49882\,\bigl(H_{4}\,\psi\bigr)\bigl(1,49882\bigr)\,\Big| 
= 1,05638 \leq \dfrac{18}{17}$,\vspace*{1.1ex}
\item\label{EWD_2_2_08_BWj}
$||\,H_{5}\,\psi\,|| = \Big|\,\bigl(H_{5}\,\psi\bigr)\bigl(0,61671\bigr)\,\Big| = 
2,30711 \leq \dfrac{7}{3}$,\vspace*{1.9ex}
\item\label{EWD_2_2_08_BWk}
$||\,x\,(H_{5}\,\psi)(x)\,|| = ||\,1,97963\,(H_{5}\,\psi)(1,97963)\,|| = 1,94588 \leq 2$,\vspace*{1.9ex}
\item\label{EWD_2_2_08_BWl}
$||\,H_{6}\,\psi\,|| = 15\,\psi(0) = \dfrac{15}{\sqrt{2\,\pi}} = 5,98413 \leq 6$,\vspace*{0.3ex}
\item\label{EWD_2_2_08_BWm}
$||\,x\,(H_{6}\,\psi)(x)\,|| = ||\,1,27145\,(H_{6}\,\psi)(1,27145)\,|| = 5,14696 \leq \dfrac{26}{5}$,\vspace*{1.7ex}
\item\label{EWD_2_2_08_BWn}
$||\,x\,(H_{7}\,\psi)(x)\,|| = ||\,1,71316\,(H_{7}\,\psi)(1,71316)\,|| = 12,83583 \leq 13$,\vspace*{1.7ex}
\item\label{EWD_2_2_08_BWo}
$||\,H_{8}\,\psi\,|| = 105\,\psi(0) = \dfrac{105}{\sqrt{2\,\pi}} = 41,88894 \leq 42$.
\end{enumerate}\vspace*{3ex}

\section{Introduction of several smooth functions}\label{EWD_Kap2_Sec3}

In chapter \ref{EWD_Kap1} we did not start from the expression $F_{n}(z) - \Phi(z)$ when applying Stein's method, but from $E\bigl(\,q(S_{n})\,\bigr) - \Phi(q)$, where $q$ was a {''smoother''} function than $q = 1_{(- \infty,\, z\,]}$.
Since we have to proceed analogously in the next chapter, where we consider linear rank statistics,
we will introduce the smooth functions\index{function!smooth} required there in this section
(or partially repeat them from the chapter \ref{EWD_Kap1}). Subsequently we prove some properties of these
smooth functions.\\[2.8ex]
For $x, z \in \mathbb{R}$, $k, m \in \mathbb{N}$ and $\lambda > 0$ denote in the following\\[2.5ex]
\refstepcounter{DSBcount}
\label{EWD_2_3_01}
\text{\hspace*{-0.8ex}(\theDSBcount)}
\hspace*{4ex}\index{function!smooth!$p_{z,k}^{\lambda}$}
\begin{tabular}[t]{@{}l@{}}
$\displaystyle{p_{z,0}^{\lambda}(x) =
\left\{
\begin{array}{ll@{}}
1& \hspace*{2ex}
\text{for}\ x \leq z,\\[1.8ex]
1 - \dfrac{1}{\lambda}\,\bigl(\,x-z\,\bigr)& \hspace*{2ex}
\text{for}\ z \leq x \leq z + \lambda,\\[2ex]
0 & \hspace*{2ex}
\text{for}\ z + \lambda \leq x\,.
\end{array}  \right.}$\\[8.5ex]
$\displaystyle{p_{z,k}^{\lambda}(x) = |x|^{k}\,p_{z,0}^{\lambda}(x)}$.
\end{tabular}\\[4ex]
\hspace*{16ex}\begin{tikzpicture}[domain=0:14]
\draw (0,3) -- (8.5, 3);
\draw (0,6) -- (8.5, 6);
\draw[dotted] (2.5,3.2) -- (2.5,6);
\draw (2.5,2.8) -- (2.5,3.2);
\draw (4,2.8) -- (4,3.2);
\draw[red,thick] (0,6) -- (2.5,6); 
\draw[red,thick] (2.5,6) -- (4,3);
\draw[red,thick] (4,3) -- (8.5,3);
\node[] at (-0.2,3) {$0$};
\node[] at (-0.2,6) {$1$};
\node[] at (2.5,2.5) {$z$};
\node[] at (4.03,2.5) {$z + \lambda$};
\node[red] at (4.1,4.5) {$p_{z,0}^{\lambda}$};
\end{tikzpicture}\\[0.5ex]
\refstepcounter{Figcount}
\label{EWDFig_2_3_01}
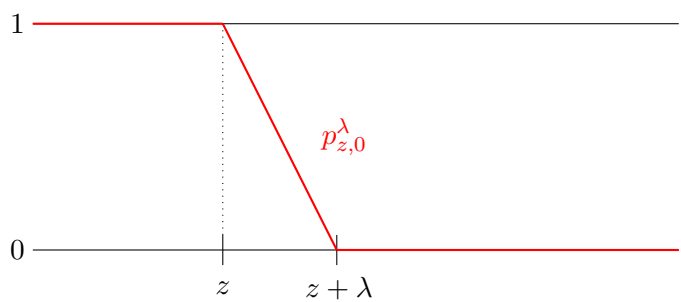
\hspace*{16ex}\textbf{Fig. \theFigcount:} Sketch of the {''linear''} function $p_{z,0}^{\lambda}$\\[3.5ex]
We use the functions $p_{z,k}^{\lambda}$ for 
{''Berry-Ess\'een proofs''}.\\[1.5ex] 
Furthermore,\\[3.5ex]
\refstepcounter{DSBcount}
\label{EWD_2_3_02}
\text{\hspace*{-0.8ex}(\theDSBcount)}
\hspace*{4ex}\index{function!smooth!$q_{z}^{\lambda,m}$}
\begin{tabular}[t]{@{}l@{}}
$\displaystyle{q_{z}^{\lambda,0}(x) =
\left\{
\begin{array}{ll@{}}
1& \hspace*{2ex}
\text{for}\ x \leq z,\\[1.8ex]
1 - \dfrac{1}{2\,\lambda^2}\,\bigl(\,x - z\,\bigr)^2& \hspace*{2ex}
\text{for}\ z \leq x \leq z + \lambda,\\[2.3ex]
\dfrac{1}{2\,\lambda^2}\,\bigl(\,z + 2\,\lambda - x\,\bigr)^2& \hspace*{2ex}
\text{for}\ z + \lambda \leq x \leq z + 2\,\lambda,\\[2.8ex]
0 & \hspace*{2ex}
\text{for}\ z + 2\,\lambda \leq x\,.
\end{array}  \right.}$\\[12ex]
$\displaystyle{q_{z}^{\lambda,m}(x) = x^{m}\,q_{z}^{\lambda,0}(x)}$.\\[3ex]
\end{tabular}\\[4ex]
\hspace*{16ex}\begin{tikzpicture}[domain=0:14]
\draw (0,3) -- (8.5, 3);
\draw (0,6) -- (8.5, 6);
\draw[dotted] (2.5,3.2) -- (2.5,6);
\draw[dotted] (4,3.2) -- (4,4.5);
\draw[dotted] (0,4.5) -- (4,4.5);
\draw (2.5,2.8) -- (2.5,3.2);
\draw (4,2.8) -- (4,3.2);
\draw (5.5,2.8) -- (5.5,3.2);
\draw[red,thick] (0,6) -- (2.5,6); 
\draw[red,thick] (2.5,6) .. controls (3.25,6) and (3.25,6) .. (4,4.5);
\draw[red,thick] (4,4.5) .. controls (4.75,3) and (4.75,3) .. (5.5,3);
\draw[red,thick] (5.5,3) -- (8.5,3);
\node[] at (-0.2,3) {$0$};
\node[] at (-0.2,4.5) {$\dfrac{1}{2}$};
\node[] at (-0.2,6) {$1$};
\node[] at (2.5,2.5) {$z$};
\node[] at (4.03,2.5) {$z + \lambda$};
\node[] at (5.64,2.5) {$z + 2\,\lambda$};
\node[red] at (4.8,4.5) {$q_{z}^{\lambda,0}$};
\end{tikzpicture}\\[0.5ex]
\refstepcounter{Figcount}
\label{EWDFig_2_3_02}
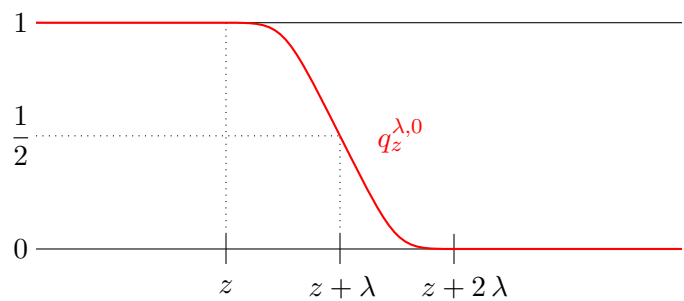
\hspace*{16ex}\textbf{Fig. \theFigcount:} Sketch of the {''quadratic''} funktion $q_{z}^{\lambda,0}$\\[3.5ex]
We need the functions $q_{z}^{\lambda,m}$ to establish first-order Edgeworth expansions.
In the iid case already considered in chapter \ref{EWD_Kap1}, we used $\lambda = \dfrac{1}{\sqrt{n}}$ 
(cf. (\ref{EWD_1_2_01})).\\[1.5ex]
Furthermore,\\[3.5ex]
\refstepcounter{DSBcount}
\label{EWD_2_3_03}
\text{\hspace*{-0.8ex}(\theDSBcount)}
\hspace*{4ex}\index{function!smooth!$r_{z}^{\lambda}$}
\begin{tabular}[t]{@{}l@{}}
$\displaystyle{r_{z}^{\lambda}(x) =
\left\{
\begin{array}{ll@{}}
1& \hspace*{2ex}
\text{for}\ x \leq z,\\[1ex]
1 - \dfrac{1}{6\,\lambda^3}\,\bigl(\,x - z\,\bigr)^3& \hspace*{2ex}
\text{for}\ z \leq x \leq z + \lambda,\\[2.8ex]
\dfrac{1}{24\,\lambda^3}\,\Bigl\{\,\bigl(\,2\,(x - z) - 
3\,\lambda\,\bigr)^3\\[2ex] 
\hspace*{9ex}-\ 18\,\lambda^2\,\bigl(\,x - z\,\bigr) + 39\,\lambda^3\,
\Bigr\}& \hspace*{2ex}
\text{for}\ z + \lambda\leq x \leq z + 2\,\lambda,\\[2.5ex]
\dfrac{1}{6\,\lambda^3}\,\bigl(\,z + 3\,\lambda - x\,\bigr)^3& \hspace*{2ex}
\text{for}\ z + 2\,\lambda \leq x \leq z + 3\,\lambda,\\[2.5ex]
0 & \hspace*{2ex}
\text{for}\ z + 3\,\lambda \leq x\,.
\end{array}  \right.}$
\end{tabular}\\[6ex]
\hspace*{16ex}\begin{tikzpicture}[domain=0:14]
\draw (0,3) -- (10, 3);
\draw (0,6) -- (10, 6);
\draw[dotted] (2.5,3.2) -- (2.5,6);
\draw[dotted] (4,3.2) -- (4,5.5);
\draw[dotted] (4.75,2) -- (4.75,4.5);
\draw[dotted] (5.5,3.2) -- (5.5,3.5);
\draw[dotted] (-0.5,3.5) -- (5.5,3.5);
\draw[dotted] (0,4.5) -- (4.75,4.5);
\draw[dotted] (-0.5,5.5) -- (4,5.5);
\draw (2.5,2.8) -- (2.5,3.2);
\draw (4,2.8) -- (4,3.2);
\draw (5.5,2.8) -- (5.5,3.2);
\draw (7,2.8) -- (7,3.2);
\draw[red,thick] (2.5,6) .. controls (3.4,6) and (3.75,5.71) .. (4,5.5);
\draw[red,thick] (4,5.5) .. controls (4.25,5.29) and (4.7,4.58) .. (4.75,4.5);
\draw[red,thick] (4.75,4.5) .. controls (4.8,4.42) and (5.25,3.71) .. (5.5,3.5);
\draw[red,thick] (5.5,3.5) .. controls (5.75,3.29) and (6.1,3) .. (7,3);
\draw[red,thick] (7,3) -- (10,3);
\draw[red,thick] (0,6) -- (2.5,6); 
\node[] at (-0.2,3) {$0$};
\node[] at (-0.7,3.5) {$\dfrac{1}{6}$};
\node[] at (-0.2,4.5) {$\dfrac{1}{2}$};
\node[] at (-0.7,5.5) {$\dfrac{5}{6}$};
\node[] at (-0.2,6) {$1$};
\node[] at (2.5,2.5) {$z$};
\node[] at (4.03,2.5) {$z + \lambda$};
\node[] at (4.91,1.72) {$z + \frac{3}{2}\,\lambda$};
\node[] at (5.64,2.5) {$z + 2\,\lambda$};
\node[] at (7.14,2.5) {$z + 3\,\lambda$};
\node[red] at (5.3,4.5) {$r_{z}^{\lambda}$};
\end{tikzpicture}\\[0.5ex]
\refstepcounter{Figcount}
\label{EWDFig_2_3_03}
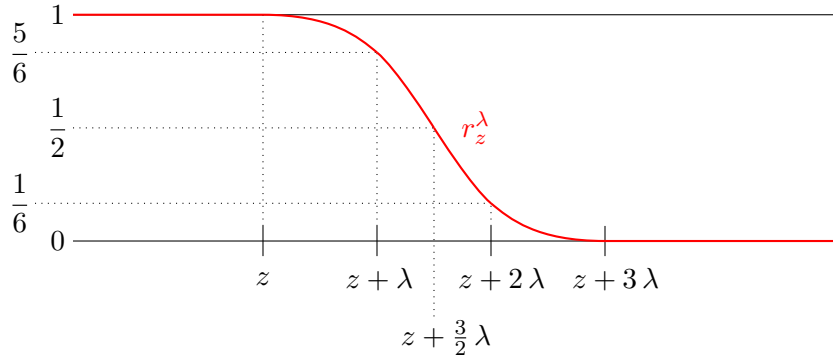
\hspace*{16ex}\textbf{Fig. \theFigcount:} Sketch of the {''cubic''} function $r_{z}^{\lambda}$\\[3.5ex]
We need the functions $r_{z}^{\lambda}$ to establish second-order Edgeworth expansions.\\[1.5ex]
In addition, we define\\[2.8ex]
\refstepcounter{DSBcount}
\label{EWD_2_3_04}
\text{\hspace*{-0.8ex}(\theDSBcount)}
\hspace*{4ex}
\begin{tabular}[t]{@{}l@{\hspace*{0.8ex}}c@{\hspace*{0.8ex}}l@{}}
$d_{z,k}^{\lambda}(x)$&$=$&$\displaystyle{\psi(x)^{-1}\,\int\limits_{-\,\infty}^{x}
\Bigl(\,p_{z,k}^{\lambda}(y) - \Phi\bigl(\,p_{z,k}^{\lambda}\,\bigr)\,\Bigr)\, \psi(y)\,dy}$,\\[3.5ex]
$f_{z}^{\lambda, m}(x)$&$=$&$\displaystyle{\psi(x)^{-1}\,\int\limits_{-\,\infty}^{x}
\Bigl(\,q_{z}^{\lambda,m}(y) - \Phi\bigl(\,q_{z}^{\lambda,m}\,\bigr)\,\Bigr)\, \psi(y)\,dy}$,\\[3.5ex]
$g_{z}^{\lambda}(x)$&$=$&$\displaystyle{\psi(x)^{-1}\,\int\limits_{-\,\infty}^{x}
\Bigl(\,r_{z}^{\lambda}(y) - \Phi\bigl(\,r_{z}^{\lambda}\,\bigr)\,\Bigr)\, \psi(y)\,dy}$.
\end{tabular}\index{function!smooth!$d_{z,k}^{\lambda}$, $f_{z}^{\lambda, m}$, $g_{z}^{\lambda}$}\\[2.5ex]
However, in cases where there is clarity about certain indices, these will be omitted in the future.

\newpage

The following results are of interest for the functions $p_{z,k}^{\lambda}$ and $d_{z,k}^{\lambda}$:\\[4ex]
\refstepcounter{DSBcount}
\label{EWD_2_3_05}
\textbf{\hspace*{-0.8ex}\theDSBcount\ Lemma}\\[0.8ex]
Let $x, y, z \in \mathbb{R}$ and $\lambda > 0$. Then
\begin{enumerate}
\item\label{EWD_2_3_05_BWa}
$\Big|\,p_{z,k}^{\lambda}(x+y) - p_{z,k}^{\lambda}(x) \,\Big|$\\[2ex]
$\displaystyle{\leq \left\{
\begin{array}{ll@{}}
|y|\,\dfrac{1}{\lambda}\,\displaystyle{\int\limits_{0}^{1}\,1_{(\,z,\,z\, + \lambda\,]}(x + s\,y)\,ds}& \hspace*{3ex}
\text{for}\ k = 0,\\[3.5ex]
|y|\,\biggl\{\,\dfrac{|x|^{k}}{\lambda}\,\displaystyle{\int\limits_{0}^{1}\,1_{(\,z,\,z\, + \lambda\,]}(x + s\,y)\,ds +
k\,2^{k-2}\,\Bigl(\,|x|^{k-1} + |y|^{k-1}\,\Bigr)\,\biggr\}}& \hspace*{3ex}
\text{otherwise}.
\end{array}  \right.}$
\item\label{EWD_2_3_05_BWb}
$\Big|\,\bigl(d_{z,0}^{\lambda}\bigr)'(x+y) - (d_{z,0}^{\lambda})'(x) \,\Big|$\\[2ex]
$\displaystyle{\leq\,|y|\,\biggl\{\,\dfrac{\sqrt{2\pi}}{4} + |x| + 
\dfrac{1}{\lambda}\,\int\limits_{0}^{1}\,1_{(\,z,\,z\, + \lambda\,]}(x + s\,y)\,ds\,\biggr\}}$. 
\item\label{EWD_2_3_05_BWc}
$\Big|\,\bigl(d_{z,1}^{\lambda}\bigr)'(x+y) - (d_{z,1}^{\lambda})'(x) \,\Big|$\\[2ex]
$\displaystyle{\leq\,|y|\,\biggl\{\,2 +  
\dfrac{|x|}{\lambda}\,\int\limits_{0}^{1}\,1_{(z,\, z+\lambda\,]}(x + s\,y)\,ds\biggr\}
+ |x|\,|\,d_{z,1}^{\lambda}(x+y) - d_{z,1}^{\lambda}(x)\,|}$\\[2ex]
$\displaystyle{\leq\,|y|\,\biggl\{\,2 + \sqrt{\dfrac{2}{\pi}}\,|x| + x^2 + |x\,y| + 
\dfrac{|x|}{\lambda}\,\int\limits_{0}^{1}\,1_{(\,z,\,z\, + \lambda\,]}(x + s\,y)\,ds\,\biggr\}}$.
\item\label{EWD_2_3_05_BWd}
Let $k \in \mathbb{N}_{0}$. Then there exists a constant $C(k) > 0$ such that\\[2.5ex]
$\Big|\,\bigl(d_{z,k}^{\lambda}\bigr)'(x+y) - (d_{z,k}^{\lambda})'(x) \,\Big|$\\[2ex]
$\displaystyle{\leq\,C(k)\,|y|\,\biggl\{\,1 + |x|^{k+1} + \bigl(\,1 + |x|\,\bigr)\,|y|^{k} + 
\dfrac{|x|^{k}}{\lambda}\,\int\limits_{0}^{1}\,1_{(\,z,\,z\, + \lambda\,]}(x + s\,y)\,ds\,\biggr\}}$.
\end{enumerate}
\vspace*{1.5ex}
\textbf{Proof:}\\[0.8ex]
Part \ref{EWD_2_3_05_BWa}) is obtained by the following calculation:\\[2.5ex]
\hspace*{12.1ex}\begin{tabular}{@{}c@{\hspace*{0.8ex}}l@{}}
&$\Big|\,p_{z,k}^{\lambda}(x+y) - p_{z,k}^{\lambda}(x) \,\Big|$\\[3ex]
$\leq$&$|x|^{k}\,\Big|\,p_{z,0}^{\lambda}(x+y) - p_{z,0}^{\lambda}(x) \,\Big|
+ \Big|\,|x+y|^{k} - |x|^{k}\,\Big|\,\Big|\,p_{z,0}^{\lambda}(x+y)\,\Big|$\\[2.6ex]
$\leq$&$\displaystyle{|x|^{k}\,\dfrac{|y|}{\lambda}\,\int\limits_{0}^{1}\,1_{(\,z,\,z\, + \lambda\,]}(x + s\,y)\,ds
+ \Big|\,|x+y|^{k} - |x|^{k}\,\Big|}$.
\end{tabular}\\[3ex]
In the case of $k \in \mathbb{N}$ the following holds for the second summand\\[2ex]
\hspace*{12.1ex}\begin{tabular}{@{}l@{\hspace*{0.8ex}}c@{\hspace*{0.8ex}}l@{}}
$\Big|\,|x+y|^{k} - |x|^{k}\,\Big|$&$\leq$&$|y|\,\sup\limits_{0 < \theta < 1}\,k\,|x + \theta\,y|^{k-1}$\\[3ex]
&$\leq$&$|y|\,k\,\Bigl(\,|x| + |y|\,\Bigr)^{k-1}$\\[3ex]
&$\leq$&$|y|\,k\,2^{k-2}\,\Bigl(\,|x|^{k-1} + |y|^{k-1}\,\Bigr)$.
\end{tabular}\\[2.5ex]
The last inequality is clear for $k = 1$ and $k = 2$. For $k > 2$ this inequality follows from 
H{\"o}lder's inequality\index{H{\"o}lder's inequality!for finite sequences using length $\nu$} 
(\ref{EWD_0_1_05}) with $\nu = 2$ and $p = k - 1$.\\[2.8ex]
To prove the parts \ref{EWD_2_3_05_BWb}) - \ref{EWD_2_3_05_BWd}), we use Stein's equation\index{Stein's equation}\\[2ex]
\hspace*{12.1ex}$\bigl(d_{z,k}^{\lambda}\bigr)'(x) = x\,d_{z,k}^{\lambda}(x) + p_{z,k}^{\lambda}(x) - 
\Phi\bigl(\,p_{z,k}^{\lambda}\,\bigr)$\\[2.5ex]
from which follows\\[2ex]
\hspace*{12.1ex}\begin{tabular}{@{}c@{\hspace*{0.8ex}}l@{}}
&$\Big|\,\bigl(d_{z,k}^{\lambda}\bigr)'(x+y) - (d_{z,k}^{\lambda})'(x) \,\Big|$\\[3ex]
$\leq$&$|y|\,\Big|\,d_{z,k}^{\lambda}(x+y)\,\Big| +
|x|\,\Big|\,d_{z,k}^{\lambda}(x+y) - d_{z,k}^{\lambda}(x) \,\Big| + 
\Big|\,p_{z,k}^{\lambda}(x+y) - p_{z,k}^{\lambda}(x) \,\Big|$\\[3ex]
$\leq$&$|y|\,\biggl\{\Big|\,d_{z,k}^{\lambda}(x+y)\,\Big| +
|x|\,\sup\limits_{0 < \theta < 1}\,\Big|\,\bigl(d_{z,k}^{\lambda}\bigr)'(x + \theta\,y)\,\Big|\,\biggr\} + 
\Big|\,p_{z,k}^{\lambda}(x+y) - p_{z,k}^{\lambda}(x) \,\Big|$
\end{tabular}\\[3ex]
Using part \ref{EWD_2_3_05_BWa}) and Corollary \ref{EWD_2_1_12} with $|z|^{l} \leq 1 + |z|^{r}$
for $0 < l < r$, $z \in \mathbb{R}$ then yields the assertion.\hfill$\Box$\\[2.8ex]
Next, we examine the functions $q_{z}^{\lambda,m}$ and $f_{z}^{\lambda,m}$.\\[4ex]
\refstepcounter{DSBcount}
\label{EWD_2_3_06}
\textbf{\hspace*{-0.8ex}\theDSBcount\ Lemma}\\[0.8ex]
Let $x, y, z \in \mathbb{R}$ and $\lambda > 0$. Then
\begin{enumerate}
\item\label{EWD_2_3_06_BWa}
$\displaystyle{\Big|\,q_{z}^{\lambda,m}(x+y) - q_{z}^{\lambda,m}(x) \,\Big|\,
\leq\, \left\{
\begin{array}{ll@{}}
|y|\,\dfrac{1}{\lambda}\,\displaystyle{\int\limits_{0}^{1}\,1_{(\,z,\,z\, + 2\,\lambda\,]}(x + s\,y)\,ds}& \hspace*{3ex}
\text{for}\ m = 0,\\[3.5ex]
|y|\,\biggl\{\,1 + \dfrac{|x|}{\lambda}\,\displaystyle{\int\limits_{0}^{1}\,1_{(\,z,\,z\, + 2\,\lambda\,]}(x + s\,y)\,ds
\,\biggr\}}& \hspace*{3ex}
\text{for}\ m = 1.
\end{array}  \right.}$\vspace*{1ex}
\item\label{EWD_2_3_06_BWb}
$\displaystyle{\bigl(q_{z}^{\lambda,0}\bigr)'(x+y) - (q_{z}^{\lambda,0})'(x)
= y\,\dfrac{1}{\lambda^2}\,\int\limits_{0}^{1}\,
\Bigl(\,1_{(\,z\, + \lambda,\,z\, + 2\,\lambda\,]} - 1_{(\,z,\,z\, + \lambda\,]}\,\Bigr)(x + s\,y)\,ds}$.\vspace*{1ex}
\item\label{EWD_2_3_06_BWc}
$\displaystyle{\bigl(q_{z}^{\lambda,1}\bigr)'(x+y) - (q_{z}^{\lambda,1})'(x)
= y\,\dfrac{x}{\lambda^2}\,\int\limits_{0}^{1}\,
\Bigl(\,1_{(\,z\, + \lambda,\,z\, + 2\,\lambda\,]} - 1_{(\,z,\,z\, + \lambda\,]}\,\Bigr)(x + s\,y)\,ds
+ R_{z}^{\lambda,1}(x,y)}$,\\[2ex]
\hspace*{5.6ex}where\hspace*{1.4ex}
$\displaystyle{\big|\,R_{z}^{\lambda,1}(x,y)\,\big|
\leq |y|\,\dfrac{1}{\lambda}\,\biggl\{\,1_{(\,z,\,z\, + 2\,\lambda\,]}(x + y) +
\int\limits_{0}^{1}\,1_{(\,z,\,z\, + 2\,\lambda\,]}(x + s\,y)\,ds
\,\biggr\}}$.\vspace*{1ex}
\item\label{EWD_2_3_06_BWd}
$\displaystyle{\bigg|\,\Bigl(\,\bigl(f_{z}^{\lambda,0}\bigr)'' - \bigl(q_{z}^{\lambda,0}\bigr)'\,\Bigr)(x+y) - 
\Bigl(\,\bigl(f_{z}^{\lambda,0}\bigr)'' - \bigl(q_{z}^{\lambda,0}\bigr)'\,\Bigr)(x)\,\bigg|}$\\[2ex]
$\displaystyle{\leq\,|y|\,\biggl\{\,3 + \dfrac{\sqrt{2\pi}}{4}\,|x| + x^2 
+ \dfrac{|x|}{\lambda}\,\displaystyle{\int\limits_{0}^{1}\,1_{(\,z,\,z\, + 2\,\lambda\,]}(x + s\,y)\,ds}\,\biggr\}}$\\[2ex]
$\displaystyle{\leq\,|y|\,\biggl\{\,\dfrac{11}{3} + \dfrac{5}{3}\,x^2 
+ \dfrac{|x|}{\lambda}\,\displaystyle{\int\limits_{0}^{1}\,1_{(\,z,\,z\, + 2\,\lambda\,]}(x + s\,y)\,ds}\,\biggr\}}$.\vspace*{1ex}
\item\label{EWD_2_3_06_BWe}
$\displaystyle{\bigg|\,\Bigl(\,\bigl(f_{z}^{\lambda,1}\bigr)'' - \bigl(q_{z}^{\lambda,1}\bigr)'\,\Bigr)(x+y) - 
\Bigl(\,\bigl(f_{z}^{\lambda,1}\bigr)'' - \bigl(q_{z}^{\lambda,1}\bigr)'\,\Bigr)(x)\,\bigg|}$\\[2ex]
$\displaystyle{\leq\,|y|\,\biggl\{\,\dfrac{33}{5} + \dfrac{32}{5}\,|x|^3 + (3 + x^2)\,|y| 
+ \dfrac{x^2}{\lambda}\,\displaystyle{\int\limits_{0}^{1}\,1_{(\,z,\,z\, + 2\,\lambda\,]}(x + s\,y)\,ds}\,\biggr\}}$.
\end{enumerate}
\vspace*{1.5ex}
\textbf{Proof:}\\[0.8ex]
By differentiating $q_{z}^{\lambda,0}(x)$ we obtain\\[2ex]
\refstepcounter{DSBcount}
\label{EWD_2_3_07}
\text{\hspace*{-0.8ex}(\theDSBcount)}
\hspace*{4ex}
$\displaystyle{\bigl(q_{z}^{\lambda,0}\bigr)'(x) =
\left\{
\begin{array}{ll@{}}
0& \hspace*{2ex}
\text{for}\ x \leq z\hspace*{1.4ex}\text{or}\hspace*{1.4ex}z + 2\,\lambda \leq x,\\[1.8ex]
-\, \dfrac{1}{\lambda^2}\,\bigl(\,x - z\,\bigr)& \hspace*{2ex}
\text{for}\ z \leq x \leq z + \lambda,\\[2.3ex]
-\, \dfrac{1}{\lambda^2}\,\bigl(\,z + 2\,\lambda - x\,\bigr)& \hspace*{2ex}
\text{for}\ z + \lambda \leq x \leq z + 2\,\lambda,\\[2.8ex]
\end{array}  \right.}$\\[2.5ex]
and therefore\\[2.5ex]
\refstepcounter{DSBcount}
\label{EWD_2_3_08}
\text{\hspace*{-0.8ex}(\theDSBcount)}
\hspace*{4ex}
$\displaystyle{\Big|\,\bigl(q_{z}^{\lambda,0}\bigr)'(x)\,\Big| \leq
\dfrac{1}{\lambda}\,1_{(\,z,\,z\, + 2\,\lambda\,]}(x)}$
\hspace*{4ex}for all $x \in \mathbb{R}$.\\[3ex] 
From this estimate, part \ref{EWD_2_3_06_BWa}) follows analogously to part \ref{EWD_2_3_05_BWa}) 
of the last Lemma \ref{EWD_2_3_05}.\\[2.8ex] 
Further differentiation of $\bigl(q_{z}^{\lambda,0}\bigr)'(x)$ at the points 
$x \not= z, z + \lambda, z + 2\,\lambda$ yields\\[2ex]
\refstepcounter{DSBcount}
\label{EWD_2_3_09}
\text{\hspace*{-0.8ex}(\theDSBcount)}
\hspace*{4ex}
$\displaystyle{\bigl(q_{z}^{\lambda,0}\bigr)''(x) =
\left\{
\begin{array}{ll@{}}
0& \hspace*{2ex}
\text{for}\ x < z\hspace*{1.4ex}\text{or}\hspace*{1.4ex}z + 2\,\lambda < x,\\[1.8ex]
-\, \dfrac{1}{\lambda^2}& \hspace*{2ex}
\text{for}\ z < x < z + \lambda,\\[2.3ex]
\dfrac{1}{\lambda^2}& \hspace*{2ex}
\text{for}\ z + \lambda < x < z + 2\,\lambda,\\[2.8ex]
\end{array}  \right.}$\\[2.5ex]
from which we can immediately see the validity of part \ref{EWD_2_3_06_BWb}).\\[2.8ex]
In addition, the product rule gives us\\[2ex] 
\hspace*{12.1ex}\begin{tabular}{@{}c@{\hspace*{0.8ex}}l@{}}
&$\bigl(q_{z}^{\lambda,1}\bigr)'(x+y) - (q_{z}^{\lambda,1})'(x)$\\[2.5ex]
$=$&$\displaystyle{x\,\Bigl[\,\bigl(q_{z}^{\lambda,0}\bigr)'(x+y) - (q_{z}^{\lambda,0})'(x)\,\Bigr] 
+  R_{z}^{\lambda,1}(x,y)}$,
\end{tabular}\\[2.5ex]
\hspace*{12.1ex}where\hspace*{1.4ex}
$\displaystyle{R_{z}^{\lambda,1}(x,y) = y\,\bigl(q_{z}^{\lambda,0}\bigr)'(x+y) + 
\Bigl[\,q_{z}^{\lambda,0}(x+y) - q_{z}^{\lambda,0}(x)\,\Bigr]}$.\\[2.5ex] 
If we now use the parts \ref{EWD_2_3_06_BWb}) and \ref{EWD_2_3_06_BWa}), and (\ref{EWD_2_3_08}), 
we get part \ref{EWD_2_3_06_BWc}).\\[2.8ex]
To derive the parts \ref{EWD_2_3_06_BWd}) and \ref{EWD_2_3_06_BWe}), 
Stein's equation\index{Stein's equation}\\[2ex]
\refstepcounter{DSBcount}
\label{EWD_2_3_10}
\text{\hspace*{-0.8ex}(\theDSBcount)}
\hspace*{2.8ex}
$f'(x) = x\,f(x) + q(x) - \Phi(q)$\\[2.5ex]
is used (with the abbreviations $f = f_{z}^{\lambda,m}$ and $q = q_{z}^{\lambda,m}$). 
Differentiation of (\ref{EWD_2_3_10}) gives\\[2ex]
\refstepcounter{DSBcount}
\label{EWD_2_3_11}
\text{\hspace*{-0.8ex}(\theDSBcount)}
\hspace*{2.8ex}
\begin{tabular}[t]{@{}l@{\hspace*{0.8ex}}c@{\hspace*{0.8ex}}l@{}}
$f''(x)$&$=$&$f(x) + x\,f'(x) + q'(x)$\\[2.5ex]
&$=$&$\bigl(\,1 + x^2\,\bigr)\,f(x) + x\,\bigl(\,q(x) - \Phi(q)\,\bigr) + q'(x)$
\end{tabular}\\[2.5ex]
and therefore\\[2ex]
\refstepcounter{DSBcount}
\label{EWD_2_3_12}
\text{\hspace*{-0.8ex}(\theDSBcount)}
\hspace*{2.8ex}
\begin{tabular}[t]{@{}c@{\hspace*{0.8ex}}l@{}}
&$\displaystyle{\Big|\,\bigl(\,f'' - q'\,\bigr)(x+y) - \bigl(\,f'' - q'\,\bigr)(x)\,\Big|}$\\[2.8ex]
$\leq$&$\displaystyle{\bigl(\,1 + x^2\,\bigr)\,\big|\,f(x+y) - f(x)\,\big| +
\big|\,2\,x\,y + y^2\,\big|\,\big|\,f(x+y)\,\big| + |y|\,\big|\,q(x+y) - \Phi(q)\,\big|}$
\end{tabular}\\[2.5ex]
\hspace*{12.1ex}\begin{tabular}[t]{@{}c@{\hspace*{0.8ex}}l@{}}
&\hspace*{5.5ex}$\displaystyle{+\ |x|\,\Big|\,\bigl(\,q(x+y) - \Phi(q)\,\bigr) - 
\bigl(\,q(x) - \Phi(q)\,\bigr)\,\Big|}$\\[2.5ex]
$\leq$&$\displaystyle{|y|\,\biggl\{\,\bigl(\,1 + x^2\,\bigr)\,
\sup\limits_{0 < \theta < 1}\,\big|\,f'(x + \theta\,y)\,\big| + 
\big|\,(x+y)\,f(x+y)\,\big| + |x|\,\big|\,f(x+y)\,\big|}$\\[2.5ex]
&\hspace*{5.5ex}$\displaystyle{+\ \big|\,q(x+y) - \Phi(q)\,\big|\,\biggr\} + |x|\,\big|\,q(x+y) - q(x)\,\big|}$.
\end{tabular}\\[2.8ex]
In the case $m = 0$ we now obtain from (\ref{EWD_2_3_12}) because of 
Corollary \ref{EWD_2_1_12}, \ref{EWD_2_1_12_BWa}) and 
Corollary \ref{EWD_2_1_17}, \ref{EWD_2_1_17_BWa})\\[2ex]
\hspace*{12.1ex}\begin{tabular}{@{}c@{\hspace*{0.8ex}}l@{}}
$\leq$&$\displaystyle{|y|\,\biggl\{\,\bigl(\,1 + x^2\,\bigr)\,1 + 1 + \dfrac{\sqrt{2\pi}}{4}\,|x| + 1\,\,\biggr\} 
+ |x|\,\big|\,q(x+y) - q(x)\,\big|}$.
\end{tabular}\\[2.5ex]
Using part \ref{EWD_2_3_06_BWa}) and the inequality $|x| \leq 1 + x^2$, 
part \ref{EWD_2_3_06_BWd}) follows from this.\\[2.8ex]
In the case $m = 1$, on the other hand, we get from (\ref{EWD_2_3_12}) because of 
Corollary \ref{EWD_2_1_17}, \ref{EWD_2_1_17_BWb}) and \ref{EWD_2_1_17_BWc}), and (\ref{EWD_2_1_07})\\[2ex]
\hspace*{12.1ex}\begin{tabular}{@{}c@{\hspace*{0.8ex}}l@{}}
$\leq$&$\displaystyle{|y|\,\biggl\{\,\bigl(\,1 + x^2\,\bigr)\,\bigl(\,|x| + |y| + \dfrac{1}{\sqrt{2\,\pi}}\,\bigr) 
+ \bigl(\,|x| + |y|\,\bigr)\,1 + |x|\,1 + \bigl(\,|x| + |y|\,\bigr) + \sqrt{\dfrac{2}{\pi}}\,\biggr\}}$\\[2.5ex]
&\hspace*{5.5ex}$\displaystyle{+\ |x|\,\big|\,q(x+y) - q(x)\,\big|}$\\[2.5ex]
$\leq$&$\displaystyle{|y|\,\biggl\{\,\dfrac{6}{5} + 4\,|x| + |x|^3 + \dfrac{2}{5}\,x^2 + 3\,|y| + x^2\,|y|\,\biggr\} 
+ |x|\,\big|\,q(x+y) - q(x)\,\big|}$.
\end{tabular}\\[3ex]
Using part \ref{EWD_2_3_06_BWa}) and the inequalities $|x| \leq 1 + |x|^3$ and $x^2 \leq 1 + |x|^3$, 
part \ref{EWD_2_3_06_BWe}) follows from this.\hfill$\Box$\\[4ex]
Finally, we look at the functions $r_{z}^{\lambda}$ and $g_{z}^{\lambda}$.\\[4ex]
\refstepcounter{DSBcount}
\label{EWD_2_3_13}
\textbf{\hspace*{-0.8ex}\theDSBcount\ Lemma}\\[0.8ex]
Let $x, y, z \in \mathbb{R}$ and $\lambda > 0$. Then
\begin{enumerate}
\item\label{EWD_2_3_13_BWa}
$\displaystyle{\Big|\,r_{z}^{\lambda}(x+y) - r_{z}^{\lambda}(x) \,\Big|\,
\leq\, \dfrac{3}{4}\,|y|\,\dfrac{1}{\lambda}\,
\int\limits_{0}^{1}\,1_{(\,z,\,z\, + 3\,\lambda\,]}(x + s\,y)\,ds}$.\vspace*{1ex}
\item\label{EWD_2_3_13_BWb}
$\displaystyle{\bigl(r_{z}^{\lambda}\bigr)''(x+y) - (r_{z}^{\lambda})''(x)}$\\[2ex]
$\displaystyle{=\,y\,\dfrac{1}{\lambda^3}\,\int\limits_{0}^{1}\,
\Bigl(\,- 1_{(\,z\, + 2\,\lambda,\,z\, + 3\,\lambda\,]} + 2 \cdot
1_{(\,z\, + \lambda,\,z\, + 2\,\lambda\,]} - 1_{(\,z,\,z\, + \lambda\,]}\,\Bigr)(x + s\,y)\,ds}$.\vspace*{1ex}
\item\label{EWD_2_3_13_BWc}
$\displaystyle{\int\limits_{z}^{z + 3\,\lambda}\,\bigl(r_{z}^{\lambda}\bigr)'(x)\,dx = -1}$,\hspace*{4ex}    
$\displaystyle{\int\limits_{z}^{z + 3\,\lambda}\,\dfrac{x - z}{\lambda}\,\bigl(r_{z}^{\lambda}\bigr)'(x)\,dx 
= -\dfrac{3}{2}}$,\\[2ex]  
$\displaystyle{\int\limits_{z}^{z + 3\,\lambda}\,\dfrac{(x - z)\,(x - z - \lambda)}
{2\,\lambda^2}\,\bigl(r_{z}^{\lambda}\bigr)'(x)\,dx 
= -\dfrac{1}{2}}$.\vspace*{1ex}    
\item\label{EWD_2_3_13_BWd}
$\bigl(r_{z}^{\lambda}\bigr)'(x) =  \dfrac{3}{4}\,\dfrac{1}{\lambda}\,\bigl(\,h_{1}(x) - h_{2}(x)\,\bigr)$
\hspace*{3ex}with $h_{1}, h_{2} \in \mathcal{H}$.\vspace*{1ex}           
\item\label{EWD_2_3_13_BWe}
$\displaystyle{\bigg|\,\Bigl(\,\bigl(g_{z}^{\lambda}\bigr)''' - \bigl(r_{z}^{\lambda}\bigr)''\,\Bigr)(x+y) - 
\Bigl(\,\bigl(g_{z}^{\lambda}\bigr)''' - \bigl(r_{z}^{\lambda}\bigr)''\,\Bigr)(x)
- x\,\Bigl(\,\bigl(r_{z}^{\lambda}\bigr)'(x+y) - \bigl(r_{z}^{\lambda}\bigr)'(x)\,\Bigr)\,\bigg|}$\\[2ex]
\hspace*{-3.55ex}$\displaystyle{\leq\,|y|\,\biggl\{\,\dfrac{35}{3} + \dfrac{32}{3}\,|x|^3 + 3\,|y| + 
\dfrac{3}{4}\,\dfrac{1}{\lambda}\,1_{(\,z,\,z\, + 3\,\lambda\,]}(x + y) + \dfrac{3}{4}\,(x^2 + 2)\,
\dfrac{1}{\lambda}\,\displaystyle{\int\limits_{0}^{1}\,1_{(\,z,\,z\, + 3\,\lambda\,]}(x + s\,y)\,ds}\,\biggr\}}$.
\end{enumerate}
\vspace*{1.5ex}
\textbf{Proof:}\\[0.8ex]
We obtain by differentiation (cf. some parallels to the proof of Lemma \ref{EWD_2_3_06}):\\[2.5ex]
\refstepcounter{DSBcount}
\label{EWD_2_3_14}
\text{\hspace*{-0.8ex}(\theDSBcount)}
\hspace*{4ex}
$\displaystyle{\bigl(r_{z}^{\lambda}\bigr)'(x) =
\left\{
\begin{array}{ll@{}}
0& \hspace*{2ex}
\text{for}\ x \leq z\hspace*{1.4ex}\text{or}\hspace*{1.4ex}z + 3\,\lambda \leq x,\\[1.8ex]
-\, \dfrac{1}{2\,\lambda^3}\,\bigl(\,x - z\,\bigr)^2& \hspace*{2ex}
\text{for}\ z \leq x \leq z + \lambda,\\[2.8ex]
\dfrac{1}{4\,\lambda^3}\,\Bigl\{\,\bigl(\,2\,(x - z) - 
3\,\lambda\,\bigr)^2 - 3\,\lambda^2\,\Bigr\}& \hspace*{2ex}
\text{for}\ z + \lambda \leq x \leq z + 2\,\lambda,\\[2.8ex]
-\, \dfrac{1}{2\,\lambda^3}\,\bigl(\,z + 3\,\lambda - x\,\bigr)^2& \hspace*{2ex}
\text{for}\ z + 2\,\lambda \leq x \leq z + 3\,\lambda,\\[2.8ex]
\end{array}  \right.}$\\[2.5ex]
and\\[2ex]
\refstepcounter{DSBcount}
\label{EWD_2_3_15}
\text{\hspace*{-0.8ex}(\theDSBcount)}
\hspace*{4ex}
$\displaystyle{\bigl(r_{z}^{\lambda}\bigr)''(x) =
\left\{
\begin{array}{ll@{}}
0& \hspace*{2ex}
\text{for}\ x \leq z\hspace*{1.4ex}\text{or}\hspace*{1.4ex}z + 3\,\lambda \leq x,\\[1.8ex]
-\, \dfrac{1}{\lambda^3}\,\bigl(\,x - z\,\bigr)& \hspace*{2ex}
\text{for}\ z \leq x \leq z + \lambda,\\[2.8ex]
\dfrac{1}{\lambda^3}\,\bigl(\,2\,(x - z) - 
3\,\lambda\,\bigr)& \hspace*{2ex}
\text{for}\ z + \lambda \leq x \leq z + 2\,\lambda,\\[2.8ex]
\dfrac{1}{\lambda^3}\,\bigl(\,z + 3\,\lambda - x\,\bigr)& \hspace*{2ex}
\text{for}\ z + 2\,\lambda \leq x \leq z + 3\,\lambda.\\[2.8ex]
\end{array}  \right.}$\\[2.5ex]
In particular, we have\\[2.5ex]
\refstepcounter{DSBcount}
\label{EWD_2_3_16}
\text{\hspace*{-0.8ex}(\theDSBcount)}
\hspace*{4ex}
$\bigl(r_{z}^{\lambda}\bigr)'\bigl(\,z + \dfrac{3}{2}\,\lambda) = - \dfrac{3}{4}\,
\dfrac{1}{\lambda} \leq \bigl(r_{z}^{\lambda}\bigr)'(x) \leq 0$
\hspace*{4ex}for all $x \in [\,z,\,z\, + 3\,\lambda\,]$,\\[2.5ex]
\refstepcounter{DSBcount}
\label{EWD_2_3_17}
\text{\hspace*{-0.8ex}(\theDSBcount)}
\hspace*{4ex}
\begin{tabular}{@{}l@{}}
$\bigl(r_{z}^{\lambda}\bigr)''(x) \leq 0$
\hspace*{4ex}for all $x \in [\,z,\,z\, + \dfrac{3}{2}\,\lambda\,]$\hspace*{4ex}and\\[2ex]
$\bigl(r_{z}^{\lambda}\bigr)''(x) \geq 0$
\hspace*{4ex}for all $x \in [\,z\, + \dfrac{3}{2}\,\lambda,\,z\, + 3\,\lambda\,]$.
\end{tabular}\\[3ex]
From (\ref{EWD_2_3_16}) and (\ref{EWD_2_3_18}) we immediately get\\[2.5ex]
\refstepcounter{DSBcount}
\label{EWD_2_3_18}
\text{\hspace*{-0.8ex}(\theDSBcount)}
\hspace*{4ex}
$\displaystyle{\Big|\,\bigl(r_{z}^{\lambda}\bigr)'(x)\,\Big| \leq \dfrac{3}{4}\,
\dfrac{1}{\lambda}\,1_{(\,z,\,z\, + 3\,\lambda\,]}(x)}$
\hspace*{4ex}for all $x \in \mathbb{R}$\\[3ex] 
and thus part \ref{EWD_2_3_13_BWa}).\\[2.8ex]
Further differentiation of $\bigl(r_{z}^{\lambda}\bigr)''(x)$ at the points $x \not= z, z + \lambda,
z + 2\,\lambda, z + 3\,\lambda$ gives\\[2ex]
\refstepcounter{DSBcount}
\label{EWD_2_3_19}
\text{\hspace*{-0.8ex}(\theDSBcount)}
\hspace*{4ex}
$\displaystyle{\bigl(r_{z}^{\lambda}\bigr)'''(x) =
\left\{
\begin{array}{ll@{}}
0& \hspace*{2ex}
\text{for}\ x < z\hspace*{1.4ex}\text{or}\hspace*{1.4ex}z + 3\,\lambda < x,\\[1.8ex]
-\, \dfrac{1}{\lambda^3}& \hspace*{2ex}
\text{for}\ z < x < z + \lambda,\\[2.3ex]
\dfrac{2}{\lambda^3}& \hspace*{2ex}
\text{for}\ z + \lambda < x < z + 2\,\lambda,\\[2.8ex]
-\, \dfrac{1}{\lambda^3}& \hspace*{2ex}
\text{for}\ z + 2\,\lambda < x < z + 3\,\lambda,\\[2.8ex]
\end{array}  \right.}$\\[2.5ex]
from which we can at once deduce the validity of part \ref{EWD_2_3_13_BWb}).\\[2.8ex]
Part \ref{EWD_2_3_13_BWc}), on the other hand, is obtained by simple calculations.
Furthermore, part \ref{EWD_2_3_13_BWd}) is valid, since we get because of (\ref{EWD_2_3_16}) 
and (\ref{EWD_2_3_17}):\\[2ex]
\refstepcounter{DSBcount}
\label{EWD_2_3_20}
\text{\hspace*{-0.8ex}(\theDSBcount)}
\hspace*{4ex}
\begin{tabular}{@{}l@{}}
$h_{1}(x) = 1_{(\,- \infty,\,z\, + \frac{3}{2}\,\lambda\,]}(x)\, + \dfrac{4}{3}\,\lambda\,
1_{(\,- \infty,\,z\, + \frac{3}{2}\,\lambda\,]}(x)\,\bigl(r_{z}^{\lambda}\bigr)'(x) \in \mathcal{H}$
\hspace*{4ex}and\\[3ex]
$h_{2}(x) = 1_{(\,- \infty,\,z\, + \frac{3}{2}\,\lambda\,]}(x)\, - \dfrac{4}{3}\,\lambda\,
1_{(\,z\, + \frac{3}{2}\,\lambda,\,+\infty\,)}(x)\,\bigl(r_{z}^{\lambda}\bigr)'(x) \in \mathcal{H}$.
\end{tabular}\\[2.8ex] 
To derive part \ref{EWD_2_3_13_BWe}), we start as in the proof of the 
last Lemma \ref{EWD_2_3_06} (cf. (\ref{EWD_2_3_12})) from Stein's equation\index{Stein's equation}\\[2ex]
\refstepcounter{DSBcount}
\label{EWD_2_3_21}
\text{\hspace*{-0.8ex}(\theDSBcount)}
\hspace*{4ex}
$g'(x) = x\,g(x) + r(x) - \Phi(r)$\\[2.5ex]
(with the abbreviations $g = g_{z}^{\lambda}$ and $r = r_{z}^{\lambda}$).
By differentiating both sides first follows\\[2ex]
\refstepcounter{DSBcount}
\label{EWD_2_3_22}
\text{\hspace*{-0.8ex}(\theDSBcount)}
\hspace*{4ex}
$g''(x) = g(x) + x\,g'(x) + r'(x)$.\\[2.5ex]
Differentiating again and then substituting (\ref{EWD_2_3_22}) and (\ref{EWD_2_3_21}) yields\\[2.5ex]
\refstepcounter{DSBcount}
\label{EWD_2_3_23}
\text{\hspace*{-0.8ex}(\theDSBcount)}
\hspace*{4ex}
\begin{tabular}[t]{@{}l@{\hspace*{0.8ex}}c@{\hspace*{0.8ex}}l@{}}
$g'''(x)$&$=$&$2\,g'(x) + x\,g''(x) + r''(x)$\\[2ex]
&$=$&$2\,g'(x) + x\,g(x) + x^2\,g'(x) + x\,r'(x) + r''(x)$\\[2ex]
&$=$&$r''(x) + x\,r'(x) + (x^3 + 3\,x)\,g(x) + (x^2 + 2)\,\bigl(\,r(x) - \Phi(r)\,\bigr)$.
\end{tabular}\\[2.5ex]
This gives us\\[3ex]
\hspace*{12.1ex}\begin{tabular}{@{}c@{\hspace*{0.8ex}}l@{}}
&$\displaystyle{\Big|\,\bigl(g''' - r''\,\bigr)(x+y) - 
\bigl(\,g''' - r''\,\bigr)(x)
- x\,\bigl(\,r'(x+y) - r'(x)\,\bigr)\,\Big|}$\\[2.5ex]
$\leq$&$|y|\,\big|\,r'(x+y)\,\big| + \big|\,(x+y)^3\,g(x+y) - x^3\,g(x)\,\big|
+ 3\,\big|\,(x+y)\,g(x+y) - x\,g(x)\,\big|$\\[2.5ex]
&$+\,\big|\,(x+y)^2\,r(x+y) - x^2\,r(x)\,\big| + 2\,\big|\,r(x+y) - r(x)\,\big|
+ \big|\,(x+y)^2 - x^2\,\big|\,\Phi(r)$.
\end{tabular}\\[2.5ex]
The various summands can now be estimated using (\ref{EWD_2_3_18}) and Corollary \ref{EWD_2_1_12}, \ref{EWD_2_1_12_BWa}) 
and Corollary \ref{EWD_2_1_17}, \ref{EWD_2_1_17_BWa}) and part \ref{EWD_2_3_13_BWa}):
\renewcommand{\theenumi}{\arabic{enumi}}
\begin{enumerate}
\item
$\displaystyle{|y|\,\big|\,r'(x + y)\,\big| \leq \dfrac{3}{4}\,|y|\,
\dfrac{1}{\lambda}\,1_{(\,z,\,z\, + 3\,\lambda\,]}(x + y)}$,
\item
\begin{tabular}[t]{@{}c@{\hspace*{0.8ex}}l@{}}
&$\big|\,(x+y)^3\,g(x+y) - x^3\,g(x)\,\big|$\\[2ex]
$\leq$&$|x|^3\,\big|\,g(x+y) - g(x)\,\big| + \big|\,\bigl(\,(x+y)^3 - x^3\,\bigr)\,g(x+y)\,\big|$\\[2ex]
$\leq$&$\displaystyle{|y|\,\biggl\{\,|x|^3\,
\sup\limits_{0 < \theta < 1}\,\big|\,g'(x + \theta\,y)\,\big| 
+ \big|\,(2\,x + y)\,(x + y)\,g(x+y)\,\big| + |\,x^2\,g(x+y)\,|\,\biggr\}}$\\[2ex]
$\leq$&$\displaystyle{|y|\,\biggl\{\,|x|^3 + 2\,|x| + |y| + \dfrac{\sqrt{2\pi}}{4}\,x^2\,\biggr\}}$,
\end{tabular}
\item
\begin{tabular}[t]{@{}c@{\hspace*{0.8ex}}l@{}}
&$3\,\big|\,(x+y)\,g(x+y) - x\,g(x)\,\big|$\\[2ex]
$\leq$&$3\,|x|\,\big|\,g(x+y) - g(x)\,\big| + 3\,\big|\,y\,g(x+y)\,\big|$\\[2ex]
$\leq$&$\displaystyle{|y|\,\biggl\{\,3\,|x|\,
\sup\limits_{0 < \theta < 1}\,\big|\,g'(x + \theta\,y)\,\big| 
+ 3\,\big|\,g(x+y)\,|\,\biggr\}}$
\end{tabular}\\[2ex]
\begin{tabular}[t]{@{}c@{\hspace*{0.8ex}}l@{}}
$\leq$&$\displaystyle{|y|\,\biggl\{\,3\,|x|\,
+ 3\,\dfrac{\sqrt{2\pi}}{4}\,\biggr\}}$,
\end{tabular}
\item
\begin{tabular}[t]{@{}c@{\hspace*{0.8ex}}l@{}}
&$\big|\,(x+y)^2\,r(x+y) - x^2\,r(x)\,\big|$\\[2ex]
$\leq$&$x^2\,\big|\,r(x+y) - r(x)\,\big| + \bigl(\,2\,|x|\,|y| + y^2\,\bigr)\,r(x+y)$\\[2ex]
$\leq$&$\displaystyle{|y|\,\biggl\{\dfrac{3}{4}\,x^2\,\dfrac{1}{\lambda}\,
\int\limits_{0}^{1}\,1_{(\,z,\,z\, + 3\,\lambda\,]}(x + s\,y)\,ds + 2\,|x| + |y|\,\biggr\}}$,
\end{tabular}
\item
$\displaystyle{2\,\big|\,r(x+y) - r(x)\,\big| \leq \dfrac{3}{2}\,|y|\,\dfrac{1}{\lambda}\,
\int\limits_{0}^{1}\,1_{(\,z,\,z\, + 3\,\lambda\,]}(x + s\,y)\,ds}$,
\item
\begin{tabular}[t]{@{}l@{\hspace*{0.8ex}}c@{\hspace*{0.8ex}}l@{}}
$\big|\,(x+y)^2 - x^2\,\big|\,\Phi(r)$&$\leq$&$\big|\,2\,x\,y + y^2\,\big|\,1
\leq |y|\,\biggl\{\,2\,|x| + |y|\,\biggr\}$.
\end{tabular}
\end{enumerate}\vspace*{1ex}
\renewcommand{\theenumi}{\alph{enumi}}
By adding the various estimates and using $|x| \leq 1 + |x|^3$ and $x^2 \leq 1 + |x|^3$,
the assertion follows.
\hspace*{1ex}\hfill$\Box$\vspace*{2ex}

\section{Some results on multiple differences and interpolating polynomials}\label{EWD_Kap2_Sec4}

In this section, some simple results on multiple differences\index{difference!multiple} 
and interpolating polynomials\index{polynomial!interpolating}
are summarized. It should be noted that these results are partly taken from the paper of Bickel and
Robinson \cite{10.1214/aop/1176993873}.\\[4ex]
\refstepcounter{DSBcount}
\label{EWD_2_4_01}
\textbf{\hspace*{-0.8ex}\theDSBcount\ Lemma}\\[0.8ex]
Let $k \in \mathbb{N}$, $\lambda > 0$ and $F\,:\,\mathbb{R} \rightarrow \mathbb{R}$ be a function 
such that there exists a constant $C > 0$ so that for all $z \in \mathbb{R}$ and
$z \leq x \leq z + (k+1)\,\lambda$:\\[2ex]
\hspace*{12.1ex}$\Big|\,F(x) - P_{\lambda}^{k}(x;z,F)\,\Big|\, \leq\, 
C\,\Bigl(\,\lambda^{k+1} + (x-z)^{k+1}\,\Bigr)$.\\[2.5ex]
Then follows for all $0 \leq y \leq \lambda$:\\[2ex]
\hspace*{12.1ex}$\big|\big|\,\Delta_{y}^{k + 1}F\,\big|\big|\, \leq\, 
D\,\Bigl(\,\lambda^{k+1} + y^{k+1}\,\Bigr)$,\\[2.5ex]
where
\hspace*{2ex}$\displaystyle{D\,=\,C\,\Bigl(\, 1 + \sum\limits_{j = 1}^{k+1}\,\dbinom{k+1}{j}\,j^{k+1}\,\Bigr)}$\hspace*{2ex}
can be selected.\\[4ex]
\textbf{Proof:}\\[0.8ex]
\begin{tabular}[t]{@{}l@{\hspace*{0.8ex}}c@{\hspace*{0.8ex}}l@{\hspace*{-11.8ex}}r@{}}
$\Big|\,\Delta_{y}^{k+1}F(z)\,\Big|$&$=$&
$\displaystyle{\Big|\,\sum\limits_{j = 0}^{k+1}\,(-1)^{k+1-j}\, \dbinom{k+1}{j}\, 
F(z + j\,y)\,\Big|}$&(cf. (\ref{EWD_0_1_03}))\\[4ex]
&$\leq$&$\displaystyle{\bigg|\,\sum\limits_{j = 0}^{k+1}\,(-1)^{k+1-j}\, \dbinom{k+1}{j}\, 
P_{\lambda}^{k}(z + j\,y;z,F)\,\bigg| +
\sum\limits_{j = 0}^{k+1}\,\dbinom{k+1}{j}\,
C\,\Bigl(\,\lambda^{k+1} + (j\,y)^{k+1}\,\Bigr)}$\\[4ex]
&$\leq$&
$\displaystyle{\bigg|\,\sum\limits_{j = 0}^{k+1}\,(-1)^{k+1-j}\, \dbinom{k+1}{j}\,
\Bigl(\,F(z)\, +\, \sum\limits_{s = 1}^{k}\,\dfrac{1}{s!}\,\Delta_{\lambda}^{s}F(z)\, \prod\limits_{i = 1}^{s}\,
\bigl(\,j\,\dfrac{y}{\lambda} - i + 1\,\bigr)\,\Bigr)\,\bigg|}$
\end{tabular}\\[4ex]
\hspace*{14.5ex}\begin{tabular}[t]{@{}l@{\hspace*{0.8ex}}c@{\hspace*{0.8ex}}l@{\hspace*{21.6ex}}r@{}}
&&$\displaystyle{+\,C\,\Bigl(\,1 + \sum\limits_{j = 1}^{k+1}\,\dbinom{k+1}{j}\,j^{k+1}\,\Bigr)\,
\Bigl(\,\lambda^{k+1} + y^{k+1}\,\Bigr)}$&(cf. (\ref{EWD_0_1_04})).
\end{tabular}\\[4.5ex]
Since\\[1ex]
\hspace*{12.1ex}$\displaystyle{\sum\limits_{j = 0}^{k+1}\,(-1)^{k+1-j}\, \dbinom{k+1}{j}\,j^l\, =\, 0}$\hspace*{4ex}
for all\hspace*{1.4ex}$l = 0,1,\ldots,k$,\\[2.5ex] 
the first sum disappears, and thus the assertion follows.
\hspace*{1ex}\hfill$\Box$\\[3ex]
Furthermore, for sufficiently smooth functions $F$ we have\\[4ex]
\refstepcounter{DSBcount}
\label{EWD_2_4_02}
\textbf{\hspace*{-0.8ex}\theDSBcount\ Lemma}\\[0.8ex]
Let $k \in \mathbb{N}$, $\lambda > 0$ and $F\,:\,\mathbb{R} \rightarrow \mathbb{R}$ be a function.
\begin{enumerate}
\item\label{EWD_2_4_02_BWa}
If $F$ is $k+1$ times differentiable with $\big|\big|\,F^{(k+1)}\,\big|\big| < \infty$, 
then\\[2ex]
\hspace*{6.3ex}
$\displaystyle{\Big|\,F(x) - P_{\lambda}^{k}(x;z,F)\,\Big|\, \leq\, \dfrac{2^k}{k+1}\,
\big|\big|\,F^{(k+1)}\,\big|\big|\,
\Bigl(\,\lambda^{k+1} + |x-z|^{k+1}\,\Bigr)}$
\hspace*{3ex}for all $z, x \in \mathbb{R}$.
\item\label{EWD_2_4_02_BWb}
If $F$ is $k$ times differentiable with $\big|\big|\,F^{(k)}\,\big|\big| < \infty$, then\\[2ex]
\hspace*{6.3ex}
$\displaystyle{\big|\big|\,\Delta_{y}^{k}F\,\big|\big|\, \leq\,  
\big|\big|\,F^{(k)}\,\big|\big|\,\,\big|\,y\,\big|^k}$
\hspace*{3ex}for all $y \in \mathbb{R}$.
\item\label{EWD_2_4_02_BWc}
If $F$ is twice differentiable with $\big|\big|\,F'\,\big|\big| < \infty$ and 
$\big|\big|\,x\,F''(x)\,\big|\big| < \infty$, then\\[2ex]
\hspace*{6.3ex}
$\displaystyle{\big|\big|\,z\,\Delta_{y}^{2}F(z)\,\big|\big|\, \leq\,  
\Bigl(\,4\,\big|\big|\,F'\,\big|\big| + \big|\big|\,x\,F''(x)\,\big|\big|\,\Bigr)\,\big|\,y\,\big|^2}$
\hspace*{3ex}for all $y \in \mathbb{R}$.
\end{enumerate}
\vspace*{1.5ex}
\textbf{Proof:}
\begin{enumerate}
\item
According to Stoer and Bulirsch \cite{stoer2002introduction}, Theorem (2.1.4.1), page 49, holds:\\[1ex]
There exists a number $\zeta$ from the smallest interval $I$ containing $z$, 
$z + \lambda$, $z + 2\,\lambda$,\ldots, $z + k\,\lambda$ and $x$ such that\\[2ex]
\hspace*{6.3ex}
$\displaystyle{F(x) - P_{\lambda}^{k}(x;z,F) = \dfrac{F^{(k+1)}(\zeta)}{(k+1)!}\,(x-z)\,(x-z-\lambda)
\,\cdots\,(x-z-k\,\lambda)}$.\\[2.5ex]
It follows that\\[2ex]
\hspace*{6.3ex}
\begin{tabular}[t]{@{}l@{\hspace*{0.8ex}}c@{\hspace*{0.8ex}}l@{}}
$\Big|\,F(x) - P_{\lambda}^{k}(x;z,F)\,\Big|$&$\leq$&
$\dfrac{||F^{(k+1)}||}{(k+1)!}\,|x-z|\,|x-z-\lambda|
\,\ldots\,|x-z-k\,\lambda|$\\[3ex]
&$\leq$&
$\dfrac{1}{k+1}\,||F^{(k+1)}||\,\Bigl(\,\lambda + |x-z|\,\Bigr)^{k+1}$
\end{tabular}\\[3ex]
\hspace*{28.3ex}\begin{tabular}[t]{@{}l@{\hspace*{0.8ex}}c@{\hspace*{0.8ex}}l@{}}
&$\leq$&
$\dfrac{2^k}{k+1}\,
||\,F^{(k+1)}\,||\,\Bigl(\,\lambda^{k+1} + |x-z|^{k+1}\,\Bigr)$.
\end{tabular}\\[3ex]
For the last inequality, H{\"o}lder's inequality\index{H{\"o}lder's inequality!for finite sequences using length $\nu$} 
(\ref{EWD_0_1_05}) with $\nu = 2$ and $p = k + 1$ was applied.
\item
For $k = 1$ the assertion follows directly from the mean value theorem\index{mean value theorem}.\\[1ex] 
If we now assume the validity of the assertion for $k - 1$,
we obtain the following by reapplying the mean value theorem\index{mean value theorem} to $\Delta_{y}^{k-1}F$\\[2.5ex] 
\hspace*{6.3ex}
\begin{tabular}[t]{@{}l@{\hspace*{0.8ex}}c@{\hspace*{0.8ex}}l@{}}
$\Big|\,\Delta_{y}^{k}F(z)\,\Big|$&$=$&
$\Big|\,\Delta_{y}^{k-1}F(z+y) - \Delta_{y}^{k-1}F(z)\,\Big|$\\[3ex]
&$\leq$&$\big|\,y\,\big|\,\,\big|\big|\bigl(\,\Delta_{y}^{k-1}F\,\bigr)'\big|\big|$\\[3ex]
&$\leq$&$\big|\,y\,\big|^k\,\,\big|\big|F^{(k)}\big|\big|\,$.
\end{tabular}\\[3ex]
since\hspace*{1.4ex}$\bigl(\,\Delta_{y}^{k-1}F\,\bigr)' = \Delta_{y}^{k-1}(F)'$\hspace*{1.4ex}(c.f. (\ref{EWD_0_1_03})).
\item
We get\\[2ex]
\hspace*{6.3ex}
\begin{tabular}[t]{@{}l@{\hspace*{0.8ex}}c@{\hspace*{0.8ex}}l@{}}
$\big|\,z\,\Delta_{y}^{2}F(z)\,\big|$&$=$&$\big|\,z\,F(z + 2\,y) - 2\,z\,F(z + y) + z\,F(z)\,\big|$\\[3ex]
&$\leq$&$2\,\big|\,y\,\big|\,\big|\,F(z + 2\,y) - F(z + y)\,\big|$\\[3ex]
&&$+\,\,\Big|\,(z + 2\,y)\,F(z + 2\,y) - 2\,(z+y)\,F(z+y) + z\,F(z)\,\Big|$\\[3ex]
&$\leq$&$2\,\big|\,y\,\big|^2\,\big|\big|\,F'\,\big|\big|$ 
\end{tabular}\\[2.5ex] 
\hspace*{18ex}
\begin{tabular}[t]{@{}l@{\hspace*{0.8ex}}c@{\hspace*{0.8ex}}l@{}}
&&$\displaystyle{+\,\,\Big|\,\int\limits_{z}^{z+y}\, \bigl[\,F(s+y) - F(s) + 
(s+y)\,F'(s+y) - s\,F'(s)\,\bigr]\,ds\,\Big|}$\\[4ex]
&$=$&$2\,\big|\,y\,\big|^2\,\big|\big|\,F'\,\big|\big|$\\[2.5ex] 
&&$\displaystyle{+\,\,\Big|\,\int\limits_{z}^{z+y}\,\Bigl(\,\int\limits_{s}^{s+y}\,F'(t)\,dt\,\Bigr)\,ds + 
\int\limits_{z}^{z+y}\,\Bigl(\,\int\limits_{s}^{s+y}\,F'(t) + t\,F''(t)\,dt\,\Bigr)\,ds\,\Big|}$\\[4ex]
&$\leq$&$3\,\big|\,y\,\big|^2\,\big|\big|\,F'\,\big|\big| + 
\big|\,y\,\big|^2\,\Bigl(\,\big|\big|\,F'\,\big|\big| + \big|\big|\,x\,F''(x)\,\big|\big|\,\Bigr)$\\[3ex] 
&$=$&$\big|\,y\,\big|^2\,\Bigl(\,4\,\big|\big|\,F'\,\big|\big| + \big|\big|\,x\,F''(x)\,\big|\big|\,\Bigr)$.
\end{tabular}\\[-2.4ex]
\hspace*{1ex}\hfill$\Box$
\end{enumerate}
                                
\rehead{Edgeworth expansions of linear rank statistics}    
\chapter[Edgeworth expansions for linear rank statistics using Stein's method]{Edgeworth expansions for linear rank statistics using Stein's method}\label{EWD_Kap3}

\section{Descriptions and results}\label{EWD_Kap3_Sec1}

Let $A = (a_{ij})$ be an $n{\times}n-$matrix of real numbers.\index{matrix} 
We first introduce some abbreviations and notations in connection with this matrix.\\[2.8ex]
We define\\[2ex]
\hspace*{12.1ex}$\displaystyle{a_{i\boldsymbol{.}} = \dfrac{1}{n}\,\sum\limits_{j=1}^{n}\,a_{ij}}$,\hspace*{1.4ex}
$\displaystyle{a_{\boldsymbol{.}j} = \dfrac{1}{n}\,\sum\limits_{i=1}^{n}\,a_{ij}}$,\hspace*{1.4ex}
$\displaystyle{a_{\boldsymbol{.}\boldsymbol{.}} = \dfrac{1}{n^2}\,\sum\limits_{i, j=1}^{n}\,a_{ij}}$\\[2.5ex]
and\index{matrix!$a_{i\boldsymbol{.}}$, $a_{\boldsymbol{.}j}$, $a_{\boldsymbol{.}\boldsymbol{.}}$, $\mu_{A}$}\\[2ex]
\hspace*{12.1ex}$\displaystyle{\mu_{A} = n\,a_{\boldsymbol{.}\boldsymbol{.}}}$\hspace*{1ex},\\[2ex] 
\hspace*{12.1ex}$\displaystyle{\sigma_{\!A}^2 =
\left\{
\begin{array}{ll@{}}
0& \hspace*{4ex}
\text{for}\ n = 1,\\[2ex]
\displaystyle{\dfrac{1}{n-1}\,\sum\limits_{i,j=1}^{n}\,\bigl(\,a_{ij} 
- a_{i\boldsymbol{.}} - a_{\boldsymbol{.}j} +  a_{\boldsymbol{.}\boldsymbol{.}}\,\bigr)^2}& \hspace*{4ex}
\text{for}\ n \geq 2.
\end{array}  \right.}$\index{matrix!$\sigma_{\hspace*{-0.3ex}A}$}\\[2.5ex]
In the following, we will only consider matrices with $\sigma_{\!A} > 0$ and thus $n \geq 2$.
This allows us to introduce the corresponding standardized matrix\index{matrix!standardized} 
$\bm\hat{A} = (\bm\hat{a}_{ij})$, whose
coefficients are\index{matrix!standardized!$\bm\hat{A}$}
\index{matrix!standardized!$\bm\hat{a}_{ij}$}\\[2ex]
\hspace*{12.1ex}$\displaystyle{\bm\hat{a}_{ij} = 
\dfrac{1}{\sigma_{\!A}}\,\bigl(\,a_{ij} 
- a_{i\boldsymbol{.}} - a_{\boldsymbol{.}j} +  a_{\boldsymbol{.}\boldsymbol{.}}\,\bigr)}$.\\[2.5ex]
We get\\[2ex]
\refstepcounter{DSBcount}
\label{EWD_3_1_01}
\text{\hspace*{-0.8ex}(\theDSBcount)}
\hspace*{4ex}
$\displaystyle{\bm\hat{a}_{i\boldsymbol{.}} = 
\bm\hat{a}_{\boldsymbol{.}j} = 
\bm\hat{a}_{\boldsymbol{.}\boldsymbol{.}} = 0}$\hspace{1.4ex}for all $i$, $j$\hspace*{1.4ex};\hspace*{1.4ex}
$\mu_{\bm\hat{A}} = 0$\hspace*{1.4ex}and\hspace*{1.4ex}$\sigma_{\!\bm\hat{A}}^2 = 1$.\\[2.5ex]
In addition, we define\index{matrix!moment}\\[2ex]
\hspace*{12.1ex}$\displaystyle{\beta_{A} = \sum\limits_{i,j = 1}^{n}|\,\bm\hat{a}_{ij}\,|^3}$,\hspace*{1.4ex}
$\displaystyle{\delta_{A} = \sum\limits_{i,j = 1}^{n}|\,\bm\hat{a}_{ij}\,|^4}$,\hspace*{1.4ex}
$\displaystyle{\eta_{A} = \sum\limits_{i,j = 1}^{n}|\,\bm\hat{a}_{ij}\,|^5}$,
\index{matrix!moment!$\beta_{A}$, $\delta_{A}$, $\eta_{A}$}\\[2.5ex]
\hspace*{12.1ex}$\displaystyle{D_{\!A} = \biggl(\,\dfrac{\delta_{A}}{n}\,\biggr)^{1/2}}$,\hspace*{1.4ex}
$\displaystyle{E_{\!A} = \biggl(\,\dfrac{\eta_{A}}{n}\,\biggr)^{1/3}}$
\index{matrix!moment!$D_{\hspace*{-0.3ex}A}$, $E_{\hspace*{-0.3ex}A}$}\\[2.5ex]
and\\[2ex]
\hspace*{12.1ex}\begin{tabular}[t]{@{}l@{\hspace*{0.8ex}}c@{\hspace*{0.8ex}}l@{}}
$M(l,A)$&$=$&set of all $(n-l)\,{\times}\,(n-l)\,-\,$matrices, which can be\\[1ex]
&&obtained from $A$ by cancelling $l$ rows and $l$ columns,\\[2ex]
$N(l,A)$&$=$&$\displaystyle{\bigcup\,\Big\{\,M(r,A)\, :\, 0 \leq r \leq \min\{\,l, n - 1\,\}\,\Big\}}$.
\end{tabular}
\index{matrix!submatrix}\index{matrix!submatrix!$M(l,A)$, $N(l,A)$}\\[3ex]
If $\pi$ is uniformly distributed on the set $\mathscr{P}_{n}$ 
\index{permutation sets!$\mathscr{P}_{n}$} of permutations of $\{1,\ldots,n\}$,
we denote by\\[2ex]
\hspace*{12.1ex}$\displaystyle{T_{A} = \sum\limits_{i=1}^n a_{i\pi(i)}.}$\\[2.5ex]
the \textbf{linear rank statistic}\index{rank statistic}\index{rank statistic!linear} for the matrix $A$ 
\index{rank statistic!$T_{A}$, $\mathscr{T}_{A}$} 
(under the null hypothesis $H_{0}$).\index{null hypothesis $H_{0}$}\\[2.8ex]
According to H{\'a}jek, \v{S}id{\'a}k and Sen \cite{sidak1999theory}, Theorem 1 in section 3.3.1, page 57, we have\\[2.5ex]
\refstepcounter{DSBcount}
\label{EWD_3_1_02}
\text{\hspace*{-0.8ex}(\theDSBcount)}
\hspace*{4ex}
$E(T_{A}) = \mu_{A}$\hspace*{1.4ex}and\hspace*{1.4ex}$\text{Var}(T_{A}) = \sigma_{\!A}^2$.\\[2.5ex]
Thus, for the related standardized rank statistic\index{rank statistic!standardized} holds\\[3ex]
\hspace*{12.1ex}$\displaystyle{\mathscr{T}_{A} = \dfrac{T_{A} - E(T_{A})}{\sqrt{\text{Var}(T_{A})}}
= \dfrac{T_{A} - \mu_{A}}{\sigma_{\!A}} = T_{\bm\hat{A}}}$.\\[3ex]
Furthermore, $F_{A}$ and $\mathscr{F}_{\!A}$ define the distribution function\index{rank statistic!distribution function}\index{rank statistic!distribution function!$F_{A}$, $\mathscr{F}_{A}$} of 
$T_{A}$ and $\mathscr{T}_{A}$, respectively.\\[2.8ex]
The starting point for this chapter is the following estimate of the distance of $\mathscr{F}_{\!A}$ to the distribution function $\Phi$ of the standard normal distribution. This result is from Bolthausen
\cite{Bolthausen1984}, Theorem, page 380, and states:\\[4ex]
\refstepcounter{DSBcount}
\label{EWD_3_1_03}
\textbf{\hspace*{-0.8ex}\theDSBcount\ Theorem}\index{Theorem!for linear rank statistics!Bolthausen}\\[0.8ex]
There exists an absolute constant $\mathcal{K}_{1} > 0$ such that for all matrices $A$ satisfying $\sigma_{\!A} > 0$\\[2.5ex]
\refstepcounter{DSBcount}
\label{EWD_3_1_04}
\text{\hspace*{-0.8ex}(\theDSBcount)}
\hspace*{4ex}
$||\,\mathscr{F}_{\!A} - \Phi\,|| \leq \mathcal{K}_{1}\,\dfrac{\beta_{A}}{n}$.\\[4ex]
If $\bigl(A^{(n)}\bigr)_{n\, \geq\, 2}$ is a sequence of $n{\times}n-$matrices, this theorem gives
in many typical applications the convergence rate $||\,\mathscr{F}_{\!A^{(n)}} - \Phi\,|| = 
\mathcal{O}\biggl(\,\dfrac{1}{\sqrt{n}}\,\biggr)$.\\[4ex]
\refstepcounter{DSBcount}
\label{EWD_3_1_05}
\textbf{\hspace*{-0.8ex}\theDSBcount\ Remark}\\[0.8ex]
Three decades after the publication of Theorem \ref{EWD_3_1_03} Chen and Fang \cite{10.3150/13-BEJ569} proved that 
$\mathcal{K}_{1} \leq 451$ can be chosen.
Shortly afterwards Th\`{a}nh \cite{Thanh2013} showed that even $\mathcal{K}_{1} \leq 90$ is possible.\\[3.3ex]
In the following, we will consider first and second order Edgeworth expansions\index{Edgeworth expansion} 
of $\mathscr{F}_{\!A}$ and prove similar results for the distance between $\mathscr{F}_{\!A}$ and these expansions.
We define the first and second order Edgeworth expansions of $\mathscr{F}_{\!A}$ as follows\\[3ex]
\refstepcounter{DSBcount}
\label{EWD_3_1_06}
\text{\hspace*{-0.8ex}(\theDSBcount)}
\hspace*{4ex}\index{Edgeworth expansion!$e_{1,A}$, $e_{2,A}$}
$e_{1,A}(x) = \Phi(x) - \psi(x)\, \dfrac{\lambda_{1,A}}{6}\, (x^2-1)$\hspace*{4ex}and\\[2ex]
\refstepcounter{DSBcount}
\label{EWD_3_1_07}
\text{\hspace*{-0.8ex}(\theDSBcount)}
\hspace*{4ex}
$e_{2,A}(x) = \Phi(x) - \psi(x)\, \biggl\{\, \dfrac{\lambda_{1,A}}{6}\, (x^2-1) + \dfrac{\lambda_{2,A}}{24}\, (x^3-3 x)
+ \dfrac{\lambda_{1,A}^{2}}{72}\, (x^5-10x^3 + 15x)\,\biggr\}$,\\[2.5ex]
where\\[1.5ex]
\refstepcounter{DSBcount}
\label{EWD_3_1_08}
\text{\hspace*{-0.8ex}(\theDSBcount)}
\hspace*{4ex}\index{Edgeworth expansion!$\lambda_{1,A}$, $\lambda_{2,A}$}
$\begin{array}{@{}l@{\hspace*{0.8ex}}c@{\hspace*{0.8ex}}l}
\lambda_{1,A}&=&\displaystyle{\dfrac{1}{n}\sum\limits_{i,j=1}^n \bm\hat{a}_{ij}^{3}}\ \ \ \ \text{and}
\end{array}$\\[2ex]
\refstepcounter{DSBcount}
\label{EWD_3_1_09}
\text{\hspace*{-0.8ex}(\theDSBcount)}
\hspace*{4ex}
$\begin{array}{@{}l@{\hspace*{0.8ex}}c@{\hspace*{0.8ex}}l}
\lambda_{2,A}&=&\displaystyle{\dfrac{1}{n}\sum\limits_{i,j=1}^n \bm\hat{a}_{ij}^{4}\, +\, \dfrac{3}{n}\, -\, 
\dfrac{3}{n^2} \sum\limits_{i,j,k=1}^n \bigl(\,\bm\hat{a}_{ij}^{2} \bm\hat{a}_{ik}^{2} + 
\bm\hat{a}_{ij}^{2} \bm\hat{a}_{kj}^{2}\,\bigr)}\ .
\end{array}$\\[3.5ex]
For $e_{1,A}$ we obtain the following result:\\[3.5ex]
\refstepcounter{DSBcount}
\label{EWD_3_1_10}
\textbf{\hspace*{-0.8ex}\theDSBcount\ Theorem}\index{Theorem!for linear rank statistics!main results}\\[0.8ex]
There exist absolute constants $\mathcal{K}_{2} > 0$ and $\mathcal{K}_{3} > 0$ such that for all matrices $A$ 
satisfying $\sigma_{\!A} > 0$ and the condition\\[2.2ex]
\refstepcounter{DSBcount}
\label{EWD_3_1_11}
\text{\hspace*{-0.8ex}(\theDSBcount)}
\hspace*{2.8ex}
\begin{tabular}[t]{@{}l@{}}
there exists a constant $\mathcal{C}_{1} > 0$ such that\\[1.5ex]
\hspace*{4ex}
$\big|\,\Delta_{y}^{2}F_{B}(z)\,\big| \leq \mathcal{C}_{1}\,\bigl(\,D_{\!A}^2 + y^2\,\bigr)$\\[1.8ex]
for all $z \in \mathbb{R}$, $0 \leq y \leq D_{\!A}$ and $B \in N(8,\bm\hat{A})$,
\end{tabular}\\[0.5ex]
we have\\[1.5ex]
\refstepcounter{DSBcount}
\label{EWD_3_1_12}
\text{\hspace*{-0.8ex}(\theDSBcount)}
\hspace*{2.8ex}
$||\,\mathscr{F}_{\!A} - e_{1,A}\,|| \leq \bigl(\,\mathcal{K}_{2}\,\mathcal{C}_{1} + 
\mathcal{K}_{3}\,\bigr)\,D_{\!A}^2$.\\[4ex]
The corresponding result for $e_{2,A}$ is:\\[3.5ex]
\refstepcounter{DSBcount}
\label{EWD_3_1_13}
\textbf{\hspace*{-0.8ex}\theDSBcount\ Theorem}\index{Theorem!for linear rank statistics!main results}\\[0.8ex]
There exist absolute constants $\mathcal{K}_{4} > 0$, $\mathcal{K}_{5} > 0$ and $\mathcal{K}_{6} > 0$ such that
for all matrices $A$ satisfying $\sigma_{\!A} > 0$ and the conditions\\[2.2ex]
\refstepcounter{DSBcount}
\label{EWD_3_1_14}
\text{\hspace*{-0.8ex}(\theDSBcount)}
\hspace*{2.8ex}
\begin{tabular}[t]{@{}l@{}}
there exists a constant $\mathcal{C}_{2} > 0$ such that\\[1.5ex]
\hspace*{4ex}
$\displaystyle{\Big|\,F_{B}(x) - P_{E_{\!A}}^{2}(x;z,F_{B})\,\Big| \leq 
\mathcal{C}_{2}\, \,\Bigl(\,E_{\!A}^{3} + 
(x - z)^{3}\,\Bigr)}$\\[1.8ex]
for all $z \in \mathbb{R}$, $z \leq x \leq z + 3\,E_{\!A}$ and $B \in N(16,\bm\hat{A})$,
\end{tabular}\\[2.2ex]
\refstepcounter{DSBcount}
\label{EWD_3_1_15}
\text{\hspace*{-0.8ex}(\theDSBcount)}
\hspace*{2.8ex}
\begin{tabular}[t]{@{}l@{}}
there exists a constant $\mathcal{C}_{3} > 0$ such that\\[1.5ex]
\hspace*{4ex}
$(|z| + 1)\,\big|\,\Delta_{y}^{2}F_{B}(z)\,\big| \leq \mathcal{C}_{3}\,\bigl(\,E_{\!A}^2 + y^2\,\bigr)$\\[1.8ex]
for all $z \in \mathbb{R}$, $0 \leq y \leq E_{\!A}$ and $B \in N(16,\bm\hat{A})$,
\end{tabular}\\[0.5ex]
we have\\[1.5ex]
\refstepcounter{DSBcount}
\label{EWD_3_1_16}
\text{\hspace*{-0.8ex}(\theDSBcount)}
\hspace*{2.8ex}
$||\,\mathscr{F}_{\!A} - e_{2,A}\,|| \leq \bigl(\,\mathcal{K}_{4}\,\mathcal{C}_{2} + \mathcal{K}_{5}\,\mathcal{C}_{3} +
\mathcal{K}_{6}\,\bigr)\,E_{\!A}^3$.\\[3.8ex]
\refstepcounter{DSBcount}
\label{EWD_3_1_17}
\textbf{\hspace*{-0.8ex}\theDSBcount\ Remarks}
\begin{enumerate}
\item\label{EWD_3_1_17_BWa}
Although the Theorems \ref{EWD_3_1_10} and \ref{EWD_3_1_13}, as Theorem \ref{EWD_3_1_03} above, are formulated for 
fixed $A$ and $n$, their real purpose is to apply them to sequences of $n{\times}n-$matrices $A^{(n)}$.
One then usually obtains (cf. \ref{EWD_3_1_17_BWb}) below) 
a statement about the rate of convergence of $||\,\mathscr{F}_{\!A^{(n)}} - e_{1,A^{(n)}}\,||$ 
and $||\,\mathscr{F}_{\!A^{(n)}} - e_{2,A^{(n)}}\,||$ to zero.
\item\label{EWD_3_1_17_BWb}
In many typical applications, $D_{\!A}^2$ is of order $\dfrac{1}{n}$ and 
$E_{\!A}^3$ is of order $\dfrac{1}{n^{3/2}}$ (cf. e.g. Lemma \ref{EWD_3_1_18}, \ref{EWD_3_1_18_BWa}) 
and Lemma \ref{EWD_4_3_06}, \ref{EWD_4_3_06_BWc})).\\[1.5ex]
However, there are also sequences $A^{(n)}$ for which $D_{\!A^{(n)}}^2 \not\rightarrow 0$ and 
$E_{\!A^{(n)}}^3 \not\rightarrow 0$. We give a simple example:\\[2ex]
\hspace*{13.2ex}$A^{(n)} = \left(
\begin{array}{rrrcc}
1&-1&\ \ \,0&\cdots&0\\
-1&1&\ \ \,0&\cdots&0\\
1&-1&\ \ \,0&\cdots&0\\
-1&1&\ \ \,0&\cdots&0\\
\vdots&\vdots&\ \ \,\vdots&\cdots&\vdots\\
\end{array} \right)$.\\[2.5ex]
If $n$ is even, then $\sigma_{\!A^{(n)}}^2 = \dfrac{2\,n}{n-1}$ and thus\\[2ex]
\hspace*{13.2ex}$\displaystyle{D_{\!A^{(n)}}^2 = \dfrac{1}{2}\,\Bigl(\,\dfrac{n-1}{n}\,\Bigr)^2, \ \ \
E_{\!A^{(n)}}^3 = \dfrac{1}{2^{\,3/2}}\,\Bigl(\,\dfrac{n-1}{n}\,\Bigr)^{5/2}}$.\vspace*{1.5ex}
\item\label{EWD_3_1_17_BWc}
In chapter \ref{EWD_Kap4} we will apply the two theorems to matrices of the form $a_{ij} = e_{i}\,d_{j}$.
In doing so, we specify conditions that are easy to prove and under which the conditions 
(\ref{EWD_3_1_11}) or (\ref{EWD_3_1_14}) and (\ref{EWD_3_1_15}) are fulfilled
(but with {''constants''} $\mathcal{C}_{1}$, $\mathcal{C}_{2}$ and $\mathcal{C}_{3}$, which grow slightly with $n$).
\end{enumerate}\vspace*{0.5ex}
The following proof of the Theorems \ref{EWD_3_1_10} and \ref{EWD_3_1_13} is divided into seven sections. 
The next five sections (i.e. \ref{EWD_Kap3_Sec2} - \ref{EWD_Kap3_Sec6}) contain numerous preparations. 
These ensure that we can prove Theorem \ref{EWD_3_1_10} in section \ref{EWD_Kap3_Sec7} and 
Theorem \ref{EWD_3_1_13} in section \ref{EWD_Kap3_Sec8} using Stein's method.\\[2.8ex]
In this section \ref{EWD_Kap3_Sec1} we will next show some simple, but repeatedly needed statements about quantities 
that were introduced at the beginning of this section. 
After that we will make some statements about the necessity of the conditions 
(\ref{EWD_3_1_11}), (\ref{EWD_3_1_14}) and (\ref{EWD_3_1_15}).\\[4ex]  
\refstepcounter{DSBcount}
\label{EWD_3_1_18}
\textbf{\hspace*{-0.8ex}\theDSBcount\ Lemma}
\begin{enumerate}
\item\label{EWD_3_1_18_BWa}
If $2 < k$, then\\[2ex]
\hspace*{13.2ex}$\displaystyle{\sum\limits_{i, j = 1}^{n} |\,\bm\hat{a}_{ij}\,|^k\,\geq\,
\dfrac{1}{2^{k/2}n^{(k/2)-2}}}$\hspace*{1ex}.\\[2ex]
In particular,\\[2ex]
$\displaystyle{\beta_{A} \geq \Bigl(\,\dfrac{n}{8}\,\Bigr)^{1/2}}$,\hspace*{2.5ex}
$\displaystyle{\delta_{A} \geq \dfrac{1}{4}}$,\hspace*{2.5ex}
$\displaystyle{\eta_{A} \geq \Bigl(\,\dfrac{1}{32\,n}\,\Bigr)^{1/2}}$\hspace*{3ex}and\hspace*{3ex} 
$D_{\!A} \geq \dfrac{1}{2\,\sqrt{n}}$,\hspace*{2.5ex}
$E_{\!A} \geq \dfrac{1}{2^{5/6}\,\sqrt{n}}$.\vspace*{1ex}
\item\label{EWD_3_1_18_BWb}
If $2 < k < r$, then\\[2ex]
\hspace*{13.2ex}$\displaystyle{\Bigl(\,\dfrac{1}{n}\,
\sum\limits_{i, j = 1}^{n} |\,\bm\hat{a}_{ij}\,|^k\,\Bigr)^{1/(k-2)}\,\leq\,
\Bigl(\,\dfrac{1}{n}\,
\sum\limits_{i, j = 1}^{n} |\,\bm\hat{a}_{ij}\,|^r\,\Bigr)^{1/(r-2)}}$\hspace*{1ex}.\\[1ex]
In particular,
\hspace*{2ex}$\dfrac{\beta_{A}}{n} \leq D_{\!A} \leq E_{\!A}$.\vspace*{1ex}
\item\label{EWD_3_1_18_BWc}
$\displaystyle{\sum\limits_{i, j, k = 1}^{n} \bm\hat{a}_{ij}^2\,\bm\hat{a}_{ik}^2\,\leq\,
n\sum\limits_{i, j = 1}^{n} \bm\hat{a}_{ij}^4}$\hspace*{2ex}and
\hspace*{2ex}$\displaystyle{\sum\limits_{i, j, k = 1}^{n} \bm\hat{a}_{ij}^2\,\bm\hat{a}_{kj}^2\,\leq\,
n\sum\limits_{i, j = 1}^{n} \bm\hat{a}_{ij}^4}$\hspace*{1ex}.\vspace*{1ex}
\item\label{EWD_3_1_18_BWd}
$\displaystyle{|\,\lambda_{1,A}\,|\,\leq\, \dfrac{\beta_{A}}{n}}$
\hspace*{3ex}and\hspace*{3ex}
$\displaystyle{|\,\lambda_{2,A}\,|\,\leq\,\dfrac{3}{n} + 7\,\dfrac{\delta_{A}}{n}
\,\leq\,19\,\dfrac{\delta_{A}}{n}}$.\vspace*{1ex}
\item\label{EWD_3_1_18_BWe}
$||\,e_{2,A} - e_{1,A}\,|| \leq \dfrac{1}{2}\,\dfrac{\delta_{A}}{n} \,\Bigl(\,= \dfrac{1}{2}\,D_{\!A}^2\,\Bigr)$
\hspace*{3ex}and\hspace*{3ex}
$||\,e'_{2,A} - e'_{1,A}\,|| \leq \dfrac{11}{10}\,\dfrac{\delta_{A}}{n}$,\\[3ex]
$||\,x\,e''_{2,A} - x\,e''_{1,A}\,|| \leq \dfrac{9}{5}\,\dfrac{\delta_{A}}{n}$.\vspace*{1ex}
\item\label{EWD_3_1_18_BWf}
$||\,\mathscr{F}_{\!A} - e_{1,A}\,|| 
\leq ||\,\mathscr{F}_{\!A} - \Phi\,|| + \dfrac{1}{6}\,\dfrac{\beta_{A}}{n}\,||\,H_{2}\,\psi\,||
\leq 1 + \dfrac{1}{15}\,\dfrac{\beta_{A}}{n}$.\vspace*{1ex}
\item\label{EWD_3_1_18_BWg}
$||\,e_{1,A}\,||
\leq ||\,\Phi\,|| + \dfrac{1}{6}\,\dfrac{\beta_{A}}{n}\,||\,H_{2}\,\psi\,||
\leq 1 + \dfrac{1}{15}\,\dfrac{\beta_{A}}{n}$
\hspace*{3ex}and\hspace*{3ex}
$||\,e'_{1,A}\,|| \leq \dfrac{2}{5} + \dfrac{1}{10}\,\dfrac{\beta_{A}}{n}$,\\[3ex]
$||\,x\,e''_{1,A}(x)\,|| \leq \dfrac{3}{10} + \dfrac{1}{5}\,\dfrac{\beta_{A}}{n}$
\hspace*{3ex}and\hspace*{3ex}
$||\,e''_{1,A}\,|| \leq \dfrac{1}{4} + \dfrac{1}{5}\,\dfrac{\beta_{A}}{n}$.\vspace*{1ex}
\item\label{EWD_3_1_18_BWh}
\begin{tabular}[t]{@{}l@{\hspace*{0.8ex}}c@{\hspace*{0.8ex}}l@{}}
$||\,\mathscr{F}_{\!A} - e_{2,A}\,||$&$\leq$&
$1 + \dfrac{1}{15}\,\dfrac{\beta_{A}}{n} + \dfrac{1}{2}\,\dfrac{\delta_{A}}{n}$.
\end{tabular}\vspace*{1ex}
\item\label{EWD_3_1_18_BWi}
$||\,e_{2,A}\,|| \leq 1 + \dfrac{1}{15}\,\dfrac{\beta_{A}}{n} +
\dfrac{1}{2}\,\dfrac{\delta_{A}}{n}$
\hspace*{3ex}and\hspace*{3ex}
$||\,e'_{2,A}\,|| \leq \dfrac{2}{5} + \dfrac{1}{10}\,\dfrac{\beta_{A}}{n} + \dfrac{11}{10}\,\dfrac{\delta_{A}}{n}$,\\[3ex]
$||\,x\,e''_{2,A}(x)\,|| \leq \dfrac{3}{10} + \dfrac{1}{5}\,\dfrac{\beta_{A}}{n} + \dfrac{9}{5}\,\dfrac{\delta_{A}}{n}$
\hspace*{3ex}and\hspace*{3ex}
$||\,e'''_{2,A}\,|| \leq \dfrac{2}{5} + \dfrac{7}{18}\,\dfrac{\beta_{A}}{n} +
\dfrac{16}{3}\,\dfrac{\delta_{A}}{n}$.
\item\label{EWD_3_1_18_BWj}
Let\hspace*{2ex}\index{Edgeworth expansion!$e_{1,A}^1$, $e_{2,A}^m$} 
$\displaystyle{e_{1,A}^1(z) = \int\limits_{- \infty}^{z} x\,e_{1,A}'(x)\,dx = 
- \psi(z) - z^3\,\psi(z)\,\dfrac{\lambda_{1,A}}{6}}$.\hspace*{2ex}
Then\hspace*{2ex} 
$||\,e_{1,A}^1\,|| \leq \dfrac{2}{5} + \dfrac{7}{90}\,\dfrac{\beta_{A}}{n}$.
\end{enumerate}
\vspace*{3.5ex}
\textbf{Proof:}
\begin{enumerate}
\item
Since we assume $\sigma_{\!A} > 0$ and thus $n \geq 2$, we have $n \leq 2\,(n-1)$.
By applying $1 = \sigma_{\!\hat{A}}^2$ (cf. (\ref{EWD_3_1_01})) and 
H{\"o}lder's inequality\index{H{\"o}lder's inequality!for finite sequences using length $\nu$} 
(\ref{EWD_0_1_05}) with $\nu = n^2$ and \mbox{\rule[0ex]{0ex}{4.4ex}$p = \dfrac{k}{2} > 1$} we obtain\\[2ex]
\hspace*{13.2ex}
\begin{tabular}{@{}l@{\hspace*{0.8ex}}c@{\hspace*{0.8ex}}l@{}}
$\displaystyle{1 = \dfrac{1}{n-1}\,\sum\limits_{i,j=1}^{n} |\hat{a}_{ij}|^2}$
&$\leq$&$\displaystyle{\dfrac{1}{n-1}\,\Bigl(\,\sum\limits_{i,j=1}^{n} \bigl(\,|\hat{a}_{ij}|^2\,\bigr)^{k/2}
\,\Bigr)^{2/k} \cdot \bigl(\,n^2\,\bigr)^{(k-2)/k}}$\\[2.5ex]
&$=$&$\displaystyle{\dfrac{1}{n-1}\,\Bigl(\,\sum\limits_{i,j=1}^{n} |\hat{a}_{ij}|^k
\,\Bigr)^{2/k} \cdot n^{2(k-2)/k}}$.
\end{tabular}\\[2ex]
Exponentiating this with $p = \dfrac{k}{2}$ and then applying $n-1 \geq \dfrac{n}{2}$ leads to\\[2ex]
\refstepcounter{DSBcount}
\label{EWD_3_1_19}
\text{\hspace*{-0.8ex}(\theDSBcount)}
\hspace*{4ex}
$\displaystyle{\sum\limits_{i,j=1}^{n} |\hat{a}_{ij}|^k \geq
\dfrac{(n-1)^{k/2}}{n^{k-2}} \geq
\dfrac{n^{k/2}}{2^{k/2}n^{k-2}} = \dfrac{1}{2^{k/2}n^{(k/2)-2}}}$.\vspace*{1ex}
\item
Because of
\mbox{\rule[-4ex]{0ex}{7ex}$k = r\,\dfrac{k - 2}{r - 2} + 2\,\dfrac{r - k}{r - 2}$},
H{\"o}lder's inequality\index{H{\"o}lder's inequality!for finite sequences} with $p = \dfrac{r - 2}{k - 2}$,
$q = \dfrac{r - 2}{r - k}$ and $1 = \sigma_{\!\hat{A}}^2$ (cf. (\ref{EWD_3_1_01})) gives\\[2ex]
\hspace*{5ex}
\begin{tabular}[t]{@{}l@{\hspace*{0.8ex}}c@{\hspace*{0.8ex}}l@{\hspace*{0.8ex}}c@{\hspace*{0.8ex}}l@{}}
$\displaystyle{\sum\limits_{i, j = 1}^{n} |\,\bm\hat{a}_{ij}\,|^k}$&$\leq$&
$\displaystyle{\biggl(\,\sum\limits_{i, j = 1}^{n} |\,\bm\hat{a}_{ij}\,|^r\,\biggr)^{\frac{k-2}{r-2}} \cdot 
\biggl(\,\sum\limits_{i, j = 1}^{n} \bm\hat{a}_{ij}^2\,\biggr)^{\frac{r-k}{r-2}}}$
&$\leq$&$\displaystyle{\biggl(\,\sum\limits_{i, j = 1}^{n} |\,\bm\hat{a}_{ij}\,|^r\,\biggr)^{\frac{k-2}{r-2}} \cdot
n^{\frac{r-k}{r-2}}}$\hspace*{1ex}.
\end{tabular}\\[2.5ex]
The assertion follows from this due to $\dfrac{1}{k-2} -\dfrac{1}{r-2} = \dfrac{r-k}{(k-2)\,(r-2)}$.
\item
A further application of H{\"o}lder's inequality\index{H{\"o}lder's inequality!for finite sequences using length $\nu$} 
(\ref{EWD_0_1_05}) with $\nu = n$ and $p = 2$ gives\\[2ex]
\hspace*{13.2ex}
$\displaystyle{\sum\limits_{i, j, k = 1}^{n} \bm\hat{a}_{ij}^2\,\bm\hat{a}_{ik}^2\,=\,
\sum\limits_{i = 1}^{n}\,\biggl(\,\sum\limits_{j = 1}^{n} \bm\hat{a}_{ij}^2\,\biggr)^2\,\leq\,
n\sum\limits_{i, j = 1}^{n} \bm\hat{a}_{ij}^4}$.\\[2.5ex]
The second inequality follows analogously.\vspace*{1ex}
\item
Using part \ref{EWD_3_1_18_BWc}) and part \ref{EWD_3_1_18_BWa}) we obtain\\[2ex]
\hspace*{5ex}
$\displaystyle{|\,\lambda_{2,A}\,| 
\,\leq\, \dfrac{\delta_{A}}{n} + \dfrac{3}{n} + \dfrac{3}{n^2}\,2\,n\,\delta_{A}
\,=\, \dfrac{3}{n} + 7\,\dfrac{\delta_{A}}{n} 
\,\leq\, 12\,\dfrac{\delta_{A}}{n} + 7\,\dfrac{\delta_{A}}{n} 
\,=\, 19\,\dfrac{\delta_{A}}{n}}$.\vspace*{1ex}
\item
Using part \ref{EWD_3_1_18_BWd}), Lemma \ref{EWD_2_2_08}, \ref{EWD_2_2_08_BWf}) and \ref{EWD_2_2_08_BWj}), and then 
part \ref{EWD_3_1_18_BWb}) we obtain\\[2ex]
\hspace*{5ex}
\begin{tabular}[t]{@{}l@{\hspace*{0.8ex}}c@{\hspace*{0.8ex}}l@{}}
$||\,e_{2,A} - e_{1,A}\,||$
&$=$&$||\,\dfrac{\lambda_{2,A}}{24}\,\psi(x)\,H_{3}(x) + 
\dfrac{\lambda_{1,A}^{2}}{72}\,\psi(x)\,H_{5}(x)\,||$\\[3.5ex]
&$\leq$&$\displaystyle{\dfrac{1}{24}\,|\lambda_{2,A}| \cdot
||\,\psi(x)\,H_{3}(x)\,|| + 
\dfrac{1}{72}\,|\lambda_{1,A}|^2 \cdot ||\,\psi(x)\,H_{5}(x)\,||}$\\[3.5ex]
&$\leq$&$\displaystyle{\dfrac{19}{24} \, \dfrac{\delta_{\!A}}{n} \, \dfrac{5}{9} +
\dfrac{1}{72} \, \Bigl(\,\dfrac{\beta_{\!A}}{n}\,\Bigr)^2} \, \dfrac{7}{3}
\leq \dfrac{1}{2}\,\dfrac{\delta_{\!A}}{n}$.
\end{tabular}\\[2.5ex]
Analogously, using Lemma \ref{EWD_2_2_08}, \ref{EWD_2_2_08_BWh}) and \ref{EWD_2_2_08_BWl}),
and Lemma \ref{EWD_2_2_08}, \ref{EWD_2_2_08_BWk}) and \ref{EWD_2_2_08_BWn}) gives\\[2.5ex]
\hspace*{5ex}
$||\,e'_{2,A} - e'_{1,A}\,|| \leq \displaystyle{\dfrac{19}{24} \, \dfrac{\delta_{\!A}}{n} \, \dfrac{6}{5} +
\dfrac{1}{72} \, \Bigl(\,\dfrac{\beta_{\!A}}{n}\,\Bigr)^2} \, 6
\leq \dfrac{11}{10}\,\dfrac{\delta_{\!A}}{n}$,\\[2ex]
\hspace*{5ex}
$||\,x\,e''_{2,A} - x\,e''_{1,A}\,|| \leq \displaystyle{\dfrac{19}{24} \, \dfrac{\delta_{\!A}}{n} \, 2 +
\dfrac{1}{72} \, \Bigl(\,\dfrac{\beta_{\!A}}{n}\,\Bigr)^2} \, 13
\leq \dfrac{9}{5}\,\dfrac{\delta_{\!A}}{n}$.\vspace*{1ex}
\item
The second inequality follows from Lemma \ref{EWD_2_2_08}, \ref{EWD_2_2_08_BWe}).\vspace*{1ex}
\item
The estimates follow because of Lemma \ref{EWD_2_2_08}, \ref{EWD_2_2_08_BWe}), \ref{EWD_2_2_08_BWf}), 
\ref{EWD_2_2_08_BWi}), \ref{EWD_2_2_08_BWh}) for the second summands (in this order) and
due to Lemma \ref{EWD_2_2_08}, \ref{EWD_2_2_08_BWa}), \ref{EWD_2_2_08_BWc}), \ref{EWD_2_2_08_BWb}) 
for the first summands.\vspace*{1ex}
\item
This follows from the parts \ref{EWD_3_1_18_BWf}) and \ref{EWD_3_1_18_BWe}).\vspace*{1ex}
\item
With the exception of the estimate for $||\,e'''_{2,A}\,||$, this follows from the
parts \ref{EWD_3_1_18_BWg}) and \ref{EWD_3_1_18_BWe}). 
To estimate $||\,e'''_{2,A}\,||$, we use Lemma \ref{EWD_2_2_08}, \ref{EWD_2_2_08_BWe}), \ref{EWD_2_2_08_BWj}),
\ref{EWD_2_2_08_BWl}), \ref{EWD_2_2_08_BWo}) and get\\[2.5ex]
\hspace*{5ex}
\begin{tabular}[b]{@{}l@{\hspace*{0.8ex}}c@{\hspace*{0.8ex}}l@{}}
$||\,e'''_{2,A}\,||$&$\leq$&$||\,H_{2}\,\psi\,|| + \dfrac{1}{6}\,\dfrac{\beta_{A}}{n}\,||\,H_{5}\,\psi\,||
+ \dfrac{19}{24}\,\dfrac{\delta_{A}}{n}\,||\,H_{6}\,\psi\,|| + 
\dfrac{1}{72}\,\Bigl(\,\dfrac{\beta_{A}}{n}\,\Bigr)^2\,||\,H_{8}\,\psi\,||$
\end{tabular}\\[3ex]
\hspace*{12.4ex}
\begin{tabular}[b]{@{}l@{\hspace*{0.8ex}}c@{\hspace*{0.8ex}}l@{}}
&$\leq$&$\dfrac{2}{5} + \dfrac{1}{6}\,\dfrac{\beta_{A}}{n}\,\dfrac{7}{3}
+ \dfrac{19}{24}\,\dfrac{\delta_{A}}{n}\,6 + \dfrac{1}{72}\,\dfrac{\delta_{A}}{n}\,42$
\end{tabular}\\[3.5ex]
\hspace*{12.4ex}
\begin{tabular}[b]{@{}l@{\hspace*{0.8ex}}c@{\hspace*{0.8ex}}l@{}}
&$=$&$\dfrac{2}{5} + \dfrac{7}{18}\,\dfrac{\beta_{A}}{n} + \dfrac{16}{3}\,\dfrac{\delta_{A}}{n}$.
\end{tabular}\vspace*{1ex}
\item
The estimate follows from Lemma \ref{EWD_2_2_08}, \ref{EWD_2_2_08_BWa}) and \ref{EWD_2_2_08_BWd}).
\hspace*{1ex}\hfill$\Box$
\end{enumerate}
\vspace*{2.5ex}
We now show that from (\ref{EWD_3_1_12}) we obtain the condition (\ref{EWD_3_1_11}) 
for $B = \bm\hat{A}$ with a different constant. In the same way, we can derive 
the conditions (\ref{EWD_3_1_14}) and (\ref{EWD_3_1_15}) 
for $B = \bm\hat{A}$ with other constants from (\ref{EWD_3_1_16}).\\[0.5ex] 
However, these proofs do not apply to the submatrices
$B \in N(8,\bm\hat{A})$ (or $B \in N(16,\bm\hat{A})$).\\[4ex]
\refstepcounter{DSBcount}
\label{EWD_3_1_20}
\textbf{\hspace*{-0.8ex}\theDSBcount\ Proposition}
\begin{enumerate}
\item\label{EWD_3_1_20_BWa}
If $||\,\mathscr{F}_{\!A} - e_{1,A}\,|| \leq K\,D_{\!A}^{2}$, then\\[2.3ex]
\refstepcounter{DSBcount}
\label{EWD_3_1_21}
\text{\hspace*{-0.8ex}(\theDSBcount)}
\hspace*{4ex}
$\big|\,\Delta_{y}^{2}\mathscr{F}_{\!A}(z)\,\big| \leq \bigl(\,4\,K + 1\,\bigr)\,\bigl(\,D_{\!A}^2 + y^2\,\bigr)$
\hspace*{2ex}for all $z, y \in \mathbb{R}$.
\item\label{EWD_3_1_20_BWb}
If $||\,\mathscr{F}_{\!A} - e_{2,A}\,|| \leq K\,E_{\!A}^{3}$, then\\[2.3ex]
\refstepcounter{DSBcount}
\label{EWD_3_1_22}
\text{\hspace*{-0.8ex}(\theDSBcount)}
\hspace*{4ex}
$\Big|\,\mathscr{F}_{\!A}(x) - P_{E_{\!A}}^{2}(x;z,\mathscr{F}_{\!A})\,\Big| \leq 
\bigl(\,8\,K + 6\,\bigr)\,\Bigl(\,E_{\!A}^{3} + |x - z|^{3}\,\Bigr)$
\hspace*{2ex}for all $z, x \in \mathbb{R}$.
\item\label{EWD_3_1_20_BWc}
If $||\,\mathscr{F}_{\!A} - e_{2,A}\,|| \leq K\,E_{\!A}^{3}$, then\\[2.3ex]
\refstepcounter{DSBcount}
\label{EWD_3_1_23}
\text{\hspace*{-0.8ex}(\theDSBcount)}
\hspace*{4ex}
$\big|\,z\,\Delta_{y}^{2}\mathscr{F}_{\!A}(z)\,\big| \leq \bigl(\,16\,K + 11\,\bigr)\,\bigl(\,E_{\!A}^2 + y^2\,\bigr)$
\hspace*{2ex}for all $z \in \mathbb{R}$, $0 \leq y \leq E_{\!A}$.
\item\label{EWD_3_1_20_BWd}
If $||\,\mathscr{F}_{\!A} - e_{2,A}\,|| \leq K\,E_{\!A}^{3}$, then\\[2.3ex]
\refstepcounter{DSBcount}
\label{EWD_3_1_24}
\text{\hspace*{-0.8ex}(\theDSBcount)}
\hspace*{4ex}
$||\,\mathscr{F}_{\!A} - e_{1,A}\,|| \leq \Bigl(\,K + \dfrac{16}{15}\,\Bigr)\,E_{\!A}^{2}$\\[2.5ex]
and thus\\[2ex]
\refstepcounter{DSBcount}
\label{EWD_3_1_25}
\text{\hspace*{-0.8ex}(\theDSBcount)}
\hspace*{4ex}
$\big|\,\Delta_{y}^{2}\mathscr{F}_{\!A}(z)\,\big| \leq 
\biggl(\,4\,\Bigl(\,K + \dfrac{16}{15}\,\Bigr) + 1\,\biggr)\,\bigl(\,E_{\!A}^2 + y^2\,\bigr)$
\hspace*{2ex}for all $z, y \in \mathbb{R}$.
\end{enumerate}
\vspace*{2.5ex}
\textbf{Proof:}
\begin{enumerate}
\item
Due to Lemma \ref{EWD_2_4_02}, \ref{EWD_2_4_02_BWb}) (for $e_{1,A}$ and $k = 2$) and 
Lemma \ref{EWD_3_1_18}, \ref{EWD_3_1_18_BWg}) (for $||\,e''_{1,A}\,||$) we get\\[2ex]
\hspace*{13.2ex}\begin{tabular}[t]{@{}l@{\hspace*{0.8ex}}c@{\hspace*{0.8ex}}l@{}}
$\big|\,\Delta_{y}^{2}\mathscr{F}_{\!A}(z)\,\big|$&$\leq$&
$4\,K\,D_{\!A}^{2} + \big|\,\Delta_{y}^{2}e_{1,A}(z)\,\big|$\\[2.5ex]
&$\leq$&$4\,K\,D_{\!A}^{2} + \Bigl(\,\dfrac{1}{4} + \dfrac{1}{5}\,\dfrac{\beta_{A}}{n}\,\Bigr)\,y^2$.
\end{tabular}\\[1.7ex]
Now, if \mbox{\rule[-3ex]{0ex}{6ex}$\dfrac{\beta_{A}}{n} \leq 1$}, then (\ref{EWD_3_1_21}) holds.
On the other hand, if \mbox{\rule[-3ex]{0ex}{6ex}$1 < \dfrac{\beta_{A}}{n}$},
then (\ref{EWD_3_1_21}) follows using Lemma \ref{EWD_3_1_18}, \ref{EWD_3_1_18_BWb})\\[2ex] 
\hspace*{13.2ex}$\big|\,\Delta_{y}^{2}\mathscr{F}_{\!A}(z)\,\big| 
\leq 1 \leq \Bigl(\,\dfrac{\beta_{A}}{n}\,\Bigr)^2 \leq D_{\!A}^2$.\vspace*{0.5ex}
\item
Due to Lemma \ref{EWD_2_4_02}, \ref{EWD_2_4_02_BWa}) (for $e_{2,A}$ and $k = 2$) and 
Lemma \ref{EWD_3_1_18}, \ref{EWD_3_1_18_BWi}) (for $||\,e'''_{2,A}\,||$) we get\\[2.3ex]
\begin{tabular}[t]{@{}l@{\hspace*{0.8ex}}c@{\hspace*{0.8ex}}l@{}}
$\Big|\,\mathscr{F}_{\!A}(x) - P_{E_{\!A}}^{2}(x;z,\mathscr{F}_{\!A})\,\Big|$&
$\leq$&$\Big|\,\mathscr{F}_{\!A}(x) - e_{2,A}(x)\,\Big| + 
\Big|\,e_{2,A}(x) - P_{E_{\!A}}^{2}(x;z,e_{2,A})\,\Big|$\\[2.8ex] 
&&$+\, \Big|\,P_{E_{\!A}}^{2}(x;z,e_{2,A}) - P_{E_{\!A}}^{2}(x;z,\mathscr{F}_{\!A})\,\Big|$
\hspace*{7.8ex}(cf. \cite{10.1214/aop/1176993873}, page 502)\\[3.3ex]
&$\leq$&$K\,E_{\!A}^{3} + \dfrac{4}{3}\,\Bigl(\,\dfrac{2}{5} + \dfrac{7}{18}\,\dfrac{\beta_{A}}{n}
+ \dfrac{16}{3}\,\dfrac{\delta_{A}}{n}\,\Bigr)\,\Bigl(\,E_{\!A}^{3} + |x - z|^{3}\,\Bigr)$\\[3.3ex]
&&$+\,K\,E_{\!A}^{3} + 2\,K\,E_{\!A}^{3}\,\dfrac{|x - z|}{E_{\!A}} + 4\,K\,E_{\!A}^{3}\,
\dfrac{|x - z|}{E_{\!A}}\,\dfrac{|x - z - E_{\!A}|}{2\,E_{\!A}}$\\[3.3ex]
&$\leq$&$K\,E_{\!A}^{3} + \Bigl(\,\dfrac{8}{15} + \dfrac{14}{27}\,\dfrac{\beta_{A}}{n}
+ \dfrac{64}{9}\,\dfrac{\delta_{A}}{n}\,\Bigr)\,\Bigl(\,E_{\!A}^{3} + |x - z|^{3}\,\Bigr)$\\[3.3ex]
&&$+\,K\,E_{\!A}^{3} + 2\,K\,E_{\!A}^{2}\,|x - z| + 2\,K\,E_{\!A}\,
\Bigl(\,|x - z|^2 + |x - z|\,E_{\!A}\,\Bigr)$\\[3.3ex]
&$\leq$&$\Bigl(\,8\,K + \dfrac{8}{15} + \dfrac{14}{27}\,\dfrac{\beta_{A}}{n} + \dfrac{64}{9}\,\dfrac{\delta_{A}}{n}    
\,\Bigr)\,\Bigl(\,E_{\!A}^{3} + |x - z|^{3}\,\Bigr)$.
\end{tabular}\\[3.5ex]
For the last inequality, the following formula was used:\\[2ex]
\hspace*{13.2ex}$a^2\,b \leq \bigl(\,\text{max}\{a, b\}\,\bigr)^3 \leq a^3 + b^3$
\hspace*{4ex}for $a, b \geq 0$.\\[2.5ex]
In the case of \mbox{\rule[-3ex]{0ex}{6ex}$E_{\!A}^{3} \leq \dfrac{1}{2}$}, from which 
\mbox{\rule[-3ex]{0ex}{6ex}$\dfrac{\beta_{A}}{n} \leq \Bigl(\,\dfrac{1}{2}\,\Bigr)^{1/3}$} 
and \mbox{\rule[-3ex]{0ex}{6ex}$\dfrac{\delta_{A}}{n} \leq \Bigl(\,\dfrac{1}{2}\,\Bigr)^{2/3}$} follow 
(cf. Lemma \ref{EWD_3_1_18}, \ref{EWD_3_1_18_BWb})), we thus obtain
(\ref{EWD_3_1_22}).\\[3ex]
On the other hand, in the case of \mbox{\rule[-2.5ex]{0ex}{5ex}$E_{\!A}^{3} > \dfrac{1}{2}$}, 
we get due to the inequalities $|x| \leq 1 + |x|^3$ and $x^2 \leq 1 + |x|^3$:\\[2.5ex]
\begin{tabular}[t]{@{}l@{\hspace*{0.8ex}}c@{\hspace*{0.8ex}}l@{}}
$\Big|\,\mathscr{F}_{\!A}(x) - P_{E_{\!A}}^{2}(x;z,\mathscr{F}_{\!A})\,\Big|$&
$\leq$&$\Big|\,\mathscr{F}_{\!A}(x) - \mathscr{F}_{\!A}(z)\,\Big| + \dfrac{|x - z|}{E_{\!A}}
+ \dfrac{1}{2}\,\dfrac{|x - z|^2 + |x - z|\,E_{\!A}}{E_{\!A}^2}$\\[3ex]
&$\leq$&$3\,\biggl(\,1 + \dfrac{|x - z|^{3}}{E_{\!A}^{3}}\,\biggr)$\\[3ex]
&$\leq$&$6\,\Bigl(\,E_{\!A}^{3} + |x - z|^{3}\,\Bigr)$.
\end{tabular}\vspace*{0.5ex}
\item
Because of Lemma \ref{EWD_2_4_02}, \ref{EWD_2_4_02_BWc}) (for $e_{2,A}$) holds\\[2.3ex]
\begin{tabular}[t]{@{}l@{\hspace*{0.8ex}}c@{\hspace*{0.8ex}}l@{}}
$\big|\,z\,\Delta_{y}^{2}\mathscr{F}_{\!A}(z)\,\big|$
&$\leq$&$|z|\,\big|\,\mathscr{F}_{\!A}(z + 2 y) - e_{2,A}(z + 2 y)\,\big| + 
2\,|z|\,\big|\,\mathscr{F}_{\!A}(z + y) - e_{2,A}(z + y)\,\big|$
\end{tabular}\\[2.5ex]
\begin{tabular}[t]{@{}l@{\hspace*{0.8ex}}c@{\hspace*{0.8ex}}l@{}}
\hspace*{13ex}&&$+\,\,|z|\,\big|\,\mathscr{F}_{\!A}(z) - e_{2,A}(z)\,\big|
+ \big|\,z\,\Delta_{y}^{2}e_{2,A}(z)\,\big|$\\[2.5ex]
&$\leq$&$4\,|z|\,K\,E_{\!A}^{3} + 
\Bigl(\,4\,\big|\big|\,e_{2,A}'\,\big|\big| + \big|\big|\,x\,e_{2,A}''(x)\,\big|\big|\,\Bigr)\,y^2$.
\end{tabular}\\[2.5ex]
We now assume 
\mbox{\rule[-1.2ex]{0ex}{3.2ex}$|z| \leq \dfrac{4}{E_{\!A}}$} and $E_{\!A} \leq 1$. 
Using Lemma \ref{EWD_3_1_18}, \ref{EWD_3_1_18_BWb}) and \ref{EWD_3_1_18_BWi}), we 
then obtain\\[2.3ex]
\hspace*{13.2ex}\begin{tabular}[t]{@{}l@{\hspace*{0.8ex}}c@{\hspace*{0.8ex}}l@{}}
$\big|\,z\,\Delta_{y}^{2}\mathscr{F}_{\!A}(z)\,\big|$
&$\leq$&$4\,\dfrac{4}{E_{\!A}}\,K\,E_{\!A}^{3} +  
\Bigl(\,4\,\dfrac{8}{5} + \dfrac{23}{10}\,\Bigr)\,y^2$\\[2.5ex]
&$\leq$&$16\,K\,E_{\!A}^{2} + 9\,y^2$.
\end{tabular}\\[2.8ex]
For the further calculation, let $z \in \mathbb{R}$ and $0 \leq y \leq E_{\!A}$ be general again:\\[3ex]
\refstepcounter{DSBcount}
\label{EWD_3_1_26}
\text{\hspace*{-0.8ex}(\theDSBcount)}
\hspace*{4ex}
\begin{tabular}{@{}l@{\hspace*{0.8ex}}c@{\hspace*{0.8ex}}l@{}}
$\big|\,z\,\Delta_{y}^{2}\mathscr{F}_{\!A}(z)\,\big|$&
$=$&$\Big|\,z\,P\Bigl(\,z+y < \mathscr{T}_{A} \leq z + 2 y\,\Bigr) - 
z\,P\Bigl(\,z < \mathscr{T}_{A} \leq z + y\,\Bigr)\,\Big|$\\[2.5ex]
&$\leq$&$|z|\,P\Bigl(\,z+y < \mathscr{T}_{A} \leq z + 2 y\,\Bigr) + |z|\,P\Bigl(\,z < 
\mathscr{T}_{A} \leq z + y\,\Bigr)$\\[2.5ex] 
&$\leq$&$|z|\,P\Bigl(\,z - 2\, E_{\!A} \leq \mathscr{T}_{A} \leq z + 2\, E_{\!A}\,\Bigr)$\\[2.5ex]
&$\leq$&$|z|\,P\Bigl(\,|z| - 2\, E_{\!A} \leq |\mathscr{T}_{A}|\,\Bigr)$\\[2.5ex]
&$=$&$|z|\,P\Bigl(\,|z| \leq |\mathscr{T}_{A}| + 2\, E_{\!A} \,\Bigr)$.
\end{tabular}\\[2ex]
We now suppose \mbox{\rule[-3ex]{0ex}{6.5ex}$|z| \geq \dfrac{4}{E_{\!A}}$} and $E_{\!A} \leq 1$. 
An application of Markov's inequality\index{Markov's inequality} with the function $h(x) = |x|^3$ then yields\\[3ex]
\hspace*{13.2ex}\begin{tabular}{@{}l@{\hspace*{0.8ex}}c@{\hspace*{0.8ex}}l@{\hspace*{6.9ex}}r@{}}
$\big|\,z\,\Delta_{y}^{2}\mathscr{F}_{\!A}(z)\,\big|$&
$\leq$&$|z|\,P\Bigl(\,|z| \leq |\mathscr{T}_{A}| + 2\,\Bigr)$&(since $E_{\!A} \leq 1$)\\[2.5ex]
&$\leq$&$|z|\,\dfrac{1}{|z|^3}\,E\bigl(\,\big|\,|\mathscr{T}_{A}| + 2\,\big|^3\,\bigr)$
&(Markov's inequality\index{Markov's inequality})\\[2.5ex]
&$\leq$&$\dfrac{1}{|z|^2}\,E\bigl(\,\bigl(\,|\mathscr{T}_{A}| + 2\,\bigr)^4\,\bigr)^{3/4}$
&(H{\"o}lder's inequality\index{H{\"o}lder's inequality!for random variables})\\[2.5ex]
&$\leq$&$\dfrac{1}{|z|^2}\,E\bigl(\,8\,\bigl(\,\mathscr{T}_{A}^4 + 16\,\bigr)\,\bigr)^{3/4}$
&(\index{H{\"o}lder's inequality!for finite sequences using length $\nu$}(\ref{EWD_0_1_05}) with $\nu = 2$, $p = 4$)\\[2.5ex]
&$=$&$\dfrac{1}{|z|^2}\,\Bigl(\,8\,E\bigl(\,\mathscr{T}_{A}^4\,\bigr) + 128\,\Bigr)^{3/4}$
\end{tabular}\\[3ex]
\hspace*{13.2ex}\begin{tabular}{@{}l@{\hspace*{0.8ex}}c@{\hspace*{0.8ex}}l@{\hspace*{23.3ex}}r@{}}
\hspace*{13ex}&$\leq$&$\dfrac{E_{\!A}^2}{4^2}\,928^{3/4}$&
$\bigl(\,\,|z| \geq \dfrac{4}{E_{\!A}}\ \Rightarrow\ \dfrac{E_{\!A}}{4} \geq
\dfrac{1}{|z|}\,\,\bigr)$\\[2.2ex]
&$\leq$&$11\,E_{\!A}^2$.
\end{tabular}\\[2.5ex]
For the penultimate inequality, we used the estimate
$E\bigl(\,\mathscr{T}_{A}^4\,\bigr) \leq 100$. 
This estimate follows from Proposition \ref{EWD_3_2_02}, \ref{EWD_3_2_02_BWd}) together with 
\mbox{\rule[-3ex]{0ex}{7ex}$\dfrac{\delta_{\!A}}{n} \leq E_{\!A}^2 \leq 1$} 
(cf. Lemma \ref{EWD_3_1_18}, \ref{EWD_3_1_18_BWb})).\\[1.5ex]
It remains to consider the case $E_{\!A} \geq 1$.
We use the inequality (\ref{EWD_3_1_26}) again:\\[2.5ex]
\hspace*{13.2ex}
\begin{tabular}{@{}l@{\hspace*{0.8ex}}c@{\hspace*{0.8ex}}l@{}}
$\big|\,z\,\Delta_{y}^{2}\mathscr{F}_{\!A}(z)\,\big|$&
$\leq$&$|z|\,P\Bigl(\,|z| \leq |\mathscr{T}_{A}| + 2\,E_{\!A}\,\Bigr)$\\[2.5ex] 
&$\leq$&$E\bigl(\,|\mathscr{T}_{A}| + 2\,E_{\!A}\,\bigr)$\\[2.5ex]
&$\leq$&$E\bigl(\,\mathscr{T}_{A}^2\,\bigr)^{1/2} + 2\,E_{\!A}$\\[2.5ex]
&$=$&$1 + 2\,E_{\!A} \leq 3\,E_{\!A}^2$.
\end{tabular}
\item
Because of Lemma \ref{EWD_3_1_18}, \ref{EWD_3_1_18_BWe}) and \ref{EWD_3_1_18_BWb}) we get:\\[2.5ex]
\hspace*{13.2ex}\begin{tabular}{@{}l@{\hspace*{0.8ex}}c@{\hspace*{0.8ex}}l@{}}
$||\,\mathscr{F}_{\!A} - e_{1,A}\,||$&$\leq$&
$||\,\mathscr{F}_{\!A} - e_{2,A}\,|| + ||\,e_{2,A} - e_{1,A}\,||$\\[2.5ex]
&$\leq$&$K\,E_{\!A}^{3} + \dfrac{1}{2}\, D_{\!A}^2$\\[2.5ex] 
&$\leq$&$\Bigl(\,K\,E_{\!A} +  \dfrac{1}{2}\,\Bigr)\, E_{\!A}^2$.
\end{tabular}\\[2.5ex]
If $E_{\!A} \leq 1$, then (\ref{EWD_3_1_24}) holds.
On the other hand, if $1 < E_{\!A}$, then because of Lemma \ref{EWD_3_1_18}, \ref{EWD_3_1_18_BWf}) 
and \ref{EWD_3_1_18_BWb}) we obtain:\\[2.5ex]
\hspace*{13.2ex}\begin{tabular}{@{}l@{\hspace*{0.8ex}}c@{\hspace*{0.8ex}}l@{}}
$||\,\mathscr{F}_{\!A} - e_{1,A}\,||$&$\leq$&$1 + \dfrac{1}{15}\,\dfrac{\beta_{\!A}}{n}$\\[2.5ex] 
&$\leq$&$1 + \dfrac{1}{15}\,E_{\!A} \leq \dfrac{16}{15}\,E_{\!A}^2$.
\end{tabular}\\[2.5ex]
The statement (\ref{EWD_3_1_25}) now follows from (\ref{EWD_3_1_24}) analogously 
to the proof in part \ref{EWD_3_1_20_BWa}).
\hfill$\Box$
\end{enumerate}
\vspace*{3.5ex}
\refstepcounter{DSBcount}
\label{EWD_3_1_27}
\textbf{\theDSBcount\ Remark}\\[0.8ex]
The question arises whether it is sufficient in the conditions of the two Theorems \ref{EWD_3_1_10}
and \ref{EWD_3_1_13} to require the validity of the inequalities only for $\bm\hat{A}$ instead of 
$B \in N(8,\bm\hat{A})$ or $B \in N(16,\bm\hat{A})$.
I~assume that this is the case. However, I cannot provide any proof of this.
                              
\section{Higher moments of linear rank statistics}\label{EWD_Kap3_Sec2}

At first, we consider the third and fourth moment of $\mathscr{T}_{A}$.\\[4ex]
\refstepcounter{DSBcount}
\label{EWD_3_2_01}
\textbf{\hspace*{-0.8ex}\theDSBcount\ Theorem}\index{Theorem!for linear rank statistics!3rd and 4th moment}
\begin{enumerate}
\item\label{EWD_3_2_01_BWa}
$\displaystyle{E\bigl(\,\mathscr{T}_{A}^3\,\bigr) = \dfrac{n}{(n-1)\,(n-2)}
\sum\limits_{i,j = 1}^{n}\,\bm\hat{a}_{ij}^3}$\hspace*{4ex}for $n \geq 3$,\vspace*{1ex}
\item\label{EWD_3_2_01_BWb}
\begin{tabular}[t]{@{}l@{\hspace*{0.8ex}}c@{\hspace*{0.8ex}}l@{}}
$E\bigl(\,\mathscr{T}_{A}^4\,\bigr)$&$=$&
$\displaystyle{3\,\dfrac{(n^2 - 3\,n +1)\,(n-1)}{n\,(n-2)\,(n-3)}}$\\[2.5ex]
&&$\displaystyle{-\, \dfrac{3}{(n-2)\,(n-3)}\,\sum\limits_{i, j, k = 1}^{n}\bigl(\,\bm\hat{a}_{ij}^2\,\bm\hat{a}_{ik}^2
+ \bm\hat{a}_{ij}^2\,\bm\hat{a}_{kj}^2\,\bigr)}$\\[3ex]
&&$\displaystyle{+\, \dfrac{n\,(n+1)}{(n-1)\,(n-2)\,(n-3)}\,
\sum\limits_{i, j = 1}^{n}\bm\hat{a}_{ij}^4}$\\[3ex]
&&$\displaystyle{+\, \dfrac{6}{n\,(n-1)\,(n-2)\,(n-3)}\,
\sum\limits_{i, j, k, h = 1}^{n}
\bm\hat{a}_{ik}\,\bm\hat{a}_{ih}\,\bm\hat{a}_{jk}\,\bm\hat{a}_{jh}}$\hspace*{4ex}for $n \geq 4$.
\end{tabular}
\end{enumerate}
\vspace*{3.5ex}
\textbf{Proof:}\\[0.8ex]
H{\'a}jek, \v{S}id{\'a}k and Sen \cite{sidak1999theory}, Problem 27 to Section 3.3, page 90, 
taking into account\\[0.5ex] 
$\displaystyle{\sum\limits_{i,j = 1}^{n}\bm\hat{a}_{ij}^2 = n - 1}$ (cf. (\ref{EWD_3_1_01})).
\hspace*{1ex}\hfill$\Box$\\[3ex]
A somewhat more manageable formula of the third and fourth moment of $\mathscr{T}_{A}$, 
which is sufficient for our purposes, is contained in the following proposition.\\[4ex]
\refstepcounter{DSBcount}
\label{EWD_3_2_02}
\textbf{\hspace*{-0.8ex}\theDSBcount\ Proposition}
\begin{enumerate}
\item\label{EWD_3_2_02_BWa}
$\displaystyle{E\bigl(\,\mathscr{T}_{A}^3\,\bigr) = \dfrac{1}{n}
\sum\limits_{i,j = 1}^{n}\,\bm\hat{a}_{ij}^3 + R_{3,A}}$,\\[1.5ex]
\hspace*{15ex}where $|R_{3,A}| \leq \dfrac{11}{n^2}\,\beta_{A}$\hspace{4ex}for $n \geq 3$.
\item\label{EWD_3_2_02_BWb}
$\displaystyle{E\bigl(\,\mathscr{T}_{A}^4\,\bigr) = 3 + \dfrac{3}{n}
- \dfrac{3}{n^2}\sum\limits_{i, j, k = 1}^{n}\bigl(\,\bm\hat{a}_{ij}^2\,\bm\hat{a}_{ik}^2
+ \bm\hat{a}_{ij}^2\,\bm\hat{a}_{kj}^2\,\bigr) + \dfrac{1}{n}
\sum\limits_{i, j = 1}^{n}\bm\hat{a}_{ij}^4 + R_{4,A}}$,\\[1.5ex]
\hspace*{15ex}where $|R_{4,A}| \leq \dfrac{336}{n^2}\,\delta_{A}$\hspace{4ex}for $n \geq 4$.
\item\label{EWD_3_2_02_BWc}
$\displaystyle{\Big|\,E\bigl(\,\mathscr{T}_{A}^3\,\bigr)\,\Big| \leq 5\,\dfrac{\beta_{A}}{n}}$
\hspace{4ex}for $n \geq 3$.
\item\label{EWD_3_2_02_BWd}
$\displaystyle{\Big|\,E\bigl(\,\mathscr{T}_{A}^4 - 3\,\bigr)\,\Big| \leq 97\,\dfrac{\delta_{A}}{n}}$
\hspace{4ex}for $n \geq 4$.
\end{enumerate}
\vspace*{3.5ex}
\textbf{Proof:}\\[0.8ex]
Part \ref{EWD_3_2_02_BWa}) follows immediately from the previous 
Theorem \ref{EWD_3_2_01}, \ref{EWD_3_2_01_BWa}) and the inequality\\[2.5ex]
\refstepcounter{DSBcount}
\label{EWD_3_2_03}
\text{\hspace*{-0.8ex}(\theDSBcount)}
\hspace*{4ex}
$\displaystyle{\Big|\,\dfrac{n}{(n-1)\,(n-2)} - \dfrac{1}{n}\,\Big| \leq \dfrac{11}{n^2}}$
\hspace{4ex}for $n \geq 3$.\\[3ex]
This inequality (\ref{EWD_3_2_03}) is derived from the following property of the
polynomial $P_{1}$\\[2ex] 
\hspace*{12.1ex}$P_{1}(x) = 8\,x^2 - 31\,x +22 \geq 0$\hspace*{4ex}for $x \geq 2,9395$\hspace*{1ex},\\[2ex]
from which we get for $n \geq 3$\\[2ex]
\hspace*{12.1ex}\begin{tabular}[t]{@{}c@{\hspace*{1.6ex}}l@{}}
&$11\,n^3 - 33\,n^2 + 22\,n \geq 3\,n^3 - 2\,n^2$\\[2ex]
$\Leftrightarrow$&$\dfrac{11}{n^2} \geq \dfrac{3\,n - 2}{n\,(n-1)\,(n-2)} = \dfrac{n}{(n-1)\,(n-2)} - \dfrac{1}{n}$.
\end{tabular}\\[2.5ex]
Furthermore, since $\dfrac{11}{n} \leq 4$ is valid for $n \geq 3$, part \ref{EWD_3_2_02_BWc}) is a simple consequence of part \ref{EWD_3_2_02_BWa}).\\[2.8ex]
For part \ref{EWD_3_2_02_BWb}) we first need\\[2.5ex]
\refstepcounter{DSBcount}
\label{EWD_3_2_04}
\text{\hspace*{-0.8ex}(\theDSBcount)}
\hspace*{4ex}
$\displaystyle{\delta_{A} \geq \Bigl(\,\dfrac{n-1}{n}\,\Bigr)^2 \geq \dfrac{9}{16}}$
\hspace{4ex}for $n \geq 4$\\[2.5ex]
(cf. (\ref{EWD_3_1_19}) in the proof of Lemma \ref{EWD_3_1_18}, \ref{EWD_3_1_18_BWa})). The difference between 
the first two summands of the formula from part \ref{EWD_3_2_02_BWb}) 
and the first summand of Theorem \ref{EWD_3_2_01}, \ref{EWD_3_2_01_BWb}) is then\\[2ex]
\hspace*{12.1ex}\begin{tabular}[t]{@{}c@{\hspace*{0.8ex}}l@{}}
&$\displaystyle{\bigg|\,3\,\dfrac{(n^2 - 3\,n +1)\,(n-1)}{n\,(n-2)\,(n-3)} - 
\Bigl(\,3 + \dfrac{3}{n}\,\Bigr)\,\bigg|}$\\[3.5ex]
$=$&$3\,\dfrac{3\,n - 7}{n\,(n-2)\,(n-3)} \leq \dfrac{30}{n^2} \leq \dfrac{54}{n^2}\,\delta_{A}$
\hspace{4ex}for $n \geq 4$.
\end{tabular}\\[3ex]
The last inequality follows from (\ref{EWD_3_2_04}). The penultimate inequality is derived
analogously to above from the following property of the polynomial $P_{2}$\\[2ex] 
\hspace*{12.1ex}$P_{2}(x) = 7\,x^2 - 43\,x + 60 \geq 0$\hspace*{4ex}for $x \geq 4$\hspace*{1ex},\\[2ex]
from which we get for $n \geq 4$\\[2ex]
\hspace*{14.4ex}\begin{tabular}[t]{@{}c@{\hspace*{1.6ex}}l@{}}
&$10\,n^3 - 50\,n^2 + 60\,n \geq 3\,n^3 - 7\,n^2$\\[2ex]
\end{tabular}

\newpage

\hspace*{12.1ex}\begin{tabular}[t]{@{}c@{\hspace*{1.6ex}}l@{}}
$\Leftrightarrow$&$\dfrac{10}{n^2} \geq \dfrac{3\,n - 7}{n\,(n-2)\,(n-3)}$.
\end{tabular}\\[2.5ex]
Furthermore, for the difference between the third summand of the formula from part \ref{EWD_3_2_02_BWb}) 
and the second summand of Theorem \ref{EWD_3_2_01}, \ref{EWD_3_2_01_BWb}), 
the following holds due to Lemma \ref{EWD_3_1_18}, \ref{EWD_3_1_18_BWc})\\[2ex]
\hspace*{12.1ex}\begin{tabular}[t]{@{}c@{\hspace*{0.8ex}}l@{}}
&$\displaystyle{\bigg|\,\Bigl(\,\dfrac{3}{(n-2)\,(n-3)} - \dfrac{3}{n^2}\,\Bigr)
\sum\limits_{i, j, k = 1}^{n}\bigl(\,\bm\hat{a}_{ij}^2\,\bm\hat{a}_{ik}^2
+ \bm\hat{a}_{ij}^2\,\bm\hat{a}_{kj}^2\,\bigr)\,\bigg|}$\\[3.8ex]
$\leq$&$\displaystyle{\bigg|\,\Bigl(\,\dfrac{3}{(n-2)\,(n-3)} - \dfrac{3}{n^2}\,\Bigr)\,\bigg|\,2\,n\,\delta_{A}}$\\[3.2ex]
$=$&$6\,\dfrac{5\,n - 6}{n\,(n-2)\,(n-3)}\,\delta_{A} \leq \dfrac{168}{n^2}\,\delta_{A}$
\hspace{4ex}for $n \geq 4$.
\end{tabular}\\[3ex]
The last inequality is derived again analogously to above from\\[2ex] 
\hspace*{12.1ex}$P_{3}(x) = 23\,x^2 - 134\,x + 168 \geq 0$\hspace*{4ex}for $x \geq 4$\hspace*{1ex},\\[2ex]
from which we get for $n \geq 4$\\[2ex]
\hspace*{12.1ex}\begin{tabular}[t]{@{}c@{\hspace*{1.6ex}}l@{}}
&$28\,n^3 - 140\,n^2 + 168\,n \geq 5\,n^3 - 6\,n^2$\\[2ex]
$\Leftrightarrow$&$\dfrac{28}{n^2} \geq \dfrac{5\,n - 6}{n\,(n-2)\,(n-3)}$.
\end{tabular}\\[2.5ex]
For the difference of the terms with the sums of the fourth powers holds\\[2ex]
\hspace*{12.1ex}\begin{tabular}[t]{@{}c@{\hspace*{0.8ex}}l@{}}
&$\displaystyle{\bigg|\,\Bigl(\,\dfrac{n\,(n+1)}{(n-1)\,(n-2)\,(n-3)} - \dfrac{1}{n}\,\Bigr)
\sum\limits_{i, j = 1}^{n}\bm\hat{a}_{ij}^4\,\bigg|}$\\[3.8ex]
$=$&$\dfrac{7\,n^2 - 11\,n + 6}{n\,(n-1)\,(n-2)\,(n-3)}\,\delta_{A} \leq \dfrac{50}{n^2}\,\delta_{A}$
\hspace{4ex}for $n \geq 4$.
\end{tabular}\\[3ex]
To justify the last inequality, we consider\\[2ex] 
\hspace*{12.1ex}$P_{4}(x) = 43\,x^3 - 289\,x^2 + 544\,x - 300 \geq 0$\hspace*{4ex}for 
$x \geq 3,9864$\hspace*{1ex},\\[2ex]
from which we get for $n \geq 4$\\[2ex]
\hspace*{12.1ex}\begin{tabular}[t]{@{}c@{\hspace*{1.6ex}}l@{}}
&$50\,n^4 - 300\,n^3 + 550\,n^2 - 300\,n \geq 7\,n^4 - 11\,n^3 + 6\,n^2$\\[2ex]
$\Leftrightarrow$&$\dfrac{50}{n^2} \geq \dfrac{7\,n^2 - 11\,n + 6}{n\,(n-1)\,(n-2)\,(n-3)}$.
\end{tabular}\\[2.8ex]
Finally, the fourth summand of the Theorem \ref{EWD_3_2_01}, \ref{EWD_3_2_01_BWb}) disappears 
completely in the term $R_{4,A}$, since\\[2ex]
\hspace*{12.1ex}\begin{tabular}[t]{@{}l@{\hspace*{0.8ex}}c@{\hspace*{0.8ex}}l@{\hspace*{6.8ex}}r@{}}
$\displaystyle{\bigg|\,\sum\limits_{i, j, k, h = 1}^{n}
\bm\hat{a}_{ik}\,\bm\hat{a}_{ih}\,\bm\hat{a}_{jk}\,\bm\hat{a}_{jh}\,\bigg|}$
&$\leq$&$\displaystyle{\sum\limits_{k, h = 1}^{n}
\Bigl(\,\sum\limits_{i = 1}^{n}|\bm\hat{a}_{ik}\,\bm\hat{a}_{ih}|\,\Bigr)^2}$\\[3.5ex]
&$\leq$&$\displaystyle{\sum\limits_{k, h = 1}^{n}\biggr[\,
\Bigl(\,\sum\limits_{i = 1}^{n}\bm\hat{a}_{ik}^2\,\Bigr)\,
\Bigl(\,\sum\limits_{j = 1}^{n}\bm\hat{a}_{jh}^2\,\Bigr)\,\biggr]}$
&(H{\"o}lder's inequality\index{H{\"o}lder's inequality!for finite sequences})\\[3.5ex] 
&$=$&$\displaystyle{\Bigl(\,\sum\limits_{i, k = 1}^{n}\bm\hat{a}_{ik}^2\,\Bigr)\,
\Bigl(\,\sum\limits_{j, h = 1}^{n}\bm\hat{a}_{jh}^2\,\Bigr)}$\\[3.2ex] 
&$=$&$(n - 1)^2$&(cf. (\ref{EWD_3_1_01}))
\end{tabular}\\[2ex]
and therefore\\[2.5ex]
\hspace*{12.1ex}\begin{tabular}[t]{@{}c@{\hspace*{0.8ex}}l@{}}
&$\displaystyle{\bigg|\,\dfrac{6}{n\,(n-1)\,(n-2)\,(n-3)}\,
\sum\limits_{i, j, k, h = 1}^{n}
\bm\hat{a}_{ik}\,\bm\hat{a}_{ih}\,\bm\hat{a}_{jk}\,\bm\hat{a}_{jh}\,\bigg|}$\\[3.8ex]
$\leq$&$\dfrac{6\,(n-1)}{n\,(n-2)\,(n-3)} \leq \dfrac{36}{n^2} \leq \dfrac{64}{n^2}\,\delta_{A}$
\hspace{4ex}for $n \geq 4$.
\end{tabular}\\[3ex]
The last inequality follows from (\ref{EWD_3_2_04}). The penultimate inequality is deduced from\\[2ex]
\hspace*{12.1ex}$P_{5}(x) = 30\,x^2 - 174\,x + 216 \geq 0$\hspace*{4ex}for $x \geq 4$\hspace*{1ex},\\[2ex]
from which we get for $n \geq 4$\\[2ex]
\hspace*{12.1ex}\begin{tabular}[t]{@{}c@{\hspace*{1.6ex}}l@{}}
&$36\,n^3 - 180\,n^2 + 216\,n \geq 6\,n^3 - 6\,n^2$\\[2ex]
$\Leftrightarrow$&$\dfrac{36}{n^2} \geq \dfrac{6\,(n-1)}{n\,(n-2)\,(n-3)}$.
\end{tabular}\\[3ex]
In summary, we obtain part \ref{EWD_3_2_02_BWb}), since\\[2.8ex]
\hspace*{12.1ex}$|R_{4,A}| \leq \dfrac{54}{n^2}\,\delta_{A} + \dfrac{168}{n^2}\,\delta_{A} +
\dfrac{50}{n^2}\,\delta_{A} + \dfrac{64}{n^2}\,\delta_{A} = \dfrac{336}{n^2}\,\delta_{A}$.\\[2.5ex]
Last but not least, part \ref{EWD_3_2_02_BWd}) is derived 
from part \ref{EWD_3_2_02_BWb}) by the following calculations:\\[2.5ex]
\hspace*{12.1ex}\begin{tabular}[t]{@{}l@{\hspace*{1.6ex}}c@{\hspace*{1.6ex}}l@{\hspace*{-4.2ex}}r@{}}
$\displaystyle{\Big|\,E\bigl(\,\mathscr{T}_{A}^4 - 3\,\bigr)\,\Big|}$
&$\leq$&$|\,\lambda_{2,A}\,| + |\,R_{4,A}\,|$&(cf. (\ref{EWD_3_1_09}))\\[2.5ex]
&$\leq$&$\dfrac{3}{n} + 7\,\dfrac{\delta_{A}}{n} + \dfrac{336}{n^2}\,\delta_{A}$
&(cf. Lemma \ref{EWD_3_1_18}, \ref{EWD_3_1_18_BWd}))\\[3.3ex]
&$\leq$&$\dfrac{16}{3}\,\dfrac{\delta_{A}}{n} + 7\,\dfrac{\delta_{A}}{n}
+ 84\,\dfrac{\delta_{A}}{n} \leq 97\,\dfrac{\delta_{A}}{n}$
\hspace*{4ex}for $n \geq 4$.
\end{tabular}\\[2.8ex]
For the penultimate inequality, we used (\ref{EWD_3_2_04}) and 
\mbox{\rule[0ex]{0ex}{4ex}$\dfrac{336}{n} \leq 84$} for $n \geq 4$.
\hspace*{1ex}\hfill$\Box$\\[3ex]
The inequalities from parts \ref{EWD_3_2_02_BWc}) and \ref{EWD_3_2_02_BWd}) of Proposition \ref{EWD_3_2_02}
can be generalized in a certain way to higher moments as follows:\\[4ex]
\refstepcounter{DSBcount}
\label{EWD_3_2_05}
\textbf{\hspace*{-0.8ex}\theDSBcount\ Proposition}\\[0.8ex]
Let $k \in \mathbb{N}\setminus\{1,2\}$. Then there exists a constant $C(k) > 0$ such that 
for all matrices $A$ satisfying $\sigma_{\!A} > 0$\\[2ex]
\refstepcounter{DSBcount}
\label{EWD_3_2_06}
\text{\hspace*{-0.8ex}(\theDSBcount)}
\hspace*{4ex}
$\displaystyle{\Big|\,E\bigl(\,\mathscr{T}_{A}^k\,\bigr) - \Phi(x^k)\,\Big| \leq
C(k)\,\Bigl(\,\dfrac{\beta_{A}}{n} 
+ \dfrac{1}{n}\,\sum\limits_{i,j = 1}^{n}|\bm\hat{a}_{ij}|^k\,\Bigr)}$.\\[3.5ex]
\textbf{Proof:}\\[0.8ex]
If we proceed as in Hoeffding \cite{Hoeffding1951}, page 563, from (29), to page 566, 
we obtain the following result:\\[2.8ex]
The difference $\Big|\,E\bigl(\,\mathscr{T}_{A}^k\,\bigr) - \Phi(x^k)\,\Big|$ 
can be estimated by $c_{1}(k)$ many summands, 
each of which can in turn be estimated 
by an expression of the form\\[2.5ex]
\hspace*{12.1ex}$c_{2}(k)\,\dfrac{1}{n}$
\hspace*{2ex}or\hspace*{2ex}
$\displaystyle{c_{3}(k)
\,\Bigl(\,\dfrac{1}{n}\,\sum\limits_{i,j = 1}^{n}|\bm\hat{a}_{ij}|^{r_{1}}\,\Bigr)\, \cdot \ldots \cdot
\,\Bigl(\,\dfrac{1}{n}\,\sum\limits_{i,j = 1}^{n}|\bm\hat{a}_{ij}|^{r_{s}}\,\Bigr)}$.\\[2.5ex]
Here, $r_{1} + \ldots + r_{s} = k$; $r_{1} \geq 2, \ldots , r_{s} \geq 2$ and $r_{t} > 2$ 
for \textbf{at least} one $t \in \{1,\ldots,s\}$.\\[2.8ex]
Because of $\beta_{A} \geq \dfrac{1}{2}$ for $n \geq 2$ (cf. Lemma \ref{EWD_3_1_18}, \ref{EWD_3_1_18_BWa})), 
it follows that\\[2.5ex]
\hspace*{12.1ex}$\displaystyle{c_{2}(k)\,\dfrac{1}{n} \leq c_{4}(k)\,\Bigl(\,\dfrac{\beta_{A}}{n} 
+ \dfrac{1}{n}\,\sum\limits_{i,j = 1}^{n}|\bm\hat{a}_{ij}|^k\,\Bigr)}$
\hspace*{4ex}with $c_{4}(k) = 2\,c_{2}(k)$.\\[2.5ex]
Furthermore, if $v_{1},\ldots,v_{u}$ are those of the exponents $r_{1},\ldots,r_{s}$ that are 
\textbf{really greater} than 2, and we denote their sum by $v = v_{1} + \ldots + v_{u} \leq k$, then\\[2.5ex] 
\hspace*{7ex}\begin{tabular}[t]{@{}c@{\hspace*{0.8ex}}l@{\hspace*{-14.3ex}}r@{}}
&$\displaystyle{c_{3}(k)
\,\Bigl(\,\dfrac{1}{n}\,\sum\limits_{i,j = 1}^{n}|\bm\hat{a}_{ij}|^{r_{1}}\,\Bigr)\, \cdot \ldots \cdot
\,\Bigl(\,\dfrac{1}{n}\,\sum\limits_{i,j = 1}^{n}|\bm\hat{a}_{ij}|^{r_{s}}\,\Bigr)}$\\[3ex]
$\leq$&$\displaystyle{c_{3}(k)
\,\Bigl(\,\dfrac{1}{n}\,\sum\limits_{i,j = 1}^{n}|\bm\hat{a}_{ij}|^{v_{1}}\,\Bigr)\, \cdot \ldots \cdot
\,\Bigl(\,\dfrac{1}{n}\,\sum\limits_{i,j = 1}^{n}|\bm\hat{a}_{ij}|^{v_{u}}\,\Bigr)}$
&$\displaystyle{\Bigl(\,\text{since }\sum\limits_{i,j = 1}^{n}|\bm\hat{a}_{ij}|^{2}
= n - 1 \leq n\Bigr)}$\\[3.5ex]
$\leq$&$\displaystyle{c_{3}(k)
\,\Bigl(\,\dfrac{1}{n}\,\sum\limits_{i,j = 1}^{n}|
\bm\hat{a}_{ij}|^{v - 2\,(u-1)}\,\Bigr)^{\textstyle{\frac{v_{1} - 2}{v - 2u}}}
\, \cdot \ldots \cdot
\,\Bigl(\,\dfrac{1}{n}\,\sum\limits_{i,j = 1}^{n}|
\bm\hat{a}_{ij}|^{v - 2\,(u-1)}\,\Bigr)^{\textstyle{\frac{v_{u} - 2}{v - 2u}}}}$\\[2.5ex]
&&(cf. Lemma \ref{EWD_3_1_18}, \ref{EWD_3_1_18_BWb}))
\end{tabular}\\[2ex]
\hspace*{7ex}\begin{tabular}[t]{@{}c@{\hspace*{0.8ex}}l@{}}
$=$&$\displaystyle{c_{3}(k)
\,\Bigl(\,\dfrac{1}{n}\,\sum\limits_{i,j = 1}^{n}|
\bm\hat{a}_{ij}|^{v - 2\,(u-1)}\,\Bigr)}$
\hspace*{15.2ex}(since $(v_{1} - 2) + \ldots + (v_{u} - 2) = v - 2u$)\\[3ex]
$\leq$&$\displaystyle{c_{3}(k)\,\Bigl(\,\dfrac{\beta_{A}}{n} 
+ \dfrac{1}{n}\,\sum\limits_{i,j = 1}^{n}|\bm\hat{a}_{ij}|^k\,\Bigr)}$
\hspace*{25.8ex}(since $3 \leq v - 2\,(u-1) \leq k$).
\end{tabular}\\[-5ex]
\hspace*{1ex}\hfill$\Box$\\[0.5ex]

\section[The submatrices $B \in M(l,\hat{A})$]
{The submatrices {\boldmath $B \in M(l,\hat{A})$}}\label{EWD_Kap3_Sec3}

In this paragraph, we mainly analyze the expected value $\mu_{B}$ and the variance 
$\sigma_{\!B}^2$ for matrices $B \in M(l,\bm\hat{A})$.\index{matrix!submatrix} 
We will prove that these two values for certain matrices $A$ do not deviate {''too far''} 
from $\mu_{\bm\hat{A}}$ and $\sigma_{\!\bm\hat{A}}^2$.\\[2.8ex]
The following result serves as preparation for this.\\[4ex]
\refstepcounter{DSBcount}
\label{EWD_3_3_01}
\textbf{\hspace*{-0.8ex}\theDSBcount\ Lemma}\\[0.8ex]
If $B \in M(l,\bm\hat{A})$, then
\begin{enumerate}
\item\label{EWD_3_3_01_BWa}
$|\,\mu_{B}| \leq \dfrac{l^2}{n-l} + \dfrac{\beta_{A}}{n-l}$;
\item\label{EWD_3_3_01_BWb}
\begin{tabular}[t]{@{}c@{\hspace*{0.8ex}}l@{}}
&$\displaystyle{\Big|\,\sigma_{\!B}^2 - \dfrac{1}{n - l - 1}\sum\limits_{i,j=1}^{n}\bm\hat{a}_{ij}^2\,\Big|}$\\[3ex]
$\leq$&$\displaystyle{\dfrac{1}{n-l-1}\,\bigg\{\,\sum\limits_{\textstyle{(i \in Z_{l})\,\vee\,(j \in S_{l})}}
\bm\hat{a}_{ij}^2 + 
5\,\Bigl(\,\dfrac{l^2}{n-l} + \dfrac{\beta_{A}}{n-l}\,\Bigr)^2 +
10\,l^2 + 10\,l\,\dfrac{\beta_{A}}{n-l}\,\bigg\}}$ ,
\end{tabular}\\[2ex]
where $Z_{l} = \{\,z_{1},\ldots,z_{l}\,\}$ ($S_{l} =\{\,s_{1},\ldots,s_{l}\,\}$) 
denotes the set of indices of the \textbf{cancelled} rows (columns).
\index{matrix!submatrix!$Z_{l}$, $S_{l}$}\\[1ex] 
In addition, the symbols $\vee$ denote the logical {''or''} and $\wedge$ the logical {''and''}.
\end{enumerate}
\vspace*{1.5ex}
\textbf{Proof:}\\[0.8ex]
We have\\[1ex]
\refstepcounter{DSBcount}
\label{EWD_3_3_02}
\text{\hspace*{-0.8ex}(\theDSBcount)}
\hspace*{4ex}
\begin{tabular}[t]{@{}l@{\hspace*{0.8ex}}c@{\hspace*{0.8ex}}l@{\hspace*{12.9ex}}r@{}}
$|\,b_{\boldsymbol{.}\boldsymbol{.}}\,|$&$=$&
$\displaystyle{\bigg|\,\dfrac{1}{(n-l)^2}\sum\limits_{\textstyle{(i \in Z_{l})\,\wedge\,(j \in S_{l})}}
\bm\hat{a}_{ij}\,\bigg|}$
&(since $\bm\hat{a}_{i\boldsymbol{.}} = 0 = \bm\hat{a}_{\boldsymbol{.}j}$)\\[4.5ex]
&$\leq$&$\displaystyle{\dfrac{1}{(n-l)^2}\sum\limits_{\textstyle{(i \in Z_{l})\,\wedge\,(j \in S_{l})}}
\bigl(\,1 + |\bm\hat{a}_{ij}|^3\,\bigr)}$
&(since $|\bm\hat{a}_{ij}| \leq 1 + |\bm\hat{a}_{ij}|^3$)\\[4.5ex]
&$\leq$&$\dfrac{l^2}{(n-l)^2} + \dfrac{\beta_{A}}{(n-l)^2}$ .
\end{tabular}\\[2.5ex]
From this follows \ref{EWD_3_3_01_BWa}). Analogously to (\ref{EWD_3_3_02}) we obtain\\[2ex]
\refstepcounter{DSBcount}
\label{EWD_3_3_03}
\text{\hspace*{-0.8ex}(\theDSBcount)}
\hspace*{4ex}
\begin{tabular}[t]{@{}l@{\hspace*{0.8ex}}c@{\hspace*{0.8ex}}l@{\hspace*{5.5ex}}r@{}}
$\displaystyle{\sum\limits_{\textstyle{i \not\in Z_{l}}}b_{i\boldsymbol{.}}^2}$&$=$&
$\displaystyle{\sum\limits_{\textstyle{i \not\in Z_{l}}}\dfrac{1}{(n-l)^2}
\,\Bigl(\,\sum\limits_{\textstyle{j \in S_{l}}}\bm\hat{a}_{ij}\,\Bigr)^2}$&
(since $\bm\hat{a}_{i\boldsymbol{.}} = 0$)\\[4.5ex]
&$\leq$&$\displaystyle{\sum\limits_{\textstyle{i \not\in Z_{l}}}\dfrac{1}{(n-l)^2}
\,l\,\sum\limits_{\textstyle{j \in S_{l}}}\bm\hat{a}_{ij}^2}$&
(\index{H{\"o}lder's inequality!for finite sequences using length $\nu$}(\ref{EWD_0_1_05}) with $\nu = l$, $p = 2$)\\[4.5ex]
&$\leq$&$\displaystyle{\dfrac{l}{(n-l)^2}\sum\limits_{\textstyle{(i \not\in Z_{l})\,\,\wedge\,(j \in S_{l})}}
\bigl(\,1 + |\bm\hat{a}_{ij}|^3\,\bigr)}$&
(since $\bm\hat{a}_{ij}^2 \leq 1 + |\bm\hat{a}_{ij}|^3$)\\[4.5ex]
&$\leq$&$\dfrac{l^2}{n-l} + l\,\dfrac{\beta_{A}}{(n-l)^2}$
\end{tabular}\\[1.5ex]
and\\[1.5ex]
\refstepcounter{DSBcount}
\label{EWD_3_3_04}
\text{\hspace*{-0.8ex}(\theDSBcount)}
\hspace*{4ex}
$\displaystyle{\sum\limits_{\textstyle{j \not\in S_{l}}}b_{\boldsymbol{.}j}^2
\leq \dfrac{l^2}{n-l} + l\,\dfrac{\beta_{A}}{(n-l)^2}}$ .\\[2ex]
Also is valid\\[1.5ex]
\refstepcounter{DSBcount}
\label{EWD_3_3_05}
\text{\hspace*{-0.8ex}(\theDSBcount)}
\hspace*{4ex}
\begin{tabular}[t]{@{}c@{\hspace*{0.8ex}}l@{\hspace*{7ex}}l@{}}
&$\displaystyle{\sum\limits_{\textstyle{(i \not\in Z_{l})\,\wedge\,(j \not\in S_{l})}}
b_{ij}\,\bigl(\,b_{\boldsymbol{.}\boldsymbol{.}} - b_{i\boldsymbol{.}}
- b_{\boldsymbol{.}j}\,\bigr)}$\\[4.5ex]
$=$&$\displaystyle{b_{\boldsymbol{.}\boldsymbol{.}}\sum\limits_{\textstyle{(i \not\in Z_{l})\,\wedge\,(j \not\in S_{l})}}
b_{ij}\, - \sum\limits_{\textstyle{i \not\in Z_{l}}}b_{i\boldsymbol{.}}
\sum\limits_{\textstyle{j \not\in S_{l}}}b_{ij}\, -
\sum\limits_{\textstyle{j \not\in S_{l}}}b_{\boldsymbol{.}j}
\sum\limits_{\textstyle{i \not\in Z_{l}}}b_{ij}}$\\[5ex]
$=$&$\displaystyle{(n - l)^2\,b_{\boldsymbol{.}\boldsymbol{.}}^2 
\,-\, (n-l)\sum\limits_{\textstyle{i \not\in Z_{l}}}b_{i\boldsymbol{.}}^2
\,-\, (n-l)\sum\limits_{\textstyle{j \not\in S_{l}}}b_{\boldsymbol{.}j}^2}$ .
\end{tabular}\\[2ex]
By reapplying H{\"o}lder's inequality\index{H{\"o}lder's inequality!for finite sequences using length $\nu$} 
(\ref{EWD_0_1_05}) with $\nu = 3$ and $p = 2$ we further obtain\\[2.2ex]
\refstepcounter{DSBcount}
\label{EWD_3_3_06}
\text{\hspace*{-0.8ex}(\theDSBcount)}
\hspace*{4ex}
$\bigl(\,b_{\boldsymbol{.}\boldsymbol{.}} - b_{i\boldsymbol{.}}
- b_{\boldsymbol{.}j}\,\bigr)^2 \leq 3\,\bigl(\,b_{\boldsymbol{.}\boldsymbol{.}}^2 + b_{i\boldsymbol{.}}^2
+ b_{\boldsymbol{.}j}^2\,\bigr)$
\hspace*{4ex}for all $i \not\in Z_{l}$ and $j \not\in S_{l}$.\\[2ex]
Now, using (\ref{EWD_3_3_05}) and (\ref{EWD_3_3_06}), we get\\[2ex]
\begin{tabular}[t]{@{}c@{\hspace*{0.8ex}}l@{}}
&$\displaystyle{\Big|\,\sigma_{\!B}^2 - \dfrac{1}{n-l-1}\sum\limits_{i,j =1}^n\bm\hat{a}_{ij}^2\,\Big|}$\\[3ex]
$\leq$&$\displaystyle{\dfrac{1}{n-l-1}\,\bigg\{\sum\limits_{\textstyle{(i \in Z_{l})\,\vee\,(j \in S_{l})}}
\bm\hat{a}_{ij}^2 + 2\,\Big|\,\sum\limits_{\textstyle{(i \not\in Z_{l})\,\wedge\,(j \not\in S_{l})}}
b_{ij}\,\bigl(\,b_{\boldsymbol{.}\boldsymbol{.}} - b_{i\boldsymbol{.}}
- b_{\boldsymbol{.}j}\,\bigr)\,\Big|}$\\[4.5ex]
&\hspace*{9.8ex}$\displaystyle{+\,\sum\limits_{\textstyle{(i \not\in Z_{l})\,\wedge\,(j \not\in S_{l})}}
\bigl(\,b_{\boldsymbol{.}\boldsymbol{.}} - b_{i\boldsymbol{.}}
- b_{\boldsymbol{.}j}\,\bigr)^2\biggr\}}$\\[4.5ex]
$\leq$&$\displaystyle{\dfrac{1}{n-l-1}\,\bigg\{\sum\limits_{\textstyle{(i \in Z_{l})\,\vee\,(j \in S_{l})}}
\bm\hat{a}_{ij}^2 + 5\,\biggl(\,\displaystyle{(n - l)^2\,b_{\boldsymbol{.}\boldsymbol{.}}^2 + 
(n-l)\sum\limits_{\textstyle{i \not\in Z_{l}}}b_{i\boldsymbol{.}}^2
+ (n-l)\sum\limits_{\textstyle{j \not\in S_{l}}}b_{\boldsymbol{.}j}^2}\,\biggr)\biggr\}}$. 
\end{tabular}\\[3ex]
This and (\ref{EWD_3_3_02}), (\ref{EWD_3_3_03}), (\ref{EWD_3_3_04}) finally give 
part \ref{EWD_3_3_01_BWb}).
\hspace*{1ex}\hfill$\Box$\\[3ex] 
A direct consequence of Lemma \ref{EWD_3_3_01} is:\\[4ex]
\refstepcounter{DSBcount}
\label{EWD_3_3_07}
\textbf{\hspace*{-0.8ex}\theDSBcount\ Corollary}\\[0.8ex]
There exists a $0 < \epsilon_{0} \leq 1$ and an $n_{0} \in \mathbb{N}$ such that 
for all matrices $A$ satisfying $\sigma_{\!A} > 0$, $\dfrac{\beta_{A}}{n} \leq \epsilon_{0}$ 
and $n \geq n_{0}$ holds:\\[2ex]
If $1 \leq l \leq 17$ and $B \in M(l,\bm\hat{A})$, then\\[2ex]
\refstepcounter{DSBcount}
\label{EWD_3_3_08}
\text{\hspace*{-0.8ex}(\theDSBcount)}
\hspace*{4ex}
$|\,\mu_{B}\,| \leq 1$,\\[2ex]
\refstepcounter{DSBcount}
\label{EWD_3_3_09}
\text{\hspace*{-0.8ex}(\theDSBcount)}
\hspace*{4ex}
$\big|\,\sigma_{\!B}^2 - 1\,\big| \leq \dfrac{1}{3}$.\\[3.5ex]
\textbf{Proof:}\\[0.8ex]
If $\epsilon_{0}$ is sufficiently small and $n_{0}$ is sufficiently large, then (\ref{EWD_3_3_08}) 
follows from part \ref{EWD_3_3_01_BWa}) of the previous Lemma \ref{EWD_3_3_01}.\\[2.8ex] 
Furthermore, an application of 
H{\"o}lder's inequality\index{H{\"o}lder's inequality!for finite sequences using length $\nu$} 
(\ref{EWD_0_1_05}) with $\nu \leq 2\,l\,n$ and $p = \dfrac{3}{2}$ gives\\[2ex]
\hspace*{12.1ex}\begin{tabular}[t]{@{}l@{\hspace*{0.8ex}}c@{\hspace*{0.8ex}}l@{}}
$\displaystyle{\sum\limits_{\textstyle{(i \in Z_{l})\,\vee\,(j \in S_{l})}}
\bm\hat{a}_{ij}^2}$&$\leq$&
$\displaystyle{(2\,l\,n)^{1/3}\,\Bigl(\,\sum\limits_{\textstyle{(i \in Z_{l})\,\vee\,(j \in S_{l})}}
|\bm\hat{a}_{ij}|^3\,\Bigr)^{2/3}}$\\[4.5ex]
&$\leq$&$4\,n\,\Bigl(\,\dfrac{\beta_{A}}{n}\,\Bigr)^{2/3}$
\hspace*{4ex}for $1 \leq l \leq 17$.
\end{tabular}\\[2.8ex]
If we insert this estimate into part \ref{EWD_3_3_01_BWb}) of the previous Lemma \ref{EWD_3_3_01}) 
and use $\displaystyle{\sum\limits_{i,j=1}^{n}\bm\hat{a}_{ij}^2 = n - 1}$, we further obtain\\[2.3ex]
\hspace*{12.1ex}\begin{tabular}[t]{@{}l@{\hspace*{0.8ex}}c@{\hspace*{0.8ex}}l@{}}
$\big|\,\sigma_{\!B}^2 - 1\,\big|$&$\leq$&
$\displaystyle{\Big|\,\sigma_{\!B}^2 - \dfrac{1}{n - l - 1}\sum\limits_{i,j=1}^{n}\bm\hat{a}_{ij}^2\,\Big| +
\Big|\,\dfrac{1}{n - l - 1}\sum\limits_{i,j=1}^{n}\bm\hat{a}_{ij}^2 - 1\,\Big|}$\\[4ex]
&$\leq$&
$\displaystyle{\dfrac{1}{n - l - 1}\,\bigg\{\,4\,n\,\Bigl(\,\dfrac{\beta_{A}}{n}\,\Bigr)^{2/3} + 
5\,\Bigl(\,\dfrac{l^2}{n-l} + \dfrac{\beta_{A}}{n-l}\,\Bigr)^2 +
10\,l^2 + 10\,l\,\dfrac{\beta_{A}}{n-l}\,\bigg\}}$\\[4ex] 
&&\hspace*{25.5ex}$\displaystyle{+\, \Big|\,\dfrac{n - 1}{n - l - 1} - 1\,\Big|}$ ,
\end{tabular}\\[3.5ex]
Thus we get (\ref{EWD_3_3_09}) for sufficiently small $\epsilon_{0}$ and sufficiently large $n_{0}$.
\hspace*{1ex}\hfill$\Box$\\[3ex]
For matrices $B \in M(l,\bm\hat{A})$ satisfying $\sigma_{\!B}^2 \geq \dfrac{2}{3}$ the following  
very useful relationship between $\bm\hat{B}$ and $\bm\hat{A}$ holds:\\[4ex]
\refstepcounter{DSBcount}
\label{EWD_3_3_10}
\textbf{\hspace*{-0.8ex}\theDSBcount\ Lemma}\\[0.8ex]
Let $k \in \mathbb{N}$ and $B \in M(l,\bm\hat{A})$ satisfying $\sigma_{\!B}^2 \geq \dfrac{2}{3}$. Then\\[2ex]
\refstepcounter{DSBcount}
\label{EWD_3_3_11}
\text{\hspace*{-0.8ex}(\theDSBcount)}
\hspace*{4ex}
$\displaystyle{\sum\limits_{\textstyle{(i \not\in Z_{l})\,\wedge\,(j \not\in S_{l})}}
|\bm\hat{b}_{ij}|^k \leq 5^k \sum\limits_{i,j = 1}^{n} |\bm\hat{a}_{ij}|^k}$.\\[3.5ex]
\textbf{Proof:}\\[0.8ex]
At first, an application of H{\"o}lder's inequality\index{H{\"o}lder's inequality!for finite sequences using length $\nu$}
(\ref{EWD_0_1_05}) with $\nu = n - l$ and $p = k$ gives\\[2ex]
\hspace*{2ex}\begin{tabular}[t]{@{}l@{\hspace*{0.8ex}}c@{\hspace*{0.8ex}}l@{}}
$\displaystyle{(n - l)\sum\limits_{\textstyle{i \not\in Z_{l}}}|b_{i\boldsymbol{.}}|^k}$&$\leq$&
$\displaystyle{\dfrac{1}{(n-l)^{k-1}}\sum\limits_{\textstyle{i \not\in Z_{l}}}
\Bigr(\,\sum\limits_{\textstyle{j \not\in S_{l}}}|\bm\hat{a}_{ij}|\,\Bigl)^k}$\\[4.5ex]
&$\leq$&$\displaystyle{\dfrac{1}{(n-l)^{k-1}}\sum\limits_{\textstyle{i \not\in Z_{l}}}
\Bigr(\,(n-l)^{k-1}\sum\limits_{\textstyle{j \not\in S_{l}}}|\bm\hat{a}_{ij}|^k\,\Bigl)
\,\,\leq \sum\limits_{i,j = 1}^{n} |\bm\hat{a}_{ij}|^k}$.
\end{tabular}\\[2.5ex]
Analogously follows\\[1.5ex]
\hspace*{11ex}$\displaystyle{(n - l)\sum\limits_{\textstyle{j \not\in S_{l}}}|b_{\boldsymbol{.}j}|^k
\,\leq \sum\limits_{i,j = 1}^{n} |\bm\hat{a}_{ij}|^k}$
\hspace*{4ex}and\hspace*{4ex}
$\displaystyle{(n - l)^2\,
|b_{\boldsymbol{.}\boldsymbol{.}}|^k
\,\leq \sum\limits_{i,j = 1}^{n} |\bm\hat{a}_{ij}|^k}$.\\[2.5ex]
A further application of H{\"o}lder's inequality\index{H{\"o}lder's inequality!for finite sequences using length $\nu$} 
(\ref{EWD_0_1_05}) with $\nu = 4$ and $p = k$ finally leads to\\[2.5ex]
\hspace*{2ex}\begin{tabular}[t]{@{}l@{\hspace*{0.8ex}}c@{\hspace*{0.8ex}}l@{}}
$\displaystyle{\sum\limits_{\textstyle{(i \not\in Z_{l})\,\wedge\,(j \not\in S_{l})}}
|\bm\hat{b}_{ij}|^k}$&$\leq$&
$\displaystyle{\Bigl(\,\dfrac{3}{2}\,\Bigr)^{k/2}\,4^{k-1}
\sum\limits_{\textstyle{(i \not\in Z_{l})\,\wedge\,(j \not\in S_{l})}}
\Bigl(\,|b_{ij}|^k + |b_{i\boldsymbol{.}}|^k 
+ |b_{\boldsymbol{.}j}|^k + |b_{\boldsymbol{.}\boldsymbol{.}}|^k\,\Bigr)}$\\[3ex] 
&$\leq$&$\displaystyle{5^k \sum\limits_{i,j = 1}^{n} |\bm\hat{a}_{ij}|^k}$
\hspace*{22.2ex}(since $\sqrt{\dfrac{3}{2}}\,4 = 4,8990 \leq 5$).
\end{tabular}\\[-5ex]
\hspace*{1ex}\hfill$\Box$\\[0.5ex]

\section[The truncated matrices $A'$ and $\bar{A}$]
{The truncated matrices {\boldmath $A'$} and {\boldmath $\bar{A}$}}\label{EWD_Kap3_Sec4}

In the following, many results can only be shown for matrices whose elements are bounded in 
absolute value by one. However, in order to be able to prove the two Theorems \ref{EWD_3_1_10} 
and \ref{EWD_3_1_13} for general matrices, in this section, we will first create the 
truncated matrices $A'$ and $\bm\bar{A}$ from $\bm\hat{A}$ 
and then establish some useful connections between these matrices.\\[2.8ex]
Let\\[2ex]
\hspace*{12.1ex}$\displaystyle{a_{ij}' =
\left\{
\begin{array}{ll@{}}
\bm\hat{a}_{ij}& \hspace*{3ex}
\text{if}\ \ \ |\bm\hat{a}_{ij}|\, \leq\, \dfrac{1}{2},\\[2ex]
0& \hspace*{3ex}
\text{if}\ \ \ |\bm\hat{a}_{ij}|\, >\, \dfrac{1}{2}.
\end{array}  \right.}$\index{matrix!truncated}
\index{matrix!truncated!$a_{ij}'$, $d_{ij}$, $\bm\bar{a}_{ij}$}\\[2.5ex]
In addition, let\\[2ex]
\hspace*{12.1ex}$\Gamma = \bigl\{(i,j)\,:\,|\bm\hat{a}_{ij}|\, >\, \dfrac{1}{2}\,\bigr\}$ ,\\[1.5ex]
\hspace*{12.1ex}$d_{ij} = a_{i\boldsymbol{.}}' + a_{\boldsymbol{.}j}' - a_{\boldsymbol{.}\boldsymbol{.}}'$ ,\\[1.5ex]
\hspace*{12.1ex}$\bm\bar{a}_{ij} = \bm\hat{a}_{ij}' = \dfrac{1}{\sigma_{\!A'}}\,(a_{ij}' - d_{ij})$
\index{matrix!truncated!$\Gamma$}\\[2.5ex]
and $A' = (a_{ij}')$, $D = (d_{ij})$, $\bm\bar{A} = (\bm\bar{a}_{ij})$ the corresponding
$n{\times}n-$matrices.\index{matrix!truncated!$A'$, $D$, $\bm\bar{A}$} 
We will first list some elementary results.\\[4ex]
\refstepcounter{DSBcount}
\label{EWD_3_4_01}
\textbf{\hspace*{-0.8ex}\theDSBcount\ Lemma}\\[0.8ex]
Let $k, r \in \mathbb{N}$ and\\[2ex]
\hspace*{12.1ex}$\displaystyle{A[k] = 2^k\sum\limits_{i,j =1}^n|\bm\hat{a}_{ij}|^k}$.
\index{matrix!truncated!$A[k]$}\\[1.5ex]
Then
\begin{enumerate}
\item\label{EWD_3_4_01_BWa}
$|\Gamma| \leq A[k]$,
\item\label{EWD_3_4_01_BWb}
$P(\mathscr{T}_{A} \not= T_{A'}) \leq \dfrac{|\Gamma|}{n} \leq \dfrac{A[k]}{n}$,
\item\label{EWD_3_4_01_BWc}
$\displaystyle{\sum\limits_{(i,j)\, \in\, \Gamma}|\bm\hat{a}_{ij}|^r \leq 2^{-r}\,A[k]}$
\hspace*{3ex}for $r \leq k$,
\item\label{EWD_3_4_01_BWd}
$\displaystyle{\sum\limits_{i=1}^n |a_{i\boldsymbol{.}}'| \leq \dfrac{A[k]}{2\,n}}$,\hspace*{2ex}
$\displaystyle{\sum\limits_{j=1}^n |a_{\boldsymbol{.}j}'| \leq \dfrac{A[k]}{2\,n}}$
\hspace*{2ex}and\hspace*{2ex}
$\displaystyle{|a_{\boldsymbol{.}\boldsymbol{.}}'| \leq \dfrac{A[k]}{2\,n^2}}$,
\item\label{EWD_3_4_01_BWe}
$\displaystyle{\sum\limits_{i,j=1}^n |d_{ij}| \leq \dfrac{3}{2}\,A[k]}$,
\item\label{EWD_3_4_01_BWf}
$|d_{ij}|\leq \dfrac{3}{2}\,\text{min}\Bigl\{\,1, \dfrac{A[k]}{n}\,\Bigr\}$
\hspace*{3ex}for all $1 \leq i, j \leq n$,
\item\label{EWD_3_4_01_BWg}
$|\mu_{A'}| \leq \dfrac{A[k]}{2\,n}$,
\item\label{EWD_3_4_01_BWh}
$\big|\,\sigma_{\!A'}^2 - 1\,\big| \leq 8\,\dfrac{A[k]}{n}$
\hspace*{3ex}for $k \geq 2$.
\end{enumerate}
\textbf{Proof:}\\[0.8ex]
This proof is based on arguments of Bolthausen \cite{Bolthausen1984}, page 382.
\begin{enumerate}
\item
$\displaystyle{\sum\limits_{i,j =1}^n|\bm\hat{a}_{ij}|^k \geq
\sum\limits_{(i,j)\, \in\, \Gamma}|\bm\hat{a}_{ij}|^k \geq \dfrac{1}{2^k}\,|\Gamma|}$.                                              
\item
For $i \in \{1,\ldots,n\}$ we define the following $i$-section of $\Gamma$:\\[2ex]
\hspace*{12.1ex}$\Gamma_{i} = \bigl\{\,j \in \{1,\ldots,n\}\,:\,(i,j) \in \Gamma\,\bigr\}$.\\[2ex]
Since $\pi(i)$ is \textbf{equally distributed} on $\{1,\ldots,n\}$ for all $i$, we obtain\\[2ex] 
\begin{tabular}[t]{@{}l@{\hspace*{0.8ex}}c@{\hspace*{0.8ex}}l@{\hspace*{37.1ex}}r@{}}
$P(\mathscr{T}_{A} \not= T_{A'})$&$=$&$P(T_{\bm\hat{A}} \not= T_{\raisebox{-2.1pt}{$\scriptstyle{\!A'}$}})$\\[1.3ex]
&$\leq$&$\displaystyle{P\Bigl(\,\sum\limits_{i = 1}^n 1_{\Gamma}(i,\pi(i)) \geq 1\,\Bigr)}$\\[3.2ex]
&$\leq$&$\displaystyle{E\Bigl(\,\sum\limits_{i = 1}^n 1_{\Gamma}(i,\pi(i))\,\Bigr)}$\\[3.2ex]
&$=$&$\displaystyle{\sum\limits_{i = 1}^n P\Bigl(\,\pi(i) \in \Gamma_{i}\,\Bigr)}$\\[3.2ex]
&$=$&$\displaystyle{\sum\limits_{i = 1}^n \dfrac{|\Gamma_{i}|}{n}} = \dfrac{|\Gamma|}{n} \leq \dfrac{A[k]}{n}$
&(cf. part \ref{EWD_3_4_01_BWa})).
\end{tabular}
\item
\begin{tabular}[t]{@{}l@{\hspace*{0.8ex}}c@{\hspace*{0.8ex}}l@{\hspace*{20ex}}r@{}}
$\displaystyle{\sum\limits_{(i,j)\, \in\, \Gamma}|\bm\hat{a}_{ij}|^r}$&$\leq$&
$\displaystyle{|\Gamma|^{(k-r)/k}\,\Bigl(\,\sum\limits_{i,j =1}^n|\bm\hat{a}_{ij}|^k\,\Bigr)^{r/k}}$
&(\index{H{\"o}lder's inequality!for finite sequences using length $\nu$}(\ref{EWD_0_1_05}) with $\nu = |\Gamma|$, $p = \dfrac{k}{r}$)\\[2ex] 
&$\leq$&$2^{-r}\,A[k]$
&(cf. part \ref{EWD_3_4_01_BWa})).
\end{tabular}
\item
Using the notation $\Gamma_{i}$ from the proof of part \ref{EWD_3_4_01_BWb}), we obtain\\[2ex]
\begin{tabular}[t]{@{}l@{\hspace*{0.8ex}}c@{\hspace*{0.8ex}}l@{\hspace*{26.8ex}}r@{}}
$\displaystyle{\sum\limits_{i =1}^n |a_{i\boldsymbol{.}}'|}$
&$=$&$\displaystyle{\sum\limits_{i =1}^n |a_{i\boldsymbol{.}}' - 
\bm\hat{a}_{\raisebox{-1.4pt}{$\scriptstyle{i\boldsymbol{.}}$}}|}$
&(since $\bm\hat{a}_{i\boldsymbol{.}} = 0$ (cf. (\ref{EWD_3_1_01}))\\[3.5ex]
&$=$&$\displaystyle{\sum\limits_{i =1}^n \Big|\,\dfrac{1}{n}\sum\limits_{j =1}^n a_{ij}' - 
\dfrac{1}{n}\sum\limits_{j =1}^n \bm\hat{a}_{ij}\,\Big|}$\\[3.5ex]
&$\leq$&$\displaystyle{\dfrac{1}{n}\sum\limits_{i =1}^n \sum\limits_{j\,\in\,\Gamma_{i}}|\bm\hat{a}_{ij}|}$\\[3.5ex]
&$=$&$\displaystyle{\dfrac{1}{n}\sum\limits_{(i,j)\, \in\, \Gamma} |\bm\hat{a}_{ij}| \leq \dfrac{A[k]}{2\,n}}$
&(cf. part \ref{EWD_3_4_01_BWc}) with $r = 1$).
\end{tabular}\\[2ex]
Similarly, $\displaystyle{\sum\limits_{j =1}^n |a_{\boldsymbol{.}j}'| \leq \dfrac{A[k]}{2\,n}}$ and
$|a_{\boldsymbol{.}\boldsymbol{.}}'| \leq \dfrac{A[k]}{2\,n^2}$.
\item
Due to part \ref{EWD_3_4_01_BWd}) is valid\hspace*{3ex}
\begin{tabular}[t]{@{}l@{\hspace*{0.8ex}}c@{\hspace*{0.8ex}}l@{}}
$\displaystyle{\sum\limits_{i,j=1}^n |d_{ij}|}$&$\leq$&
$\displaystyle{n\sum\limits_{i=1}^n |a_{i\boldsymbol{.}}'| +
n\sum\limits_{j=1}^n |a_{\boldsymbol{.}j}'| + n^2\,|a_{\boldsymbol{.}\boldsymbol{.}}'|
\leq \dfrac{3}{2}\,A[k]}$. 
\end{tabular}
\item
\begin{tabular}[t]{@{}l@{\hspace*{3ex}}l@{}}
We have $|a_{ij}'| \leq \dfrac{1}{2}$ and thus&
$|d_{ij}| \leq |a_{i\boldsymbol{.}}'| + |a_{\boldsymbol{.}j}'| + |a_{\boldsymbol{.}\boldsymbol{.}}'|
\leq \dfrac{1}{2} + \dfrac{1}{2} + \dfrac{1}{2} = \dfrac{3}{2}$.\\[2ex]
Due to part \ref{EWD_3_4_01_BWd}) also holds&
$|d_{ij}| \leq |a_{i\boldsymbol{.}}'| + |a_{\boldsymbol{.}j}'| + |a_{\boldsymbol{.}\boldsymbol{.}}'|
\leq \dfrac{A[k]}{2\,n} + \dfrac{A[k]}{2\,n} + \dfrac{A[k]}{2\,n^2} 
\leq \dfrac{3}{2}\,\dfrac{A[k]}{n}$.
\end{tabular}
\item
Due to part \ref{EWD_3_4_01_BWd}) is valid\hspace*{3ex}
$|\mu_{\raisebox{-2.1pt}{$\scriptstyle{\!A'}$}}| =
|n\,a_{\boldsymbol{.}\boldsymbol{.}}'| \leq  \dfrac{A[k]}{2\,n}$.
\item
\begin{tabular}[t]{@{}l@{\hspace*{0.8ex}}c@{\hspace*{0.8ex}}l@{\hspace*{-1.7ex}}r@{}}
$\big|\,\sigma_{\raisebox{-2.5pt}{$\scriptstyle{\!A'}$}}^2 - 1\,\big|$&$=$&
$\big|\,\sigma_{\raisebox{-2.5pt}{$\scriptstyle{\!A'}$}}^2 - \sigma_{\!\bm\hat{A}}^2\,\big|$
&(since $\sigma_{\!\bm\hat{A}}^2 = 1$ (cf. (\ref{EWD_3_1_01}))\\[2ex]
&$=$&$\displaystyle{\dfrac{1}{n-1}\,\Big|\,\sum\limits_{i,j = 1}^n(a_{ij}' - d_{ij})^2 - \bm\hat{a}_{ij}^2\,\Big|}$\\[3ex]
&$=$&$\displaystyle{\dfrac{1}{n-1}\,\Big|\,\sum\limits_{i,j = 1}^n
\bigl(\,a_{ij}'^{\,2} - \bm\hat{a}_{ij}^2\,\bigr) - 2\,a_{ij}'\,d_{ij} + d_{ij}^2\,\Big|}$\\[3ex]
&$\leq$&$\displaystyle{\dfrac{1}{n-1}\,\Bigl(\,\sum\limits_{(i,j)\, \in\, \Gamma}\bm\hat{a}_{ij}^2
+ \dfrac{5}{2}\sum\limits_{i,j = 1}^n|d_{ij}|\,\Bigr)}$
&(since $|a_{ij}'| \leq \dfrac{1}{2}$, $|d_{ij}| \leq \dfrac{3}{2}$ (cf. part \ref{EWD_3_4_01_BWf}))\\[3.5ex]
&$\leq$&$4\,\dfrac{A[k]}{n-1} \leq 8\,\dfrac{A[k]}{n}$
&(cf. part \ref{EWD_3_4_01_BWc}), part \ref{EWD_3_4_01_BWe}) and $n \geq 2$).\hspace*{1.2ex}$\Box$
\end{tabular}
\end{enumerate}
\vspace*{3ex}
An important consequence of Lemma \ref{EWD_3_4_01} is:\\[4ex]
\refstepcounter{DSBcount}
\label{EWD_3_4_02}
\textbf{\hspace*{-0.8ex}\theDSBcount\ Corollary}\\[0.8ex]
There exists a $0 < \epsilon_{0} \leq 1$ such that for all matrices $A$ satisfying 
$\sigma_{\!A} > 0$ and $\dfrac{\beta_{A}}{n} \leq \epsilon_{0}$ holds:\\[2ex]
\refstepcounter{DSBcount}
\label{EWD_3_4_03}
\text{\hspace*{-0.8ex}(\theDSBcount)}
\hspace*{4ex}
$\big|\,\sigma_{\!A'}^2 - 1\,\big| \leq \dfrac{1}{3}$,\\[2ex]
\refstepcounter{DSBcount}
\label{EWD_3_4_04}
\text{\hspace*{-0.8ex}(\theDSBcount)}
\hspace*{4ex}
$|\bm\bar{a}_{ij}| \leq 1$
\hspace*{4ex}for all $1 \leq i, j \leq n$.\\[3.5ex]
\textbf{Proof:}\\[0.8ex]
Let \mbox{\rule[-3ex]{0ex}{6ex}$\epsilon_{0} = \dfrac{1}{192}$}. 
Then we obtain (\ref{EWD_3_4_03}) by applying Lemma \ref{EWD_3_4_01}, \ref{EWD_3_4_01_BWh}) with $k = 3$.
From this and because of Lemma \ref{EWD_3_4_01}, \ref{EWD_3_4_01_BWf}) with $k = 3$, we also get (\ref{EWD_3_4_04}):\\[2ex]
\hspace*{12.1ex}$|\bm\bar{a}_{ij}| \leq
\dfrac{1}{\sigma_{\!A'}}\,\bigl(\,|a_{ij}'| + |d_{ij}|\,\bigr) \leq \sqrt{\dfrac{3}{2}}\,
\Bigl(\,\dfrac{1}{2} + 12\,\dfrac{\beta_{A}}{n}\,\Bigr) \leq \sqrt{\dfrac{3}{2}}\,\dfrac{9}{16} \leq 1$.
\hspace*{1ex}\hfill$\Box$\\[3ex]
Now, under the condition \mbox{\rule[-2.5ex]{0ex}{5.5ex}$\sigma_{\!A'}^2 \geq \dfrac{2}{3}$}, 
which follows from (\ref{EWD_3_4_03}), 
we can derive some useful relations between the matrices $\bm\hat{A}$ and $\bm\bar{A}$.
The first result is an analogue of Lemma \ref{EWD_3_3_10}.\\[3.5ex]
\refstepcounter{DSBcount}
\label{EWD_3_4_05}
\textbf{\hspace*{-0.8ex}\theDSBcount\ Lemma}\\[0.8ex]
Let $k \in \mathbb{N}$. Then, for each matrix $A$ satisfying $\sigma_{\!A'}^2 \geq \dfrac{2}{3}$,\\[2ex]
\hspace*{12.1ex}$\displaystyle{\sum\limits_{i, j = 1}^{n} |\bm\bar{a}_{ij}|^k \leq
5^k \sum\limits_{i,j = 1}^{n} |\bm\hat{a}_{ij}|^k}$.\\[3.5ex]
\textbf{Proof:}\\[0.8ex]
We use $|a_{ij}'| \leq |\bm\hat{a}_{ij}|$ for all $1 \leq i, j \leq n$ and
proceed in the same way as in the proof of Lemma \ref{EWD_3_3_10} (with $l = 0$).
\hspace*{1ex}\hfill$\Box$\\[3ex]
Since the expansions $e_{1,A}$ and $e_{2,A}$ contain the terms 
\mbox{\rule[-3.8ex]{0ex}{7ex}$\displaystyle{\dfrac{1}{n}\sum\limits_{i,j = 1}^{n} \bm\hat{a}_{ij}^{3}}$,} 
\mbox{\rule[-3.8ex]{0ex}{7ex}$\displaystyle{\dfrac{1}{n}\sum\limits_{i,j = 1}^{n} \bm\hat{a}_{ij}^{4}}$,} 
\mbox{\rule[-3.8ex]{0ex}{7ex}$\displaystyle{\Bigl(\,\dfrac{1}{n}\sum\limits_{i,j = 1}^{n} \bm\hat{a}_{ij}^{3}\,\Bigr)^2}$,} etc., the following estimates are also of interest for the further procedure.\\[4ex]
\refstepcounter{DSBcount}
\label{EWD_3_4_06}
\textbf{\hspace*{-0.8ex}\theDSBcount\ Lemma}\\[0.8ex]
Let $k \in \mathbb{N}$. Then there exist positive constants $C_{1}(k)$, $C_{2}(k)$, $C_{3}(k)$ and $C_{4}(k)$ 
such that for all matrices $A$ satisfying \mbox{\rule[-2ex]{0ex}{6ex}$\sigma_{\!A'}^2 \geq \dfrac{2}{3}$} holds:
\begin{enumerate}
\item\label{EWD_3_4_06_BWa}
$\displaystyle{\bigg|\,\dfrac{1}{n}\sum\limits_{i, j = 1}^{n}\bm\hat{a}_{ij}^3 -   
\dfrac{1}{n}\sum\limits_{i, j = 1}^{n}\bm\bar{a}_{ij}^3\,\bigg| \leq C_{1}(k)\,
\dfrac{1}{n}\sum\limits_{i, j = 1}^{n}|\bm\hat{a}_{ij}|^k}$
\hspace*{4ex}for $k \geq 3$.\vspace*{0.5ex}
\item\label{EWD_3_4_06_BWb}
$\displaystyle{\bigg|\,\dfrac{1}{n}\sum\limits_{i, j = 1}^{n}\bm\hat{a}_{ij}^4 -   
\dfrac{1}{n}\sum\limits_{i, j = 1}^{n}\bm\bar{a}_{ij}^4\,\bigg| \leq C_{2}(k)\,
\dfrac{1}{n}\sum\limits_{i, j = 1}^{n}|\bm\hat{a}_{ij}|^k}$
\hspace*{4ex}for $k \geq 4$.\vspace*{0.5ex}
\item\label{EWD_3_4_06_BWc}
$\displaystyle{\bigg|\,\dfrac{1}{n^2}\sum\limits_{i, j, r = 1}^{n}\bm\hat{a}_{ij}^2\,\bm\hat{a}_{ir}^2 -   
\dfrac{1}{n^2}\sum\limits_{i, j, r = 1}^{n}\bm\bar{a}_{ij}^2\,\bm\bar{a}_{ir}^2\,\bigg| \leq C_{3}(k)\,
\dfrac{1}{n}\sum\limits_{i, j = 1}^{n}|\bm\hat{a}_{ij}|^k}$
\hspace*{4ex}for $k \geq 4$.\\[2ex]
The same is true for $\displaystyle{\sum\limits_{i, j, r = 1}^{n}
\bm\hat{a}_{ij}^2\,\bm\hat{a}_{rj}^2}$.
\item\label{EWD_3_4_06_BWd}
$\displaystyle{\bigg|\,\Bigl(\,\dfrac{1}{n}\sum\limits_{i, j = 1}^{n}\bm\hat{a}_{ij}^3\,\Bigr)^2 -   
\Bigl(\,\dfrac{1}{n}\sum\limits_{i, j = 1}^{n}\bm\bar{a}_{ij}^3\,\Bigr)^2\,\bigg| \leq C_{4}(k)\,
\dfrac{1}{n}\sum\limits_{i, j = 1}^{n}|\bm\hat{a}_{ij}|^k}$
\hspace*{4ex}for $k \geq 4$.
\end{enumerate}
\vspace*{2.5ex}
\textbf{Proof:}
\begin{enumerate}
\item
Let $k \geq 3$ be fixed. Then an application of the triangle inequality gives\\[2ex]
\refstepcounter{DSBcount}
\label{EWD_3_4_07}
\text{\hspace*{-0.8ex}(\theDSBcount)}
\hspace*{4ex}
\begin{tabular}[t]{@{}c@{\hspace*{0.8ex}}l@{}}
&$\displaystyle{\bigg|\,\dfrac{1}{n}\sum\limits_{i, j = 1}^{n}\bm\hat{a}_{ij}^3 -   
\dfrac{1}{n}\sum\limits_{i, j = 1}^{n}\bm\bar{a}_{ij}^3\,\bigg|}$\\[3.5ex]
$\leq$&$\displaystyle{\Big|\,1 - \dfrac{1}{\sigma_{\!A'}^3}\,\Big|\,\dfrac{\beta_{A}}{n} 
+ \dfrac{1}{\sigma_{\!A'}^3}\,\bigg|\,\dfrac{1}{n}\sum\limits_{i, j = 1}^{n}\bm\hat{a}_{ij}^3 -   
\dfrac{1}{n}\sum\limits_{i, j = 1}^{n}\bigl(\,a_{ij}' - d_{ij}\,\bigr)^3\,\bigg|}$.
\end{tabular}\\[2.5ex]
To estimate the first summand of (\ref{EWD_3_4_07}), we first consider\\[2.8ex]
\hspace*{2ex}\begin{tabular}[t]{@{}l@{\hspace*{0.8ex}}c@{\hspace*{0.8ex}}l@{\hspace*{7.7ex}}r@{}}
$\displaystyle{\Big|\,1 - \dfrac{1}{\sigma_{\!A'}^3}\,\Big|}$&$\leq$&
$\displaystyle{\dfrac{3}{2}\,\text{max}\Bigl\{\,\dfrac{1}{\sigma_{\!A'}^5}, 1\,\Bigr\} 
\big|\,\sigma_{\!A'}^2 - 1\,\big|}$
&(mean value theorem\index{mean value theorem} for $f(x) = \dfrac{1}{x^{3/2}}$)\\[3ex]
&$\leq$&$\displaystyle{\Bigl(\,\dfrac{3}{2}\,\Bigr)^{7/2}\,2^{k+2}
\,\dfrac{1}{n}\sum\limits_{i, j = 1}^{n}|\bm\hat{a}_{ij}|^{k-1}}$
&(due to Lemma \ref{EWD_3_4_01}, \ref{EWD_3_4_01_BWh}) and $\sigma_{\!A'}^2 \geq \dfrac{2}{3}$).
\end{tabular}\\[2.5ex]
We claim that, using the definition $c_{1}(k) = \Bigl(\,\dfrac{3}{2}\,\Bigr)^{7/2}\,2^{k+2}$, we get\\[2ex]
\hspace*{2ex}$\displaystyle{\Big|\,1 - \dfrac{1}{\sigma_{\!A'}^3}\,\Big|\,\dfrac{\beta_{A}}{n}
\leq c_{1}(k)\,\dfrac{1}{n}\sum\limits_{i, j = 1}^{n}|\bm\hat{a}_{ij}|^{k}}$.\\[1ex]
This inequality is obtained in the case $k = 3$ because of
\mbox{\rule[-3.8ex]{0ex}{7ex}$\displaystyle{\dfrac{1}{n}\sum\limits_{i, j = 1}^{n}|\bm\hat{a}_{ij}|^{k-1} 
= \dfrac{n-1}{n} \leq 1}$} 
and in the case $k > 3$ by a double application of Lemma \ref{EWD_3_1_18}, \ref{EWD_3_1_18_BWb}), 
which yields\\[2.5ex]
\refstepcounter{DSBcount}
\label{EWD_3_4_08}
\text{\hspace*{-0.8ex}(\theDSBcount)}
\hspace*{4ex}
\begin{tabular}[t]{@{}c@{\hspace*{0.8ex}}l@{}}
&$\displaystyle{\Bigl(\,\dfrac{1}{n}\sum\limits_{i, j = 1}^{n}|\bm\hat{a}_{ij}|^{k-1}\,\Bigr)
\,\Bigl(\,\dfrac{1}{n}\sum\limits_{i, j = 1}^{n}|\bm\hat{a}_{ij}|^{3}\,\Bigr)}$\\[3ex]
$\leq$&
$\displaystyle{\Bigl(\,\dfrac{1}{n}
\sum\limits_{i, j = 1}^{n} |\,\bm\hat{a}_{ij}\,|^k\,\Bigr)^{(k -3)/(k-2)}
\,\Bigl(\,\dfrac{1}{n}
\sum\limits_{i, j = 1}^{n} |\,\bm\hat{a}_{ij}\,|^k\,\Bigr)^{1/(k-2)} 
= \dfrac{1}{n}\sum\limits_{i, j = 1}^{n} |\,\bm\hat{a}_{ij}\,|^k}$.
\end{tabular}\\[2.5ex]
To estimate the second summand of (\ref{EWD_3_4_07}), we calculate\\[2ex]
\hspace*{2ex}\begin{tabular}[t]{@{}c@{\hspace*{0.8ex}}l@{\hspace*{-11ex}}r@{}}
&$\displaystyle{\Big|\,\sum\limits_{i, j = 1}^{n}\bm\hat{a}_{ij}^3 -   
\sum\limits_{i, j = 1}^{n}\bigl(\,a_{ij}' - d_{ij}\,\bigr)^3\,\Big|}$\\[3.5ex]
$=$&$\displaystyle{\Big|\,\sum\limits_{i, j = 1}^{n}\bigl(\,\bm\hat{a}_{ij}^3 - a_{ij}'^{\,3}\,\bigr) + 
3\,a_{ij}'^{\,2}\,d_{ij} - 3\,a_{ij}'\,d_{ij}^2 + d_{ij}^3\,\Big|}$\\[3ex]
$\leq$&$\displaystyle{\sum\limits_{(i,j)\, \in\, \Gamma}|\bm\hat{a}_{ij}|^3 + 
\dfrac{21}{4}\,\sum\limits_{i, j = 1}^{n}|d_{ij}|}$
&(since $|a_{ij}'| \leq \dfrac{1}{2}$ and $|d_{ij}| \leq \dfrac{3}{2}$, 
cf. Lemma \ref{EWD_3_4_01}, \ref{EWD_3_4_01_BWf}))\\[3ex]
$\leq$&$\displaystyle{\Bigl(\,2^{k-3} + \dfrac{63}{8}\,2^k\,\Bigr)\,\sum\limits_{i, j = 1}^{n}|\bm\hat{a}_{ij}|^k}$
&\hspace*{13ex}(due to Lemma \ref{EWD_3_4_01}, \ref{EWD_3_4_01_BWc}) and \ref{EWD_3_4_01_BWe}))\\[3ex]
$=$&$\displaystyle{2^{k+3}\sum\limits_{i, j = 1}^{n}|\bm\hat{a}_{ij}|^k}$.
\end{tabular}
\item
We obtain the claimed estimate analogously to part \ref{EWD_3_4_06_BWa}) via\\[2ex]
\hspace*{2ex}\begin{tabular}[t]{@{}l@{\hspace*{0.8ex}}c@{\hspace*{0.8ex}}l@{}}
$\displaystyle{\Big|\,1 - \dfrac{1}{\sigma_{\!A'}^4}\,\Big|}$&$\leq$&
$\displaystyle{c_{2}(k)\,\dfrac{1}{n}\sum\limits_{i, j = 1}^{n}|\bm\hat{a}_{ij}|^{k-2}}$
\hspace*{4ex}for $k \geq 4$.
\end{tabular}
\item
Using Lemma \ref{EWD_3_1_18}, \ref{EWD_3_1_18_BWc}) we get analogously to (\ref{EWD_3_4_07})\\[2ex]
\hspace*{2ex}\begin{tabular}[t]{@{}c@{\hspace*{0.8ex}}l@{}}
&$\displaystyle{\bigg|\,\dfrac{1}{n^2}\sum\limits_{i, j, r = 1}^{n}\bm\hat{a}_{ij}^2\,\bm\hat{a}_{ir}^2 -   
\dfrac{1}{n^2}\sum\limits_{i, j, r = 1}^{n}\bm\bar{a}_{ij}^2\,\bm\bar{a}_{ir}^2\,\bigg|}$\\[3.5ex]
$\leq$&$\displaystyle{\Big|\,1 - \dfrac{1}{\sigma_{\!A'}^4}\,\Big|\,\dfrac{\delta_{A}}{n} + 
\dfrac{1}{\sigma_{\!A'}^4}\bigg|\,\dfrac{1}{n^2}\sum\limits_{i, j, r = 1}^{n}\bm\hat{a}_{ij}^2\,\bm\hat{a}_{ir}^2 -   
\dfrac{1}{n^2}\sum\limits_{i, j, r = 1}^{n}
\bigl(\,a_{ij}' - d_{ij}\,\bigr)^2\,\bigl(\,a_{ir}' - d_{ir}\,\bigr)^2\,\bigg|}$.
\end{tabular}\\[3ex]
The second summand of this is estimated as follows. Let $k \geq 2$, then\\[2ex]
\hspace*{2ex}\begin{tabular}[t]{@{}c@{\hspace*{0.8ex}}l@{\hspace*{2.1ex}}r@{}}
&$\displaystyle{\dfrac{1}{n}\,\bigg|\,\sum\limits_{i, j, r = 1}^{n}\bm\hat{a}_{ij}^2\,\bm\hat{a}_{ir}^2 -   
\sum\limits_{i, j, r = 1}^{n} \bigl(\,a_{ij}' - d_{ij}\,\bigr)^2\,\bigl(\,a_{ir}' - d_{ir}\,\bigr)^2\,\bigg|}$\\[3.5ex]
$\leq$&$\displaystyle{\dfrac{1}{n}\,\bigg|\,\sum\limits_{i, j, r = 1}^{n}
\bigl(\,\bm\hat{a}_{ij}^2\,\bm\hat{a}_{ir}^2 - a_{ij}'^{\,2}\,a_{ir}'^{\,2}\,\bigr)\,\bigg|\, +\,
\dfrac{85}{8}\sum\limits_{i,j = 1}^{n}|d_{ij}|}$
&(since $|a_{ij}'| \leq \dfrac{1}{2}$ and $|d_{ij}| \leq \dfrac{3}{2}$)\\[3.5ex]
$\leq$&$\displaystyle{\dfrac{2}{n}\sum\limits_{\begin{array}{@{}c@{}}\\[-5ex]\scriptstyle{i, j, r = 1}\\[-1.5ex]\scriptstyle{(i,j)\, \in\, \Gamma}\end{array}}^{n}\! \bm\hat{a}_{ij}^2\,\bm\hat{a}_{ir}^2\, +\,
\dfrac{85}{8}\sum\limits_{i,j = 1}^{n}|d_{ij}|}$\\[3.5ex]
$\leq$&$\displaystyle{2\sum\limits_{(i,j)\, \in\, \Gamma}^{n}\! \bm\hat{a}_{ij}^2\, +\,
\dfrac{85}{8}\sum\limits_{i,j = 1}^{n}|d_{ij}|}$
&(since $\displaystyle{\sum\limits_{l,r = 1}^{n}\bm\hat{a}_{lr}^2 = n - 1 \leq n}$)\\[3.5ex]
$\leq$&$\displaystyle{\Bigl(\,2^{k-1} + \dfrac{255}{16}\,2^k\,\Bigr)\,\sum\limits_{i, j = 1}^{n}|\bm\hat{a}_{ij}|^k}$
&(due to Lemma \ref{EWD_3_4_01}, \ref{EWD_3_4_01_BWc}) and \ref{EWD_3_4_01_BWe}))\\[3ex]
$=$&$\displaystyle{\dfrac{263}{16}\,2^k\sum\limits_{i, j = 1}^{n}|\bm\hat{a}_{ij}|^k}$.
\end{tabular}
\item
At first, we remark that\\[1.5ex]
\hspace*{2ex}$|\,a^2 - b^2\,| = |\,a + b\,|\,|\,a - b\,| = |\,2\,a - (a - b)\,|\,|\,a - b\,|
\leq \bigl(\,2\,|\,a\,| + |\,a - b\,|\,\bigr)\,|\,a - b\,|$\\[1.5ex]
holds for all $a, b \in \mathbb{R}$. 
We use this inequality and then part \ref{EWD_3_4_06_BWa}) twice, so that we get\\[2.5ex]
\hspace*{2ex}\begin{tabular}[t]{@{}l@{\hspace*{0.8ex}}c@{\hspace*{0.8ex}}l@{}}
$\displaystyle{\bigg|\,\Bigl(\,\dfrac{1}{n}\sum\limits_{i, j = 1}^{n}\bm\hat{a}_{ij}^3\,\Bigr)^2 -   
\Bigl(\,\dfrac{1}{n}\sum\limits_{i, j = 1}^{n}\bm\bar{a}_{ij}^3\,\Bigr)^2\,\bigg|}$
&$\leq$&$\displaystyle{\Bigl(\,2 + C_{1}(3)\,\Bigr)\,\dfrac{\beta_{A}}{n}\,C_{1}(k-1)\,
\dfrac{1}{n}\sum\limits_{i, j = 1}^{n}|\bm\hat{a}_{ij}|^{k-1}}$\\[3ex]
&$\leq$&$\displaystyle{\Bigl(\,2 + C_{1}(3)\,\Bigr)\,C_{1}(k-1)\,\dfrac{1}{n}
\sum\limits_{i, j = 1}^{n}|\bm\hat{a}_{ij}|^{k}}$\hspace*{5.5ex}for $k \geq 4$.
\end{tabular}\\[2.5ex]
For the last inequality, we again used the estimate from (\ref{EWD_3_4_08}).
\hspace*{1ex}\hfill$\Box$
\end{enumerate}
Next, we show a small analytic lemma that we need to derive our most important result 
about the relations between the matrices $\bm\hat{A}$ and $\bm\bar{A}$.\\[4ex]
\refstepcounter{DSBcount}
\label{EWD_3_4_09}
\textbf{\hspace*{-0.8ex}\theDSBcount\ Lemma}\\[0.8ex]
Let $g$ be differentiable and let an $r \in \mathbb{R}$ be given such that\\[2.5ex]
\refstepcounter{DSBcount}
\label{EWD_3_4_10}
\text{\hspace*{-0.8ex}(\theDSBcount)}
\hspace*{4ex}
$\displaystyle{\left\{
\begin{array}{ll@{}}
g(z) = r + g(-z)& \hspace*{2ex}
\text{for all}\ z \in \mathbb{R}\ \ \ \ \ \text{or}\\[1ex]
g(z) = r - g(-z)& \hspace*{2ex}
\text{for all}\ z \in \mathbb{R}\,.
\end{array}  \right.}$\\[2.5ex]
Then
\begin{enumerate}
\item\label{EWD_3_4_09_BWa}
$\sup\limits_{z \in \mathbb{R}} \big|\,g(z+q) - g(z)\,\big| \leq |q|\,||g'||$
\hspace*{4ex}for $q \in \mathbb{R}$,
\item\label{EWD_3_4_09_BWb}
$\sup\limits_{z \in \mathbb{R}} \big|\,g(p\,z) - g(z)\,\big| 
\leq (p-1)\,||id_{\mathbb{R}}\,g'||$
\hspace*{4ex}for $p \geq 1$,
\item\label{EWD_3_4_09_BWc}
$\sup\limits_{z \in \mathbb{R}} \big|\,g(p\,z) - g(z)\,\big| \leq 
\Bigl(\,\dfrac{1}{p} - 1\,\Bigr)\,||id_{\mathbb{R}}\,g'||$
\hspace*{4ex}for $0 < p \leq 1$.
\end{enumerate}
\vspace*{2.5ex}
\textbf{Proof:}\\[0.8ex]
Part \ref{EWD_3_4_09_BWa}) is a trivial conclusion from the mean value theorem\index{mean value theorem} 
and part \ref{EWD_3_4_09_BWc}) 
follows from part \ref{EWD_3_4_09_BWb}) using\\[2ex]  
\hspace*{12.1ex}$\sup\limits_{z \in \mathbb{R}} \big|\,g(p\,z) - g(z)\,\big| =
\sup\limits_{z \in \mathbb{R}} \big|\,g(z) - g(\dfrac{z}{p})\,\big|$.\\[2.5ex]
Part \ref{EWD_3_4_09_BWb}) therefore remains to be verified. Because of (\ref{EWD_3_4_10}), 
$z > 0$ can be assumed without loss of generality. A further application of the 
mean value theorem\index{mean value theorem} then yields\\[2ex]
\hspace*{12.1ex}\begin{tabular}[t]{@{}l@{\hspace*{0.8ex}}c@{\hspace*{0.8ex}}l@{}}
$\big|\,g(p\,z) - g(z)\,\big|$&$=$&$\big|\,g(z + (p-1)\,z) - g(z)\,\big|$\\[2ex]
&$\leq$&$(p-1)\,z\,\sup\limits_{z\, \leq\, t\, \leq\, (p-1)\,z}|\,g'(t)\,|$\\[2ex]
&$\leq$&$(p-1)\,\sup\limits_{z\, \leq\, t\, \leq\, (p-1)\,z}|\,t\,g'(t)\,|$\\[2.5ex]
&$\leq$&$(p-1)\,||id_{\mathbb{R}}\,g'||$.
\end{tabular}\\[-2.4ex]
\hspace*{1ex}\hfill$\Box$\\[4ex]
\refstepcounter{DSBcount}
\label{EWD_3_4_11}
\textbf{\hspace*{-0.8ex}\theDSBcount\ Example}\\[0.8ex]
The condition (\ref{EWD_3_4_10}) is fulfilled, for example, by the function $g(z) = \Phi(z)$ 
and all its derivatives (cf. (\ref{EWD_2_2_01})).\\[3ex]
\refstepcounter{DSBcount}
\label{EWD_3_4_12}
\textbf{\hspace*{-0.8ex}\theDSBcount\ Conclusion}\\[0.8ex]
If a differentiable function $g$ satisfies the condition (\ref{EWD_3_4_10}) 
and $\sigma_{\!A'}^2 \geq \dfrac{2}{3}$ is valid, then\\[2ex]
\hspace*{12.1ex}\begin{tabular}[t]{@{}l@{\hspace*{0.8ex}}c@{\hspace*{0.8ex}}l@{}}
$\Big|\,g\Bigl(\,\dfrac{z - \mu_{A'}}{\sigma_{\!A'}}\,\Bigr) - g(z)\,\Big|$&$\leq$&
$\Big|\,g(z - \mu_{A'}) - g(z)\,\Big| + \Big|\,g\Bigl(\,\dfrac{z - \mu_{A'}}{\sigma_{\!A'}}\,\Bigr) 
- g(z - \mu_{A'})\,\Big|$\\[2ex]
&$\leq$&$|\,\mu_{A'}\,|\,||g'|| + \sqrt{\dfrac{3}{2}}\,\Big|\,\sigma_{\!A'}^2 - 1\,\Big|\,||id_{\mathbb{R}}\,g'||$
\hspace*{8ex}for all $z \in \mathbb{R}$. 
\end{tabular}\\[2ex]
To prove the last inequality, we consider the two cases $\sigma_{\!A'} \leq 1$ and $\sigma_{\!A'} \geq 1$.\\[1ex] 
In the case $\sigma_{\!A'} \leq 1$, Lemma \ref{EWD_3_4_09}, \ref{EWD_3_4_09_BWb}) gives 
the following factor of $||id_{\mathbb{R}}\,g'||$\\[1.5ex]
\hspace*{12.1ex}$\Bigl(\,\dfrac{1}{\sigma_{\!A'}} - 1\,\Bigr) = \dfrac{1}{\sigma_{\!A'}}\,
\Bigl(\,1 - \sigma_{\!A'}\,\Bigr) \leq \sqrt{\dfrac{3}{2}}\,
\Bigl(\,1 - \sigma_{\!A'}\,\Bigr) \leq \sqrt{\dfrac{3}{2}}\,
\Bigl(\,1 - \sigma_{\!A'}^2\,\Bigr)$.\\[2ex]
In the case $\sigma_{\!A'} \geq 1$, on the other hand, Lemma \ref{EWD_3_4_09}, \ref{EWD_3_4_09_BWc})
gives\\[2ex]
\hspace*{12.1ex}$\Bigl(\,\sigma_{\!A'} - 1\,\Bigr) \leq \Bigl(\,\sigma_{\!A'}^2 - 1\,\Bigr)
\leq \sqrt{\dfrac{3}{2}}\,\Bigl(\,\sigma_{\!A'}^2 - 1\,\Bigr)$.
\hspace*{1ex}\hfill$\Box$\\[3.5ex]
By using the Conclusion \ref{EWD_3_4_12} we obtain the following important result.\\[3.5ex]
\refstepcounter{DSBcount}
\label{EWD_3_4_13}
\textbf{\hspace*{-0.8ex}\theDSBcount\ Proposition}\\[0.8ex]
There exist constants $C_{1} > 0$ and $C_{2} > 0$ such that for all matrices $A$ satisfying 
\mbox{\rule[-2.5ex]{0ex}{5ex}$\sigma_{\!A'}^2 \geq \dfrac{2}{3}$} and for all functions 
$h \in \mathcal{H}$ (cf. (\ref{EWD_0_1_01}) for the definition of $\mathcal{H}$) 
holds:\index{function!$\mathcal{H}$}\\[2ex]
\refstepcounter{DSBcount}
\label{EWD_3_4_14}
\text{\hspace*{-0.8ex}(\theDSBcount)}
\hspace*{4ex}
\begin{tabular}[t]{@{}c@{\hspace*{0.8ex}}l@{}}
&$\displaystyle{\bigg|\,E\bigl(h(\mathscr{T}_{A})\bigr) - \int\limits_{\mathbb{R}}\,h(x)\,e_{1,A}'(x)\,dx\,\bigg|}$\\[2.5ex]
$\leq$&$\displaystyle{\bigg|\,E\Bigl(h\bigl(\sigma_{\!A'}\,T_{\bm\bar{A}} + \mu_{A'}\bigr)\Bigr) - 
\int\limits_{\mathbb{R}}\,h\bigl(\sigma_{\!A'}\,x + \mu_{A'}\bigr)\,e_{1,\bm\bar{A}}'(x)\,dx\,\bigg| 
+ C_{1}\,D_{\!A}^2}$,
\end{tabular}\\[2ex]
\refstepcounter{DSBcount}
\label{EWD_3_4_15}
\text{\hspace*{-0.8ex}(\theDSBcount)}
\hspace*{4ex}
\begin{tabular}[t]{@{}c@{\hspace*{0.8ex}}l@{}}
&$\displaystyle{\bigg|\,E\bigl(h(\mathscr{T}_{A})\bigr) - \int\limits_{\mathbb{R}}\,h(x)\,e_{2,A}'(x)\,dx\,\bigg|}$\\[2.5ex]
$\leq$&$\displaystyle{\bigg|\,E\Bigl(h\bigl(\sigma_{\!A'}\,T_{\bm\bar{A}} + \mu_{A'}\bigr)\Bigr) - 
\int\limits_{\mathbb{R}}\,h\bigl(\sigma_{\!A'}\,x + \mu_{A'}\bigr)\,e_{2,\bm\bar{A}}'(x)\,dx\,\bigg| 
+ C_{2}\,E_{\!A}^3}$.
\end{tabular}\\[3ex]
\textbf{Proof:}\\[0.8ex]
To prove (\ref{EWD_3_4_14}) it is sufficient to show\\[2ex]
\hspace*{2ex}$\displaystyle{\bigg|\,E\bigl(h(\mathscr{T}_{A})\bigr) - 
E\Bigl(h\bigl(\sigma_{\!A'}\,T_{\bm\bar{A}} + \mu_{A'}\bigr)\Bigr) - 
\int\limits_{\mathbb{R}}\,h(x)\,e_{1,A}'(x)\,dx +
\int\limits_{\mathbb{R}}\,h\bigl(\sigma_{\!A'}\,x + 
\mu_{A'}\bigr)\,e_{1,\bm\bar{A}}'(x)\,dx\,\bigg|}$\\
\hspace*{1ex}\hfill$\leq C_{1}\,D_{\!A}^2\ $\\[2.5ex]
for the functions $h = 1_{(- \infty,\, z\,]}$, $z \in \mathbb{R}$. Once this has been shown, we can move on to the 
convex combinations\index{function!convex combinations} 
of these functions and then to their limits 
via majorized convergence\index{integration!using majorized convergence}.\\[2.8ex]
Using the notation $H_{2}(x) = x^2 - 1$ (see section \ref{EWD_Kap2_Sec2}), we obtain\\[2.5ex]
\hspace*{3.5ex}\begin{tabular}[t]{@{}c@{\hspace*{0.8ex}}l@{}}
&$\displaystyle{\bigg|\,P\bigl(\,\mathscr{T}_{A} \leq z\,\bigr)
- P\Bigl(\,T_{\bm\bar{A}} \leq \dfrac{z - \mu_{A'}}{\sigma_{\,A'}}\,\Bigr) - e_{1,A}(z) + 
e_{1,\bm\bar{A}}\Bigl(\,\dfrac{z - \mu_{A'}}{\sigma_{\,A'}}\,\Bigr)\,\bigg|}$\\[3ex]
$=$&$\displaystyle{\bigg|\,P\bigl(\,\mathscr{T}_{A} \leq z\,\bigr)
- P\bigl(\,T_{A'} \leq z\,\bigr) - e_{1,A}(z) + 
e_{1,\bm\bar{A}}\Bigl(\,\dfrac{z - \mu_{A'}}{\sigma_{\,A'}}\,\Bigr)\,\bigg|}$
\hspace*{6.1ex}$\Bigl(\,\text{since }  T_{\bm\bar{A}} = \dfrac{T_{A'} - \mu_{A'}}{\sigma_{\,A'}} \,\Bigr)$\\[3ex]
$\leq$&$\displaystyle{P(\mathscr{T}_{A} \not= T_{A'}) + 
\bigg|\,\Phi(z) - \Phi\Bigl(\,\dfrac{z - \mu_{A'}}{\sigma_{\,A'}}\,\Bigr)\,\bigg|}$\\[3ex] 
&$\displaystyle{+\,
\bigg|\,H_{2}(z)\,\psi(z) - H_{2}\Bigl(\,\dfrac{z - 
\mu_{A'}}{\sigma_{\,A'}}\,\Bigr)\,\psi\Bigl(\,\dfrac{z - \mu_{A'}}{\sigma_{\,A'}}\,\Bigr)\,\bigg|\,
\dfrac{\beta_{A}}{6\,n}}$\\[3ex]
&$\displaystyle{+\,\dfrac{1}{6}\,||\,H_{2}\,\psi\,||\,\bigg|\,\dfrac{1}{n}\sum\limits_{i, j = 1}^{n}\bm\hat{a}_{ij}^3 -   
\dfrac{1}{n}\sum\limits_{i, j = 1}^{n}\bm\bar{a}_{ij}^3\,\bigg|}$\\[3.5ex]
$\leq$&$c_{1}\,\dfrac{\delta_{A}}{n} + c_{2}\,\Bigl(\,\dfrac{\beta_{A}}{n}\,\Bigr)^2 \leq C_{1}\,D_{\!A}^2$.
\end{tabular}\\[3ex]
The penultimate inequality is obtained by applying of
\begin{itemize}
\item
Lemma \ref{EWD_3_4_01}, \ref{EWD_3_4_01_BWb}) with $k = 4$, 
and Lemma \ref{EWD_3_4_06}, \ref{EWD_3_4_06_BWa}) with $k = 4$,
\item
Conclusion \ref{EWD_3_4_12} with $g(z) = \Phi(z)$ and $g(z) = H_{2}(z)\,\psi(z)$,
\item
Lemma \ref{EWD_3_4_01}, \ref{EWD_3_4_01_BWg}) and \ref{EWD_3_4_01_BWh})
with $k = 3$ and $k = 4$, and Lemma \ref{EWD_2_2_08}, \ref{EWD_2_2_08_BWa}), \ref{EWD_2_2_08_BWb}), 
\ref{EWD_2_2_08_BWe}), \ref{EWD_2_2_08_BWf}) and \ref{EWD_2_2_08_BWg}). 
\end{itemize}
The last inequality, on the other hand, follows from Lemma \ref{EWD_3_1_18}, \ref{EWD_3_1_18_BWb}).\\[2.8ex]
The proof of (\ref{EWD_3_4_15}) is completely analogous and finally yields the inequalities\\[2ex]
\hspace*{12.1ex}$\ldots \leq c_{3}\,\dfrac{\eta_{A}}{n} + 
c_{4}\,\dfrac{\delta_{A}}{n}\,\dfrac{\beta_{A}}{n} +
c_{5}\,\Bigl(\,\dfrac{\beta_{A}}{n}\,\Bigr)^3 
\leq C_{2}\,E_{\!A}^3$.\\[3ex]
The penultimate inequality in the estimate for (\ref{EWD_3_4_15}) is obtained by applying of
\begin{itemize}
\item
Lemma \ref{EWD_3_4_01}, \ref{EWD_3_4_01_BWb}) with $k = 5$, 
and Lemma \ref{EWD_3_4_06}, \ref{EWD_3_4_06_BWa}), \ref{EWD_3_4_06_BWb}), \ref{EWD_3_4_06_BWc}) 
and \ref{EWD_3_4_06_BWd}) with $k = 5$,
\item
Conclusion \ref{EWD_3_4_12} with $g(z) = \Phi(z)$, $g(z) = H_{2}(z)\,\psi(z)$, 
$g(z) = H_{3}(z)\,\psi(z)$ and\\ 
$g(z) = H_{5}(z)\,\psi(z)$ 
(cf. section \ref{EWD_Kap2_Sec2} for the definition of the Hermite polynomials $H_{n}$),
\item
Lemma \ref{EWD_3_4_01}, \ref{EWD_3_4_01_BWg}) and \ref{EWD_3_4_01_BWh})
with $k = 3$, $k = 4$ and $k = 5$,
\item
Lemma \ref{EWD_3_1_18}, \ref{EWD_3_1_18_BWd}), and
Lemma \ref{EWD_2_2_08}, \ref{EWD_2_2_08_BWa}), \ref{EWD_2_2_08_BWb}), 
\ref{EWD_2_2_08_BWe}), \ref{EWD_2_2_08_BWf}), \ref{EWD_2_2_08_BWg}),
\ref{EWD_2_2_08_BWh}), \ref{EWD_2_2_08_BWi}),
\ref{EWD_2_2_08_BWj}), \ref{EWD_2_2_08_BWl}) and \ref{EWD_2_2_08_BWm}).
\end{itemize}
The last inequality, on the other hand, follows again from Lemma \ref{EWD_3_1_18}, \ref{EWD_3_1_18_BWb}).
\hspace*{1ex}\hfill$\Box$\\[4ex]
The last part of this section deals with the conditions of the Theorems \ref{EWD_3_1_10} and \ref{EWD_3_1_13}. 
The aim of the following two propositions is to transfer the inequalities of (\ref{EWD_3_1_11}) and (\ref{EWD_3_1_15}) 
\linebreak
from the matrices $B \in N(l,\bm\hat{A})$ to the 
corresponding {''truncated''} matrices $\BTB \in N(l,\bm\bar{A})$.\\[2.8ex]
For the remainder of this section, let $B \in N(l,\bm\hat{A})$ be fixed
and $Z_{l} = \{\,z_{1},\ldots,z_{l}\,\}$ ($S_{l} =\{\,s_{1},\ldots,s_{l}\,\}$) 
the set of indices of the cancelled rows (columns).\index{matrix!submatrix!$Z_{l}$, $S_{l}$}
Furthermore,\\[2ex]
\hspace*{12.2ex}\begin{tabular}[t]{@{}l@{\hspace*{0.8ex}}l@{}}
$\BTB$&is the $(n-l){\times}(n-l)-$matrix, which is obtained from $\bm\bar{A}$ by cancelling the\\[1ex]
&the same rows and columns as for $B$,\\[2ex]
$B'$&is the $(n-l){\times}(n-l)-$matrix, which is obtained from $A'$ by cancelling the\\[1ex]
&the same rows and columns as for $B$ (also $B'$ = truncated matrix $B$), and\\[2ex]
$\theta_{B'}$&$\displaystyle{=\,\sum\limits_{\textstyle{i \not\in Z_{l}}} a_{i\boldsymbol{.}}' +
\sum\limits_{\textstyle{j \not\in S_{l}}} a_{\boldsymbol{.}j}' - (n - l)\,a_{\boldsymbol{.}\boldsymbol{.}}'}$ .
\end{tabular}
\index{matrix!submatrix}\index{matrix!submatrix!$\BTB$, $B'$, $\theta_{B'}$}\\[2.8ex]
It follows that\\[2ex]
\refstepcounter{DSBcount}
\label{EWD_3_4_16}
\text{\hspace*{-0.8ex}(\theDSBcount)}
\hspace*{3.8ex}$T_{\SBTB} = \dfrac{T_{B'} - \theta_{B'}}{\sigma_{\!A'}}$.\\[3ex]
The next proposition provides a relation between the 
\textit{r}th differences\index{difference!\textit{r}th} of $F_{B}$ and $F_{\SBTB}$.\\[4ex]
\refstepcounter{DSBcount}
\label{EWD_3_4_17}
\textbf{\hspace*{-0.8ex}\theDSBcount\ Proposition}\\[0.8ex]
Let $\dfrac{n}{n - l} \leq 10$, 
$\sigma_{\!A'} > 0$ and $k, r \in \mathbb{N}$.\\[2ex]
Then there exists a constant $C > 0$, which depends only on $k$ and $r$, such that\\[2ex]
\hspace*{12.1ex}$\displaystyle{\big|\big|\,\Delta_{y}^{r}F_{\SBTB}\,\big|\big| \leq
\big|\big|\,\Delta_{\sigma_{\!A'}y}^{r}F_{B}\,\big|\big|
+ \dfrac{C}{n}\sum\limits_{i, j = 1}^{n}|\bm\hat{a}_{ij}|^k}$
\hspace*{4ex}for all $y \in \mathbb{R}$.\\[4ex]
\textbf{Proof:}\\[0.8ex]
Due to Lemma \ref{EWD_3_4_01}, \ref{EWD_3_4_01_BWb}) ($B'$ = truncated matrix $B$) 
we have\\[2ex]
\refstepcounter{DSBcount}
\label{EWD_3_4_18}
\text{\hspace*{-0.8ex}(\theDSBcount)}
\hspace*{3.8ex}$\displaystyle{P(\,T_{B} \not= T_{B'}\,) \leq \dfrac{|\Gamma|}{n-l}
\leq \dfrac{10}{n}\,2^k\,\sum\limits_{i, j = 1}^{n}|\bm\hat{a}_{ij}|^k}$.\\[2.5ex]
Thus, for all $z \in \mathbb{R}$,\\[2.5ex]
\hspace*{12.1ex}\begin{tabular}[t]{@{}l@{\hspace*{0.8ex}}c@{\hspace*{0.8ex}}l@{}}
$\big|\,\Delta_{y}^{r}F_{\SBTB}(z)\,\big|$&$=$&
$\big|\,\Delta_{\sigma_{\!A'}y}^{r}F_{B'}(\sigma_{\!A'}\,z + \theta_{B'})\,\big|$
\hspace*{27.9ex}(due to (\ref{EWD_3_4_16}))\\[3ex]
&$\leq$&$\big|\,\Delta_{\sigma_{\!A'}y}^{r}F_{B}(\sigma_{\!A'}\,z + \theta_{B'})\,\big| +
2^{r}\,P(\,T_{B} \not= T_{B'}\,)$\\[3ex]
&$\leq$&$\displaystyle{\big|\big|\,\Delta_{\sigma_{\!A'}y}^{r}F_{B}\,\big|\big| +
\dfrac{C}{n}\sum\limits_{i, j = 1}^{n}|\bm\hat{a}_{ij}|^k}$,\hspace*{7ex}where $C = 10 \cdot 2^{k+r}$.
\end{tabular}\\[-4.8ex] 
\hspace*{1ex}\hfill$\Box$\\[6ex]
Finally, with more effort, we prove a similar relation between 
$z\,\Delta_{y}^{2}F_{\SBTB}(z)$ and 
$z\,\Delta_{\sigma_{\!A'}y}^{2}F_{B}(z)$.\\[4ex]
\refstepcounter{DSBcount}
\label{EWD_3_4_19}
\textbf{\hspace*{-0.8ex}\theDSBcount\ Proposition}\\[0.8ex]
Let $\dfrac{n}{n - l} \leq 10$, 
$|\mu_{B}| \leq 1$, $\big|\,\sigma_{\!B}^2 - 1\,\big| \leq \dfrac{1}{3}$
and $\big|\,\sigma_{\!A'}^2 - 1\,\big| \leq \dfrac{1}{3}$.\\[2ex]
Then there exists a constant $C > 0$ (which does not depend on anything!) such that\\[2ex]
\hspace*{12.1ex}$\displaystyle{\big|\big|\,z\,\Delta_{y}^{2}F_{\SBTB}(z)\,\big|\big| \leq
2\,\big|\big|\,z\,\Delta_{\sigma_{\!A'}y}^{2}F_{B}(z)\,\big|\big|
+ C\,\bigl(\,E_{\!A}^2 + y^2\,\bigr)}$
\hspace*{4ex}for all $y \geq 0$.\\[4ex]
\textbf{Proof:}\\[0.8ex]
First, we show the existence of a constant $C_{1}$ such that\\[2.5ex]
\refstepcounter{DSBcount}
\label{EWD_3_4_20}
\text{\hspace*{-0.8ex}(\theDSBcount)}
\hspace*{3.8ex}$\displaystyle{\Big|\,E\bigl(\,T_{B'}\,1_{\{\,T_{B'}\, \leq\, z\,\}}\,\bigr) -
E\bigl(\,T_{B}\,1_{\{\,T_{B}\, \leq\, z\,\}}\,\bigr)\,\Big| \leq C_{1}\,E_{\!A}^2}$
\hspace*{4ex}for all $z \in \mathbb{R}$.\\[2.5ex]
To prove this, we use the triangle inequality and obtain\\[2.5ex]
\refstepcounter{DSBcount}
\label{EWD_3_4_21}
\text{\hspace*{-0.8ex}(\theDSBcount)}
\hspace*{1.2ex}\begin{tabular}[t]{@{}c@{\hspace*{0.8ex}}l@{}}
&$\displaystyle{\Big|\,E\bigl(\,T_{B'}\,1_{\{\,T_{B'}\, \leq\, z\,\}}\,\bigr) -
E\bigl(\,T_{B}\,1_{\{\,T_{B}\, \leq\, z\,\}}\,\bigr)\,\Big|}$\\[3ex]
$\leq$&$\displaystyle{\Big|\,E\bigl(\,T_{B'}\,1_{\{\,T_{B'}\, \leq\, z\,\}}\,\bigr) -
E\bigl(\,T_{B}\,1_{\{\,T_{B'}\, \leq\, z\,\}}\,\bigr)\,\Big| + 
\Big|\,E\bigl(\,T_{B}\,1_{\{\,T_{B'}\, \leq\, z\,\}}\,\bigr) -
E\bigl(\,T_{B}\,1_{\{\,T_{B}\, \leq\, z\,\}}\,\bigr)\,\Big|}$\\[3ex]
$\leq$&$E\bigl(\,|\,T_{B'} - T_{B}\,|\,\bigr) + E\bigl(\,|\,T_{B}\,|\,1_{\{\,T_{B'}\, \not=\, T_{B}\,\}}\,\bigr)$\\[3ex]
$\leq$&$E\bigl(\,|\,T_{B'} - T_{B}\,|\,\bigr) + 
\sigma_{\!B}\,E\bigl(\,|\,\mathscr{T}_{B}\,|\,1_{\{\,T_{B'}\, \not=\, T_{B}\,\}}\,\bigr) +
|\,\mu_{B}\,|\,P\bigl(\,T_{B'} \not= T_{B}\,\bigr)$.
\end{tabular}\\[3.5ex]
The first summand of (\ref{EWD_3_4_21}) can be estimated as follows:\\[2ex]
\hspace*{10ex}\begin{tabular}[t]{@{}l@{\hspace*{0.8ex}}c@{\hspace*{0.8ex}}l@{\hspace*{9.3ex}}r@{}}
$\displaystyle{E\bigl(\,|\,T_{B'} - T_{B}\,|\,\bigr)}$&$\leq$&
$\displaystyle{\sum\limits_{\textstyle{i \not\in Z_{l}}} E\bigl(\,\Big|\,b_{i\pi(i)}' - b_{i\pi(i)}\,\Big|\,\bigr)}$\\[4ex]
&$=$&$\displaystyle{\dfrac{1}{n-l}\,\sum\limits_{\begin{array}{@{}c@{}}\\[-5ex]
\textstyle{(i \not\in Z_{l})\,\wedge\,(j \not\in S_{l})}\\[-0.5ex]\textstyle{(i,j)\, \in\, \Gamma}\end{array}}|b_{ij}|}$&
(since $P(\pi(i) = j) = \dfrac{1}{n-l}$)
\end{tabular}\\[3ex]
\hspace*{10ex}\begin{tabular}[t]{@{}l@{\hspace*{0.8ex}}c@{\hspace*{0.8ex}}l@{\hspace*{24.9ex}}r@{}}
\hspace*{16.05ex}&$\leq$&$\displaystyle{\dfrac{10}{n}\,\sum\limits_{\textstyle{(i,j)\, \in\, \Gamma}}|\bm\hat{a}_{ij}|}$&
(since $\dfrac{n}{n - l} \leq 10$)\\[4.5ex]
&$\leq$&$c_{1}\,\dfrac{\delta_{A}}{n}$&
(cf. Lemma \ref{EWD_3_4_01}, \ref{EWD_3_4_01_BWc}))\\[3ex]
&$\leq$&$c_{1}\,E_{\!A}^2$&
(cf. Lemma \ref{EWD_3_1_18}, \ref{EWD_3_1_18_BWb})).
\end{tabular}\\[3ex]
Furthermore, we estimate\\[2ex]
\hspace*{13.4ex}\begin{tabular}[t]{@{}l@{\hspace*{0.8ex}}c@{\hspace*{0.8ex}}l@{\hspace*{14.7ex}}r@{}}
$E\bigl(\,|\mathscr{T}_{B}|^3\,\bigr)^{1/3}$&$\leq$&$E\bigl(\,|\mathscr{T}_{B}|^4\,\bigr)^{1/4}$&
(H{\"o}lder's inequality\index{H{\"o}lder's inequality!for random variables})\\[3ex]
&$\leq$&$\Bigl(\,3 + 97\,\dfrac{\delta_{B}}{n-l}\,\Bigr)^{1/4}$&
(cf. Proposition \ref{EWD_3_2_02}, \ref{EWD_3_2_02_BWd}))\\[3ex]
&$\leq$&$\Bigl(\,3 + 97\,\dfrac{\delta_{B}}{n-l}\,\Bigr)^{1/3}$&
(since $\Bigl(\,3 + 97\,\dfrac{\delta_{B}}{n-l}\,\Bigr) \geq 1$)\\[3ex]
&$\leq$&$c_{2} + c_{3}\,\Bigl(\,\dfrac{\delta_{A}}{n}\,\Bigr)^{1/3}$&
(cf. Lemma \ref{EWD_3_3_10} and $\sigma_{\!B}^2 \geq \dfrac{2}{3}$).
\end{tabular}\\[3ex]
We note that for the last inequality we used the formula\\[2ex]
\hspace*{14.8ex}$(a + b\,x)^{1/3} \leq a^{1/3} + b^{1/3}\,x^{1/3}$
\hspace*{4ex}for all $a, b, x \geq 0$.\\[2.5ex]
For the second summand of (\ref{EWD_3_4_21}) we thus obtain:\\[2ex]
\hspace*{2ex}\begin{tabular}[t]{@{}l@{\hspace*{0.8ex}}c@{\hspace*{0.8ex}}l@{\hspace*{-27.1ex}}r@{}}
$\sigma_{\!B}\,E\bigl(\,|\,\mathscr{T}_{B}\,|\,1_{\{\,T_{B'}\, \not=\, T_{B}\,\}}\,\bigr)$&$\leq$&
$\sigma_{\!B}\,E\bigl(\,|\mathscr{T}_{B}|^3\,\bigr)^{1/3}\,P(\,T_{B} \not= T_{B'}\,)^{2/3}$&
(H{\"o}lder's ineq.\index{H{\"o}lder's inequality!for random variables}, $p = 3$, $q = \dfrac{3}{2}$)\\[3ex]
&$\leq$&$\sqrt{\dfrac{4}{3}}\,\biggl(\,c_{2} + c_{3}\,\Bigl(\,\dfrac{\delta_{A}}{n}\,\Bigr)^{1/3}\,\biggr)\,
P(\,T_{B} \not= T_{B'}\,)^{2/3}$&
(since $\sigma_{\!B}^2 \leq \dfrac{4}{3}$)\\[3ex]
&$=$&$\sqrt{\dfrac{4}{3}}\,c_{2}\,P(\,T_{B} \not= T_{B'}\,)^{2/3} + 
\sqrt{\dfrac{4}{3}}\,c_{3}\,\Bigl(\,\dfrac{\delta_{A}}{n}\,\Bigr)^{1/3}\,
P(\,T_{B} \not= T_{B'}\,)^{2/3}$\\[3ex]
&$\leq$&$c_{4}\,\Bigl(\,\dfrac{\eta_{A}}{n}\,\Bigr)^{2/3} + c_{5}\,\dfrac{\delta_{A}}{n}$&
(cf. (\ref{EWD_3_4_18}) for $k = 5$ and $k =4$)\\[3ex]
&$\leq$&$(c_{4} + c_{5})\,E_{\!A}^2$&
(cf. Lemma \ref{EWD_3_1_18}, \ref{EWD_3_1_18_BWb})).
\end{tabular}\\[3ex]
Finally, the estimate of the third summand of (\ref{EWD_3_4_21}) is\\[2.1ex]
\hspace*{6.5ex}\begin{tabular}[t]{@{}l@{\hspace*{0.8ex}}c@{\hspace*{0.8ex}}l@{\hspace*{28.8ex}}r@{}}
$|\,\mu_{B}\,|\,P\bigl(\,T_{B'} \not= T_{B}\,\bigr)$&$\leq$&$P\bigl(\,T_{B'} \not= T_{B}\,\bigr)$
&(since $|\mu_{B}| \leq 1$)\\[3ex]
&$\leq$&$c_{6}\,\dfrac{\delta_{A}}{n} \leq c_{6}\,E_{\!A}^2$
&(cf. (\ref{EWD_3_4_18}) for $k = 4$).
\end{tabular}\\[3ex]
As a consequence of (\ref{EWD_3_4_20}), we now get a constant $C_{2}$ such that\\[2.2ex]
\refstepcounter{DSBcount}
\label{EWD_3_4_22}
\text{\hspace*{-0.8ex}(\theDSBcount)}
\hspace*{3.8ex}\begin{tabular}{@{}l@{}}
$\Big|\,z\,P\bigl(\,T_{B'} \in I\,\bigr) - z\,P\bigl(\,T_{B} \in I\,\bigr)\,\Big|
\leq C_{2}\,\bigl(\,E_{\!A}^2 + a^2 + b^2\,\bigr)$\\[2ex]
for all $I = (\,z + a,\,z+a+b\,]$, where $z \in \mathbb{R}$ and $a, b \geq 0$.
\end{tabular}\\[2.5ex]
We proof this as follows:\\[2.1ex]
\hspace*{7.5ex}\begin{tabular}[t]{@{}c@{\hspace*{0.8ex}}l@{\hspace*{-42.3ex}}r@{}}
&$\Big|\,z\,P\bigl(\,T_{B'} \in I\,\bigr) - z\,P\bigl(\,T_{B} \in I\,\bigr)\,\Big|$\\[3ex]
$\leq$&$\Big|\,E\bigl(\,T_{B'}\,1_{\{\,T_{B'}\, \in\, I\,\}}\,\bigr) -
E\bigl(\,T_{B}\,1_{\{\,T_{B}\, \in\, I\,\}}\,\bigr)\,\Big| +
\bigl(\,a + b\,\bigr)\,\Bigl(\,P\bigl(\,T_{B'} \in I\,\bigr) + P\bigl(\,T_{B} \in I\,\bigr)\,\Bigr)$\\[3ex]
$\leq$&$2\,C_{1}\,E_{\!A}^2 + \bigl(\,a + b\,\bigr)\,\Bigl(\,2\,P\bigl(\,T_{B} \in I\,\bigr) 
+ P(\,T_{B} \not= T_{B'}\,)\,\Bigr)$&
(due to (\ref{EWD_3_4_20}))\\[3ex]
$\leq$&$2\,C_{1}\,E_{\!A}^2 + \bigl(\,a + b\,\bigr)\,\biggl(\,4\,\mathcal{K}_{1}\,\dfrac{\beta_{B}}{n-l}
+ \dfrac{2}{\sqrt{2\,\pi}}\,\dfrac{b}{\sigma_{\!B}} + c_{7}\,\dfrac{\beta_{A}}{n}\,\biggr)$\\[1.5ex]
&&(cf. Theorem \ref{EWD_3_1_03} for $\mathscr{T}_{B}$ and (\ref{EWD_3_4_18}))\\[2ex]
$\leq$&$C_{2}\,\bigl(\,E_{\!A}^2 + a^2 + b^2\,\bigr)$&
(cf. Lemma \ref{EWD_3_3_10}, $\sigma_{\!B}^2 \geq \dfrac{2}{3}$ and $\dfrac{\beta_{A}}{n} \leq E_{\!A}$).
\end{tabular}\\[3ex]
Since we additionally have\\[2.1ex]
\refstepcounter{DSBcount}
\label{EWD_3_4_23}
\text{\hspace*{-0.8ex}(\theDSBcount)}
\hspace*{3.8ex}$\displaystyle{|\theta_{B'}| \leq c_{8}\,\dfrac{\delta_{A}}{n} \leq c_{8}\,E_{\!A}^2}$,
\hfill(cf. Lemma \ref{EWD_3_4_01}, \ref{EWD_3_4_01_BWd})),\hspace*{0.1ex}\\[2.5ex]
we obtain for all $z \in \mathbb{R}$ and $y \in \mathbb{R}$ (with the
abbreviation $z' = \sigma_{\!A'}\,z + \theta_{B'}$):\\[2.2ex]
\hspace*{7.5ex}\begin{tabular}[t]{@{}l@{\hspace*{0.8ex}}c@{\hspace*{0.8ex}}l@{\hspace*{-20.5ex}}r@{}}
$\big|\,z\,\Delta_{y}^{2}F_{\SBTB}(z)\,\big|$&$\leq$&
$\dfrac{1}{\sigma_{\!A'}}\,\Bigl(\,\big|\,z'\,\Delta_{\sigma_{\!A'}y}^{2}F_{B'}(z')\,\big| 
+ |\theta_{B'}|\,\Bigr)$&
(due to (\ref{EWD_3_4_16}) and $\big|\,\Delta_{y}^{2}F_{\SBTB}(z)\,\big| \leq 1$)\\[3ex]
&$\leq$&$\dfrac{1}{\sigma_{\!A'}}\,\Bigl(\,\big|\,z'\,\Delta_{\sigma_{\!A'}y}^{2}F_{B}(z')\,\big| 
+ C_{2}\,\bigl(\,2\,E_{\!A}^2 + 3\,\sigma_{\!A'}^2\,y^2\,\bigr) + c_{8}\,E_{\!A}^2\,\Bigr)$\\[3ex]
&$\leq$&$2\,\big|\big|\,z\,\Delta_{\sigma_{\!A'}y}^{2}F_{B}(z)\,\big|\big|
+ C\,\bigl(\,E_{\!A}^2 + y^2\,\bigr)$.
\end{tabular}\\[3ex]
For the penultimate inequality, (\ref{EWD_3_4_22}) was applied twice, 
once with $a = 0$, $b = \sigma_{\!A'}y\ (\geq 0)$ and 
once with $a = b = \sigma_{\!A'}y\ (\geq 0)$. Furthermore, (\ref{EWD_3_4_23}) was used.\\[2ex]
For the last inequality, we used 
the condition $\dfrac{2}{3} \leq \sigma_{\!A'}^2 \leq \dfrac{4}{3}$.\hspace*{1ex}\hfill$\Box$\\[5ex]
\refstepcounter{DSBcount}
\label{EWD_3_4_24}
\textbf{\hspace*{-0.8ex}\theDSBcount\ Remark}\\[0.8ex]
Note, that the condition \mbox{\rule[-2.8ex]{0ex}{4ex}$\dfrac{n}{n - l} \leq 10$} is fulfilled in our relevant cases
$l \leq 8 < 10 \leq n$ (cf. section \ref{EWD_Kap3_Sec7}) and $l \leq 16 < 18 \leq n$
(cf. section \ref{EWD_Kap3_Sec8}).                              
\section{Construction of several random permutations}\label{EWD_Kap3_Sec5}

In order to apply Stein's method, as we used it in chapter \ref{EWD_Kap1}, 
to linear rank statistics, we must, among other things, find an analogue to the independence of 
$S_{n-1}^{n} = (\,1/\sqrt{n}\,)\,(\,X_{1} + \ldots + X_{n-1}\,)$ and $X_{n}$ 
(see section \ref{EWD_Kap1_Sec3}).\\[2.8ex] 
This is done in the following by generating five specially constructed 
random permutations\index{random permutation} $\pi_{1}$, $\pi_{2}$, $\pi_{3}$, $\pi_{4}$, $\pi_{5}$
(see (\ref{EWD_3_5_12}) and (\ref{EWD_3_5_13}) together with Table \ref{EWDTab_3_5_01})
and subsequent by defining the differences $\Delta T_{i} = T_{i+1} - T_{i}$, where the $T_{i}$ are the rank statistics 
belonging to these random permutations $\pi_{i}$ (see the definitions in (\ref{EWD_3_5_21})). 
The role of $S_{n-1}^{n}$ is then played by the $T_{i}$ and that of $X_{n}$ by the 
$\Delta T_{j}$, $j \geq i+1$ (cf. Lemma \ref{EWD_3_5_23}).\\[2.8ex]
The existing independencies between the $T_{i}$ and $\Delta T_{j}$ can be 
explained and motivated to some extent using the important Table \ref{EWDTab_3_5_02}.
This table shows the random indices\index{random index} that are obtained when the random permutations 
$\pi_{1},\ldots,\pi_{5}$ are applied to the random indices $I_{1},\ldots,I_{16}$.
It is essential to note that the $I_{1},\ldots,I_{16}$
are \textbf{in\-de\-pen\-dent} of the $\pi_{1},\ldots,\pi_{5}$ (cf. (\ref{EWD_3_5_12})
and Lemma \ref{EWD_3_5_15}, \ref{EWD_3_5_15_BWa})).
Under $\pi_{1}$ the obtained random indices are completely given by 
the $J_{1},\ldots,J_{16}$ (cf. (\ref{EWD_3_5_14})).
Under $\pi_{2},\ldots,\pi_{5}$ some of these obtained random indices are no longer specified explicitly, as they are 
not used explicitly in these cases and would also be very complex.\\[2.8ex]
In the following, we look at the Table \ref{EWDTab_3_5_02} from \textbf{right to left} and consider the measurability statements of Lemma \ref{EWD_3_5_22}, i.e. 
$\Delta T_{j} \in \sigma(I_{1},\ldots,J_{1},\ldots)$\index{random permutation!$\sigma(X)$, $f \in \sigma(X)$}
(briefly, $\Delta T_{j} \hookrightarrow 
(I_{1},\ldots,J_{1},\ldots)$).\index{random permutation!$\Delta T_{j} \hookrightarrow (I_{1},\ldots,J_{1},\ldots)$}
First of all, according to the above remarks,
\begin{itemize}
\item
$\pi_{5}$ and $I_{1}$ are independent, but \textbf{not} $\pi_{5}$ and $(\,I_{1},\,J_{1}\,)$.
\end{itemize}
Therefore, we swap $J_{1}$ and $J_{2}$ in the first two lines of the Table \ref{EWDTab_3_5_02} when we move 
from $\pi_{5}$ to $\pi_{4}$. We get:
\begin{itemize}
\item
$T_{4} \hookrightarrow \pi_{4}$ 
and $(\,I_{1},\,J_{1}\,) = (\,I_{1},\,\pi_{4}(\,I_{2}\,)\,) \hookleftarrow a_{I_{1}J_{1}}$ are independent.
\end{itemize}
For $\pi_{3}$ the blocks $(\,J_{1},\,J_{2}\,)$ and $(\,J_{3},\,J_{4}\,)$ are swapped in 
Table \ref{EWDTab_3_5_02}, so that
\begin{itemize}
\item
$T_{3} \hookrightarrow \pi_{3}$ and $(\,I_{1}, I_{2}, J_{1}, J_{2}\,) =
(\,I_{1}, I_{2},\,\pi_{3}(\,I_{3}\,), \,\pi_{3}(\,I_{4}\,)\,) \hookleftarrow \Delta T_{4}$ are independent.
\end{itemize}
The blocks that are swapped now become larger and larger. To be more precisely, the block size doubles 
with each additional step. For $\pi_{2}$ the blocks
$(\,J_{1},\ldots,\,J_{4}\,)$ and $(\,J_{5},\ldots,\,J_{8}\,)$ are 
swapped in Table \ref{EWDTab_3_5_02}, so that
\vspace*{1.2ex}

\pagebreak

\begin{itemize}
\item
$T_{2} \hookrightarrow \pi_{2}$ and $(\,I_{1},.\,.\,,\, I_{4}, J_{1},.\,.\,,\,J_{4}\,) =
(\,I_{1},.\,.\,,\, I_{4},\,\pi_{2}(\,I_{5}\,),.\,.\,,\,\pi_{2}(\,I_{8}\,)\,) \hookleftarrow \Delta T_{3}$\\ 
are independent.
\end{itemize}
Finally, for $\pi_{1}$ the blocks $(\,J_{1},\ldots,\,J_{8}\,)$ and
$(\,J_{9},\ldots,\,J_{16}\,)$ are swapped in Table \ref{EWDTab_3_5_02}, so that
\begin{itemize}
\item
$T_{1} \hookrightarrow \pi_{1}$ and $(\,I_{1},.\,.\,,\, I_{8}, J_{1},.\,.\,,\,J_{8}\,) =
(\,I_{1},.\,.\,,\, I_{8},\,\pi_{1}(\,I_{9}\,),.\,.\,,\,\pi_{1}(\,I_{16}\,)\,) \hookleftarrow \Delta T_{2}$\\ 
are independent.
\end{itemize}\vspace*{1.5ex}
After these considerations, we start with the formal, explicit construction 
(in reverse order, i.e. $\pi_{1}$ is defined first, then $\pi_{2}$ \ldots and finally $\pi_{5}$).\\[2.8ex]
Let
\index{random permutation!$N$, $M_{16}$}
\index{random permutation!$\underline{i}$}
\index{random permutation!$i_{1},\ldots,i_{16}$}\\[2ex]
\refstepcounter{DSBcount}
\label{EWD_3_5_01}
\text{\hspace*{-0.8ex}(\theDSBcount)} 
\hspace*{4ex}
$\begin{array}{@{}l@{\hspace*{0.8ex}}c@{\hspace*{0.8ex}}
l@{\hspace*{0.6ex}}l@{\hspace*{3ex}}l@{\hspace*{1.6ex}}
l@{\hspace*{1.6ex}}l@{\hspace*{3ex}}l@{}}
N&=&\multicolumn{5}{@{\hspace*{0.2ex}}l@{}}{\bigl\{\,1,\ldots,n\,\bigr\}\ \ \text{and}}\\[1ex]
M_{16}&=&\biggl\{&\multicolumn{5}{@{\hspace*{0.2ex}}l@{}}
{\underline{i} = (\,i_{1},\ldots,i_{16}\,)\,\in\,N^{16}\,:\,\text{$\underline{i}$ satisfies 
the equivalences:}}\\[1ex]
&&&\text{(t1)}&i_{1} = i_{2}&\Leftrightarrow&i_{3} = i_{4},\index{random permutation!properties!(t1)}\\[1.2ex]
&&&\text{(u1)}&i_{1} = i_{2}&\Leftrightarrow&i_{7} = i_{8},\\[1.2ex]
&&&\text{(u2)}&i_{3} = i_{4}&\Leftrightarrow&i_{5} = i_{6},\\[1.2ex]
&&&\text{(u3)}&i_{l} = i_{k}&\Leftrightarrow&i_{l+6} = 
i_{k+2}&\text{for}\ 1 \leq l \leq 2,\ \ 3 \leq k \leq 4,\index{random permutation!properties!(u1), (u2), (u3)}\\[1.2ex]
&&&\text{(v1)}&i_{l} = i_{k}&\Leftrightarrow&i_{l+12} = i_{k+12}&\text{for}\ 1 \leq l \leq 4,\ \ 1 \leq k \leq 4,\\[1.2ex]
&&&\text{(v2)}&i_{l} = i_{k}&\Leftrightarrow&i_{l+4} = i_{k+4}&\text{for}\ 5 \leq l \leq 8,\ \ 5 \leq k \leq 8,\\[1ex]
&&&\text{(v3)}&i_{l} = i_{k}&\Leftrightarrow&i_{l+12} = i_{k+4}&\text{for}\ 1 \leq l \leq 4,\ \ 5 \leq k \leq 8
\ \,\biggr\}.\index{random permutation!properties!(v1), (v2), (v3)}
\end{array}$\\[2.8ex]
For each $\underline{i} = (\,i_{1},\ldots,i_{16}\,) \in M_{16}$ 
we \textbf{\underline{fix}} once and for all permutations 
$v(\underline{i})$, $u(\underline{i})$, $t(\underline{i})$ and $s(\underline{i})$ 
\linebreak
of $N$ 
with properties as described in Table \ref{EWDTab_3_5_01} on the following page.\\[3.5ex]
\refstepcounter{DSBcount}
\label{EWD_3_5_02}
\textbf{\hspace*{-0.8ex}\theDSBcount\ Remarks on the Table \ref{EWDTab_3_5_01}}
\begin{enumerate}
\item\label{EWD_3_5_02_BWa}
The mapping\\[2ex]
\hspace*{6.3ex}\begin{tabular}[t]{@{}l@{\hspace*{3ex}}l@{}}
$M_{16} \ni \underline{i} \longrightarrow v(\underline{i})$&is well defined beause 
of (v1), (v2) and (v3); and\index{random permutation!properties!(v1), (v2), (v3)}\\[1ex]
$M_{16} \ni \underline{i} \longrightarrow u(\underline{i})$&is well defined beause 
of (u1), (u2) and (u3); and\index{random permutation!properties!(u1), (u2), (u3)}\\[1ex]
$M_{16} \ni \underline{i} \longrightarrow t(\underline{i})$&is well defined beause 
of (t1).\index{random permutation!properties!(t1)}
\end{tabular}\vspace*{1.3ex}
\item\label{EWD_3_5_02_BWb}
As can be seen from the Table \ref{EWDTab_3_5_01}, we have not defined all values of 
$v(\underline{i})$, $u(\underline{i})$ and $t(\underline{i})$ \textbf{explicitly}.
The reason for this is that, on the one hand, we do not need these explicit 
definitions and, 
on the other hand, we would have to distinguish many cases for a complete explicit definition. 
The following two cases are mentioned as an example:\\[2ex] 
\hspace*{6.3ex}\begin{tabular}[t]{@{}l@{\hspace*{3ex}}c@{\hspace*{3ex}}l@{}}
$i_{8} = i_{9}$&$\Rightarrow$&$[v(\underline{i})](i_{9}) = [v(\underline{i})](i_{8}) = i_{12}$,\hspace*{2ex}but\\[1ex]
$i_{8} \not= i_{9}$&$\Rightarrow$&$[v(\underline{i})](i_{9}) \not= [v(\underline{i})](i_{8}) = i_{12}$.
\end{tabular}
\end{enumerate}
\vspace*{2ex}
\refstepcounter{Tabcount}
\label{EWDTab_3_5_01}
\hspace*{12.1ex}\textbf{Table \theTabcount:} Definition of the permutations $v(\underline{i})$, $u(\underline{i})$,
$t(\underline{i})$ and $s(\underline{i})$
\index{random permutation!$v(\underline{i})$, $u(\underline{i})$, $t(\underline{i})$, $s(\underline{i})$}\\[3ex]
\hspace*{0.2ex}\begin{tabular}{@{}|c|c|c|c|c|@{}} \hline
\hspace*{11.2ex}&
\hspace*{6.3ex}$v(\underline{i})$\hspace*{6.3ex}&
\hspace*{6.3ex}$u(\underline{i})$\hspace*{6.3ex}&
\hspace*{6.3ex}$t(\underline{i})$\hspace*{6.3ex}&
\hspace*{6.3ex}$s(\underline{i})$\hspace*{6.3ex}\\ \hline
$i_{1}$&$i_{13}$&$i_{7}$&$i_{4}$&$i_{2}$\\
$i_{2}$&$i_{14}$&$i_{8}$&$i_{3}$&$i_{1}$\\ \hline
$i_{3}$&$i_{15}$&$i_{5}$&Values&\\
$i_{4}$&$i_{16}$&$i_{6}$&$\in \{\,i_{1},\ldots,i_{4}\,\}$&\\ \cline{1-4}
$i_{5}$&$i_{9}$&&\multicolumn{2}{|c|}{}\\
$i_{6}$&$i_{10}$&Values&\multicolumn{2}{|c|}{}\\
$i_{7}$&$i_{11}$&$\in \{\,i_{1},\ldots,i_{8}\,\}$&\multicolumn{2}{|c|}{}\\
$i_{8}$&$i_{12}$&&\multicolumn{2}{|c|}{}\\ \cline{1-3}
$i_{9}$&&\multicolumn{3}{|c|}{}\\
$i_{10}$&&\multicolumn{3}{|c|}{}\\
$i_{11}$&&\multicolumn{3}{|c|}{}\\
$i_{12}$&Values&\multicolumn{3}{|c|}{Remain}\\
$i_{13}$&$\in \{\,i_{1},\ldots,i_{16}\,\}$&\multicolumn{3}{|c|}{fixed}\\
$i_{14}$&&\multicolumn{3}{|c|}{}\\
$i_{15}$&&\multicolumn{3}{|c|}{}\\
$i_{16}$&&\multicolumn{3}{|c|}{}\\ \cline{1-2}
$N \setminus \{\,i_{1},\ldots,i_{16}\,\}$&\multicolumn{4}{|c|}{}\\ \hline
\end{tabular}\\[4.5ex]
For our considerations we need further useful notations:
\index{random permutation!$M_{4}$, $M_{8}$}\index{random permutation!$\overline{i}$, $\bm{\breve}{i}$}\\[2.5ex]
\refstepcounter{DSBcount}
\label{EWD_3_5_03}
\text{\hspace*{-0.8ex}(\theDSBcount)} 
\hspace*{4ex}
$\begin{array}{@{}l@{\hspace*{0.8ex}}c@{\hspace*{0.8ex}}l@{\hspace*{0.6ex}}
l@{\hspace*{0.6ex}}l@{\hspace*{3ex}}l@{}}
M_{4}&=&\biggl\{&\bm{\breve}{i} = (\,i_{1},\ldots,i_{4}\,)\,\in\,N^{4}\,:\
&\text{$\bm{\breve}{i}$ satisfies
the equivalence (t1)\index{random permutation!properties!(t1)} from (\ref{EWD_3_5_01})}\,\biggr\},\\[2.5ex]
M_{8}&=&\biggl\{&\overline{i} = (\,i_{1},\ldots,i_{8}\,)\,\in\,N^{8}\,:
&\text{$\overline{i}$ satisfies
the equivalences}\\[1ex] 
&&&&\text{(t1), (u1), (u2) and 
(u3)\index{random permutation!properties!(u1), (u2), (u3)}\index{random permutation!properties!(t1)} 
from (\ref{EWD_3_5_01})}\,\biggr\}.
\end{array}$\\[2.8ex]
Moreover, for each $\underline{i} \in M_{16}$ let\\[2.5ex]
\refstepcounter{DSBcount}
\label{EWD_3_5_04}
\text{\hspace*{-0.8ex}(\theDSBcount)} 
\hspace*{4ex}
$\gamma(\underline{i}) = \big|\,\{i_{1}, i_{2}\}\,\big| = 
\big|\,\{i_{3}, i_{4}\}
\,\big|$,\index{random permutation!$\gamma(\underline{i})$, $\theta(\underline{i})$, $\mu(\underline{i})$}\\[2.5ex]
\refstepcounter{DSBcount}
\label{EWD_3_5_05}
\text{\hspace*{-0.8ex}(\theDSBcount)} 
\hspace*{4ex}
$\theta(\underline{i}) = \big|\,\{i_{1}, i_{2}, i_{3}, i_{4}\}\,\big| = 
\big|\,\{i_{5}, i_{6}, i_{7}, i_{8}\}\,\big|$,\\[2.5ex]
\refstepcounter{DSBcount}
\label{EWD_3_5_06}
\text{\hspace*{-0.8ex}(\theDSBcount)} 
\hspace*{4ex}
$\mu(\underline{i}) = \big|\,\{\,i_{1}, i_{2}, i_{3}, i_{4},i_{5}, i_{6}, i_{7}, i_{8}\,\}\,\big|
= \big|\,\{\,i_{9}, i_{10}, i_{11}, i_{12},
i_{13}, i_{14}, i_{15}, i_{16}\,\}\,\big|$.\\[2.5ex]
The second {''$=$''} in (\ref{EWD_3_5_04}) is valid 
due to (t1).\index{random permutation!properties!(t1)} 
Similarly, the second {''$=$''} in (\ref{EWD_3_5_05}) is 
from (u1), (u2), (u3)\index{random permutation!properties!(u1), (u2), (u3)}  
and the second {''$=$''} in (\ref{EWD_3_5_06}) is 
from (v1), (v2), (v3).\index{random permutation!properties!(v1), (v2), (v3)}\\[2.8ex] 
Furthermore, let $\underline{I} = (\,I_{1},\ldots,I_{16}\,)$ be a random element on $M_{16}$ whose distribution
can be described as 
follows\index{random permutation!$\underline{I}$}\index{random permutation!$I_{1},\ldots,I_{16}$}:\\[2.3ex] 
\refstepcounter{DSBcount}
\label{EWD_3_5_07}
\text{\hspace*{-0.8ex}(\theDSBcount)} 
\hspace*{16.7ex}
$(\,I_{1}, I_{2}\,)$ is \textbf{uniformly distributed} on $N^2$.\\[1.5ex]
\refstepcounter{DSBcount}
\label{EWD_3_5_08}
\text{\hspace*{-0.8ex}(\theDSBcount)} 
\hspace*{6.2ex}
\begin{tabular}[t]{@{}c@{}}
For each $(\,i_{1},i_{2}\,)\, \in\, N^2$,\\[1.5ex]
$\bigl(\,I_{3},\,I_{4}\,\bigr)\,\Big|\,\bigl(\,I_{1},\,I_{2}\,\bigr) = \bigl(\,i_{1},\,i_{2}\,\bigr)$\\[1.5ex]
is \textbf{uniformly distributed} on 
$\bigl\{\,(\,\eta_{3},\eta_{4}\,)\,:\,(\,i_{1},i_{2},\eta_{3},\eta_{4}\,)\, \in\, M_{4}\,\bigr\}$.
\end{tabular}\\[1.5ex]
\refstepcounter{DSBcount}
\label{EWD_3_5_09}
\text{\hspace*{-0.8ex}(\theDSBcount)} 
\hspace*{12.7ex}
\begin{tabular}[t]{@{}c@{}}
For each $(\,i_{1},i_{2},i_{3},i_{4}\,)\, \in\, M_{4}$,\\[1.5ex]
$\bigl(\,I_{5},\,I_{6},\,I_{7},\,I_{8}\,\bigr)\,\Big|\,\bigl(\,I_{1},\,I_{2},\,I_{3},\,I_{4}\,\bigr)
= \bigl(\,i_{1},\,i_{2},\,i_{3},\,i_{4}\,\bigr)$\\[1.5ex]
is \textbf{uniformly distributed} on\\[1.5ex] 
$\bigl\{\,(\,\eta_{5},\eta_{6},\,\eta_{7},\eta_{8}\,)\,:
\,(\,i_{1},i_{2},,i_{3},i_{4},\eta_{5},\eta_{6},\eta_{7},\eta_{8}\,)\, \in\, M_{8}\,\bigr\}$.
\end{tabular}\\[1.5ex]
\refstepcounter{DSBcount}
\label{EWD_3_5_10}
\text{\hspace*{-0.8ex}(\theDSBcount)} 
\hspace*{8.2ex}
\begin{tabular}[t]{@{}c@{}}
For each $(\,i_{1},i_{2},\ldots,i_{8}\,)\, \in\, M_{8}$,\\[1.5ex]
$\bigl(\,I_{9},\,I_{10},\,\ldots,\,I_{16}\,\bigr)\,\Big|
\,\bigl(\,I_{1},\,I_{2},\,\ldots,\,I_{8}\,\bigr) = \bigl(\,i_{1},\,i_{2},\,\ldots,\,i_{8}\,\bigr)$\\[1.5ex]
is \textbf{uniformly distributed} on\\[1.5ex] 
$\bigl\{\,(\,\eta_{9},\eta_{10},\,\ldots,\,\eta_{16}\,)\,:
\,(\,i_{1},i_{2},\,\ldots,i_{8},\eta_{9},\eta_{10},\,\ldots,\eta_{16}\,)\, \in\, M_{16}\,\bigr\}$.
\end{tabular}\\[2.8ex]
To summarize, the multiplication rule\index{multiplication rule} 
together with (\ref{EWD_3_5_04}), (\ref{EWD_3_5_05}) 
and (\ref{EWD_3_5_06}) yields\\[2ex]
\refstepcounter{DSBcount}
\label{EWD_3_5_11}
\text{\hspace*{-0.8ex}(\theDSBcount)} 
\hspace*{4ex}
$P\Bigl(\,\underline{I} = \underline{i}\,\Bigr) =
\dfrac{\bigl( n - \mu(\underline{i}) \bigr)!}{n!} \cdot
\dfrac{\bigl( n - \theta(\underline{i}) \bigr)!}{n!} \cdot
\dfrac{\bigl( n - \gamma(\underline{i}) \bigr)!}{n!} \cdot \dfrac{1}{n^2}$
\hspace*{4ex}for each $\underline{i} \in M_{16}$.
\index{random permutation!distributions!$P\bigl(\,\underline{I} = \underline{i}\,\bigr)$}\\[2.8ex]
Additionally, let\\[2.3ex]
\refstepcounter{DSBcount}
\label{EWD_3_5_12}
\text{\hspace*{-0.8ex}(\theDSBcount)} 
\hspace*{4ex}
\begin{tabular}[t]{@{}l@{}}
$\pi_{1}$ be a random permutation that is \textbf{uniformly distributed}\\[0.5ex] 
on the set $\mathscr{P}_{n}$\index{permutation sets!$\mathscr{P}_{n}$} 
of permutations of $N$ and \textbf{independent} of $\underline{I}$.
\end{tabular}\\[2.8ex] 
Next, we compose $\pi_{1}$ and the random permutations $v(\underline{I})$, $u(\underline{I})$, 
$t(\underline{I})$, $s(\underline{I})$. In detail:\\[2.3ex]
\refstepcounter{DSBcount}
\label{EWD_3_5_13}
\text{\hspace*{-0.8ex}(\theDSBcount)} 
\hspace*{2.8ex}
\begin{tabular}[t]{@{}l@{\hspace*{3ex}}l@{\hspace*{3ex}}l@{\hspace*{3ex}}l@{}}
$\pi_{2} = \pi_{1} \circ v(\underline{I})$,&
$\pi_{3} = \pi_{2} \circ u(\underline{I})$,&
$\pi_{4} = \pi_{3} \circ t(\underline{I})$,&
$\pi_{5} = \pi_{4} \circ s(\underline{I})$.
\end{tabular}
\index{random permutation!$\pi_{1},\ldots,\pi_{5}$}\\[2.8ex]
Finally, we define:
\index{random permutation!$J_{1},\ldots,J_{16}$}\\[2.3ex]
\refstepcounter{DSBcount}
\label{EWD_3_5_14}
\text{\hspace*{-0.8ex}(\theDSBcount)} 
\hspace*{2.8ex}
\begin{tabular}{@{}l@{\hspace*{3ex}}l@{\hspace*{0.4ex}}l@{\hspace*{3ex}}l@{}}
$J_{1} = \pi_{1}(I_{9})$,&$J_{2} = \pi_{1}(I_{10})$,&\ldots ,&$J_{8} = \pi_{1}(I_{16})$,\\[1ex]
$J_{9} = \pi_{1}(I_{1})$,&$J_{10} = \pi_{1}(I_{2})$,&\ldots ,&$J_{16} = \pi_{1}(I_{8})$.
\end{tabular}\\[2.8ex]
Using the definition of $v$, $u$, $t$ and $s$, we see 
in the following Table \ref{EWDTab_3_5_02} how $\pi_{1},\ldots,\pi_{5}$ map 
$I_{1},\ldots,I_{16}$.\\[2.8ex]
\refstepcounter{Tabcount}
\label{EWDTab_3_5_02}
\hspace*{18ex}\textbf{Tabelle \theTabcount:} Values of $I_{1},\ldots,I_{16}$ under $\pi_{1},\ldots,\pi_{5}$\\[3ex]
\begin{tabular}{@{}|c|c|c|c|c|c|@{}} \hline
\hspace*{11.2ex}&
\hspace*{4.5ex}$\pi_{1}$\hspace*{4.5ex}&
\hspace*{4.5ex}$\pi_{2}$\hspace*{4.5ex}&
\hspace*{4.5ex}$\pi_{3}$\hspace*{4.5ex}&
\hspace*{4.5ex}$\pi_{4}$\hspace*{4.5ex}&
\hspace*{4.5ex}$\pi_{5}$\hspace*{4.5ex}\\ \hline
$I_{1}$&$J_{9}$&$J_{5}$&$J_{3}$&$J_{2}$&$J_{1}$\\
$I_{2}$&$J_{10}$&$J_{6}$&$J_{4}$&$J_{1}$&$J_{2}$\\ \cline{4-6}
&&&&Random&\\
$I_{3}$&$J_{11}$&$J_{7}$&$J_{1}$&variables&Same\\
$I_{4}$&$J_{12}$&$J_{8}$&$J_{2}$&$\in \sigma(I_{1},\ldots,I_{4},$&as $\pi_{4}$\\ 
&&&&$\ \ \ \ \,\,\,\,J_{1},\ldots,J_{4})$&\\ \cline{3-6}
$I_{5}$&$J_{13}$&$J_{1}$&Random&&\\
$I_{6}$&$J_{14}$&$J_{2}$&variables&Same&Same\\
$I_{7}$&$J_{15}$&$J_{3}$&$\in \sigma(I_{1},\ldots,I_{8},$&as $\pi_{3}$&as $\pi_{3}$\\
$I_{8}$&$J_{16}$&$J_{4}$&$\ \ \ \ \,\,\,\,J_{1},\ldots,J_{8})$&&\\ \cline{2-6}
$I_{9}$&$J_{1}$&&&&\\ 
$I_{10}$&$J_{2}$&&&&\\
$I_{11}$&$J_{3}$&Random&&&\\
$I_{12}$&$J_{4}$&variables&Same&Same&Same\\
$I_{13}$&$J_{5}$&$\in \sigma(I_{1},\ldots,I_{16},$&as $\pi_{2}$&as $\pi_{2}$&as $\pi_{2}$\\
$I_{14}$&$J_{6}$&$\ \ \ \ \,\,\,\,J_{1},\ldots,J_{16})$&&&\\
$I_{15}$&$J_{7}$&&&&\\
$I_{16}$&$J_{8}$&&&&\\ \hline
\end{tabular}\\[4.5ex]
As usual we define that $\sigma(X)$ is the $\sigma$-field generated by the random vector $X$ 
and $f \in \sigma(X)$ means that $f$ is measurable relative to $\sigma(X)$.
\index{random permutation!$\sigma(X)$, $f \in \sigma(X)$}\\[2.8ex] 
Further important results about the above random permutations and 
random indices\index{random index} are:\\[4ex] 
\refstepcounter{DSBcount}
\label{EWD_3_5_15}
\textbf{\hspace*{-0.8ex}\theDSBcount\ Lemma}
\begin{enumerate}
\item\label{EWD_3_5_15_BWa}
$\pi_{1}$, $\pi_{2}$, $\pi_{3}$, $\pi_{4}$, $\pi_{5}$ 
are uniformly distributed on $\mathscr{P}_{n}$\index{permutation sets!$\mathscr{P}_{n}$}
and independent of $\underline{I} = (\,I_{1}, \ldots, I_{16}\,)$.
\item\label{EWD_3_5_15_BWb}
$\pi_{1}$ and $(\,I_{1}, \ldots, I_{8}, J_{1}, \ldots, J_{8}\,)$ are independent.
\item\label{EWD_3_5_15_BWc}
$\pi_{2}$ and $(\,I_{1}, \ldots, I_{4}, J_{1}, \ldots, J_{4}\,)$ are independent.
\item\label{EWD_3_5_15_BWd}
$\pi_{3}$ and $(\,I_{1}, I_{2}, J_{1}, J_{2}\,)$ are independent.
\item\label{EWD_3_5_15_BWe}
$\pi_{4}$ and $(\,I_{1}, J_{1}\,)$ are independent.
\item\label{EWD_3_5_15_BWf}
$\pi_{5}$ and $I_{1}$ are independent.
\item\label{EWD_3_5_15_BWg}
$I_{1},\ldots,I_{16}$ are uniformly distributed on $N$.
\item\label{EWD_3_5_15_BWh}
$\bigl(\,I_{l}, \pi_{k}(I_{l})\,\bigr)$ are uniformly distributed on $N^2$ for all $1 \leq l \leq 16$, 
$1 \leq k \leq 5$.
\end{enumerate}
\vspace*{3.5ex}
\textbf{Proof:}
\begin{enumerate}
\item
Since $\pi_{1}$ is uniformly distributed on the set $\mathscr{P}_{n}$\index{permutation sets!$\mathscr{P}_{n}$} 
and since $\pi_{1}$ is also independent of $\underline{I}$, (cf. (\ref{EWD_3_5_12})), 
it follows for all $\pi \in \mathscr{P}_{n}$ and $\underline{i} \in M_{16}$ that\\[2ex]
\hspace*{7.2ex}\begin{tabular}[t]{@{}l@{\hspace*{0.8ex}}c@{\hspace*{0.8ex}}l@{}}
$P\bigl(\,\pi_{2} = \pi,\, \underline{I} = \underline{i}\, \bigr)$&$=$&
$P\bigl(\,\pi_{1} \circ v(\underline{i}) = \pi,\, \underline{I} = \underline{i}\, \bigr)$
\hspace*{7.9ex}(cf. (\ref{EWD_3_5_13}): $\pi_{2} = \pi_{1} \circ v(\underline{I})$)\\[1.2ex]
&$=$&$P\bigl(\,\pi_{1} = \pi \circ v(\underline{i})^{-1},\, \underline{I} = \underline{i}\, \bigr)$\\[1.2ex]
&$=$&$P\bigl(\,\pi_{1} = \pi \circ v(\underline{i})^{-1}\,\bigr)\,P\bigl(\, \underline{I} 
= \underline{i}\, \bigr)$\\[1.2ex]
&$=$&$\dfrac{1}{n!}\,P\bigl(\, \underline{I} = \underline{i}\, \bigr)$.
\end{tabular}\\[2.5ex]
Thus summation over all $\underline{i} \in M_{16}$ gives first\\[2ex]
\hspace*{7.2ex}$P\bigl(\,\pi_{2} = \pi\,\bigr) = \dfrac{1}{n!}$\\[2.5ex]
and then the independence of $\pi_{2}$ and $\underline{I}$.\\[1.5ex]
From this we obtain the assertion for $\pi_{3}$, then for $\pi_{4}$ 
and finally for $\pi_{5}$ in a completely analogous way.
\item
Since $\underline{I}$ and $\pi_{1}$ are independent (cf. (\ref{EWD_3_5_12})), we get
for all $\underline{i} = (\,i_{1},\ldots,i_{16}\,)\,\in\,M_{16}$ and $\pi \in \mathscr{P}_{n}$:\\[2ex]
\hspace*{7.2ex}\begin{tabular}[t]{@{}c@{\hspace*{0.8ex}}l@{}}
&$P\Bigl(\,\bigl(\,I_{1}, \ldots, I_{8}, J_{1}, \ldots, J_{8}\,\bigr) = \underline{i},\,\pi_{1} = \pi\,\Bigr)$\\[1.4ex]
$=$&$P\Bigl(\,\bigl(\,I_{1}, \ldots, I_{8}, \pi(I_{9}), \ldots, \pi(I_{16})\,\bigr) = \underline{i},
\,\pi_{1} = \pi\,\Bigr)$
\hspace*{19.5ex}(cf. (\ref{EWD_3_5_14}))\\[1.4ex]
$=$&$P\Bigl(\,\underline{I} = \bigl(\,i_{1}, \ldots, i_{8}, \pi^{-1}(i_{9}), \ldots, \pi^{-1}(i_{16})\,\bigr),
\,\pi_{1} = \pi\,\Bigr)$\\[1.4ex]
$=$&$P\Bigl(\,\underline{I} = \bigl(\,i_{1}, \ldots, i_{8}, \pi^{-1}(i_{9}), \ldots, \pi^{-1}(i_{16})\,\bigr)\,\Bigr)\,
P\Bigl(\,\pi_{1} = \pi\,\Bigr)$\\[1.4ex]
$=$&$P\Bigl(\,\underline{I} = \underline{i}\,\Bigr)\,
P\Bigl(\,\pi_{1} = \pi\,\Bigr)$.
\end{tabular}\\[2.5ex]
The last equation is valid, since we can use (\ref{EWD_3_5_10}) for\\[2ex]
\hspace*{3ex}$\underline{i} = \bigl(\,i_{1}, \ldots, i_{8}, i_{9}, \ldots, i_{16}\,\bigr)\in M_{16}$ and
$\underline{\xi} = 
\bigl(\,i_{1}, \ldots, i_{8}, \pi^{-1}(i_{9}), \ldots, \pi^{-1}(i_{16})\,\bigr) \in M_{16}$.\\[2.5ex]
Remember that $\pi^{-1} \in \mathscr{P}_{n}$\index{permutation sets!$\mathscr{P}_{n}$} 
is bijective and thus $\underline{\xi} \in M_{16}$
due to (\ref{EWD_3_5_01}), (v1), (v2) and (v3).\index{random permutation!properties!(v1), (v2), (v3)}
Now summation over all $\pi \in \mathscr{P}_{n}$ gives first\\[2ex]
\refstepcounter{DSBcount}
\label{EWD_3_5_16}
\text{\hspace*{-0.8ex}(\theDSBcount)} 
\hspace*{4ex}
$P\Bigl(\,\bigl(\,I_{1}, \ldots, I_{8}, J_{1}, \ldots, J_{8}\,\bigr) = \underline{i}\,\Bigr)
= P\Bigl(\,\underline{I} = \underline{i}\,\Bigr)$\\[2.5ex]
and then the assertion. 
\item
The proof of part \ref{EWD_3_5_15_BWc}) is completely analogous to that of part \ref{EWD_3_5_15_BWb}), if we use
\begin{itemize}
\item
$\pi_{2}$ instead of $\pi_{1}$, 
\item
$J_{1} = \pi_{2}(I_{5})$, $J_{2} = \pi_{2}(I_{6})$, $J_{3} = \pi_{2}(I_{7})$, $J_{4} = \pi_{2}(I_{8})$
from Table \ref{EWDTab_3_5_02},
\item
$(\,I_{1}, \ldots, I_{8}\,)$ and $\pi_{2}$ are independent due to part \ref{EWD_3_5_15_BWa}),
\item
(\ref{EWD_3_5_09}) instead of (\ref{EWD_3_5_10}) for\\[2ex]
\hspace*{7.2ex}$\bigl(\,i_{1}, \ldots, i_{4}, i_{5}, \ldots, i_{8}\,\bigr)\in M_{8}$ and
$\bigl(\,i_{1}, \ldots, i_{4}, \pi^{-1}(i_{5}), \ldots, \pi^{-1}(i_{8})\,\bigr) \in M_{8}$,\vspace*{1.5ex}
\item
$\bigl(\,i_{1}, \ldots, i_{4}, \pi^{-1}(i_{5}), \ldots, \pi^{-1}(i_{8})\,\bigr) \in M_{8}$ due to
(\ref{EWD_3_5_01}), (u1), (u2) and (u3).\index{random permutation!properties!(u1), (u2), (u3)}
\end{itemize}
\item
The proof of part \ref{EWD_3_5_15_BWd}) is completely analogous to that of part \ref{EWD_3_5_15_BWb}), if we use
\begin{itemize}
\item
$\pi_{3}$ instead of $\pi_{1}$, 
\item
$J_{1} = \pi_{3}(I_{3})$, $J_{2} = \pi_{3}(I_{4})$ from Table \ref{EWDTab_3_5_02},
\item
$(\,I_{1}, \ldots, I_{4}\,)$ and $\pi_{3}$ are independent due to part \ref{EWD_3_5_15_BWa}),
\item
(\ref{EWD_3_5_08}) instead of (\ref{EWD_3_5_10}) for\\[2ex]
\hspace*{7.2ex}$\bigl(\,i_{1}, i_{2}, i_{3}, i_{4}\,\bigr)\in M_{4}$ and
$\bigl(\,i_{1}, i_{2}, \pi^{-1}(i_{3}), \pi^{-1}(i_{4})\,\bigr) \in M_{4}$,\vspace*{1.5ex}
\item
$\bigl(\,i_{1}, i_{2}, \pi^{-1}(i_{3}), \pi^{-1}(i_{4})\,\bigr) \in M_{4}$
due to
(\ref{EWD_3_5_01}), (t1).\index{random permutation!properties!(t1)}
\end{itemize}
\item
The proof of part \ref{EWD_3_5_15_BWe}) is completely analogous to that of part \ref{EWD_3_5_15_BWb}), if we use
\begin{itemize}
\item
$\pi_{4}$ instead of $\pi_{1}$, 
\item
$J_{1} = \pi_{4}(I_{2})$ from Table \ref{EWDTab_3_5_02},
\item
$(\,I_{1}, I_{2}\,)$ and $\pi_{4}$ are independent due to part \ref{EWD_3_5_15_BWa}),
\item
(\ref{EWD_3_5_07}) instead of (\ref{EWD_3_5_10}) for\\[2ex]
\hspace*{7.2ex}$\bigl(\,i_{1}, i_{2}\,\bigr) \in N^2$ and
$\bigl(\,i_{1}, \pi^{-1}(i_{2})\,\bigr) \in N^2$.
\end{itemize}
\item
$I_{1}$ and $\pi_{5}$ are independent due to part \ref{EWD_3_5_15_BWa}).
\item
Trivially, we have for all $1 \leq l, k \leq 16$ and $i_{l}, i_{k} \in N$:\\[2ex]
\hspace*{13.5ex}$i_{l} = i_{k}\ \Leftrightarrow\ i_{l} + 1 = i_{k} + 1\ (\text{mod}\ n)
\ \Leftrightarrow\ i_{l} - 1 = i_{k} - 1$\ (\text{mod}\ n).\\[2ex] 
Therefore the mapping\\[2ex]
\hspace*{13.6ex}$\tau\,:\,M_{16}\, \longrightarrow\, M_{16}\,:\,
\underline{i}\, \longrightarrow\, \underline{i} + (1,1,\ldots,1)\ (\text{mod}\ n)$\\[2ex]
is \textbf{bijective} with the inverse mapping\\[2ex]
\hspace*{13.5ex}$\tau^{-1}\,:\,M_{16}\, \longrightarrow\, M_{16}\,:\,
\underline{i}\, \longrightarrow\, \underline{i} - (1,1,\ldots,1)\ (\text{mod}\ n)$.\\[2ex]
Since also\\[2ex]
\hspace*{13.5ex}$\gamma(\underline{i}) = \gamma\bigl(\tau(\underline{i})\bigr)$, 
$\theta(\underline{i}) = \theta\bigl(\tau(\underline{i})\bigr)$ 
and $\mu(\underline{i}) = \mu\bigl(\tau(\underline{i})\bigr)$,\\[2ex] 
we get using (\ref{EWD_3_5_11}):\\[2ex]
\refstepcounter{DSBcount}
\label{EWD_3_5_17}
\text{\hspace*{-0.8ex}(\theDSBcount)} 
\hspace*{4ex}
$P\bigl(\,\underline{I} = \underline{i}\,\bigr) =
P\bigl(\,\underline{I} = \tau(\underline{i})\,\bigr)$
\hspace*{3ex}for all $\underline{i} \in M_{16}$.\\[2ex]
From this we immediately obtain for all $1 \leq l \leq 16$ and $1 \leq r \leq n - 1$:\\[2ex]
\hspace*{13.5ex}\begin{tabular}[t]{@{}l@{\hspace*{0.8ex}}c@{\hspace*{0.8ex}}l@{\hspace*{26.8ex}}r@{}}
$P(I_{l} = r)$&$=$&$\displaystyle{\sum\limits_{\begin{array}{@{}c@{}}\\[-5ex]
\scriptstyle{\underline{i}\, \in\, M_{16}}\\[-1.5ex]\scriptstyle{i_{l}\, =\, r}\end{array}}
P\bigl(\,\underline{I} = \underline{i}\,\bigr)}$\\[3.5ex]
&$=$&$\displaystyle{\sum\limits_{\begin{array}{@{}c@{}}\\[-5ex]
\scriptstyle{\underline{i}\, \in\, M_{16}}\\[-1.5ex]\scriptstyle{i_{l}\, =\, r}\end{array}}
P\bigl(\,\underline{I} = \tau(\underline{i})\,\bigr)}$&(cf. (\ref{EWD_3_5_17}))\\[3.5ex]
&$=$&$\displaystyle{\sum\limits_{\begin{array}{@{}c@{}}\\[-5ex]
\scriptstyle{\underline{i}\, \in\, M_{16}}\\[-1.5ex]\scriptstyle{i_{l}\, =\, r + 1}\end{array}}
P\bigl(\,\underline{I} = \underline{i}\,\bigr)}$&
$\Bigl(\,\underline{i} = \tau\bigl(\,\tau^{-1}(\underline{i})\,\bigr)\,\Bigr)$\\[3.5ex]
&$=$&$\displaystyle{P(I_{l} = r + 1)}$.
\end{tabular}
\item
We get for all $(i,j) \in N^2$:\\[2.3ex]
\hspace*{7.2ex}\begin{tabular}[t]{@{}l@{\hspace*{0.8ex}}c@{\hspace*{0.8ex}}l@{\hspace*{3.3ex}}r@{}}
$P\Bigl(\,\bigl(\,I_{l}, \pi_{k}(I_{l})\,\bigr) = (\,i,j\,)\,\Bigr)$&
$=$&$\displaystyle{\sum\limits_{\pi\, \in\, \mathscr{P}_{n}}}
P\Bigl(\,\bigl(\,I_{l}, \pi(I_{l})\,\bigr) = (\,i,j\,),\,\pi_{k} = \pi\,\Bigr)$\\[3.2ex]
&$=$&$\displaystyle{\sum\limits_{\begin{array}{@{}c@{}}\\[-5ex]
\scriptstyle{\pi\, \in\, \mathscr{P}_{n}}\\[-1.5ex]\scriptstyle{\pi(i)\, =\, j}\end{array}}
P\Bigl(\,I_{l} = i,\,\pi_{k} = \pi\,\Bigr)}$\\[1.6ex]
&$=$&$\displaystyle{\sum\limits_{\begin{array}{@{}c@{}}\\[-5ex]
\scriptstyle{\pi\, \in\, \mathscr{P}_{n}}\\[-1.5ex]\scriptstyle{\pi(i)\, =\, j}\end{array}}
P\Bigl(\,I_{l} = i\,\Bigr)\,P\Bigl(\,\pi_{k} = \pi\,\Bigr)}$&(cf. part \ref{EWD_3_5_15_BWa}))\\[1.6ex]
&$=$&$(n-1)!\,\dfrac{1}{n}\,\dfrac{1}{n!} = \dfrac{1}{n^2}$.
\end{tabular}\\[3ex]
For the last equation, we used that $I_l$ is uniformly distributed on $N$ (see part \ref{EWD_3_5_15_BWg}))
and $\pi_{k}$ is uniformly distributed on $\mathscr{P}_{n}$\index{permutation sets!$\mathscr{P}_{n}$} 
(see part \ref{EWD_3_5_15_BWa})).
\hspace*{1ex}\hfill$\Box$
\end{enumerate}
\vspace*{3ex}
In the proof of Lemma \ref{EWD_3_5_15}, which we have just given, we also proved 
equalities of distributions (see e.g. (\ref{EWD_3_5_16})), which we 
list again separately in the following corollary.\\[2.8ex]
In addition, the proof structure for the part \ref{EWD_3_5_15_BWh}) of this lemma can also be used 
for other statements.\\[4ex] 
\refstepcounter{DSBcount}
\label{EWD_3_5_18}
\textbf{\hspace*{-0.8ex}\theDSBcount\ Corollary}\begin{enumerate}
\item\label{EWD_3_5_18_BWa}
$(\,I_{1}, \ldots, I_{8}, I_{9}, \ldots, I_{16}\,)$ and $(\,I_{1}, \ldots, I_{8}, J_{1}, \ldots, J_{8}\,)$
have the same distribution.
\item\label{EWD_3_5_18_BWb}
$(\,I_{1}, \ldots, I_{4}, I_{5}, \ldots, I_{8}\,)$ and $(\,I_{1}, \ldots, I_{4}, J_{1}, \ldots, J_{4}\,)$
have the same distribution.
\item\label{EWD_3_5_18_BWc}
$(\,I_{1}, I_{2}, I_{3}, I_{4}\,)$ and $(\,I_{1}, I_{2}, J_{1}, J_{2}\,)$
have the same distribution.
\item\label{EWD_3_5_18_BWd}
Let $\underline{i} = \bigl(\,i_{1},\ldots, i_{16}\,\bigr) \in M_{16}$ and 
$\underline{j} = \bigl(\,j_{1},\ldots, j_{16}\,\bigr) \in M_{16}$\index{random permutation!$\underline{j}$} 
\index{random permutation!$j_{1},\ldots,j_{16}$}have the properties\\[1.8ex]
\refstepcounter{DSBcount}
\label{EWD_3_5_19}
\text{\hspace*{-0.8ex}(\theDSBcount)} 
\hspace*{4ex}
$\begin{array}[t]{@{}l@{\hspace*{3ex}}l@{\hspace*{1.6ex}}
l@{\hspace*{1.6ex}}l@{\hspace*{3ex}}l@{}}
\text{(w1)}&i_{l} = i_{k}&\Leftrightarrow&j_{l+8} = j_{k+8}&\text{for}\ 1 \leq l \leq 8,\ \ 1 \leq k \leq 8,\\[1.2ex]
\text{(w2)}&i_{l} = i_{k}&\Leftrightarrow&j_{l-8} = j_{k-8}&\text{for}\ 9 \leq l \leq 16,\ \ 9 \leq k \leq 16,\\[1ex]
\text{(w3)}&i_{l} = i_{k}&\Leftrightarrow&j_{l+8} = j_{k-8}&\text{for}\ 1 \leq l \leq 8,\ \ 9 \leq k \leq 16.
\end{array}$\index{random permutation!properties!(w1), (w2), (w3)}\\[2.2ex]
Furthermore, let $\underline{J} = (\,J_{1},\ldots,J_{16}\,)$\index{random permutation!$\underline{J}$}
and $\rho(\underline{i}) = \big|\,\{\,i_{1}, i_{2}, \ldots, 
i_{15}, i_{16}\,\}\,\big|$. Then\index{random permutation!$\rho(\underline{i})$}\\[3ex]
\refstepcounter{DSBcount}
\label{EWD_3_5_20}
\text{\hspace*{-0.8ex}(\theDSBcount)} 
\hspace*{4ex}
\begin{tabular}[t]{@{}l@{\hspace*{0.8ex}}c@{\hspace*{0.8ex}}l@{}}
$P\Bigl(\,\underline{I} = \underline{i},\,\underline{J} = \underline{j}\,\Bigr)$&$=$&$ 
\dfrac{\bigl( n - \rho(\underline{i}) \bigr)!}{n!} \cdot
\dfrac{\bigl( n - \mu(\underline{i}) \bigr)!}{n!} \cdot
\dfrac{\bigl( n - \theta(\underline{i}) \bigr)!}{n!} \cdot
\dfrac{\bigl( n - \gamma(\underline{i}) \bigr)!}{n!} \cdot \dfrac{1}{n^2}$.
\end{tabular}
\index{random permutation!distributions!$P\bigl(\,\underline{I} = \underline{i},\,\underline{J} = \underline{j}\,\bigr)$}
\end{enumerate}
\newpage
\textbf{Proof:}\\[2ex]
\hspace*{2.4ex}a) - c)\hspace*{1.4ex}See proof of Lemma \ref{EWD_3_5_15}, especially part \ref{EWD_3_5_15_BWb}), but also the parts \ref{EWD_3_5_15_BWc}) and \ref{EWD_3_5_15_BWd}).
\begin{enumerate}
\addtocounter{enumi}{+3}
\item
We obtain for $\underline{i} = \bigl(\,i_{1},\ldots, i_{16}\,\bigr) \in M_{16}$ and 
$\underline{j} = \bigl(\,j_{1},\ldots, j_{16}\,\bigr) \in M_{16}$ with the properties (w1), (w2) 
and (w3):\index{random permutation!properties!(w1), (w2), (w3)}\\[2.3ex]
\hspace*{3ex}\begin{tabular}[t]{@{}c@{\hspace*{0.8ex}}l@{}}
&$P\Bigl(\,\bigl(\,I_{1}, \ldots, I_{8}, I_{9}, \ldots, I_{16}\,\bigr) = 
\bigl(\,i_{1}, \ldots, i_{8}, i_{9}, \ldots, i_{16}\,\bigr),$\\[1.4ex]
&\hspace*{9.2ex}$\bigl(\,J_{1}, \ldots, J_{8}, J_{9}, \ldots, J_{16}\,\bigr) = 
\bigl(\,j_{1}, \ldots, j_{8}, j_{9}, \ldots, j_{16}\,\bigr)\,
\Bigr)$\\[2.5ex]
$=$&$\displaystyle{\sum\limits_{\pi\, \in\, \mathscr{P}_{n}}
P\Bigl(\,\bigl(\,I_{1}, \ldots, I_{8}, I_{9}, \ldots, I_{16}\,\bigr) = 
\bigl(\,i_{1}, \ldots,i_{8}, i_{9}, \ldots, i_{16}\,\bigr),}$\\[1.4ex]
&\hspace*{9.2ex}$\bigl(\,\pi(I_{9}), \ldots, \pi(I_{16}), \pi(I_{1}), \ldots, \pi(I_{8})\,\bigr) = 
\bigl(\,j_{1}, \ldots, j_{8}, j_{9}, \ldots, j_{16}\,\bigr),\,\pi_{1} = \pi\,\Bigr)$\\[1.4ex]
&\hspace*{71.7ex}(cf. (\ref{EWD_3_5_14}))\\[1.4ex]
$=$&$\displaystyle{\sum\limits_{\begin{array}{@{}c@{}}\\[-5ex]
\scriptstyle{\pi\, \in\, \mathscr{P}_{n}}\\[-1.5ex]
\scriptstyle{\pi(i_{9})\, =\, j_{1},\ldots,\,\pi(i_{16})\, =\, j_{8}}\\[-1.5ex]
\scriptstyle{\pi(i_{1})\, =\, j_{9},\ldots,\,\pi(i_{8})\, =\, j_{16}}
\end{array}}
P\Bigl(\,\bigl(\,I_{1}, \ldots, I_{16}\,\bigr) = 
\bigl(\,i_{1},\ldots, i_{16}\,\bigr),\,\pi_{1} = \pi\,\Bigr)}$\\[2.5ex]
$=$&$\displaystyle{\sum\limits_{\begin{array}{@{}c@{}}\\[-5ex]
\scriptstyle{\pi\, \in\, \mathscr{P}_{n}}\\[-1.5ex]
\scriptstyle{\pi(i_{9})\, =\, j_{1},\ldots,\,\pi(i_{16})\, =\, j_{8}}\\[-1.5ex]
\scriptstyle{\pi(i_{1})\, =\, j_{9},\ldots,\,\pi(i_{8})\, =\, j_{16}}
\end{array}}
P\Bigl(\,\bigl(\,I_{1}, \ldots, I_{16}\,\bigr) = 
\bigl(\,i_{1},\ldots, i_{16}\,\bigr)\,\Bigr)\,P\Bigl(\,\pi_{1} = \pi\,\Bigr)}$
\hspace*{5.1ex}(cf. (\ref{EWD_3_5_12}))\\[9.5ex]
$=$&$(n - \rho(\underline{i}))!\,P\Bigl(\,\bigl(\,I_{1}, \ldots, I_{16}\,\bigr) = 
\bigl(\,i_{1},\ldots, i_{16}\,\bigr)\,\Bigr)\,\dfrac{1}{n!}$.
\end{tabular}\\[3ex]
By using (\ref{EWD_3_5_11}) we then get the assertion (\ref{EWD_3_5_20}).
\hspace*{1ex}\hfill$\Box$
\end{enumerate}
\vspace*{3ex}
Next, we define
\index{random permutation!$T_{i}$, $\Delta T_{i}$}\\[2.5ex]
\refstepcounter{DSBcount}
\label{EWD_3_5_21}
\text{\hspace*{-0.8ex}(\theDSBcount)}
\hspace*{4ex}
\begin{tabular}{@{}r@{\hspace*{0.8ex}}c@{\hspace*{0.8ex}}l@{\hspace*{3ex}}l@{}}
$T_{i}$&$=$&$\displaystyle{\sum\limits_{j = 1}^{n} a_{j\pi_{i}(j)}}$&for $i = 1, 2, 3, 4, 5$;\\[3.5ex]
$\Delta T_{i}$&$=$&$T_{i+1} - T_{i}$&for $i = 1, 2, 3, 4$.
\end{tabular}\\[3.5ex]
According to our construction, the following representations and measurability statements are valid for
the $\Delta T_{i}$ with regard to the random indices $I_{1},\ldots,I_{16}$ 
and $J_{1},\ldots,J_{16}$.\index{random permutation!$\sigma(X)$, $f \in \sigma(X)$}\index{random index}\\[4ex]
\refstepcounter{DSBcount}
\label{EWD_3_5_22}
\textbf{\hspace*{-0.8ex}\theDSBcount\ Lemma}
\begin{enumerate}
\item\label{EWD_3_5_22_BWa}
$\displaystyle{\Delta T_{1} = \sum\limits_{j = 1}^{16} \Bigl(\,a_{I_{j}\pi_{2}(I_{j})} -
a_{I_{j}\pi_{1}(I_{j})}\,\Bigr) \in \sigma\bigl(\,I_{1},\ldots,I_{16},J_{1},\ldots,J_{16}\,\bigr)}$,
\item\label{EWD_3_5_22_BWb}
$\displaystyle{\Delta T_{2} = \sum\limits_{j = 1}^{8} \Bigl(\,a_{I_{j}\pi_{3}(I_{j})} -
a_{I_{j}\pi_{2}(I_{j})}\,\Bigr) \in \sigma\bigl(\,I_{1},\ldots,I_{8},J_{1},\ldots,J_{8}\,\bigr)}$,
\item\label{EWD_3_5_22_BWc}
$\displaystyle{\Delta T_{3} = \sum\limits_{j = 1}^{4} \Bigl(\,a_{I_{j}\pi_{4}(I_{j})} -
a_{I_{j}\pi_{3}(I_{j})}\,\Bigr) \in \sigma\bigl(\,I_{1},\ldots,I_{4},J_{1},\ldots,J_{4}\,\bigr)}$,
\item\label{EWD_3_5_22_BWd}
$\displaystyle{\Delta T_{4} = a_{I_{1}J_{1}} + a_{I_{2}J_{2}} - a_{I_{1}J_{2}} - a_{I_{2}J_{1}}
\in \sigma\bigl(\,I_{1}, I_{2}, J_{1}, J_{2}\,\bigr)}$.
\end{enumerate}
\vspace*{3.5ex}
\textbf{Proof:}\\[0.8ex]
The representations of the $\Delta T_{i}$, $i = 1, 2, 3, 4$, follow from the fact that
$v(\underline{i})$ leaves the numbers outside of $\{\,i_{1},\ldots,i_{16}\,\}$ fixed,
$u(\underline{i})$ leaves the numbers outside of $\{\,i_{1},\ldots,i_{8}\,\}$ fixed,
$t(\underline{i})$ leaves the numbers outside of $\{\,i_{1},\ldots,i_{4}\,\}$ fixed, and
$s(\underline{i})$ leaves the numbers outside of $\{\,i_{1}, i_{2}\,\}$ fixed
(see Table \ref{EWDTab_3_5_01}).\\[2.8ex]
The measurability statements, on the other hand, can be found 
in the other Table \ref{EWDTab_3_5_02}.\hspace*{1ex}\hfill$\Box$\\[4ex]
The following lemma summarizes the statements of the Lemmata \ref{EWD_3_5_15} and \ref{EWD_3_5_22}.
This gives us the necessary independence statements for the application of Stein's method in the following sections \ref{EWD_Kap3_Sec6}, \ref{EWD_Kap3_Sec7} and \ref{EWD_Kap3_Sec8}.\\[4ex]
\refstepcounter{DSBcount}
\label{EWD_3_5_23}
\textbf{\hspace*{-0.8ex}\theDSBcount\ Lemma}
\begin{enumerate}
\item\label{EWD_3_5_23_BWa}
$T_{4}$ and $a_{I_{1}J_{1}}$ are independent.
\item\label{EWD_3_5_23_BWb}
$T_{3}$ and $\bigl(\,a_{I_{1}J_{1}},\,\Delta T_{4}\,\bigr)$ are independent.
\item\label{EWD_3_5_23_BWc}
$T_{2}$ and $\bigl(\,a_{I_{1}J_{1}},\,\Delta T_{4},\,\Delta T_{3}\,\bigr)$ are independent.
\item\label{EWD_3_5_23_BWd}
$T_{1}$ and $\bigl(\,a_{I_{1}J_{1}},\,\Delta T_{4},\,\Delta T_{3},\,\Delta T_{2}\,\bigr)$ are independent.
\end{enumerate}
\vspace*{3.5ex}
\textbf{Proof:}
\begin{enumerate}
\item
Lemma \ref{EWD_3_5_15}, \ref{EWD_3_5_15_BWe}).
\item
Lemma \ref{EWD_3_5_15}, \ref{EWD_3_5_15_BWd}) and Lemma \ref{EWD_3_5_22}, \ref{EWD_3_5_22_BWd}).
\item
Lemma \ref{EWD_3_5_15}, \ref{EWD_3_5_15_BWc}) and Lemma \ref{EWD_3_5_22}, \ref{EWD_3_5_22_BWc}), \ref{EWD_3_5_22_BWd}).
\item
Lemma \ref{EWD_3_5_15}, \ref{EWD_3_5_15_BWb}) and Lemma \ref{EWD_3_5_22}, \ref{EWD_3_5_22_BWb}), \ref{EWD_3_5_22_BWc}), 
\ref{EWD_3_5_22_BWd}).
\hspace*{1ex}\hfill$\Box$
\end{enumerate}
\refstepcounter{DSBcount}
\label{EWD_3_5_24}
\textbf{\theDSBcount\ Remark}\\[0.8ex]
Illustration of the dependencies between the $T_{i}$ and $\Delta T_{j}$:\\[2ex]
\begin{tabular}{l@{\hspace*{0.8ex}}c@{\hspace*{0.8ex}}l@{}}
$T_{A} = T_{5}$&$=$&\colorbox{brilliantlavender}{$T_{4}$} $+$ \colorbox{lightskyblue}{$\Delta T_{4}$}\\[1.5ex]
&$=$&\colorbox{brilliantlavender}{$T_{3}$} $+$ \colorbox{lightskyblue}{$\Delta T_{3}$} 
$+$ \colorbox{brilliantlavender}{$\Delta T_{4}$}\\[1.5ex]
&$=$&\colorbox{brilliantlavender}{$T_{2}$} $+$ \colorbox{lightskyblue}{$\Delta T_{2}$} 
$+$ \colorbox{brilliantlavender}{$\Delta T_{3} + \Delta T_{4}$}\\[1.5ex]
&$=$&\colorbox{brilliantlavender}{$T_{1}$} $+$ \colorbox{lightskyblue}{$\Delta T_{1}$} 
$+$ \colorbox{brilliantlavender}{$\Delta T_{2} + \Delta T_{3} + \Delta T_{4}$}
\end{tabular}
\hspace*{14.5ex}
\begin{tabular}{@{}c@{}}
For each line:\\[1ex] 
same colour = independent.
\end{tabular}\\[4ex]
The following lemma contains some useful and repeatedly used results when applying Stein's method 
to the random variables $T_{i}$, $\Delta T_{j}$, $a_{I_{1}J_{1}}$, etc.\\[4ex]
\refstepcounter{DSBcount}
\label{EWD_3_5_25}
\textbf{\hspace*{-0.8ex}\theDSBcount\ Lemma}
\begin{enumerate}
\item\label{EWD_3_5_25_BWa}
Let $f : \mathbb{R} \rightarrow \mathbb{R}$ be a function.
Then\\[2ex]
\hspace*{9.1ex}$E\bigl(\,T_{A}\,f(T_{A})\,\bigr) = E\bigl(\,T_{5}\,f(T_{5})\,\bigr)
= n\,E\bigl(\,a_{I_{1}J_{1}}\,f(T_{5})\,\bigr)$.
\item\label{EWD_3_5_25_BWb}
Let $r \in \mathbb{N}$ and $l_{1}, l_{2},\ldots,l_{r} \in \{\,1, 2, \ldots, 16\,\}$ and 
$k_{1}, k_{2},\ldots,k_{r} \in \{\,1, 2, \ldots, 5\,\}$. Then\\[2ex]
\hspace*{9.1ex}$\displaystyle{n\,E\Bigl(\,\big|\,a_{I_{l_{1}}\pi_{k_{1}}(I_{l_{1}})}\,\big|\,  
\big|\,a_{I_{l_{2}}\pi_{k_{2}}(I_{l_{2}})}\,\big|\,\ldots\,
\big|\,a_{I_{l_{r}}\pi_{k_{r}}(I_{l_{r}})}\,\big|\,\,\Bigr) 
\leq \dfrac{1}{n}\,\sum\limits_{i, j=1}^{n}\,|\,a_{ij}\,|^r}$.
\item\label{EWD_3_5_25_BWc}
Let $e_{1}, e_{2}, e_{3}, e_{4} \in \mathbb{N}_{0}$ and $e = 1 + e_{1} + e_{2} + e_{3} + e_{4}$. Then\\[2ex]
\hspace*{9.1ex}$\displaystyle{n\,E\Bigl(\,\big|\,a_{I_{1}J_{1}}\,\big|\,\big|\,\Delta T_{4}\,\big|^{e_{4}}\,\big|\,
\big|\,\Delta T_{3}\,\big|^{e_{3}}\,\big|\,\big|\,\Delta T_{2}\,\big|^{e_{2}}\,\big|\,
\big|\,\Delta T_{1}\,\big|^{e_{1}}\,\Bigr)
\leq C\,\dfrac{1}{n}\,\sum\limits_{i, j=1}^{n}\,|\,a_{ij}\,|^e}$,\\[2.5ex]
where $C = 4^{e_{4}}\,8^{e_{3}}\,16^{e_{2}}\,32^{e_{1}} \leq 32^{e - 1}$.
\end{enumerate}
\vspace*{3.5ex}
\textbf{Proof:}
\begin{enumerate}
\item
\hspace*{4ex}\begin{tabular}[t]{@{}l@{\hspace*{0.8ex}}l@{\hspace*{13.7ex}}r@{}}
&$n\,E\bigl(\,a_{I_{1}J_{1}}\,f(\,T_{5}\,)\,\bigr)$\\[2.5ex]
$=$&$n\,E\bigl(\,a_{I_{1}\pi_{5}(I_{1})}\,f(\,T_{5}\,)\,\bigr)$
&(cf. Table \ref{EWDTab_3_5_02}: $J_{1} = \pi_{5}(I_{1})$)\\[2ex] 
$=$&$n\,\displaystyle{\sum\limits_{i = 1}^{n} P(I_{1} = i)\,E\bigl(\,a_{i\pi_{5}(i)}\,f(\,T_{5}\,)\,
\big|\,I_{1} = i\,\bigr)}$
&(law of total probability)\index{law of total probability}\\[3.5ex]
$=$&$n\,\displaystyle{\sum\limits_{i = 1}^{n} P(I_{1} = i)\,E\bigl(\,a_{i\pi_{5}(i)}\,f(\,T_{5}\,)\,\bigr)}$
&(cf. Lemma \ref{EWD_3_5_15}, \ref{EWD_3_5_15_BWf}))
\end{tabular}\\[4ex]
\hspace*{4ex}\begin{tabular}[t]{@{}l@{\hspace*{0.8ex}}l@{\hspace*{34.9ex}}r@{}}
$=$&$n\,\displaystyle{\sum\limits_{i = 1}^{n} \dfrac{1}{n}\,E\bigl(\,a_{i\pi_{5}(i)}\,f(\,T_{5}\,)\,\bigr)}$
&(cf. Lemma \ref{EWD_3_5_15}, \ref{EWD_3_5_15_BWg}))\\[4ex]
$=$&$E\bigl(\,T_{5}\,f(\,T_{5}\,)\,\bigr)$
&($E$ linear and (\ref{EWD_3_5_21}))\\[3ex]
$=$&$E\bigl(\,T_{A}\,f(\,T_{A}\,)\,\bigr)$
&(cf. Lemma \ref{EWD_3_5_15}, \ref{EWD_3_5_15_BWa})).
\end{tabular}\vspace*{1ex}
\item
An application of H{\"o}lder's inequality\index{H{\"o}lder's inequality!$r$ terms}
\index{H{\"o}lder's inequality!for random variables} with $r$ terms and of Lemma \ref{EWD_3_5_15}, 
\ref{EWD_3_5_15_BWh}) gives\\[2.5ex]
\hspace*{4ex}\begin{tabular}[t]{@{}c@{\hspace*{0.8ex}}l@{}}
&$\displaystyle{n\,E\Bigl(\,\big|\,a_{I_{l_{1}}\pi_{k_{1}}(I_{l_{1}})}\,\big|\,
\big|\,a_{I_{l_{2}}\pi_{k_{2}}(I_{l_{2}})}\,\big|\,\ldots\,
\big|\,a_{I_{l_{r}}\pi_{k_{r}}(I_{l_{r}})}\,\big|\,\,\Bigr)}$\\[3.5ex]
$\leq$&$\displaystyle{n\,E\Bigl(\,\big|\,a_{I_{l_{1}}\pi_{k_{1}}(I_{l_{1}})}\,\big|^r\,\Bigr)^{1/r}\,
E\Bigl(\,\big|\,a_{I_{l_{2}}\pi_{k_{2}}(I_{l_{2}})}\,\big|^r\,\Bigr)^{1/r}\,
\ldots  
E\Bigl(\,\big|\,a_{I_{l_{r}}\pi_{k_{r}}(I_{l_{r}})}\,\big|^r\,\Bigr)^{1/r}}$\\[3.5ex]
$=$&$\displaystyle{n\,\Bigl(\,\dfrac{1}{n^2}\,\sum\limits_{i, j=1}^{n}\,|\,a_{ij}\,|^r
\,\Bigr)^{1/r\,+\,1/r\,+\,\ldots\,+\,1/r}
= \dfrac{1}{n}\,\sum\limits_{i, j=1}^{n}\,|\,a_{ij}\,|^r}$.
\end{tabular}\vspace*{1ex}
\item
$J_{1} = \pi_{5}(I_{1})$, $J_{2} = \pi_{5}(I_{2})$, $J_{1} = \pi_{4}(I_{2})$ and
$J_{2} = \pi_{4}(I_{1})$ (cf. Table \ref{EWDTab_3_5_02}).
Thus, all $\Delta T_{i}$ for $i = 1, 2, 3, 4$ can be represented as a sum of terms 
$\pm a_{I_{l}\pi_{k}(I_{l})}$ (cf. Lemma \ref{EWD_3_5_22}).\\[2.8ex] 
Consequently, the term in this part \ref{EWD_3_5_25_BWc}) can be estimated 
by a sum of terms from part \ref{EWD_3_5_25_BWb}). The constant $C$ indicates the number of these summands.
\hspace*{1ex}\hfill$\Box$
\end{enumerate}
\vspace*{3.5ex}
The following result is important for the correct determination of the Edgeworth expansions.\\[4ex]
\refstepcounter{DSBcount}
\label{EWD_3_5_26}
\textbf{\hspace*{-0.8ex}\theDSBcount\ Lemma}\\[0.8ex]
If $A = \bm\hat{A}$, then
\begin{enumerate}
\item\label{EWD_3_5_26_BWa}
$E\bigl(\,a_{I_{1}J_{1}}\,\bigr) = 0$,
\item\label{EWD_3_5_26_BWb}
$n\,E\bigl(\,a_{I_{1}J_{1}}\,\Delta T_{4}\,\bigr) = 1$,
\item\label{EWD_3_5_26_BWc}
$n\,\biggl\{\,E\bigl(\,a_{I_{1}J_{1}}\,\Delta T_{4}\,\Delta T_{3}\,\bigr) +
E\Bigl(\,a_{I_{1}J_{1}}\,\dfrac{(\Delta T_{4})^2}{2}\,\Bigr)\,\biggr\} = \dfrac{1}{2}\,E(T_{\!A}^3)$,
\item\label{EWD_3_5_26_BWd}
$n\,\biggl\{\,E\bigl(\,a_{I_{1}J_{1}}\,\Delta T_{4}\,\Delta T_{3}\,\Delta T_{2}\,\bigr) +
E\Bigl(\,a_{I_{1}J_{1}}\,\Delta T_{4}\,\dfrac{(\Delta T_{3})^2}{2}\,\Bigr)$\\[2.5ex]
\hspace*{3ex}$+\ E\Bigl(\,a_{I_{1}J_{1}}\,\dfrac{(\Delta T_{4})^2}{2}\,\Delta T_{3}\,\Bigr)
+ E\Bigl(\,a_{I_{1}J_{1}}\,\dfrac{(\Delta T_{4})^3}{6}\,\Bigr)\,\biggr\} = 
\dfrac{1}{6}\,\bigl(\,E(T_{\!A}^4) - 3\,\bigr)$.
\end{enumerate}
\vspace*{3.5ex}

\pagebreak

\textbf{Proof:}
\begin{enumerate}
\item
Since $J_{1} = \pi_{5}(I_{1})$ (cf. Table \ref{EWDTab_3_5_02}) and since $\bigl(\,I_{1}, \pi_{5}(I_{1})\,\bigr)$ is uniformly distributed on $N^{2}$ (cf. Lemma \ref{EWD_3_5_15}, \ref{EWD_3_5_15_BWh})), we obtain, 
because of $A = \bm\hat{A}$ and (\ref{EWD_3_1_01}), i.e.
$\displaystyle{\bm\hat{a}_{\boldsymbol{.}\boldsymbol{.}} = 0}\,\,$,\\[2ex]
\hspace*{4ex}$\displaystyle{E\bigl(\,a_{I_{1}J_{1}}\,\bigr) = 
\dfrac{1}{n^2}\,\sum\limits_{i, j=1}^{n}\,a_{ij} = 0}$.
\item
\hspace*{4ex}\begin{tabular}[t]{@{}c@{\hspace*{0.8ex}}l@{\hspace*{13.4ex}}r@{}}
&$E\bigl(\,a_{I_{1}J_{1}}\,\Delta T_{4}\,\bigr)$\\[2.5ex]
$=$&$E\bigl(\,a_{I_{1}J_{1}}\,T_{5}\,\bigr) - E\bigl(\,a_{I_{1}J_{1}}\,T_{4}\,\bigr)$
&(cf. (\ref{EWD_3_5_21}): $\Delta T_{4} = T_{5} - T_{4}$)\\[2.5ex]
$=$&$\displaystyle{E\bigl(\,a_{I_{1}J_{1}}\,T_{5}\,\bigr)
- E\bigl(\,a_{I_{1}J_{1}}\,\bigr)\,E\bigl(\,T_{4}\,\bigr)}$&(cf. Lemma \ref{EWD_3_5_23}, \ref{EWD_3_5_23_BWa}))\\[2.5ex]
$=$&$\displaystyle{\dfrac{1}{n}\,E\bigl(\,T_{A}^2\,\bigr)
- E\bigl(\,a_{I_{1}J_{1}}\,\bigr)\,E\bigl(\,T_{A}\,\bigr)}$&(cf. Lemma \ref{EWD_3_5_25}, \ref{EWD_3_5_25_BWa})
with $f(x) = x$)\\[2.5ex]
$=$&$\dfrac{1}{n}$&(cf. (\ref{EWD_3_1_01}): $\mu_{\bm\hat{A}} = 0$, $\sigma_{\!\bm\hat{A}}^2 = 1$).
\end{tabular}
\item
\hspace*{4ex}\begin{tabular}[t]{@{}c@{\hspace*{0.8ex}}l@{\hspace*{-8.8ex}}r@{}}
&$E\bigl(\,a_{I_{1}J_{1}}\,\Delta T_{4}\,\Delta T_{3}\,\bigr) +
E\Bigl(\,a_{I_{1}J_{1}}\,\dfrac{(\Delta T_{4})^2}{2}\,\Bigr)$\\[2.5ex]
$=$&$E\bigl(\,a_{I_{1}J_{1}}\,T_{5}\,T_{4}\,\bigr) - E\bigl(\,a_{I_{1}J_{1}}\,T_{4}^2\,\bigr)
- E\bigl(\,a_{I_{1}J_{1}}\,\Delta T_{4}\,T_{3}\,\bigr)$\\[2.5ex]
&$+\ \dfrac{1}{2}\,E\bigl(\,a_{I_{1}J_{1}}\,T_{5}^2\,\bigr) - E\bigl(\,a_{I_{1}J_{1}}\,T_{5}\,T_{4}\,\bigr)
+ \dfrac{1}{2}\,E\bigl(\,a_{I_{1}J_{1}}\,T_{4}^2\,\bigr)$\\[2.5ex]
$=$&$-\ \dfrac{1}{2}\,E\bigl(\,a_{I_{1}J_{1}}\,\bigr)\,E\bigl(\,T_{4}^2\,\bigr)
- E\bigl(\,a_{I_{1}J_{1}}\,\Delta T_{4}\,\bigr)\,E\bigl(\,T_{3}\,\bigr)$\\[2.5ex]
&$+\ \dfrac{1}{2}\,E\bigl(\,a_{I_{1}J_{1}}\,T_{5}^2\,\bigr)$
&(cf. Lemma \ref{EWD_3_5_23}, \ref{EWD_3_5_23_BWa}) and \ref{EWD_3_5_23_BWb}))\\[2.5ex]
$=$&$\dfrac{1}{2}\,E\bigl(\,a_{I_{1}J_{1}}\,T_{5}^2\,\bigr)$
&(cf. part \ref{EWD_3_5_26_BWa}) and (\ref{EWD_3_1_01}): $\mu_{\bm\hat{A}} = 0$)\\[2.5ex]
$=$&$\dfrac{1}{n}\,\dfrac{1}{2}\,E(T_{\!A}^3)$
&(cf. Lemma \ref{EWD_3_5_25}, \ref{EWD_3_5_25_BWa}) with $f(x) = x^2$).
\end{tabular}\vspace*{1.5ex}
\item
We obtain for the first two summands:\\[2ex]
\hspace*{4ex}\begin{tabular}[t]{@{}c@{\hspace*{0.8ex}}l@{}}
&$E\bigl(\,a_{I_{1}J_{1}}\,\Delta T_{4}\,\Delta T_{3}\,\Delta T_{2}\,\bigr) +
E\Bigl(\,a_{I_{1}J_{1}}\,\Delta T_{4}\,\dfrac{(\Delta T_{3})^2}{2}\,\Bigr)$\\[2.5ex]
$=$&$E\bigl(\,a_{I_{1}J_{1}}\,\Delta T_{4}\,T_{4}\,T_{3}\,\bigr) -
E\bigl(\,a_{I_{1}J_{1}}\,\Delta T_{4}\,T_{3}^2\,\bigr) - 
E\bigl(\,a_{I_{1}J_{1}}\,\Delta T_{4}\,\Delta T_{3}\,T_{2}\,\bigr)$\\[2.5ex]
&$+\ \dfrac{1}{2}\,E\bigl(\,a_{I_{1}J_{1}}\,T_{5}\,T_{4}^2\,\bigr)
- \dfrac{1}{2}\,E\bigl(\,a_{I_{1}J_{1}}\,T_{4}^3\,\bigr)
- E\bigl(\,a_{I_{1}J_{1}}\,\Delta T_{4}\,T_{4}\,T_{3}\,\bigr)
+ \dfrac{1}{2}\,E\bigl(\,a_{I_{1}J_{1}}\,\Delta T_{4}\,T_{3}^2\,\bigr)$
\end{tabular}\\[2.5ex]
\hspace*{4ex}\begin{tabular}[t]{@{}c@{\hspace*{0.8ex}}l@{\hspace*{-22.4ex}}r@{}}
$=$&$-\ \dfrac{1}{2}\,E\bigl(\,a_{I_{1}J_{1}}\,\Delta T_{4}\,\bigr)\,E\bigl(\,T_{3}^2\,\bigr)
- E\bigl(\,a_{I_{1}J_{1}}\,\Delta T_{4}\,\Delta T_{3}\,\bigr)\,E\bigl(\,T_{2}\,\bigr)$\\[2.5ex]
&$+\ \dfrac{1}{2}\,E\bigl(\,a_{I_{1}J_{1}}\,T_{5}\,T_{4}^2\,\bigr)
- \dfrac{1}{2}\,E\bigl(\,a_{I_{1}J_{1}}\,\bigr)\,E\bigl(\,T_{4}^3\,\bigr)$
&(cf. Lemma \ref{EWD_3_5_23}, \ref{EWD_3_5_23_BWa}), \ref{EWD_3_5_23_BWb}) and \ref{EWD_3_5_23_BWc}))\\[2.5ex]
$=$&$\dfrac{1}{2}\,E\bigl(\,a_{I_{1}J_{1}}\,T_{5}\,T_{4}^2\,\bigr) - \dfrac{1}{2}\,\dfrac{1}{n}$
&(cf. part \ref{EWD_3_5_26_BWa}) and \ref{EWD_3_5_26_BWb}), and
(\ref{EWD_3_1_01}): $\mu_{\bm\hat{A}} = 0$, $\sigma_{\!\bm\hat{A}}^2 = 1$).
\end{tabular}\\[4ex]
For the third summand we get:\\[2.2ex]
\hspace*{4ex}\begin{tabular}[t]{@{}c@{\hspace*{0.8ex}}l@{}}
\hspace*{1.8ex}&$E\Bigl(\,a_{I_{1}J_{1}}\,\dfrac{(\Delta T_{4})^2}{2}\,\Delta T_{3}\,\Bigr)$
\end{tabular}\\[2.5ex]
\hspace*{4ex}\begin{tabular}[t]{@{}c@{\hspace*{0.8ex}}l@{\hspace*{2.4ex}}r@{}}
$=$&$E\Bigl(\,a_{I_{1}J_{1}}\,\dfrac{(\Delta T_{4})^2}{2}\,T_{4}\,\Bigr) - 
E\Bigl(\,a_{I_{1}J_{1}}\,\dfrac{(\Delta T_{4})^2}{2}\,\Bigr)\,E\bigl(\,T_{3}\,\bigr)$
&(cf. Lemma \ref{EWD_3_5_23}, \ref{EWD_3_5_23_BWb}))\\[2.5ex]
$=$&$\dfrac{1}{2}\,E\bigl(\,a_{I_{1}J_{1}}\,T_{5}^2\,T_{4}\,\bigr)
- E\bigl(\,a_{I_{1}J_{1}}\,T_{5}\,T_{4}^2\,\bigr)
+ \dfrac{1}{2}\,E\bigl(\,a_{I_{1}J_{1}}\,T_{4}^3\,\bigr)$
&(cf. (\ref{EWD_3_1_01}): $\mu_{\bm\hat{A}} = 0$)\\[2.5ex]
$=$&$\dfrac{1}{2}\,E\bigl(\,a_{I_{1}J_{1}}\,T_{5}^2\,T_{4}\,\bigr)
- E\bigl(\,a_{I_{1}J_{1}}\,T_{5}\,T_{4}^2\,\bigr)
+ \dfrac{1}{2}\,E\bigl(\,a_{I_{1}J_{1}}\,\bigr)\,E\bigl(\,T_{4}^3\,\bigr)$
&(Lemma \ref{EWD_3_5_23}, \ref{EWD_3_5_23_BWa}))\\[2.5ex]
$=$&$\dfrac{1}{2}\,E\bigl(\,a_{I_{1}J_{1}}\,T_{5}^2\,T_{4}\,\bigr)
- E\bigl(\,a_{I_{1}J_{1}}\,T_{5}\,T_{4}^2\,\bigr)$
&(cf. part \ref{EWD_3_5_26_BWa})).
\end{tabular}\\[4ex]
For the fourth summand we get analogously to above:\\[2.2ex]
\hspace*{4ex}\begin{tabular}[t]{@{}c@{\hspace*{0.8ex}}l@{\hspace*{1.7ex}}r@{}}
&$E\Bigl(\,a_{I_{1}J_{1}}\,\dfrac{(\Delta T_{4})^3}{6}\,\Bigr)$\\[2.5ex]
$=$&$\dfrac{1}{6}\,\biggl(\,E\Bigl(\,a_{I_{1}J_{1}}\,\bigl(\,T_{5}^3 - 3\,T_{5}^2\,T_{4} + 
3\,T_{5}\,T_{4}^2\,\bigr)\,\Bigr) - E\bigl(\,a_{I_{1}J_{1}}\,\bigr)\,E\bigl(\,T_{4}^3\,\bigr)\,\biggr)$
&(Lemma \ref{EWD_3_5_23}, \ref{EWD_3_5_23_BWa}))\\[2.5ex]
$=$&$\dfrac{1}{6}\,E\bigl(\,a_{I_{1}J_{1}}\,T_{5}^3\,\bigr)
- \dfrac{1}{2}\,E\bigl(\,a_{I_{1}J_{1}}\,T_{5}^2\,T_{4}\,\bigr)
+ \dfrac{1}{2}\,E\bigl(\,a_{I_{1}J_{1}}\,T_{5}\,T_{4}^2\,\bigr)$
&(cf. part \ref{EWD_3_5_26_BWa})).
\end{tabular}\\[3.5ex]
Lastly, a summary of the three partial results gives:\\[2ex]
\hspace*{4ex}\begin{tabular}[t]{@{}l@{\hspace*{0.8ex}}c@{\hspace*{0.8ex}}l@{\hspace*{2.5ex}}r@{}}
$\biggl\{\,\ldots\ldots\ldots\ldots\,\biggr\}$&$=$&
$\dfrac{1}{6}\,E\bigl(\,a_{I_{1}J_{1}}\,T_{5}^3\,\bigr) - \dfrac{1}{2}\,\dfrac{1}{n}$\\[2.5ex]
&$=$&$\dfrac{1}{n}\,\dfrac{1}{6}\,\bigl(\,E(T_{\!A}^4) - 3\,\bigr)$
&(cf. Lemma \ref{EWD_3_5_25}, \ref{EWD_3_5_25_BWa}) with $f(x) = x^3$). $\Box$
\end{tabular}
\end{enumerate}
\vspace*{2.5ex}
Finally, we consider conditional distributions on the permutations of $N$.\\[4ex]
\refstepcounter{DSBcount}
\label{EWD_3_5_27}
\textbf{\hspace*{-0.8ex}\theDSBcount\ Lemma}
\begin{enumerate}
\item\label{EWD_3_5_27_BWa}
Let $(\,i_{1},\ldots,i_{4},\,j_{1},\ldots,j_{4}\,) \in M_{8}$, $\bm{\breve}{i} = (\,i_{i},\ldots,i_{4}\,)$
and $\theta(\bm{\breve}{i}) = \big|\,\{\,i_{1}, i_{2}, i_{3}, i_{4}\,\}\,\big|$.\\[1ex]
Also, let $\pi \in \mathscr{R}_{n} = \Bigl\{\,\nu \in \mathscr{P}_{n}\,:\,\nu(i_{1}) = j_{3},\,
\nu(i_{2}) = j_{4},\,\nu(i_{3}) = j_{1},\,\nu(i_{4}) = j_{2}\,\Bigr\}$.
\index{permutation sets!$\mathscr{R}_{n}$, $\mathscr{S}_{n}$, $\mathscr{W}_{n}$}
\index{random permutation!$\overline{i}$, $\bm{\breve}{i}$}\\[1ex] 
Then\\[2ex]
\hspace*{1.2ex}$P\Bigl(\,\pi_{3} = \pi\,\Big|\,\,\bigl(\,I_{1},\ldots, I_{4}\,\bigr) = 
(\,i_{1},\ldots,i_{4}\,),\,\bigl(\,J_{1},\ldots, J_{4}\,\bigr) = 
(\,j_{1},\ldots,j_{4}\,)\,\Bigr) = \dfrac{1}{\bigl(n - \theta(\bm{\breve}{i})\bigr)!}$,\\[2ex]
i.e. $\pi_{3}$, given $\bigl(\,I_{1},\ldots, I_{4}\,\bigr) = 
(\,i_{1},\ldots,i_{4}\,)$ and $\bigl(\,J_{1},\ldots, J_{4}\,\bigr) = 
(\,j_{1},\ldots,j_{4}\,)$, is uniformly distributed on $\mathscr{R}_{n}$.
\item\label{EWD_3_5_27_BWb}
Let $(\,i_{1},\ldots,i_{8},\,j_{1},\ldots,j_{8}\,) \in M_{16}$, $\overline{i} = (\,i_{i}, \ldots, i_{8}\,)$
and $\mu(\overline{i}) = \big|\,\{\,i_{1}, i_{2}, \ldots, i_{7}, i_{8}\,\}\,\big|$.\\[1ex]
Also, let $\pi \in \mathscr{S}_{n} = \Bigl\{\,\nu \in \mathscr{P}_{n}\,:\,\nu(i_{1}) = j_{5},\,\ldots,
\nu(i_{4}) = j_{8},\,\nu(i_{5}) = j_{1},\,\ldots,\,\nu(i_{8}) = j_{4}\,\Bigr\}$.
\index{random permutation!$\overline{i}$, $\bm{\breve}{i}$}
\index{permutation sets!$\mathscr{R}_{n}$, $\mathscr{S}_{n}$, $\mathscr{W}_{n}$}
\index{random permutation!$\overline{i}$, $\bm{\breve}{i}$}\\[1ex]
Then\\[2ex]
\hspace*{1.2ex}$P\Bigl(\,\pi_{2} = \pi\,\Big|\,\,\bigl(\,I_{1},\ldots, I_{8}\,\bigr) = 
(\,i_{1},\ldots,i_{8}\,),\,\bigl(\,J_{1},\ldots, J_{8}\,\bigr) = 
(\,j_{1},\ldots,j_{8}\,)\,\Bigr) = \dfrac{1}{\bigl(n - \mu(\overline{i})\bigr)!}$,\\[2ex]
i.e. $\pi_{2}$, given $\bigl(\,I_{1},\ldots, I_{8}\,\bigr) = 
(\,i_{1},\ldots,i_{8}\,)$ and $\bigl(\,J_{1},\ldots, J_{8}\,\bigr) = 
(\,j_{1},\ldots,j_{8}\,)$, is uniformly distributed on $\mathscr{S}_{n}$.
\item\label{EWD_3_5_27_BWc}
Let $\underline{i} = \bigl(\,i_{1},\ldots, i_{16}\,\bigr) \in M_{16}$ 
and $\underline{j} = \bigl(\,j_{1},\ldots, j_{16}\,\bigr) \in M_{16}$ 
with the properties (w1), 
\linebreak
(w2) and (w3) from (\ref{EWD_3_5_19}).\index{random permutation!properties!(w1), (w2), (w3)}
Also, let $\rho(\underline{i}) = \big|\,\{\,i_{1}, i_{2}, \ldots, i_{15}, i_{16}\,\}\,\big|$ and\\[1ex]
$\pi \in \mathscr{W}_{n} = \Bigl\{\,\nu \in \mathscr{P}_{n}\,:\,\nu(i_{1}) = j_{9},\,\ldots,
\nu(i_{8}) = j_{16},\,\nu(i_{9}) = j_{1},\,\ldots,\,\nu(i_{16}) = j_{8}\,\Bigr\}$.
\index{permutation sets!$\mathscr{R}_{n}$, $\mathscr{S}_{n}$, $\mathscr{W}_{n}$}
\index{random permutation!$\overline{i}$, $\bm{\breve}{i}$}\\[1ex]
Then\\[2ex]
\hspace*{1.2ex}$P\Bigl(\,\pi_{1} = \pi\,\Big|\,\,\bigl(\,I_{1},\ldots, I_{16}\,\bigr) = 
(\,i_{1},\ldots,i_{16}\,),\,\bigl(\,J_{1},\ldots, J_{16}\,\bigr) = 
(\,j_{1},\ldots,j_{16}\,)\,\Bigr) = \dfrac{1}{\bigl(n - \rho(\underline{i})\bigr)!}$,\\[2ex]
i.e. $\pi_{1}$, given $\bigl(\,I_{1},\ldots, I_{16}\,\bigr) = 
(\,i_{1},\ldots,i_{16}\,)$ and $\bigl(\,J_{1},\ldots, J_{16}\,\bigr) = 
(\,j_{1},\ldots,j_{16}\,)$, is uniformly distributed on $\mathscr{W}_{n}$.
\end{enumerate}
\vspace*{3.5ex}
\textbf{Proof:}
\begin{enumerate}
\item
We get for $\pi \in \mathscr{R}_{n}$:\\[2ex]
\begin{tabular}{@{}c@{\hspace*{0.8ex}}l@{}}
&$P\Bigl(\,\pi_{3} = \pi\,\Big|\,\,\bigl(\,I_{1},\ldots, I_{4}\,\bigr) = 
(\,i_{1},\ldots,i_{4}\,),\,\bigl(\,J_{1},\ldots, J_{4}\,\bigr) = 
(\,j_{1},\ldots,j_{4}\,)\,\Bigr)$\\[2.5ex]
$=$&$\dfrac{P\Bigl(\,\pi_{3} = \pi,\,\bigl(\,I_{1},\!.., I_{4}\,\bigr) = 
(\,i_{1},\!..,i_{4}\,),\,\bigl(\,\pi_{3}(I_{3}),\, \pi_{3}(I_{4}),\, \pi_{3}(I_{1}),\, \pi_{3}(I_{2})\,\bigr) = 
(\,j_{1},\, j_{2},\, j_{3},\, j_{4}\,)\Bigr)}
{P\Bigl(\,\bigl(\,I_{1},\ldots, I_{4}\,\bigr) = 
(\,i_{1},\ldots,i_{4}\,),\,\bigl(\,J_{1},\ldots, J_{4}\,\bigr) = 
(\,j_{1},\ldots,j_{4}\,)\,\Bigr)}$\\[3ex]
&\hspace*{65.4ex}(cf. $\pi_{3}$ in Table \ref{EWDTab_3_5_02})\\[1ex]
$=$&$\dfrac{P\Bigl(\,\pi_{3} = \pi,\,\bigl(\,I_{1},\ldots, I_{4}\,\bigr) = 
(\,i_{1},\ldots,i_{4}\,),\,\bigl(\,\pi(i_{3}),\, \pi(i_{4}),\, \pi(i_{1}),\, \pi(i_{2})\,\bigr) = 
(\,j_{1},\, j_{2},\, j_{3},\, j_{4}\,)\,\Bigr)}
{P\Bigl(\,\bigl(\,I_{1},\ldots, I_{4}\,\bigr) = 
(\,i_{1},\ldots,i_{4}\,),\,\bigl(\,J_{1},\ldots, J_{4}\,\bigr) = 
(\,j_{1},\ldots,j_{4}\,)\,\Bigr)}$
\end{tabular}\\
\begin{tabular}{@{}c@{\hspace*{0.8ex}}l@{}}
$=$&$\dfrac{P\Bigl(\,\pi_{3} = \pi,\,\bigl(\,I_{1},\ldots, I_{4}\,\bigr) = 
(\,i_{1},\ldots,i_{4}\,)\,\Bigr)}
{P\Bigl(\,\bigl(\,I_{1},\ldots, I_{4}\,\bigr) = 
(\,i_{1},\ldots,i_{4}\,),\,\bigl(\,J_{1},\ldots, J_{4}\,\bigr) = 
(\,j_{1},\ldots,j_{4}\,)\,\Bigr)}$\hspace*{17.1ex}($\pi \in \mathscr{R}_{n}$)\\[4ex]
$=$&$\dfrac{1}{n!}\,\dfrac{P\Bigl(\,\bigl(\,I_{1},\ldots, I_{4}\,\bigr) = 
(\,i_{1},\ldots,i_{4}\,)\,\Bigr)}
{P\Bigl(\,\bigl(\,I_{1},\ldots, I_{4}\,\bigr) = 
(\,i_{1},\ldots,i_{4}\,),\,\bigl(\,J_{1},\ldots, J_{4}\,\bigr) = 
(\,j_{1},\ldots,j_{4}\,)\,\Bigr)}$\hspace*{2.9ex}(Lemma \ref{EWD_3_5_15}, \ref{EWD_3_5_15_BWa})))
\end{tabular}\\
\begin{tabular}{@{}c@{\hspace*{0.8ex}}l@{}}
$=$&$\dfrac{1}{n!}\,\dfrac{n!}{\bigl(n - \theta(\bm{\breve}{i})\bigr)!} 
= \dfrac{1}{\bigl(n - \theta(\bm{\breve}{i})\bigr)!}$
\hspace*{19.4ex}(cf. (\ref{EWD_3_5_11}) and Corollary \ref{EWD_3_5_18}, \ref{EWD_3_5_18_BWb})).
\end{tabular}\vspace*{1ex}
\item
Completely analogous to part \ref{EWD_3_5_27_BWa}), where now Corollary \ref{EWD_3_5_18}, \ref{EWD_3_5_18_BWa}) 
is used instead of Corollary \ref{EWD_3_5_18}, \ref{EWD_3_5_18_BWb}).
\item
Also completely analogous to part \ref{EWD_3_5_27_BWa}), but now for the penultimate {''$=$''} in the 
denominator Corollary \ref{EWD_3_5_18}, \ref{EWD_3_5_18_BWd}) is used instead of 
Corollary \ref{EWD_3_5_18}, \ref{EWD_3_5_18_BWb}).
\hspace*{1ex}\hfill$\Box$
\end{enumerate}
\vspace*{2ex}

\section[A modification of Theorem 3.1.3]{A modification of Theorem \ref{EWD_3_1_03}}\label{EWD_Kap3_Sec6}

In the next two sections we need an estimate of\\[2ex]
\hspace*{12.1ex}$\displaystyle{\sup\limits_{z\, \in\, \mathbb{R}}\, 
\bigg|\,E\Bigl(\,|\,\mathscr{T}_{A}|^k\, 1_{(- \infty,\, z\,]}(\,\mathscr{T}_{A})\,\Bigr) 
- \Phi\Bigl(\,|x|^k\, 1_{(- \infty,\, z\,]}(x)\,\Bigr)\,\bigg|}$\\[2.5ex] 
for $k = 1$ and $k = 2$ to prove the two Theorems \ref{EWD_3_1_10} and \ref{EWD_3_1_13}.
The aim of this section is therefore to provide this estimate.
Since the corresponding proof can be carried out for each $k \in \mathbb{N}$ without much more work, 
we will consider this more general framework in the 
\linebreak
following.\\[2.8ex] 
The methods used in the proof are almost identical to those in Bolthausen \cite{Bolthausen1984}.
However, it was not possible to transfer the truncation technique used there (see in particular page 383, 
inequality (3.6)) to the situation here, which led in the estimate obtained to the high power $k + 4$ of the term
\mbox{\rule[-3ex]{0ex}{7.3ex}$\displaystyle{\dfrac{1}{n}\,\sum\limits_{i,\,j}|\,\bm\hat{a}_{ij}\,|^{k+4}}$}. 
For the further procedure, however, this is irrelevant. The estimate is as follows:\\[4ex]
\refstepcounter{DSBcount}
\label{EWD_3_6_01}
\textbf{\hspace*{-0.8ex}\theDSBcount\ Theorem}\index{Theorem!for linear rank statistics!further results}\\[0.8ex]
Let $k \in \mathbb{N}$ be fixed.\\[1.5ex] 
Then there exists a constant $\mathcal{K}_{7}(k) > 0$  such 
that for all matrices $A$ satisfying $\sigma_{\!A} > 0$\\[2.5ex]
\refstepcounter{DSBcount}
\label{EWD_3_6_02}
\text{\hspace*{-0.8ex}(\theDSBcount)}
\hspace*{4ex}
$\displaystyle{\sup\limits_{z\, \in\, \mathbb{R}}\, 
\bigg|\,E\Bigl(\,|\,\mathscr{T}_{A}|^k\, 1_{(- \infty,\, z\,]}(\,\mathscr{T}_{A})\,\Bigr) 
- \Phi\Bigl(\,|x|^k\, 1_{(- \infty,\, z\,]}(x)\,\Bigr)\,\bigg|}$\\[2ex]
\hspace*{13.2ex}$\displaystyle{\leq \mathcal{K}_{7}(k)\,\biggl(\, \dfrac{\beta_{A}}{n} +
\dfrac{1}{n}\sum\limits_{i,j=1}^n |\,\bm\hat{a}_{ij}\,|^{k+4} \,\biggr)}$.\\[3.5ex]
\textbf{Proof:}\\[0.8ex]
At first, it should be briefly noted that this proof has a similar structure to the proof of 
Theorem \ref{EWD_1_4_03}.\\[2.8ex]
In the following, $k \in \mathbb{N}$ is assumed to be \textbf{fixed}. Therefore the \textbf{dependence} 
of the constants $c_{1}, c_{2},\ldots$ of $k$ is \textbf{no} longer emphasized in this proof.\\[2.8ex]
Furthermore, we define for $\gamma, \lambda > 0$ and $n \in \mathbb{N} \setminus \{1\}$:\\[2ex]
\hspace*{4ex}\begin{tabular}[t]{@{}l@{\hspace*{0.8ex}}c@{\hspace*{0.8ex}}l@{}}
$\mathscr{M}(n, \gamma)$\hspace*{0.8ex}&$=$&set of $n{\times}n-$matrices $A$ satisfying $\sigma_{\!A} > 0$,
\end{tabular}\\
\hspace*{4ex}\begin{tabular}[t]{@{}l@{\hspace*{0.8ex}}c@{\hspace*{0.8ex}}l@{}}
&&$a_{ij} = \bm\hat{a}_{ij}$ for $1 \leq i,j \leq n\,\,$ and 
$\,\,\displaystyle{\beta_{A} + \sum\limits_{i,j=1}^n |\,\bm\hat{a}_{ij}\,|^{k+4} \leq \gamma}$;\\[2.8ex]
$\delta(\lambda, \gamma, n)$&$=$&$\sup\,\Bigl\{\,\Big|\,E\bigl(p_{z,k}^{\lambda}(\,T_{A})\bigr) - 
\Phi\bigl(p_{z,k}^{\lambda}\bigr)\,\Big|\,:\,z \in \mathbb{R},\,
A \in \mathscr{M}(n, \gamma)\,\Bigr\}$;\\[2.8ex]
$\delta(\gamma, n)$&$=$&$\sup\,\Bigl\{\,\Big|\,E\bigl(\,|\,T_{A}|^k\, 1_{(- \infty,\, z\,]}(\,T_{A})\,\bigr) - 
\Phi\bigl(\,|x|^k\, 1_{(- \infty,\, z\,]}(x)\,\bigr)\,\Big|\,:\,z \in \mathbb{R},\,
A \in \mathscr{M}(n, \gamma)\,\Bigr\}$.
\end{tabular}\\[2.5ex]
The functions $p_{z,k}^{\lambda}$ (cf. (\ref{EWD_2_3_01})) and 
$d_{z,k}^{\lambda}$ (cf. (\ref{EWD_2_3_04})) are defined as in section \ref{EWD_Kap2_Sec3}.
Moreover, let $\sup\,\emptyset = 0$.\\[2.8ex]
If\\[1.5ex]
\hspace*{12.1ex}$\displaystyle{c_{1} = \sup_{x \in \mathbb{R}}\, |x|^k\,\psi(x) = 
\dfrac{1}{\sqrt{2\,\pi}}\,\Bigl(\,\dfrac{k}{e}\,\Bigr)^{k/2}}$,\\[1.8ex]
then\\[2ex]
\refstepcounter{DSBcount}
\label{EWD_3_6_03}
\text{\hspace*{-0.8ex}(\theDSBcount)}
\hspace*{4ex}
$\delta(\gamma, n) \leq \delta(\lambda, \gamma, n) + \dfrac{c_{1}}{2}\,\lambda$
\hfill(analogous to (\ref{EWD_1_4_04})).\\[2.5ex]
To show the assertion $\delta(\gamma, n) \leq C(k)\,\dfrac{\gamma}{n}$, we first estimate the expression\\[2ex] 
\hspace*{12.1ex}$\Big|\,E\bigl(p_{z,k}^{\lambda}(\,T_{A})\bigr) - 
\Phi\bigl(p_{z,k}^{\lambda}\bigr)\,\Big|$\\[2.5ex]
for matrices $A \in \mathscr{M}(n, \gamma)$ using Stein's method and then apply (\ref{EWD_3_6_03}). 
For the time being, these matrices $A$ should also have the following property:\\[2ex]
\refstepcounter{DSBcount}
\label{EWD_3_6_04}
\text{\hspace*{-0.8ex}(\theDSBcount)}
\hspace*{4ex}
\begin{tabular}[t]{@{}l@{}}
$\beta_{A} \leq \epsilon_{0}\,n$ and $n \geq n_{0}$, where $0 < \epsilon_{0} \leq 1$ and $n_{0} \in \mathbb{N}$ 
are chosen so that\\[0.7ex] 
(\ref{EWD_3_3_08}) and (\ref{EWD_3_3_09}) of Corollary \ref{EWD_3_3_07} are valid. However, we require at\\[0.7ex] 
least $n_{0} \geq 6$ (so that $n_{0} - 4 \geq 2$ and thus $\mathscr{M}(n_{0} - 4, \gamma)$ is defined).
\end{tabular}\\[2.5ex]
Furthermore, let $z \in \mathbb{R}$ and $\lambda > 0$ be fixed in the following considerations, 
so that we can set $p = p_{z,k}^{\lambda}$ and $d = d_{z,k}^{\lambda}$ for simplification.\\[2.8ex] 
Using these abbreviations and $T_{5} = T_{4} + \Delta T_{4}$ and
$T_{4} = T_{3} + \Delta T_{3}$ (cf. (\ref{EWD_3_5_21})), we obtain:\\[2.3ex]
\refstepcounter{DSBcount}
\label{EWD_3_6_05}
\text{\hspace*{-0.8ex}(\theDSBcount)}
\hspace*{4ex}
\begin{tabular}[t]{@{}c@{\hspace*{0.8ex}}l@{}}
&$E\bigl(\,T_{A}\,d(\,T_{A}\,)\,\bigr)$\\[2.5ex]
$=$&$n\,E\bigl(\,a_{I_{1}J_{1}}\,d(\,T_{5}\,)\,\bigr)$
\hspace*{40ex}(cf. Lemma \ref{EWD_3_5_25}, \ref{EWD_3_5_25_BWa}))\\[1.5ex]
$=$&$\displaystyle{n\,E\bigl(\,a_{I_{1}J_{1}}\,d(\,T_{4}\,)\,\bigr) + 
n\,E\biggl(\,a_{I_{1}J_{1}}\,\Delta T_{4}\,\int\limits_{0}^{1}\,d'\bigl(\,T_{4} + t\,\Delta T_{4}\,\bigr)
\,dt\,\biggr)}$\\[4ex]
$=$&$\displaystyle{n\,E\bigl(\,a_{I_{1}J_{1}}\,d(\,T_{4}\,)\,\bigr) + 
n\,E\bigl(\,a_{I_{1}J_{1}}\,\Delta T_{4}\,d'(\,T_{3}\,)\,\bigr)}$\\[2ex] 
&$\displaystyle{+\ n\,E\biggl(\,a_{I_{1}J_{1}}\,\Delta T_{4}\,\int\limits_{0}^{1}\,
\Bigl(\,d'\bigl(\,T_{3} + \Delta T_{3} + t\,\Delta T_{4}\,\bigr)
- d'(\,T_{3}\,)\,\Bigr)\,dt\,\biggr)}$\\[2ex]
$=$&$\displaystyle{E\bigl(\,d'(T_{3}\,)\,\bigr) + n\,E\biggl(\,a_{I_{1}J_{1}}\,\Delta T_{4}\,\int\limits_{0}^{1}\,
\Bigl(\,d'\bigl(\,T_{3} + \Delta T_{3} + t\,\Delta T_{4}\,\bigr)
- d'(\,T_{3}\,)\,\Bigr)\,dt\,\biggr)}$.
\end{tabular}\\[3ex]
The last equation follows because of Lemma \ref{EWD_3_5_23}, \ref{EWD_3_5_23_BWa}) and 
\ref{EWD_3_5_23_BWb}), and Lemma \ref{EWD_3_5_26}, \ref{EWD_3_5_26_BWa}) and \ref{EWD_3_5_26_BWb}).\\[2.8ex]
From this calculation and Lemma \ref{EWD_2_3_05}, \ref{EWD_2_3_05_BWd}) and 
H{\"o}lder's inequality (\ref{EWD_0_1_05}) with $\nu = 2$ and 
$p = k+1$\index{H{\"o}lder's inequality!for finite sequences using length $\nu$} 
we further get:\\[2.5ex]
\refstepcounter{DSBcount}
\label{EWD_3_6_06}
\text{\hspace*{-0.8ex}(\theDSBcount)}
\hspace*{4ex}
\begin{tabular}[t]{@{}c@{\hspace*{0.8ex}}l@{}}
&$\Big|\,E\bigl(p(\,T_{A})\bigr) - 
\Phi\bigl(p\bigr)\,\Big|$\\[2.5ex]
$=$&$\Big|\,E\bigl(d'(\,T_{A})\bigr) - 
E\bigl(\,T_{A}\,d(\,T_{A}\,)\,\bigr)\,\Big|$
\hspace*{32.3ex}(Stein's equation\index{Stein's equation})\\[2.5ex]
$=$&$\Big|\,E\bigl(d'(\,T_{3})\bigr) - 
E\bigl(\,T_{A}\,d(\,T_{A}\,)\,\bigr)\,\Big|$\\[1.5ex]
$\leq$&$\displaystyle{n\,E\biggl(\,\big|\,a_{I_{1}J_{1}}\,\big|\,\big|\,\Delta T_{4}\,\big|\,
\int\limits_{0}^{1}\,\Big|\,d'\bigl(\,T_{3} + \Delta T_{3} + t\,\Delta T_{4}\,\bigr)
- d'(\,T_{3}\,)\,\Big|\,dt\,\biggr)}$
\hspace*{6.2ex}(cf. (\ref{EWD_3_6_05}))\\[3ex]
$\leq$&$\displaystyle{c_{2}\,n\,E\biggl(\,\big|\,a_{I_{1}J_{1}}\,\big|\,\big|\,\Delta T_{4}\,\big|\,
\int\limits_{0}^{1}\,\big|\,\Delta T_{3} +  t\,\Delta T_{4}\,\big|\,\biggl[\,1 + \big|\,T_{3}\,\big|^{k + 1}}$
\end{tabular}\\[3ex]
\hspace*{12.1ex}\begin{tabular}[t]{@{}c@{\hspace*{0.8ex}}l@{}}
&\hspace*{24.9ex}
$\displaystyle{+\ \Bigl(\,1 + \big|\,T_{3}\,\big|\,\Bigr)\,\big|\,\Delta T_{3} +  t\,\Delta T_{4}\,\big|^k}$\\[2.5ex]
&\hspace*{24.9ex}
$\displaystyle{+\ \dfrac{\big|\,T_{3}\,\big|^k}{\lambda}\,
\int\limits_{0}^{1}\,1_{(\,z,\,z\, + \lambda\,]}\bigl(\,T_{3} + s\,\Delta T_{3} + s\,t\,\Delta T_{4}\,\bigr)
\,ds\,\biggr]\,dt\,\biggl)}$
\end{tabular}\\
\hspace*{12.1ex}\begin{tabular}[t]{@{}c@{\hspace*{0.8ex}}l@{}}
$\leq$&$\displaystyle{c_{2}\Biggl\{\,n\,E\biggl(\,\big|\,a_{I_{1}J_{1}}\,\big|\,\big|\,\Delta T_{4}\,\big|\,
\Bigl(\,\big|\,\Delta T_{3}\,\big| + \big|\,\Delta T_{4}\,\big|\,\Bigr)\,\biggr)}$\\[3ex]
&\hspace*{1.7ex}$\displaystyle{+\ n\,E\biggl(\,\big|\,a_{I_{1}J_{1}}\,\big|\,\big|\,\Delta T_{4}\,\big|\,
\Bigl(\,\big|\,\Delta T_{3}\,\big| + \big|\,\Delta T_{4}\,\big|\,\Bigr)\,\big|\,T_{3}\,\big|^{k + 1}\,\biggr)}$\\[3ex]
&\hspace*{1.7ex}$\displaystyle{+\ n\,E\biggl(\,\big|\,a_{I_{1}J_{1}}\,\big|\,\big|\,\Delta T_{4}\,\big|\,
2^k\,\Bigl(\,\big|\,\Delta T_{3}\,\big|^{k+1} + \big|\,\Delta T_{4}\,\big|^{k+1}\,\Bigr)\,
\Bigl(\,1 + \big|\,T_{3}\,\big|\,\Bigr)\,\biggr)}$\\[3ex]
&\hspace*{1.7ex}$\displaystyle{+\ n\,E\biggl(\,\big|\,a_{I_{1}J_{1}}\,\big|\,\big|\,\Delta T_{4}\,\big|\,
\Bigl(\,\big|\,\Delta T_{3}\,\big| + \big|\,\Delta T_{4}\,\big|\,\Bigr)\,
\dfrac{\big|\,T_{3}\,\big|^k}{\lambda}}$\\[3ex]
&\hspace*{24.9ex}$\displaystyle{\boldsymbol{\cdot}\ \int\limits_{0}^{1}\,\int\limits_{0}^{1}\,
1_{(\,z,\,z\, + \lambda\,]}\bigl(\,T_{3} + s\,\Delta T_{3} + s\,t\,\Delta T_{4}\,\bigr)
\,ds\,dt\,\biggl)\,\Biggr\}}$\\[4ex]
$=$&$c_{2}\,\bigl(\,A_{1} + A_{2} + A_{3} + A_{4}\,\bigr)$.
\end{tabular}\\[2.7ex]
To estimate $A_{1}$ we use Lemma \ref{EWD_3_5_25}, \ref{EWD_3_5_25_BWc}) and obtain\\[2.5ex]
\refstepcounter{DSBcount}
\label{EWD_3_6_07}
\text{\hspace*{-0.8ex}(\theDSBcount)}
\hspace*{4ex}
$A_{1} \leq c_{3}\,\dfrac{\beta_{A}}{n} \leq c_{3}\,\dfrac{\gamma}{n}$.\\[2.2ex]
Because of (cf. H{\"o}lder's inequality (\ref{EWD_0_1_05}) with $\nu = 2$ and 
$p = k+1$\index{H{\"o}lder's inequality!for finite sequences using length $\nu$})\\[1.5ex]
\hspace*{12.1ex}$\big|\,T_{3}\,\big|^{k + 1} = \big|\,T_{2} + \Delta T_{2}\,\big|^{k + 1}
\leq 2^k\,\Bigl(\,\big|\,T_{2}\,\big|^{k + 1} + \big|\,\Delta T_{2}\,\big|^{k + 1}\,\Bigr)$\\[2.2ex]
and because of Lemma \ref{EWD_3_5_23}, \ref{EWD_3_5_23_BWc}) and Lemma \ref{EWD_3_5_25}, \ref{EWD_3_5_25_BWc}) 
we get analogously\\[2.2ex]
\hspace*{12.1ex}$\displaystyle{A_{2} \leq c_{4}\,\Bigl(\,\dfrac{\beta_{A}}{n}\,
E\bigl(\,\big|\,T_{A}\,\big|^{k + 1}\,\bigr)
+ \dfrac{1}{n}\sum\limits_{i,j=1}^n |\,a_{ij}\,|^{k+4}\,\Bigr)}$.\\[1.5ex]
We estimate the $(k+1)$th absolute moment of $T_{A}$. Since $E\bigl(\,T_{A}^{2}\,\bigr) = 1$, we
may assume $k \geq 2$. Using Proposition \ref{EWD_3_2_05}, we get:\\[2.2ex]
\refstepcounter{DSBcount}
\label{EWD_3_6_08}
\text{\hspace*{-0.8ex}(\theDSBcount)}
\hspace*{4ex}
\begin{tabular}[t]{@{}l@{\hspace*{0.8ex}}c@{\hspace*{0.8ex}}l@{}}
$E\bigl(\,\big|\,T_{A}\,\big|^{k + 1}\,\bigr)$&
$\leq$&$\displaystyle{\left\{
\begin{array}{l@{\hspace*{0.8ex}}ll@{}}
&\displaystyle{\Phi(x^{k+1}) + c_{5}\,
\Bigl(\,\dfrac{\beta_{A}}{n} 
+ \dfrac{1}{n}\,\sum\limits_{i,j = 1}^{n}|\,a_{ij}\,|^{k+1}\,\Bigr)}
& \hspace*{1.5ex}
\text{for}\ k \geq 3\ \text{odd},\\[2.5ex]
&\displaystyle{\Phi(x^{k+2}) + c_{6}\,
\Bigl(\,\dfrac{\beta_{A}}{n} 
+ \dfrac{1}{n}\,\sum\limits_{i,j = 1}^{n}|\,a_{ij}\,|^{k+2}\,\Bigr)}
& \hspace*{1.5ex}
\text{for}\ k \geq 2\ \text{even},
\end{array}  \right.}$\\[7ex]
&$\leq$&$c_{7} + c_{8}\,\dfrac{\gamma}{n}$.
\end{tabular}\\[2.5ex]
The first case {''$k \geq 3\ \text{odd}$''} follows directly from Proposition \ref{EWD_3_2_05}. For the 
second case {''$k \geq 2\ \text{even}$''} we use additionally
H{\"o}lder's inequality\index{H{\"o}lder's inequality!for random variables}\\[2.5ex]
\hspace*{12.1ex}$E\bigl(\,\big|\,T_{A}\,\big|^{k + 1}\,\bigr) 
\leq E\bigl(\,\big|\,T_{A}\,\big|^{k + 2}\,\bigr)^{(k+1)/(k+2)}
\leq E\bigl(\,\big|\,T_{A}\,\big|^{k + 2}\,\bigr)$\hspace*{3ex}for $k \geq 2$\\[2.5ex] 
and then apply Proposition \ref{EWD_3_2_05} to $E\bigl(\,\big|\,T_{A}\,\big|^{k + 2}\,\bigr)$.
For the second inequality of (\ref{EWD_3_6_08}), 
the inequality $|x|^{k+1} \leq |x|^3 + |x|^{k+4}$ is used for odd $k \geq 3$ 
and $|x|^{k+2} \leq |x|^3 + |x|^{k+4}$ for even $k \geq 2$.\\[2.8ex]
For the sake of completeness, we set $c_{7} = 1$, $c_{8} = 0$ for $k = 1$. Due to $\epsilon_{0} \leq 1$ and thus 
\mbox{\rule[0ex]{0ex}{3.8ex}$\dfrac{\beta_{A}}{n} \leq \epsilon_{0} \leq 1$} in (\ref{EWD_3_6_04}), it follows that\\[2ex]
\refstepcounter{DSBcount}
\label{EWD_3_6_09}
\text{\hspace*{-0.8ex}(\theDSBcount)}
\hspace*{4ex}
$\displaystyle{A_{2} \leq 
c_{4}\,\biggl(\,c_{7}\,\dfrac{\beta_{A}}{n} + c_{8}\,\dfrac{\gamma}{n}
+ \dfrac{1}{n}\sum\limits_{i,j=1}^n |\,a_{ij}\,|^{k+4}\,\biggr)
\leq c_{9}\,\dfrac{\gamma}{n}}$.\\[2ex]
Because of\\[1.5ex]
\hspace*{12.1ex}$\big|\,T_{3}\,\big| \leq \big|\,T_{2}\,\big| + \big|\,\Delta T_{2}\,\big|$
\hspace*{1.4ex}and\hspace*{1.4ex}
$E\bigl(\,\big|\,T_{2}\,\big|\,\bigr) \leq E\bigl(\,\big|\,T_{2}\,\big|^2\,\bigr)^{1/2} = \sigma_{\!A} = 1$\\[2.5ex]
and again because of Lemma \ref{EWD_3_5_23}, \ref{EWD_3_5_23_BWc}) and
Lemma \ref{EWD_3_5_25}, \ref{EWD_3_5_25_BWc}) also holds\\[2.5ex]
\refstepcounter{DSBcount}
\label{EWD_3_6_10}
\text{\hspace*{-0.8ex}(\theDSBcount)}
\hspace*{2.8ex}
$\displaystyle{A_{3} \leq c_{10}\,\biggl(\,
\dfrac{1}{n}\sum\limits_{i,j=1}^n |\,a_{ij}\,|^{k+3} + 
\dfrac{1}{n}\sum\limits_{i,j=1}^n |\,a_{ij}\,|^{k+4}\,\biggr) \leq 2\,c_{10}\,\dfrac{\gamma}{n}}$.\\[2.5ex]
For the second inequality, the inequality $|x|^{k+3} \leq |x|^3 + |x|^{k+4}$ for $k \in \mathbb{N}$ is used.\\[2.8ex]
To estimate $A_{4}$, we first have to estimate the term\\[2ex]
\refstepcounter{DSBcount}
\label{EWD_3_6_11}
\text{\hspace*{-0.8ex}(\theDSBcount)}
\hspace*{2.8ex}
\begin{tabular}[t]{@{}l@{}l@{}}
$\displaystyle{A_{5} = E\biggl(\,\big|\,T_{3}\,\big|^k\,
1_{(\,z,\,z\, +\, \lambda\,]}\bigl(\,T_{3} + s\,\Delta T_{3} + s\,t\,\Delta T_{4}\,\bigr)
\,\Big|}$&$\,\bigl(\,I_{1},\ldots, I_{4}\,\bigr) = 
(\,i_{1},\ldots,i_{4}\,),$\\[1.5ex] 
&$\,\bigl(\,J_{1},\ldots, J_{4}\,\bigr) = 
(\,j_{1},\ldots,j_{4}\,)\,\biggr)$
\end{tabular}\\[2ex]
for all  $0 \leq s, t \leq 1$ and $(\,i_{1},\ldots,i_{4},\,j_{1},\ldots,j_{4}\,) \in M_{8}$.\\[2.8ex] 
Because of Lemma \ref{EWD_3_5_27}, \ref{EWD_3_5_27_BWa}) is $T_{3}$, 
given $\bigl(\,I_{1},\ldots, I_{4}\,\bigr) = (\, i_{1},\ldots,i_{4}\,)$ and 
$\bigl(\,J_{1},\ldots, J_{4}\,\bigr) = (\,j_{1},\ldots,j_{4}\,)$, distributed as\\[2ex]
\refstepcounter{DSBcount}
\label{EWD_3_6_12}
\text{\hspace*{-0.8ex}(\theDSBcount)}
\hspace*{2.8ex}
$\displaystyle{\sum\limits_{(i,\,j)\, \in\, S} a_{ij} + T_{B}}$,\\[2.5ex]
where\\[2ex]
\refstepcounter{DSBcount}
\label{EWD_3_6_13}
\text{\hspace*{-0.8ex}(\theDSBcount)}
\hspace*{2.8ex}
$S = \Bigl\{\,(\,i_{1},\,j_{3}\,),\,(\,i_{2},\,j_{4}\,),\,(\,i_{3},\,j_{1}\,),\,(\,i_{4},\,j_{2}\,)\,\Bigr\}$\\[2.5ex]
and\\[2ex]
\refstepcounter{DSBcount}
\label{EWD_3_6_14}
\text{\hspace*{-0.8ex}(\theDSBcount)}
\hspace*{2.8ex}
\begin{tabular}[t]{@{}l@{}}
$B$ is the $(n - l) \times (n - l)-$matrix obtained from $A$ by \textbf{cancelling} the\\[0.7ex]
$l = \big|\,\{\,i_{1}, i_{2}, i_{3}, i_{4}\,\}\,\big|$ rows $i_{1},\ldots,i_{4}$ and the $l$ columns
$j_{1},\ldots,j_{4}$.
\end{tabular}\\[2.5ex]
Using $\Delta T_{3},\, \Delta T_{4} \in \sigma\bigl(\,I_{1},\ldots,I_{4},J_{1},\ldots,J_{4}\,\bigr)$
(cf. Lemma \ref{EWD_3_5_22}, \ref{EWD_3_5_22_BWc}) and \ref{EWD_3_5_22_BWd}))
and H{\"o}lder's in\-equality (\ref{EWD_0_1_05}) with $\nu = 5$ and 
$p = k$\index{H{\"o}lder's inequality!for finite sequences using length $\nu$}, we now get\\[2ex]
\hspace*{12.1ex}$A_{5} \leq 5^{k-1}\,\biggl\{\,\Bigl(\,|a_{i_{1}j_{3}}|^k + |a_{i_{2}j_{4}}|^k +
|a_{i_{3}j_{1}}|^k + |a_{i_{4}j_{2}}|^k\,\Bigr)\,\alpha(\lambda, A) + \beta(\lambda, A)\,\biggr\}$,\\[2.5ex]
where\\[2ex]
\hspace*{12.1ex}\begin{tabular}[t]{@{}l@{\hspace*{0.8ex}}c@{\hspace*{0.8ex}}l@{}}
$\alpha(\lambda, A)$&$=$&$\sup\,\Bigl\{\,P\bigl(\,z < T_{B} \leq z + \lambda\,\bigr)\,:
\,z \in \mathbb{R}, B \in M(l,A), 1 \leq l \leq 4\,\Bigr\}$,\\[2ex]
$\beta(\lambda, A)$&$=$&$\sup\,\biggl\{\,E\bigl(\,|\,T_{B}\,|^k\,1_{(\,z,\,z\, +\, \lambda\,]}(\,T_{B}\,)\,\bigr)\,:
\,z \in \mathbb{R}, B \in M(l,A), 1 \leq l \leq 4\,\Bigr\}$.
\end{tabular}\\[2.5ex]
To estimate these two quantities, we utilize $|\,\mu_{B}\,| \leq 1$ and 
$\big|\,\sigma_{\!B}^2 - 1\,\big| \leq \dfrac{1}{3}$ for $B \in M(l,A)$, 
$1 \leq l \leq 4$, 
which is valid according to (\ref{EWD_3_6_04}).\\[2.8ex] 
Thus we obtain for $B \in M(l,A)$, $1 \leq l \leq 4$:\\[2.5ex]
\refstepcounter{DSBcount}
\label{EWD_3_6_15}
\text{\hspace*{-0.8ex}(\theDSBcount)}
\hspace*{2.8ex}
\begin{tabular}[t]{@{}c@{\hspace*{0.8ex}}l@{\hspace*{1.8ex}}r@{}}
&$\sup\limits_{z\, \in\, \mathbb{R}}\,P\bigl(\,z < T_{B} \leq z + \lambda\,\bigr)$\\[3ex]
$\leq$&$\sup\limits_{z\, \in\, \mathbb{R}}\,P\bigl(\,z < \mathscr{T}_{B} \leq z + 2\,\lambda\,\bigr)$
&(since $T_{B} = \sigma_{B}\,\mathscr{T}_{B} + \mu_{B}$ and $\dfrac{1}{\sigma_{B}} \leq 2$)\\[3ex]
$\leq$&$2\,\Bigl(\,\mathcal{K}_{1}\,\dfrac{\beta_{B}}{n - l} + \dfrac{\lambda}{\sqrt{2\,\pi}}\,\Bigr)$
&(cf. Theorem \ref{EWD_3_1_03})\\[3ex]
$\leq$&$c_{11}\,\Bigl(\,\dfrac{\beta_{A}}{n} + \lambda\,\Bigr)$
&(cf. Lemma \ref{EWD_3_3_10} and $l \leq 4 < 6 \leq n\, \Rightarrow\, \dfrac{n}{n - l} \leq 3$)\\[3ex]
$\leq$&$c_{11}\,\Bigl(\,\dfrac{\gamma}{n} + \lambda\,\Bigr)$
&(since $A \in \mathscr{M}(n, \gamma)$)
\end{tabular}\\[2.5ex]
and\\[2ex]
\refstepcounter{DSBcount}
\label{EWD_3_6_16}
\text{\hspace*{-0.8ex}(\theDSBcount)}
\hspace*{2.8ex}
\begin{tabular}[t]{@{}c@{\hspace*{0.8ex}}l@{\hspace*{-0.7ex}}r@{}}
&$\sup\limits_{z\, \in\, \mathbb{R}}\,E\bigl(\,|\,T_{B}\,|^k\,1_{(\,z,\,z\, +\, \lambda\,]}(\,T_{B}\,)\,\bigr)$\\[3ex]
$\leq$&$\sup\limits_{z\, \in\, \mathbb{R}}\,E\bigl(\,|\,T_{B}\,|^k\,1_{(\,z,\,z + 2\,\lambda\,]}
(\,\mathscr{T}_{B}\,)\,\bigr)$
&(since $T_{B} = \sigma_{B}\,\mathscr{T}_{B} + \mu_{B}$ and $\dfrac{1}{\sigma_{B}} \leq 2$)\\[3ex]
$\leq$&$2^{k-1}\,\Bigl(\,\dfrac{3}{2}\,\Bigr)^k\,
\sup\limits_{z\, \in\, \mathbb{R}}\,E\bigl(\,|\,\mathscr{T}_{B}\,|^k\,1_{(\,z,\,z\, +\, 2\,\lambda\,]}
(\,\mathscr{T}_{B}\,)\,\bigr)$
&($T_{B} = \sigma_{B}\,\mathscr{T}_{B} + \mu_{B}$ and $\sigma_{B} \leq \dfrac{3}{2}$)
\end{tabular}\\
\hspace*{12.1ex}\begin{tabular}[t]{@{}c@{\hspace*{0.8ex}}l@{\hspace*{8.1ex}}r@{}}
&$+\ 2^{k-1}\,|\,\mu_{B}\,|^k\,\sup\limits_{z\, \in\, \mathbb{R}}\,P\bigl(\,z < \mathscr{T}_{B} 
\leq z + 2\,\lambda\,\bigr)$
&(and (\ref{EWD_0_1_05}) with $\nu = 2$, 
$p = k$)\index{H{\"o}lder's inequality!for finite sequences using length $\nu$}\\[3ex]
$\leq$&$\displaystyle{3^k\,\delta\Bigl(\,\beta_{B} + 
\sum\limits_{i \not= i_{m},\,j \not= j_{m}}|\,\bm\hat{b}_{ij}\,|^{k+4}\,,\, n-l\,\Bigr)}$\\[4ex]
&$+\ 3^k\,\lambda\,\sup\limits_{x\, \in\, \mathbb{R}}\,|x|^k\,\psi(x) + 
2^{k-1}\,c_{11}\,\Bigl(\,\dfrac{\gamma}{n} + \lambda\,\Bigr)$
&(cf. (\ref{EWD_3_6_15}) and $|\,\mu_{B}\,| \leq 1$)\\[4ex]
$\leq$&$3^k\,\delta\bigl(\,c_{12}\,\gamma,\, n - l\,\bigr) + c_{13}\,\Bigl(\,\dfrac{\gamma}{n} + \lambda\,\Bigr)$
&(cf. Lemma \ref{EWD_3_3_10}).
\end{tabular}\\[2.5ex]
It follows that\\[3ex]
\hspace*{4ex}$A_{5} \leq c_{14}\,\biggl\{\,\Bigl(\,|a_{i_{1}j_{3}}|^k + |a_{i_{2}j_{4}}|^k +
|a_{i_{3}j_{1}}|^k + |a_{i_{4}j_{2}}|^k + 1\,\Bigr)\,\Bigl(\,\dfrac{\gamma}{n} + \lambda\,\Bigr)
+ \max\limits_{1\, \leq\, l\, \leq\, 4}\,\delta\bigl(\,c_{12}\,\gamma,\, n - l\,\bigr)\,\biggr\}$\\[3ex]
and thus by reusing the inequality $|x|^{k+3} \leq |x|^3 + |x|^{k+4}$ for $k \in \mathbb{N}$\\[3ex]
\refstepcounter{DSBcount}
\label{EWD_3_6_17}
\text{\hspace*{-0.8ex}(\theDSBcount)}
\hspace*{2.8ex}
$A_{4} \leq c_{15}\,\dfrac{\gamma}{n}\,\biggl[\,1 + \dfrac{1}{\lambda}\,\dfrac{\gamma}{n}
+ \dfrac{1}{\lambda}\,\max\limits_{1\, \leq\, l\, \leq\, 4}\,\delta\bigl(\,c_{12}\,\gamma,\, n - l\,\bigr)\,\biggr]$.\\[3ex]
Overall, (\ref{EWD_3_6_07}), (\ref{EWD_3_6_09}), (\ref{EWD_3_6_10}) and 
(\ref{EWD_3_6_17}) yield the existence of a constant $c_{16}$ with\\[3ex]
\refstepcounter{DSBcount}
\label{EWD_3_6_18}
\text{\hspace*{-0.8ex}(\theDSBcount)}
\hspace*{2.8ex}
$\sup\limits_{z\, \in\, \mathbb{R}}\,\Big|\,E\bigl(p_{z,k}^{\lambda}(\,T_{A})\bigr) - 
\Phi\bigl(p_{z,k}^{\lambda}\bigr)\,\Big| \leq
c_{16}\,\dfrac{\gamma}{n}\,\biggl[\,1 + \dfrac{1}{\lambda}\,\dfrac{\gamma}{n}
+ \dfrac{1}{\lambda}\,\max\limits_{1\, \leq\, l\, \leq\, 4}\,\delta\bigl(\,c_{12}\,\gamma,\, n - l\,\bigr)\,\biggr]$\\[3ex]
for all matrices $A \in \mathscr{M}(n, \gamma)$, which have the property (\ref{EWD_3_6_04}).
For the remaining matrices $A \in \mathscr{M}(n, \gamma)$, on the other hand, we can show that\\[2.5ex] 
\refstepcounter{DSBcount}
\label{EWD_3_6_19}
\text{\hspace*{-0.8ex}(\theDSBcount)}
\hspace*{2.8ex}
\begin{tabular}[t]{@{}c@{\hspace*{0.8ex}}l@{}}
&$\sup\limits_{z\, \in\, \mathbb{R}}\,\Big|\,E\bigl(\eta(\,T_{A})\bigr) - 
\Phi\bigl(\eta\bigr)\Big|$,\hspace*{3ex} 
where $\eta(x) = p_{z,k}^{\lambda}(x)$ or $\eta(x) = |x|^k\, 1_{(- \infty,\, z\,]}(x)$\\[3ex]
$\leq$&$E\bigl(\,|\,T_{A}\,|^k\,\bigr) + \Phi\bigl(\,|x|^k\,\bigr)$\\[2.5ex]
$\leq$&$c_{17} + c_{18}\,\dfrac{\gamma}{n}$
\hspace*{32.1ex}(analogous to (\ref{EWD_3_6_08}) and cf. (\ref{EWD_2_1_07})).
\end{tabular}\\[2.5ex] 
In the case of $\beta_{A} > \epsilon_{0}\,n$, we further obtain\\[3ex]
\refstepcounter{DSBcount}
\label{EWD_3_6_20}
\text{\hspace*{-0.8ex}(\theDSBcount)}
\hspace*{2.8ex}
$\leq \dfrac{c_{17}}{\epsilon_{0}}\,\dfrac{\beta_{A}}{n} + c_{18}\,\dfrac{\gamma}{n}
\leq \Bigl(\,\dfrac{c_{17}}{\epsilon_{0}} + c_{18}\,\Bigr)\,\dfrac{\gamma}{n}$.\\[2.5ex]
Moreover, in the case $n < n_{0}$ we get because of 
$\dfrac{\beta_{A}}{n} \geq \dfrac{1}{\sqrt{8\,n}}$ (cf. Lemma \ref{EWD_3_1_18}, \ref{EWD_3_1_18_BWa})):\\[2.5ex]
\refstepcounter{DSBcount}
\label{EWD_3_6_21}
\text{\hspace*{-0.8ex}(\theDSBcount)}
\hspace*{2.8ex}
$\leq c_{17}\,\sqrt{8\,n_{0}}\,\dfrac{\beta_{A}}{n_{0}} + c_{18}\,\dfrac{\gamma}{n}
\leq c_{17}\,\sqrt{8\,n_{0}}\,\dfrac{\beta_{A}}{n} + c_{18}\,\dfrac{\gamma}{n}
\leq \Bigl(\,c_{17}\,\sqrt{8\,n_{0}} + c_{18}\,\Bigr)\,\dfrac{\gamma}{n}$.\\[3ex]
Combining the estimates (\ref{EWD_3_6_18}), (\ref{EWD_3_6_19}), (\ref{EWD_3_6_20}), (\ref{EWD_3_6_21}) 
with (\ref{EWD_3_6_03}), we obtain a constant $c_{19}$ such that\\[3ex]  
\hspace*{12.1ex}$\delta(\gamma, n) \leq c_{19}\,\biggl\{\,\dfrac{\gamma}{n}\,\biggl[\,1 
+ \dfrac{1}{\lambda}\,\dfrac{\gamma}{n} + \dfrac{1}{\lambda}\,\max\limits_{1\, \leq\, l\, \leq\, 4}\,
\delta\bigl(\,c_{12}\,\gamma,\, n - l\,\bigr)\,\biggr] + \lambda\,\biggr\}$
\hspace*{3ex}for all $n \geq 6$.\\[2.5ex]
If we now choose $\lambda = 4\,c_{19}\,c_{12}\,\dfrac{\gamma}{n}$, we get\\[2.5ex]
\hspace*{12.1ex}$\delta(\gamma, n) \leq c_{20}\,\dfrac{\gamma}{n} + \dfrac{1}{4\,c_{12}}\,
\max\limits_{1\, \leq\, l\, \leq\, 4}\,\delta\bigl(\,c_{12}\,\gamma,\, n - l\,\bigr)$
\hspace*{3ex}for all $n \geq 6$.\\[2.5ex]
Since the inequality $n \leq 3\,(n - l)$ is valid for $n \geq 6$ and $l \leq 4$, we 
conclude\\[2.5ex]
\hspace*{12.1ex}$\dfrac{n}{\gamma}\,\delta(\gamma, n) \leq c_{20} + \dfrac{3}{4}\,
\max\limits_{1\, \leq\, l\, \leq\, 4}\,
\dfrac{n - l}{c_{12}\,\gamma}\,\delta\bigl(\,c_{12}\,\gamma,\, n - l\,\bigr)$
\hspace*{3ex}for all $n \geq 6$\\[2.5ex]
and therefore\\[2.5ex]
\refstepcounter{DSBcount}
\label{EWD_3_6_22}
\text{\hspace*{-0.8ex}(\theDSBcount)}
\hspace*{2.8ex}
$\sup\limits_{\gamma > 0}\,\dfrac{n}{\gamma}\,\delta(\gamma, n) \leq c_{20} + \dfrac{3}{4}\,
\max\limits_{1\, \leq\, l\, \leq\, 4}\,\sup\limits_{\gamma > 0}\,
\dfrac{n - l}{\gamma}\,\delta\bigl(\,\gamma,\, n - l\,\bigr)$
\hspace*{3ex}for all $n \geq 6$.\\[3ex]
A repeated application of this inequality (\ref{EWD_3_6_22}) then gives\\[2.5ex]
\refstepcounter{DSBcount}
\label{EWD_3_6_23}
\text{\hspace*{-0.8ex}(\theDSBcount)}
\hspace*{2.8ex}
$\sup\limits_{\gamma > 0}\,\dfrac{n}{\gamma}\,\delta(\gamma, n) \leq 
\displaystyle{\sum\limits_{s = 0}^{n - 6}\Bigl(\,\dfrac{3}{4}\,\Bigr)^s\,c_{20} +
\dfrac{3}{4}\,\max\limits_{2\, \leq\, r\, \leq\, 5}\,\sup\limits_{\gamma > 0}\,
\dfrac{r}{\gamma}\,\delta\bigl(\,\gamma,\, r\,\bigr)}$
\hspace*{3ex}for all $n \geq 6$.\\[3ex]
If we now use the summation formula of the geometric series\index{geometric series}, we further get\\[2.5ex]
\hspace*{12.1ex}\begin{tabular}{@{}l@{\hspace*{0.8ex}}c@{\hspace*{0.8ex}}l@{}}
$\sup\limits_{\gamma > 0}\,\dfrac{n}{\gamma}\,\delta(\gamma, n)$&$\leq$&
$4\,c_{20} + \dfrac{3}{4}\,\max\limits_{2\, \leq\, r\, \leq\, 5}\,\sup\limits_{\gamma > 0}\,
\dfrac{r}{\gamma}\,\delta\bigl(\,\gamma,\, r\,\bigr)$\\[3.5ex]
&$\leq$&$4\,c_{20} + \Bigl(\,3\,\sqrt{3}\,c_{17} + \dfrac{3}{4}\,c_{18}\,\Bigr) < \infty$
\hspace*{3ex}for all $n \geq 2$.
\end{tabular}\\[3.5ex]
For the last inequality we proceeded as in the calculations in (\ref{EWD_3_6_19}) and (\ref{EWD_3_6_21}) 
with $n_{0} = 6$ and $\eta(x) = |x|^k\, 1_{(- \infty,\, z\,]}(x)$.
\hspace*{1ex}\hfill$\Box$\\[4ex]
\refstepcounter{DSBcount}
\label{EWD_3_6_24}
\textbf{\hspace*{-0.8ex}\theDSBcount\ Remark}\\[0.8ex]
For $k = 1$, a slightly different proof of Theorem \ref{EWD_3_6_01} 
with the remainder term $c\,\beta_{A}/n$ instead of $c\,\bigl(\,\beta_{A}/n + \eta_{A}/n\,\bigr)$
is sketched in \cite{10.1214/aos/1176347258}, Proposition 5.7.\\[4ex] 
The following analogous result to Theorem \ref{EWD_3_6_01}, in which the terms $|\,\mathscr{T}_{A}|^k$ and $|x|^k$
are replaced by the terms $\mathscr{T}_{A}^k$ and $x^k$, 
can be derived from this theorem without much effort.\\[4ex]
\refstepcounter{DSBcount}
\label{EWD_3_6_25}
\textbf{\hspace*{-0.8ex}\theDSBcount\ Corollary}\\[0.8ex]
Let $k \in \mathbb{N}$ be fixed.\\[1.5ex] 
Then there exists a constant $\mathcal{K}_{8}(k) > 0$ such 
that for all matrices $A$ satisfying $\sigma_{\!A} > 0$\\[2.5ex]
\refstepcounter{DSBcount}
\label{EWD_3_6_26}
\text{\hspace*{-0.8ex}(\theDSBcount)}
\hspace*{4ex}
$\displaystyle{\sup\limits_{z\, \in\, \mathbb{R}}\, 
\bigg|\,E\Bigl(\,\mathscr{T}_{A}^k\, 1_{(- \infty,\, z\,]}(\,\mathscr{T}_{A})\,\Bigr) 
- \Phi\Bigl(\,x^k\, 1_{(- \infty,\, z\,]}(x)\,\Bigr)\,\bigg|}$\\[2ex]
\hspace*{13.2ex}$\displaystyle{\leq \mathcal{K}_{8}(k)\,\biggl(\, \dfrac{\beta_{A}}{n} +
\dfrac{1}{n}\sum\limits_{i,j=1}^n |\,\bm\hat{a}_{ij}\,|^{k+4} \,\biggr)}$.\\[3.5ex]
\textbf{Proof:}\\[0.8ex]
We consider the two cases $z \leq 0$ and $z \geq 0$:\\[2ex]
\begin{tabular}{@{}l@{\hspace*{2ex}}c@{\hspace*{0.8ex}}l@{}}
$z \leq 0$:&\multicolumn{2}{@{}l@{}}{$E\Bigl(\,\mathscr{T}_{A}^k\, 1_{(- \infty,\, z\,]}(\,\mathscr{T}_{A})\,\Bigr)
= (-1)^k\,E\Bigl(\,|\,\mathscr{T}_{A}|^k\, 1_{(- \infty,\, z\,]}(\,\mathscr{T}_{A})\,\Bigr)$}.\\[2ex]
$z \geq 0$:&&$E\Bigl(\,\mathscr{T}_{A}^k\, 1_{(- \infty,\, z\,]}(\,\mathscr{T}_{A})\,\Bigr)$\\[2ex]
&$=$&$E\Bigl(\,\mathscr{T}_{A}^k\, 1_{(- \infty,\, 0\,]}(\,\mathscr{T}_{A})\,\Bigr)
+ E\Bigl(\,\mathscr{T}_{A}^k\, 1_{(\,0,\, z\,]}(\,\mathscr{T}_{A})\,\Bigr)$\\[2ex]
&$=$&$E\Bigl(\,|\,\mathscr{T}_{A}|^k\, 1_{(- \infty,\, z\,]}(\,\mathscr{T}_{A})\,\Bigr)
+ \bigl(\,(-1)^k - 1\,\bigr)\,E\Bigl(\,|\,\mathscr{T}_{A}|^k\, 1_{(- \infty,\, 0\,]}(\,\mathscr{T}_{A})\,\Bigr)$.
\end{tabular}\\[2.5ex]
Since we can proceed completely analogously with the term $\Phi\Bigl(\,x^k\, 1_{(- \infty,\, z\,]}(x)\,\Bigr)$, 
the assertion follows from Theorem \ref{EWD_3_6_01}.
\hspace*{1ex}\hfill$\Box$\\[4ex]
Finally, a slightly more general version of Corollary \ref{EWD_3_6_25} and Theorem \ref{EWD_3_1_03} 
is also necessary for the further procedure.\\[4ex]
\refstepcounter{DSBcount}
\label{EWD_3_6_27}
\textbf{\hspace*{-0.8ex}\theDSBcount\ Corollary}\\[0.8ex]
Let $k \in \mathbb{N}_{0}$ be fixed, and let $\mathcal{H}$ be defined 
according to (\ref{EWD_0_1_01})\index{function!$\mathcal{H}$}. Then\\[2.5ex]
\refstepcounter{DSBcount}
\label{EWD_3_6_28}
\text{(\theDSBcount)}\\[-2.5ex]
\hspace*{1ex}\hfill$\displaystyle{\sup\limits_{h\, \in\, \mathcal{H}}\, 
\bigg|\,E\Bigl(\,\mathscr{T}_{A}^k\, h(\,\mathscr{T}_{A})\,\Bigr) 
- \Phi\Bigl(\,x^k\, h(x)\,\Bigr)\,\bigg|
\, \leq\, 
\left\{
\begin{array}{l@{\hspace*{3ex}}l@{}}
\mathcal{K}_{1}\,\dfrac{\beta_{A}}{n}
&\text{for}\ k = 0,\\[3.5ex]
\mathcal{K}_{8}(k)\,\biggl(\, \dfrac{\beta_{A}}{n} +
\displaystyle{\dfrac{1}{n}\sum\limits_{i,j=1}^n |\,\bm\hat{a}_{ij}\,|^{k+4}} \,\biggr)
&\text{for}\ k \geq 1.
\end{array}  \right.}$\\[4ex] 
\textbf{Proof:}\\[0.8ex]
Proceed as in the proof of Corollary \ref{EWD_1_4_08}. 
\hspace*{1ex}\hfill$\Box$
                           
\section[Edgeworth expansions of first order]{Edgeworth expansions of first order}\label{EWD_Kap3_Sec7}

In this section we prove the Theorem \ref{EWD_3_1_10} and a similar version with {''first moment''}. 
The latter is needed in the next section to prove the Theorem \ref{EWD_3_1_13}. 
The main work to achieve these two goals is done by proving the following theorem.\\[3.5ex]
\refstepcounter{DSBcount}
\label{EWD_3_7_01}
\textbf{\hspace*{-0.8ex}\theDSBcount\ Theorem}\index{Theorem!for linear rank statistics!further results}\\[0.8ex]
Let $m \in \{\,0,\,1\,\}$, $G_{\!A} \in \{\,D_{\!A},\,E_{\!A}\,\}$\index{matrix!moment!$G_{\hspace*{-0.3ex}A}$}\index{Edgeworth expansion!$e_{1,A}^1$, $e_{2,A}^m$} be fixed and\\[2.2ex]
\hspace*{12.1ex}$\displaystyle{e_{1,A}^m(z) = \int\limits_{- \infty}^{z} x^m\,e_{1,A}'(x)\,dx}$
\hfill(cf. Lemma \ref{EWD_3_1_18}, \ref{EWD_3_1_18_BWj}) for $m = 1$).\\[2.2ex]
Then there exist constants $\mathfrak{K}_{1}(m) > 0$ and $\mathfrak{K}_{2}(m) > 0$ such that
for all matrices $A$ satisfying $\sigma_{\!A} > 0$ and the conditions\\[2.2ex]
\refstepcounter{DSBcount}
\label{EWD_3_7_02}
\text{\hspace*{-0.8ex}(\theDSBcount)}
\hspace*{4ex}
\mbox{\rule[0ex]{0ex}{3.8ex}$\sigma_{\!A'}^2 \geq \dfrac{2}{3}$},\\[2ex]
\refstepcounter{DSBcount}
\label{EWD_3_7_03}
\text{\hspace*{-0.8ex}(\theDSBcount)}
\hspace*{4ex}
\begin{tabular}[t]{@{}l@{}}
there exists a constant $\mathfrak{L}(m) > 0$ such that\\[1.5ex]
\hspace*{4ex}
$\bigl(\,|z|^m + 1\,\bigr)\,\big|\,\Delta_{y}^{2}F_{\SBTB}(z)\,\big| 
\leq \mathfrak{L}(m)\,\bigl(\,G_{\!A}^2 + y^2\,\bigr)$\\[1.8ex]
for all $z \in \mathbb{R}$, $0 \leq y \leq \dfrac{G_{\!A}}{\sigma_{\!A'}}$ 
and $\BTB \in N(8,\bm\bar{A})$,
\end{tabular}\\[0.5ex]
we have\\[2ex]
\refstepcounter{DSBcount}
\label{EWD_3_7_04}
\text{\hspace*{-0.8ex}(\theDSBcount)}
\hspace*{4ex}
$\displaystyle{\sup\limits_{z\, \in\, \mathbb{R}}\, 
\bigg|\,E\Bigl(\,T_{\bm\bar{A}\,}^m\, 1_{(- \infty,\, z\,]}(\,T_{\bm\bar{A}}\,)\,\Bigr) 
- e_{1,\bm\bar{A}}^{m}(z)\,\bigg| \leq \Bigl(\,\mathfrak{K}_{1}(m)\,\mathfrak{L}(m) 
+ \mathfrak{K}_{2}(m)\,\Bigr)\,G_{\!A}^2}$.\\[5.8ex]
The case $m = 0$ and $G_{\!A} = D_{\!A}$ leads to the proof of Theorem \ref{EWD_3_1_10} 
by the considerations in the following subsection \ref{EWD_Kap3_Sec7_1}.\\[2ex]
The case $m = 1$ and $G_{\!A} = E_{\!A}$ is important for the next section. There we also need
a slightly more general version of Theorem \ref{EWD_3_7_01}.\\[3.5ex]
\refstepcounter{DSBcount}
\label{EWD_3_7_05}
\textbf{\hspace*{-0.8ex}\theDSBcount\ Corollary}\\[0.8ex]
Suppose that the conditions and notations of the previous Theorem \ref{EWD_3_7_01} apply.
Furthermore, let $\mathcal{H}$\index{function!$\mathcal{H}$} be defined according to (\ref{EWD_0_1_01}).
Then\\[3ex]
\refstepcounter{DSBcount}
\label{EWD_3_7_06}
\text{\hspace*{-0.8ex}(\theDSBcount)}
\hspace*{4ex}
$\displaystyle{\sup\limits_{h\, \in\, \mathcal{H}}\,
\bigg|\,E\Bigl(\,T_{\bm\bar{A}\,}^m\, h(\,T_{\bm\bar{A}}\,)\,\Bigr) 
- \int\limits_{\mathbb{R}} x^m\,h(x)\,e_{1,\bm\bar{A}}'(x)\,dx\,\bigg| \leq \Bigl(\,\mathfrak{K}_{1}(m)\,\mathfrak{L}(m) 
+ \mathfrak{K}_{2}(m)\,\Bigr)\,G_{\!A}^2}$.\\[3.5ex]
\textbf{Proof:}\\[0.8ex]
Proceed as in the proof of Corollary \ref{EWD_1_4_08}. 
\hspace*{1ex}\hfill$\Box$\vspace*{1ex}

\subsection[Proof of Theorem \ref{EWD_3_1_10} using Theorem \ref{EWD_3_7_01}]
{Proof of Theorem \ref{EWD_3_1_10} using Theorem \ref{EWD_3_7_01}}\label{EWD_Kap3_Sec7_1}

According to Corollary \ref{EWD_3_4_02}, there exists a $0 < \epsilon_{0} \leq 1$ 
such that the inequalities \mbox{\rule[0ex]{0ex}{3.8ex}$\dfrac{2}{3} \leq \sigma_{\!A'}^2 \leq \dfrac{4}{3}$}
hold for all matrices $A$ satisfying $\sigma_{\!A} > 0$ and $\beta_{A} \leq \epsilon_{0}\,n$.\\[2.8ex]
We thus prove the Theorem \ref{EWD_3_1_10} by considering the following two cases:\\[4ex]
\textbf{1. Case:} $\beta_{A} \leq \epsilon_{0}\,n$ and $n \geq n_{0} = 10$.\\[2.8ex]
First, we derive condition (\ref{EWD_3_7_03}) for $m = 0$ and $G_{\!A} = D_{\!A}$
from condition (\ref{EWD_3_1_11}). 
To do this, we use Proposition \ref{EWD_3_4_17} (with $r = 2$ and $k = 4$) and obtain:\\[2.5ex]
\hspace*{12.1ex}\begin{tabular}{@{}l@{\hspace*{0.8ex}}c@{\hspace*{0.8ex}}l@{}}
$\displaystyle{\big|\big|\,\Delta_{y}^{2}F_{\SBTB}\,\big|\big|}$&$\leq$&$\displaystyle{
\big|\big|\,\Delta_{\sigma_{\!A'}y}^{2}F_{B}\,\big|\big| + c_{1}\,D_{\!A}^2}$\\[2.8ex]
&$\leq$&$\Bigl(\,\dfrac{4}{3}\,\mathcal{C}_{1} + c_{1}\,\Bigr)\,\bigl(\,D_{\!A}^2 + y^2\,\bigr)$
\hspace*{4ex}for all $0 \leq y \leq \dfrac{D_{\!A}}{\sigma_{\!A'}}$ 
and $\BTB \in N(8,\bm\bar{A})$.
\end{tabular}\\[2.2ex]
So because of $|z|^0 + 1 = 2$, we get $\mathfrak{L}(0) = 2\,\Bigl(\,\dfrac{4}{3}\,\mathcal{C}_{1} + c_{1}\,\Bigr)$.\\[3ex]
Furthermore, we obtain from (\ref{EWD_3_4_14}), 
applied to the functions $h = 1_{(\,- \infty,\,z\,]}, z \in \mathbb{R}$, that\\[2.5ex]
\hspace*{12.1ex}$||\,\mathscr{F}_{\!A} - e_{1,A}\,|| \leq 
||\,\mathscr{F}_{\bm\bar{A}} - e_{1,\bm\bar{A}}\,|| + c_{2}\,D_{\!A}^2$.\\[2.5ex]
Therefore, the inequality (\ref{EWD_3_1_12}) of Theorem \ref{EWD_3_1_10} follows directly from
the inequality (\ref{EWD_3_7_04}) with $m = 0$, $G_{\!A} = D_{\!A}$ and 
$\mathfrak{L}(0) = 2\,\Bigl(\,\dfrac{4}{3}\,\mathcal{C}_{1} + c_{1}\,\Bigr)$.\\[4ex]
\textbf{2. Case:} $\beta_{A} > \epsilon_{0}\,n$ or $n < n_{0} = 10$.\\[2.8ex]
Because of Lemma \ref{EWD_3_1_18}, \ref{EWD_3_1_18_BWb}), we have\\[2ex]
\refstepcounter{DSBcount}
\label{EWD_3_7_07}
\text{\hspace*{-0.8ex}(\theDSBcount)}
\hspace*{4ex}
$\beta_{A} > \epsilon_{0}\,n$
\hspace*{2ex}$\Rightarrow$\hspace*{2ex}
$1 \leq \dfrac{1}{\epsilon_{0}}\,\dfrac{\beta_{A}}{n} \leq 
\dfrac{1}{\epsilon_{0}}\,D_{\!A}$\\[2.5ex]
and because of Lemma \ref{EWD_3_1_18}, \ref{EWD_3_1_18_BWa}) also\\[2ex]
\refstepcounter{DSBcount}
\label{EWD_3_7_08}
\text{\hspace*{-0.8ex}(\theDSBcount)}
\hspace*{4ex}
$n < n_{0}$
\hspace*{2ex}$\Rightarrow$\hspace*{2ex}
$1 \leq \dfrac{n_{0}}{n} \leq \dfrac{n_{0}}{n}\,4\,\delta_{\!A} = 4\,n_{0}\,D_{\!A}^2$.\\[2ex]
Thus, it follows that\\[2ex]
\refstepcounter{DSBcount}
\label{EWD_3_7_09}
\text{\hspace*{-0.8ex}(\theDSBcount)}
\hspace*{4ex}
\begin{tabular}[t]{@{}l@{\hspace*{0.8ex}}c@{\hspace*{0.8ex}}l@{\hspace*{-4.4ex}}r@{}}
$||\,\mathscr{F}_{\!A} - e_{1,A}\,||$&$\leq$&$1 + \dfrac{1}{15}\,\dfrac{\beta_{A}}{n}$
&(cf. Lemma \ref{EWD_3_1_18}, \ref{EWD_3_1_18_BWf}))\\[2.5ex] 
&$\leq$&$1 + \dfrac{1}{15}\,D_{\!A}$
&(cf. Lemma \ref{EWD_3_1_18}, \ref{EWD_3_1_18_BWb}))\\[2.5ex] 
&$\leq$&$\biggl[\,\max\,\Bigl\{\,\dfrac{1}{\epsilon_{0}^2},\,4\,n_{0}\,\Bigr\} 
+ \dfrac{1}{15}\,\max\,\Bigl\{\,\dfrac{1}{\epsilon_{0}},\,\sqrt{\rule[0ex]{0ex}{2.2ex}{4\,n_{0}}}\,\Bigr\}\,
\biggr]\,D_{\!A}^2$.
\end{tabular}\\[3ex]
A sufficiently large $\mathcal{K}_{3}$ in Theorem \ref{EWD_3_1_10} therefore 
ensures the validity of this theorem in the 2nd case as well.\hspace*{1ex}\hfill$\Box$\\[2.5ex]
We need the rest of this section for the

\subsection[Proof of Theorem \ref{EWD_3_7_01}]{Proof of Theorem \ref{EWD_3_7_01}}\label{EWD_Kap3_Sec7_2}

Due to Corollary \ref{EWD_3_4_02}, Lemma \ref{EWD_3_4_05} (for $\beta_{\bm\bar{A}}
\leq 125\,\beta_{A}$) and Corollary \ref{EWD_3_3_07} 
there exists a $0 < \epsilon_{0} \leq 1$ and an $n_{0} \geq 10$
such that for all matrices $A$ satisfying
$\sigma_{\!A} > 0$, $\beta_{A} \leq \epsilon_{0}\,n$ and $n \geq n_{0}$ holds:\\[2.5ex]
\refstepcounter{DSBcount}
\label{EWD_3_7_10}
\text{\hspace*{-0.8ex}(\theDSBcount)}
\hspace*{4ex}
\begin{tabular}{@{}l@{}}
$\displaystyle{
\left\{
\begin{array}{l@{\hspace*{3ex}}l@{}}
\big|\,\sigma_{\!A'}^2 - 1\,\big| \leq \dfrac{1}{3},\\[2.5ex]
|\bm\bar{a}_{ij}| \leq 1&\text{for}\ 1 \leq i, j \leq n,\\[2.5ex]
|\,\mu_{\SBTB}\,| \leq 1&\text{for}\ \BTB \in N(8,\bm\bar{A}),\\[2.8ex]
\big|\,\sigma_{\SBTB}^2 - 1\,\big| \leq \dfrac{1}{3}&\text{for}\ \BTB \in N(8,\bm\bar{A}).
\end{array}  \right.}$
\end{tabular}\\[2.5ex]
A suitable $\epsilon_{0}$ would be, for example, 
\mbox{\rule[-3.5ex]{0ex}{6ex}$\epsilon_{0} = \min\Bigl\{\,\epsilon_{0, 1},\,\dfrac{\epsilon_{0, 2}}{125}\,\Bigr\}$}, 
where $\epsilon_{0, 1}$ is determined according to Corollary \ref{EWD_3_4_02} and 
$\epsilon_{0, 2}$ according to Corollary \ref{EWD_3_3_07}.\\[2.8ex]
We now assume\\[2.5ex]
\refstepcounter{DSBcount}
\label{EWD_3_7_11}
\text{\hspace*{-0.8ex}(\theDSBcount)}
\hspace*{4ex}
$\beta_{A} \leq \epsilon_{0}\,n$
\hspace*{2ex}and\hspace*{2ex}
$n \geq n_{0}$\\[2.5ex]
for this proof without loss of generality, so that (\ref{EWD_3_7_10}) is valid.\\[2.8ex] 
If, on the other hand, $\beta_{A} > \epsilon_{0}\,n$ or $n < n_{0}$, we proceed
as in the 2nd case of the proof of Theorem \ref{EWD_3_1_10} above. We obtain
for matrices $A$ satisfying $\sigma_{\!A} > 0$ and 
\mbox{\rule[0ex]{0ex}{3.8ex}$\sigma_{\!A'}^2 \geq \dfrac{2}{3}$}:\\[2ex]
\hspace*{4ex}\begin{tabular}[t]{@{}c@{\hspace*{0.8ex}}l@{\hspace*{-14.1ex}}r@{}}
\hspace*{1.8ex}&$\displaystyle{\sup\limits_{z\, \in\, \mathbb{R}}\, 
\bigg|\,E\Bigl(\,T_{\bm\bar{A}\,}^m\, 1_{(- \infty,\, z\,]}(\,\bm\bar{A}\,)\,\Bigr) 
- e_{1,\bm\bar{A}}^{m}(z)\,\bigg|}$\\[3.5ex]
$\leq$&$\displaystyle{E\bigl(\,\bigl|\,T_{\bm\bar{A}\,}\,\bigr|^m\,\bigr) + ||\,e_{1,\bm\bar{A}}^{m}\,||}$\\[2.4ex] 
$\leq$&$c_{1}(m) + c_{2}(m)\,\dfrac{\beta_{\bm\bar{A}}}{n}$
&(cf. Lemma \ref{EWD_3_1_18}, \ref{EWD_3_1_18_BWg}), \ref{EWD_3_1_18_BWj}) and 
$E\bigl(\,\bigl|\,T_{\bm\bar{A}\,}\,\bigr|\,\bigr) \leq 1$)\\[2.4ex]
$\leq$&$c_{1}(m) + c_{3}(m)\,\dfrac{\beta_{A}}{n}$
&(cf. Lemma \ref{EWD_3_4_05} and (\ref{EWD_3_7_02}), i.e. $\sigma_{\!A'}^2 \geq \dfrac{2}{3}$)\\[3ex]
$\leq$&$\biggl[\,c_{1}(m)\,\max\,\Bigl\{\,\dfrac{1}{\epsilon_{0}^2},\,4\,n_{0}\,\Bigr\} 
+ c_{3}(m)\,\max\,\Bigl\{\,\dfrac{1}{\epsilon_{0}},\,\sqrt{\rule[0ex]{0ex}{2.2ex}{4\,n_{0}}}\,\Bigr\}\,
\biggr]\,G_{\!A}^2$
&(cf. Lemma \ref{EWD_3_1_18}, \ref{EWD_3_1_18_BWa}), \ref{EWD_3_1_18_BWb})).
\end{tabular}\\[3ex]
A sufficiently large choice of $\mathfrak{K}_{2}(m)$ thus ensures the validity 
of the Theorem \ref{EWD_3_7_01} for such matrices as well.\\[2.8ex]
We further note that we use the definition\\[2.4ex]
\refstepcounter{DSBcount}
\label{EWD_3_7_12}
\text{\hspace*{-0.8ex}(\theDSBcount)}
\hspace*{2.8ex}
$\lambda = \dfrac{G_{\!A}}{\sigma_{\!A'}}$.\\[2.4ex]
for this proof. It follows from (\ref{EWD_3_7_10}) that\\[2.4ex]
\refstepcounter{DSBcount}
\label{EWD_3_7_13}
\text{\hspace*{-0.8ex}(\theDSBcount)}
\hspace*{2.8ex}
$\sqrt{\dfrac{3}{4}}\,G_{\!A} \leq \lambda \leq \sqrt{\dfrac{3}{2}}\,G_{\!A}$.\\[2.5ex]
Since $\lambda$ is therefore fixed in this proof, we will suppress in the following for the frequently used functions 
$q_{z}^{\lambda,m}$ (see (\ref{EWD_2_3_02})) and $f_{z}^{\lambda,m}$ (see (\ref{EWD_2_3_04}))
the dependence on $\lambda$ and set\\[2.4ex]
\refstepcounter{DSBcount}
\label{EWD_3_7_14}
\text{\hspace*{-0.8ex}(\theDSBcount)}
\hspace*{2.8ex}
$q_{z, m} = q_{z}^{\lambda,m}$
\hspace*{2ex}and\hspace*{2ex}
$f_{z, m}  = f_{z}^{\lambda,m}$.\index{function!smooth!$q_{z, m}$, $f_{z, m}$}\\[2.4ex]
After these preliminary remarks, we now turn to the central part of the proof. 
As the proof of Theorem \ref{EWD_1_1_03}, this is divided into 4 parts or subsections.\vspace*{1.2ex}

\subsubsection[Transition to smooth functions]{Transition to smooth functions}\label{EWD_Kap3_Sec7_2_1}

In order to apply Stein's method in the next subsection \ref{EWD_Kap3_Sec7_2_2}, we first proceed analogously to 
chapter \ref{EWD_Kap1}, section \ref{EWD_Kap1_Sec2} and replace in the expression\\[2ex]
\hspace*{12.1ex}$\displaystyle{\sup\limits_{z\, \in\, \mathbb{R}}\, 
\bigg|\,E\Bigl(\,T_{\bm\bar{A}\,}^m\, 1_{(- \infty,\, z\,]}(\,T_{\bm\bar{A}}\,)\,\Bigr) 
- \int\limits_{- \infty}^{z} x^m\,e_{1,\bm\bar{A}}'(x)\,dx\,\bigg|}$\\[2.3ex]
the functions $x^m\,1_{(\,- \infty,\,z\,]}(x)$, $z \in \mathbb{R}$, by smoother functions.
The functions $q_{z, m}$, $z \in \mathbb{R}$, prove to be suitable. We have\\[4ex]
\refstepcounter{DSBcount}
\label{EWD_3_7_15}
\textbf{\hspace*{-0.8ex}\theDSBcount\ Proposition}\\[1.5ex]
\hspace*{12.1ex}\begin{tabular}[t]{@{}c@{\hspace*{0.8ex}}l@{}}
&$\displaystyle{\sup\limits_{z\, \in\, \mathbb{R}}\, 
\bigg|\,E\Bigl(\,T_{\bm\bar{A}\,}^m\, 1_{(- \infty,\, z\,]}(\,T_{\bm\bar{A}}\,)\,\Bigr) 
- e_{1,\bm\bar{A}}^{m}(z)\,\bigg|}$\\[3ex]
$\leq$&$\displaystyle{\sup\limits_{z\, \in\, \mathbb{R}}\, 
\bigg|\,E\Bigl(\,q_{z, m}(\,T_{\bm\bar{A}}\,)\,\Bigr) 
- \int\limits_{\mathbb{R}} q_{z,m}(x)\,e_{1,\bm\bar{A}}'(x)\,dx\,\bigg| 
+ \Bigl(\,\mathfrak{L}(m) + C(m)\,\Bigr)\,G_{\!A}^2}$.
\end{tabular}\\[3.5ex]
\textbf{Proof:}
\begin{enumerate}
\item
Integration by parts\index{integration!by parts for Lebesgue-Stieltjes integrals} 
for Lebesgue-Stieltjes integrals (cf. e.g. \cite{athreya2006measure}, Theorem 5.2.3)
and use of $\dfrac{d}{dx} x^m = m$ for $m = 0$ and $m = 1$ gives\\[2.5ex]
\begin{tabular}{@{\hspace*{3.7ex}}l@{\hspace*{0.8ex}}l@{}}
\hspace*{1.9ex}&$\bigg|\,E\Bigl(\,T_{\bm\bar{A}\,}^m\, 1_{(- \infty,\, z\,]}(\,T_{\bm\bar{A}}\,)\,\Bigr)
- E\Bigl(\,q_{z - \lambda, m}(\,T_{\bm\bar{A}}\,)\,\Bigr)\,\bigg|$
\end{tabular}\\[3ex]
\begin{tabular}{@{\hspace*{3.7ex}}l@{\hspace*{0.8ex}}l@{}}
$=$&$\displaystyle{\bigg|\,\int\limits_{\textstyle{\bigl(\,z - \lambda,\,z\,\bigr]}}
x^m\,\dfrac{\bigl(\,x - (\,z - \lambda\,)\,\bigr)^2}{2\,\lambda^2}\,dF_{\bm\bar{A}}(x)\,\, - 
\int\limits_{\textstyle{\bigl(\,z,\,z + \lambda\,\bigr]}}
x^m\,\dfrac{\bigl(\,(\,z + \lambda\,) - x\,\bigr)^2}{2\,\lambda^2}\,dF_{\bm\bar{A}}(x)\,\bigg|}$\\[4.5ex]
$=$&$\displaystyle{\bigg|\,z^m\,F_{\bm\bar{A}}(z)\, - \int\limits_{z - \lambda}^{z}
x^m\,\dfrac{\bigl(\,x - (\,z - \lambda\,)\,\bigr)}{\lambda^2}\,F_{\bm\bar{A}}(x)\,dx\, -
\int\limits_{z}^{z + \lambda}
x^m\,\dfrac{\bigl(\,(\,z + \lambda\,) - x\,\bigr)}{\lambda^2}\,F_{\bm\bar{A}}(x)\,dx}$\\[3.5ex]
&$\displaystyle{-\ m \int\limits_{z - \lambda}^{z}
\dfrac{\bigl(\,x - (\,z - \lambda\,)\,\bigr)^2}{2\,\lambda^2}\,F_{\bm\bar{A}}(x)\,dx\, +\,
m \int\limits_{z}^{z + \lambda}
\dfrac{\bigl(\,(\,z + \lambda\,) - x\,\bigr)^2}{2\,\lambda^2}\,F_{\bm\bar{A}}(x)\,dx\,\bigg|}$.
\end{tabular}\\[3.5ex]
We now substitute\index{integration!by substitution} in the first and third integral $y = z - x$, 
and in the second and fourth integral $y = x - z$. Because of
\mbox{\rule[0ex]{0ex}{6.5ex}$\displaystyle{\int\limits_{0}^{\lambda} \dfrac{\lambda - y}{\lambda^2}\,dy = \dfrac{1}{2}}$}
we then obtain\\[0.8ex]
\begin{tabular}{@{\hspace*{3.7ex}}l@{\hspace*{0.8ex}}l@{}}
$=$&$\displaystyle{\bigg|\,\int\limits_{0}^{\lambda} \dfrac{\lambda - y}{\lambda^2}\,
\Bigl\{\,(z + y)^m\,F_{\bm\bar{A}}(z+y) - 2\,z^m\,F_{\bm\bar{A}}(z) + (z-y)^m\,F_{\bm\bar{A}}(z-y)\,\Bigr\}\,dy}$\\[3ex]
&$\displaystyle{-\ m\,\int\limits_{0}^{\lambda} \dfrac{\bigl(\,\lambda - y\,\bigr)^2}{2\,\lambda^2}\,
\Bigl\{\,F_{\bm\bar{A}}(z+y) - F_{\bm\bar{A}}(z-y)\,\Bigr\}\,dy\,\bigg|}$\\[3.5ex]
$=$&$\displaystyle{\bigg|\,\int\limits_{0}^{\lambda} \dfrac{\lambda - y}{\lambda^2}\,
(z-y)^m\,\Delta_{y}^{2}F_{\bm\bar{A}}(z - y)\,dy\,
+ m\,\int\limits_{0}^{\lambda} \dfrac{\lambda - y}{\lambda^2}\,2\,y\,
\Bigl\{\,F_{\bm\bar{A}}(z+y) - F_{\bm\bar{A}}(z)\,\Bigr\}\,dy}$\\[3.5ex]
\hspace*{1.9ex}&$\displaystyle{-\ m\,\int\limits_{0}^{\lambda} \dfrac{\bigl(\,\lambda - y\,\bigr)^2}{2\,\lambda^2}\,
\Bigl\{\,F_{\bm\bar{A}}(z+y) - F_{\bm\bar{A}}(z-y)\,\Bigr\}\,dy\,\bigg|}$.
\end{tabular}\\[3ex]
Next, we apply condition (\ref{EWD_3_7_03}) to the first integral and 
Theorem \ref{EWD_3_1_03} twice to each of the other two integrals.
We get\\[2.5ex]
\begin{tabular}{@{\hspace*{3.7ex}}l@{\hspace*{0.8ex}}l@{}}
$\leq$&$\displaystyle{\int\limits_{0}^{\lambda} \dfrac{\lambda - y}{\lambda^2}\,
\mathfrak{L}(m)\,\bigl(\,G_{\!A}^2 + y^2\,\bigr)\,dy\,
+ m\,\int\limits_{0}^{\lambda} \dfrac{\lambda - y}{\lambda^2}\,2\,y\,
\Bigl\{\,2\,\mathcal{K}_{1}\,\dfrac{\beta_{\bm\bar{A}}}{n} + \dfrac{y}{\sqrt{2\,\pi}}\,\Bigr\}\,dy}$\\[3.5ex]
&$\displaystyle{+\ m\,\int\limits_{0}^{\lambda} \dfrac{\bigl(\,\lambda - y\,\bigr)^2}{2\,\lambda^2}\,
\Bigl\{\,2\,\mathcal{K}_{1}\,\dfrac{\beta_{\bm\bar{A}}}{n} + \dfrac{2\,y}{\sqrt{2\,\pi}}\,\Bigr\}\,dy}$\\[4ex]
$=$&$\displaystyle{\dfrac{1}{2}\,\mathfrak{L}(m)\,G_{\!A}^2\, +\, \mathfrak{L}(m)\,\dfrac{\lambda^2}{12}\,
+\, m\,\biggl[\,\,\dfrac{2}{3}\,\mathcal{K}_{1}\,\lambda\,\dfrac{\beta_{\bm\bar{A}}}{n}\, 
+\, \dfrac{1}{\sqrt{2\,\pi}}\,\dfrac{\lambda^2}{6}\, 
+\, \dfrac{1}{3}\,\mathcal{K}_{1}\,\lambda\,\dfrac{\beta_{\bm\bar{A}}}{n}\,
+\, \dfrac{1}{\sqrt{2\,\pi}}\,\dfrac{\lambda^2}{12}\,\,\biggr]}$\\[3.5ex]
$\leq$&$\Bigl(\,\mathfrak{L}(m) + c_{1}(m)\,\Bigr)\,G_{\!A}^2$.
\end{tabular}\\[2.5ex]
For the last inequality, we used $\lambda \leq \sqrt{\dfrac{3}{2}}\,G_{\!A}$ (cf. (\ref{EWD_3_7_13})) 
and $\dfrac{\beta_{\bm\bar{A}}}{n} \leq 125\,\dfrac{\beta_{A}}{n} \leq 125\,G_{\!A}$ \mbox{\rule[0ex]{0ex}{3.2ex}(cf.
Lemma \ref{EWD_3_4_05} and Lemma \ref{EWD_3_1_18}, \ref{EWD_3_1_18_BWb})).}
\item
Furthermore, we have\\[2.5ex]
\begin{tabular}{@{\hspace*{3.7ex}}l@{\hspace*{0.8ex}}l@{}}
&$\displaystyle{\bigg|\,e_{1,\bm\bar{A}}^{m}(z) 
- \int\limits_{\mathbb{R}} q_{z - \lambda,m}(x)\,e_{1,\bm\bar{A}}'(x)\,dx\,\bigg|}$\\[3.5ex]
$=$&$\displaystyle{\bigg|\,\int\limits_{\mathbb{R}} x^{m}\,1_{(- \infty,\, z\,]}(x)\,e_{1,\bm\bar{A}}'(x)\,dx 
- \int\limits_{\mathbb{R}} q_{z - \lambda,m}(x)\,e_{1,\bm\bar{A}}'(x)\,dx\,\bigg|}$.
\end{tabular}\\[3ex]
We estimate this as in part a) with $e_{1,\bm\bar{A}}(x)$ instead of $F_{\bm\bar{A}}(x)$.
Using Lemma \ref{EWD_2_4_02}, \ref{EWD_2_4_02_BWc}) and \ref{EWD_2_4_02_BWb}) then yields\\[2.5ex]
\begin{tabular}{@{\hspace*{3.7ex}}l@{\hspace*{0.8ex}}l@{}}
$\leq$&$\displaystyle{\int\limits_{0}^{\lambda} \dfrac{\lambda - y}{\lambda^2}\,
\biggl(\,m\,\Bigl(\,4\,\big|\big|\,e_{1,\bm\bar{A}}'\,\big|\big| + \big|\big|\,x\,
e_{1,\bm\bar{A}}''(x)\,\big|\big|\,\Bigr) + 
\bigl(\,1 - m\,\bigr)\,\big|\big|\,e_{1,\bm\bar{A}}''\,\big|\big|\,\biggr)\,y^2\,dy}$\\[3.5ex]
&$\displaystyle{+\ m\,\int\limits_{0}^{\lambda} \dfrac{\lambda - y}{\lambda^2}\,2\,y\,
\big|\big|\,e_{1,\bm\bar{A}}'\,\big|\big|\,y\,dy 
+ m\,\int\limits_{0}^{\lambda} \dfrac{\bigl(\,\lambda - y\,\bigr)^2}{2\,\lambda^2}\,
\big|\big|\,e_{1,\bm\bar{A}}'\,\big|\big|\,2\,y\,dy}$\\[4ex]
$=$&$\displaystyle{\biggl(\,m\,\Bigl(\,4\,\big|\big|\,e_{1,\bm\bar{A}}'\,\big|\big| + \big|\big|\,x\,
e_{1,\bm\bar{A}}''(x)\,\big|\big|\,\Bigr) + 
\bigl(\,1 - m\,\bigr)\,\big|\big|\,e_{1,\bm\bar{A}}''\,\big|\big|\,\biggr)\,\dfrac{\lambda^2}{12}}$\\[3.5ex]
&$\displaystyle{+\ m\,\biggl(\,\,\big|\big|\,e_{1,\bm\bar{A}}'\,\big|\big|\,\dfrac{\lambda^2}{6}\,+\,
\big|\big|\,e_{1,\bm\bar{A}}'\,\big|\big|\,\dfrac{\lambda^2}{12}\,\,\biggr)}$\\[3.5ex]
$\leq$&$c_{2}(m)\,G_{\!A}^2$.
\end{tabular}\\[2.5ex]
For the last inequality, we used $\lambda \leq \sqrt{\dfrac{3}{2}}\,G_{\!A}$ (cf. (\ref{EWD_3_7_13}))
and Lemma \ref{EWD_3_1_18}, \ref{EWD_3_1_18_BWg}) together with
$\dfrac{\beta_{\bm\bar{A}}}{n} \leq 125\,\dfrac{\beta_{A}}{n} \leq 125$ 
(cf. Lemma \ref{EWD_3_4_05} and (\ref{EWD_3_7_11}) and $\epsilon_{0} \leq 1$).
\item
The assertion now follows by putting the parts a) and b) together.\hspace*{1ex}\hfill$\Box$
\end{enumerate}\vspace*{1ex}

\subsubsection[Application of Stein's method]{Application of Stein's method}\label{EWD_Kap3_Sec7_2_2}

In this part, we carry out the central part of the proof of Theorem \ref{EWD_3_7_01} using Stein's method.
Similar to chapter \ref{EWD_Kap1}, section \ref{EWD_Kap1_Sec3} (cf. in particular (\ref{EWD_1_3_03})), we apply Stein's equation\index{Stein's equation}\\[2.8ex]
\refstepcounter{DSBcount}
\label{EWD_3_7_16}
\text{\hspace*{-0.8ex}(\theDSBcount)}
\hspace*{2.8ex}
$f'_{z,m}(x) - x\,f_{z,m}(x) = q_{z,m}(x) - \Phi(q_{z,m})$
\hspace*{4ex}for $x \in \mathbb{R}$\\[2.5ex]
to the expression $E\bigl(\,q_{z, m}(\,T_{\bm\bar{A}}\,)\,\bigr) - \Phi\bigl(\,q_{z, m}\,\bigr)$ 
and then expand the functions $f_{z,m}$ and $f'_{z,m}$ according to 
Taylor's theorem\index{Taylor's theorem}.\\[2.8ex]
In these expansions, the statistics $T_{1}, \ldots, T_{5}$ (cf. (\ref{EWD_3_5_21})) come into play, 
which will always refer to the matrix $\bm\bar{A}$ in this and the following parts of the proof.
With this agreement we 
\linebreak
have\\[4ex]
\refstepcounter{DSBcount}
\label{EWD_3_7_17}
\textbf{\hspace*{-0.8ex}\theDSBcount\ Proposition}\\[0.8ex]
For every $z \in \mathbb{R}$,\\[2ex]
\hspace*{12.1ex}$\displaystyle{E\bigl(\,q_{z, m}(\,T_{\bm\bar{A}}\,)\,\bigr) - \Phi\bigl(\,q_{z, m}\,\bigr)
+ \dfrac{1}{2}\,E\bigl(\,T_{\bm\bar{A}}^3\,\bigr)\,E\bigl(\,T_{\bm\bar{A}}\,f'_{z,m}(\,T_{\bm\bar{A}}\,)\,\bigr)
= R(z,\,m)}$,\\[2.5ex]
where\\[2ex]
\begin{tabular}{@{\hspace*{7ex}}l@{\hspace*{0.8ex}}c@{\hspace*{0.8ex}}l@{}}
$R(z,\,m)$&$=$&$\displaystyle{\dfrac{1}{2}\,E\bigl(\,T_{\bm\bar{A}}^3\,\bigr)\,n\,
E\biggl(\,\bm\bar{a}_{I_{1}J_{1}}\,\Delta T_{4}\,\int\limits_{0}^{1}
\Bigl(\,f''_{z,m}(\,T_{3} + \Delta T_{3} + t\,\Delta T_{4}\,) - f''_{z,m}(\,T_{3}\,)\,\Bigr)\,dt\,\biggr)}$\\[3.5ex]
&&$\displaystyle{-\ n\,
E\biggl(\,\bm\bar{a}_{I_{1}J_{1}}\,\Delta T_{4}\,\Delta T_{3}\int\limits_{0}^{1}
\Bigl(\,f''_{z,m}(\,T_{2} + \Delta T_{2} + t\,\Delta T_{3}\,) - f''_{z,m}(\,T_{2}\,)\,\Bigr)\,dt\,\biggr)}$\\[3.5ex]
&&$\displaystyle{-\ n\,
E\biggl(\,\bm\bar{a}_{I_{1}J_{1}}\,\bigl(\,\Delta T_{4}\,\bigr)^2\,\int\limits_{0}^{1} (1-t)\,
\Bigl(\,f''_{z,m}(\,T_{3} + \Delta T_{3} + t\,\Delta T_{4}\,) - f''_{z,m}(\,T_{3}\,)\,\Bigr)\,dt\,\biggr)}$.
\end{tabular}\\[3.5ex]
\textbf{Proof:}\\[0.8ex]
In the following, let $m \in \{\,0,\,1\,\}$ and $z \in \mathbb{R}$ be fixed. 
We can therefore set $q = q_{z,m}$ and $f = f_{z,m}$.
Using (\ref{EWD_3_7_16}), we then get\\[2ex]
\hspace*{12.1ex}$E\bigl(\,q(\,T_{\bm\bar{A}}\,)\,\bigr) - \Phi\bigl(\,q\,\bigr) =
E\bigl(\,f'(\,T_{\bm\bar{A}}\,)\,\bigr) - E\bigl(\,T_{\bm\bar{A}}\,f(\,T_{\bm\bar{A}}\,)\,\bigr)$.\\[2.5ex]
We will now look at the last term on its own. 
Because of Lemma \ref{EWD_3_5_25}, \ref{EWD_3_5_25_BWa}) we have\\[2ex]
\hspace*{12.1ex}$E\bigl(\,T_{\bm\bar{A}}\,f(\,T_{\bm\bar{A}}\,)\,\bigr) =
n\,E\bigl(\,\bm\bar{a}_{I_{1}J_{1}}\,f(\,T_{5}\,)\,\bigr)$.\\[2.5ex]
Next, we use Taylor's theorem\index{Taylor's theorem} and the sums $T_{5} = T_{4} + \Delta T_{4}$ and 
$T_{4} = T_{3} + \Delta T_{3}$ (cf. (\ref{EWD_3_5_21})). We expand $f$ about $T_{4}$ and obtain\\[2ex]
\hspace*{12.1ex}\begin{tabular}[t]{@{}c@{\hspace*{0.8ex}}l@{}}
$=$&$\displaystyle{n\,E\bigl(\,\bm\bar{a}_{I_{1}J_{1}}\,f(\,T_{4}\,)\,\bigr) + 
n\,E\bigl(\,\bm\bar{a}_{I_{1}J_{1}}\,\Delta T_{4}\,f'(\,T_{4}\,)\,\bigr)}$
\end{tabular}\\[2ex]
\hspace*{12.1ex}\begin{tabular}[t]{@{}c@{\hspace*{0.8ex}}l@{}}
&$\displaystyle{+\ n\,E\biggl(\,\bm\bar{a}_{I_{1}J_{1}}\,\bigl(\,\Delta T_{4}\,\bigr)^2\,
\int\limits_{0}^{1}\,(1 - t)\,f''(\,T_{4} + t\,\Delta T_{4}\,)\,dt\,\biggr)}$\\[3.5ex]
$=$&$\displaystyle{n\,E\bigl(\,\bm\bar{a}_{I_{1}J_{1}}\,f(\,T_{4}\,)\,\bigr) + 
n\,E\bigl(\,\bm\bar{a}_{I_{1}J_{1}}\,\Delta T_{4}\,f'(\,T_{4}\,)\,\bigr) +
n\,E\Bigl(\,\bm\bar{a}_{I_{1}J_{1}}\,\dfrac{\bigl(\,\Delta T_{4}\,\bigr)^2}{2}\,f''(\,T_{3}\,)\,\Bigr)}$\\[2.5ex]
&$\displaystyle{+\ n\,E\biggl(\,\bm\bar{a}_{I_{1}J_{1}}\,\bigl(\,\Delta T_{4}\,\bigr)^2\,
\int\limits_{0}^{1}\,(1 - t)\,\Bigl(\,f''(\,T_{3} + \Delta T_{3} + t\,\Delta T_{4}\,)
- f''(\,T_{3}\,)\,\Bigr)\,dt\,\biggr)}$.
\end{tabular}\\[2.7ex]
For the first and third summand we get, because of Lemma \ref{EWD_3_5_23}, \ref{EWD_3_5_23_BWa}) and 
\ref{EWD_3_5_23_BWb}), and Lemma \ref{EWD_3_5_26}, \ref{EWD_3_5_26_BWa}),\\[2ex]
\hspace*{12.1ex}$n\,E\bigl(\,\bm\bar{a}_{I_{1}J_{1}}\,f(\,T_{4}\,)\,\bigr) =
n\,E\bigl(\,\bm\bar{a}_{I_{1}J_{1}}\,\bigr)\,E\bigl(\,f(\,T_{4}\,)\,\bigr) = 0$,\\[2.5ex]
\hspace*{12.1ex}$n\,E\Bigl(\,\bm\bar{a}_{I_{1}J_{1}}\,\dfrac{\bigl(\,\Delta T_{4}\,\bigr)^2}{2}\,f''(\,T_{3}\,)\,\Bigr)
= n\,E\Bigl(\,\bm\bar{a}_{I_{1}J_{1}}\,\dfrac{\bigl(\,\Delta T_{4}\,\bigr)^2}{2}\,\Bigr)
\,E\bigl(\,f''(\,T_{3}\,)\,\bigr)$.\\[2.7ex]
In the second summand, however, we expand $f'$ about $T_{3}$ and use the sums 
$T_{4} = T_{3} + \Delta T_{3}$ and $T_{3} = T_{2} + \Delta T_{2}$ (cf. (\ref{EWD_3_5_21}))\\[2ex]
\hspace*{12.1ex}\begin{tabular}[t]{@{}c@{\hspace*{0.8ex}}l@{}}
&$\displaystyle{n\,E\bigl(\,\bm\bar{a}_{I_{1}J_{1}}\,\Delta T_{4}\,f'(\,T_{4}\,)\,\bigr)}$\\[2.5ex]
$=$&$\displaystyle{n\,E\bigl(\,\bm\bar{a}_{I_{1}J_{1}}\,\Delta T_{4}\,f'(\,T_{3}\,)\,\bigr)}$\\[1.7ex]
&$\displaystyle{+\ n\,E\biggl(\,\bm\bar{a}_{I_{1}J_{1}}\,\Delta T_{4}\,\Delta T_{3}\,
\int\limits_{0}^{1}\,f''(\,T_{3} + t\,\Delta T_{3}\,)\,dt\,\biggr)}$\\[4ex]
$=$&$\displaystyle{n\,E\bigl(\,\bm\bar{a}_{I_{1}J_{1}}\,\Delta T_{4}\,f'(\,T_{3}\,)\,\bigr)
+ n\,E\bigl(\,\bm\bar{a}_{I_{1}J_{1}}\,\Delta T_{4}\,\Delta T_{3}\,f''(\,T_{2}\,)\,\bigr)}$
\end{tabular}\\[2ex]
\hspace*{12.1ex}\begin{tabular}[t]{@{}c@{\hspace*{0.8ex}}l@{}}
&$\displaystyle{+\ n\,
E\biggl(\,\bm\bar{a}_{I_{1}J_{1}}\,\Delta T_{4}\,\Delta T_{3}\int\limits_{0}^{1}
\Bigl(\,f''(\,T_{2} + \Delta T_{2} + t\,\Delta T_{3}\,) - f''(\,T_{2}\,)\,\Bigr)\,dt\,\biggr)}$\\[4ex]
$=$&$\displaystyle{E\bigl(\,f'(\,T_{3}\,)\,\bigr)
+ n\,E\bigl(\,\bm\bar{a}_{I_{1}J_{1}}\,\Delta T_{4}\,\Delta T_{3}\,\bigr)\,E\bigl(\,f''(\,T_{2}\,)\,\bigr)}$\\[1.7ex]
&$\displaystyle{+\ n\,
E\Bigl(\,\bm\bar{a}_{I_{1}J_{1}}\,\Delta T_{4}\,\Delta T_{3}\int\limits_{0}^{1}
\Bigl(\,f''(\,T_{2} + \Delta T_{2} + t\,\Delta T_{3}\,) - f''(\,T_{2}\,)\,\Bigr)\,dt\,\Bigr)}$.
\end{tabular}\\[2.5ex]
For the last equation, Lemma \ref{EWD_3_5_23}, \ref{EWD_3_5_23_BWb}) and 
\ref{EWD_3_5_23_BWc}), and Lemma \ref{EWD_3_5_26}, \ref{EWD_3_5_26_BWb}) were applied.\\[2.8ex]
Overall, taking into account Lemma \ref{EWD_3_5_15}, \ref{EWD_3_5_15_BWa}) and 
Lemma \ref{EWD_3_5_26}, \ref{EWD_3_5_26_BWc}), we have\\[2ex]
\refstepcounter{DSBcount}
\label{EWD_3_7_18}
\text{\hspace*{-0.8ex}(\theDSBcount)}
\hspace*{2.8ex}
\begin{tabular}[t]{@{}c@{\hspace*{0.8ex}}l@{}}
&$E\bigl(\,T_{\bm\bar{A}}\,f(\,T_{\bm\bar{A}}\,)\,\bigr)$\\[2.5ex]
$=$&$E\bigl(\,f'(\,T_{\bm\bar{A}}\,)\,\bigr) 
+ E\bigl(\,f''(\,T_{\bm\bar{A}}\,)\,\bigr)\,\dfrac{1}{2}\,E\bigl(\,T_{\bm\bar{A}}^3\,\bigr)$\\[2ex]
&$\displaystyle{+\ n\,
E\biggl(\,\bm\bar{a}_{I_{1}J_{1}}\,\Delta T_{4}\,\Delta T_{3}\int\limits_{0}^{1}
\Bigl(\,f''(\,T_{2} + \Delta T_{2} + t\,\Delta T_{3}\,) - f''(\,T_{2}\,)\,\Bigr)\,dt\,\biggr)}$\\[3.5ex]
&$\displaystyle{+\ n\,E\biggl(\,\bm\bar{a}_{I_{1}J_{1}}\,\bigl(\,\Delta T_{4}\,\bigr)^2\,
\int\limits_{0}^{1}\,(1 - t)\,\Bigl(\,f''(\,T_{3} + \Delta T_{3} + t\,\Delta T_{4}\,)
- f''(\,T_{3}\,)\,\Bigr)\,dt\,\biggr)}$.
\end{tabular}\\[2.5ex]
The degree of the derivative of $f$ in the term $E\bigl(\,f''(\,T_{\bm\bar{A}}\,)\,\bigr)$ 
will now be reduced again. 
If we argue as in (\ref{EWD_3_6_05}) with $f'$ instead of $d$ and 
$\bm\bar{A}$ instead of $A$, we get\\[2ex]
\hspace*{12.1ex}\begin{tabular}[t]{@{}c@{\hspace*{0.8ex}}l@{}}
&$E\bigl(\,T_{\bm\bar{A}}\,f'(\,T_{\bm\bar{A}}\,)\,\bigr)$\\[1.5ex]
$=$&$\displaystyle{E\bigl(\,f''(T_{\bm\bar{A}}\,)\,\bigr) 
+ n\,E\biggl(\,\bm\bar{a}_{I_{1}J_{1}}\,\Delta T_{4}\,\int\limits_{0}^{1}\,
\Bigl(\,f''(\,T_{3} + \Delta T_{3} + t\,\Delta T_{4}\,)
- f''(\,T_{3}\,)\,\Bigr)\,dt\,\biggr)}$.
\end{tabular}\\[2.5ex]
By inserting the corresponding expression for $E\bigl(\,f''(T_{\bm\bar{A}}\,)\,\bigr)$ 
into the equation (\ref{EWD_3_7_18}), we obtain the assertion of the proposition.\hspace*{1ex}\hfill$\Box$\vspace*{0.5ex}

\subsubsection[Deducing the correct expansion]{Deducing the correct expansion}\label{EWD_Kap3_Sec7_2_3}

This subsection deals with the term
$\dfrac{1}{2}\,E\bigl(\,T_{\bm\bar{A}}^3\,\bigr)\,E\bigl(\,T_{\bm\bar{A}}\,f'_{z,m}(T_{\bm\bar{A}})\,\bigr)$. 
We prove the following\\[3.5ex]
\refstepcounter{DSBcount}
\label{EWD_3_7_19}
\textbf{\hspace*{-0.8ex}\theDSBcount\ Proposition}\\[2ex]
\hspace*{12.1ex}$\displaystyle{\sup\limits_{z\, \in\, \mathbb{R}}\, 
\bigg|\,\dfrac{1}{2}\,E\bigl(\,T_{\bm\bar{A}}^3\,\bigr)\,E\bigl(\,T_{\bm\bar{A}}\,f'_{z,m}(\,T_{\bm\bar{A}}\,)\,\bigr)
- \dfrac{1}{2}\,\lambda_{1,\bm\bar{A}}
\,\Phi\bigl(\,x\,f'_{z,m}(x)\,\bigr)\,\bigg|
\leq C(m)\,G_{\!A}^2}$.\\[3.5ex]
This gives us the correct Edgeworth expansion in Proposition \ref{EWD_3_7_17} 
except for an error of size $C(m)\,G_{\!A}^2$, since according to Lemma \ref{EWD_2_2_05}, \ref{EWD_2_2_05_BWc}) 
and (\ref{EWD_2_2_02}) holds:\\[2ex]
\hspace*{12.1ex}$\displaystyle{\Phi\bigl(\,x\,f'_{z,m}(x)\,\bigr) =
\dfrac{1}{3}\,\int\limits_{\mathbb{R}} q_{z,m}(x)\,\Bigl[\,(x^2-1)\,\psi(x)\,\Bigr]'\,dx}$.\\[4ex]
\textbf{Proof of Proposition \ref{EWD_3_7_19}:}\\[0.8ex]
Let $z \in \mathbb{R}$ be fixed.
An application of the triangle inequality and of Lemma \ref{EWD_3_1_18}, \ref{EWD_3_1_18_BWd}) yields\\[2ex]
\hspace*{12.1ex}\begin{tabular}[t]{@{}c@{\hspace*{0.8ex}}l@{}}
\hspace*{1.9ex}&$\displaystyle{\bigg|\,E\bigl(\,T_{\bm\bar{A}}^3\,\bigr)\,
E\bigl(\,T_{\bm\bar{A}}\,f'_{z,m}(\,T_{\bm\bar{A}}\,)\,\bigr)
- \lambda_{1,\bm\bar{A}}
\,\Phi\bigl(\,x\,f'_{z,m}(x)\,\bigr)\,\bigg|}$\\[3.5ex]
$\leq$&$\displaystyle{\bigg|\,E\bigl(\,T_{\bm\bar{A}}^3\,\bigr) - \lambda_{1,\bm\bar{A}}\,\bigg|\,
E\bigl(\,\big|\,T_{\bm\bar{A}}\,f'_{z,m}(\,T_{\bm\bar{A}}\,)\,\big|\,\bigr)}$\\[3.5ex]
&$+\ \dfrac{\beta_{\bm\bar{A}}}{n}\,\bigg|\,E\bigl(\,T_{\bm\bar{A}}\,f'_{z,m}(\,T_{\bm\bar{A}}\,)\,\bigr)
- \Phi\bigl(\,x\,f'_{z,m}(x)\,\bigr)\,\bigg|$\\[3.5ex] 
$=$&$A_{1} + A_{2}$.
\end{tabular}\\[3ex]
Because of Proposition \ref{EWD_3_2_02}, \ref{EWD_3_2_02_BWa})
and because of $|f'_{z,0}(x)| \leq 1$ (cf. Corollary \ref{EWD_2_1_12}, \ref{EWD_2_1_12_BWa})) and
\mbox{\rule[0ex]{0ex}{3.8ex}$|f'_{z,1}(x)| \leq |x| + \dfrac{1}{\sqrt{2\,\pi}}$} 
(cf. Corollary \ref{EWD_2_1_17}, \ref{EWD_2_1_17_BWc})), we obtain for $A_{1}$:\\[2ex]
\hspace*{12.1ex}\begin{tabular}[t]{@{}l@{\hspace*{0.8ex}}c@{\hspace*{0.8ex}}l@{\hspace*{10.6ex}}r@{}} 
$A_{1}$&$\leq$&$\dfrac{c_{1}(m)}{n^2}\,\beta_{\bm\bar{A}}$\\[2.7ex]
&$\leq$&$\dfrac{c_{2}(m)}{\sqrt{n}}\,\dfrac{\beta_{A}}{n}$
&(cf. Lemma \ref{EWD_3_4_05}; $\sigma_{\!A'}^2 \geq \dfrac{2}{3}$ and $\sqrt{n} \leq n$)\\[2.7ex]
&$\leq$&$c_{3}(m)\,\Bigl(\,\dfrac{\beta_{A}}{n}\,\Bigr)^2 \leq c_{3}(m)\,G_{\!A}^2$
&(cf. Lemma \ref{EWD_3_1_18},  \ref{EWD_3_1_18_BWa}) and \ref{EWD_3_1_18_BWb})).
\end{tabular}\\[2.8ex]
To estimate $A_{2}$, the following remains to be shown due to $\dfrac{\beta_{\bm\bar{A}}}{n}\leq 125\,G_{\!A}$:\\[2ex]
\refstepcounter{DSBcount}
\label{EWD_3_7_20}
\text{\hspace*{-0.8ex}(\theDSBcount)}
\hspace*{2.8ex}
$\bigg|\,E\bigl(\,T_{\bm\bar{A}}\,f'_{z,m}(\,T_{\bm\bar{A}}\,)\,\bigr)
- \Phi\bigl(\,x\,f'_{z,m}(x)\,\bigr)\,\bigg| \leq c_{4}(m)\,\dfrac{\beta_{\bm\bar{A}}}{n}$.\\[3ex]
In the case $m = 0$ we use (\ref{EWD_3_7_16}) and obtain $f'_{z,0}(x) = h_{1}(x) - h_{2}(x)$ with\\[2.3ex]
\hspace*{12.1ex}$h_{1}(x) = q_{z,0}(x) \in \mathcal{H}$\hspace*{2ex}and\hspace*{2ex}
$h_{2}(x) = -\,\,\bigl(\,x\,f_{z,0}(x) - \Phi(q_{z,0})\,\bigr) \in \mathcal{H}$
\index{function!$\mathcal{H}$}\\[2.8ex]
(cf. Lemma \ref{EWD_2_1_13}, \ref{EWD_2_1_13_BWa})). 
If we also take into account that we have
$|\bm\bar{a}_{ij}| \leq 1$ for all $1 \leq i, j \leq n$ (cf. (\ref{EWD_3_7_10})) and therefore the inequality\\[2ex]
\refstepcounter{DSBcount}
\label{EWD_3_7_21}
\text{\hspace*{-0.8ex}(\theDSBcount)}
\hspace*{2.8ex}
$\displaystyle{\dfrac{1}{n}\sum\limits_{i,j=1}^n |\,\bm\bar{a}_{ij}\,|^{k+4} 
\leq \dfrac{\beta_{\bm\bar{A}}}{n}}$,\\[2.5ex]
then we get (\ref{EWD_3_7_20}) from Corollary \ref{EWD_3_6_27}, 
applied with $k = 1$ to $h_{1}$, $h_{2}$.\\[2.8ex]
In the case $m = 1$ we also use (\ref{EWD_3_7_16}) and obtain
$f'_{z,1}(x) = x\,h_{3}(x) - x\,h_{4}(x) - h_{5}(x)$ with\\[2.3ex]
\hspace*{12.1ex}$h_{3}(x) = q_{z,0}(x) \in \mathcal{H}$, $h_{4}(x) = - f_{z,1}(x) \in \mathcal{H}$
\hspace*{2ex}and\hspace*{2ex}
$h_{5}(x) = \Phi(q_{z,1}) \in \mathcal{H}$\index{function!$\mathcal{H}$}\\[2.8ex]
(cf. Lemma \ref{EWD_2_1_13}, \ref{EWD_2_1_13_BWb})).
Thus (\ref{EWD_3_7_20}) follows from Corollary \ref{EWD_3_6_27}, applied with $k = 2$
to $h_{3}$, $h_{4}$ and with $k = 1$ to $h_{5}$.
\hspace*{1ex}\hfill$\Box$\vspace*{1ex}

\subsubsection[The estimate of $R(z,\,m)$]{The estimate of {\boldmath $R(z,\,m)$}}
\label{EWD_Kap3_Sec7_2_4}

In this subsection we finally show for the $R(z,\,m)$ from Proposition \ref{EWD_3_7_17}:\\[4ex]
\refstepcounter{DSBcount}
\label{EWD_3_7_22}
\textbf{\hspace*{-0.8ex}\theDSBcount\ Proposition}\\[1.5ex]
\hspace*{12.1ex}$\displaystyle{\sup\limits_{z\, \in\, \mathbb{R}}\, 
\big|\,R(z,\,m)\,\big|\, \leq\, \Bigl(\,C_{1}(m)\,\mathfrak{L}(m) + C_{2}(m)\,\Bigr)\,G_{\!A}^2}$.\\[3.5ex]
\textbf{Proof:}\\[0.8ex]
In the following, let $\overline{I} = (\,I_{1},\ldots,I_{8}\,)$, 
$\bm{\breve}{I} = (\,I_{1},\ldots,I_{4}\,)$ and analogously 
$\overline{i} = (\,i_{1},\ldots,i_{8}\,)$ $\in M_{8}$, 
$\bm{\breve}{i} = (\,i_{1},\ldots,i_{4}\,) \in M_{4}$.\index{random permutation!$\overline{i}$, $\bm{\breve}{i}$}
The symbols 
$\overline{J}$, $\bm{\breve}{J}$, $\overline{j}$ 
and $\bm{\breve}{j}$ are defined in the same 
way.\index{random permutation!$\overline{I}$, $\bm{\breve}{I}$, $\overline{J}$, $\bm{\breve}{J}$}
\index{random permutation!$\overline{j}$, $\bm{\breve}{j}$}\\[2.8ex]
To prove Proposition \ref{EWD_3_7_22}, we show that there exist 
constants $c_{1}(m)$, $c_{2}(m)$, $c_{3}(m)$ and $c_{4}(m)$ such that\\[3.5ex]
\refstepcounter{DSBcount}
\label{EWD_3_7_23}
\text{\hspace*{-0.8ex}(\theDSBcount)}
\hspace*{2.8ex}
\begin{tabular}[t]{@{}l@{}}
$\displaystyle{\bigg|\,E\Bigl(\,f''_{z,m}\bigl(\,T_{2} + \Delta T_{2} + t\,\Delta T_{3}\,\bigr)
- f''_{z,m}\bigl(\,T_{2}\,\bigr)\,\Big|\,\overline{I} = \overline{i},\, \overline{J} = \overline{j}\,\Bigr)\,\bigg|}$\\[3ex]
$\displaystyle{\leq\ \Bigl(\,c_{1}(m)\,\mathfrak{L}(m) + c_{2}(m)\,\Bigr)\,
E\Bigl(\,\big|\,\Delta T_{2}\,\big| + \big|\,\Delta T_{3}\,\big|
\,\Big|\,\overline{I} = \overline{i},\, \overline{J} = \overline{j}\,\Bigr)}$\\[3ex]
for all $z \in \mathbb{R}$, $0 \leq t \leq 1$ and $\overline{i}, \overline{j} \in M_{8}$ 
with $(\overline{i},\, \overline{j}) \in M_{16}$,
\end{tabular}\\[3.5ex]
\refstepcounter{DSBcount}
\label{EWD_3_7_24}
\text{\hspace*{-0.8ex}(\theDSBcount)}
\hspace*{2.8ex}
\begin{tabular}[t]{@{}l@{}}
$\displaystyle{\bigg|\,E\Bigl(\,f''_{z,m}\bigl(\,T_{3} + \Delta T_{3} + t\,\Delta T_{4}\,\bigr)
- f''_{z,m}\bigl(\,T_{3}\,\bigr)\,\Big|\,
\bm{\breve}{I} = \bm{\breve}{i},\, 
\bm{\breve}{J} = \bm{\breve}{j}\,\Bigr)\,\bigg|}$\\[3ex]
$\displaystyle{\leq\ \Bigl(\,c_{3}(m)\,\mathfrak{L}(m) + c_{4}(m)\,\Bigr)\,
E\Bigl(\,\big|\,\Delta T_{3}\,\big| + \big|\,\Delta T_{4}\,\big|
\,\Big|\,
\bm{\breve}{I} = \bm{\breve}{i},\, 
\bm{\breve}{J} = \bm{\breve}{j}\,\Bigr)}$\\[3ex]
for all $z \in \mathbb{R}$, $0 \leq t \leq 1$ and 
$\bm{\breve}{i}, \bm{\breve}{j} \in M_{4}$ 
with $(\bm{\breve}{i},\, \bm{\breve}{j}) \in M_{8}$.
\end{tabular}\\[3.5ex]
Because of Lemma \ref{EWD_3_5_25}, \ref{EWD_3_5_25_BWc}) and Proposition \ref{EWD_3_2_02}, \ref{EWD_3_2_02_BWc})
and\\[2.5ex]
\hspace*{12.1ex}$\Bigl(\,\dfrac{\beta_{\bm\bar{A}}}{n}\,\Bigr)^2 \leq 
c_{5}\, \Bigl(\,\dfrac{\beta_{A}}{n}\,\Bigr)^2 \leq c_{5}\,G_{\!A}^2$
\hspace*{2ex}and\hspace*{2ex}
$\dfrac{\delta_{\bm\bar{A}}}{n} \leq 
c_{6}\,\dfrac{\delta_{A}}{n} = c_{6}\,D_{\!A}^2 \leq c_{6}\,G_{\!A}^2$\\[2.5ex]
(cf. Lemma \ref{EWD_3_4_05} and Lemma \ref{EWD_3_1_18}, \ref{EWD_3_1_18_BWb})), 
Proposition \ref{EWD_3_7_22} then follows from (\ref{EWD_3_7_23}) and (\ref{EWD_3_7_24}).\\[2.8ex]
Since the proofs of (\ref{EWD_3_7_23}) and (\ref{EWD_3_7_24}) are very similar, 
we will only deal with the case (\ref{EWD_3_7_23}) as an example.\\[2.8ex]
For that we fix in the following the quantities $m \in \{\,0,\,1\,\}$, $z \in \mathbb{R}$, $0 \leq t \leq 1$ and $\overline{i}, \overline{j} \in M_{8}$ with $(\overline{i},\, \overline{j})$ $\in M_{16}$. 
We can therefore set $q = q_{z,m}$ and $f = f_{z,m}$.\\[2.8ex]
Because of Lemma \ref{EWD_3_5_27}, \ref{EWD_3_5_27_BWb}) is $T_{2}$, 
given $\overline{I} = \overline{i}$ and $\overline{J} = \overline{j}$, distributed as\\[2ex]
\refstepcounter{DSBcount}
\label{EWD_3_7_25}
\text{\hspace*{-0.8ex}(\theDSBcount)}
\hspace*{2.8ex}
$\displaystyle{\sum\limits_{(i,\,j)\, \in\, S} \bm\bar{a}_{ij} + T_{\SBTB}}$,\\[2.5ex]
where\\[2ex]
\refstepcounter{DSBcount}
\label{EWD_3_7_26}
\text{\hspace*{-0.8ex}(\theDSBcount)}
\hspace*{2.8ex}
$S = \Bigl\{\,(\,i_{s},\,j_{s+4}\,)\,:\,1\, \leq\, s\, \leq\, 4\,\Bigr\}\,
\cup\, \Bigl\{\,(\,i_{s},\,j_{s-4}\,)\,:\,5\, \leq\, s\, \leq\, 8\,\Bigr\}$\\[2.5ex]
and\\[2ex]
\refstepcounter{DSBcount}
\label{EWD_3_7_27}
\text{\hspace*{-0.8ex}(\theDSBcount)}
\hspace*{2.8ex}
\begin{tabular}[t]{@{}l@{}}
$\BTB$ is the $(n - l) \times (n - l)-$matrix obtained from $\bm\bar{A}$ by \textbf{cancelling} the\\[0.7ex]
$l = \big|\,\{\,i_{1},\ldots, i_{8}\,\}\,\big|$ rows $i_{1},\ldots,i_{8}$ and the $l$ columns
$j_{1},\ldots,j_{8}$.
\end{tabular}\\[3ex]
With the abbreviations\\[2ex]
\hspace*{12.1ex}\begin{tabular}[t]{@{}l@{\hspace*{0.8ex}}c@{\hspace*{0.8ex}}l@{\hspace*{15.4ex}}l@{}}
$a$&$=$&$\displaystyle{\sum\limits_{(i,\,j)\, \in\, S} \bm\bar{a}_{ij}}$,\\[3.5ex]
$p$&$=$&$\displaystyle{E\Bigl(\,\Delta T_{2}\,\Big|\,\overline{I} = \overline{i}, \overline{J} = \overline{j}\,\Bigr)}$
&($\Delta T_{2} \in \sigma\bigl(\,\overline{I},\,\overline{J}\,\bigr)$, 
cf. Lemma \ref{EWD_3_5_22}, \ref{EWD_3_5_22_BWb})),\\[2.5ex]
$q$&$=$&$\displaystyle{E\Bigl(\,\Delta T_{3}\,\Big|\,\overline{I} = \overline{i}, \overline{J} = \overline{j}\,\Bigr)}$
&($\Delta T_{3} \in \sigma\bigl(\,\overline{I},\,\overline{J}\,\bigr)$,
cf. Lemma \ref{EWD_3_5_22}, \ref{EWD_3_5_22_BWc}))
\end{tabular}\\[3ex]
remains to be shown:\\[2.5ex]
\refstepcounter{DSBcount}
\label{EWD_3_7_28}
\text{\hspace*{-0.8ex}(\theDSBcount)}
\hspace*{2.8ex}
\begin{tabular}[t]{@{}c@{\hspace*{0.8ex}}l@{}}
&$\displaystyle{\bigg|\,E\Bigl(\,f''\bigl(\,T_{\SBTB} + a + p + t\,q\,\bigr)
- f''\bigl(\,T_{\SBTB} + a\,\bigr)\,\Bigr)\,\bigg|}$\\[3ex]
$\leq$&$\displaystyle{\Bigl(\,c_{7}(m)\,\mathfrak{L}(m) + c_{8}(m)\,\Bigr)\,
\Bigl(\,|\,p\,| + |\,q\,|\,\Bigr)}$.
\end{tabular}\\[3ex]
If we estimate the left-hand side of the inequality (\ref{EWD_3_7_28}) using Lemma \ref{EWD_2_3_06}, 
\ref{EWD_2_3_06_BWb}) - \ref{EWD_2_3_06_BWe}) 
(with \mbox{\rule[0ex]{0ex}{2.8ex}$x = T_{\SBTB} + a$} and 
\mbox{\rule[0ex]{0ex}{2.8ex}$y = p + t\,q$)}, we get\\[2ex]
\hspace*{10.9ex}\begin{tabular}[t]{@{}c@{\hspace*{0.8ex}}l@{}}
&$\displaystyle{\bigg|\,E\Bigl(\,f''\bigl(\,T_{\SBTB} + a + p + t\,q\,\bigr)
- f''\bigl(\,T_{\SBTB} + a\,\bigr)\,\Bigr)\,\bigg|}$\\[3ex]
$\leq$&$\displaystyle{E\Bigl(\,\Big|\,\bigl(\,f'' - q'\,\bigr)\,\bigl(\,T_{\SBTB} + a + p + t\,q\,\bigr)
- \bigl(\,f'' - q'\,\bigr)\,\bigl(\,T_{\SBTB} + a\,\bigr)\,\Big|\,\Bigr)}$
\end{tabular}\\[3ex]
\hspace*{12.7ex}\begin{tabular}[t]{@{}c@{\hspace*{0.8ex}}l@{}}
&$\displaystyle{+\ \bigg|\,E\Bigl(\,q'\bigl(\,T_{\SBTB} + a + p + t\,q\,\bigr)
- q'\bigl(\,T_{\SBTB} + a\,\bigr)\,\Bigr)\,\bigg|}$
\end{tabular}\\[3ex]
\hspace*{10.9ex}\begin{tabular}[t]{@{}c@{\hspace*{0.8ex}}l@{}}
$\leq$&$\displaystyle{c_{9}(m)\,\Bigl(\,|\,p\,| + |\,q\,|\,\Bigr)
\,\Biggl\{\,1 + E\Bigl(\,\big|\,T_{\SBTB} + a\,\big|^{2 + m}\,\Bigr)
+ m\,\Bigl(\,|\,p\,| + |\,q\,|\,\Bigr)\,E\Bigl(\,1 + \bigl(\,T_{\SBTB} + a\,\bigr)^2\,\Bigr)}$\\[3.5ex]
&$\displaystyle{+\ \dfrac{1}{\lambda}\,\int\limits_{0}^{1}\,E\Bigl(\,\big|\,T_{\SBTB} + a\,\big|^{1 + m}\,
1_{\textstyle{(\,z,\,z\, + 2\,\lambda\,]}}\bigl(\,T_{\SBTB} + a + s\,p + s\,t\,q\,\bigr)\,\Bigr)\,ds}$\\[4.5ex]
&$\displaystyle{+\ \dfrac{1}{\lambda^2}\,\int\limits_{0}^{1}\,\bigg|\,E\Bigl(\,\bigl(\,T_{\SBTB} + a\,\bigr)^{m}\,
\bigl(\,1_{\textstyle{(\,z + \lambda,\,z\, + 2\,\lambda\,]}} - 1_{\textstyle{(\,z,\,z\, + \lambda\,]}}\,\bigr)}$\\[1.5ex]
&\hspace*{45ex}$\displaystyle{\bigl(\,T_{\SBTB} + a + s\,p + s\,t\,q\,\bigr)\,\Bigr)\,\bigg|\,ds}$\\[2.5ex]
&$\displaystyle{+\ \dfrac{m}{\lambda}\,\biggl(\,P\bigl(\,z < T_{\SBTB} + a + p + t\,q \leq z\, 
+ 2\,\lambda\,\bigr)}$\\[1.5ex]
&\hspace*{22.5ex}$\displaystyle{+\ \int\limits_{0}^{1}\,P\bigl(\,z < T_{\SBTB} + a + s\,p 
+ s\,t\,q \leq z\, + 2\,\lambda\,\bigr)\,ds
\,\biggr)\,\Biggr\}}$\\[4ex]
$=$&$\displaystyle{c_{9}(m)\,\Bigl(\,|\,p\,| + |\,q\,|\,\Bigr)
\,\Bigl(\,1 + A_{1} + A_{2} + A_{3} + A_{4} + A_{5}\,\Bigr)}$.
\end{tabular}\\[3.5ex] 
It therefore remains to estimate the $A_{i}$, $1 \leq i \leq 5$.\\[2.8ex]
Due to $|\bm\bar{a}_{ij}| \leq 1$ for all $1 \leq i, j \leq n$, we first obtain\\[2ex]
\hspace*{12.2ex}$|\,a\,| \leq 8$,
\hspace*{2ex}$|\,p\,| \leq 16$
\hspace*{2ex}and\hspace*{2ex}
$|\,q\,| \leq 8$.\\[2.5ex]
Furthermore, since $|\,\mu_{\SBTB}\,| \leq 1$ and $\dfrac{2}{3} \leq \sigma_{\SBTB}^2 \leq \dfrac{4}{3}$
(cf. (\ref{EWD_3_7_10})), we also get\\[2ex]
\hspace*{12.2ex}\begin{tabular}[t]{@{}r@{\hspace*{0.8ex}}c@{\hspace*{0.8ex}}l@{\hspace*{15ex}}l@{}}
$E\bigl(\,T_{\SBTB}^2\,\bigr)$&$=$&$ \sigma_{\SBTB}^2 + \mu_{\SBTB}^2 \leq \dfrac{7}{3}$\hspace*{2ex}and
\end{tabular}\\[2.5ex]
\refstepcounter{DSBcount}
\label{EWD_3_7_29}
\text{\hspace*{-0.8ex}(\theDSBcount)}
\hspace*{3ex}
\begin{tabular}[t]{@{}r@{\hspace*{0.8ex}}c@{\hspace*{0.8ex}}l@{\hspace*{14ex}}r@{}}
$E\bigl(\,\big|\,T_{\SBTB}\,\big|^3\,\bigr)$&$=$&
$E\bigl(\,\big|\,\sigma_{\SBTB}\,\mathscr{T}_{\SBTB} + \mu_{\SBTB}\,\big|^3\,\bigr)$\\[2.5ex]
&$\leq$&$4\,\big|\,\sigma_{\SBTB}\,\big|^3\,E\bigl(\,\big|\,\mathscr{T}_{\SBTB}\,\big|^3\,\bigr)
+ 4\,\big|\,\mu_{\SBTB}\,\big|^3$
&(\index{H{\"o}lder's inequality!for finite sequences using length $\nu$}(\ref{EWD_0_1_05}) 
with $\nu = 2$, $p = 3$)\\[2.5ex]
&$\leq$&$4\,\Bigl(\,\dfrac{4}{3}\,\Bigr)^{3/2} E\bigl(\,\big|\,\mathscr{T}_{\SBTB}\,\big|^4\,\bigr)^{3/4} + 4$
&(H{\"o}lder's inequality\index{H{\"o}lder's inequality!for random variables})\\[2.5ex]
&$\leq$&$7\,\Bigl(\,97\,\dfrac{\delta_{\SBTB}}{n - l} + 3\,\Bigr)^{3/4} + 4$
&(cf. Proposition \ref{EWD_3_2_02}, \ref{EWD_3_2_02_BWd})).
\end{tabular}\\[2.8ex]
Next, we use the estimate\\[2ex]
\refstepcounter{DSBcount}
\label{EWD_3_7_30}
\text{\hspace*{-0.8ex}(\theDSBcount)}
\hspace*{2.9ex}
\begin{tabular}[t]{@{}r@{\hspace*{0.8ex}}c@{\hspace*{0.8ex}}l@{\hspace*{38.7ex}}r@{}}
$\dfrac{\delta_{\SBTB}}{n - l}$&$\leq$&$5^4\,\dfrac{\delta_{\bm\bar{A}}}{n - l}$&
(cf. Lemma \ref{EWD_3_3_10}; $\sigma_{\SBTB}^2 \geq \dfrac{2}{3}$)
\end{tabular}\\[2ex]
\hspace*{12.2ex}\begin{tabular}[t]{@{}r@{\hspace*{0.8ex}}c@{\hspace*{0.8ex}}l@{\hspace*{22.5ex}}r@{}}
\hspace*{5.5ex}&$\leq$&$10 \cdot 5^4 \cdot \dfrac{\delta_{\bm\bar{A}}}{n}$
&(cf. Remark \ref{EWD_3_4_24} $\Rightarrow$ $\dfrac{n}{n - l} \leq 10$)\\[2ex]
&$\leq$&$2 \cdot 5^5 \cdot \dfrac{\beta_{\bm\bar{A}}}{n}$
&($|\bm\bar{a}_{ij}| \leq 1$ for all $1 \leq i, j \leq n$ $\Rightarrow$ 
$\delta_{\bm\bar{A}} \leq \beta_{\bm\bar{A}}$)\\[2ex]
&$\leq$&$2 \cdot 5^8 \cdot \dfrac{\beta_{A}}{n}$&(cf. Lemma \ref{EWD_3_4_05}; $\sigma_{\!A'}^2 \geq \dfrac{2}{3}$).
\end{tabular}\\[1.2ex]
Because of (\ref{EWD_3_7_11}) and $\epsilon_{0} \leq 1$, we now conclude that 
\mbox{\rule[-2.5ex]{0ex}{3.5ex}$\dfrac{\beta_{A}}{n} \leq 1$} and thus
$E\bigl(\,\big|\,T_{\SBTB}\,\big|^3\,\bigr) \leq c_{10}$.
Therefore, a repeated application of H{\"o}lder's inequality (\ref{EWD_0_1_05}) with $\nu = 2$ and 
$p = 2$ ($p = 3$)\index{H{\"o}lder's inequality!for finite sequences using length $\nu$}
shows that $A_{1} \leq c_{11}$ and $A_{2} \leq c_{12}$.\\[2.8ex]
Furthermore, we obtain for all $w \in \{\,0, 1, 2\,\}$ and $b \in \mathbb{R}$\\[2ex]
\refstepcounter{DSBcount}
\label{EWD_3_7_31}
\text{\hspace*{-0.8ex}(\theDSBcount)}
\hspace*{2.8ex}
\begin{tabular}[t]{@{}c@{\hspace*{0.8ex}}l@{\hspace*{0.8ex}}l@{}}
\hspace*{1.9ex}&$E\Bigl(\,\big|\,T_{\SBTB} + a\,\big|^{w}\,
1_{\textstyle{(\,z,\,z\, + 2\,\lambda\,]}}\bigl(\,T_{\SBTB} + b\,\bigr)\,\Bigr)$\\[2.5ex]
$=$&$E\Bigl(\,\big|\,\sigma_{\SBTB}\,\mathscr{T}_{\SBTB} + \mu_{\SBTB} + a\,\big|^{w}\,
1_{\textstyle{(\,z,\,z\, + 2\,\lambda\,]}}\bigl(\,\sigma_{\SBTB}\,\mathscr{T}_{\SBTB} 
+ \mu_{\SBTB} + b\,\bigr)\,\Bigr)$\\[2.5ex]
$\leq$&$3^{w-1}\,\big|\,\sigma_{\SBTB}\,\big|^{w}\,E\Bigl(\,\big|\,\mathscr{T}_{\SBTB}\,\big|^{w}\,
1_{\textstyle{(\,z,\,z\, + 2\,\lambda\,]}}\bigl(\,\sigma_{\SBTB}\,\mathscr{T}_{\SBTB} 
+ \mu_{\SBTB} + b\,\bigr)\,\Bigr)$\\[2.5ex]
&\hspace*{15.8ex}$+\ 3^{w-1}\,\Bigl(\,\big|\,\mu_{\SBTB}\,\big|^w + \big|\,a\,\big|^{w}\,\Bigr)\,E\Bigl(\,
1_{\textstyle{(\,z,\,z\, + 2\,\lambda\,]}}\bigl(\,\sigma_{\SBTB}\,\mathscr{T}_{\SBTB} 
+ \mu_{\SBTB} + b\,\bigr)\,\Bigr)$\\[2.5ex]
$\leq$&$3^{w-1}\,\big|\,\sigma_{\SBTB}\,\big|^{w}\,\Bigl(\,c_{13}(w)\,\dfrac{\beta_{A}}{n} 
+ \dfrac{1}{\sqrt{2\,\pi}}\,\dfrac{2\,\lambda}{\sigma_{\SBTB}}\,\Bigr)$\\[2.5ex]
&\hspace*{27.1ex}$+\ 3^{w-1}\,\Bigl(\,\big|\,\mu_{\SBTB}\,\big|^w + \big|\,a\,\big|^{w}\,\Bigr)\, 
\,\Bigl(\,c_{13}(0)\,\dfrac{\beta_{A}}{n} + \dfrac{1}{\sqrt{2\,\pi}}\,
\dfrac{2\,\lambda}{\sigma_{\SBTB}}\,\Bigr)$\\[2.5ex]
$\leq$&$c_{14}(w)\,\Bigl\{\,\dfrac{\beta_{A}}{n} + \lambda\,\Bigr\}$\\[3ex]
$\leq$&$c_{15}(w)\,\lambda$.
\end{tabular}\\[2.5ex]
The first inequality is clear for $w = 0$ and $w = 1$. For $w = 2$ it 
follows from H{\"o}lder's inequality (\ref{EWD_0_1_05}) with $\nu = 3$ and 
$p = 2$.\index{H{\"o}lder's inequality!for finite sequences using length $\nu$}\\[1ex]
The second inequality is a consequence of Theorem \ref{EWD_3_6_01}, Theorem \ref{EWD_3_1_03} 
and argumentations similar to those in (\ref{EWD_3_7_30}). In (\ref{EWD_3_7_30})
we estimate a term using the exponent $e = 4$ (i.e. $\delta_{\SBTB}$). However, similar estimates also hold
for the exponents $e = 3$ and $e = w + 4$.\\[1.5ex]
The third inequality uses again $|\,\mu_{\SBTB}\,| \leq 1$ and $\dfrac{2}{3} \leq \sigma_{\SBTB}^2 \leq \dfrac{4}{3}$
(cf. (\ref{EWD_3_7_10})), and $|\,a\,| \leq 8$.\\[1ex]
The fourth and last inequality follows from \mbox{\rule[-2.5ex]{0ex}{3.5ex}$\dfrac{\beta_{A}}{n} \leq G_{\!A}$}
(cf. Lemma \ref{EWD_3_1_18}, \ref{EWD_3_1_18_BWb})) and 
$G_{\!A} \leq \sqrt{\dfrac{4}{3}}\,\lambda$ (cf. (\ref{EWD_3_7_13})).\\[1ex]
Therefore, we also have $A_{3} \leq c_{16}$ and $A_{5} \leq c_{17}$.\\[2.8ex]
Finally, we estimate $A_{4}$. For the case $m = 0$ we get\\[2ex]
\refstepcounter{DSBcount}
\label{EWD_3_7_32}
\text{\hspace*{-0.8ex}(\theDSBcount)}
\hspace*{2.8ex}
\begin{tabular}[t]{@{}c@{\hspace*{0.8ex}}l@{\hspace*{-30.4ex}}r@{}}
&$\displaystyle{\bigg|\,E\Bigl(\,\bigl(\,1_{\textstyle{(\,z + \lambda,\,z\, + 2\,\lambda\,]}} - 
1_{\textstyle{(\,z,\,z\, + \lambda\,]}}\,\bigr)\,
\bigl(\,T_{\SBTB} + a + s\,p + s\,t\,q\,\bigr)\,\Bigr)\,\bigg|}$\\[3ex]
$\leq$&$\displaystyle{\sup\limits_{\varrho\, \in\, \mathbb{R}}\,
\,\big|\,\Delta_{\lambda}^{2}F_{\SBTB}(\varrho)\,\big|}$
&($\varrho = z - a - s\,p - s\,t\,q$)\\[3ex]
$\leq$&$\mathfrak{L}(0)\,\bigl(\,G_{\!A}^2 + \lambda^2\,\bigr)$
&(cf. condition (\ref{EWD_3_7_03}), fulfilled due to (\ref{EWD_3_7_12})).
\end{tabular}\\[2.5ex]
Thus, because of $G_{\!A}^2 \leq \dfrac{4}{3}\,\lambda^2$ (cf. (\ref{EWD_3_7_13})), we have\\[2ex]
\hspace*{12.1ex}$A_{4} \leq \dfrac{7}{3}\,\mathfrak{L}(0)$.\\[3.5ex]
For the case $m = 1$, on the other hand,\\[2ex]
\refstepcounter{DSBcount}
\label{EWD_3_7_33}
\text{\hspace*{-0.8ex}(\theDSBcount)}
\hspace*{5.4ex}
$\displaystyle{\bigg|\,E\Bigl(\,T_{\SBTB}\,
\bigl(\,1_{\textstyle{(\,z + \lambda,\,z\, + 2\,\lambda\,]}} - 
1_{\textstyle{(\,z,\,z\, + \lambda\,]}}\,\bigr)\,
\bigl(\,T_{\SBTB} + a + s\,p + s\,t\,q\,\bigr)\,\Bigr)\,\bigg|}$\\[2.5ex]
has to be estimated. To do this, we define again $\varrho = z - a - s\,p - s\,t\,q$
and obtain then by applying integration by parts\index{integration!by parts for Lebesgue-Stieltjes integrals} 
for Lebesgue-Stieltjes integrals (cf. e.g. \cite{athreya2006measure}, Theorem 5.2.3)\\[3ex]
\hspace*{12.1ex}\begin{tabular}[t]{@{}c@{\hspace*{0.8ex}}l@{}}
$=$&$\displaystyle{\bigg|\,\int\limits_{\textstyle{\bigl(\,\varrho + \lambda,\,\varrho + 2\,\lambda\,\bigr]}}
\!\!x\,\,\,dF_{\SBTB}(x)\,\, - 
\int\limits_{\textstyle{\bigl(\,\varrho,\,\varrho + \lambda\,\bigr]}}
\!\!x\,\,\,dF_{\SBTB}(x)\,\bigg|}$\\[6ex]
$=$&$\displaystyle{\bigg|\,\Delta_{\lambda}^{2}\bigl(\,\varrho\,F_{\SBTB}(\varrho)\,\bigr) - 
\int\limits_{\varrho + \lambda}^{\varrho + 2\,\lambda} F_{\SBTB}(x)\,dx + 
\int\limits_{\varrho}^{\varrho + \lambda} F_{\SBTB}(x)\,dx
\,\bigg|}$
\end{tabular}\\[5ex]
\hspace*{12.1ex}\begin{tabular}[t]{@{}c@{\hspace*{0.8ex}}l@{\hspace*{-1ex}}l@{}}
$\leq$&$\displaystyle{\bigg|\,\varrho\,\Delta_{\lambda}^{2} F_{\SBTB}(\varrho) 
+ 2\,\lambda\,F_{\SBTB}(\,\varrho + 2\,\lambda\,) - 2\,\lambda\,F_{\SBTB}(\,\varrho + \lambda\,)}$\\[1.8ex] 
&$\displaystyle{-\ \int\limits_{0}^{\lambda} F_{\SBTB}(\,v + \varrho + \lambda\,)\,dv + 
\int\limits_{0}^{\lambda} F_{\SBTB}(\,v + \varrho\,)\,dv
\,\bigg|}$&(subst.\index{integration!by substitution} $x = v + \varrho + t\,\lambda$, $t = 0\,; 1$)\\[4.5ex]
$\leq$&$\displaystyle{\big|\,\varrho\,\big|\,\Big|\,\Delta_{\lambda}^{2}F_{\SBTB}(\varrho)\,\Big|
+ 2\,\lambda\,\Big|\,F_{\SBTB}(\,\varrho + 2\,\lambda\,) - F_{\SBTB}(\,\varrho + \lambda\,)\,\Big|}$\\[1.8ex] 
&$\displaystyle{+\ \int\limits_{0}^{\lambda} \Big|\,F_{\SBTB}(\,v + \varrho + \lambda\,) 
- F_{\SBTB}(\,v + \varrho\,)\,\Big|\,dv}$.
\end{tabular}\\[3ex]
If we now use condition (\ref{EWD_3_7_03}) for the first summand and Theorem \ref{EWD_3_1_03} 
as in (\ref{EWD_3_7_31}) for the other two summands, we further get\\[2.2ex]
\hspace*{12.1ex}\begin{tabular}[t]{@{}c@{\hspace*{0.8ex}}l@{\hspace*{10ex}}l@{}}
$\leq$&$\mathfrak{L}(1)\,\bigl(\,G_{\!A}^2 + \lambda^2\,\bigr)
+ 3\,\lambda\,\Bigl\{\,c_{13}(0)\,\dfrac{\beta_{A}}{n} 
+ \dfrac{1}{\sqrt{2\,\pi}}\,
\dfrac{\lambda}{\sigma_{\SBTB}}\,\Bigr\}$\\[2.5ex]
$\leq$&$\Bigl(\,\dfrac{7}{3}\,\mathfrak{L}(1) + c_{18}\,\Bigr)\,\lambda^2$.
\end{tabular}\\[1.2ex]
For the last inequality, again $\dfrac{2}{3} \leq \sigma_{\SBTB}^2$ and 
\mbox{\rule[-2.5ex]{0ex}{3.5ex}$\dfrac{\beta_{A}}{n} \leq G_{\!A}$}
(cf. Lemma \ref{EWD_3_1_18}, \ref{EWD_3_1_18_BWb})) and 
$G_{\!A} \leq \sqrt{\dfrac{4}{3}}\,\lambda$ (cf. (\ref{EWD_3_7_13})) were applied.\\[2.8ex]
By using the triangle inequality, the estimates of (\ref{EWD_3_7_33}) and (\ref{EWD_3_7_32}) 
(with $\mathfrak{L}(1)$ instead of $\mathfrak{L}(0)$), and 
$|\,a\,| \leq 8$, we obtain $A_{4} \leq c_{19}\,\mathfrak{L}(1) + c_{20}$
for the case $m = 1$.\hspace*{1ex}\hfill$\Box$\vspace*{2ex}

\section[Edgeworth expansions of second order]{Edgeworth expansions of second order}\label{EWD_Kap3_Sec8}

The aim of this section is to prove Theorem \ref{EWD_3_1_13}.\\[2.8ex] 
To do this, we choose analogously to the proof of Theorem \ref{EWD_3_7_01} (cf. (\ref{EWD_3_7_10})) 
a $0 < \epsilon_{0} \leq 1$ and an $n_{0} \geq 18$ such that that for all matrices $A$ satisfying
$\sigma_{\!A} > 0$, $\beta_{A} \leq \epsilon_{0}\,n$ and $n \geq n_{0}$ holds:\\[2.5ex]
\refstepcounter{DSBcount}
\label{EWD_3_8_01}
\text{\hspace*{-0.8ex}(\theDSBcount)}
\hspace*{4ex}
\begin{tabular}{@{}l@{}}
$\displaystyle{
\left\{
\begin{array}{l@{\hspace*{3ex}}l@{}}
\big|\,\sigma_{\!A'}^2 - 1\,\big| \leq \dfrac{1}{3},\\[2.5ex]
|\bm\bar{a}_{ij}| \leq 1&\text{for}\ 1 \leq i, j \leq n,\\[2.5ex]
|\,\mu_{Q}\,| \leq 1&\text{for}\ Q \in N(16, \bm\hat{A}) \cup N(16,\bm\bar{A}),\\[2.8ex]
\big|\,\sigma_{\!Q}^2 - 1\,\big| \leq \dfrac{1}{3}&\text{for}\ Q \in N(16, \bm\hat{A}) \cup N(16,\bm\bar{A}).
\end{array}  \right.}$
\end{tabular}\\[2.5ex]
As in the argumentation to (\ref{EWD_3_7_10}), this follows from Corollary \ref{EWD_3_4_02}, 
Lemma \ref{EWD_3_4_05} and Corol-
\linebreak
lary \ref{EWD_3_3_07}.\\[2.8ex]
In this section, therefore, we now assume without loss of generality that\\[2.5ex]
\refstepcounter{DSBcount}
\label{EWD_3_8_02}
\text{\hspace*{-0.8ex}(\theDSBcount)}
\hspace*{4ex}
$\dfrac{\beta_{A}}{n} \leq E_{\!A} \leq \epsilon_{0}$
\hspace*{2ex}and\hspace*{2ex}
$n \geq n_{0}$\hfill(cf. Lemma \ref{EWD_3_1_18}, \ref{EWD_3_1_18_BWb})).\\[2.5ex]
If, on the other hand, $E_{\!A} > \epsilon_{0}$ or $n < n_{0}$, we first proceed
analogously to (\ref{EWD_3_7_07}) and (\ref{EWD_3_7_08}):\\[2.5ex] 
\refstepcounter{DSBcount}
\label{EWD_3_8_03}
\text{\hspace*{-0.8ex}(\theDSBcount)}
\hspace*{4ex}
$E_{\!A} > \epsilon_{0}$
\hspace*{2ex}$\Rightarrow$\hspace*{2ex}
$1 \leq \dfrac{E_{\!A}}{\epsilon_{0}}$\hspace*{7ex}and\\[2.5ex]
\refstepcounter{DSBcount}
\label{EWD_3_8_04}
\text{\hspace*{-0.8ex}(\theDSBcount)}
\hspace*{4ex}
$n < n_{0}$
\hspace*{2ex}$\Rightarrow$\hspace*{2ex}
$1 \leq \dfrac{n_{0}}{n} \leq \dfrac{n_{0}}{n}\,4\,\delta_{\!A} = 4\,n_{0}\,D_{\!A}^2 
\leq 4\,n_{0}\,E_{\!A}^2$
\hfill(cf. Lemma \ref{EWD_3_1_18}, \ref{EWD_3_1_18_BWa}), \ref{EWD_3_1_18_BWb})).\\[2.5ex]
From this we get similar to (\ref{EWD_3_7_09}) for matrices $A$ 
with $E_{\!A} > \epsilon_{0}$ or $n < n_{0}$:\\[2.5ex]
\refstepcounter{DSBcount}
\label{EWD_3_8_05}
\text{\hspace*{-0.8ex}(\theDSBcount)}
\hspace*{4ex}
\begin{tabular}[t]{@{}c@{\hspace*{0.8ex}}l@{\hspace*{-14.4ex}}r@{}}
&$||\,\mathscr{F}_{\!A} - e_{2,A}\,||$\\[2ex]
$\leq$&$1 + \dfrac{1}{15}\,\dfrac{\beta_{A}}{n}
+ \dfrac{1}{2}\,\dfrac{\delta_{A}}{n}$&
(cf. Lemma \ref{EWD_3_1_18}, \ref{EWD_3_1_18_BWh}))\\[2.5ex] 
$\leq$&$1 + \dfrac{1}{15}\,E_{\!A}
+ \dfrac{1}{2}\,E_{\!A}^2$&
(cf. Lemma \ref{EWD_3_1_18}, \ref{EWD_3_1_18_BWb}))\\[2.5ex]
$\leq$&$\biggl[\,\max\,\Bigl\{\,\dfrac{1}{\epsilon_{0}^3},\,8\,\sqrt{\rule[0ex]{0ex}{2ex}{n_{0}^3}}\,\Bigr\} 
+ \dfrac{1}{15}\,\max\,\Bigl\{\,\dfrac{1}{\epsilon_{0}^2},\,4\,n_{0}\,\Bigr\} 
+ \dfrac{1}{2}\,\max\,\Bigl\{\,\dfrac{1}{\epsilon_{0}},\,2\,\sqrt{\rule[0ex]{0ex}{2ex}{n_{0}}}\,\Bigr\}\,
\biggr]\,E_{\!A}^3$.
\end{tabular}\\[3ex]
A sufficiently large $\mathcal{K}_{6}$ in Theorem \ref{EWD_3_1_13} therefore 
ensures the validity of this theorem in the cases
$E_{\!A} > \epsilon_{0}$ and $n < n_{0}$ as well.\\[2.8ex]
Due to its length, we divide the following proof of Theorem \ref{EWD_3_1_13} again into four parts, 
\linebreak
similar to the proof of Theorem \ref{EWD_3_7_01} (or Theorem \ref{EWD_1_1_03} for iid 
random variables).\vspace*{2ex}

\subsection[Transition to smooth functions and replacing $A$ with $\bar{A}$]
{Transition to smooth functions and replacing {\boldmath $A$} with {\boldmath $\bar{A}$}}\label{EWD_Kap3_Sec8_1}

In order to apply Stein's method in the next section \ref{EWD_Kap3_Sec8_2}, we start analogously to
subsection \ref{EWD_Kap3_Sec7_2_1} (or section \ref{EWD_Kap1_Sec2}) and replace in the expression\\[2ex]
\hspace*{12.1ex}$\displaystyle{\sup\limits_{z\, \in\, \mathbb{R}}\, 
\bigg|\,E\Bigl(\,1_{(- \infty,\, z\,]}(\,\mathscr{T}_{A}\,)\,\Bigr) 
- \int\limits_{\mathbb{R}} 1_{(\,- \infty,\,z\,]}(x)\,e_{2,A}'(x)\,dx\,\bigg|}$\\[2.5ex]
the functions $1_{(\,- \infty,\,z\,]}(x)$, $z \in \mathbb{R}$, by smoother functions.
The functions $r_{z}^{\eta}$, $z \in \mathbb{R}$ (cf. (\ref{EWD_2_3_03})) with\\[2ex]
\refstepcounter{DSBcount}
\label{EWD_3_8_06}
\text{\hspace*{-0.8ex}(\theDSBcount)}
\hspace*{4ex}
$\eta = E_{\!A}$\\[2ex]
prove to be suitable. We have\\[4ex]
\refstepcounter{DSBcount}
\label{EWD_3_8_07}
\textbf{\hspace*{-0.8ex}\theDSBcount\ Proposition}\\[2ex]
\hspace*{12.1ex}$\displaystyle{||\,\mathscr{F}_{\!A} - e_{2,A}\,||
\leq 3\,\sup\limits_{z\, \in\, \mathbb{R}}\, 
\bigg|\,E\Bigl(\,r_{z}^{\eta}(\,\mathscr{T}_{A}\,)\,\Bigr) 
- \int\limits_{\mathbb{R}} r_{z}^{\eta}(x)\,e_{2,A}'(x)\,dx\,\bigg| 
+ \Bigl(\,C_{1}\,\mathcal{C}_{2} + C_{2}\,\Bigr)\,E_{\!A}^3}$.\\[3.5ex]
\textbf{Proof:}\\[0.8ex]
In this proof, let $r_{z} = r_{z}^{\eta}$. Then for all $z \in \mathbb{R}$ holds\\[2ex]
\hspace*{12.1ex}\begin{tabular}[t]{@{}c@{\hspace*{0.8ex}}l@{\hspace*{-13.4ex}}l@{}}
&$\displaystyle{\dfrac{1}{2}\,\Big|\,\bigl(\,\mathscr{F}_{\!A} - e_{2,A}\,\bigr)(\,z + 2\,\eta\,) + 
\bigl(\,\mathscr{F}_{\!A} - e_{2,A}\,\bigr)(\,z + \eta\,)\,\Big|}$\\[3ex]
$=$&$\displaystyle{\Big|\,-\ \bigl(\,\mathscr{F}_{\!A} - e_{2,A}\,\bigr)(\,z\,)
- \dfrac{3}{2}\,\Delta_{\eta}^{1}\bigl(\,\mathscr{F}_{\!A} - e_{2,A}\,\bigr)(\,z\,)
- \dfrac{1}{2}\,\Delta_{\eta}^{2}\bigl(\,\mathscr{F}_{\!A} - e_{2,A}\,\bigr)(\,z\,)\,\Big|}$\\[2.5ex]
$=$&$\displaystyle{\bigg|\,\int\limits_{z}^{z\, +\, 3\,\eta}
r_{z}'(x)\,P_{\eta}^{2}\bigl(\,x;\,z,\,\mathscr{F}_{\!A} - e_{2,A}\,\bigr)\,dx\,\bigg|}$
&(cf. Lemma \ref{EWD_2_3_13}, \ref{EWD_2_3_13_BWc}))\\[4ex]
$\leq$&$\displaystyle{\bigg|\,\int\limits_{z}^{z\, +\, 3\,\eta}
r_{z}'(x)\,\Bigl(\,P_{\eta}^{2}\bigl(\,x;\,z,\,\mathscr{F}_{\!A}\,\bigr) - \mathscr{F}_{\!A}(x)\,\Bigr)\,dx\,\bigg|}$\\[4ex]
&$\displaystyle{+\ \bigg|\,\int\limits_{z}^{z\, +\, 3\,\eta}
r_{z}'(x)\,\Bigl(\,\mathscr{F}_{\!A}(x) - e_{2,A}(x)\,\Bigr)\,dx\,\bigg|}$\\[4ex]
&$\displaystyle{+\ \bigg|\,\int\limits_{z}^{z\, +\, 3\,\eta}
r_{z}'(x)\,\Bigl(\,e_{2,A}(x) - P_{\eta}^{2}\bigl(\,x;\,z,\,e_{2,A}\,\bigr)\,\Bigr)\,dx\,\bigg|}$\\[4ex]
$\leq$&$\displaystyle{\sup\limits_{z\, \in\, \mathbb{R}}\,
\bigg|\,\int\limits_{\mathbb{R}} r_{z}'(x)\,\mathscr{F}_{\!A}(x)\,dx - 
\int\limits_{\mathbb{R}} r_{z}'(x)\,e_{2,A}(x)\,dx\,\bigg|}$\\[4ex]
&$\displaystyle{+\ \Bigl(\,28\,\mathcal{C}_{2} + \dfrac{4}{3} \cdot 28\,||\,e'''_{2,A}\,||\,\Bigr)\,\eta^3\,
\int\limits_{z}^{z\, +\, 3\,\eta} \big|\,r_{z}'(x)\,\big|\,dx}$.
\end{tabular}\\[3.5ex]
For the last inequality, we used condition (\ref{EWD_3_1_14}) together with (\ref{EWD_3_8_06}) 
and Lemma \ref{EWD_2_4_02}, \ref{EWD_2_4_02_BWa}) for $F = e_{2,A}$.
Furthermore, we exploited that $r_{z}'$ is identical to zero 
outside of $[\,z,\,z + 3\,\eta\,]$ (cf. (\ref{EWD_2_3_14})).\\[2.8ex]
The next step is\\[2.5ex]
\hspace*{12.1ex}\begin{tabular}[t]{@{}c@{\hspace*{0.8ex}}l@{\hspace*{-10.5ex}}l@{}}
$\leq$&$\displaystyle{\sup\limits_{z\, \in\, \mathbb{R}}\,
\bigg|\,E\Bigl(\,r_{z}(\,\mathscr{T}_{A}\,)\,\Bigr) - 
\int\limits_{\mathbb{R}} r_{z}(x)\,e'_{2,A}(x)\,dx\,\bigg|
+ \Bigl(\,28\,\mathcal{C}_{2} + 229\,\Bigr)\,\eta^3} =: R$
\end{tabular}\\[3.5ex]
We obtained the first summand by an
integration by parts\index{integration!by parts for Lebesgue-Stieltjes integrals} 
for Lebesgue-Stieltjes integrals (cf. e.g. \cite{athreya2006measure}, Theorem 5.2.3)
using $r_{z}(+ \infty) = 0$, $\mathscr{F}_{\!A}(- \infty) = 0$ 
and $e_{2,A}(- \infty) = 0$.\\[2.8ex]
For the second summand, we used the first integral of Lemma \ref{EWD_2_3_13}, \ref{EWD_2_3_13_BWc}) again, 
taking into account the fact that $r_{z}'$ is negative on the entire interval $[\,z,\,z + 3\,\eta\,]$.
In addition, we utilized $||\,e'''_{2,A}\,|| \leq 6,13$ (cf. Lemma \ref{EWD_3_1_18}, 
\ref{EWD_3_1_18_BWi}) and \ref{EWD_3_1_18_BWb}), and $E_{\!A} \leq \epsilon_{0} \leq 1$).\\[2.8ex]
Using the formula $\Delta_{y}^{k+1}G(x) = \Delta_{y}^{k}G(x + y) - \Delta_{y}^{k}G(x)$,
we now get from\\[2ex]
\hspace*{12.1ex}$\displaystyle{\dfrac{1}{2}\,\Big|\,\bigl(\,\mathscr{F}_{\!A} - e_{2,A}\,\bigr)(\,z + \eta\,) + 
\bigl(\,\mathscr{F}_{\!A} - e_{2,A}\,\bigr)(\,z\,)\,\Big| \leq R}$
\hspace*{4ex}for all $z \in \mathbb{R}$\\[2.5ex]
the inequalities\\[2.5ex]
\hspace*{12.1ex}$\displaystyle{\dfrac{1}{2}\,\Big|\,\Delta_{\eta}^{1}\bigl(\,\mathscr{F}_{\!A} - e_{2,A}\,\bigr)
(\,z + \eta\,) + 
\Delta_{\eta}^{1}\bigl(\,\mathscr{F}_{\!A} - e_{2,A}\,\bigr)(\,z\,)\,\Big| \leq 2\,R}$
\hspace*{4ex}and\\[3ex]
\hspace*{12.1ex}$\displaystyle{\dfrac{1}{2}\,\Big|\,\Delta_{\eta}^{2}\bigl(\,\mathscr{F}_{\!A} - e_{2,A}\,\bigr)
(\,z + \eta\,) + 
\Delta_{\eta}^{2}\bigl(\,\mathscr{F}_{\!A} - e_{2,A}\,\bigr)(\,z\,)\,\Big| \leq 4\,R}$
\hspace*{4ex}for all $z \in \mathbb{R}$.\\[3ex]
If we then apply the formula\\[2ex]
\hspace*{12.1ex}\begin{tabular}[t]{@{}l@{\hspace*{0.8ex}}c@{\hspace*{0.8ex}}l@{}}
$\big|\,G(x)\,\big|$&$=$&$\dfrac{1}{2}\,\big|\,G(x) + G(x + \eta) + G(x) - G(x + \eta)\,\big|$\\[2ex]
&$\leq$&$\dfrac{1}{2}\,\big|\,G(x + \eta) + G(x)\,\big| + 
\dfrac{1}{2}\,\big|\,G(x + \eta) - G(x)\,\big|$
\end{tabular}\\[2.5ex]
successively to the functions $G = \mathscr{F}_{\!A} - e_{2,A}$, 
$G = \Delta_{\eta}^{1}\bigl(\,\mathscr{F}_{\!A} - e_{2,A}\,\bigr)$
and $G =$ $\Delta_{\eta}^{2}\bigl(\,\mathscr{F}_{\!A} - e_{2,A}\,\bigr)$, we obtain\\[2ex]
\hspace*{12.1ex}\begin{tabular}[t]{@{}l@{\hspace*{0.8ex}}c@{\hspace*{0.8ex}}l@{}}
$\Big|\,\bigl(\,\mathscr{F}_{\!A} - e_{2,A}\,\bigr)(\,z\,)\,\Big|$&$\leq$&
$R + \dfrac{1}{2}\,\Big|\,\Delta_{\eta}^{1}\bigl(\,\mathscr{F}_{\!A} - e_{2,A}\,\bigr)(\,z\,)\,\Big|$\\[2.5ex]
&$\leq$&
$2\,R + \dfrac{1}{4}\,\Big|\,\Delta_{\eta}^{2}\bigl(\,\mathscr{F}_{\!A} - e_{2,A}\,\bigr)(\,z\,)\,\Big|$\\[2.5ex]
&$\leq$&
$3\,R + \dfrac{1}{8}\,\Big|\,\Delta_{\eta}^{3}\bigl(\,\mathscr{F}_{\!A} - e_{2,A}\,\bigr)(\,z\,)\,\Big|$\\[2.5ex]
&$\leq$&
$3\,R + \dfrac{1}{8}\,\Big|\,\Delta_{\eta}^{3}\,\mathscr{F}_{\!A}(z)\,\Big| +
\dfrac{1}{8}\,\Big|\,\Delta_{\eta}^{3}\,e_{2,A}(z)\,\Big|$.
\end{tabular}\\[2.8ex]
Using condition (\ref{EWD_3_1_14}) and Lemma \ref{EWD_2_4_01} for the second summand and
Lemma \ref{EWD_2_4_02}, \ref{EWD_2_4_02_BWb}) and $||\,e'''_{2,A}\,|| \leq 6,13$ 
(cf. reasoning above) for the third summand, we finally get the proposition.\hspace*{1ex}\hfill$\Box$\\[4ex]
For the further proof, we must now replace the matrix $A$ 
with the truncated matrix $\bm\bar{A}$. The following two propositions serve this purpose.
Before doing so, it should be noted that in the entire section \ref{EWD_Kap3_Sec8} 
we always use the definition\\[2.8ex]
\refstepcounter{DSBcount}
\label{EWD_3_8_08}
\text{\hspace*{-0.8ex}(\theDSBcount)}
\hspace*{4ex}
$\lambda = \dfrac{E_{\!A}}{\sigma_{\!A'}} = \dfrac{\eta}{\sigma_{\!A'}}$.\\[2.2ex]
We then get because of \mbox{\rule[0ex]{0ex}{3.8ex}$\dfrac{2}{3} \leq \sigma_{\!A'}^2 \leq \dfrac{4}{3}$}
(cf. (\ref{EWD_3_8_01}) and (\ref{EWD_3_8_02}))\\[2.4ex]
\refstepcounter{DSBcount}
\label{EWD_3_8_09}
\text{\hspace*{-0.8ex}(\theDSBcount)}
\hspace*{4ex}
$\sqrt{\dfrac{3}{4}}\,E_{\!A} \leq \lambda \leq \sqrt{\dfrac{3}{2}}\,E_{\!A}$.\\[2.5ex]
We point out the analogies between (\ref{EWD_3_8_08}) and (\ref{EWD_3_7_12}) as well as
(\ref{EWD_3_8_09}) and (\ref{EWD_3_7_13}).\\[2.8ex]
With these terms $\lambda$ and $\eta$ we now have:\\[4ex]
\refstepcounter{DSBcount}
\label{EWD_3_8_10}
\textbf{\hspace*{-0.8ex}\theDSBcount\ Proposition}\\[2ex]
\hspace*{12.1ex}\begin{tabular}[t]{@{}c@{\hspace*{0.8ex}}l@{}}
&$\displaystyle{\sup\limits_{z\, \in\, \mathbb{R}}\, 
\bigg|\,E\Bigl(\,r_{z}^{\eta}(\,\mathscr{T}_{A}\,)\,\Bigr) 
- \int\limits_{\mathbb{R}} r_{z}^{\eta}(x)\,e_{2,A}'(x)\,dx\,\bigg|}$\\[3ex]
$\leq$&$\displaystyle{\sup\limits_{z\, \in\, \mathbb{R}}\, 
\bigg|\,E\Bigl(\,r_{z}^{\lambda}(\,T_{\bm\bar{A}}\,)\,\Bigr) 
- \int\limits_{\mathbb{R}} r_{z}^{\lambda}(x)\,e_{2,\bm\bar{A}}'(x)\,dx\,\bigg| 
+ C\,E_{\!A}^3}$.
\end{tabular}\\[3.5ex]
\textbf{Proof:}\\[0.8ex]
The proposition follows from (\ref{EWD_3_4_15}) with $h = r_{z}^{\eta}$ taking into account\\[2ex] 
\refstepcounter{DSBcount}
\label{EWD_3_8_11}
\text{\hspace*{-0.8ex}(\theDSBcount)}
\hspace*{2.8ex}
$r_{z}^{\eta}(a\,x + b) = r_{(z-b)/a}^{\eta/a}(x)$\hspace*{3ex}for all 
$z, x, b \in \mathbb{R}$ and $a > 0$.\\[2.5ex] 
Since $a = \sigma_{\!A'}$ 
in (\ref{EWD_3_4_15}), we get $\dfrac{\eta}{a} = \lambda$.
\hspace*{1ex}\hfill$\Box$\\[4ex]
The conditions (\ref{EWD_3_1_14}) and (\ref{EWD_3_1_15}) must also be transferred to the matrix $\bm\bar{A}$.
In the following, a weaker version of (\ref{EWD_3_1_14}) is sufficient.\\[4ex]
\refstepcounter{DSBcount}
\label{EWD_3_8_12}
\textbf{\hspace*{-0.8ex}\theDSBcount\ Proposition}
\begin{enumerate}
\item\label{EWD_3_8_12_BWa}
From (\ref{EWD_3_1_14}) follows:\\[2.5ex]
\refstepcounter{DSBcount}
\label{EWD_3_8_13}
\text{\hspace*{-0.8ex}(\theDSBcount)}
\hspace*{2.8ex}
$\big|\,\Delta_{y}^{3}F_{\SBTB}(z)\,\big| 
\leq \Bigl(\,C_{1}\,\mathcal{C}_{2} + C_{2}\,\Bigr)\,\,\bigl(\,E_{\!A}^3 + y^3\,\bigr)$\\[1.8ex]
\hspace*{12.1ex}for all $z \in \mathbb{R}$, $0 \leq y \leq \lambda$ 
and $\BTB \in N(16,\bm\bar{A})$.
\item\label{EWD_3_8_12_BWb}
From (\ref{EWD_3_1_15}) follows:\\[2.5ex]
\refstepcounter{DSBcount}
\label{EWD_3_8_14}
\text{\hspace*{-0.8ex}(\theDSBcount)}
\hspace*{2.8ex}
$\bigl(\,|z| + 1\,\bigr)\,\big|\,\Delta_{y}^{2}F_{\SBTB}(z)\,\big| 
\leq \Bigl(\,C_{3}\,\mathcal{C}_{3} + C_{4}\,\Bigr)\,\,\bigl(\,E_{\!A}^2 + y^2\,\bigr)$\\[1.8ex]
\hspace*{12.1ex}for all $z \in \mathbb{R}$, $0 \leq y \leq \lambda$ 
and $\BTB \in N(16,\bm\bar{A})$.
\end{enumerate}
\vspace*{2ex}
\textbf{Proof:}
\begin{enumerate}
\item
This follows from Lemma \ref{EWD_2_4_01} and Proposition \ref{EWD_3_4_17} (with $r = 3$ and $k = 5$).
\item
We obtain this part from Proposition \ref{EWD_3_4_17} (with $r = 2$ and $k = 4$), 
Lemma \ref{EWD_3_1_18}, \ref{EWD_3_1_18_BWb}) (for $D_{\!A}^2 \leq E_{\!A}^2$) 
and Proposition \ref{EWD_3_4_19} (only for this Proposition \ref{EWD_3_4_19} we need
from (\ref{EWD_3_8_01}) the inequalities $|\,\mu_{Q}\,| \leq 1$ and 
$\big|\,\sigma_{\!Q}^2 - 1\,\big| \leq \dfrac{1}{3}$ for $Q \in N(16, \bm\hat{A})$).
\hspace*{1ex}\hfill$\Box$
\end{enumerate}
\vspace*{2ex}
Given (\ref{EWD_3_8_13}) and (\ref{EWD_3_8_14}), we now have to estimate the expression\\[2.5ex]
\refstepcounter{DSBcount}
\label{EWD_3_8_15}
\text{\hspace*{-0.8ex}(\theDSBcount)}
\hspace*{2.8ex}
$\displaystyle{\sup\limits_{z\, \in\, \mathbb{R}}\, 
\bigg|\,E\Bigl(\,r_{z}^{\lambda}(\,T_{\bm\bar{A}}\,)\,\Bigr) 
- \int\limits_{\mathbb{R}} r_{z}^{\lambda}(x)\,e_{2,\bm\bar{A}}'(x)\,dx\,\bigg|}$.\\[2.5ex]
In these estimates, we suppress the $\lambda$ of the functions $r_{z}^{\lambda}$ (cf. (\ref{EWD_2_3_03})) 
and $g_{z}^{\lambda}$ (cf. (\ref{EWD_2_3_04})) similar to (\ref{EWD_3_7_14}) 
in the last section \ref{EWD_Kap3_Sec7}. We set:\\[2ex]
\refstepcounter{DSBcount}
\label{EWD_3_8_16}
\text{\hspace*{-0.8ex}(\theDSBcount)}
\hspace*{2.8ex}
$r_{z} = r_{z}^{\lambda}$
\hspace*{2ex}and\hspace*{2ex}
$g_{z}  = g_{z}^{\lambda}$.\index{function!smooth!$r_{z}$, $g_{z}$}\vspace*{1.5ex}

\subsection[Application of Stein's method]{Application of Stein's method}\label{EWD_Kap3_Sec8_2}

In order to estimate (\ref{EWD_3_8_15}), we apply, as in the last section \ref{EWD_Kap3_Sec7} 
(see in particular (\ref{EWD_3_7_16})), Stein's equation\index{Stein's equation}\\[2ex]
\refstepcounter{DSBcount}
\label{EWD_3_8_17}
\text{\hspace*{-0.8ex}(\theDSBcount)}
\hspace*{2.8ex}
$g_{z}'(x) - x\,g_{z}(x) = r_{z}(x) - \Phi(r_{z})$
\hspace*{4ex}for $x \in \mathbb{R}$\\[2.5ex]
to the expression $E\bigl(\,r_{z}(\,T_{\bm\bar{A}}\,)\,\bigr) - \Phi\bigl(\,r_{z}\,\bigr)$ 
and then expand the functions $g_{z}$ and $g'_{z}$ according to Taylor's theorem\index{Taylor's theorem}.\\[2.8ex]
In these expansions, the statistics $T_{1}, \ldots, T_{5}$ (cf. (\ref{EWD_3_5_21})) come into play, 
which, as in the proof of Theorem \ref{EWD_3_7_01}, will always refer to the matrix $\bm\bar{A}$.
With this agreement we have\\[4ex]

\pagebreak

\refstepcounter{DSBcount}
\label{EWD_3_8_18}
\textbf{\theDSBcount\ Proposition}\\[0.8ex]
For every $z \in \mathbb{R}$,\\[2ex]
\hspace*{12.1ex}\begin{tabular}{@{}l@{}}
$\displaystyle{E\bigl(\,r_{z}(\,T_{\bm\bar{A}}\,)\,\bigr) - \Phi\bigl(\,r_{z}\,\bigr)
+ \dfrac{1}{2}\,E\bigl(\,T_{\bm\bar{A}}^3\,\bigr)\,E\bigl(\,T_{\bm\bar{A}}\,g'_{z}(T_{\bm\bar{A}})\,\bigr)}$\\[2.5ex]
$\displaystyle{+\ \biggl\{\,\dfrac{1}{6}\,\Bigl(\,E\bigl(\,T_{\bm\bar{A}}^4\,\bigr) - 3\,\Bigr) - 
\dfrac{1}{4}\,E\bigl(\,T_{\bm\bar{A}}^3\,\bigr)^2\,\biggr\}\,E\bigl(\,T_{\bm\bar{A}}\,g''_{z}(T_{\bm\bar{A}})\,\bigr)
= R(z)}$
\end{tabular}\\[2.5ex]
where\\[2ex]
\begin{tabular}{@{\hspace*{4ex}}l@{\hspace*{0.8ex}}c@{\hspace*{0.8ex}}l@{}}
$R(z)$&$=$&$\displaystyle{\biggl\{\,\dfrac{1}{6}\,\Bigl(\,E\bigl(\,T_{\bm\bar{A}}^4\,\bigr) - 3\,\Bigr) - 
\dfrac{1}{4}\,E\bigl(\,T_{\bm\bar{A}}^3\,\bigr)^2\,\biggr\}\,n\,
E\biggl(\,\bm\bar{a}_{I_{1}J_{1}}\,\Delta T_{4}\,\int\limits_{0}^{1}
\Bigl(\,g'''_{z}(\,T_{3} + \Delta T_{3} + t\,\Delta T_{4}\,)}$\\[2ex] 
&&\hspace*{63ex}$\displaystyle{-\ g'''_{z}(\,T_{3}\,)\,\Bigr)\,dt\,\biggr)}$\\[3.5ex]
&&$\displaystyle{+\ \dfrac{1}{2}\,E\bigl(\,T_{\bm\bar{A}}^3\,\bigr)\,n\,
E\biggl(\,\bm\bar{a}_{I_{1}J_{1}}\,\Delta T_{4}\,\Delta T_{3}\,\int\limits_{0}^{1}
\Bigl(\,g'''_{z}(\,T_{2} + \Delta T_{2} + t\,\Delta T_{3}\,) - g'''_{z}(\,T_{2}\,)\,\Bigr)\,dt\,\biggr)}$\\[3.5ex]
&&$\displaystyle{+\ \dfrac{1}{2}\,E\bigl(\,T_{\bm\bar{A}}^3\,\bigr)\,n\,
E\biggl(\,\bm\bar{a}_{I_{1}J_{1}}\,\bigl(\,\Delta T_{4}\,\bigr)^2\,\int\limits_{0}^{1} (1-t)\,
\Bigl(\,g'''_{z}(\,T_{3} + \Delta T_{3} + t\,\Delta T_{4}\,) - g'''_{z}(\,T_{3}\,)\,\Bigr)\,dt\,\biggr)}$\\[3.5ex]
&&$\displaystyle{-\ n\,
E\biggl(\,\bm\bar{a}_{I_{1}J_{1}}\,\dfrac{\bigl(\,\Delta T_{4}\,\bigr)^3}{2}\,\int\limits_{0}^{1} (1-t)^2\,
\Bigl(\,g'''_{z}(\,T_{3} + \Delta T_{3} + t\,\Delta T_{4}\,) - g'''_{z}(\,T_{3}\,)\,\Bigr)\,dt\,\biggr)}$\\[3.5ex]
&&$\displaystyle{-\ n\,
E\biggl(\,\bm\bar{a}_{I_{1}J_{1}}\,\dfrac{\bigl(\,\Delta T_{4}\,\bigr)^2}{2}\,\Delta T_{3}\int\limits_{0}^{1}
\Bigl(\,g'''_{z}(T_{2} + \Delta T_{2} + t\,\Delta T_{3}\,) - g'''_{z}(\,T_{2}\,)\,\Bigr)\,dt\,\biggr)}$\\[3.5ex]
&&$\displaystyle{-\ n\,
E\biggl(\,\bm\bar{a}_{I_{1}J_{1}}\,\Delta T_{4}\,\bigl(\,\Delta T_{3}\,\bigr)^2\,\int\limits_{0}^{1} (1-t)\,
\Bigl(\,g'''_{z}(\,T_{2} + \Delta T_{2} + t\,\Delta T_{3}\,) - g'''_{z}(\,T_{2}\,)\,\Bigr)\,dt\,\biggr)}$\\[3.5ex]
&&$\displaystyle{-\ n\,
E\biggl(\,\bm\bar{a}_{I_{1}J_{1}}\,\Delta T_{4}\,\Delta T_{3}\,\Delta T_{2}\,\int\limits_{0}^{1}
\Bigl(\,g'''_{z}(\,T_{1} + \Delta T_{1} + t\,\Delta T_{2}\,) - g'''_{z}(\,T_{1}\,)\,\Bigr)\,dt\,\biggr)}$.
\end{tabular}\\[3.5ex]
\textbf{Proof:}\\[0.8ex]
In the following, let $z \in \mathbb{R}$ be fixed and $r = r_{z}$, $g = g_{z}$. 
Using (\ref{EWD_3_8_17}), we then get\\[2ex]
\hspace*{12.1ex}$E\bigl(\,r(\,T_{\bm\bar{A}}\,)\,\bigr) - \Phi\bigl(\,r\,\bigr) =
E\bigl(\,g'(\,T_{\bm\bar{A}}\,)\,\bigr) - E\bigl(\,T_{\bm\bar{A}}\,g(\,T_{\bm\bar{A}}\,)\,\bigr)$.\\[2.5ex]
We will now look at the last term on its own. 
Because of Lemma \ref{EWD_3_5_25}, \ref{EWD_3_5_25_BWa}) we have\\[2ex]
\hspace*{12.1ex}$E\bigl(\,T_{\bm\bar{A}}\,g(\,T_{\bm\bar{A}}\,)\,\bigr) =
n\,E\bigl(\,\bm\bar{a}_{I_{1}J_{1}}\,g(\,T_{5}\,)\,\bigr)$.\\[2.5ex]
Next, we use Taylor's theorem\index{Taylor's theorem} and the sums $T_{5} = T_{4} + \Delta T_{4}$ and 
$T_{4} = T_{3} + \Delta T_{3}$ (cf. (\ref{EWD_3_5_21})). We expand $g$ about $T_{4}$ and obtain\\[2ex]
\refstepcounter{DSBcount}
\label{EWD_3_8_19}
\text{\hspace*{-0.8ex}(\theDSBcount)}
\hspace*{2.8ex}
\begin{tabular}{@{}c@{\hspace*{0.8ex}}l@{}}
$=$&$\displaystyle{n\,E\bigl(\,\bm\bar{a}_{I_{1}J_{1}}\,g(\,T_{4}\,)\,\bigr) + 
n\,E\Bigl(\,\bm\bar{a}_{I_{1}J_{1}}\,\Delta T_{4}\,g'(\,T_{4}\,)\,\Bigr) +
n\,E\Bigl(\,\bm\bar{a}_{I_{1}J_{1}}\,\dfrac{\bigl(\,\Delta T_{4}\,\bigr)^2}{2}\,g''(\,T_{4}\,)\,\Bigr)}$\\[2ex]
&$\displaystyle{+\ n\,E\biggl(\,\bm\bar{a}_{I_{1}J_{1}}\,\dfrac{\bigl(\,\Delta T_{4}\,\bigr)^3}{2}\,
\int\limits_{0}^{1} (1 - t)^2\,g'''(\,T_{4} + t\,\Delta T_{4}\,)\,dt\,\biggr)}$\\[3ex]
$=$&$\displaystyle{n\,E\bigl(\,\bm\bar{a}_{I_{1}J_{1}}\,g(\,T_{4}\,)\,\bigr) + 
n\,E\Bigl(\,\bm\bar{a}_{I_{1}J_{1}}\,\Delta T_{4}\,g'(\,T_{4}\,)\,\Bigr) +
n\,E\Bigl(\,\bm\bar{a}_{I_{1}J_{1}}\,\dfrac{\bigl(\,\Delta T_{4}\,\bigr)^2}{2}\,g''(\,T_{4}\,)\,\Bigr)}$\\[3ex]
&$\displaystyle{+\
n\,E\Bigl(\,\bm\bar{a}_{I_{1}J_{1}}\,\dfrac{\bigl(\,\Delta T_{4}\,\bigr)^3}{6}\,g'''(\,T_{3}\,)\,\Bigr)}$\\[2.5ex]
&$\displaystyle{+\ n\,E\biggl(\,\bm\bar{a}_{I_{1}J_{1}}\,\dfrac{\bigl(\,\Delta T_{4}\,\bigr)^3}{2}\,
\int\limits_{0}^{1} (1 - t)^2\,\Bigl(\,g'''(\,T_{3} + \Delta T_{3} + t\,\Delta T_{4}\,)
- g'''(\,T_{3}\,)\,\Bigr)\,dt\,\biggr)}$.
\end{tabular}\\[3.5ex]
For the first and fourth summand of (\ref{EWD_3_8_19}) we get, because of
Lemma \ref{EWD_3_5_23}, \ref{EWD_3_5_23_BWa}) and \ref{EWD_3_5_23_BWb}), and
Lemma \ref{EWD_3_5_26}, \ref{EWD_3_5_26_BWa}),\\[2ex]
\hspace*{12.1ex}$n\,E\bigl(\,\bm\bar{a}_{I_{1}J_{1}}\,g(\,T_{4}\,)\,\bigr) =
n\,E\bigl(\,\bm\bar{a}_{I_{1}J_{1}}\,\bigr)\,E\bigl(\,g(\,T_{4}\,)\,\bigr) = 0$,\\[2.5ex]
\hspace*{12.1ex}$n\,E\Bigl(\,\bm\bar{a}_{I_{1}J_{1}}\,\dfrac{\bigl(\,\Delta T_{4}\,\bigr)^3}{6}\,g'''(\,T_{3}\,)\,\Bigr)
= n\,E\Bigl(\,\bm\bar{a}_{I_{1}J_{1}}\,\dfrac{\bigl(\,\Delta T_{4}\,\bigr)^3}{6}\,\Bigr)
\,E\bigl(\,g'''(\,T_{3}\,)\,\bigr)$.\\[3ex]
In the third summand of (\ref{EWD_3_8_19}), however, we expand $g''$ about $T_{3}$ and use the sums 
$T_{4} = T_{3} + \Delta T_{3}$ and $T_{3} = T_{2} + \Delta T_{2}$ (cf. (\ref{EWD_3_5_21}))\\[2ex]
\hspace*{12.1ex}\begin{tabular}[t]{@{}c@{\hspace*{0.8ex}}l@{}}
&$\displaystyle{n\,E\Bigl(\,\bm\bar{a}_{I_{1}J_{1}}\,
\dfrac{\bigl(\,\Delta T_{4}\,\bigr)^2}{2}\,g''(\,T_{4}\,)\,\Bigr)}$\\[2.5ex]
$=$&$\displaystyle{n\,E\Bigl(\,\bm\bar{a}_{I_{1}J_{1}}\,
\dfrac{\bigl(\,\Delta T_{4}\,\bigr)^2}{2}\,g''(\,T_{3}\,)\,\Bigr)}$\\[2.5ex]
&$\displaystyle{+\ n\,E\biggl(\,\bm\bar{a}_{I_{1}J_{1}}\,\dfrac{\bigl(\,\Delta T_{4}\,\bigr)^2}{2}\,\Delta T_{3}\,
\int\limits_{0}^{1} g'''(\,T_{3} + t\,\Delta T_{3}\,)\,dt\,\biggr)}$\\[3.5ex]
$=$&$\displaystyle{n\,E\Bigl(\,\bm\bar{a}_{I_{1}J_{1}}\,
\dfrac{\bigl(\,\Delta T_{4}\,\bigr)^2}{2}\,g''(\,T_{3}\,)\,\Bigr)
+ n\,E\Bigl(\,\bm\bar{a}_{I_{1}J_{1}}\,
\dfrac{\bigl(\,\Delta T_{4}\,\bigr)^2}{2}\,\Delta T_{3}\,g'''(\,T_{2}\,)\,\Bigr)}$\\[2.5ex]
&$\displaystyle{+\ n\,
E\biggl(\,\bm\bar{a}_{I_{1}J_{1}}\,\dfrac{\bigl(\,\Delta T_{4}\,\bigr)^2}{2}\,\Delta T_{3}\int\limits_{0}^{1}
\Bigl(\,g'''(\,T_{2} + \Delta T_{2} + t\,\Delta T_{3}\,) - g'''(\,T_{2}\,)\,\Bigr)\,dt\,\biggr)}$
\end{tabular}\\[3.5ex]
\hspace*{12.1ex}\begin{tabular}[t]{@{}c@{\hspace*{0.8ex}}l@{}}
$=$&$\displaystyle{n\,E\Bigl(\,\bm\bar{a}_{I_{1}J_{1}}\,
\dfrac{\bigl(\,\Delta T_{4}\,\bigr)^2}{2}\,\Bigr)\,E\bigl(\,g''(\,T_{3}\,)\,\bigr)
+ n\,E\Bigl(\,\bm\bar{a}_{I_{1}J_{1}}\,
\dfrac{\bigl(\,\Delta T_{4}\,\bigr)^2}{2}\,\Delta T_{3}\,\Bigr)\,E\bigl(\,g'''(\,T_{2}\,)\,\bigr)}$\\[2.5ex]
&$\displaystyle{+\ n\,
E\biggl(\,\bm\bar{a}_{I_{1}J_{1}}\,\dfrac{\bigl(\,\Delta T_{4}\,\bigr)^2}{2}\,\Delta T_{3}\int\limits_{0}^{1}
\Bigl(\,g'''(\,T_{2} + \Delta T_{2} + t\,\Delta T_{3}\,) - g'''(\,T_{2}\,)\,\Bigr)\,dt\,\biggr)}$.
\end{tabular}\\[3ex]
For the last equation, Lemma \ref{EWD_3_5_23}, \ref{EWD_3_5_23_BWb}) and 
\ref{EWD_3_5_23_BWc}) were applied.\\[2.8ex]
In the second summand of (\ref{EWD_3_8_19}) we expand $g'$ also about $T_{3}$ 
and use the sums $T_{4} = T_{3} + \Delta T_{3}$ and 
$T_{3} = T_{2} + \Delta T_{2}$ (cf. (\ref{EWD_3_5_21})) again:\\[2ex]
\refstepcounter{DSBcount}
\label{EWD_3_8_20}
\text{\hspace*{-0.8ex}(\theDSBcount)}
\hspace*{2.8ex}
\begin{tabular}{@{}c@{\hspace*{0.8ex}}l@{}}
&$\displaystyle{n\,E\Bigl(\,\bm\bar{a}_{I_{1}J_{1}}\,\Delta T_{4}\,g'(\,T_{4}\,)\,\Bigr)}$\\[3ex]
$=$&$\displaystyle{n\,E\Bigl(\,\bm\bar{a}_{I_{1}J_{1}}\,\Delta T_{4}\,g'(\,T_{3}\,)\,\Bigr)
+ n\,E\Bigl(\,\bm\bar{a}_{I_{1}J_{1}}\,\Delta T_{4}\,\Delta T_{3}\,g''(\,T_{3}\,)\,\Bigr)}$\\[2ex]
&$\displaystyle{+\ n\,E\biggl(\,\bm\bar{a}_{I_{1}J_{1}}\,\Delta T_{4}\,\bigl(\,\Delta T_{3}\,\bigr)^2\,
\int\limits_{0}^{1} (1-t)\,g'''(\,T_{3} + t\,\Delta T_{3}\,)\,dt\,\biggr)}$\\[4ex]
$=$&$\displaystyle{n\,E\Bigl(\,\bm\bar{a}_{I_{1}J_{1}}\,\Delta T_{4}\,g'(\,T_{3}\,)\,\Bigr)
+ n\,E\Bigl(\,\bm\bar{a}_{I_{1}J_{1}}\,\Delta T_{4}\,\Delta T_{3}\,g''(\,T_{3}\,)\,\Bigr)}$\\[3ex]
&$\displaystyle{+\ n\,E\Bigl(\,\bm\bar{a}_{I_{1}J_{1}}\,\Delta T_{4}\,
\dfrac{\bigl(\,\Delta T_{3}\,\bigr)^2}{2}\,g'''(\,T_{2}\,)\,\Bigr)}$\\[2.5ex]
&$\displaystyle{+\ n\,E\biggl(\,\bm\bar{a}_{I_{1}J_{1}}\,\Delta T_{4}\,\bigl(\,\Delta T_{3}\,\bigr)^2\,
\int\limits_{0}^{1} (1-t)\,\Bigl(\,g'''(\,T_{2} + \Delta T_{2} + t\,\Delta T_{3}\,)
- g'''(T_{2})\Bigr)\,dt\,\biggr)}$.
\end{tabular}\\[3ex]
For the first and third summand of (\ref{EWD_3_8_20}) we get, because of Lemma \ref{EWD_3_5_23}, \ref{EWD_3_5_23_BWb}) 
and \ref{EWD_3_5_23_BWc}), and Lemma \ref{EWD_3_5_26}, \ref{EWD_3_5_26_BWb}),\\[2ex]
\hspace*{12.1ex}$n\,E\Bigl(\,\bm\bar{a}_{I_{1}J_{1}}\,\Delta T_{4}\,g'(\,T_{3}\,)\,\Bigr)
= n\,E\Bigl(\,\bm\bar{a}_{I_{1}J_{1}}\,\Delta T_{4}\,\Bigr)\,E\bigl(\,g'(\,T_{3}\,)\,\bigr)
= E\bigl(\,g'(\,T_{3}\,)\,\bigr)$,\\[2.5ex]
\hspace*{12.1ex}$n\,E\Bigl(\,\bm\bar{a}_{I_{1}J_{1}}\,\Delta T_{4}\,
\dfrac{\bigl(\,\Delta T_{3}\,\bigr)^2}{2}\,g'''(\,T_{2}\,)\,\Bigr)
= n\,E\Bigl(\,\bm\bar{a}_{I_{1}J_{1}}\,\Delta T_{4}\,
\dfrac{\bigl(\,\Delta T_{3}\,\bigr)^2}{2}\,\Bigr)\,E\bigl(\,g'''(\,T_{2}\,)\,\bigr)$.\\[3ex]
Finally, in the second summand of (\ref{EWD_3_8_20}), we expand $g''$ about $T_{2}$ 
and use the sums $T_{3} = T_{2} + \Delta T_{2}$ and 
$T_{2} = T_{1} + \Delta T_{1}$ (cf. (\ref{EWD_3_5_21}))\\[2ex]
\hspace*{12.1ex}\begin{tabular}[t]{@{}c@{\hspace*{0.8ex}}l@{}}
&$\displaystyle{n\,E\Bigl(\,\bm\bar{a}_{I_{1}J_{1}}\,\Delta T_{4}\,\Delta T_{3}\,g''(\,T_{3}\,)\,\Bigr)}$\\[3.2ex]
$=$&$\displaystyle{n\,E\Bigl(\,\bm\bar{a}_{I_{1}J_{1}}\,\Delta T_{4}\,\Delta T_{3}\,g''(\,T_{2}\,)\,\Bigr)}$
\end{tabular}\\[2.3ex]
\hspace*{12.1ex}\begin{tabular}[t]{@{}c@{\hspace*{0.8ex}}l@{}}
&$\displaystyle{+\ n\,E\biggl(\,\bm\bar{a}_{I_{1}J_{1}}\,\Delta T_{4}\,\Delta T_{3}\,\Delta T_{2}\,
\int\limits_{0}^{1} g'''(\,T_{2} + t\,\Delta T_{2}\,)\,dt\,\biggr)}$\\[4ex]
$=$&$\displaystyle{n\,E\Bigl(\,\bm\bar{a}_{I_{1}J_{1}}\,\Delta T_{4}\,\Delta T_{3}\,g''(\,T_{2}\,)\,\Bigr)
+ n\,E\Bigl(\,\bm\bar{a}_{I_{1}J_{1}}\,\Delta T_{4}\,\Delta T_{3}\,\Delta T_{2}\,g'''(\,T_{1}\,)\,\Bigr)}$\\[3ex]
&$\displaystyle{+\ n\,E\biggl(\,\bm\bar{a}_{I_{1}J_{1}}\,\Delta T_{4}\,\Delta T_{3}\,\Delta T_{2}\,
\int\limits_{0}^{1} \Bigl(\,g'''(\,T_{1} + \Delta T_{1} 
+ t\,\Delta T_{2}\,) - g'''(\,T_{1}\,)\,\Bigr)\,dt\,\biggr)}$\\[4ex]
$=$&$\displaystyle{n\,E\Bigl(\,\bm\bar{a}_{I_{1}J_{1}}\,\Delta T_{4}\,\Delta T_{3}\,\Bigr)\,
E\bigl(\,g''(\,T_{2}\,)\,\bigr)
+ n\,E\Bigl(\,\bm\bar{a}_{I_{1}J_{1}}\,\Delta T_{4}\,\Delta T_{3}\,\Delta T_{2}\,\Bigr)\,
E\bigl(\,g'''(\,T_{1}\,)\,\bigr)}$\\[3ex]
&$\displaystyle{+\ n\,E\biggl(\,\bm\bar{a}_{I_{1}J_{1}}\,\Delta T_{4}\,\Delta T_{3}\,\Delta T_{2}\,
\int\limits_{0}^{1} \Bigl(\,g'''(\,T_{1} + \Delta T_{1} + t\,\Delta T_{2}\,) - g'''(\,T_{1}\,)\,\Bigr)\,dt\,\biggr)}$.
\end{tabular}\\[3ex]
For the last equation, Lemma \ref{EWD_3_5_23}, \ref{EWD_3_5_23_BWc}) and 
\ref{EWD_3_5_23_BWd}) were applied.\\[2.8ex]
We summarize the previous results taking into account Lemma \ref{EWD_3_5_15}, \ref{EWD_3_5_15_BWa}):\\[2ex]
\hspace*{12.1ex}\begin{tabular}[t]{@{}c@{\hspace*{0.8ex}}l@{}}
&$E\bigl(\,T_{\bm\bar{A}}\,g(\,T_{\bm\bar{A}}\,)\,\bigr)$\\[2.5ex]
$=$&$E\bigl(\,g'(\,T_{\bm\bar{A}}\,)\,\bigr)
+ E\bigl(\,g''(\,T_{\bm\bar{A}}\,)\,\bigr)\,\biggl\{\,
n\,E\Bigl(\,\bm\bar{a}_{I_{1}J_{1}}\,\Delta T_{4}\,\Delta T_{3}\,\Bigr) +
n\,E\Bigl(\,\bm\bar{a}_{I_{1}J_{1}}\,
\dfrac{\bigl(\,\Delta T_{4}\,\bigr)^2}{2}\,\Bigr)\,\biggr\}$\\[2.5ex]
&$+\ E\bigl(\,g'''(\,T_{\bm\bar{A}}\,)\,\bigr)\,\biggl\{\,
n\,E\Bigl(\,\bm\bar{a}_{I_{1}J_{1}}\,\Delta T_{4}\,\Delta T_{3}\,\Delta T_{2}\,\Bigr) +
n\,E\Bigl(\,\bm\bar{a}_{I_{1}J_{1}}\,\Delta T_{4}\,
\dfrac{\bigl(\,\Delta T_{3}\,\bigr)^2}{2}\,\Bigr)$\\[2.5ex] 
&\hspace*{19ex}$+\ n\,E\Bigl(\,\bm\bar{a}_{I_{1}J_{1}}\,
\dfrac{\bigl(\,\Delta T_{4}\,\bigr)^2}{2}\,\Delta T_{3}\,\Bigr) +
n\,E\Bigl(\,\bm\bar{a}_{I_{1}J_{1}}\,\dfrac{\bigl(\,\Delta T_{4}\,\bigr)^3}{6}\,\Bigr)\,\biggr\}$\\[2.5ex]
&$+\ N(z)$,
\end{tabular}\\[2.5ex]
where $N(z)$ is the sum of the last four terms of $R(z)$ (without minus sign!).
An application of Lemma \ref{EWD_3_5_26}, \ref{EWD_3_5_26_BWc}) and \ref{EWD_3_5_26_BWd}) further yields\\[2.5ex]
\refstepcounter{DSBcount}
\label{EWD_3_8_21}
\text{\hspace*{-0.8ex}(\theDSBcount)}
\hspace*{2.8ex}
\begin{tabular}[t]{@{}l@{\hspace*{0.8ex}}c@{\hspace*{0.8ex}}l@{}}
$E\bigl(\,r(\,T_{\bm\bar{A}}\,)\,\bigr) - \Phi\bigl(\,r\,\bigr)$&$+$&$\dfrac{1}{2}\,E\bigl(\,T_{\bm\bar{A}}^3\,\bigr)
\,E\bigl(\,g''(\,T_{\bm\bar{A}}\,)\,\bigr)$\\[2.5ex]
&$+$&$\biggl\{\,\dfrac{1}{6}\,\bigl(\,E(T_{\bm\bar{A}}^4) - 3\,\bigr)\,\biggr\}\,E\bigl(\,g'''(\,T_{\bm\bar{A}}\,)\,\bigr)
= - N(z)$.
\end{tabular}\\[3ex]
The degree of the derivative of $g''$ in the term $E\bigl(\,g''(\,T_{\bm\bar{A}}\,)\,\bigr)$
and of $g'''$ in the term $E\bigl(\,g'''(\,T_{\bm\bar{A}}\,)\,\bigr)$ 
will now be reduced again.\\[2.8ex]
If we argue as for (\ref{EWD_3_7_18}) with $g'$ instead of $f$, we get\\[2ex]
\hspace*{12.1ex}\begin{tabular}[t]{@{}c@{\hspace*{0.8ex}}l@{}}
\hspace*{1.9ex}&$E\bigl(\,T_{\bm\bar{A}}\,g'(\,T_{\bm\bar{A}}\,)\,\bigr)$\\[2.5ex]
$=$&$E\bigl(\,g''(\,T_{\bm\bar{A}}\,)\,\bigr) 
+ E\bigl(\,g'''(\,T_{\bm\bar{A}}\,)\,\bigr)\,\dfrac{1}{2}\,E\bigl(\,T_{\bm\bar{A}}^3\,\bigr)$\\[2ex]
&$\displaystyle{+\ n\,
E\biggl(\,\bm\bar{a}_{I_{1}J_{1}}\,\Delta T_{4}\,\Delta T_{3}\int\limits_{0}^{1}
\Bigl(\,g'''(\,T_{2} + \Delta T_{2} + t\,\Delta T_{3}\,) - g'''(\,T_{2}\,)\,\Bigr)\,dt\,\biggr)}$
\end{tabular}\\
\hspace*{12.1ex}\begin{tabular}[t]{@{}c@{\hspace*{0.8ex}}l@{}}
&$\displaystyle{+\ n\,E\biggl(\,\bm\bar{a}_{I_{1}J_{1}}\,\bigl(\,\Delta T_{4}\,\bigr)^2\,
\int\limits_{0}^{1}\,(1 - t)\,\Bigl(\,g'''(\,T_{3} + \Delta T_{3} + t\,\Delta T_{4}\,)
- g'''(\,T_{3}\,)\,\Bigr)\,dt\,\biggr)}$.
\end{tabular}\\[3ex]
Next, we insert the formula just derived for $E\bigl(\,g''(\,T_{\bm\bar{A}}\,)\,\bigr)$ 
into the third term of (\ref{EWD_3_8_21}) and obtain\\[2.5ex]
\refstepcounter{DSBcount}
\label{EWD_3_8_22}
\text{\hspace*{-0.8ex}(\theDSBcount)}
\hspace*{2.8ex}
\begin{tabular}[t]{@{}c@{\hspace*{0.8ex}}l@{}}
$\dfrac{1}{2}\,E\bigl(\,T_{\bm\bar{A}}^3\,\bigr)
\,E\bigl(\,g''(\,T_{\bm\bar{A}}\,)\,\bigr) =$&
$\dfrac{1}{2}\,E\bigl(\,T_{\bm\bar{A}}^3\,\bigr)\,E\bigl(\,T_{\bm\bar{A}}\,g'(T_{\bm\bar{A}})\,\bigr)$\\[2.5ex] 
&$-\ \Bigl(\,\dfrac{1}{2}\,E\bigl(\,T_{\bm\bar{A}}^3\,\bigr)\,\Bigr)^2\,E\bigl(\,g'''(\,T_{\bm\bar{A}}\,)\,\bigr)
- P_{1}(z)$,
\end{tabular}\\[3ex]
where $P_{1}(z)$ is the sum of the second and third term of $R(z)$.\\[2.8ex] 
If we further argue as in (\ref{EWD_3_6_05}) with $g''$ instead of $d$ 
and $\bm\bar{A}$ instead of $A$, and use Lemma \ref{EWD_3_5_15}, \ref{EWD_3_5_15_BWa}) again, we get\\[2ex]
\hspace*{12.1ex}\begin{tabular}[t]{@{}c@{\hspace*{0.8ex}}l@{}}
&$E\bigl(\,T_{\bm\bar{A}}\,g''(T_{\bm\bar{A}})\,\bigr)$\\[2ex]
$=$&$\displaystyle{E\bigl(\,g'''(\,T_{\bm\bar{A}}\,)\,\bigr) 
+ n\,E\biggl(\,\bm\bar{a}_{I_{1}J_{1}}\,\Delta T_{4}\,\int\limits_{0}^{1}\,
\Bigl(\,g'''\bigl(\,T_{3} + \Delta T_{3} + t\,\Delta T_{4}\,\bigr)
- g'''(\,T_{3}\,)\,\Bigr)\,dt\,\biggr)}$.
\end{tabular}\\[3ex]
From this it finally follows for the two terms from (\ref{EWD_3_8_21}) and (\ref{EWD_3_8_22}), 
which contain the value $E\bigl(\,g'''(\,T_{\bm\bar{A}}\,)\,\bigr)$, that\\[2.5ex]
\hspace*{12.1ex}\begin{tabular}[t]{@{}c@{\hspace*{0.8ex}}l@{}}
&$\biggl\{\,\dfrac{1}{6}\,\Bigl(\,E\bigl(\,T_{\bm\bar{A}}^4\,\bigr) - 3\,\Bigr) - 
\dfrac{1}{4}\,E\bigl(\,T_{\bm\bar{A}}^3\,\bigr)^2\,\biggr\}\,E\bigl(\,g'''(\,T_{\bm\bar{A}}\,)\,\bigr)$\\[2.5ex]
$=$&$\biggl\{\,\dfrac{1}{6}\,\Bigl(\,E\bigl(\,T_{\bm\bar{A}}^4\,\bigr) - 3\,\Bigr) - 
\dfrac{1}{4}\,E\bigl(\,T_{\bm\bar{A}}^3\,\bigr)^2\,\biggr\}\,E\bigl(\,T_{\bm\bar{A}}\,g''(T_{\bm\bar{A}})\,\bigr)
- P_{2}(z)$,
\end{tabular}\\[3ex]
where $P_{2}(z)$ is the first term of $R(z)$.\\[2.8ex]
Overall, we have thus proved the assertion of Proposition \ref{EWD_3_8_18}.\hspace*{1ex}\hfill$\Box$\vspace*{2ex}

\subsection[Deducing the correct expansion]{Deducing the correct expansion}\label{EWD_Kap3_Sec8_3}

This subsection deals with the two terms\\[2.5ex]
\hspace*{5.5ex}$\dfrac{1}{2}\,E\bigl(\,T_{\bm\bar{A}}^3\,\bigr)\,E\bigl(\,T_{\bm\bar{A}}\,g'_{z}(T_{\bm\bar{A}})\,\bigr)$
\hspace*{2ex}and\hspace*{2ex}
$\biggl\{\,\dfrac{1}{6}\,\Bigl(\,E\bigl(\,T_{\bm\bar{A}}^4\,\bigr) - 3\,\Bigr) - 
\dfrac{1}{4}\,E\bigl(\,T_{\bm\bar{A}}^3\,\bigr)^2\,\biggr\}\,
E\bigl(\,T_{\bm\bar{A}}\,g''_{z}(T_{\bm\bar{A}})\,\bigr)$.\\[3ex]
We prove the following\\[4ex]
\refstepcounter{DSBcount}
\label{EWD_3_8_23}
\textbf{\hspace*{-0.8ex}\theDSBcount\ Proposition}\\[0.8ex]
Let $\lambda_{1,\bm\bar{A}}$ be defined according to (\ref{EWD_3_1_08}) 
and $\lambda_{2,\bm\bar{A}}$ according to (\ref{EWD_3_1_09}) with $\bm\bar{A}$ instead of $A$ or $\bm\hat{A}$. Then:
\begin{enumerate}
\item\label{EWD_3_8_23_BWa}
\hspace*{4ex}\begin{tabular}[t]{@{}c@{\hspace*{0.8ex}}l@{}}
$\displaystyle{\sup\limits_{z\, \in\, \mathbb{R}}}$&$\displaystyle{ 
\bigg|\,\dfrac{1}{2}\,E\bigl(\,T_{\bm\bar{A}}^3\,\bigr)\,E\bigl(\,T_{\bm\bar{A}}\,g_{z}'(\,T_{\bm\bar{A}}\,)\,\bigr)
- \dfrac{1}{2}\,\lambda_{1,\bm\bar{A}}
\,\Phi\bigl(\,x\,g_{z}'(x)\,\bigr)}$\\[3.5ex]
&$\displaystyle{\ +\ \dfrac{1}{12}\,\lambda_{1,\bm\bar{A}}^{2}\,
\int\limits_{\mathbb{R}} x\,g_{z}'(x)\,(\,3\,x - x^3\,)\,\psi(x)\,dx \,\bigg|
\leq \Bigl(\,C_{1}\,\mathcal{C}_{3} + C_{2}\,\Bigr)\,E_{\!A}^3}$.
\end{tabular}
\item\label{EWD_3_8_23_BWb}
\hspace*{4ex}\begin{tabular}[t]{@{}c@{\hspace*{0.8ex}}l@{}}
$\displaystyle{\sup\limits_{z\, \in\, \mathbb{R}}}$&$\displaystyle{ 
\bigg|\,\biggl\{\,\dfrac{1}{6}\,\Bigl(\,E\bigl(\,T_{\bm\bar{A}}^4\,\bigr) - 3\,\Bigr) - 
\dfrac{1}{4}\,E\bigl(\,T_{\bm\bar{A}}^3\,\bigr)^2\,\biggr\}\,E\bigl(\,T_{\bm\bar{A}}\,
g''_{z}(T_{\bm\bar{A}})\,\bigr)}$\\[3.5ex]
&$\displaystyle{\ -\ \biggl\{\,\dfrac{1}{6}\,\lambda_{2,\bm\bar{A}} - 
\dfrac{1}{4}\,\lambda_{1,\bm\bar{A}}^{2}\,\biggr\}\,
\Phi\bigl(\,x\,g_{z}''(x)\,\bigr)\,\bigg|
\leq \Bigl(\,C_{3}\,\mathcal{C}_{3} + C_{4}\,\Bigr)\,E_{\!A}^3}$.
\end{tabular}
\end{enumerate}
\vspace*{3.5ex}
Using Lemma \ref{EWD_2_2_05}, \ref{EWD_2_2_05_BWc}), \ref{EWD_2_2_05_BWd}) and \ref{EWD_2_2_05_BWe}), 
we obtain the correct Edgeworth expansion in Proposition \ref{EWD_3_8_18} except for an error of size 
$\Bigl(\,\bigl(\,C_{1} + C_{3}\,\bigr)\,\mathcal{C}_{3} +
\bigl(\,C_{2} + C_{4}\,\bigr)\,\Bigr)\,E_{\!A}^3$.\\[4ex]
\textbf{Proof of Proposition \ref{EWD_3_8_23}:}
\begin{enumerate}
\item
Let $z \in \mathbb{R}$ be fixed.
We apply the triangle inequality and Lemma \ref{EWD_3_1_18}, \ref{EWD_3_1_18_BWd}), and get\\[2ex]
\begin{tabular}[t]{@{}c@{\hspace*{0.8ex}}l@{}}
&$\displaystyle{\bigg|\,E\bigl(\,T_{\bm\bar{A}}^3\,\bigr)\,E\bigl(\,T_{\bm\bar{A}}\,g_{z}'(\,T_{\bm\bar{A}}\,)\,\bigr)
- \lambda_{1,\bm\bar{A}}
\,\Phi\bigl(\,x\,g_{z}'(x)\,\bigr)
+ \dfrac{1}{6}\,\lambda_{1,\bm\bar{A}}^{2}\,
\int\limits_{\mathbb{R}} x\,g_{z}'(x)\,(\,3\,x - x^3\,)\,\psi(x)\,dx \,\bigg|}$\\[3.5ex]
$=$&$\displaystyle{\bigg|\,E\bigl(\,T_{\bm\bar{A}}^3\,\bigr)\,E\bigl(\,T_{\bm\bar{A}}\,g_{z}'(\,T_{\bm\bar{A}}\,)\,\bigr)
- \lambda_{1,\bm\bar{A}}\,
\int\limits_{\mathbb{R}} x\,g_{z}'(x)\,e_{1,\bm\bar{A}}'(x)\,dx \,\bigg|}$\\[3.5ex]
$\leq$&$\displaystyle{\bigg|\,E\bigl(\,T_{\bm\bar{A}}^3\,\bigr) - \lambda_{1,\bm\bar{A}}\,\bigg|\, 
E\Bigl(\,\big|\,T_{\bm\bar{A}}\,g_{z}'(\,T_{\bm\bar{A}}\,)\,\big|\,\Bigr)}$\\[3.5ex] 
&$\displaystyle{+\
\dfrac{\beta_{\bm\bar{A}}}{n}\,\bigg|\,E\bigl(\,T_{\bm\bar{A}}\,g_{z}'(\,T_{\bm\bar{A}}\,)\,\bigr)  
- \int\limits_{\mathbb{R}} x\,g_{z}'(x)\,e_{1,\bm\bar{A}}'(x)\,dx \,\bigg|}$\\[3.5ex]
$=$&$A_{1} + A_{2}$.
\end{tabular}\\[3.5ex]
Because of Proposition \ref{EWD_3_2_02}, \ref{EWD_3_2_02_BWa})
and because of $|\,g_{z}'\,| \leq 1$ (cf. Corollary \ref{EWD_2_1_12}, \ref{EWD_2_1_12_BWa})), 
we obtain for $A_{1}$:\\[2ex]
\hspace*{12.1ex}\begin{tabular}[t]{@{}l@{\hspace*{0.8ex}}c@{\hspace*{0.8ex}}l@{\hspace*{22.2ex}}r@{}} 
$A_{1}$&$\leq$&$\dfrac{11}{n^2}\,\beta_{\bm\bar{A}}$\\[2.5ex]
&$\leq$&$\dfrac{c_{1}}{n}\,\dfrac{\beta_{A}}{n}$
&(cf. Lemma \ref{EWD_3_4_05}; $\sigma_{\!A'}^2 \geq \dfrac{2}{3}$)\\[2.5ex]
\hspace*{2.8ex}&$\leq$&$c_{2}\,\Bigl(\,\dfrac{\beta_{A}}{n}\,\Bigr)^3 \leq c_{2}\,E_{\!A}^3$
&(cf. Lemma \ref{EWD_3_1_18}, \ref{EWD_3_1_18_BWa}) and \ref{EWD_3_1_18_BWb})).
\end{tabular}\\[2.5ex]
To estimate $A_{2}$, the following remains to be shown due to $\dfrac{\beta_{\bm\bar{A}}}{n}\leq 125\,E_{\!A}$:\\[2ex]
\refstepcounter{DSBcount}
\label{EWD_3_8_24}
\text{\hspace*{-0.8ex}(\theDSBcount)}
\hspace*{2.8ex}
$\displaystyle{\bigg|\,\,E\bigl(\,T_{\bm\bar{A}}\,g_{z}'(\,T_{\bm\bar{A}}\,)\,\bigr)  
- \int\limits_{\mathbb{R}} x\,g_{z}'(x)\,e_{1,\bm\bar{A}}'(x)\,dx \,\bigg| \leq 
\Bigl(\,c_{3}\,\mathcal{C}_{3} + c_{4}\,\Bigr)\,E_{\!A}^2}$.\\[2.5ex]
For this proof we use (\ref{EWD_3_8_17}) and obtain $g_{z}'(x) = h_{1}(x) - h_{2}(x)$ with\\[2ex]
\refstepcounter{DSBcount}
\label{EWD_3_8_25}
\text{\hspace*{-0.8ex}(\theDSBcount)}
\hspace*{2.8ex}
$h_{1}(x) = r_{z}(x) \in \mathcal{H}$
\hspace*{2ex}and\hspace*{2ex}
$h_{2}(x) = -\,\,\bigl(\,x\,g_{z}(x) - \Phi(r_{z})\,\bigr) \in \mathcal{H}$\index{function!$\mathcal{H}$}\\[2.5ex]
(cf. Lemma \ref{EWD_2_1_13}, \ref{EWD_2_1_13_BWa})). 
Thus (\ref{EWD_3_8_24}) follows from Corollary \ref{EWD_3_7_05} with $m = 1$ and $G_{\!A} = E_{\!A}$.
The condition (\ref{EWD_3_7_03}) of Corollary \ref{EWD_3_7_05} is fulfilled due to (\ref{EWD_3_8_14}).
\item
Let $z \in \mathbb{R}$ be fixed.
A further application of the triangle inequality and of 
Lemma \ref{EWD_3_1_18}, \ref{EWD_3_1_18_BWd}) yields\\[2ex]
\refstepcounter{DSBcount}
\label{EWD_3_8_26}
\text{\hspace*{-0.8ex}(\theDSBcount)}
\hfill
\begin{tabular}[t]{@{}c@{\hspace*{0.8ex}}l@{}}
&$\displaystyle{ 
\bigg|\,\biggl\{\,\dfrac{1}{6}\,\Bigl(\,E\bigl(\,T_{\bm\bar{A}}^4\,\bigr) - 3\,\Bigr) - 
\dfrac{1}{4}\,E\bigl(\,T_{\bm\bar{A}}^3\,\bigr)^2\,\biggr\}\,E\bigl(\,T_{\bm\bar{A}}\,
g''_{z}(T_{\bm\bar{A}})\,\bigr)}$\\[3.5ex]
&$\displaystyle{-\ \biggl\{\,\dfrac{1}{6}\,\lambda_{2,\bm\bar{A}} - 
\dfrac{1}{4}\,\lambda_{1,\bm\bar{A}}^{2}\,\biggr\}\,
\Phi\bigl(\,x\,g_{z}''(x)\,\bigr)\,\bigg|}$\\[3.5ex]
$\leq$&$\displaystyle{ 
\bigg|\,\biggl[\,\Bigl\{\,\dfrac{1}{6}\,\Bigl(\,E\bigl(\,T_{\bm\bar{A}}^4\,\bigr) - 3\,\Bigr) -
\dfrac{1}{4}\,E\bigl(\,T_{\bm\bar{A}}^3\,\bigr)^2\,\Bigr\} - \Bigl\{\,
\dfrac{1}{6}\,\lambda_{2,\bm\bar{A}} - \dfrac{1}{4}\,\lambda_{1,\bm\bar{A}}^{2}\,\Bigr\}\,\biggr]
\,E\bigl(\,T_{\bm\bar{A}}\,g''_{z}(T_{\bm\bar{A}})\,\bigr)\,\bigg|}$\\[3.5ex]
&$\displaystyle{+\ \bigg|\,\biggl\{\,\dfrac{1}{6}\,\lambda_{2,\bm\bar{A}} - 
\dfrac{1}{4}\,\lambda_{1,\bm\bar{A}}^{2}\,\biggr\}\,
\Bigl(\,E\bigl(\,T_{\bm\bar{A}}\,
g''_{z}(T_{\bm\bar{A}})\,\bigr) - \Phi\bigl(\,x\,g_{z}''(x)\,\bigr)\,\Bigr)\,\bigg|}$\\[3.5ex]
$\leq$&$\displaystyle{\biggl(\,\dfrac{1}{6}\,\Big|\,E\bigl(\,T_{\bm\bar{A}}^4\,\bigr) - 3
- \lambda_{2,\bm\bar{A}}\,\Big|\,+\,\dfrac{1}{4}\,
\Big|\,E\bigl(\,T_{\bm\bar{A}}^3\,\bigr)^2 - \lambda_{1,\bm\bar{A}}^{2}\,\Big|\,\biggr)\,
E\bigl(\,\big|\,T_{\bm\bar{A}}\,g''_{z}(T_{\bm\bar{A}})\,\big|\,\bigr)}$\\[3.5ex]
&$\displaystyle{+\ \biggl(\,\dfrac{19}{6}\,\dfrac{\delta_{\bm\bar{A}}}{n} + \dfrac{1}{4}\,
\Bigl(\,\dfrac{\beta_{\bm\bar{A}}}{n}\,\Bigr)^2\,\biggr)\,
\Big|\,E\bigl(\,T_{\bm\bar{A}}\,g''_{z}(T_{\bm\bar{A}})\,\bigr) - \Phi\bigl(\,x\,g_{z}''(x)\,\bigr)\,\Big|}$.
\end{tabular}\\[3.5ex]
Because of Proposition \ref{EWD_3_2_02}, \ref{EWD_3_2_02_BWb}), we now have\\[2ex]
\refstepcounter{DSBcount}
\label{EWD_3_8_27}
\text{\hspace*{-0.8ex}(\theDSBcount)}
\hspace*{2.8ex}
$\dfrac{1}{6}\,\Big|\,E\bigl(\,T_{\bm\bar{A}}^4\,\bigr) - 3
- \lambda_{2,\bm\bar{A}}\,\Big| = \dfrac{1}{6}\,\Big|\,R_{4,\bm\bar{A}}\,\Big| \leq \dfrac{1}{6}\,
 \dfrac{336}{n^2}\,\delta_{\bm\bar{A}} = \dfrac{56}{n}\,\dfrac{\delta_{\bm\bar{A}}}{n}$.\\[2.5ex]
Furthermore, if we apply the inequality\\[1.5ex]
\hspace*{12.1ex}$|\,a^2 - b^2\,| = |\,a - b\,|\,|\,a + b\,| 
\leq |\,a - b\,|\,\bigl(\,|\,a\,| + |\,b\,|\,\bigr)$\hspace*{4ex}for $a, b \in \mathbb{R}$,\\[1.5ex]
and use Proposition \ref{EWD_3_2_02}, \ref{EWD_3_2_02_BWa}) and \ref{EWD_3_2_02_BWc}),
and once more Lemma \ref{EWD_3_1_18}, \ref{EWD_3_1_18_BWd}),
we get\\[2.2ex]
\refstepcounter{DSBcount}
\label{EWD_3_8_28}
\text{\hspace*{-0.8ex}(\theDSBcount)}
\hspace*{2.8ex}
\begin{tabular}[t]{@{}l@{\hspace*{0.8ex}}c@{\hspace*{0.8ex}}l@{}} 
$\dfrac{1}{4}\,\Big|\,E\bigl(\,T_{\bm\bar{A}}^3\,\bigr)^2 - \lambda_{1,\bm\bar{A}}^{2}\,\Big|$&
$\leq$&$\dfrac{1}{4}\,\Big|\,E\bigl(\,T_{\bm\bar{A}}^3\,\bigr) - \lambda_{1,\bm\bar{A}}\,\Big|\,
\Bigl(\,\big|\,E\bigl(\,T_{\bm\bar{A}}^3\,\bigr)\,\big| + \big|\,\lambda_{1,\bm\bar{A}}\,\big|\,\Bigr)$\\[3ex]
&$\leq$&$\dfrac{1}{4}\,\dfrac{11}{n^2}\,\beta_{\bm\bar{A}}\,
\Bigl(\,5\,\dfrac{\beta_{\bm\bar{A}}}{n} + \dfrac{\beta_{\bm\bar{A}}}{n}
\,\Bigr)$\\[3ex] 
&$\leq$&$\dfrac{17}{n}\,\Bigl(\,\dfrac{\beta_{\bm\bar{A}}}{n}\,\Bigr)^2$.
\end{tabular}\\[2.5ex]
We now insert (\ref{EWD_3_8_27}) and (\ref{EWD_3_8_28}) into (\ref{EWD_3_8_26}).
Because of Lemma \ref{EWD_3_4_05} and Lemma \ref{EWD_3_1_18}, \ref{EWD_3_1_18_BWb}),
we then obtain for (\ref{EWD_3_8_26}) further\\[2.2ex]
\refstepcounter{DSBcount}
\label{EWD_3_8_29}
\text{\hspace*{-0.8ex}(\theDSBcount)}
\hspace*{2.8ex}
\begin{tabular}[t]{@{}c@{\hspace*{0.8ex}}l@{}}
$\leq$&$\displaystyle{\biggl(\,\dfrac{56}{n}\,\dfrac{\delta_{\bm\bar{A}}}{n} + 
\dfrac{17}{n}\,\Bigl(\,\dfrac{\beta_{\bm\bar{A}}}{n}\,\Bigr)^2\,\biggr)\,
E\bigl(\,\big|\,T_{\bm\bar{A}}\,g''_{z}(T_{\bm\bar{A}})\,\big|\,\bigr)}$\\[3.5ex]
&$\displaystyle{+\ \biggl(\,\dfrac{19}{6}\,\dfrac{\delta_{\bm\bar{A}}}{n} + \dfrac{1}{4}\,
\Bigl(\,\dfrac{\beta_{\bm\bar{A}}}{n}\,\Bigr)^2\,\biggr)\,
\Big|\,E\bigl(\,T_{\bm\bar{A}}\,
g''_{z}(T_{\bm\bar{A}})\,\bigr) - \Phi\bigl(\,x\,g_{z}''(x)\,\bigr)\,\Big|}$\\[3.5ex]
$\leq$&$\displaystyle{\Bigl(\,\dfrac{c_{5}}{n}\,E_{\!A}^2\,\Bigr)\,
E\bigl(\,\big|\,T_{\bm\bar{A}}\,g''_{z}(T_{\bm\bar{A}})\,\big|\,\bigr)
+ \Bigl(\,c_{6}\,E_{\!A}^2\,\Bigr)\,\Big|\,E\bigl(\,T_{\bm\bar{A}}\,
g''_{z}(T_{\bm\bar{A}})\,\bigr) - \Phi\bigl(\,x\,g_{z}''(x)\,\bigr)\,\Big|}$\\[3ex]
$=$&$A_{3} + A_{4}$.
\end{tabular}\\[2.5ex]
To estimate $A_{3}$, $|\,g_{z}''\,|$ must be estimated. To do this, 
we differentiate the equation (\ref{EWD_3_8_17}) on both sides and obtain\\[2ex]
\refstepcounter{DSBcount}
\label{EWD_3_8_30}
\text{\hspace*{-0.8ex}(\theDSBcount)}
\hspace*{2.8ex}
$g_{z}''(x) = g_{z}(x) + x\,g_{z}'(x) + r_{z}'(x)$
\hspace*{4ex}for $x \in \mathbb{R}$.\\[2ex]
Due to Corollary \ref{EWD_2_1_12}, \ref{EWD_2_1_12_BWa}) and due to
inequality (\ref{EWD_2_3_18}) and (\ref{EWD_3_8_09}), we then get\\[2ex]
\hspace*{12.1ex}$|\,g_{z}''(x)\,| \leq \dfrac{\sqrt{2\pi}}{4} + |\,x\,| + \sqrt{\dfrac{3}{4}}\,\dfrac{1}{E_{\!A}}$
\hspace*{4ex}for $x \in \mathbb{R}$.\\[2ex]
From this follows\\[2ex]
\hspace*{12.1ex}$E\bigl(\,\big|\,T_{\bm\bar{A}}\,g''_{z}(T_{\bm\bar{A}})\,\big|\,\bigr)
\leq \Bigl(\,\dfrac{2}{3} + \dfrac{7}{8}\,\dfrac{1}{E_{\!A}}\,\Bigr)\,E\bigl(\,\big|\,T_{\bm\bar{A}}\,\big|\,\bigr)
+ E\bigl(\,\big|\,T_{\bm\bar{A}}\,\big|^2\,\bigr)
\leq  \dfrac{5}{3} + \dfrac{7}{8}\,\dfrac{1}{E_{\!A}}$\\[2.2ex]
and therefore\\[2ex]
\hspace*{12.1ex}\begin{tabular}[t]{@{}l@{\hspace*{0.8ex}}c@{\hspace*{0.8ex}}l@{\hspace*{12.3ex}}r@{}}
$A_{3}$&$\leq$&$\dfrac{c_{5}}{n}\,E_{\!A}^2\,\Bigl(\,\dfrac{5}{3} + \dfrac{7}{8}\,\dfrac{1}{E_{\!A}}\,\Bigr)$\\[3.5ex]
&$\leq$&$\dfrac{5}{3}\,\dfrac{c_{5}}{\sqrt{n}}\,E_{\!A}^2
+ \dfrac{7}{8}\,\dfrac{c_{5}}{n}\,E_{\!A} 
$&(since $\sqrt{n} \leq n$)\\[3.5ex]
&$\leq$&$\dfrac{5}{3}\,2^{5/6}\,c_{5}\,E_{\!A}^3
+ \dfrac{7}{8}\,2^{5/3}\,c_{5}\,E_{\!A}^3 
 \leq 6\,c_{5}\,E_{\!A}^3$&(cf. Lemma \ref{EWD_3_1_18}, \ref{EWD_3_1_18_BWa})).
\end{tabular}\\[2.5ex]
To estimate $A_{4}$, we show\\[2ex]
\refstepcounter{DSBcount}
\label{EWD_3_8_31}
\text{\hspace*{-0.8ex}(\theDSBcount)}
\hspace*{2.8ex}
$\Big|\,E\bigl(\,T_{\bm\bar{A}}\,
g''_{z}(T_{\bm\bar{A}})\,\bigr) - \Phi\bigl(\,x\,g_{z}''(x)\,\bigr)\,\Big|
\leq \Bigl(\,c_{7}\,\mathcal{C}_{3} + c_{8}\,\Bigr)\,E_{\!A}$.\\[2.5ex]
To do this, we use the following representation of $x\,g_{z}''(x)$, where $x \in \mathbb{R}$:\\[2.2ex]
$x\,g_{z}''(x) = x\,\Bigl(\,g_{z}(x) + x\,g_{z}'(x) + r_{z}'(x)\,\Bigr)$\hspace*{38.8ex}(cf. (\ref{EWD_3_8_30}))\\[2.5ex]
\hspace*{4.5ex}\begin{tabular}[t]{@{}c@{\hspace*{0.8ex}}l@{\hspace*{-12.9ex}}r@{}}
$=$&$x\,g_{z}(x) 
+ x^2\,\bigl(\,x\,g_{z}(x) + r_{z}(x) - \Phi(r_{z})\,\bigr) + x\,r_{z}'(x)$&(cf. (\ref{EWD_3_8_17}))\\[2.5ex]
$=$&$\Phi(r_{z}) - \Bigl(\,-\,\bigl(\,x\,g_{z}(x) - \Phi(r_{z})\,\bigr)\,\Bigr)
- x^2\,\Bigl(\,-\,\bigl(\,x\,g_{z}(x) - \Phi(r_{z})\,\bigr)\,\Bigr) + x^2\,r_{z}(x) + x\,r_{z}'(x)$.
\end{tabular}\\[3ex]
We note that as in (\ref{EWD_3_8_25}), $\Bigl(\,-\,\,\bigl(\,x\,g_{z}(x) - \Phi(r_{z})\,\bigr)\,\Bigr) \in \mathcal{H}$
and $r_{z}(x) \in \mathcal{H}$.\index{function!$\mathcal{H}$}\\[2.5ex]
Thus, due to Corollary \ref{EWD_3_6_27} for $k = 0$ and $k = 2$,
and due to (\ref{EWD_3_7_21}) and $\dfrac{\beta_{\bm\bar{A}}}{n} \leq 125\,E_{\!A}$,
it suffices to carry out the following estimate for the proof of (\ref{EWD_3_8_31}):\\[2ex]
\refstepcounter{DSBcount}
\label{EWD_3_8_32}
\text{\hspace*{-0.8ex}(\theDSBcount)}
\hspace*{2.5ex}
\begin{tabular}{@{}l@{\hspace*{0.8ex}}c@{\hspace*{0.8ex}}l@{}}
$\Big|\,E\bigl(\,T_{\bm\bar{A}}\,
r'_{z}(T_{\bm\bar{A}})\,\bigr) - \Phi\bigl(\,x\,r_{z}'(x)\,\bigr)\,\Big|$&
$\leq$&$\displaystyle{\Big|\,E\bigl(\,T_{\bm\bar{A}}\,
r'_{z}(T_{\bm\bar{A}})\,\bigr) - \int\limits_{\mathbb{R}} x\,r_{z}'(x)\,e_{1,\bm\bar{A}}'\,dx\,\Big|}$\\[4ex]
&&$\displaystyle{\ +\ \dfrac{\beta_{\bm\bar{A}}}{6\,n}\,\int\limits_{\mathbb{R}}  
\big|\,r_{z}'(x)\,\big|\,\Big|\,x\,(\,3\,x - x^3\,)\,\psi(x)\,\Big|\,dx}$.
\end{tabular}\\[2ex]
We now utilize some properties of $r'_{z}$ from Lemma \ref{EWD_2_3_13}, \ref{EWD_2_3_13_BWd}) 
and from (\ref{EWD_2_3_18}). Because of Corollary \ref{EWD_3_7_05} 
(with $m = 1$ and $G_{\!A} = E_{\!A}$), which we can use due to (\ref{EWD_3_8_14}), 
and because of Lemma \ref{EWD_2_2_08}, \ref{EWD_2_2_08_BWg}), we then further obtain\\[2.2ex]
\hspace*{12.1ex}\begin{tabular}[t]{@{}c@{\hspace*{0.8ex}}l@{}}
$\leq$&
$2\,\Bigl(\,\dfrac{3}{4}\,\dfrac{1}{\lambda}\,\Bigr)
\,\Bigl(\,c_{9}\,\mathcal{C}_{3} + c_{10}\,\Bigr)\,E_{\!A}^2 + 
\dfrac{\beta_{\bm\bar{A}}}{6\,n}\,\bigl(\,(\,z + 3\,\lambda\,) - z\,\bigr)
\,\Bigl(\,\dfrac{3}{4}\,\dfrac{1}{\lambda}\,\Bigr)\,\dfrac{1}{2}$\\[2.5ex]
$\leq$&$\Bigl(\,c_{11}\,\mathcal{C}_{3} + c_{12}\,\Bigr)\,E_{\!A}$.
\end{tabular}\\[1.7ex]
For the last inequality, we used (\ref{EWD_3_8_09}), i.e.
$\dfrac{3}{4}\,\dfrac{1}{\lambda} \leq \sqrt{\dfrac{3}{4}}\,\dfrac{1}{E_{\!A}}$, and
again $\dfrac{\beta_{\bm\bar{A}}}{n} \leq 125\,E_{\!A}$.\\[2ex]
In total, (\ref{EWD_3_8_31}) and thus part \ref{EWD_3_8_23_BWb}) is proven.\hspace*{1ex}\hfill$\Box$
\end{enumerate}

\subsection[The estimate of $R(z)$]{The estimate of {\boldmath $R(z)$}}\label{EWD_Kap3_Sec8_4}

In this subsection we finally show for the $R(z)$ from Proposition \ref{EWD_3_8_18}:\\[4ex]
\refstepcounter{DSBcount}
\label{EWD_3_8_33}
\textbf{\hspace*{-0.8ex}\theDSBcount\ Proposition}\\[0.8ex]
\hspace*{12.1ex}$\displaystyle{\sup\limits_{z\, \in\, \mathbb{R}}\,\big|\,R(z)\,\big|\, 
\leq\, \Bigl(\,C_{1}\,\mathcal{C}_{2} + C_{2}\,\mathcal{C}_{3} + C_{3}\,\Bigr)\,E_{\!A}^3}$.\\[3.5ex]
\textbf{Proof:}\\[0.8ex]
In the following, let $\underline{I} = (\,I_{1},\ldots,I_{16}\,)$,\index{random permutation!$\underline{I}$}
\index{random permutation!$I_{1},\ldots,I_{16}$}$\overline{I} = (\,I_{1},\ldots,I_{8}\,)$, 
$\bm{\breve}{I} = (\,I_{1},\ldots,I_{4}\,)$
\index{random permutation!$\overline{I}$, $\bm{\breve}{I}$, $\overline{J}$, $\bm{\breve}{J}$}and 
analogously $\underline{i} = (\,i_{1},\ldots,i_{16}\,)$ $\in M_{16}$,\index{random permutation!$\underline{i}$}
\index{random permutation!$i_{1},\ldots,i_{16}$}$\overline{i} = (\,i_{1},\ldots,i_{8}\,)$ $\in M_{8}$, 
$\bm{\breve}{i} = (\,i_{1},\ldots,i_{4}\,) \in M_{4}$.\index{random permutation!$\overline{i}$, $\bm{\breve}{i}$}
The symbols $\underline{J}$\index{random permutation!$\underline{J}$}, 
$\overline{J}$, $\bm{\breve}{J}$, 
\linebreak
$\underline{j}$\index{random permutation!$\underline{j}$}, 
$\overline{j}$ and $\bm{\breve}{j}$\index{random permutation!$\overline{j}$, $\bm{\breve}{j}$} 
are defined in the same way.\\[2.8ex]
To prove Proposition \ref{EWD_3_8_33}, we show that there exist constants $c_{1},\ldots,c_{9}$ such that\\[3.5ex]
\refstepcounter{DSBcount}
\label{EWD_3_8_34}
\text{\hspace*{-0.8ex}(\theDSBcount)}
\hspace*{2.8ex}
\begin{tabular}[t]{@{}l@{}}
$\displaystyle{\bigg|\,E\Bigl(\,g'''_{z}\bigl(\,T_{1} + \Delta T_{1} + t\,\Delta T_{2}\,\bigr)
- g'''_{z}\bigl(\,T_{1}\,\bigr)\,\Big|\,\underline{I} = 
\underline{i},\, \underline{J} = \underline{j}\,\Bigr)\,\bigg|}$\\[3ex]
$\displaystyle{\leq\ \Bigl(\,c_{1}\,\mathcal{C}_{2} + c_{2}\,\mathcal{C}_{3} + c_{3}\,\Bigr)\,
E\Bigl(\,\big|\,\Delta T_{1}\,\big| + \big|\,\Delta T_{2}\,\big|
\,\Big|\,\underline{I} = \underline{i},\, \underline{J} = \underline{j}\,\Bigr)}$
\end{tabular}\\[3ex]
\hspace*{12.1ex}\begin{tabular}[t]{@{}l@{}}
for all $z \in \mathbb{R}$, $0 \leq t \leq 1$ and $\underline{i}, \underline{j} \in M_{16}$,
which have the properties (w1), (w2), (w3)\index{random permutation!properties!(w1), (w2), (w3)}
\end{tabular}\\
\hspace*{12.1ex}\begin{tabular}[t]{@{}l@{}}
from (\ref{EWD_3_5_19}),
\end{tabular}\\[3.5ex]
\refstepcounter{DSBcount}
\label{EWD_3_8_35}
\text{\hspace*{-0.8ex}(\theDSBcount)}
\hspace*{2.8ex}
\begin{tabular}[t]{@{}l@{}}
$\displaystyle{\bigg|\,E\Bigl(\,g'''_{z}\bigl(\,T_{2} + \Delta T_{2} + t\,\Delta T_{3}\,\bigr)
- g'''_{z}\bigl(\,T_{2}\,\bigr)\,\Big|\,\overline{I} = \overline{i},\, \overline{J} = \overline{j}\,\Bigr)\,\bigg|}$\\[3ex]
$\displaystyle{\leq\ \Bigl(\,c_{4}\,\mathcal{C}_{2} + c_{5}\,\mathcal{C}_{3} + c_{6}\,\Bigr)\,
E\Bigl(\,\big|\,\Delta T_{2}\,\big| + \big|\,\Delta T_{3}\,\big|
\,\Big|\,\overline{I} = \overline{i},\, \overline{J} = \overline{j}\,\Bigr)}$\\[3ex]
for all $z \in \mathbb{R}$, $0 \leq t \leq 1$ and $\overline{i}, \overline{j} \in M_{8}$ 
with $(\overline{i},\, \overline{j}) \in M_{16}$,
\end{tabular}\\[3.5ex]
\refstepcounter{DSBcount}
\label{EWD_3_8_36}
\text{\hspace*{-0.8ex}(\theDSBcount)}
\hspace*{2.8ex}
\begin{tabular}[t]{@{}l@{}}
$\displaystyle{\bigg|\,E\Bigl(\,g'''_{z}\bigl(\,T_{3} + \Delta T_{3} + t\,\Delta T_{4}\,\bigr)
- g'''_{z}\bigl(\,T_{3}\,\bigr)\,\Big|\,
\bm{\breve}{I} = \bm{\breve}{i},\, 
\bm{\breve}{J} = \bm{\breve}{j}\,\Bigr)\,\bigg|}$\\[3ex]
$\displaystyle{\leq\ \Bigl(\,c_{7}\,\mathcal{C}_{2} + c_{8}\,\mathcal{C}_{3} + c_{9}\,\Bigr)\,
E\Bigl(\,\big|\,\Delta T_{3}\,\big| + \big|\,\Delta T_{4}\,\big|
\,\Big|\,
\bm{\breve}{I} = \bm{\breve}{i},\, 
\bm{\breve}{J} = \bm{\breve}{j}\,\Bigr)}$\\[3ex]
for all $z \in \mathbb{R}$, $0 \leq t \leq 1$ and 
$\bm{\breve}{i}, \bm{\breve}{j} \in M_{4}$ 
with $(\bm{\breve}{i},\, \bm{\breve}{j}) \in M_{8}$.
\end{tabular}\\[3.5ex]
Because of Lemma \ref{EWD_3_5_25}, \ref{EWD_3_5_25_BWc}) and Proposition \ref{EWD_3_2_02}, \ref{EWD_3_2_02_BWc})
and \ref{EWD_3_2_02_BWd}), and because of\\[2.5ex]
\hspace*{12.1ex}$\dfrac{\delta_{\bm\bar{A}}}{n}\,\dfrac{\beta_{\bm\bar{A}}}{n} \leq
c_{10}\,\dfrac{\delta_{A}}{n}\,\dfrac{\beta_{A}}{n} \leq c_{10}\,E_{\!A}^3$
\hspace*{2ex}and\hspace*{2ex}
$\Bigl(\,\dfrac{\beta_{\bm\bar{A}}}{n}\,\Bigr)^3 \leq
c_{11}\,\Bigl(\,\dfrac{\beta_{A}}{n}\,\Bigr)^3 \leq c_{11}\,E_{\!A}^3$\\[2.5ex]
(cf. Lemma \ref{EWD_3_4_05} and Lemma \ref{EWD_3_1_18}, \ref{EWD_3_1_18_BWb})),
Proposition \ref{EWD_3_8_33} then follows from (\ref{EWD_3_8_34}), (\ref{EWD_3_8_35}) and (\ref{EWD_3_8_36}).\\[2.8ex]
Since the proofs of (\ref{EWD_3_8_34}), (\ref{EWD_3_8_35}) and (\ref{EWD_3_8_36}) are very similar, 
we will only deal with the case (\ref{EWD_3_8_34}) as an example.\\[2.8ex]
For that we fix in the following the quantities $z \in \mathbb{R}$, $0 \leq t \leq 1$ 
and $\underline{i}, \underline{j} \in M_{16}$, which have the properties (w1), (w2), (w3) 
from (\ref{EWD_3_5_19}).\index{random permutation!properties!(w1), (w2), (w3)}\\[2.8ex]
Because of Lemma \ref{EWD_3_5_27}, \ref{EWD_3_5_27_BWc}) is $T_{1}$, given $\underline{I} = \underline{i}$ and 
$\underline{J} = \underline{j}$, distributed as\\[2ex]
\refstepcounter{DSBcount}
\label{EWD_3_8_37}
\text{\hspace*{-0.8ex}(\theDSBcount)}
\hspace*{2.8ex}
$\displaystyle{\sum\limits_{(i,\,j)\, \in\, S} \bm\bar{a}_{ij} + T_{\SBTB}}$,\\[2.5ex]
where\\[2ex]
\refstepcounter{DSBcount}
\label{EWD_3_8_38}
\text{\hspace*{-0.8ex}(\theDSBcount)}
\hspace*{2.8ex}
$S = \Bigl\{\,(\,i_{s},\,j_{s+8}\,)\,:\,1\, \leq\, s\, \leq\, 8\,\Bigr\}\,
\cup\, \Bigl\{\,(\,i_{s},\,j_{s-8}\,)\,:\,9\, \leq\, s\, \leq\, 16\,\Bigr\}$\\[2.5ex]
and\\[2ex]
\refstepcounter{DSBcount}
\label{EWD_3_8_39}
\text{\hspace*{-0.8ex}(\theDSBcount)}
\hspace*{2.8ex}
\begin{tabular}[t]{@{}l@{}}
$\BTB$ is the $(n - l) \times (n - l)-$matrix obtained from $\bm\bar{A}$ by \textbf{cancelling} the\\[0.7ex]
$l = \big|\,\{\,i_{1},\ldots, i_{16}\,\}\,\big|$ rows $i_{1},\ldots,i_{16}$ and the $l$ columns
$j_{1},\ldots,j_{16}$.
\end{tabular}\\[3ex]
With the abbreviations\\[2ex]
\hspace*{12.1ex}\begin{tabular}[t]{@{}l@{\hspace*{0.8ex}}c@{\hspace*{0.8ex}}l@{\hspace*{15.5ex}}l@{}}
$a$&$=$&$\displaystyle{\sum\limits_{(i,\,j)\, \in\, S} \bm\bar{a}_{ij}}$,\\[3.5ex]
$p$&$=$&$\displaystyle{E\Bigl(\,\Delta T_{1}\,\Big|\,\underline{I} = \underline{i}, \underline{J} = \underline{j}\,\Bigr)}$
&($\Delta T_{1} \in \sigma\bigl(\,\underline{I},\,\underline{J}\,\bigr)$,
cf. Lemma \ref{EWD_3_5_22}, \ref{EWD_3_5_22_BWa})),\\[2.5ex]
$q$&$=$&$\displaystyle{E\Bigl(\,\Delta T_{2}\,\Big|\,\underline{I} = \underline{i}, \underline{J} = \underline{j}\,\Bigr)}$
&($\Delta T_{2} \in \sigma\bigl(\,\underline{I},\,\underline{J}\,\bigr)$,
cf. Lemma \ref{EWD_3_5_22}, \ref{EWD_3_5_22_BWb}))
\end{tabular}\\[3ex]
remains to be shown:\\[2.5ex]
\refstepcounter{DSBcount}
\label{EWD_3_8_40}
\text{\hspace*{-0.8ex}(\theDSBcount)}
\hspace*{2.8ex}
\begin{tabular}[t]{@{}l@{}}
$\displaystyle{\bigg|\,E\Bigl(\,g_{z}'''\bigl(\,T_{\SBTB} + a + p + t\,q\,\bigr)
- g_{z}'''\bigl(\,T_{\SBTB} + a\,\bigr)\,\Bigr)\,\bigg|}$\\[3ex]
$\displaystyle{\leq\ \Bigl(\,c_{12}\,\mathcal{C}_{2} + c_{13}\,\mathcal{C}_{3} + c_{14}\,\Bigr)\,
\Bigl(\,|\,p\,| + |\,q\,|\,\Bigr)}$.
\end{tabular}\\[3ex]
If we estimate the left-hand side of the inequality (\ref{EWD_3_8_40}) using Lemma \ref{EWD_2_3_13}, 
\ref{EWD_2_3_13_BWb}) and \ref{EWD_2_3_13_BWe}) 
\linebreak
(with \mbox{\rule[0ex]{0ex}{2.8ex}$x = T_{\SBTB} + a$} and 
\mbox{\rule[0ex]{0ex}{2.8ex}$y = p + t\,q$)}, we get\\[2ex]
\hspace*{12.1ex}\begin{tabular}[t]{@{}c@{\hspace*{0.8ex}}l@{}}
&$\displaystyle{\bigg|\,E\Bigl(\,g_{z}'''\bigl(\,T_{\SBTB} + a + p + t\,q\,\bigr)
- g_{z}'''\bigl(\,T_{\SBTB} + a\,\bigr)\,\Bigr)\,\bigg|}$\\[3ex]
$\leq$&$\displaystyle{E\biggl(\,\bigg|\,\bigl(\,g_{z}''' - r_{z}''\,\bigr)\,\bigl(\,T_{\SBTB} + a + p + t\,q\,\bigr)
- \bigl(\,g_{z}''' - r_{z}''\,\bigr)\,\bigl(\,T_{\SBTB} + a\,\bigr)}$\\[3ex]
&$\displaystyle{-\ \bigl(\,T_{\SBTB} + a\,\bigr)\,\Bigl(\,r_{z}'\bigl(\,T_{\SBTB} + a + p + t\,q\,\bigr)
- r_{z}'\bigl(\,T_{\SBTB} + a\,\bigr)\,\Bigr)\,\bigg|\,\biggr)}$\\[3ex]
&$\displaystyle{+\ \bigg|\,E\Bigl(\,r_{z}''\bigl(\,T_{\SBTB} + a + p + t\,q\,\bigr)
- r_{z}''\bigl(\,T_{\SBTB} + a\,\bigr)\,\Bigr)\,\bigg|}$\\[3ex]
&$\displaystyle{+\ \bigg|\,E\biggl(\,\bigl(\,T_{\SBTB} + a\,\bigr)\,\Bigl(\,r_{z}'\bigl(\,T_{\SBTB} + a + p + t\,q\,\bigr)
- r_{z}'\bigl(\,T_{\SBTB} + a\,\bigr)\,\Bigr)\,\biggr)\,\bigg|}$\\[3ex]
$\leq$&$\displaystyle{c_{15}\,\Bigl(\,|\,p\,| + |\,q\,|\,\Bigr)
\,\Biggl\{\,1 + E\Bigl(\,\big|\,T_{\SBTB} + a\,\big|^{3}\,\Bigr) + \Bigl(\,|\,p\,| + |\,q\,|\,\Bigr)}$\\[4ex]
&$\displaystyle{+\ \dfrac{1}{\lambda}\,P\Bigl(\,z\, <\, T_{\SBTB} + a + p + t\,q\, \leq\, z + 3\,\lambda\,\Bigr)}$\\[3ex]
&$\displaystyle{+\ \dfrac{1}{\lambda}\,\int\limits_{0}^{1}\,E\Bigl(\,\Bigl[\,(\,T_{\SBTB} + a\,)^{2} + 2\,\Bigr]\,
1_{\textstyle{(\,z,\,z\, + 3\,\lambda\,]}}\bigl(\,T_{\SBTB} + a + s\,p + s\,t\,q\,\bigr)\,\Bigr)\,ds}$
\end{tabular}\\
\hspace*{12.1ex}\begin{tabular}[t]{@{}c@{\hspace*{0.8ex}}l@{}}
\hspace*{1.8ex}&$\displaystyle{+\ \dfrac{1}{\lambda^3}\,\int\limits_{0}^{1}\,\bigg|\,E\Bigl(\,
\bigl(\,1_{\textstyle{(\,z + 2\,\lambda,\,z\, + 3\,\lambda\,]}}
- 2 \cdot 1_{\textstyle{(\,z + \lambda,\,z\, + 2\,\lambda\,]}}
+ 1_{\textstyle{(\,z,\,z\, + \lambda\,]}}\,\bigr)}$\\[3ex]
&\hspace*{45ex}$\displaystyle{\bigl(\,T_{\SBTB} + a + s\,p + s\,t\,q\,\bigr)\,\Bigr)\,\bigg|\,ds}$\\[2ex]
&$\displaystyle{+\ \int\limits_{0}^{1}\,\bigg|\,
E\Bigl(\,\bigl(\,T_{\SBTB} + a\,\bigr)\,r_{z}''\bigl(\,T_{\SBTB} + a + s\,p + s\,t\,q\,\bigr)\,
\Bigr)\,\bigg|\,ds\,\Biggr\}}$\\[4.5ex]
$=$&$\displaystyle{c_{15}\,\Bigl(\,|\,p\,| + |\,q\,|\,\Bigr)
\,\Bigl(\,1 + A_{1} + A_{2} + A_{3} + A_{4} + A_{5} + A_{6}\,\Bigr)}$.
\end{tabular}\\[3.2ex]
As in the proof of Proposition \ref{EWD_3_7_22}, it therefore remains 
to estimate the $A_{i}$, $1 \leq i \leq 6$.\\[2.8ex]
Due to $|\bm\bar{a}_{ij}| \leq 1$ for all $1 \leq i, j \leq n$, we first obtain\\[2ex]
\hspace*{12.1ex}$|\,a\,| \leq 16$,
\hspace*{2ex}$|\,p\,| \leq 32$
\hspace*{2ex}and\hspace*{2ex}
$|\,q\,| \leq 16$,\\[2.5ex]
from which $A_{2} = |\,p\,| + |\,q\,| \leq 48$ follows.\\[2.8ex] 
If we additionally use the estimates (\ref{EWD_3_7_29}) and (\ref{EWD_3_7_30})
from the proof of Proposition \ref{EWD_3_7_22}, and if we use as in this 
proof \mbox{\rule[-1.5ex]{0ex}{5.5ex}$\dfrac{\beta_{A}}{n} \leq \epsilon_{0}$} (cf. (\ref{EWD_3_8_02})) 
and $\epsilon_{0} \leq 1$, then we get $A_{1} \leq c_{16}$.\\[1.2ex]
Moreover, if we proceed as in the proof of (\ref{EWD_3_7_31}) (with $G_{\!A} = E_{\!A}$,
$3\,\lambda$ instead of $2\,\lambda$, and $|\,a\,| \leq 16$ instead of $|\,a\,| \leq 8$)
and use \mbox{\rule[-3ex]{0ex}{7ex}$E_{\!A} \leq \sqrt{\dfrac{4}{3}}\,\lambda$} (cf. (\ref{EWD_3_8_09})), 
we also get $A_{3} \leq c_{17}$ and $A_{4} \leq c_{18}$.\\[2.8ex]
Furthermore, the relationship (\ref{EWD_3_8_13}) yields\\[2ex]
\hspace*{4ex}
$\displaystyle{\bigg|\,E\Bigl(\,
\bigl(\,1_{\textstyle{(\,z + 2\,\lambda,\,z\, + 3\,\lambda\,]}}
- 2 \cdot 1_{\textstyle{(\,z + \lambda,\,z\, + 2\,\lambda\,]}}
+ 1_{\textstyle{(\,z,\,z\, + \lambda\,]}}\,\bigr)\,\bigl(\,T_{\SBTB} + a 
+ s\,p + s\,t\,q\,\bigr)\,\Bigr)\,\bigg|}$\\[2.2ex]
\hspace*{12.1ex}\begin{tabular}[t]{@{}c@{\hspace*{0.8ex}}l@{}}
$\leq$&
$\displaystyle{\sup\limits_{\varrho\, \in\, \mathbb{R}}\,\,\big|\,\Delta_{\lambda}^{3}F_{\SBTB}(\varrho)\,\big|}$
\hspace*{40.5ex}($\varrho = z - a - s\,p - s\,t\,q$)\\[3.2ex]
$\leq$&
$\Bigl(\,c_{19}\,\mathcal{C}_{2} + c_{20}\,\Bigr)\,\bigl(\,E_{\!A}^3 + \lambda^3\,\bigr)$,
\end{tabular}\\[2.5ex]
so that due to $E_{\!A}^3 \leq \Bigl(\,\dfrac{4}{3}\,\Bigr)^{3/2}\,\lambda^3$ 
(cf. (\ref{EWD_3_8_09})) we obtain\\[2.2ex]
\hspace*{12.1ex}$A_{5} \leq 3\,\Bigl(\,c_{19}\,\mathcal{C}_{2} + c_{20}\,\Bigr)$.\\[2.5ex]
Finally, $A_{6}$ has to be estimated. Because of\\[2ex]
\refstepcounter{DSBcount}
\label{EWD_3_8_41}
\text{\hspace*{-0.8ex}(\theDSBcount)}
\hspace*{2.8ex}
$r''_{z}(x + b) = r''_{z-b}(x)$\hspace*{3ex}for all $z, x, b \in \mathbb{R}$
\hfill(cf. (\ref{EWD_3_8_11}) with $a = 1$),\\[2ex]
it is enough to show\\[2.2ex]
\refstepcounter{DSBcount}
\label{EWD_3_8_42}
\text{\hspace*{-0.8ex}(\theDSBcount)}
\hspace*{2.8ex}
$\displaystyle{\sup\limits_{\varrho\, \in\, \mathbb{R}}\,\Big|\,E\Bigl(\,\bigl(\,T_{\SBTB} + a\,\bigr)\,
r_{\varrho}''\bigl(\,T_{\SBTB}\,\bigr)\,\Bigr)\,\Big| \leq c_{21}\,\mathcal{C}_{3} + c_{22}}$.\\[2.5ex]
By using (\ref{EWD_2_3_15}) we get for all $\varrho \in \mathbb{R}$\\[2ex]
\hspace*{9.2ex}\begin{tabular}[t]{@{}c@{\hspace*{0.8ex}}l@{}}
&$\Big|\,E\Bigl(\,\bigl(\,T_{\SBTB} + a\,\bigr)\,r_{\varrho}''\bigl(\,T_{\SBTB}\,\bigr)\,\Bigr)\,\Big|$\\[2.5ex] 
$=$&$\displaystyle{\bigg|\,\dfrac{1}{\lambda^3}\,\biggl\{\,\int\limits_{\varrho}^{\varrho + \lambda}
\bigl(\,x + a\,\bigr)\,\bigl(\,\varrho - x\,\bigr)\,dF_{\SBTB}(x) + 2\, \int\limits_{\varrho + \lambda}^{\varrho + 2\,\lambda}
\bigl(\,x + a\,\bigr)\,\bigl(\,x - \varrho - \dfrac{3}{2}\,\lambda\,\bigr)\,dF_{\SBTB}(x)}$\\[4.5ex]
&$\displaystyle{+\ \int\limits_{\varrho + 2\,\lambda}^{\varrho + 3\,\lambda}
\bigl(\,x + a\,\bigr)\,\bigl(\,\varrho + 3\,\lambda - x\,\bigr)\,dF_{\SBTB}(x)\,\biggr\}\,\bigg|.}$
\end{tabular}\\[2.5ex]
Next, we apply
integration by parts\index{integration!by parts for Lebesgue-Stieltjes integrals} 
for Lebesgue-Stieltjes integrals (cf. e.g. \cite{athreya2006measure}, Theorem 5.2.3)
to all three integrals.
Because of\\[2.2ex]
\hspace*{3.3ex}\begin{tabular}[t]{@{}l@{}}
$\dfrac{d}{dx}\,\bigl(\,x + a\,\bigr)\,\bigl(\,\varrho - x\,\bigr)
= \dfrac{d}{dx}\,\bigl(\,\varrho\,x - x^2 + a\,\varrho - a\,x\,\bigr) = - 2\,x + \varrho - a$,
\end{tabular}\\[2ex]
\hspace*{3.3ex}\begin{tabular}[t]{@{}l@{}}
$\dfrac{d}{dx}\,\bigl(\,x + a\,\bigr)\,\bigl(\,x - \varrho - \dfrac{3}{2}\,\lambda\,\bigr)
= \dfrac{d}{dx}\,\bigl(\,x^2 - \varrho\,x - \dfrac{3}{2}\,\lambda\,x + a\,x - a\,\varrho - \dfrac{3}{2}\,\lambda\,a\,\bigr)
= 2\,x - \varrho - \dfrac{3}{2}\,\lambda + a$,\\[2ex]
$\dfrac{d}{dx}\,\bigl(\,x + a\,\bigr)\,\bigl(\,\varrho + 3\,\lambda - x\,)
= \dfrac{d}{dx}\,\bigl(\,\varrho\,x + 3\,\lambda\,x - x^2 + a\,\varrho + 3\,\lambda\,a - a\,x\,\bigr)
= - 2\,x + \varrho + 3\,\lambda - a$
\end{tabular}\\[3ex]
and since each of the terms $\pm (\,\varrho + \lambda + a\,)\,\lambda\,F_{\SBTB}(\,\varrho + \lambda\,)$ 
and $\pm (\,\varrho + 2\,\lambda + a\,)\,\lambda\,F_{\SBTB}(\,\varrho + 2\,\lambda\,)$ cancels each other out,
we further obtain\\[2.5ex]
\hspace*{9.2ex}\begin{tabular}[t]{@{}c@{\hspace*{0.8ex}}l@{\hspace*{-32.8ex}}r@{}}
$=$&$\displaystyle{\bigg|\,\dfrac{1}{\lambda^3}\,\biggl\{\,\int\limits_{\varrho}^{\nu + \lambda}
\bigl(\,2\,x - \varrho + a\,\bigr)\,F_{\SBTB}(x)\,dx - 2\, \int\limits_{\varrho + \lambda}^{\varrho + 2\,\lambda}
\bigl(\,2\,x - \varrho - \dfrac{3}{2}\,\lambda + a\,\bigr)\,F_{\SBTB}(x)\,dx}$\\[4.5ex]
&$\displaystyle{+\ \int\limits_{\varrho + 2\,\lambda}^{\varrho + 3\,\lambda}
\bigl(\,2\,x - \varrho - 3\,\lambda + a\,\bigr)\,F_{\SBTB}(x)\,dx\,\biggr\}\,\bigg|}$\\[4.5ex]
$=$&$\displaystyle{\bigg|\,\dfrac{1}{\lambda^3}\,\biggl\{\,\int\limits_{0}^{\lambda}
\bigl(\,2\,v + \varrho + a\,\bigr)\,F_{\SBTB}(\,v + \varrho\,)\,dv - 2\, \int\limits_{0}^{\lambda}
\bigl(\,2\,v + \varrho + \dfrac{1}{2}\,\lambda + a\,\bigr)\,F_{\SBTB}(\,v + \varrho + \lambda\,)\,dv}$\\[4.5ex]
&$\displaystyle{+\ \int\limits_{0}^{\lambda}
\bigl(\,2\,v + \varrho + \lambda + a\,\bigr)\,F_{\SBTB}(\,v + \varrho + 2\,\lambda\,)\,dv\,\biggr\}\,\bigg|}$
&(subst.\index{integration!by substitution} $x = v + \varrho + t\,\lambda$, $t = 0\,; 1\,; 2$)\\[4.5ex]
$=$&$\displaystyle{\bigg|\,\dfrac{1}{\lambda^3}\,\biggl\{\,\int\limits_{0}^{\lambda}
(\,v + \varrho\,)\,\Delta_{\lambda}^{2}F_{\SBTB}(\,v + \varrho\,)\,dv + a\,\int\limits_{0}^{\lambda}
\Delta_{\lambda}^{2}F_{\SBTB}(\,v + \varrho\,)\,dv}$\\[4.5ex]
&$\displaystyle{+\ \int\limits_{0}^{\lambda}
v\,\Bigl(\,F_{\SBTB}(\,v + \varrho + 2\,\lambda\,) - 
F_{\SBTB}(\,v + \varrho + \lambda\,)\,\Bigr)\,dv}$\\[4.5ex]
&$\displaystyle{-\ \int\limits_{0}^{\lambda}
v\,\Bigl(\,F_{\SBTB}(\,v + \varrho + \lambda\,) - 
F_{\SBTB}(\,v + \varrho\,)\,\Bigr)\,dv}$\\[4.5ex]
&$\displaystyle{+\ \lambda\,\int\limits_{0}^{\lambda}
\Bigl(\,F_{\SBTB}(\,v + \varrho + 2\,\lambda\,) - 
F_{\SBTB}(\,v + \varrho + \lambda\,)\,\Bigr)\,dv\,\biggr\}\,\bigg|}$\\[4.5ex]
$\leq$&$c_{21}\,\mathcal{C}_{3} + c_{22}$.
\end{tabular}\\[0.8ex]
For the inequality, (\ref{EWD_3_8_14}) and $|\,a\,| \leq 16$ and 
\mbox{\rule[0ex]{0ex}{4.2ex}$E_{\!A} \leq \sqrt{\dfrac{4}{3}}\,\lambda$} 
(cf. (\ref{EWD_3_8_09})) were used to estimate the first and second integral.\\[0.8ex] 
To estimate the other three integrals, \mbox{\rule[0ex]{0ex}{2.4ex}Theorem} \ref{EWD_3_1_03}
was applied as in (\ref{EWD_3_7_31}) (or (\ref{EWD_3_7_33})).\\[0.8ex]
All in all, we obtain the assertion of Proposition \ref{EWD_3_8_33} and thus the Theorem \ref{EWD_3_1_13}.
\hspace*{1ex}\hfill$\Box$
                           
\rehead{Applications for simple linear rank statistics}    
\chapter[Applications of the two theorems from the chapter \ref{EWD_Kap3}]
{Applications of the two theorems from the chapter \ref{EWD_Kap3}}\label{EWD_Kap4}

\section{Preparations}\label{EWD_Kap4_Sec1}

In this chapter, we give easily provable conditions for simple linear rank statistics
\index{rank statistic!linear!simple} under which the requirements (\ref{EWD_3_1_11}), (\ref{EWD_3_1_14}) 
and (\ref{EWD_3_1_15}) are fulfilled (but with {''constants''} $\mathcal{C}_{1}$, $\mathcal{C}_{2}$ 
and $\mathcal{C}_{3}$, which grow slightly with $n$). We then obtain the estimates 
of the Theorems \ref{EWD_3_1_10} and \ref{EWD_3_1_13} (cf. Theorem \ref{EWD_4_3_13}).\\[2.8ex]
The conditions originate from van Zwet \cite{vanZwet1982},\index{conditions of van Zwet} who used them to derive an estimate of the characteristic function of $\mathscr{F}_{\!A}$ (cf. Theorem \ref{EWD_4_3_01}).
In order to be able to deduce our requirements (\ref{EWD_3_1_11}), (\ref{EWD_3_1_14}) and (\ref{EWD_3_1_15}) 
from this result, we must derive relations between the expressions on the left-hand side of these requirements and the characteristic functions of the distribution functions used there. This is done in the following lemma.\\[2.8ex]
To formulate this lemma, we first choose a distribution function $U$ on $\mathbb{R}$ 
that has a density $u$ with the following properties:\\[2.5ex]
\refstepcounter{DSBcount}
\label{EWD_4_1_01}
\text{\hspace*{-0.8ex}(\theDSBcount)}
\hspace*{4ex}
$u$ is infinitely differentiable,\\[2.5ex]
\refstepcounter{DSBcount}
\label{EWD_4_1_02}
\text{\hspace*{-0.8ex}(\theDSBcount)}
\hspace*{4ex}
the support of $u$ is contained in the intervall $[\,-1,\,1\,]$.
\index{function!with compact support}\index{function!with compact support!$u$, $U$}\\[2.5ex]
The functions $U$ and $u$ are fixed for the entire chapter.
From the two properties of $u$ we can conclude that all derivatives of
$u(x)$ and $x\,u(x)$ are integrable and therefore holds
(cf. Feller \cite{feller1971introduction}, Chapter XV, Section 4, Lemmas 2 and 4):\\[2.5ex]
\refstepcounter{DSBcount}
\label{EWD_4_1_03}
\text{\hspace*{-0.8ex}(\theDSBcount)}
\hspace*{4ex}
\begin{tabular}[t]{@{}l@{}}
$\displaystyle{
\left\{
\begin{array}{l@{\hspace*{2.6ex}}l@{\hspace*{2.6ex}}l@{}}
\bigl|\,\hat{U}(t)\,\bigr| = \smallO(\,|\,t\,|^{-m})&\text{for }|\,t\,| \rightarrow \infty
&\text{and for all }m \in \mathbb{N},\\[2ex]
\bigl|\,\hat{U}'(t)\,\bigr| = \smallO(\,|\,t\,|^{-m})&\text{for }|\,t\,| \rightarrow \infty
&\text{and for all }m \in \mathbb{N}.
\end{array}  \right.}$
\end{tabular}\\[2.5ex]
In the following, $\hat{G}$\index{function!characteristic}\index{function!characteristic!$\hat{G}$} denotes
the characteristic function of the (distribution) function $G$.\\[2.8ex]
From (\ref{EWD_4_1_03}) we conclude\\[2.5ex]
\refstepcounter{DSBcount}
\label{EWD_4_1_04}
\text{\hspace*{-0.8ex}(\theDSBcount)}
\hspace*{4ex}
$\displaystyle{\int\limits_{\mathbb{R}} \bigl|\,t\,\bigr|^{m}\,\bigl|\,\hat{U}(t)\,\bigr|\,dt < \infty}$
\hspace*{2ex}and\hspace*{2ex}
$\displaystyle{\int\limits_{\mathbb{R}} \bigl|\,t\,\bigr|^{m}\,\bigl|\,\hat{U}'(t)\,\bigr|\,dt < \infty}$
\hspace*{3ex}for all $m \in \mathbb{N}_{0}$,\\[2.5ex]
\refstepcounter{DSBcount}
\label{EWD_4_1_05}
\text{\hspace*{-0.8ex}(\theDSBcount)}
\hspace*{4ex}
$\displaystyle{\int\limits_{\mathbb{R}} \bigl|\,\hat{U}(t)\,\bigr|^{m}\,dt < \infty}$
\hspace*{2ex}and\hspace*{2ex}
$\displaystyle{\int\limits_{\mathbb{R}} \bigl|\,\hat{U}'(t)\,\bigr|^{m}\,dt < \infty}$
\hspace*{3ex}for all $m \in \mathbb{N}$.\\[2.5ex]
The following also follows from (\ref{EWD_4_1_02})\\[2ex]
\refstepcounter{DSBcount}
\label{EWD_4_1_06}
\text{\hspace*{-0.8ex}(\theDSBcount)}
\hspace*{4ex}
$\bigl|\,\hat{U}'(t)\,\bigr| \leq 1$
\hspace*{3ex}for all $t \in \mathbb{R}$.\\[2.8ex] 
After these preparations, we now prove the lemma announced above. The proof used here is
based on a technique introduced by von Bahr \cite{vonBahr1967} (cf. in particular Section 3).\\[4ex]
\refstepcounter{DSBcount}
\label{EWD_4_1_07}
\textbf{\hspace*{-0.8ex}\theDSBcount\ Lemma}\\[0.8ex]
Let $F$ be a distribution function on $\mathbb{R}$ and
\index{function!with compact support!$U_{\theta}$}\\[2ex]
\hspace*{12.1ex}$U_{\theta}(x) = U\Bigl(\,\dfrac{x}{\theta}\,\Bigr)$
\hspace*{3ex}for $x \in \mathbb{R}$ and $\theta > 0$.\\[2.5ex]
Then
\begin{enumerate}
\item\label{EWD_4_1_07_BWa}
$\displaystyle{\big|\,\Delta_{y}^{2}F(z)\,\big|\, \leq\, \bigl(\,y^2 + \theta\,\bigr)\,\int\limits_{\mathbb{R}}
\bigl(\,1 + |\,t\,|\,\bigr)\,\big|\,\hat{F}(t)\,\hat{U}_{\theta}(t)\,\big|\,dt}$\\[2.2ex]
for all $z \in \mathbb{R}$, $y > 0$ and $\theta > 0$.
\item\label{EWD_4_1_07_BWb}
$\displaystyle{\Big|\,F(x) - P_{\lambda}^{2}(x;z,F)\,\Big|\, \leq\, \theta\,\Bigl(\,1 + 
\dfrac{|\,x - z\,|^3}{\lambda^3}\,\Bigr)\,\int\limits_{\mathbb{R}}
\bigl(\,2 + t^2\,\bigr)\,\big|\,\hat{F}(t)\,\hat{U}_{\theta}(t)\,\big|\,dt}$\\[2.2ex]
for all $x, z \in \mathbb{R}$ and $\theta \geq \lambda^3 > 0$.
\item\label{EWD_4_1_07_BWc}
If $F$ has an existing first moment, then\\[2ex]
\begin{tabular}[t]{@{}l@{\hspace*{1ex}}c@{\hspace*{1ex}}l@{}}
$\displaystyle{\big|\,z\,\Delta_{y}^{2}F(z)\,\big|}$&$\leq$&
$\displaystyle{\bigl(\,y^2 + \theta\,\bigr)\,\int\limits_{\mathbb{R}}
\bigl(\,1 + |\,t\,|\,\bigr)\,\Bigl\{\,\big|\,\hat{F}'(t)\,\hat{U}_{\theta}(t)\,\big|\,+\,
\big|\,\hat{F}(t)\,\hat{U}'_{\theta}(t)\,\big|\,\Bigr\}\,dt}$\\[3.5ex]
&&$\displaystyle{+\ 2\,\bigl(\,y^2 + \theta^2\,\bigr)\,\int\limits_{\mathbb{R}}
\big|\,\hat{F}(t)\,\hat{U}_{\theta}(t)\,\big|\,dt}$
\end{tabular}\\[2.2ex]
for all $z \in \mathbb{R}$, $y > 0$ and $\theta > 0$.
\end{enumerate}
\vspace*{3.5ex}

\pagebreak

\textbf{Proof:}\\[0.8ex]
We note in advance that\\[2ex]
\refstepcounter{DSBcount}
\label{EWD_4_1_08}
\text{\hspace*{-0.8ex}(\theDSBcount)}
\hspace*{4ex}
$\displaystyle{\bigg|\,\int\limits_{u + \alpha}^{v - \alpha} e^{-its}\,ds - 
\int\limits_{u}^{v} e^{-its}\,ds\,\bigg|
= \bigg|\,\int\limits_{v - \alpha}^{v} e^{-its}\,ds + 
\int\limits_{u}^{u + \alpha} e^{-its}\,ds\,\bigg|\, \leq\, 2\,| \alpha |}$\\[2.5ex]
is valid for all $u, v, \alpha, t \in \mathbb{R}$.\\[2.8ex]
In the following, let $F_{\theta}$ be the convolution of $F$ and $U_{\theta}$. Because of (\ref{EWD_4_1_02})
then holds\\[2ex]
\refstepcounter{DSBcount}
\label{EWD_4_1_09}
\text{\hspace*{-0.8ex}(\theDSBcount)}
\hspace*{4ex}
$F_{\theta}(x - \theta)\, \leq\, F(x)\, \leq\, F_{\theta}(x + \theta)$
\hspace*{3ex}for all $x \in \mathbb{R}$ and $\theta > 0$.\\[2.5ex]
Furthermore, $F_{\theta}$ has a density $f_{\theta}$ that is square-integrable. This follows from
the square-integrability of $\hat{f}_{\theta} = \hat{F}\,\hat{U}_{\theta}$ 
(cf. (\ref{EWD_4_1_05})) and from the Plancherel identity\index{Plancherel!Plancherel identity}  
(cf. Feller \cite{feller1971introduction}, Chapter XV, Section 3, Example c)).\\[2.8ex]
We now obtain the assertions of the lemma as follows:
\begin{enumerate}
\item
If $y \geq 2\,\theta$, then $z + y - \theta \geq z + \theta$ and thus\\[2ex]
\begin{tabular}[t]{@{}c@{\hspace*{0.8ex}}l@{\hspace*{6.4ex}}l@{}}
&$F(z + 2\,y) - 2\,F(z + y) + F(z)$\\[2.5ex]
$\leq$&$F_{\theta}(z + 2\, y + \theta) - 2\,F_{\theta}(z + y - \theta) + F_{\theta}(z + \theta)$
&(cf. (\ref{EWD_4_1_09}))\\[2.5ex]
$=$&$\Bigl(\,F_{\theta}(z + 2\, y + \theta) - F_{\theta}(z + y - \theta)\,\Bigr) - 
\Bigl(\,F_{\theta}(z + y - \theta) - F_{\theta}(z + \theta)\,\Bigr)$\\[2.5ex]
$\leq$&$\displaystyle{\bigg|\,\int\limits_{\mathbb{R}}
\Bigl(\,1_{\textstyle{(\,z + y - \theta,\,z + 2\,y + \theta\,]}}(t) -
1_{\textstyle{(\,z + \theta,\,z + y - \theta\,]}}(t)\,\Bigr)\,f_{\theta}(t)\,dt\,\bigg|}$.
\end{tabular}\\[3ex]
An application of Plancherel's theorem\index{Plancherel!Plancherel's theorem} for scalar products 
(see Bhattacharya, Rao \cite{doi:10.1137/1.9780898719895}, Theorem 4.2) further leads to\\[2.5ex] 
\begin{tabular}[t]{@{}c@{\hspace*{0.8ex}}l@{}}
$=$&$\displaystyle{\dfrac{1}{2\,\pi}\,\bigg|\,\int\limits_{\mathbb{R}}
\Bigl(\,\int\limits_{z + y - \theta}^{z + 2y + \theta} e^{-its}\,ds\, -
\int\limits_{z + \theta}^{z + y - \theta} e^{-its}\,ds\,\Bigr)\,\hat{f}_{\theta}(t)\,dt\,\bigg|}$\\[5.5ex]
$\leq$&$\displaystyle{\dfrac{1}{2\,\pi}\,\int\limits_{\mathbb{R}} \dfrac{1}{|\,t\,|}\,
\Big|\,e^{-it(z + 2y)} - 2\,e^{-it(z + y)} + e^{-itz}\,\Big|\,\big|\,\hat{f}_{\theta}(t)\,\big|\,dt\,\, +}$\\[4.5ex]
&$\displaystyle{\dfrac{1}{2\,\pi}\,\int\limits_{\mathbb{R}}\,\biggl\{\,\bigg|\,
\int\limits_{z + y - \theta}^{z + 2y + \theta} e^{-its}\,ds -
\int\limits_{z + y}^{z + 2y} e^{-its}\,ds\,\bigg|\, +\,
\bigg|\,
\int\limits_{z + \theta}^{z + y - \theta} e^{-its}\,ds -
\int\limits_{z}^{z + y} e^{-its}\,ds\,\bigg|\,\biggr\}\,\big|\,\hat{f}_{\theta}(t)\,\big|\,dt}$.
\end{tabular}\\[3ex]
For the first integral, we now apply Lemma \ref{EWD_2_4_02}, \ref{EWD_2_4_02_BWb}) for each $t \in \mathbb{R}$,
once for the function $g_{re}(w) = \cos(tw)$ and once for $g_{im}(w) = \sin(tw)$.
For the second integral, we use the inequality (\ref{EWD_4_1_08}).
Thus, together with $\hat{f}_{\theta} = \hat{F}\,\hat{U}_{\theta}$, we further obtain\\[2.5ex]
\begin{tabular}[t]{@{}c@{\hspace*{0.8ex}}l@{\hspace*{6ex}}l@{}}
$\leq$&$\displaystyle{y^2\,\int\limits_{\mathbb{R}} |\,t\,|\,\big|\,\hat{F}(t)\,\hat{U}_{\theta}(t)\,\big|\,dt
+ \theta\,\int\limits_{\mathbb{R}} \big|\,\hat{F}(t)\,\hat{U}_{\theta}(t)\,\big|\,dt}$\\[4.5ex]
$\leq$&$\displaystyle{\bigl(\,y^2 + \theta\,\bigr)\,\int\limits_{\mathbb{R}}
\bigl(\,1 + |\,t\,|\,\bigr)\,\big|\,\hat{F}(t)\,\hat{U}_{\theta}(t)\,\big|\,dt := R}$.
\end{tabular}\\[2.5ex]
Furthermore, from $y \geq 2\,\theta$ also follows $z + 2y - \theta \geq z + y + \theta$, 
so that we get by com\-plete\-ly analogous considerations\\[2ex]
\begin{tabular}[t]{@{}c@{\hspace*{0.8ex}}l@{\hspace*{-3.2ex}}r@{}}
&$-\ \Bigl(\,F(z + 2\,y) - 2\,F(z + y) + F(z)\,\Bigr)$\\[2.5ex]
$\leq$&$\Bigl(\,F_{\theta}(z + y + \theta) - F_{\theta}(z - \theta)\,\Bigr) -
\Bigl(\,F_{\theta}(z + 2\, y - \theta) - F_{\theta}(z + y + \theta)\,\Bigr)$
&(cf. (\ref{EWD_4_1_09}))\\[2.5ex]
$\leq$&$\displaystyle{\dfrac{1}{2\,\pi}\,\bigg|\,\int\limits_{\mathbb{R}}
\Bigl(\,\int\limits_{z + y + \theta}^{z + 2y - \theta} e^{-its}\,ds\, -
\int\limits_{z - \theta}^{z + y + \theta} e^{-its}\,ds\,\Bigr)\,\hat{f}_{\theta}(t)\,dt\,\bigg|}$
&(Plancherel's theorem\index{Plancherel!Plancherel's theorem}).
\end{tabular}\\[2.5ex]
The last expression can be estimated in the same way as above by $R$.\\[2.8ex]
If, on the other hand, $y < 2\,\theta$, then $z + y - \theta < z + \theta$ and thus\\[2ex]
\begin{tabular}[t]{@{}c@{\hspace*{0.8ex}}l@{\hspace*{-13.4ex}}r@{}}
&$F(z + 2\,y) - 2\,F(z + y) + F(z)$\\[2.5ex]
$\leq$&$\Bigl(\,F_{\theta}(z + 2\, y + \theta) - F_{\theta}(z + y - \theta)\,\Bigr) + 
\Bigl(\,F_{\theta}(z + \theta) - F_{\theta}(z + y - \theta)\,\Bigr)$
&(cf. (\ref{EWD_4_1_09}))\\[2.5ex]
$\leq$&$\displaystyle{\dfrac{1}{2\,\pi}\,\bigg|\,\int\limits_{\mathbb{R}}
\Bigl(\,\int\limits_{z + y - \theta}^{z + 2y + \theta} e^{-its}\,ds\, +
\int\limits_{z + y - \theta}^{z + \theta} e^{-its}\,ds\,\Bigr)\,\hat{f}_{\theta}(t)\,dt\,\bigg|}$
&(Plancherel's theorem\index{Plancherel!Plancherel's theorem})\\[4.5ex]
$=$&$\displaystyle{\dfrac{1}{2\,\pi}\,\bigg|\,\int\limits_{\mathbb{R}}
\Bigl(\,\int\limits_{z + y - \theta}^{z + 2y + \theta} e^{-its}\,ds\, -
\int\limits_{z + \theta}^{z + y - \theta} e^{-its}\,ds\,\Bigr)\,\hat{f}_{\theta}(t)\,dt\,\bigg|}$
&(reversing the integration limits).
\end{tabular}\\[2.5ex]
This expression has already been estimated in detail by $R$ above.\\[1.5ex]
Finally, from $y < 2\,\theta$ also follows $z + 2y - \theta < z + y + \theta$, 
so we get analogously to above\\[2ex]
\begin{tabular}[t]{@{}c@{\hspace*{0.8ex}}l@{\hspace*{3.1ex}}r@{}}
&$-\ \Bigl(\,F(z + 2\,y) - 2\,F(z + y) + F(z)\,\Bigr)$\\[2.5ex]
$\leq$&$\Bigl(\,F_{\theta}(z + y + \theta) - F_{\theta}(z - \theta)\,\Bigr) +
\Bigl(\,F_{\theta}(z + y + \theta) - F_{\theta}(z + 2\, y - \theta)\,\Bigr) \leq R$
&(cf. (\ref{EWD_4_1_09})).
\end{tabular}\\[2.5ex]
Overall, the validity of part \ref{EWD_4_1_07_BWa}) follows from the 4 cases considered.
\item
Let $x \geq z + \lambda$. Then\\[2ex]
\begin{tabular}[t]{@{}c@{\hspace*{0.8ex}}l@{\hspace*{-1.8ex}}r@{}}
&$F(x) - P_{\lambda}^{2}(x;z,F)$\\[2.5ex]
$\leq$&$F_{\theta}(x + \theta) - F_{\theta}(z - \theta) - 
\dfrac{x - z}{\lambda}\,\Bigl(\,F_{\theta}(z + \lambda - \theta) - F_{\theta}(z + \theta)\,\Bigr)$\\[2.5ex]
&$-\ \dfrac{1}{2}\,\dfrac{x - z}{\lambda}\,\dfrac{x - z - \lambda}{\lambda}\,
\Bigl(\,F_{\theta}(z + 2\,\lambda - \theta) - 2\,F_{\theta}(z + \lambda + \theta)
+ F_{\theta}(z - \theta)\,\Bigr)$\\[3ex]
$\leq$&$\displaystyle{\dfrac{1}{2\,\pi}\,\Bigg|\,\int\limits_{\mathbb{R}}\, \biggl\{
\,\int\limits_{z - \theta}^{x + \theta} e^{-its}\,ds\, - \, \dfrac{x - z}{\lambda}\,
\int\limits_{z + \theta}^{z + \lambda - \theta} e^{-its}\,ds\,}$\\[4.5ex]
&\hspace*{9ex}$\displaystyle{-\ \dfrac{1}{2}\,\dfrac{x - z}{\lambda}\,\dfrac{x - z - \lambda}{\lambda}\,
\Bigl(\,\int\limits_{z + \lambda + \theta}^{x + 2\,\lambda - \theta} e^{-its}\,ds\, - \,
\int\limits_{z - \theta}^{x + \lambda + \theta} e^{-its}\,ds\,\Bigr)\,\biggr\}\,\,\hat{f}_{\theta}(t)\,dt\,\Bigg|}$\\[4.5ex]
$\leq$&$\displaystyle{\dfrac{1}{2\,\pi}\,\int\limits_{\mathbb{R}} \dfrac{1}{|\,t\,|}\,
\Big|\,e^{-itx} - P_{\lambda}^{2}(x;z,e^{-it\boldsymbol{.}})\,\Big|\,\big|\,\hat{f}_{\theta}(t)\,\big|\,dt\,\, +}$\\[4.5ex]
&\hspace*{9ex}$\displaystyle{+\ \dfrac{1}{2\,\pi}\,\int\limits_{\mathbb{R}} 2\,\theta\,
\Bigl(\,1\, +\, \dfrac{|x - z|}{\lambda}\, + \,\dfrac{|x - z|}{\lambda}\,\dfrac{|x - z - \lambda|}{\lambda}\,\Bigr)\,
\big|\,\hat{f}_{\theta}(t)\,\big|\,dt}$&(cf. (\ref{EWD_4_1_08})).
\end{tabular}\\[3ex]
We now apply Lemma \ref{EWD_2_4_02}, \ref{EWD_2_4_02_BWa}) to the first integral and the inequality\\[1.5ex]
\hspace*{2.7ex}$1 + \eta + \eta\,(1 + \eta) = (1 + \eta)^2 \leq (1 + \eta)^3 \leq
4\,(1 + \eta^3)$\hspace*{2ex}for $\eta > 0$
\hfill(\index{H{\"o}lder's inequality!for finite sequences using length $\nu$}(\ref{EWD_0_1_05}), 
$\nu = 2$, $p = 3$)\\[1.5ex] 
to the second integral. This gives\\[2.5ex]
\begin{tabular}[t]{@{}c@{\hspace*{0.8ex}}l@{\hspace*{-6.4ex}}r@{}}
$\leq$&$\displaystyle{\Bigl(\,\lambda^3 + |x-z|^3\,\Bigr)\,\int\limits_{\mathbb{R}} 
t^2\,\big|\,\hat{f}_{\theta}(t)\,\big|\,dt\, +\,
2\,\theta\,\Bigl(\,1 + \dfrac{|x - z|^3}{\lambda^3}\,\Bigr)\,
\int\limits_{\mathbb{R}} \big|\,\hat{f}_{\theta}(t)\,\big|\,dt}$\\[4.5ex]
$\leq$&$\displaystyle{\theta\,\Bigl(\,1 + \dfrac{|x - z|^3}{\lambda^3}\,\Bigr)\,
\int\limits_{\mathbb{R}}
\bigl(\,2 + t^2\,\bigr)\,\big|\,\hat{F}(t)\,\hat{U}_{\theta}(t)\,\big|\,dt}$
&(cf. condition $\theta \geq \lambda^3 > 0$).
\end{tabular}\\[2.5ex]
Since we obtain these estimates completely analogously (as in part \ref{EWD_4_1_07_BWa})) 
for $x \leq z$ and
\linebreak 
$z \leq x \leq z + \lambda$ as well as
$-\ \Bigl(\,F(x) - P_{\lambda}^{2}(x;z,F)\,\Bigr)$, part \ref{EWD_4_1_07_BWb}) follows.
\item
Let $y \geq 2\,\theta$ and $z \geq 0$. Then\\[2ex]
\begin{tabular}[t]{@{}c@{\hspace*{0.8ex}}l@{\hspace*{25.4ex}}l@{}}
&$z\,\Bigl(\,F(z + 2\,y) - 2\,F(z + y) + F(z)\,\Bigr)$\\[2.5ex]
$\leq$&$z\,\Bigl(\,F_{\theta}(z + 2\, y + \theta) - 2\,F_{\theta}(z + y - \theta) + F_{\theta}(z + \theta)\,\Bigr)$
&(cf. (\ref{EWD_4_1_09}))
\end{tabular}\\[2.5ex]
\begin{tabular}[t]{@{}c@{\hspace*{0.8ex}}l@{\hspace*{-2.2ex}}r@{}}
$\leq$&$\displaystyle{\bigg|\,\int\limits_{\mathbb{R}}
z\,\Bigl(\,1_{\textstyle{(\,z + y - \theta,\,z + 2\,y + \theta\,]}}(t) -
1_{\textstyle{(\,z + \theta,\,z + y - \theta\,]}}(t)\,\Bigr)\,f_{\theta}(t)\,dt\,\bigg|}$\\[4.5ex]
$\leq$&$\displaystyle{\bigg|\,\int\limits_{\mathbb{R}}
t\,\Bigl(\,1_{\textstyle{(\,z + y - \theta,\,z + 2\,y + \theta\,]}}(t) -
1_{\textstyle{(\,z + \theta,\,z + y - \theta\,]}}(t)\,\Bigr)\,f_{\theta}(t)\,dt\,\bigg|}$\\[4.5ex]
&$\displaystyle{+\ \bigg|\,\int\limits_{z + y - \theta}^{z + 2\,y + \theta} (t - z)\,f_{\theta}(t)\,dt\,\bigg|\,
+\, \bigg|\,\int\limits_{z + \theta}^{z + y - \theta} (t - z)\,f_{\theta}(t)\,dt\,\bigg|}$\\[5ex]
$\leq$&$\displaystyle{\bigg|\,\int\limits_{\mathbb{R}}
\Bigl(\,1_{\textstyle{(\,z + y - \theta,\,z + 2\,y + \theta\,]}}(t) -
1_{\textstyle{(\,z + \theta,\,z + y - \theta\,]}}(t)\,\Bigr)\,\bigl(\,t\,f_{\theta}(t)\,\bigr)\,dt\,\bigg|}$\\[4.5ex]
&$\displaystyle{+\ (\,2\,y + \theta\,)\,
\int\limits_{z + y - \theta}^{z + 2\,y + \theta} f_{\theta}(t)\,dt\,
+\, (\,y + \theta\,)\,\int\limits_{z + \theta}^{z + y - \theta} f_{\theta}(t)\,dt}$
&($f_{\theta}$ is a density!)\\[4.5ex]
$\leq$&$\displaystyle{\dfrac{1}{2\,\pi}\,\,\bigg|\,\int\limits_{\mathbb{R}}
\Bigl(\,\int\limits_{z + y - \theta}^{z + 2y + \theta} e^{-its}\,ds\, -
\int\limits_{z + \theta}^{z + y - \theta} e^{-its}\,ds\,\Bigr)\,\widehat{\bigl(\,id_{\mathbb{R}}\,f_{\theta}\,\bigr)}(t)\,dt\,\bigg|}$\\[4.5ex]
&$\displaystyle{+\ (\,2\,y + \theta\,)\,\,\dfrac{1}{2\,\pi}\,\,\bigg|\,\int\limits_{\mathbb{R}}
\Bigl(\,\int\limits_{z + y - \theta}^{z + 2y + \theta} e^{-its}\,ds\,\,\Bigr)\,\hat{f}_{\theta}(t)\,dt\,\bigg|}$\\[4.5ex]
&$\displaystyle{+\ (\,y + \theta\,)\,\,\dfrac{1}{2\,\pi}\,\,\bigg|\,\int\limits_{\mathbb{R}}
\Bigl(\,\int\limits_{z + \theta}^{z + y - \theta} e^{-its}\,ds\,\,\Bigr)\,\hat{f}_{\theta}(t)\,dt\,\bigg|}$\\[5ex]
$\leq$&$\displaystyle{\bigl(\,y^2 + \theta\,\bigr)\,\int\limits_{\mathbb{R}}
\bigl(\,1 + |\,t\,|\,\bigr)\,\Big|\,\widehat{\bigl(\,id_{\mathbb{R}}\,f_{\theta}\,\bigr)}(t)\,\Big|\,dt}$\\[4.5ex]
&$\displaystyle{+\ \bigl\{\,(\,2\,y + \theta\,)\,(\,y + 2\,\theta\,) + (\,y + \theta\,)\,(\,y + 2\,\theta\,)\,\bigr\}\,
\dfrac{1}{2\,\pi}\,\int\limits_{\mathbb{R}} \big|\,\hat{f}_{\theta}(t)\,\big|\,dt}$.
\end{tabular}\\[2.5ex]
For the last inequality, the first summand was estimated as in the proof of part \ref{EWD_4_1_07_BWa}).
If we now use\\[2ex] 
\hspace*{7ex}$\hat{f}_{\theta} = \hat{F}\,\hat{U}_{\theta}$,\hspace*{2ex} 
$\widehat{\bigl(\,id_{\mathbb{R}}\,f_{\theta}\,\bigr)} = \dfrac{1}{i}\,(\hat{f}_{\theta})'$\hspace*{2ex}and\\[2.2ex]
\hspace*{7ex}$(\,2\,y + \theta\,)\,(\,y + 2\,\theta\,) + (\,y + \theta\,)\,(\,y + 2\,\theta\,)
= 3\,y^2 + 8\,\theta\,y + 4\,\theta^2 \leq 7\,y^2 + 8\,\theta^2$,\\[2.5ex]
the asserted estimate in part \ref{EWD_4_1_07_BWc}) follows for the considered case.\\[2.8ex] 
All other cases are derived using completely analogous considerations.\hspace*{1ex}\hfill$\Box$
\end{enumerate}

\section{A conclusion from Cram\'er's condition}\label{EWD_Kap4_Sec2}

This section is to be understood as an insertion and serves to illustrate some basic ideas of the 
further procedure in the simpler
case of independent and identically distributed (iid) random variables.
More precisely, we show how to obtain the assumption (\ref{EWD_1_1_04}) of 
Theorem \ref{EWD_1_1_03} from Cram\'er's condition. In the next section \ref{EWD_Kap4_Sec3}, 
an estimate of the characteristic function\index{conditions of van Zwet} of a simple linear rank statistic 
based on van Zwet will play the role of Cram\'er's condition.\\[2.8ex]
In this section, let $X_{i}$, $i \in \mathbb{N}$, be a sequence of independent and identically distributed
random variables with a common distribution function $F$ and\\[2ex]
\refstepcounter{DSBcount}
\label{EWD_4_2_01}
\text{\hspace*{-0.8ex}(\theDSBcount)}
\hspace*{4ex}
$E(X_{i}) = 0$
\hspace*{2ex}and\hspace*{2ex}
$E(X_{i}^2) = 1$.\\[2.5ex]
Furthermore, $F_{n}$ denotes the distribution function of 
$S_{n} = \displaystyle{\dfrac{1}{\sqrt{n}} \sum\limits_{i = 1}^n X_{i}}$ for $n \in \mathbb{N}$. Then
we get:\\[4ex]
\refstepcounter{DSBcount}
\label{EWD_4_2_02}
\textbf{\hspace*{-0.8ex}\theDSBcount\ Theorem}\index{Theorem!for the iid-case}\\[0.8ex]
Suppose that\\[2ex]
\refstepcounter{DSBcount}
\label{EWD_4_2_03}
\text{\hspace*{-0.8ex}(\theDSBcount)}
\hspace*{4ex}
$\displaystyle{\limsup\limits_{|\,t\,| \rightarrow \infty}\, \big|\,\hat{F}(t)\,\big| < 1}$
\hspace*{20ex}({''Cram\'er's condition\index{Cram\'er's condition}''}).\\[2.5ex]
Then there exists a constant $\mathcal{C}' > 0$ depending only on $F$ such that\\[2.3ex]
\refstepcounter{DSBcount}
\label{EWD_4_2_04}
\text{\hspace*{-0.8ex}(\theDSBcount)}
\hspace*{4ex}
$\displaystyle{||\Delta_{y}^{2}\, F_{n}|| \leq \mathcal{C}'\Bigl(\,\dfrac{1}{n} + y^2\, \Bigr)}$
\hspace*{2ex}for all $0 \leq y \leq \dfrac{1}{\sqrt{n}}$, $n \in \mathbb{N}$.\\[3.5ex]
\textbf{Proof:}\\[0.8ex]
At first, it should be noted that the constants $c_{1}$, $c_{2}$, $\delta$ and $\eta$ used in this proof 
depend on $F$ and nothing else.\\[2.8ex] 
According to Lemma \ref{EWD_4_1_07}, \ref{EWD_4_1_07_BWa}) with $\theta = \dfrac{1}{n}$ 
and $F = F_{n}$ remains to be shown:\\[2ex]
\hspace*{12.1ex}$\displaystyle{\int\limits_{\mathbb{R}}
\bigl(\,1 + |\,t\,|\,\bigr)\,\big|\,\hat{F}_{n}(t)\,\hat{U}_{\textstyle{\frac{1}{n}}}(t)\,\big|\,dt \leq c_{1}}$
\hspace*{2ex}for all $n \in \mathbb{N}$.\\[2.5ex]
Because of\\[2ex] 
\hspace*{12.1ex}$\big|\,\hat{F}_{n}(t)\,\hat{U}_{\textstyle{\frac{1}{n}}}(t)\,\big| \leq 1$ 
\hspace*{3ex}and\hspace*{3ex} 
$1 \leq |\,t\,|\ $ for $\ t \not\in\,\, ]-1,\,1\,[\ $,\\[2.5ex] 
and because of\\[2ex] 
\hspace*{12.1ex}$\hat{F}_{n}(t) = \Bigl(\,\hat{F}\bigl(\,\dfrac{t}{\sqrt{n}}\,\bigr)\,\Bigr)^n$ 
\hspace*{3ex}and\hspace*{3ex}
$\hat{U}_{\textstyle{\frac{1}{n}}}(t) = \hat{U}\bigl(\,\dfrac{t}{n}\,\bigr)$,\\[2.5ex] 
it is sufficient to prove the existence of a constant $c_{2}$ such that\\[3ex]
\refstepcounter{DSBcount}
\label{EWD_4_2_05}
\text{\hspace*{-0.8ex}(\theDSBcount)}
\hspace*{4ex}
$\displaystyle{\int\limits_{\mathbb{R}}
|\,t\,|\,\Big|\,\Bigl(\,\hat{F}\bigl(\,\dfrac{t}{\sqrt{n}}\,\bigr)\,\Bigr)^n\,\Big|\,
\Big|\,\hat{U}\bigl(\,\dfrac{t}{n}\,\bigr)\,\Big|\,dt \leq c_{2}}$
\hspace*{2ex}for all $n \in \mathbb{N}$.\\[2.8ex]
To do this, we first expand the characteristic function $\hat{F}$ about $0$ according to
Taylor's theorem\index{Taylor's theorem} taking into account (\ref{EWD_4_2_01}).
From this follows the existence of a $\delta > 0$, so that\\[2ex]
\hspace*{12.1ex}$\big|\,\hat{F}(t)\,\big| \leq 1 - \dfrac{t^2}{4}$
\hspace*{2ex}for $|\,t\,| \leq \delta$\\[2.5ex]
(cf. e.g. \cite{petrov1975sums}, Chapter I, Theorem 2 incl. proof, page 11).
We now split the above integral into the integration regions $|\,t\,| \leq \delta\,\sqrt{n}$ and
$|\,t\,| > \delta\,\sqrt{n}$. The estimate over the first region gives\\[2.5ex]
\hspace*{12.1ex}\begin{tabular}[t]{@{}l@{\hspace*{0.8ex}}c@{\hspace*{0.8ex}}l@{}}
$\displaystyle{\int\limits_{|\,t\,|\, \leq\, \delta\,\sqrt{n}}
|\,t\,|\,\Big|\,\hat{F}\bigl(\,\dfrac{t}{\sqrt{n}}\,\bigr)\,\Big|^n\,
\Big|\,\hat{U}\bigl(\,\dfrac{t}{n}\,\bigr)\,\Big|\,dt}$&
$\leq$&$\displaystyle{\int\limits_{|\,t\,|\, \leq\, \delta\,\sqrt{n}}
|\,t\,|\,\Bigl(\,1 - \dfrac{t^2}{4\,n}\,\Bigr)^n\,dt}$\\[4.5ex]
&$\leq$&$\displaystyle{\int\limits_{\mathbb{R}}
|\,t\,|\,\exp\Bigl(\, - \dfrac{t^2}{4}\,\Bigr)\,\,dt = 4}$.
\end{tabular}\\[3ex]
For the last inequality, $1 - x \leq e^{-x}$ for $x \in \mathbb{R}$ was used.\\[2.8ex]
For the estimate over the other region, we use that there is an $\epsilon > 0$ 
and a $\rho > 0$ due to Cram\'er's condition, so that\\[2ex]
\hspace*{12.1ex}$\big|\,\hat{F}(t)\,\big| \leq 1 - \epsilon$
\hspace*{2ex}for $|\,t\,| > \rho$.\\[2.5ex]
Since from this $\big|\,\hat{F}(t)\,\big| < 1$ follows for all $t \not= 0$ (cf. Feller \cite{feller1971introduction}, Chapter XV, Section 1, Lemma 4) and since $\big|\,\hat{F}\,\big|$ (in the case of $\delta < \rho$) is continuous on the compact interval $[\,\delta,\,\rho\,]$, there exists an $\eta > 0$ such that\\[2ex]
\hspace*{12.1ex}$\big|\,\hat{F}(t)\,\big| \leq 1 - \eta$
\hspace*{2ex}for $|\,t\,| > \delta$.\\[2.5ex]
It follows that\\[2ex]
\hspace*{12.1ex}\begin{tabular}[t]{@{}l@{\hspace*{0.8ex}}c@{\hspace*{0.8ex}}l@{}}
$\displaystyle{\int\limits_{|\,t\,|\, >\, \delta\,\sqrt{n}}
|\,t\,|\,\Big|\,\hat{F}\bigl(\,\dfrac{t}{\sqrt{n}}\,\bigr)\,\Big|^n\,
\Big|\,\hat{U}\bigl(\,\dfrac{t}{n}\,\bigr)\,\Big|\,dt}$&
$\leq$&$\displaystyle{\bigl(\,1 - \eta\,\bigr)^n\,\int\limits_{|\,t\,|\, >\, \delta\,\sqrt{n}}
|\,t\,|\,\Big|\,\hat{U}\bigl(\,\dfrac{t}{n}\,\bigr)\,\Big|\,dt}$.
\end{tabular}\\[3ex]
Now we make the substitution\index{integration!by substitution} 
$v = \dfrac{t}{n}$ and obtain furthermore\\[2ex]
\hspace*{22.1ex}\begin{tabular}[t]{@{}c@{\hspace*{0.8ex}}l@{\hspace*{10.5ex}}r@{}}
$=$&$\displaystyle{\bigl(\,1 - \eta\,\bigr)^n\,n^2\,\int\limits_{|\,v\,|\, >\, \textstyle{\frac{\delta}{\sqrt{n}}}}
|\,v\,|\,|\,\hat{U}(\,v\,)\,|\,dv}$\\[5ex]
$\leq$&$\displaystyle{\bigl(\,1 - \eta\,\bigr)^n\,n^2\,\int\limits_{\mathbb{R}}
|\,v\,|\,|\,\hat{U}(\,v\,)\,|\,dv} \rightarrow 0$
\hspace*{2ex}for $n \rightarrow \infty$&(cf. (\ref{EWD_4_1_04})).
\end{tabular}\\[2ex]
Therefore (\ref{EWD_4_2_05}) is valid.\hspace*{1ex}\hfill$\Box$\\[4ex]
If we summarize the above Theorem \ref{EWD_4_2_02} and the Theorem \ref{EWD_1_1_03}, we obtain the following
result, which goes back to Cram\'er \cite{cramer_1970}, Theorem 25 in Chapter VII
(cf. also the original articles \cite{doi:10.1080/03461238.1928.10416862} 
and \cite{doi:10.1080/03461238.1928.10416872}):\\[4ex]
\refstepcounter{DSBcount}
\label{EWD_4_2_06}
\textbf{\hspace*{-0.8ex}\theDSBcount\ Corollary}\\[0.8ex]
Suppose that\\[2ex]
\hspace*{12.1ex}$E(X_{i}^4) \leq \infty$
\hspace*{2ex}and\hspace*{2ex}
$\displaystyle{\limsup\limits_{|\,t\,| \rightarrow \infty}\, \big|\,\hat{F}(t)\,\big| < 1}$.\\[2.5ex]
Then there exists a constant $\mathcal{K}' > 0$ depending only on $F$ such that\\[2.5ex]
\hspace*{12.1ex}$\displaystyle{||F_{n} - e_{n}|| \leq \dfrac{\mathcal{K}'}{n}}$
\hspace*{2ex}for all $n \in \mathbb{N}$.
\vspace*{1ex}

\section{Use of a result from van Zwet}\label{EWD_Kap4_Sec3}

In this section we consider simple linear rank statistics.
We use a result of van Zwet\index{conditions of van Zwet} and Lemma \ref{EWD_4_1_07} to derive the
conditions (\ref{EWD_3_1_11}), (\ref{EWD_3_1_14}) and (\ref{EWD_3_1_15}) in a similar way as in the previous
section \ref{EWD_Kap4_Sec2}. However, the resulting {''constants''} $\mathcal{C}_{1}$, $\mathcal{C}_{2}$ 
and $\mathcal{C}_{3}$ will grow slightly with $n$.\\[2.8ex]
In the following $e_{1},\ldots,e_{n}$ (regression constants)\index{regression constants} 
and $d_{1},\ldots,d_{n}$ (scores)\index{scores} are sequences of real numbers and\\[2ex]
\hspace*{12.1ex}$\displaystyle{\bm\bar{e} = \dfrac{1}{n}\,\sum\limits_{i=1}^{n}\,e_{i}}$,\hspace*{4ex}
$\displaystyle{\bm\bar{d} = \dfrac{1}{n}\,\sum\limits_{j=1}^{n}\,d_{j}}$\\[2.5ex]
and $a_{ij} = e_{i}\,d_{j}$ for $i, j = 1,\ldots,n$. For the matrix $A = (a_{ij})$,
the same notations as in chapter \ref{EWD_Kap3} shall apply. 
Furthermore, let $\lambda$ be the Lebesgue measure on $\mathbb{R}$
\index{Lebesgue measure}\index{Lebesgue measure!$\lambda$}.\\[2.8ex]
The result used from van Zwet is now:\index{conditions of van Zwet}\\[4ex]
\refstepcounter{DSBcount}
\label{EWD_4_3_01}
\textbf{\hspace*{-0.8ex}\theDSBcount\ Theorem}\index{Theorem!for simple linear rank statistics!van Zwet}\\[0.8ex]
Suppose that there exist positive constants $e$, $E$, $d$, $D$ and $\delta$ such that\\[2ex]
\refstepcounter{DSBcount}
\label{EWD_4_3_02}
\text{\hspace*{-0.8ex}(\theDSBcount)}
\hspace*{4ex}
\begin{tabular}{@{}l@{}}
$\displaystyle{\sum\limits_{i=1}^{n}\,\big|\,e_{i} - \bm\bar{e}\,\big|^{r}\, \geq\, e n}$,\hspace*{4ex}
$\displaystyle{\sum\limits_{i=1}^{n}\,\big|\,e_{i} - \bm\bar{e}\,\big|^{k}\, \leq\, E n}$\\[2.5ex]
for some $k > 2$ and $0 < r < k$.
\end{tabular}\\[3ex]
\refstepcounter{DSBcount}
\label{EWD_4_3_03}
\text{\hspace*{-0.8ex}(\theDSBcount)}
\hspace*{4ex}
\begin{tabular}{@{}l@{}}
$\displaystyle{\sum\limits_{j=1}^{n}\,\big|\,d_{j} - \bm\bar{d}\,\big|^{m}\, \geq\, d n}$,\hspace*{4ex}
$\displaystyle{\sum\limits_{j=1}^{n}\,\big|\,d_{j} - \bm\bar{d}\,\big|^{s}\, \leq\, D n}$\\[2.5ex]
for some $s > 2$ and $0 < m < s$.
\end{tabular}\\[3ex]
\refstepcounter{DSBcount}
\label{EWD_4_3_04}
\text{\hspace*{-0.8ex}(\theDSBcount)}
\hspace*{4ex}
\begin{tabular}{@{}l@{}}
$\displaystyle{\lambda\Bigl(\,\bigcup_{j = 1}^{n}\,\Bigl\{\,x \in \mathbb{R}\,:\,\big|\,x - d_{j}\,\big| < \zeta\,
\Bigr\}\,\Bigr)\, \geq\, \delta n \zeta}$\\[3ex]
for some $\zeta \geq n^{-3/2}\,\log\ n$.
\end{tabular}\\[3.5ex]
Then there exist positive constants $a_{1}$, $a_{2}$, $a_{3}$ and $a_{4}$ depending only on
$e$, $E$, $d$, $D$, $\delta$ and $r$, $k$, $m$, $s$ such that\\[2ex]
\refstepcounter{DSBcount}
\label{EWD_4_3_05}
\text{\hspace*{-0.8ex}(\theDSBcount)}
\hspace*{4ex}
$\displaystyle{\big|\,\hat{F}_{\bm\hat{A}}(t)\,\big|\, \leq\,
a_{1}\,n^{- a_{2}\log n}}$
\hspace*{3ex}for\hspace*{2ex}$a_{3}\log n \leq \big|\,t\,\big| \leq a_{4}\,n^{3/2}$
\hfill(note $\hat{F}_{\bm\hat{A}} = \hat{\mathscr{F}}_{\!A}$ !).\\[3.5ex]
\textbf{Proof (by citations):}\\[0.8ex]
This theorem is a direct consequence of Theorem 2.1 and the comments from 
Section 3 of van Zwet \cite{vanZwet1982}. Some of these comments also follow 
from the next Lemma \ref{EWD_4_3_06}.\hspace*{1ex}\hfill$\Box$\\[4ex]
The estimate (\ref{EWD_4_3_05}) will play in this section \ref{EWD_Kap4_Sec3} the role that 
Cram\'er's condition (\ref{EWD_4_2_03}) played in the previous section \ref{EWD_Kap4_Sec2},
i.e. with the help of (\ref{EWD_4_3_05}) we will derive the conditions (\ref{EWD_3_1_11}), (\ref{EWD_3_1_14}) 
and (\ref{EWD_3_1_15}) for $B = \bm\hat{A}$.\\[2.8ex]
However, in order to be able to prove these conditions also for the submatrices 
$B \in N(8,\bm\hat{A})$ (or $B \in N(16,\bm\hat{A})$), 
we need the validity of the estimate (\ref{EWD_4_3_05}) for these submatrices as well.
The following proposition, which establishes this validity for a sufficient number of matrices $B$,
serves this purpose.\\[2.8ex]
Before that, we show a simple lemma that is useful for the proof of this proposition and 
the further procedure.\\[4ex]
\refstepcounter{DSBcount}
\label{EWD_4_3_06}
\textbf{\hspace*{-0.8ex}\theDSBcount\ Lemma}
\begin{enumerate}
\item\label{EWD_4_3_06_BWa}
Let $E$, $k$ and $t$ be positive constants. Then\\[2.5ex]
\begin{tabular}[t]{@{}l@{\hspace*{5.3ex}}l@{}} 
from&$\displaystyle{\sum\limits_{i=1}^{n}\,\big|\,e_{i} - \bm\bar{e}\,\big|^{k}\, \leq\, E n}$\\[4ex]
follows&$\displaystyle{\sum\limits_{i=1}^{n}\,\big|\,e_{i} - \bm\bar{e}\,\big|^{t}\, \leq\, 
E^{\,t/k}\,n^{1\, +\, \bigl((t/k)\, -\, 1\bigr)^{+}}}$.
\end{tabular}\vspace*{0.5ex}
\item\label{EWD_4_3_06_BWb}
Let $e$, $E$, $r$, $k$ and $t$ be positive constants such that $0 < r < k$. Then\\[2.5ex]
\begin{tabular}[t]{@{}l@{\hspace*{5.3ex}}l@{}} 
from&$\displaystyle{\sum\limits_{i=1}^{n}\,\big|\,e_{i} - \bm\bar{e}\,\big|^{r}\, \geq\, e n}$
\hspace*{2ex}and\hspace*{2ex}
$\displaystyle{\sum\limits_{i=1}^{n}\,\big|\,e_{i} - \bm\bar{e}\,\big|^{k}\, \leq\, E n}$\\[4ex]
follows&$\displaystyle{\sum\limits_{i=1}^{n}\,\big|\,e_{i} - \bm\bar{e}\,\big|^{t}\, \geq\, e' n}$
\end{tabular}\\[2.5ex]
\hspace*{12.1ex}with a positive constant $e'$ depending only on $e$, $E$, $r$, $k$ and $t$.\\[3ex]
In particular, we can assume $r = m = 2$ in (\ref{EWD_4_3_02}) and (\ref{EWD_4_3_03}) without 
loss of generality.\vspace*{0.5ex}
\item\label{EWD_4_3_06_BWc}
From (\ref{EWD_4_3_02}) and (\ref{EWD_4_3_03}) follows that there exist constants
$C_{1} > 0$ and $C_{2} > 0$ depending only on $e$, $E$, $d$, $D$ and $r$, $k$, $m$, $s$ such that\\[3ex]
\refstepcounter{DSBcount}
\label{EWD_4_3_07}
\text{\hspace*{-0.8ex}(\theDSBcount)}
\hspace*{4ex}
$\displaystyle{D_{\!A}^2 \leq C_{1}\,n^{-\,1\, +\, \bigl((4/k)\, -\, 1\bigr)^{+}\, 
+\, \bigl((4/s)\, -\, 1\bigr)^{+}}}$,\\[3ex]
\refstepcounter{DSBcount}
\label{EWD_4_3_08}
\text{\hspace*{-0.8ex}(\theDSBcount)}
\hspace*{4ex}
$\displaystyle{E_{\!A}^3 \leq C_{2}\,n^{-\,(3/2)\, +\, \bigl((5/k)\, -\, 1\bigr)^{+}\, 
+\, \bigl((5/s)\, -\, 1\bigr)^{+}}}$.
\end{enumerate}
\vspace*{4ex}
\textbf{Proof:} 
\begin{enumerate}
\item
\textbf{1. Case:} $t \leq k$.\\[2ex]
H{\"o}lder's inequality\index{H{\"o}lder's inequality!for finite sequences using length $\nu$} 
(\ref{EWD_0_1_05}) with $\nu = n$ and $p = \dfrac{k}{t}$ gives\\[2ex]
\hspace*{12.1ex}\begin{tabular}[t]{@{}l@{\hspace*{1ex}}c@{\hspace*{1ex}}l@{}}
$\displaystyle{\sum\limits_{i=1}^{n}\,\big|\,e_{i} - \bm\bar{e}\,\big|^{t}}$&
$\leq$&$\displaystyle{\Bigl(\,\sum\limits_{i=1}^{n}\,\big|\,e_{i} 
- \bm\bar{e}\,\big|^{k}\,\Bigr)^{t/k}\,n^{1 - (t/k)}}$\\[3ex]
&$\leq$&$E^{\,t/k}\,n$.
\end{tabular}\\[2.5ex]

\pagebreak

\textbf{2. Case:} $t \geq k$.\\[2ex]
\hspace*{12.1ex}\begin{tabular}[t]{@{}l@{\hspace*{1ex}}c@{\hspace*{1ex}}l@{}}
$\displaystyle{\sum\limits_{i=1}^{n}\,\big|\,e_{i} - \bm\bar{e}\,\big|^{t}}$&
$\leq$&$\displaystyle{\max\limits_{1\, \leq\, i\, \leq\, n}\,\big|\,e_{i} - \bm\bar{e}\,\big|^{t - k}\,
\sum\limits_{i=1}^{n}\,\big|\,e_{i} - \bm\bar{e}\,\big|^{k}}$\\[3ex]
&$\leq$&$\displaystyle{\Bigl(\,\sum\limits_{i=1}^{n}\,\big|\,e_{i} - \bm\bar{e}\,\big|^{k}\,\Bigr)^{1 + (t-k)/k}}$\\[3.5ex]
&$\leq$&$E^{\,t/k}\,n^{\,t/k}$.
\end{tabular}
\item
\textbf{1. Case:} $r \leq t$.\\[2ex]
From H{\"o}lder's inequality\index{H{\"o}lder's inequality!for finite sequences using length $\nu$} 
(\ref{EWD_0_1_05}) with $\nu = n$ and $p = \dfrac{t}{r}$ we obtain\\[2ex]
\hspace*{12.1ex}$\displaystyle{e n\, \leq\, \sum\limits_{i=1}^{n}\,\big|\,e_{i} - \bm\bar{e}\,\big|^{r}
\, \leq\, \Bigl(\,\sum\limits_{i=1}^{n}\,\big|\,e_{i} - \bm\bar{e}\,\big|^{t}\,\Bigr)^{r/t}\,n^{1 - (r/t)}}$\\[2.5ex]
and therefore\\[2ex]
\hspace*{12.1ex}$\displaystyle{e^{t/r} n\, \leq\, \sum\limits_{i=1}^{n}\,\big|\,e_{i} - \bm\bar{e}\,\big|^{t}}$.\\[2.5ex]
\textbf{2. Case:} $r > t$.\\[2ex]
H{\"o}lder's inequality\index{H{\"o}lder's inequality!for finite sequences}
with $p = \dfrac{k-t}{k-r}$ and $q = \dfrac{k-t}{r-t}$ gives\\[1.5ex]
\hspace*{12.1ex}\begin{tabular}[t]{@{}l@{\hspace*{1ex}}c@{\hspace*{1ex}}l@{}}
$e n$&$\leq$&$\displaystyle{\sum\limits_{i=1}^{n}\,\big|\,e_{i} - \bm\bar{e}\,\big|^{t/p}
\,\big|\,e_{i} - \bm\bar{e}\,\big|^{r-(t/p)}}$\\[3ex]
&$\leq$&$\displaystyle{\Bigl(\,\sum\limits_{i=1}^{n}\,\big|\,e_{i} - \bm\bar{e}\,\big|^{t}\,\Bigr)^{1/p}\,
\Bigl(\,\sum\limits_{i=1}^{n}\,\big|\,e_{i} - \bm\bar{e}\,\big|^{k}\,\Bigr)^{1/q}}$
\hspace*{9.2ex}(since $\Bigl(\,r - \dfrac{t}{p}\,\Bigr)\,q = k$)\\[3ex]
&$\leq$&$\displaystyle{\Bigl(\,\sum\limits_{i=1}^{n}\,\big|\,e_{i} - \bm\bar{e}\,\big|^{t}\,\Bigr)^{1/p}\,
(En)^{1/q}}$.
\end{tabular}\\[3ex]
Because of $\dfrac{p}{q} = p - 1$, it follows that\\[2ex]
\hspace*{12.1ex}$\displaystyle{e^{p}E^{1-p}n\, \leq\, \sum\limits_{i=1}^{n}\,\big|\,e_{i} - \bm\bar{e}\,\big|^{t}}$.
\item
Due to part \ref{EWD_4_3_06_BWb}) we can assume $r = m = 2$ without loss of generality.
Therefore we get\\[2ex]
\hspace*{12.1ex}$\displaystyle{\sigma_{\!A}^{2} = \dfrac{1}{n-1}\,
\sum\limits_{i=1}^{n}\,\bigl(\,e_{i} - \bm\bar{e}\,\bigr)^{2}\,
\sum\limits_{j=1}^{n}\,\bigl(\,d_{j} - \bm\bar{d}\,\bigr)^{2} \geq edn}$.\\[2.5ex]
Using part \ref{EWD_4_3_06_BWa}) we then obtain\\[2ex]
\hspace*{12.1ex}\begin{tabular}[t]{@{}l@{\hspace*{1ex}}c@{\hspace*{1ex}}l@{}}
$D_{\!A}^2$&$=$&$\displaystyle{\dfrac{1}{n}\,\dfrac{1}{\sigma_{\!A}^{4}}\,
\sum\limits_{i=1}^{n}\,\bigl(\,e_{i} - \bm\bar{e}\,\bigr)^{4}\,
\sum\limits_{j=1}^{n}\,\bigl(\,d_{j} - \bm\bar{d}\,\bigr)^{4}}$\\[4ex]
&$\leq$&$\dfrac{E^{\,4/k}D^{\,4/s}}{e^{\,2}\,{\,d^2}}\,n^{-\,1\, +\, \bigl((4/k)\, -\, 1\bigr)^{+}\, 
+\, \bigl((4/s)\, -\, 1\bigr)^{+}}$
\end{tabular}\\[2ex]
and analogously\\[2ex]
\hspace*{12.1ex}$\displaystyle{E_{\!A}^3\, \leq\, \dfrac{E^{\,5/k}D^{\,5/s}}{e^{\,5/2}\,d^{\,5/2}}\,n^{-\,(3/2)\, 
+\, \bigl((5/k)\, -\, 1\bigr)^{+}\, 
+\, \bigl((5/s)\, -\, 1\bigr)^{+}}}$.\hspace*{1ex}\hfill$\Box$
\end{enumerate}
\vspace*{3.5ex}
The announced proposition is now:\\[4ex]
\refstepcounter{DSBcount}
\label{EWD_4_3_09}
\textbf{\hspace*{-0.8ex}\theDSBcount\ Proposition}\\[0.8ex]
Suppose that the conditions and notations of the Theorem \ref{EWD_4_3_01} apply.\\[1.8ex]
Furthermore, we choose a $0 < \epsilon_{0} \leq 1$ and an $n_{0} \in \mathbb{N}$
such that for all matrices $A$ satisfying $\sigma_{\!A} > 0$, $\beta_{A} \leq \epsilon_{0}\,n$ 
and $n \geq n_{0}$ holds:\\[1.8ex]
\hspace*{12.1ex}$\big|\,\sigma_{\!B}^2 - 1\,\big| \leq \dfrac{1}{3}$
\hspace*{2ex}for $B \in N(17,\bm\hat{A})$\\[1.8ex]
(cf. Corollary \ref{EWD_3_3_07}). Then there exist positive constants $b_{1}$, $b_{2}$, $b_{3}$ and $b_{4}$ 
depending only on $e$, $E$, $d$, $D$, $\delta$ and $r$, $k$, $m$, $s$ such that in the case of\\[1.8ex]
\hspace*{12.1ex}$n \geq n_{1} = \max \Bigl\{\,n_{0},\,\dfrac{68}{\delta},\,19\,\Bigr\}$
\hspace*{2ex}and\hspace*{2ex}
$\beta_{A} \leq \epsilon_{0}\,n$\\[1.8ex]
we have\\[2ex]
\refstepcounter{DSBcount}
\label{EWD_4_3_10}
\text{\hspace*{-0.8ex}(\theDSBcount)}
\hspace*{2.8ex}
$\displaystyle{\big|\,\hat{F}_{B}(t)\,\big|\, \leq\,
b_{1}\,n^{- b_{2}\log n}}$
\hspace*{3ex}for\hspace*{2ex}$b_{3}\log n \leq \big|\,t\,\big| \leq b_{4}\,n^{3/2}$
\hspace*{2ex}and\hspace*{2ex}
$B \in N(17,\bm\hat{A})$.\\[3.5ex]
\textbf{Proof:}
\begin{enumerate}
\item\label{EWD_4_3_09_a}
In the first part of the proof, we show that the conditions (\ref{EWD_4_3_02}), (\ref{EWD_4_3_03}) and (\ref{EWD_4_3_04})
also hold for matrices $\sigma_{\!A}\,B$ with $B \in N(17,\bm\hat{A})$
(but with other constants $e'$, $E'$, $d'$, $D'$, $\delta'$, which depend only on $e$, $E$, $d$, $D$, $\delta$ 
and $r$, $k$, $m$, $s$, and with another value $\zeta'$, which depends only on $\zeta$).
Because of Lemma \ref{EWD_4_3_06}, \ref{EWD_4_3_06_BWb}) we can assume $r = 2 = m$ without loss of generality.\\[2.8ex]
The matrices $\sigma_{\!A}\,B$, where $B \in M(l,\bm\hat{A})$, have the form\\[2ex]
\hspace*{12.1ex}$\sigma_{\!A}\,B = (g_{i}f_{j})_{1\, \leq\, i, j\, \leq\, n - l}$\\[2.3ex] 
with $g_{i} \in \bigl\{\,e_{k} - \bm\bar{e}\,:\, k = 1,\ldots,n\,\bigr\}$ and
$f_{j} \in \bigl\{\,d_{k} - \bm\bar{d}\,:\, k = 1,\ldots,n\,\bigr\}$ (without repetitions!). 
Furthermore, we define\\[2ex]
\hspace*{12.1ex}$\displaystyle{\bm\bar{g} = \dfrac{1}{n - l}\,\sum\limits_{i=1}^{n - l}\,g_{i}}$
\hspace*{2ex}and\hspace*{2ex}
$\displaystyle{\bm\bar{f} = \dfrac{1}{n - l}\,\sum\limits_{j=1}^{n - l}\,f_{j}}$.\\[2ex]
Using these abbreviations, we now determine the items $e'$, $E'$, $d'$, $D'$, $\delta'$ and $\zeta'$.
\vspace*{1ex}
\renewcommand\labelenumii{\textbf{\theenumii)}}
\begin{enumerate}
\makeatletter
\renewcommand*\p@enumii{}
\makeatother
\item\label{EWD_4_3_09_a_i} \textbf{Computation of {\boldmath $E'$} and {\boldmath $D'$}:}\\[0.8ex] 
An application of H{\"o}lder's inequality\index{H{\"o}lder's inequality!for finite sequences using length $\nu$}
(\ref{EWD_0_1_05}) with $\nu = 2$ and $p = k$ to $\bigl(\,|\,g_{i}\,| + |\,\bm\bar{g}\,|\,\bigr)^{k}$
and with $\nu = n-l$ and $p = k$ to $\bigl(\,|\,g_{1}\,| + \ldots + |\,g_{n-l}\,|\,\bigr)^{k}$ gives\\[2ex]
\hspace*{6.9ex}\begin{tabular}[t]{@{}l@{\hspace*{1ex}}c@{\hspace*{1ex}}l@{\hspace*{-6.7ex}}l@{}}
$\displaystyle{\sum\limits_{i=1}^{n-l}\,\big|\,g_{i} - \bm\bar{g}\,\big|^{k}}$&
$\leq$&$\displaystyle{2^{k-1}\,\sum\limits_{i=1}^{n-l}\,\big|\,g_{i}\,\big|^{k}\,
+\, 2^{k-1}\,(n - l)\,\dfrac{1}{(n - l)^k}\,\Bigl(\,\sum\limits_{i=1}^{n - l}\,\big|\,g_{i}\,\big|\,\Bigr)^k}$\\[3ex]
&$\leq$&$\displaystyle{2^{k}\,\sum\limits_{i=1}^{n-l}\,\big|\,g_{i}\,\big|^{k}}$\\[3ex]
&$\leq$&$\displaystyle{2^{k}\,\sum\limits_{i=1}^{n}\,\big|\,e_{i} - \bm\bar{e}\,\big|^{k}}$\\[2.5ex]
&$\leq$&$\displaystyle{2^{k}\,En}$
&(due to (\ref{EWD_4_3_02}))\\[2ex]
&$\leq$&$\displaystyle{10\,(2^{k}\,E)\,(n - l)}$\hspace*{4ex}for $0 \leq l \leq 17$.
\end{tabular}\\[2.5ex]
For the last inequality, we used that $n \leq 10\,(n - l)$ 
holds for $n \geq 19$ and $0 \leq l \leq 17$.
Therefore $E' = 10\,(2^{k}\,E)$ is suitable.\\[1.5ex]
We proceed in the same way with the $f_{j}$ and obtain $D' = 10\,(2^{s}\,D)$.
\vspace*{1ex}
\item\label{EWD_4_3_09_a_ii} \textbf{Computation of {\boldmath $e'$} and {\boldmath $d'$}:}\\[0.8ex]
Since we assume \mbox{\rule[0ex]{0ex}{2ex}$\sigma_{\!B}^2 \geq \dfrac{2}{3}$} and since we have 
$(n - 1) \leq 18\,\Bigl(\,(n - 1) - l\,\Bigr)$ for $n \geq 19$ 
and \mbox{\rule[0ex]{0ex}{2.8ex}$0 \leq l \leq 17$}, we obtain\\[2ex]
\hspace*{6.9ex}\begin{tabular}[t]{@{}l@{\hspace*{1ex}}c@{\hspace*{1ex}}l@{}}
$\displaystyle{\sum\limits_{i=1}^{n - l}\,\bigl(\,g_{i} - \bm\bar{g}\,\bigr)^{2}\,
\sum\limits_{j=1}^{n - l}\,\bigl(\,f_{j} - \bm\bar{f}\,\bigr)^{2}}$&
$=$&$(n - l - 1)\,\sigma_{\!A}^2\,\sigma_{\!B}^2$\\[2.5ex]
&$\geq$&$\dfrac{1}{18}\,(n - 1)\,\sigma_{\!A}^2\,\dfrac{2}{3}$\\[2.5ex]
&$=$&$\displaystyle{\dfrac{1}{27}\,\sum\limits_{i=1}^{n}\,\bigl(\,e_{i} - \bm\bar{e}\,\bigr)^{2}\,
\sum\limits_{j=1}^{n}\,\bigl(\,d_{j} - \bm\bar{d}\,\bigr)^{2}}$
\end{tabular}\\
\hspace*{6.9ex}\begin{tabular}[t]{@{}l@{\hspace*{1ex}}c@{\hspace*{1ex}}l@{}}
\hspace*{28.2ex}&$\geq$&$\dfrac{1}{27}\,edn^2$\hspace*{4ex}for $0 \leq l \leq 17$.
\end{tabular}\\[2.5ex]
For the last inequality, (\ref{EWD_4_3_02}) and (\ref{EWD_4_3_03}) were used.\\[1.5ex]
It follows that\\[2ex]
\hspace*{6.9ex}\begin{tabular}[t]{@{}l@{\hspace*{1ex}}c@{\hspace*{1ex}}l@{}}
$\displaystyle{\sum\limits_{i=1}^{n - l}\,\bigl(\,g_{i} - \bm\bar{g}\,\bigr)^{2}}$&
$\geq$&$\displaystyle{\dfrac{1}{27}\,edn^2\,\Bigl(\,
\sum\limits_{j=1}^{n - l}\,\bigl(\,f_{j} - \bm\bar{f}\,\bigr)^{2}\,\Bigr)^{-1}}$\\[3ex]
&$\geq$&$\displaystyle{\dfrac{1}{27}\,edn^2\,\Bigl(\,(D')^{2/s}\,(n - l)\,\Bigr)^{-1}}$\\[3ex]
&$\geq$&$\displaystyle{\dfrac{1}{27}\,ed\,(D')^{-\,2/s}\,(n - l)}$.
\end{tabular}\\[3ex]
For the penultimate inequality, Lemma \ref{EWD_4_3_06}, \ref{EWD_4_3_06_BWa}) with
$k = s > 2$, $E = D'$ (cf. part \ref{EWD_4_3_09_a_i})) and $t = 2$ was used.\\[1.5ex]
Therefore $e' = \dfrac{1}{27}\,ed\,(D')^{-\,2/s}$ and 
analogously $d' = \dfrac{1}{27}\,ed\,(E')^{-\,2/k}$ are suitable.
\vspace*{2ex}
\item\label{EWD_4_3_09_a_iii} \textbf{Computation of {\boldmath $\delta'$} and {\boldmath $\zeta'$} :}\\[0.8ex]
The translation invariance of $\lambda$ yields\\[2ex]
\hspace*{6.9ex}\begin{tabular}[t]{@{}c@{\hspace*{1ex}}l@{}}
&$\displaystyle{\lambda\Bigl(\,\bigcup_{j = 1}^{n}\,\Bigl\{\,x \in \mathbb{R}\,:\,\big|\,x - d_{j}\,\big| < \zeta\,
\Bigr\}\,\Bigr)}$\\[3ex]
$=$&$\displaystyle{\lambda\Bigl(\,\bigcup_{j = 1}^{n}\,\Bigl\{\,x \in \mathbb{R}\,:\,
\big|\,x - (\,d_{j} - \bm\bar{d}\,)\,\big| < \zeta\,
\Bigr\}\,\Bigr)}$\\[3ex]
$\leq$&$\displaystyle{\lambda\Bigl(\,\bigcup_{j = 1}^{n-l}\,\Bigl\{\,x \in \mathbb{R}\,:\,
\big|\,x - f_{j}\,\big| < \zeta\,
\Bigr\}\,\Bigr) + l\,2\zeta}$.
\end{tabular}\\[3ex]
For $n \geq \dfrac{68}{\delta}$ it follows that\\[2ex]
\hspace*{6.9ex}\begin{tabular}[t]{@{}c@{\hspace*{1ex}}l@{\hspace*{15.7ex}}l@{}}
&$\displaystyle{\lambda\Bigl(\,\bigcup_{j = 1}^{n-l}\,\Bigl\{\,x \in \mathbb{R}\,:\,
\big|\,x - f_{j}\,\big| < \zeta\,
\Bigr\}\,\Bigr)}$\\[3ex]
$\geq$&$\displaystyle{\lambda\Bigl(\,\bigcup_{j = 1}^{n}\,\Bigl\{\,x \in \mathbb{R}\,:\,\big|\,x - d_{j}\,\big| < \zeta\,
\Bigr\}\,\Bigr) - l\,2\zeta}$\\[3ex]
$\geq$&$\delta n \zeta - l\,2\zeta$
&(because of (\ref{EWD_4_3_04}))\\[2ex]
$\geq$&$\bigl(\,\delta n - 34\,\bigr)\,\zeta$\\[2ex]
$\geq$&$\dfrac{\delta}{2}\,n\,\zeta$\hspace*{4ex}for $0 \leq l \leq 17$.
\end{tabular}\\[2.5ex]
Since\\[1.5ex]
\hspace*{6.9ex}$7\,\dfrac{\log\ n}{n^{3/2}} \geq \dfrac{\log\ (n-l)}{(n-l)^{3/2}}$
\hspace*{4ex}for $n \geq 19$ and $0 \leq l \leq 17$,\\[2ex]
we get all in all\\[1.5ex]
\hspace*{6.9ex}\begin{tabular}[t]{@{}c@{\hspace*{1ex}}l@{\hspace*{19.3ex}}l@{}}
&$\displaystyle{\lambda\Bigl(\,\bigcup_{j = 1}^{n-l}\,\Bigl\{\,x \in \mathbb{R}\,:\,
\big|\,x - f_{j}\,\big| < 7\,\zeta\,
\Bigr\}\,\Bigr)}$\\[3ex]
$\geq$&$\displaystyle{\lambda\Bigl(\,\bigcup_{j = 1}^{n-l}\,\Bigl\{\,x \in \mathbb{R}\,:\,\big|\,x - f_{j}\,\big| < \zeta\,
\Bigr\}\,\Bigr)}$\\[3ex]
$\geq$&$\dfrac{\delta}{2}\,n\,\zeta$\\[2.5ex]
$=$&$\dfrac{\delta}{14}\,n\,(7\,\zeta)$\\[2.5ex]
$\geq$&$\dfrac{\delta}{14}\,(n-l)\,(7\,\zeta)$\hspace*{4ex}for $0 \leq l \leq 17$
\end{tabular}\\[2ex]
and\\[1.5ex]
\hspace*{6.9ex}$\zeta' = 7\,\zeta \geq 7\,\dfrac{\log\ n}{n^{3/2}} \geq \dfrac{\log\ (n-l)}{(n-l)^{3/2}}$
\hspace*{4ex}for $n \geq \max \Bigl\{\,\dfrac{68}{\delta},\,19\,\Bigr\}$ and $0 \leq l \leq 17$.\\[1.7ex]
Thus, $\delta' = \dfrac{\delta}{14}$ and $\zeta' = 7\,\zeta$ are suitable.
\end{enumerate}
\item\label{EWD_4_3_09_b}
The existence of positive constants $a_{1}'$, $a_{2}'$, $a_{3}'$ and $a_{4}'$ depending only on 
$e$, $E$, $d$, $D$, $\delta$ and $r$, $k$, $m$, $s$ such that\\[2ex]
\hspace*{6.9ex}\begin{tabular}[t]{@{}l@{\hspace*{4ex}}l@{}}
$\displaystyle{\big|\,\hat{F}_{\widehat{\sigma_{\!A}\,B}}(t)\,\big|\, \leq\,
a_{1}'\,(n - l)^{- a_{2}'\log\ (n - l)}}$&
for\hspace*{2ex}$a_{3}'\log\ (n - l) \leq \big|\,t\,\big| \leq a_{4}'\,(n - l)^{3/2}$\\[1ex]
&and\hspace*{2ex}$B \in M(l,\bm\hat{A})$, $0 \leq l \leq 17$
\end{tabular}\\[2.5ex]
now follows from part \ref{EWD_4_3_09_a}) and Theorem \ref{EWD_4_3_01}. 
Using $\widehat{\sigma_{\!A}\,B} = \bm\hat{B}$ we get from this that\\[2ex]
\hspace*{6.9ex}\begin{tabular}[t]{@{}l@{\hspace*{1ex}}c@{\hspace*{0.8ex}}l@{\hspace*{1.1ex}}l@{\hspace*{-23.8ex}}r@{}}
$\big|\,\hat{F}_{B}(t)\,\big|$&$\leq$&$\Big|\,E\Bigl(\,\exp\bigl\{it(T_{B} - \mu_{B})\bigr\}\,\Bigr)\,\Big|$&
&(since $\big|\exp\bigl(\,it\mu_{B}\,\bigr)\big| \leq 1$)\\[2.5ex]
&$=$&$\big|\,\hat{F}_{\bm\hat{B}}(\sigma_{\!B}\,t)\,\big|$&
&(since $T_{\bm\hat{B}} = (T_{B} - \mu_{B})/\sigma_{\!B}$)\\[2.5ex]
&$\leq$&$\displaystyle{a_{1}'\,(n - l)^{- a_{2}'\log\ (n - l)}}$&
for\hspace*{2ex}$\dfrac{a_{3}'\log\ (n - l)}{\sigma_{\!B}} \leq \big|\,t\,\big| \leq 
\dfrac{a_{4}'\,(n - l)^{3/2}}{\sigma_{\!B}}$\\[1.5ex]
&&&and\hspace*{2ex}$B \in M(l,\bm\hat{A})$, $0 \leq l \leq 17$.
\end{tabular}\\[2.2ex]
Because of
\begin{itemize}
\item
$n \leq 10\,(n - l)$\hspace*{4ex}for $n \geq 19$ and $0 \leq l \leq 17$, 
\item
$\log\ (n -l) \geq \dfrac{1}{5}\,\log\ n$\hspace*{4ex}for $n \geq 19$ and $0 \leq l \leq 17$,\vspace*{0.5ex} 
\item
$\sqrt{\dfrac{3}{4}} \leq \dfrac{1}{\sigma_{\!B}} \leq \sqrt{\dfrac{3}{2}}$,\hspace*{4ex}and\vspace*{1ex}
\item
$f(x) = c_{1}\,x^{-\,c_{2}\,\log\,x} = c_{1}\,e^{-\,c_{2}\,(\log\,x)^2}$\hspace*{4ex}for $x > 0$
\end{itemize}
\vspace*{0.5ex}
we finally obtain Proposition \ref{EWD_4_3_09} with\\[1.5ex]
\hspace*{6.9ex}$b_{1} = a_{1}'$,\hspace*{4ex}$b_{2} = \dfrac{a_{2}'}{25}$,\hspace*{4ex}$b_{3} 
= \sqrt{\dfrac{3}{2}}\,a_{3}'$
\hspace*{2ex}and\hspace*{2ex} 
$b_{4} = \sqrt{\dfrac{3}{4}}\,\Bigl(\,\dfrac{1}{10}\,\Bigr)^{3/2}\,a_{4}'$.\hspace*{1ex}\hfill$\Box$
\end{enumerate}
\vspace*{3.5ex}
As Lemma \ref{EWD_4_1_07}, \ref{EWD_4_1_07_BWc}) shows, we need an additional estimate of the \textbf{derivative} 
of $\hat{F}_{B}$, when proving the condition (\ref{EWD_3_1_15}).\\[2.8ex]
This estimate is given in the following proposition.\\[4ex]
\refstepcounter{DSBcount}
\label{EWD_4_3_11}
\textbf{\hspace*{-0.8ex}\theDSBcount\ Proposition}\\[0.8ex]
Suppose that the conditions and notations of the previous Proposition \ref{EWD_4_3_09} apply.\\[1.8ex] 
Furthermore, let $b_{1}$, $b_{2}$, $b_{3}$ and $b_{4}$ be the same positive constants as 
in this Proposition \ref{EWD_4_3_09}.\\[1.8ex]
Then in the case of $n \geq n_{1}$ and $\beta_{A} \leq \epsilon_{0}\,n$ holds:\\[2ex]
\refstepcounter{DSBcount}
\label{EWD_4_3_12}
\text{\hspace*{-0.8ex}(\theDSBcount)}
\hspace*{2.8ex}
$\displaystyle{\big|\,\hat{F}_{B}'(t)\,\big|\, \leq\,
b_{1}\,n^{- b_{2}\log n\,+ \!\textstyle{\frac{1}{2}}}}$
\hspace*{3ex}for\hspace*{2ex}$b_{3}\log n \leq \big|\,t\,\big| \leq b_{4}\,n^{3/2}$
\hspace*{2ex}and\hspace*{2ex}
$B \in N(16,\bm\hat{A})$.\\[3.5ex]
\textbf{Proof:}\\[0.8ex]
For $B \in M(l,\bm\hat{A})$, we denote by $B_{jk}$ the matrix that we obtain from $B$ 
by cancelling the $j$th row and the $k$th column. Thus, $B_{jk} \in M(l+1,\bm\hat{A})$.
\index{matrix!submatrix}\index{matrix!submatrix!$B_{jk}$}\\[2.8ex]
So we have:\\[2ex]
\hspace*{4ex}\begin{tabular}[t]{@{}l@{\hspace*{1ex}}c@{\hspace*{1ex}}l@{\hspace*{3.1ex}}r@{}}
$\big|\,\hat{F}_{B}'(t)\,\big|$&$=$&
$\Big|\,E\Bigl(\,T_{B}\,\exp\bigl\{itT_{B}\bigr\}\,\Bigr)\,\Big|$\\[2.5ex]
&$\leq$&$\displaystyle{\sum\limits_{j = 1}^{n - l}\, \Big|\,E\Bigl(\,b_{j\pi(j)}\,
\exp\bigl\{itT_{B}\bigr\}\,\Bigr)\,\Big|}$\\[3.5ex]
&$\leq$&$\displaystyle{\dfrac{1}{n - l}\, \sum\limits_{j, k = 1}^{n-l}\, \Big|\,E\Bigl(\,b_{j\pi(j)}\,
\exp\bigl\{itT_{B}\bigr\}\,\Big|\,\pi(j) = k\,\Bigr)\,\Big|}$
&(since $P\bigl(\,\pi(j) = k\,\bigr) = \dfrac{1}{n-l}$)\\[4ex]
&$=$&$\displaystyle{\dfrac{1}{n - l}\, \sum\limits_{j, k = 1}^{n-l}\, \Big|\,E\Bigl(\,b_{jk}\,
\exp\bigl\{itb_{jk}\bigr\}\,\exp\bigl\{itT_{B_{jk}}\bigr\}\,\Bigr)\,\Big|}$.
\end{tabular}\\[3ex]
The last equation follows, since $P\bigr(\,\pi = \bm{.} \,\,\big|\,\pi(j) = k\,\bigr)$ is the 
uniform distribution on the set of permutations 
$\pi_{0} \in \mathscr{P}_{n-l}$\index{permutation sets!$\mathscr{P}_{n-l}$} with $\pi_{0}(j) = k$.
Due to $\big|\exp\bigl(\,itb_{jk}\,\bigr)\big| \leq 1$, we further get\\[2ex]
\hspace*{4ex}\begin{tabular}[t]{@{}l@{\hspace*{1ex}}c@{\hspace*{1ex}}l@{\hspace*{2.1ex}}r@{}}
\hspace*{8ex}&$\leq$&$\displaystyle{\dfrac{1}{n - l}\, \sum\limits_{j, k = 1}^{n-l}\,\big|\,b_{jk}\,\big|\,
\Big|\,E\Bigl(\,\exp\bigl\{itT_{B_{jk}}\bigr\}\,\Bigr)\,\Big|}$\\[3.5ex]
&$\leq$&$\displaystyle{b_{1}\,n^{- b_{2}\log n}\,\dfrac{1}{n - l}\,
\sum\limits_{j, k = 1}^{n-l}\,\big|\,b_{jk}\,\big|}$
&(note that $B_{jk} \in N(17,\bm\hat{A})$)\\[3.5ex]
&$\leq$&$\displaystyle{b_{1}\,n^{- b_{2}\log n\,+ \!\textstyle{\frac{1}{2}}}}$
\hspace*{7ex}for\hspace*{2ex}$b_{3}\log n \leq \big|\,t\,\big| \leq b_{4}\,n^{3/2}$.
\end{tabular}\\[3ex]
The Proposition \ref{EWD_4_3_09} was used for the penultimate inequality. For the last inequality, 
we first applied H{\"o}lder's inequality\index{H{\"o}lder's inequality!for finite sequences using length $\nu$} 
(\ref{EWD_0_1_05}) with $\nu = (n - l)^2$ and $p = 2$, and then (\ref{EWD_3_1_01}), i.e.\\[2ex]
\hspace*{12.1ex}\begin{tabular}[t]{@{}l@{\hspace*{1ex}}c@{\hspace*{1ex}}l@{}}
$\displaystyle{\dfrac{1}{n - l}\,
\sum\limits_{j, k = 1}^{n-l}\,\big|\,b_{jk}\,\big|}$&$\leq$&
$\displaystyle{\dfrac{1}{n - l}\,\Bigl(\,(n - l)^2 \,\sum\limits_{j, k = 1}^{n-l}\,\,b_{jk}^2\,\Bigr)^{1/2}}$\\[3.5ex]
&$\leq$&$\displaystyle{\Bigl(\,\sum\limits_{j, k = 1}^{n}\,\,\bm\hat{a}_{jk}^2\,\Bigr)^{1/2}\, \leq\, n^{1/2}}$.
\end{tabular}\\[-2.8ex]
\hspace*{1ex}\hfill$\Box$\\[4ex]
After all these preparations, we can now derive from the conditions of van Zwet\index{conditions of van Zwet}
the requirements (\ref{EWD_3_1_11}), (\ref{EWD_3_1_14}), (\ref{EWD_3_1_15}) 
and thus the results of the Theorems \ref{EWD_3_1_10} and \ref{EWD_3_1_13}  
(but with {''constants''} $\mathcal{C}_{1}$, $\mathcal{C}_{2}$ and $\mathcal{C}_{3}$, which grow slightly with $n$).\\[2.8ex]
We obtain the following result, which is central to this chapter:\\[4ex]
\refstepcounter{DSBcount}
\label{EWD_4_3_13}
\textbf{\hspace*{-0.8ex}\theDSBcount\ Theorem}\index{Theorem!for simple linear rank statistics!main results}\\[0.8ex]
Let $A = (a_{ij})$ be an $n{\times}n-$matrix, whose elements have the form $a_{ij} = e_{i} \, d_{j}$. 
In addition, let $\mathscr{F}_{\!A}$, $e_{1,A}$, $e_{2,A}$ and
$D_{\!A}$, $E_{\!A}$ be defined as in chapter \ref{EWD_Kap3}.\\[2.8ex]
Furthermore, suppose that there exist positive constants $e$, $E$, $d$, $D$ and $\delta$ such that\\[2ex]
(\ref{EWD_4_3_02})
\hspace*{4ex}
\begin{tabular}{@{}l@{}}
$\displaystyle{\sum\limits_{i=1}^{n}\,\big|\,e_{i} - \bm\bar{e}\,\big|^{r}\, \geq\, e n}$,\hspace*{4ex}
$\displaystyle{\sum\limits_{i=1}^{n}\,\big|\,e_{i} - \bm\bar{e}\,\big|^{k}\, \leq\, E n}$\\[2.5ex]
for some $k > 2$ and $0 < r < k$.
\end{tabular}\\[3ex]
(\ref{EWD_4_3_03})
\hspace*{4ex}
\begin{tabular}{@{}l@{}}
$\displaystyle{\sum\limits_{j=1}^{n}\,\big|\,d_{j} - \bm\bar{d}\,\big|^{m}\, \geq\, d n}$,\hspace*{4ex}
$\displaystyle{\sum\limits_{j=1}^{n}\,\big|\,d_{j} - \bm\bar{d}\,\big|^{s}\, \leq\, D n}$\\[2.5ex]
for some $s > 2$ and $0 < m < s$.
\end{tabular}\\[3ex]
(\ref{EWD_4_3_04})
\hspace*{4ex}
\begin{tabular}{@{}l@{}}
$\displaystyle{\lambda\Bigl(\,\bigcup_{j = 1}^{n}\,\Bigl\{\,x \in \mathbb{R}\,:\,\big|\,x - d_{j}\,\big| < \zeta\,
\Bigr\}\,\Bigr)\, \geq\, \delta n \zeta}$\\[3ex]
for some $\zeta \geq n^{-3/2}\,\log\ n$.
\end{tabular}
\begin{enumerate}
\item\label{EWD_4_3_13_BWa}
Then there exist constants $\mathcal{E}_{1} > 0$ and $\mathcal{E}_{2} > 0$ depending only on
$e$, $E$, $d$, $D$, $\delta$ and $r$, $k$, $m$, $s$ such that\\[2ex]
\refstepcounter{DSBcount}
\label{EWD_4_3_14}
\text{\hspace*{-0.8ex}(\theDSBcount)}
\hspace*{2.8ex}
\begin{tabular}[t]{@{}l@{\hspace*{1ex}}c@{\hspace*{1ex}}l@{}}
$||\,\mathscr{F}_{\!A} - e_{1,A}\,||$&$\leq$&$\mathcal{E}_{1}\,(\log n)^2\,D_{\!A}^2$\\[1.5ex]
&$\leq$&$\displaystyle{\mathcal{E}_{2}\,(\log n)^2\,n^{-\,1\, +\, \bigl((4/k)\, -\, 1\bigr)^{+}\, 
+\, \bigl((4/s)\, -\, 1\bigr)^{+}}}$.
\end{tabular}\vspace*{0.5ex}
\item\label{EWD_4_3_13_BWb}
Let $\epsilon > 0$. Then there exist constants $\mathcal{E}_{3} > 0$ and $\mathcal{E}_{4} > 0$ depending only on
$e$, $E$, $d$, $D$, $\delta$, $r$, $k$, $m$, $s$ and $\epsilon$ such that\\[2ex]
\refstepcounter{DSBcount}
\label{EWD_4_3_15}
\text{\hspace*{-0.8ex}(\theDSBcount)}
\hspace*{2.8ex}
\begin{tabular}[t]{@{}l@{\hspace*{1ex}}c@{\hspace*{1ex}}l@{}}
$||\,\mathscr{F}_{\!A} - e_{2,A}\,||$&$\leq$&$\mathcal{E}_{3}\,n^{\epsilon}\,E_{\!A}^3$\\[1.5ex]
&$\leq$&$\displaystyle{\mathcal{E}_{4}\,n^{-\,(3/2)\, +\, \epsilon\, +\, \bigl((5/k)\, -\, 1\bigr)^{+}\, 
+\, \bigl((5/s)\, -\, 1\bigr)^{+}}}$.
\end{tabular}
\end{enumerate}
\vspace*{3.5ex}
\textbf{Proof:}
\begin{enumerate}
\item
At first, it should be noted that we can assume\\[2ex] 
\hspace*{12.1ex}$n \geq n_{1}$ and $\beta_{A} \leq \epsilon_{0}\,n$\hspace*{2ex} 
with $n_{1}$, $\epsilon_{0}$ 
from Proposition \ref{EWD_4_3_09}\\[2.5ex]
without loss of generality. To justify this, we proceed as in (\ref{EWD_3_7_07}), (\ref{EWD_3_7_08}) 
and (\ref{EWD_3_7_09}).\\[2.8ex]
Because of Theorem \ref{EWD_3_1_10}, Lemma \ref{EWD_4_1_07}, \ref{EWD_4_1_07_BWa}) and inequality (\ref{EWD_4_3_07}), 
it is now sufficient to show the existence of a constant $c_{1}$ depending only on $e$, $E$, $d$, $D$, $\delta$ 
and $r$, $k$, $m$, $s$ such that\\[2ex]
\refstepcounter{DSBcount}
\label{EWD_4_3_16}
\text{\hspace*{-0.8ex}(\theDSBcount)}
\hspace*{3.6ex}$\displaystyle{\int\limits_{\mathbb{R}}
\bigl(\,1 + |\,t\,|\,\bigr)\,\big|\,\hat{F}_{B}(t)\,\hat{U}_{D_{\!A}^2}(t)\,\big|\,dt\, \leq\,
c_{1}\,(\log n)^2}$\hspace*{4ex}for all $B \in N(8,\bm\hat{A})$.\\[2.5ex] 
Thereby $U$ is defined as in the section \ref{EWD_Kap4_Sec1} (see in particular (\ref{EWD_4_1_01}) 
and (\ref{EWD_4_1_02})) of this chapter \ref{EWD_Kap4}. Because of\\[2ex]
\hspace*{12.1ex}$\displaystyle{\int\limits_{\mathbb{R}}
\bigl(\,1 + |\,t\,|\,\bigr)\,\big|\,\hat{F}_{B}(t)\,\hat{U}_{D_{\!A}^2}(t)\,\big|\,dt\, \leq\,
2 + 2\,\int\limits_{\mathbb{R}}
|\,t\,|\,\big|\,\hat{F}_{B}(t)\,\hat{U}_{D_{\!A}^2}(t)\,\big|\,dt}$,\\[2.5ex]
it is also sufficient to prove (\ref{EWD_4_3_16}) with $|\,t\,|$ instead of $1 + |\,t\,|$.\\[2.8ex]
Using the constants $b_{1}$, $b_{2}$, $b_{3}$, $b_{4}$ from Proposition \ref{EWD_4_3_09} and the
abbreviation $\rho = D_{\!A}^2$ we obtain\\[2ex]
\hspace*{12.1ex}\begin{tabular}[t]{@{}l@{\hspace*{1ex}}l@{}}
&$\displaystyle{\int\limits_{\mathbb{R}}
|\,t\,|\,\big|\,\hat{F}_{B}(t)\,\hat{U}_{\rho}(t)\,\big|\,dt}$\\[4ex]
$\leq$&$\displaystyle{2\,\bigg\{\,\int\limits_{0}^{b_{3} \log n}\,t\,dt\, +\,
\int\limits_{b_{3} \log n}^{b_{4} n^{3/2}}\,t\,\big|\,\hat{F}_{B}(t)\,\big|\,dt\, +\,
\int\limits_{b_{4} n^{3/2}}^{\infty}\,t\,\big|\,\hat{U}_{\rho}(t)\,\big|\,dt\,\bigg\}}$\\[4ex]
$\leq$&$\displaystyle{b_{3}^2\,(\log n)^2\, +\, b_{1}\,n^{- b_{2}\log n}\,
\Bigl[\,b_{4}^2\,n^3 - b_{3}^2\,(\log n)^2\,\Bigr]\, +\, 2\,
\int\limits_{b_{4} n^{3/2}}^{\infty}\,t\,\big|\,\hat{U}_{\rho}(t)\,\big|\,dt}$.
\end{tabular}\\[2ex]
Because of $b_{2} > 0$, it follows that\\[2ex]
\hspace*{12.1ex}\begin{tabular}[t]{@{}l@{\hspace*{1ex}}l@{}}
&$\displaystyle{\lim\limits_{n \rightarrow \infty}\,
n^{- b_{2}\log n}\,
\Bigl[\,b_{4}^2\,n^3 - b_{3}^2\,(\log n)^2\,\Bigr]}$\\[2.5ex] 
$=$&$\displaystyle{\lim\limits_{n \rightarrow \infty}\,
\exp\bigl(\,- b_{2}\,(\log n)^2\,\bigr)\,\Bigl[\,b_{4}^2\,\exp\bigl(\,3 \log n\,\bigr) 
- b_{3}^2\,(\log n)^2\,\Bigr]\, =\, 0}$,
\end{tabular}\\[2.5ex]
so that only the third summand remains to be analyzed:\\[2ex]
\hspace*{12.1ex}\begin{tabular}[t]{@{}l@{\hspace*{1ex}}l@{\hspace*{31.2ex}}r@{}}
&$\displaystyle{\int\limits_{b_{4} n^{3/2}}^{\infty}\,t\,\big|\,\hat{U}(\rho\,t)\,\big|\,dt}$
&($\hat{U}_{\rho}(t) = \hat{U}(\rho\,t)$)\\[4.5ex]
$=$&$\displaystyle{\dfrac{1}{\rho^2}\,\int\limits_{b_{4} n^{3/2} \rho}^{\infty}\,v\,\big|\,\hat{U}(v)\,\big|\,dv}$
&(substitution\index{integration!by substitution} $v = t\,\rho$).
\end{tabular}\\[2ex]
Furthermore, since the inequality $\rho \geq \dfrac{1}{4n}$ is valid due to Lemma \ref{EWD_3_1_18}, \ref{EWD_3_1_18_BWa}),
we get with the definition $c_{2} = \dfrac{1}{4}\,b_{4}$\\[1.5ex]
\hspace*{12.1ex}\begin{tabular}[t]{@{}l@{\hspace*{1ex}}l@{\hspace*{20.2ex}}r@{}}
$\leq$&$\displaystyle{16\,n^2\,\int\limits_{c_{2} n^{1/2}}^{\infty}\,v\,\big|\,\hat{U}(v)\,\big|\,dv}$\\[4.5ex]
$\leq$&$\displaystyle{16\,n^2\,
\int\limits_{c_{2} n^{1/2}}^{\infty}\,\dfrac{v^6}{c_{2}^5\,n^{5/2}}\,\big|\,\hat{U}(v)\,\big|\,dv}$\\[5ex]
$\leq$&$\displaystyle{\dfrac{16}{c_{2}^5\,n^{1/2}}\,\,
\int\limits_{\mathbb{R}}\,v^6\,\big|\,\hat{U}(v)\,\big|\,dv \rightarrow 0}$
\hspace*{3ex}as $n \rightarrow \infty$&(cf. (\ref{EWD_4_1_04})).
\end{tabular}\\[3ex]
Thus, the part \ref{EWD_4_3_13_BWa}) is proved.
\item
Analogous to the proof of part \ref{EWD_4_3_13_BWa}) we can assume\\[2ex] 
\hspace*{12.1ex}$n \geq n_{1}$ and $\,\dfrac{\beta_{A}}{n} \leq E_{\!A} \leq \epsilon_{0}$\hspace*{2ex} 
with $n_{1}$, $\epsilon_{0}$ 
from Proposition \ref{EWD_4_3_09}\\[2.5ex]
without loss of generality. To justify this, we proceed as in (\ref{EWD_3_8_03}), (\ref{EWD_3_8_04}) 
and (\ref{EWD_3_8_05}).\\[2.8ex]
Because of Theorem \ref{EWD_3_1_13}, Lemma \ref{EWD_4_1_07}, inequality (\ref{EWD_4_3_08}) and\\[2ex]
\hspace*{12.1ex}$\displaystyle{\big|\,\hat{U}_{E_{\!A}^2}'(t)\,\big| = 
E_{\!A}^2\,\big|\,\hat{U}'(E_{\!A}^2\,t)\,\big| \leq E_{\!A}^2 \leq \epsilon_{0}^2 \leq 1}$
\hspace*{4ex}and\\[2ex]
\hspace*{12.1ex}$\displaystyle{\big|\,\hat{F}_{B}'(t)\,\big| =
\Big|\,E\bigl(\,T_{B}\,\exp\bigl\{itT_{B}\bigr\}\,\bigr)\,\Big| \leq E\bigl(\,\big|\,T_{B}\,\big|\,\bigr)
\leq \bigl(\,\sigma_{\!B}^2\,\bigr)^{1/2} \leq \sqrt{\dfrac{4}{3}}}$\\[2.5ex]
\hspace*{65ex}for all  $B \in N(16,\bm\hat{A})$,\\[2ex]
it is then sufficient to show the existence of constants $c_{3}$, $c_{4}$, $c_{5}$ and $c_{6}$
depending only on $e$, $E$, $d$, $D$, $\delta$, $r$, $k$, $m$, $s$ and $\epsilon$ such that\\[2.5ex]
\refstepcounter{DSBcount}
\label{EWD_4_3_17}
\text{\hspace*{-0.8ex}(\theDSBcount)}
\hspace*{3.6ex}$\displaystyle{\int\limits_{\mathbb{R}}
|\,t\,|\,\big|\,\hat{F}_{B}'(t)\,\hat{U}_{E_{\!A}^2}(t)\,\big|\,dt\, \leq\,
c_{3}\,(\log n)^2}$,\\[3ex]
\refstepcounter{DSBcount}
\label{EWD_4_3_18}
\text{\hspace*{-0.8ex}(\theDSBcount)}
\hspace*{3.6ex}$\displaystyle{\int\limits_{\mathbb{R}}
|\,t\,|\,\big|\,\hat{F}_{B}(t)\,\hat{U}_{E_{\!A}^2}'(t)\,\big|\,dt\, \leq\,
c_{4}\,(\log n)^2}$,\\[3ex]
\refstepcounter{DSBcount}
\label{EWD_4_3_19}
\text{\hspace*{-0.8ex}(\theDSBcount)}
\hspace*{3.6ex}$\displaystyle{\int\limits_{\mathbb{R}}
|\,t\,|\,\big|\,\hat{F}_{B}(t)\,\hat{U}_{E_{\!A}^2}(t)\,\big|\,dt\, \leq\,
c_{5}\,(\log n)^2}$,\\[3ex]
\refstepcounter{DSBcount}
\label{EWD_4_3_20}
\text{\hspace*{-0.8ex}(\theDSBcount)}
\hspace*{3.6ex}$\displaystyle{\int\limits_{\mathbb{R}}
t^2\,\big|\,\hat{F}_{B}(t)\,\hat{U}_{E_{\!A}^3\, n^{\epsilon/2}}(t)\,\big|\,dt\, \leq\,
c_{6}\,(\log n)^3}$
\hspace*{4ex}for all $B \in N(16,\bm\hat{A})$.\\[3ex]
The statements (\ref{EWD_4_3_17}), (\ref{EWD_4_3_18}) and (\ref{EWD_4_3_19}) are proved analogously
to (\ref{EWD_4_3_16}) using the Propositions \ref{EWD_4_3_09} and \ref{EWD_4_3_11} and\\[2ex]
\hspace*{12.1ex}$E_{\!A}^2 \geq 2^{-5/3}\,n^{-1}$ (cf. \ref{EWD_3_1_18}, \ref{EWD_3_1_18_BWa})),
$\big|\,\hat{U}_{E_{\!A}^2}'(t)\,\big| \leq 1$ and $\big|\,\hat{F}_{B}'(t)\,\big| \leq \sqrt{\dfrac{4}{3}}$.\\[2.8ex]
It remains to verify (\ref{EWD_4_3_20}).
As in part \ref{EWD_4_3_13_BWa}) we obtain with the abbreviation $\theta = E_{\!A}^3\, n^{\epsilon/2}$:\\[2ex]
\hspace*{12.1ex}\begin{tabular}[t]{@{}l@{\hspace*{1ex}}l@{}}
&$\displaystyle{\int\limits_{\mathbb{R}}
t^2\,\big|\,\hat{F}_{B}(t)\,\hat{U}_{\theta}(t)\,\big|\,dt}$\\[5ex]
$\leq$&$\displaystyle{\dfrac{2}{3}\,b_{3}^3\,(\log n)^3\, +\, \dfrac{2}{3}\,b_{1}\,n^{- b_{2}\log n}\,
\Bigl[\,b_{4}^3\,n^{9/2} - b_{3}^3\,(\log n)^3\,\Bigr]}$
\end{tabular}\\
\hspace*{12.1ex}\begin{tabular}[t]{@{}l@{\hspace*{1ex}}l@{}}
\hspace*{1.8ex}&$\displaystyle{+\ 2\,\int\limits_{b_{4} n^{3/2}}^{\infty}\,t^2\,\big|\,\hat{U}_{\theta}(t)\,\big|\,dt}$.
\end{tabular}\\[2.5ex]
Analogous to part \ref{EWD_4_3_13_BWa}), we get for the third summand\\[2ex]
\hspace*{12.1ex}\begin{tabular}[t]{@{}l@{\hspace*{1ex}}l@{\hspace*{30.3ex}}r@{}}
&$\displaystyle{\int\limits_{b_{4} n^{3/2}}^{\infty}\,t^2\,\big|\,\hat{U}(\theta\,t)\,\big|\,dt}$
&($\hat{U}_{\theta}(t) = \hat{U}(\theta\,t)$)\\[5ex]
$=$&$\displaystyle{\dfrac{1}{\theta^3}\,\int\limits_{b_{4} n^{3/2} \theta}^{\infty}\,v^2\,\big|\,\hat{U}(v)\,\big|\,dv}$
&(substitution\index{integration!by substitution} $v = t\,\theta$).
\end{tabular}\\[2ex]
Furthermore, since the inequality $\theta \geq \dfrac{n^{\epsilon/2}}{\sqrt{32}\,n^{3/2}}$ is valid 
due to Lemma \ref{EWD_3_1_18}, \ref{EWD_3_1_18_BWa}),
we get with the definition $c_{7} = \dfrac{1}{\sqrt{32}}\,b_{4}$\\[2ex]
\hspace*{12.1ex}\begin{tabular}[t]{@{}l@{\hspace*{1ex}}l@{\hspace*{9.4ex}}r@{}}
$\leq$&$\displaystyle{\dfrac{(\sqrt{32})^3\,n^{9/2}}{n^{3\epsilon/2}}\,
\int\limits_{c_{7} n^{\epsilon/2}}^{\infty}\,v^2\,\big|\,\hat{U}(v)\,\big|\,dv}$\\[5ex]
$\leq$&$\displaystyle{\dfrac{(\sqrt{32})^3\,n^{9/2}}{n^{3\epsilon/2}}\,
\int\limits_{c_{7} n^{\epsilon/2}}^{\infty}\,\dfrac{v^{\,2\, +\, (9/\epsilon)}}
{\bigl(\,c_{7}\,n^{\epsilon/2}\,\bigr)^{9/\epsilon}}\,\big|\,\hat{U}(v)\,\big|\,dv}$\\[5.5ex]
$\leq$&$\displaystyle{\dfrac{(\sqrt{32})^3}{c_{7}^{9/\epsilon}\,n^{3\epsilon/2}}\,\,
\int\limits_{\mathbb{R}}\,|\,v\,|^{\,2\, +\, (9/\epsilon)}\,\big|\,\hat{U}(v)\,\big|\,dv \rightarrow 0}$
\hspace*{3ex}as $n \rightarrow \infty$&(cf. (\ref{EWD_4_1_04})).
\end{tabular}\\[3ex]
Thus, (\ref{EWD_4_3_20}) is proved. All together, since 
$n^{\epsilon/2} (\log n)^3 \leq n^{\epsilon}$ for sufficiently large $n$, 
the part \ref{EWD_4_3_13_BWb}) is also proved.\hspace*{1ex}\hfill$\Box$
\end{enumerate}
\vspace*{3.5ex}
\refstepcounter{DSBcount}
\label{EWD_4_3_21}
\textbf{\hspace*{-0.8ex}\theDSBcount\ Remarks}
\begin{enumerate}
\item\label{EWD_4_3_21_BWa}
According to Lemma \ref{EWD_3_1_18}, \ref{EWD_3_1_18_BWe}) follows from (\ref{EWD_4_3_15})\\[2ex]
\hspace*{12.1ex}\begin{tabular}[t]{@{}l@{\hspace*{1ex}}c@{\hspace*{1ex}}l@{}}
$||\,\mathscr{F}_{\!A} - e_{1,A}\,||$&$\leq$&
$||\,\mathscr{F}_{\!A} - e_{2,A}\,||\, +\, ||\,e_{2,A} - e_{1,A}\,||$\\[2.5ex]
&$\leq$&$\mathcal{E}_{3}\,n^{\epsilon}\,E_{\!A}^3\, +\, \dfrac{1}{2}\,D_{\!A}^2$.
\end{tabular}\\[2.5ex]
Therefore, in the case of\\[2ex]
\refstepcounter{DSBcount}
\label{EWD_4_3_22}
\text{\hspace*{-0.8ex}(\theDSBcount)}
\hspace*{2.8ex}
$\displaystyle{\Bigl[\,\bigl((5/k)\, -\, 1\bigr)^{+}\, - \,\bigl((4/k)\, -\, 1\bigr)^{+}\,\Bigr]\, +\,
\Bigl[\,\bigl((5/s)\, -\, 1\bigr)^{+}\, - \,\bigl((4/s)\, -\, 1\bigr)^{+}\,\Bigr]\, <\, \dfrac{1}{2}}$,\\[2.5ex]
we even have\\[2ex]
\refstepcounter{DSBcount}
\label{EWD_4_3_23}
\text{\hspace*{-0.8ex}(\theDSBcount)}
\hspace*{2.8ex}
$\displaystyle{||\,\mathscr{F}_{\!A} - e_{1,A}\,||\, \leq\, \mathcal{E}_{2}'\,
n^{-\,1\, +\, \bigl((4/k)\, -\, 1\bigr)^{+}\, +\, \bigl((4/s)\, -\, 1\bigr)^{+}}}$\\[2.5ex]
due to Lemma \ref{EWD_4_3_06}, \ref{EWD_4_3_06_BWc}), 
where $\mathcal{E}_{2}' > 0$ depends only on $e$, $E$, $d$, $D$, $\delta$ and $r$, $k$, $m$, $s$.
\item\label{EWD_4_3_21_BWb}
If, compared to Lemma \ref{EWD_3_1_18}, \ref{EWD_3_1_18_BWa}), we additionally have\\[2.2ex]
\hspace*{12.1ex}$E_{\!A}^3 \geq C\,n^{-\,3/2\, +\, \eta}$, where $C, \eta > 0$,\\[2.3ex]
then the proof of (\ref{EWD_4_3_20}) shows that (\ref{EWD_4_3_15}) can be strengthened as follows:\\[2ex]
\refstepcounter{DSBcount}
\label{EWD_4_3_24}
\text{\hspace*{-0.8ex}(\theDSBcount)}
\hspace*{2.8ex}
\begin{tabular}[t]{@{}l@{\hspace*{1ex}}c@{\hspace*{1ex}}l@{}}
$||\,\mathscr{F}_{\!A} - e_{2,A}\,||$&$\leq$&$\mathcal{E}_{3}'\,(\log n)^3\,E_{\!A}^3$\\[1.5ex]
&$\leq$&$\displaystyle{\mathcal{E}_{4}'\,(\log n)^3\,n^{-\,(3/2)\, +\, \bigl((5/k)\, -\, 1\bigr)^{+}\, 
+\, \bigl((5/s)\, -\, 1\bigr)^{+}}}$,
\end{tabular}\\[2.5ex]
where $\mathcal{E}_{3}' > 0$ 
and $\mathcal{E}_{4}' > 0$ depend only on $e$, $E$, $d$, $D$, $\delta$, $r$, $k$, $m$, $s$ 
and $C$, $\eta$.
\item\label{EWD_4_3_21_BWc}
We refer again to the Remark \ref{EWD_3_1_17}, \ref{EWD_3_1_17_BWa}). 
This remark also applies analogously to the Theorem \ref{EWD_4_3_13} and the two remarks above.
\item\label{EWD_4_3_21_BWd}
In the two-sample case\index{regression constants!two-sample statistic}, 
i.e. in the case $e_{1} = \ldots = e_{m} = 1$ and $e_{m+1} = \ldots = e_{n} = 0$,
the corresponding condition of (\ref{EWD_4_3_02}) is:\\[2ex]
\refstepcounter{DSBcount}
\label{EWD_4_3_25}
\text{\hspace*{-0.8ex}(\theDSBcount)}
\hspace*{2.8ex}
There exists an $\eta > 0$ such that $\eta \leq \dfrac{m}{n} \leq 1 - \eta$.\\[2.5ex]
Since (\ref{EWD_4_3_02}) follows from this even for any $0 < r < k$, the $(\log n)^2$ can always be omitted 
in the two-sample case in (\ref{EWD_4_3_14}) (cf. Remark \ref{EWD_4_3_21_BWa}) above).\\[2.8ex]
Furthermore, it is known from the literature that in the two-sample case the $n^{\epsilon}$ 
from (\ref{EWD_4_3_15}) (or $(\log n)^3$ from (\ref{EWD_4_3_24})) can either be omitted completely 
or replaced by smaller terms. 
We refer to the papers of Bickel and van Zwet (\cite{10.1214/aos/1176344305}, Corollary 2.1, 
and second remark, p. 949 centre) and Robinson (\cite{10.1214/aos/1176344306}, Theorem and, 
to weaken his condition (c), the paper \cite{10.1214/aos/1176345078}, p. 861, Remark, and p. 852, 
Condition (B), (C)).
\item\label{EWD_4_3_21_BWe}
The condition (\ref{EWD_4_3_04}) is necessary to prevent the scores $d_{1},\ldots,d_{n}$ from 
being too close together in a few points and thus giving the distribution of $\mathscr{T}_{A}$ a
grid character that is too strong (such as $S_{n}$ or $X_{i}$ in the Remark \ref{EWD_1_1_06}, 
\ref{EWD_1_1_06_BWc})).\\[2.8ex]
As an example where this happens and therefore the convergence rates of the theorem do not hold, 
the following sequence of $n{\times}n-$matrices $A^{(n)} = \bigl(\,e_{i}^{(n)}\,d_{j}^{(n)}\,\bigr)$ is given:\\[2ex]
\hspace*{6.9ex}$\displaystyle{e_{i}^{(n)} =
\left\{
\begin{array}{ll@{}}
1& \hspace*{3ex}
\text{for}\ \ \ i\,\leq\,\Bigl\lfloor\,\dfrac{n}{2}\,\Bigr\rfloor\\[2ex]
0& \hspace*{3ex}
\text{for}\ \ \ i\,>\,\Bigl\lfloor\,\dfrac{n}{2}\,\Bigr\rfloor
\end{array}  \right.}$ ,
\hspace*{6.9ex}$\displaystyle{d_{j}^{(n)} =
\left\{
\begin{array}{rl@{}}
- 1& \hspace*{3ex}
\text{for}\ \ \ j\,\leq\,\Bigl\lfloor\,\dfrac{n}{2}\,\Bigr\rfloor\\[2ex]
1& \hspace*{3ex}
\text{for}\ \ \ j\,>\,\Bigl\lfloor\,\dfrac{n}{2}\,\Bigr\rfloor
\end{array}  \right.}$ ,\\[3ex]
({''Median test''}, cf. \cite{sidak1999theory}, Section 4.1.1, page 97).
\index{regression constants!Median test}\index{scores!Median test}\\[2.8ex]
For $n \geq 2$ and $k \geq 0$ then holds\\[2ex]
\hspace*{6.9ex}$\displaystyle{\dfrac{n}{3^{k}}\, \leq\, 
\sum\limits_{i=1}^{n}\,\big|\,e_{i}^{(n)} - \bm\bar{e}^{(n)}\,\big|^{k}\, \leq\, n}$,
\hspace*{9.1ex}$\displaystyle{\Bigl(\,\dfrac{2}{3}\,\Bigr)^{k}\,n\, \leq\, 
\sum\limits_{j=1}^{n}\,\big|\,d_{j}^{(n)} - \bm\bar{d}^{(n)}\,\big|^{k}\, \leq\, 2^{k}\,n}$.\\[2.5ex]
On the other hand, however, we get (analogously to the Remark \ref{EWD_1_1_06}, \ref{EWD_1_1_06_BWc}))\\[2ex] 
\hspace*{6.9ex}\begin{tabular}[t]{@{}l@{\hspace*{1ex}}l@{\hspace*{-3.9ex}}r@{}}
&$||\,\mathscr{F}_{\!A^{(4n)}} - e_{1,A^{(4n)}}\,||$\hspace*{2ex}bzw.\hspace*{2ex}
$||\,\mathscr{F}_{\!A^{(4n)}} - e_{2,A^{(4n)}}\,||$\\[2.5ex]
$\geq$&$\dfrac{1}{2}\,P\bigl(\,\mathscr{T}_{A^{(4n)}} = 0\,\bigr)$\\[2.5ex]
$=$&$\dfrac{1}{2}\,P\bigl(\,T_{A^{(4n)}} = 0\,\bigr)$
&(since $\mu_{A^{(4n)}} = 0$ and $0 \cdot \sigma_{\!A^{(4n)}} = 0$)\\[2.5ex]
$=$&$\displaystyle{\dfrac{1}{2}\, \dfrac{\dbinom{2n}{n}^{2}}{\dbinom{4n}{2n}} \sim 
\dfrac{1}{\sqrt{2 \pi n}}}$\hspace*{4ex}for $n \rightarrow \infty$
&(Stirling's formula\index{Stirling's formula}).
\end{tabular}
\end{enumerate}
\vspace*{1.5ex}

\section{Applications to the case of exact and approximating scores}\label{EWD_Kap4_Sec4}

Let $J: (0,\,1) \rightarrow \mathbb{R}$ be an integrable function. Furthermore, we denote by
$U_{j:n}$\index{order statistic!$U_{j:n}$} the $j$th
order statistic\index{order statistic} in a random sample of size $n$
from the uniform distribution on $(0,\,1)$.\\[2.8ex]
In the following we consider the scores\index{scores}\\[2ex]
\refstepcounter{DSBcount}
\label{EWD_4_4_01}
\text{\hspace*{-0.8ex}(\theDSBcount)}
\hspace*{4ex}
$d_{j} = E\bigl(\,J(\,U_{j:n}\,)\,\bigr)$,\hspace*{4ex}$j = 1,\ldots,n$\hspace*{4ex}({''exact 
scores\index{scores!exact}''})\\[2.5ex]
and\\[2ex]
\refstepcounter{DSBcount}
\label{EWD_4_4_02}
\text{\hspace*{-0.8ex}(\theDSBcount)}
\hspace*{4ex}
$d_{j} = J\Bigl(\,\dfrac{j}{n+1}\,\Bigr)$,\hspace*{5.2ex}
$j = 1,\ldots,n$\hspace*{4ex}({''approximating 
scores\index{scores!approximating}''}).\\[2.5ex]
We investigate which conditions must be placed on $J$ in each of the two upper cases 
so that the scores $d_{j}$ fulfil the conditions (\ref{EWD_4_3_03}) and (\ref{EWD_4_3_04}).\\[2.8ex]
In addition, we analyze for which $J$ the expansions $e_{1,A}$ and $e_{2,A}$ can be replaced by the
following somewhat simpler expansions (with integrals)
\index{Edgeworth expansion!$e_{1,A}^{I}$, $e_{2,A}^{I}$}\\[3ex]
$\begin{array}{@{\hspace*{7ex}}l@{\hspace*{0.8ex}}c@{\hspace*{0.8ex}}l} 
e_{1,A}^{I}(x)&=&\Phi(x) - \psi(x)\, \dfrac{\xi_{1,A}}{6}\, (x^2-1)\ \ \ \ \text{and}\\[2.5ex]
e_{2,A}^{I}(x)&=&\Phi(x) - \psi(x)\, \biggl\{\, \dfrac{\xi_{1,A}}{6}\, (x^2-1) +   \dfrac{\xi_{2,A}}{24}\, (x^3-3 x)
+  \dfrac{\xi_{1,A}^{2}}{72}\, (x^5-10x^3 + 15x)\,\biggr\}\ ,
\end{array}$\\[3.2ex]
where\\[2ex] 
\refstepcounter{DSBcount}
\label{EWD_4_4_03}
\text{\hspace*{-0.8ex}(\theDSBcount)}
\hspace*{4ex}
$\displaystyle{\bm\hat{e}_{i} = 
\dfrac{e_{i} - \bm\bar{e}}{\sqrt{\sum\limits_{l=1}^{n}\,\big(\,e_{l} - \bm\bar{e}\,\big)^{2}\,}}}$ ,
\hspace*{5ex}$\displaystyle{\bm\hat{J} = 
\dfrac{J - \int\limits_{0}^{1} J(t)\,dt}{\sqrt{\int\limits_{0}^{1} J^2(t)\,dt - 
\Bigl(\,\int\limits_{0}^{1} J(t)\,dt\,\Bigr)^2\,}}}$
\index{matrix!standardized!$\bm\hat{e}_{i}$, $\bm\hat{J}$}\\[2.5ex]
and\\[2ex]
\refstepcounter{DSBcount}
\label{EWD_4_4_04}
\text{\hspace*{-0.8ex}(\theDSBcount)}
\hspace*{4ex}
$\displaystyle{
\xi_{1,A} = \sum\limits_{i=1}^n\, \bm\hat{e}_{i}^{3} \int\limits_{0}^{1} \bm\hat{J}^{3}(t)\,dt}$,\\[2.5ex]
\refstepcounter{DSBcount}
\label{EWD_4_4_05}
\text{\hspace*{-0.8ex}(\theDSBcount)}
\hspace*{4ex}
$\displaystyle{
\xi_{2,A} = \biggl(\,\sum\limits_{i=1}^n\, \bm\hat{e}_{i}^{4}\,\Bigl\{\,\int\limits_{0}^{1} \bm\hat{J}^{4}(t)\,dt
\, -\, 3\,\Bigr\}\,\biggr)\, - \, \dfrac{3}{n}\,\Bigl\{\,\int\limits_{0}^{1} \bm\hat{J}^{4}(t)\,dt
\, -\, 1\,\Bigr\} }$.
\index{Edgeworth expansion}\index{Edgeworth expansion!$\xi_{1,A}$, $\xi_{2,A}$}\\[3ex]
It should be mentioned that the type of conditions placed on $J$ in this section can also be 
found in Albers, Bickel and van Zwet \cite{10.1214/aos/1176343350}, 
Bickel and van Zwet \cite{10.1214/aos/1176344305} and 
Does \cite{10.1214/aos/1176346166}, \cite{Does1982}, \cite{does1982higher}.\\[2.8ex]
We start with a few simple but useful preliminary considerations.\\[4ex] 
\refstepcounter{DSBcount}
\label{EWD_4_4_06}
\textbf{\hspace*{-0.8ex}\theDSBcount\ Lemma}\\[0.8ex]
Let $\beta > 1$ and $n \in \mathbb{N}$. Then\\[2.5ex]
\refstepcounter{DSBcount}
\label{EWD_4_4_07}
\text{\hspace*{-0.8ex}(\theDSBcount)}
\hspace*{4ex}
$\displaystyle{\sum\limits_{j=1}^n\,
\biggl(\,\dfrac{j}{n+1}\,\Bigl(\,1 - \dfrac{j}{n+1}\,\Bigr)
\,\biggr)^{-\beta}\, \leq\, C(\beta)\,n^{\beta}}$,\\[2ex]
where $\displaystyle{C(\beta) = 4^{\beta}\,\sum\limits_{j = 1}^{\infty}\,\dfrac{1}{j^{\beta}}\, <\, \infty}$.\\[4ex]
\textbf{Proof:}\\[2ex] 
\hspace*{7ex}\begin{tabular}[t]{@{}l@{\hspace*{1ex}}c@{\hspace*{1ex}}l@{}}
$\displaystyle{\sum\limits_{j=1}^n\,\biggl(\,\dfrac{j}{n+1}\,\Bigl(\,1 - \dfrac{j}{n+1}\,\Bigr)
\,\biggr)^{-\beta}}$
$=$&$\displaystyle{\sum\limits_{j=1}^n\,
\biggl(\,\dfrac{(n+1)^2}{j\,(n+1-j)}
\,\biggr)^{\beta}}$
\end{tabular}\\[4ex]
\hspace*{20.1ex}\begin{tabular}[t]{@{}l@{\hspace*{1ex}}l@{}}
$=$&$\displaystyle{(n+1)^{\beta}\,\sum\limits_{j=1}^n\,
\biggl(\,\dfrac{1}{j} + \dfrac{1}{n+1-j}
\,\biggr)^{\beta}}$.
\end{tabular}\\[2.5ex]
An application of H{\"o}lder's inequality\index{H{\"o}lder's inequality!for finite sequences using length $\nu$} 
(\ref{EWD_0_1_05}) with $\nu = 2$ and $p = \beta$ further yields\\[2.5ex]
\hspace*{20.1ex}\begin{tabular}[t]{@{}l@{\hspace*{1ex}}l@{}}
$\leq$&$\displaystyle{(n+1)^{\beta}\,\sum\limits_{j=1}^n\,2^{\beta-1}
\biggl(\,\Bigl(\,\dfrac{1}{j}\,\Bigr)^{\beta} + \Bigl(\,\dfrac{1}{n+1-j}\,\Bigr)^{\beta}
\,\biggr)}$\\[4ex]
$=$&$\displaystyle{(n+1)^{\beta}\,2^{\beta-1+1}\,\sum\limits_{j=1}^n\,\dfrac{1}{j^{\beta}}}$\\[4ex]
$\leq$&$\displaystyle{(2n)^{\beta}\,2^{\beta}\,\sum\limits_{j=1}^{\infty}\,\dfrac{1}{j^{\beta}}}$
\hspace*{32.2ex}(because of $n + 1 \leq 2n$)\\[4ex]
$=$&$C(\beta)\,n^{\beta}$.
\end{tabular}\\[-3ex]
\hspace*{1ex}\hfill$\Box$\\[7ex]
The following lemma contains an important definition for this section and some basic consequences.\\[4ex]
\refstepcounter{DSBcount}
\label{EWD_4_4_08}
\textbf{\hspace*{-0.8ex}\theDSBcount\ Lemma}\\[0.8ex]
Let $g: (0,\,1) \rightarrow \mathbb{R}$ be continuously differentiable and $0 < \alpha < 1$.\\[1.5ex]
Then $g$ satisfies the \textbf{condition} {\boldmath $V_{\alpha}$}\index{condition $V_{\alpha}$},
if a $\Gamma > 0$ exists such that\\[2ex]
\refstepcounter{DSBcount}
\label{EWD_4_4_09}
\text{\hspace*{-0.8ex}(\theDSBcount)}
\hspace*{4ex}
$|\,g'(t)\,| \leq \Gamma\,\bigl(\,t\,(1-t)\,\bigr)^{-1-\alpha}$
\hspace*{4ex}for all $t \in (0,\,1)$.\\[2.5ex]
For functions $g$ satisfying the condition $V_{\alpha}$ holds:
\begin{enumerate}
\item\label{EWD_4_4_08_BWa}
There exists a constant $C > 0$ depending only on $g$, $\Gamma$ and $\alpha$ such that\\[2ex]
\refstepcounter{DSBcount}
\label{EWD_4_4_10}
\text{\hspace*{-0.8ex}(\theDSBcount)}
\hspace*{2.8ex}
$|\,g(t)\,| \leq C\,\bigl(\,t\,(1-t)\,\bigr)^{-\alpha}$
\hspace*{4ex}for all $t \in (0,\,1)$.
\item\label{EWD_4_4_08_BWb}
$g^s$ is integrable for every $s > 0$ with $\alpha < \dfrac{1}{s}\ \Leftrightarrow\ s\alpha < 1$.
\end{enumerate}
\vspace*{3.5ex}
\textbf{Proof:}\\[0.8ex]
Since the function $f(t) = \Gamma\,\bigl(\,t\,(1-t)\,\bigr)^{-1-\alpha}$ is 
axially symmetric to the line $t = 1/2$, it is sufficient to consider the intervall
$(0,\,1/2\,]$ instead of the intervall $(0,\,1)$.
\begin{enumerate}
\item
For $x \in (0,\,1/2\,]$ we have\\[1ex]
\hspace*{4ex}\begin{tabular}[t]{@{}l@{\hspace*{1ex}}c@{\hspace*{1ex}}l@{}}
$g(x) - g(1/2)$&$=$&$\displaystyle{\int\limits_{1/2}^{x} g'(t)\,dt
= \int\limits_{x}^{1/2} - g'(t)\,dt \leq \int\limits_{x}^{1/2} |\,g'(t)\,|\,dt}$\\[4ex]
&$\leq$&$\displaystyle{\Gamma\,\int\limits_{x}^{1/2} \bigl(\,t\,(1-t)\,\bigr)^{-1-\alpha}\,dt
\leq \Gamma\,\Bigl(\,\dfrac{1}{2}\,\Bigr)^{-1-\alpha}\,\int\limits_{x}^{1/2} t^{-1-\alpha}\,dt}$\\[4.5ex]
&$=$&$\displaystyle{\Gamma\,\dfrac{2^{1+\alpha}}{\alpha}\,\Bigl(\,x^{-\alpha} - 2^{\alpha}\,\Bigr)
\leq \Gamma\,\dfrac{2^{1+\alpha}}{\alpha}\,x^{-\alpha}}
\leq \displaystyle{\Gamma\,\dfrac{2^{1+\alpha}}{\alpha}\,\bigl(\,x\,(1-x)\,\bigr)^{-\alpha}}$.
\end{tabular}\\[2ex]
It follows that\\[2ex]
\hspace*{4ex}$\displaystyle{g(x) \leq \biggl[\,\Gamma\,\dfrac{2^{1+\alpha}}{\alpha} + \big|\,g(1/2)\,\big|
\,4^{-\alpha}\,\biggr]\,\bigl(\,x\,(1-x)\,\bigr)^{-\alpha} = C\,\bigl(\,x\,(1-x)\,\bigr)^{-\alpha}}$.\\[2.5ex]
Since $-g(x) \leq C\,\bigl(\,x\,(1-x)\,\bigr)^{-\alpha}$ can be derived similarly, we get
$|\,g(x)\,| \leq C\,\bigl(\,x\,(1-x)\,\bigr)^{-\alpha}$.
\item
This follows from (\ref{EWD_4_4_10}) and since 
$\displaystyle{\int\limits_{0}^{1/2}t^{- \beta}\,dt < \infty}$ for $0 < \beta < 1$.
\hspace*{1ex}\hfill$\Box$
\end{enumerate}
\vspace*{2ex}
We now consider the case of the approximating scores (cf. (\ref{EWD_4_4_02})). The following elementary lemma 
is an important help for this case.\\[4ex]
\refstepcounter{DSBcount}
\label{EWD_4_4_11}
\textbf{\hspace*{-0.8ex}\theDSBcount\ Lemma}\\[0.8ex]
Let $g: (0,\,1) \rightarrow \mathbb{R}$ be continuously differentiable and $0 < \alpha < 1$ such that
the \textbf{condition} {\boldmath $V_{\alpha}$} (cf. (\ref{EWD_4_4_09})) is satisfied for $g$.\\[2.5ex]
Then there exist constants $C_{1} > 0$ and $C_{2} > 0$ depending only on $g$ and $\alpha$ such that\\[2ex]
\refstepcounter{DSBcount}
\label{EWD_4_4_12}
\text{\hspace*{-0.8ex}(\theDSBcount)}
\hspace*{2.8ex}
$\displaystyle{\bigg|\,\sum\limits_{j=1}^n\,g\Bigl(\,\dfrac{j}{n+1}\,\Bigr)\, -\, n 
\int\limits_{0}^{1} g(t)\,dt\,\bigg|\, \leq\, C_{1}\,n^{\alpha}}$
\hspace*{8ex}for all $n \in \mathbb{N}$,\\[3ex]
\refstepcounter{DSBcount}
\label{EWD_4_4_13}
\text{\hspace*{-0.8ex}(\theDSBcount)}
\hspace*{2.8ex}
$\displaystyle{\bigg|\,\sum\limits_{j=1}^n\,\Big|\,g\Bigl(\,\dfrac{j}{n+1}\,\Bigr)\,\Big|\, -\, n 
\int\limits_{0}^{1} |\,g(t)\,|\,dt\,\bigg|\, \leq\, C_{2}\,n^{\alpha}}$
\hspace*{3.6ex}for all $n \in \mathbb{N}$.\\[4ex]
\textbf{Proof:}\\[0.8ex] 
We get for all $n \geq 2$:\\[2ex]
\hspace*{12.1ex}\begin{tabular}[t]{@{}l@{\hspace*{1ex}}l@{}}
\hspace*{1.7ex}&$\displaystyle{\bigg|\,\dfrac{1}{n}\,\sum\limits_{j=1}^n\,g\Bigl(\,\dfrac{j}{n+1}\,\Bigr)\, -\, 
\int\limits_{0}^{1} g(t)\,dt\,\bigg|}$
\end{tabular}\\[4.5ex]
\hspace*{12.1ex}\begin{tabular}[t]{@{}l@{\hspace*{1ex}}l@{}}
$\leq$&$\displaystyle{\int\limits_{0}^{1/(n+1)} |\,g(t)\,|\,dt\, +\,  
\int\limits_{n/(n+1)}^{1} |\,g(t)\,|\,dt\, +\, \dfrac{1}{n}\,\Big|\,g\Bigl(\,\dfrac{n}{n+1}\,\Bigr)\,\Big|}$\\[4.5ex]
&$\displaystyle{+\ \Big|\,\dfrac{1}{n} - \dfrac{1}{n+1}\,\Big|\,
\sum\limits_{j=1}^{n-1}\,\Big|\,g\Bigl(\,\dfrac{j}{n+1}\,\Bigr)\,\Big|\, +\,
\bigg|\,\int\limits_{1/(n+1)}^{n/(n+1)} g(t)\,dt\, 
-\, \dfrac{1}{n+1}\,\sum\limits_{j=1}^{n-1}\,g\Bigl(\,\dfrac{j}{n+1}\,\Bigr)\,\bigg|}$\\[5ex]
$=$&$S_{1} + S_{2} + S_{3} + S_{4} + S_{5}$.
\end{tabular}\\[3ex]
For the estimates of $S_{1}$ to $S_{4}$ we use (\ref{EWD_4_4_10}) from Lemma \ref{EWD_4_4_08}, \ref{EWD_4_4_08_BWa})
and for the estimate of $S_{5}$ we use (\ref{EWD_4_4_09}).\\[2.8ex]
In addition to the constant $C$ from (\ref{EWD_4_4_10}), we obtain further constants $c_{1}$, $c_{2}$, $c_{3}$ 
and $c_{4}$ depending only on $g$, $\Gamma$ and $\alpha$ such that\\[2ex]
\hspace*{12.1ex}$\displaystyle{S_{1}\, \leq\, C\,2^{\alpha} \int\limits_{0}^{1/(n+1)} t^{-\alpha}\,dt
= C\,\dfrac{2^{\alpha}}{1 - \alpha}\,\Bigl(\,\dfrac{1}{n+1}\,\Bigr)^{1-\alpha}\, \leq\, c_{1}\,n^{\alpha-1}}$,\\[1.5ex]
\hspace*{12.1ex}$\displaystyle{S_{2}\, \leq\, c_{1}\,n^{\alpha-1}}$
\hspace*{3.3ex}(analogous due to the axial symmetry of $\bigl(\,t\,(1-t)\,\bigr)^{-\alpha}$ to $t = \dfrac{1}{2}$),\\[2ex]
\hspace*{12.1ex}$\displaystyle{S_{3}\, 
\leq\, C\,\dfrac{1}{n}\,\biggl(\,\dfrac{n}{n+1}\,\Bigl(\,1 - \dfrac{n}{n+1}\,\Bigr)\,\biggr)^{-\alpha}\, 
\leq\, C\, 2^{\alpha}\,\dfrac{1}{n}\,\Bigl(\,\dfrac{1}{n+1}\,\Bigr)^{-\alpha}
\, \leq\, c_{2}\,n^{\alpha-1}}$,\\[2.5ex]
\hspace*{12.1ex}\begin{tabular}[t]{@{}l@{\hspace*{1ex}}c@{\hspace*{1ex}}l@{}}
$S_{4}$&$\leq$&$\displaystyle{C\,\dfrac{1}{n^2}\,\sum\limits_{j = 1}^{n-1}\,
\biggl(\,\dfrac{j}{n+1}\,\Bigl(\,1 - \dfrac{j}{n+1}\,\Bigr)\,\biggr)^{-\alpha}}$\\[4.5ex]
&$\leq$&$\displaystyle{C\,\dfrac{1}{n^2}\,\sum\limits_{j = 1}^{n-1}\,
\biggl(\,\dfrac{j}{n+1}\,\Bigl(\,1 - \dfrac{j}{n+1}\,\Bigr)\,\biggr)^{-1-\alpha}}$
\hspace*{2.5ex}(since $\biggl(\,\dfrac{j}{n+1}\,\Bigl(\,1 - \dfrac{j}{n+1}\,\Bigr)\,\biggr) \leq 1$)\\[4.5ex]
&$\leq$&$c_{3}\,n^{\alpha-1}$\hspace*{31.3ex}(due to Lemma \ref{EWD_4_4_06} with $\beta = 1 + \alpha$),\\[2.5ex]
$S_{5}$&$\leq$&$\displaystyle{\sum\limits_{j = 1}^{n-1}\,
\int\limits_{j/(n+1)}^{(j+1)/(n+1)}\, \Big|\,g(t)\, -\, g\Bigl(\,\dfrac{j}{n+1}\,\Bigr)\,\Big|\,dt}$\\[5ex]
&$\leq$&$\displaystyle{\dfrac{1}{n^2}\,\sum\limits_{j = 1}^{n-1}\,
\sup\Bigl\{\,\big|\,g'(t)\,\big|\,:\,\dfrac{j}{n+1} < t < \dfrac{j+1}{n+1}\,\Bigr\}}$
\hspace*{11ex}(mean value theorem\index{mean value theorem})\\[4.5ex]
&$\leq$&$\displaystyle{\dfrac{\Gamma}{n^2}\,\sum\limits_{j = 1}^{n}\,
\biggl(\,\dfrac{j}{n+1}\,\Bigl(\,1 - \dfrac{j}{n+1}\,\Bigr)\,\biggr)^{-1-\alpha}}$
\hspace*{1.4ex}$\Bigl( \bigl(\,t\,(1-t)\,\bigr)^{-1-\alpha}$ monotone for $t \lessgtr \dfrac{1}{2} \Bigr)$
\end{tabular}\\
\hspace*{12.1ex}\begin{tabular}[t]{@{}l@{\hspace*{1ex}}c@{\hspace*{1ex}}l@{}}
\hspace*{2.5ex}&$\leq$&$c_{4}\,n^{\alpha-1}$\hspace*{31.3ex}(due to Lemma \ref{EWD_4_4_06} with $\beta = 1 + \alpha$).
\end{tabular}\\[2.5ex]
Now, since we can choose $\Gamma$ as a function of $g$ and $\alpha$ by defining\\[2ex]
\refstepcounter{DSBcount}
\label{EWD_4_4_14}
\text{\hspace*{-0.8ex}(\theDSBcount)}
\hspace*{2.8ex}
$\Gamma = \Gamma(g,\alpha) := \min\Bigl\{\,\widetilde{\Gamma}\,:\,|\,g'(t)\,| 
\leq \widetilde{\Gamma}\,\bigl(\,t\,(1-t)\,\bigr)^{-1-\alpha}\ \ \text{for all}\ t \in (0,\,1)\,\Bigr\}$,
\index{condition $V_{\alpha}$!$\Gamma(g,\alpha)$}\\[2.5ex]
the estimate (\ref{EWD_4_4_12}) follows.\\[2.8ex]
The proof of (\ref{EWD_4_4_13}) is completely analogous except that the reverse triangle inequality
is additionally used for the estimate of $S_{5}$.\hspace*{1ex}\hfill$\Box$\\[4ex]
The lemma just proved is now used to give estimates for the sums from (\ref{EWD_4_3_03}) for the 
approximating scores $d_{j}$ from (\ref{EWD_4_4_02}).\\[4ex]
\refstepcounter{DSBcount}
\label{EWD_4_4_15}
\textbf{\hspace*{-0.8ex}\theDSBcount\ Corollary}\\[0.8ex]
Let $J: (0,\,1) \rightarrow \mathbb{R}$ be continuously differentiable and $0 < \alpha < 1$ such that
the \textbf{condition} {\boldmath $V_{\alpha}$} (cf. (\ref{EWD_4_4_09})) is satisfied for $J$.
In addition, let \mbox{\rule[0ex]{0ex}{4ex}$w \in \mathbb{N}$ with $\alpha < \dfrac{1}{w}$.}\\[2.5ex]
Furthermore, we define\\[2ex]
\hspace*{12.1ex}$\displaystyle{\bm\bar{J} = 
\dfrac{1}{n}\,\sum\limits_{j=1}^{n}\,J\Bigl(\,\dfrac{j}{n+1}\,\Bigr)}$.
\index{matrix!$\bm\bar{J}$}\\[0.5ex]
If \mbox{\rule[-4.5ex]{0ex}{6ex}$\displaystyle{\int\limits_{0}^{1} J(t)\,dt = 0}$}, 
then there exist constants $C_{1} > 0$ and $C_{2} > 0$
depending only on $J$, $\alpha$ and $w$ such that\\[2.3ex]
\refstepcounter{DSBcount}
\label{EWD_4_4_16}
\text{\hspace*{-0.8ex}(\theDSBcount)}
\hspace*{2.8ex}
$\displaystyle{\bigg|\,\sum\limits_{j=1}^n\,
\Bigl(\,J\Bigl(\,\dfrac{j}{n+1}\,\Bigr) - \bm\bar{J}\,\Bigr)^w\, -\, n 
\int\limits_{0}^{1} J^w(t)\,dt\,\bigg|\, \leq\, C_{1}\,n^{w\alpha}}$
\hspace*{5.2ex}for all $n \in \mathbb{N}$,\\[3ex]
\refstepcounter{DSBcount}
\label{EWD_4_4_17}
\text{\hspace*{-0.8ex}(\theDSBcount)}
\hspace*{2.8ex}
$\displaystyle{\bigg|\,\sum\limits_{j=1}^n\,
\Big|\,J\Bigl(\,\dfrac{j}{n+1}\,\Bigr) - \bm\bar{J}\,\Big|^w\, -\, n 
\int\limits_{0}^{1} |\,J^w(t)\,|\,dt\,\bigg|\, \leq\, C_{2}\,n^{w\alpha}}$
\hspace*{4ex}for all $n \in \mathbb{N}$.\\[4ex]
\textbf{Proof:}\\[0.8ex] 
At first, we note that $\,0 < w\alpha < 1\,$ and that $J^w$ satisfies the condition $V_{w\alpha}$.
This follows for $w \geq 2$ because of
$\bigl( J^w \bigr)' = w\,J^{w-1}\,J'$
from (\ref{EWD_4_4_10}) and (\ref{EWD_4_4_09}).\\[2.8ex]
To prove (\ref{EWD_4_4_16}), we therefore apply Lemma \ref{EWD_4_4_11} to $g = J^w$ 
and obtain\\[2ex]
\refstepcounter{DSBcount}
\label{EWD_4_4_18}
\text{\hspace*{-0.8ex}(\theDSBcount)}
\hspace*{2.8ex}
$\displaystyle{\bigg|\,\sum\limits_{j=1}^n\,J^w\Bigl(\,\dfrac{j}{n+1}\,\Bigr)\, -\, n 
\int\limits_{0}^{1} J^w(t)\,dt\,\bigg|\, \leq\, c_{1}\,n^{w\alpha}}$.\\[2.5ex]
Like the following constants, $c_{1}$ only depends on $J$, $\alpha$ and $w$.\\[2.8ex]
In addition, we have\\[2ex]
\refstepcounter{DSBcount}
\label{EWD_4_4_19}
\text{\hspace*{-0.8ex}(\theDSBcount)}
\hspace*{2.8ex}
\begin{tabular}[t]{@{}l@{\hspace*{1ex}}l@{}}
&$\displaystyle{\bigg|\,\sum\limits_{j=1}^n\,
\Bigl(\,J\Bigl(\,\dfrac{j}{n+1}\,\Bigr) - \bm\bar{J}\,\Bigr)^w\, -\, 
\sum\limits_{j=1}^n\,J^w\Bigl(\,\dfrac{j}{n+1}\,\Bigr)\,\bigg|}$\\[4.5ex]
$\leq$&$\displaystyle{\sum\limits_{j=1}^n\,\bigg|\,
\Bigl(\,J\Bigl(\,\dfrac{j}{n+1}\,\Bigr) - \bm\bar{J}\,\Bigr)^w\, -\, 
J^w\Bigl(\,\dfrac{j}{n+1}\,\Bigr)\,\bigg|}$.
\end{tabular}\\[3ex]
An application of the mean value theorem\index{mean value theorem} 
to $f(\theta) = (x - \theta\,y)^w$, where 
\mbox{\rule[-2ex]{0ex}{4ex}$x = J\Bigl(\,\dfrac{j}{n+1}\,\Bigr)$} 
and $y = \bm\bar{J}$, further gives\\[2.5ex]
\hspace*{12.1ex}\begin{tabular}[t]{@{}l@{\hspace*{1ex}}l@{}}
$\leq$&$\displaystyle{w\,|\,\bm\bar{J}\,|\,\sum\limits_{j=1}^n\,
\sup\limits_{0 < \theta < 1}\,\bigg|\,J\Bigl(\,\dfrac{j}{n+1}\,\Bigr) - \theta \bm\bar{J}\,\bigg|^{w-1}}$\\[4.5ex]
$\leq$&$\displaystyle{w\,|\,\bm\bar{J}\,|\,\sum\limits_{j=1}^n\,
\biggl(\,\Big|\,J\Bigl(\,\dfrac{j}{n+1}\,\Bigr)\,\Big|\, +\, |\,\bm\bar{J}\,|\,\biggr)^{w-1}}$.
\end{tabular}\\[3ex]
Now, for $w > 2$ we use H{\"o}lder's inequality\index{H{\"o}lder's inequality!for finite sequences using length $\nu$} 
(\ref{EWD_0_1_05}) with $\nu = 2$ and $p = w-1$. For the remaining cases $w = 2$ and $w = 1$, however,
the next step is clear. We obtain\\[2.5ex] 
\hspace*{12.1ex}\begin{tabular}[t]{@{}l@{\hspace*{1ex}}l@{}}
$\leq$&$\displaystyle{w\,2^{w-2}\,\biggl(\,|\,\bm\bar{J}\,|\,\sum\limits_{j=1}^n\,
\Big|\,J\Bigl(\,\dfrac{j}{n+1}\,\Bigr)\,\Big|^{w-1}\, +\, n\,|\,\bm\bar{J}\,|^w\,\biggr)}$\\[3ex]
$\leq$&$c_{2}\,\bigl(\,n^{\alpha} + n^{w\alpha - w + 1}\,\bigr)$\\[2ex]
$\leq$&$2\,c_{2}\,n^{w\alpha}$.
\end{tabular}\\[3ex]
The penultimate inequality is valid, since we get the inequalities\\[2ex]
\hspace*{12.1ex}$|\,\bm\bar{J}\,| \leq c_{3}\,n^{\alpha-1}$
\hspace*{2ex}and\hspace*{2ex}
$\displaystyle{\sum\limits_{j=1}^n\,
\Big|\,J\Bigl(\,\dfrac{j}{n+1}\,\Bigr)\,\Big|^{w-1}\, \leq\, c_{4}\,n}$\\[0.5ex]
due to $\displaystyle{\int\limits_{0}^{1} J(t)\,dt = 0}$
and Lemma \ref{EWD_4_4_11} with $g = J$ and $g = J^{w-1}$. Thus (\ref{EWD_4_4_16}) is proved.\\[2.8ex]
The proof of (\ref{EWD_4_4_17}) is completely analogous if we additionally use the 
reverse triangle inequality in the above estimate (\ref{EWD_4_4_19}).\hspace*{1ex}\hfill$\Box$\\[4ex] 
\refstepcounter{DSBcount}
\label{EWD_4_4_20}
\textbf{\hspace*{-0.8ex}\theDSBcount\ Remark}\\[0.8ex]
Corollary \ref{EWD_4_4_15} improves analogous results of Does \cite{does1982higher} (cf. e.g.
Lemmas 3.2.1, 2.2.1 and the proof of Lemma 3.4.1), where the additional condition\\[2ex]
\hspace*{12.1ex}$\limsup\limits_{t\, \rightarrow\, 0, 1}\,t\,(1-t)\,\bigg|\,\dfrac{J''(t)}{J'(t)}\,\bigg| < 2$\\[2.5ex]
is required.\\[4ex]
After the above considerations, we can now state the results for the case of the approximating scores. 
They are summarized in the following theorem.\\[4ex]
\refstepcounter{DSBcount}
\label{EWD_4_4_21}
\textbf{\hspace*{-0.8ex}\theDSBcount\ Theorem}\index{Theorem!for simple linear rank statistics!with approximating scores}\\[0.8ex]
Let $J: (0,\,1) \rightarrow \mathbb{R}$ be a nonconstant and continuously differentiable function, and\\[2ex]
\hspace*{12.1ex}$\displaystyle{n_{1}(J) = \sup \bigg\{\,n \in \mathbb{N}\,:\,\sum\limits_{j=1}^n\,
\Bigl(\,J\Bigl(\,\dfrac{j}{n+1}\,\Bigr) - \bm\bar{J}\,\Bigr)^2\, =\, 0\,\biggr\}}$.\\[2.5ex]
In addition, let $A = (a_{ij})$ be an $n{\times}n-$matrix such that\\[2ex] 
\hspace*{12.1ex}$n > n_{1}(J)$, $a_{ij} = e_{i}\,d_{j}$ and
$d_{j} = J\Bigl(\,\dfrac{j}{n+1}\,\Bigr)$.\\[2.8ex]
Furthermore, suppose that the following two conditions hold:\\[2ex]
(\ref{EWD_4_3_02})
\hspace*{4ex}
\begin{tabular}{@{}l@{}}
$\displaystyle{\sum\limits_{i=1}^{n}\,\big|\,e_{i} - \bm\bar{e}\,\big|^{r}\, \geq\, e n}$,\hspace*{4ex}
$\displaystyle{\sum\limits_{i=1}^{n}\,\big|\,e_{i} - \bm\bar{e}\,\big|^{k}\, \leq\, E n}$\\[2.5ex]
for some $k > 2$, $0 < r < k$ and $e > 0$, $E > 0$.
\end{tabular}\\[3ex]
\refstepcounter{DSBcount}
\label{EWD_4_4_22}
\text{\hspace*{-0.8ex}(\theDSBcount)}
\hspace*{2.8ex}
\begin{tabular}{@{}l@{}}
$J$ satisfies the \textbf{condition} {\boldmath $V_{\alpha}$} (cf. (\ref{EWD_4_4_09})) 
for some $0 < \alpha < \dfrac{1}{s}$,\\[1.5ex]
where $s \in \mathbb{N}$ and $s \geq 3$.
\end{tabular}\vspace*{0.5ex}
\begin{enumerate}
\item\label{EWD_4_4_21_BWa}
Then there exists a constant $\mathcal{E}_{5} > 0$ depending only on $e$, $E$, $r$, $k$ and $J$, $\alpha$, $s$
such that\\[2ex]
\refstepcounter{DSBcount}
\label{EWD_4_4_23}
\text{\hspace*{-0.8ex}(\theDSBcount)}
\hspace*{2.8ex}
\begin{tabular}[t]{@{}l@{\hspace*{1ex}}c@{\hspace*{1ex}}l@{}}
$||\,\mathscr{F}_{\!A} - e_{1,A}\,||$&$\leq$&
$\displaystyle{\mathcal{E}_{5}\,(\log n)^2\,n^{-\,1\, +\, \bigl((4/k)\, -\, 1\bigr)^{+}\, 
+\, \bigl((4/s)\, -\, 1\bigr)^{+}}}$.
\end{tabular}\vspace*{0.5ex}
\item\label{EWD_4_4_21_BWb}
Let $\epsilon > 0$. 
Then there exists a constant $\mathcal{E}_{6} > 0$ depending only on
$e$, $E$, $r$, $k$, $J$, $\alpha$, $s$ and $\epsilon$ such that\\[2ex]
\refstepcounter{DSBcount}
\label{EWD_4_4_24}
\text{\hspace*{-0.8ex}(\theDSBcount)}
\hspace*{2.8ex}
\begin{tabular}[t]{@{}l@{\hspace*{1ex}}c@{\hspace*{1ex}}l@{}}
$||\,\mathscr{F}_{\!A} - e_{2,A}\,||$&$\leq$&
$\displaystyle{\mathcal{E}_{6}\,n^{-\,(3/2)\, +\, \epsilon\, +\, \bigl((5/k)\, -\, 1\bigr)^{+}\, 
+\, \bigl((5/s)\, -\, 1\bigr)^{+}}}$.
\end{tabular}\vspace*{0.5ex}
\item\label{EWD_4_4_21_BWc}
Let $k \geq 3$.
Then there exists a constant $\mathcal{E}_{7} > 0$ depending only on 
$e$, $E$, $r$, $k$ and $J$, $\alpha$, $s$ such that\\[2ex]
\refstepcounter{DSBcount}
\label{EWD_4_4_25}
\text{\hspace*{-0.8ex}(\theDSBcount)}
\hfill
\begin{tabular}[t]{@{}l@{\hspace*{1ex}}c@{\hspace*{1ex}}l@{}}
$||\,\mathscr{F}_{\!A} - e_{1,A}^{I}\,||$&$\leq$&
$\displaystyle{\mathcal{E}_{7}\,\max \bigg\{\,n^{-\,(3/2)\, +\, 3\,\alpha}\, ,\,
(\log n)^2\,n^{-\,1\, +\, \bigl((4/k)\, -\, 1\bigr)^{+}\, 
+\, \bigl((4/s)\, -\, 1\bigr)^{+}}\,\bigg\}}$.
\end{tabular}\vspace*{0.5ex}
\item\label{EWD_4_4_21_BWd}
Let $\epsilon > 0$, $k \geq 4$ and $s \geq 4$.
Then there exists a constant $\mathcal{E}_{8} > 0$ depending only on 
$e$, $E$, $r$, $k$, $J$, $\alpha$, $s$ and $\epsilon$ such that\\[2ex]
\refstepcounter{DSBcount}
\label{EWD_4_4_26}
\text{\hspace*{-0.8ex}(\theDSBcount)}
\hfill
\begin{tabular}[t]{@{}l@{\hspace*{1ex}}c@{\hspace*{1ex}}l@{}}
$||\,\mathscr{F}_{\!A} - e_{2,A}^{I}\,||$&$\leq$&
$\displaystyle{\mathcal{E}_{8}\,\max \bigg\{\,n^{-\,(3/2)\, +\, 3\,\alpha}\, ,\,
n^{-\,(3/2)\, +\, \epsilon\, +\, \bigl((5/k)\, -\, 1\bigr)^{+}\, 
+\, \bigl((5/s)\, -\, 1\bigr)^{+}}\,\bigg\}}$.
\end{tabular}
\end{enumerate}
\vspace*{3.5ex}
\textbf{Proof:}\\[0.8ex]
\textbf{For a) and b):}\\[0.8ex]
To show \ref{EWD_4_4_21_BWa}) and \ref{EWD_4_4_21_BWb}), we use the Theorem \ref{EWD_4_3_13}.
For this, we must first prove (\ref{EWD_4_3_03}). We apply Corollary \ref{EWD_4_4_15}
to the function \mbox{\rule[0ex]{0ex}{4.5ex}$\displaystyle{\widetilde{J} = J - \int\limits_{0}^{1} J(t)\,dt}$},
once with $w = 2$ and once with $w = s$. 
Note that from (\ref{EWD_4_4_22}) follows 
\mbox{\rule[0ex]{0ex}{3.5ex}$\alpha < \dfrac{1}{s} \leq \dfrac{1}{3} < \dfrac{1}{2}$}. 
We then get:\\[2.8ex]
There exist positive constants $d$, $D$ depending only on $J$, $\alpha$ and $s$ such that\\[2ex]
\hspace*{12.1ex}$\displaystyle{\sum\limits_{j=1}^{n}\,\big|\,d_{j} - \bm\bar{d}\,\big|^{2}\, \geq\, d n}$
\hspace*{2ex}for $n > n_{1}(J)$,\hspace*{7ex}
$\displaystyle{\sum\limits_{j=1}^{n}\,\big|\,d_{j} - \bm\bar{d}\,\big|^{s}\, \leq\, D n}$
\hspace*{2ex}for $n \in \mathbb{N}$.\\[2.5ex]
Next we prove (\ref{EWD_4_3_04}). Since $J$ is not constant and continuously differentiable,
there exists an intervall $\Psi = [\,t_{1},\,t_{2}\,]$ with $0 < t_{1} < t_{2} < 1$ 
such that\\[2ex]
\hspace*{12.1ex}$\displaystyle{0 < c_{1} = \min\limits_{t \in \Psi} |\,J'(t)\,| \leq 
c_{2} = \max\biggl\{\,1,\,\max\limits_{t \in \Psi} |\,J'(t)\,|\,\biggr\} < \infty}$.\\[2ex]
Furthermore, let $n_{0} \in \mathbb{N}$ such that $t_{2} - t_{1} \geq \dfrac{1}{n_{0} + 1}$ and $n \geq n_{0}$.\\[1.5ex]
Now, for $x \in \Psi$ there exists an $j_{0} \in \{\,1,\ldots,n\,\}$ such that $\dfrac{j_{0}}{n+1} \in \Psi$ and
$\Big|\,x - \dfrac{j_{0}}{n+1}\,\Big| < \dfrac{1}{n+1}$. 
Due to the mean value theorem\index{mean value theorem} we get\\[2.2ex]
\hspace*{12.1ex}$\Big|\,J(x) - J\Bigl(\,\dfrac{j_{0}}{n+1}\,\Bigr)\,\Big|
\leq \Big|\,x - \dfrac{j_{0}}{n+1}\,\Big| \cdot c_{2} < \dfrac{c_{2}}{n+1} < \dfrac{c_{2}}{n}$\\[2.2ex]
and then with the definitions $\zeta = \dfrac{c_{2}}{n} \,\Bigl(\,\geq \dfrac{1}{n} \geq n^{-3/2}\,\log\ n\,\Bigr)$ and 
$\delta = \dfrac{c_{1} \cdot (t_{2} - t_{1})}{c_{2}}$:\\[2ex]
\hspace*{12.1ex}
$\displaystyle{\lambda\Bigl(\,\bigcup_{j = 1}^{n}\,\Bigl\{\,x \in \mathbb{R}\,:\,\big|\,x - d_{j}\,\big| < \zeta\,
\Bigr\}\,\Bigr) \geq \lambda\Bigl(\,J(\,\Psi\,)\,\Bigr)
\geq c_{1}\,(t_{2} - t_{1}) = \delta n \zeta}$\hspace*{3ex}for $n \geq n_{0}$.\\[4ex]
\textbf{For c) and d):}\\[0.8ex]
Since the proofs of \ref{EWD_4_4_21_BWc}) and \ref{EWD_4_4_21_BWd}) are completely analogous, 
we only show \ref{EWD_4_4_21_BWd}).\\[2.8ex]
Because of \ref{EWD_4_4_21_BWb}) it is sufficient to prove the following estimate:\\[2ex]
\refstepcounter{DSBcount}
\label{EWD_4_4_27}
\text{\hspace*{-0.8ex}(\theDSBcount)}
\hspace*{2.8ex}
$\displaystyle{||\,e_{2,A} - e_{2,A}^{I}\,||\, \leq\, c_{3}\,n^{-\,(3/2)\, +\, 3\,\alpha}}$
\hspace*{4ex}for $n > n_{1}(J)$.\\[2.5ex]
Like the following constants, $c_{3}$ only depends on $e$, $E$, $r$, $k$ and $J$, $\alpha$.\\[2.8ex]
According to Corollary \ref{EWD_4_4_15}, applied to $\bm\hat{J}$ (cf. (\ref{EWD_4_4_03})) 
and $w = 2$, we now obtain for\\[2ex]
\hspace*{12.1ex}$\displaystyle{b_{n} = \dfrac{1}{n-1}\,\sum\limits_{j=1}^{n}\,
\Bigl(\,\bm\hat{J}\Bigl(\,\dfrac{j}{n+1}\,\Bigr) - \bigl(\overline{\bm\hat{J}}\bigr)\,\Bigr)^2}$\\[2ex]
the inequality\\[2.2ex]
\hspace*{12.1ex}$\displaystyle{\big|\,b_{n} - 1\,\big| \leq c_{4}\,n^{2\alpha\, -\, 1}}$.\\[2ex]
From this we get\\[2ex]
\hspace*{12.1ex}\begin{tabular}[t]{@{}l@{\hspace*{1ex}}l@{}}
&$\displaystyle{\Biggg|\,\,\dfrac{\dfrac{1}{n}\,\sum\limits_{j=1}^n\,
\Bigl(\,J\Bigl(\,\dfrac{j}{n+1}\,\Bigr) - \bm\bar{J}\,\Bigr)^3}
{\Bigl\{\,\dfrac{1}{n-1}\,\sum\limits_{j=1}^n\,
\Bigl(\,J\Bigl(\,\dfrac{j}{n+1}\,\Bigr) - \bm\bar{J}\,\Bigr)^2\,\Bigr\}^{3/2}}\, -\,
\int\limits_{0}^{1} \bm\hat{J}^3(t)\,dt\,\,\Biggg|}$\\[7ex]
$=$&$\displaystyle{\biggg|\,\,\dfrac{\dfrac{1}{n}\,\sum\limits_{j=1}^n\,
\Bigl(\,\bm\hat{J}\Bigl(\,\dfrac{j}{n+1}\,\Bigr) - \bigl(\overline{\bm\hat{J}}\bigr)\,\Bigr)^3}
{b_{n}^{3/2}}\, -\,
\int\limits_{0}^{1} \bm\hat{J}^3(t)\,dt\,\,\biggg|}$\\[4ex]
$\leq$&$\displaystyle{\dfrac{1}{b_{n}^{3/2}}\,\bigg|\,\dfrac{1}{n}\,\sum\limits_{j=1}^n\,
\Bigl(\,\bm\hat{J}\Bigl(\,\dfrac{j}{n+1}\,\Bigr) - \bigl(\overline{\bm\hat{J}}\bigr)\,\Bigr)^3\, -\,
\int\limits_{0}^{1} \bm\hat{J}^3(t)\,dt\,\bigg|\, +\,
\biggl(\,\int\limits_{0}^{1} |\bm\hat{J}|^3(t)\,dt\,\biggr)\,
\bigg|\,\dfrac{1}{b_{n}^{3/2}} - 1\,\bigg|}$\\[4.5ex]
$\leq$&$c_{5}\,n^{3\alpha\, -\, 1} + c_{6}\,n^{2\alpha\, -\, 1}$
\hspace*{24.5ex}(Corollary \ref{EWD_4_4_15} for $\bm\hat{J}$ and $w = 3$)
\end{tabular}\\[2.5ex]
\hspace*{12.1ex}\begin{tabular}[t]{@{}l@{\hspace*{1ex}}l@{}}
$\leq$&$c_{7}\,n^{3\alpha\, -\, 1}$.
\end{tabular}\\[3ex]
Analogously we obtain using Corollary \ref{EWD_4_4_15} for $\bm\hat{J}$ and $w = 4 \leq s$\\[2.5ex]
\hspace*{12.1ex}\begin{tabular}[t]{@{}l@{\hspace*{1ex}}l@{}}
&$\displaystyle{\Biggg|\,\,\dfrac{\dfrac{1}{n}\,\sum\limits_{j=1}^n\,
\Bigl(\,J\Bigl(\,\dfrac{j}{n+1}\,\Bigr) - \bm\bar{J}\,\Bigr)^4}
{\Bigl\{\,\dfrac{1}{n-1}\,\sum\limits_{j=1}^n\,
\Bigl(\,J\Bigl(\,\dfrac{j}{n+1}\,\Bigr) - \bm\bar{J}\,\Bigr)^2\,\Bigr\}^{2}}\, -\,
\int\limits_{0}^{1} \bm\hat{J}^4(t)\,dt\,\,\Biggg|\, \leq\, c_{8}\,n^{4\alpha\, -\, 1}}$.
\end{tabular}\\[2.5ex]
Because of\\[2.5ex]
\hspace*{12.1ex}$\displaystyle{\sum\limits_{i=1}^n\,|\bm\hat{e}_{i}|^3\, \leq\, \dfrac{c_{9}}{n^{1/2}}}$,
\hspace*{3ex}
$\displaystyle{\sum\limits_{i=1}^n\,\bm\hat{e}_{i}^{\,4}\, \leq\, \dfrac{c_{10}}{n}}$
\hspace*{3ex}and\hspace*{3ex}
$\displaystyle{n^{4\alpha\, -\, 2}\, \leq\, n^{3\alpha\, -\, (3/2)}}$,\\[3ex]
we then get (\ref{EWD_4_4_27}).\hspace*{1ex}\hfill$\Box$\\[4ex]
\refstepcounter{DSBcount}
\label{EWD_4_4_28}
\textbf{\hspace*{-0.8ex}\theDSBcount\ Remarks}
\begin{enumerate}
\item\label{EWD_4_4_28_BWa}
The Remarks \ref{EWD_4_3_21}, \ref{EWD_4_3_21_BWa}) and \ref{EWD_4_3_21_BWb}) apply analogously 
to Theorem \ref{EWD_4_4_21}, i.e. in the case of\\[2ex]
(\ref{EWD_4_3_22})
\hspace*{2.4ex}
$\displaystyle{\Bigl[\,\bigl((5/k)\, -\, 1\bigr)^{+}\, - \,\bigl((4/k)\, -\, 1\bigr)^{+}\,\Bigr]\, +\,
\Bigl[\,\bigl((5/s)\, -\, 1\bigr)^{+}\, - \,\bigl((4/s)\, -\, 1\bigr)^{+}\,\Bigr]\, <\, \dfrac{1}{2}}$,\\[2.5ex]
we can omit the $(\log n)^2$ in (\ref{EWD_4_4_23}) and (\ref{EWD_4_4_25}).\\[2.8ex] 
Furthermore, the $n^{\epsilon}$ in (\ref{EWD_4_4_24}) and (\ref{EWD_4_4_26}) can be replaced by $(\log n)^3$, if\\[2ex] 
\hspace*{12.1ex}$E_{\!A}^3 \geq C\,n^{-\,3/2\, +\, \eta}$, where $C, \eta > 0$.\\[2.5ex]
The constants $\mathcal{E}_{6}$ and $\mathcal{E}_{8}$ then also depend on $C$ and $\eta$.
\item\label{EWD_4_4_28_BWb}
If $J$ is not constant and if $J$ fulfils the condition (\ref{EWD_4_4_22}), then follows 
$\displaystyle{\int\limits_{0}^{1} \bigl(\widetilde{J}\bigr)^2(t)\,dt > 0}$ 
(see $\widetilde{J}$ in the proof above)
and therefore $n_{1}(J) < \infty$ according to Corollary \ref{EWD_4_4_15}.
\item\label{EWD_4_4_28_BWc}
Does \cite{10.1214/aos/1176346166} (Theorem 3.1.1 and Theorem 2.1) and \cite{does1982higher} shows
(\ref{EWD_4_4_26}) with the convergence rate $\smallO(n^{-1})$ and under stronger conditions for the
regression constants $e_{i}$ and for the function $J$. 
In particular, he requires for $J$ that\\[2.5ex] 
\hspace*{12.1ex}$\limsup\limits_{t\, \rightarrow\, 0, 1}\,t\,(1-t)\,\bigg|\,\dfrac{J''(t)}{J'(t)}\,\bigg| < 2$
\hspace*{3ex}and\hspace*{3ex}
$|\,J'''(t)\,| \leq \Gamma\,\bigl(\,t\,(1-t)\,\bigr)^{-3-\alpha}$,\\[2.5ex]
where $0 < \alpha < \dfrac{1}{14}$. But the last inequality implies
(\ref{EWD_4_4_22}) with $s = 14$.
\item\label{EWD_4_4_28_BWd}
For the two most popular score functions $J(t) = t$ (Wilcoxon-Scores\index{scores!approximating!Wilcoxon}) and 
$J(t) = \Phi^{-1}(t)$ (van der Waerden-Scores\index{scores!approximating!van der Waerden}) 
there is a $\Gamma > 0$ with\\[2ex]
\hspace*{12.1ex}$|\,J'(t)\,| \leq \Gamma\,\bigl(\,t\,(1-t)\,\bigr)^{-1}$\hspace*{4ex}for all $t \in (0,\,1)$,\\[2.5ex]
so that (\ref{EWD_4_4_22}) is fulfilled for any $\alpha > 0$ and $s \in \mathbb{N}$ with $\alpha < \dfrac{1}{s}$.
\end{enumerate}
\vspace*{3.5ex}
Let us now look at the exact scores. We need the following elementary statements 
about order statistics of uniformly distributed random variables on $(0,\,1)$.\\[4ex]  
\refstepcounter{DSBcount}
\label{EWD_4_4_29}
\textbf{\hspace*{-0.8ex}\theDSBcount\ Lemma}\\[0.8ex]
Let $U_{j:n}$\index{order statistic!$U_{j:n}$} denote the $j$th
order statistic\index{order statistic} in a random sample of size $n$
from the uniform distribution on $(0,\,1)$ and let 
$g: (0,\,1) \rightarrow \mathbb{R}$ be an integrable function. Then 
\begin{enumerate}
\item\label{EWD_4_4_29_BWa}
$E\bigl(\,U_{j:n}\,\bigr) = \dfrac{j}{n+1}$\hspace*{4ex}for all $1 \leq j \leq n$,
\item\label{EWD_4_4_29_BWb}
$\displaystyle{\dfrac{1}{n}\,\sum\limits_{j=1}^n\,E\bigl(\,g(\,U_{j:n}\,)\,\bigr) =
\int\limits_{0}^{1} g(t)\,dt}$.
\end{enumerate}
\vspace*{3.5ex}
\textbf{Proof:}\\[0.8ex]
According to \cite{sidak1999theory}, Chapter 3, Section 3.1, Theorem 4, page 38 and 39, $U_{j:n}$ has the density\\[2ex]
\hspace*{12.1ex}$\displaystyle{f_{j:n}(t) = n\,\dbinom{n-1}{j-1}\,t^{j-1}\,(1-t)^{n-j}\,1_{(0,\,1)}(t)}$ 
\hspace*{4ex}for $t \in \mathbb{R}$\\[1.5ex]
and the expected value $E\bigl(\,U_{j:n}\,\bigr) = \dfrac{j}{n+1}$ for all $1 \leq j \leq n$. 
Therefore only \ref{EWD_4_4_29_BWb}) remains to be shown:\\[2ex]
\hspace*{12.1ex}\begin{tabular}[t]{@{}l@{\hspace*{1ex}}c@{\hspace*{1ex}}l@{\hspace*{7.4ex}}r@{}}
$\displaystyle{\dfrac{1}{n}\,\sum\limits_{j=1}^n\,E\bigl(\,g(\,U_{j:n}\,)\,\bigr)}$&$=$&
$\displaystyle{\dfrac{1}{n}\,\sum\limits_{j=1}^n\,n\,\int\limits_{0}^{1}\,g(t)\,\dbinom{n-1}{j-1}\,
t^{j-1}\,(1-t)^{n-j}\,dt}$\\[3.5ex]
&$=$&$\displaystyle{\int\limits_{0}^{1}\,g(t)\,\sum\limits_{i=0}^{n-1}\,\dbinom{n-1}{i}\,
t^{i}\,(1-t)^{(n-1)-i}\,dt}$&($i = j-1$)\\[3.5ex]
&$=$&$\displaystyle{\int\limits_{0}^{1}\,g(t)\,\Bigl(\,t + (1-t)\,\Bigr)^{n-1}\,dt = 
\int\limits_{0}^{1} g(t)\,dt}$.&$\Box$
\end{tabular}\\[4ex]
The following lemma, which is an analogue of (\ref{EWD_4_4_18}) for the exact scores, 
is needed for the change from $e_{1,A}$ and $e_{2,A}$
to $e_{1,A}^{I}$ and $e_{2,A}^{I}$ in the case of the exact scores.

\pagebreak

\refstepcounter{DSBcount}
\label{EWD_4_4_30}
\textbf{\hspace*{-0.8ex}\theDSBcount\ Lemma}\\[0.8ex]
Let $J: (0,\,1) \rightarrow \mathbb{R}$ be continuously differentiable and $0 < \alpha < 1$ such that
the \textbf{condition} {\boldmath $V_{\alpha}$} (cf. (\ref{EWD_4_4_09})) is satisfied for $J$.
In addition, let \mbox{\rule[0ex]{0ex}{4ex}$w \in \mathbb{N}$ with $\alpha < \dfrac{1}{w}$.}\\[2.5ex]
Then there exists a constant $C > 0$ depending only on $J$, $\alpha$ and $w$ (and all $l \leq w$) such that\\[2ex]
\refstepcounter{DSBcount}
\label{EWD_4_4_31}
\text{\hspace*{-0.8ex}(\theDSBcount)}
\hspace*{2.5ex}
$\displaystyle{\bigg|\,\sum\limits_{j=1}^n\,\Bigl(\,E\bigl(\,J(\,U_{j:n}\,)\,\bigr)\,\Bigr)^w\, -\,
n \int\limits_{0}^{1} J^w(t)\,dt\,\bigg|\, \leq\, C\,n^{w\alpha}}$
\hspace*{5.2ex}for all $n \in \mathbb{N}$.\\[4ex]
\textbf{Proof:}\\[0.8ex]
Due to Lemma \ref{EWD_4_4_29}, \ref{EWD_4_4_29_BWb}) with $g = J$, we can assume $w \geq 2$
without loss of generality. For all $1 \leq j \leq n$ we then obtain using the binomial formula\\[2ex]
\hspace*{11.8ex}\begin{tabular}[t]{@{}l@{\hspace*{1ex}}l@{}}
&$\displaystyle{E\bigl(\,J^w(\,U_{j:n}\,)\,\bigr)\, -\, \Bigl(\,E\bigl(\,J(\,U_{j:n}\,)\,\bigr)\,\Bigr)^w}$\\[2.5ex] 
$=$&$\displaystyle{\sum\limits_{l = 0}^{w-2}\,\dbinom{w}{l}\,\Bigl(\,E\bigl(\,J(\,U_{j:n}\,)\,\bigr)\,\Bigr)^l\,
E\Bigl(\,\Bigl[\,J(\,U_{j:n}\,) - E\bigl(\,J(\,U_{j:n}\,)\,\bigr)\,\Bigr]^{w-l}\,\Bigr)}$.
\end{tabular}\\[3ex]
Furthermore, we get for all $1 \leq j \leq n$ and $0 \leq l \leq w - 2$\\[2.5ex]
\hspace*{11.8ex}\begin{tabular}[t]{@{}l@{\hspace*{1ex}}l@{}}
&$\displaystyle{E\Bigl(\,\Big|\,J(\,U_{j:n}\,) - E\bigl(\,J(\,U_{j:n}\,)\,\bigr)\,\Big|^{w-l}\,\Bigr)}$\\[3.5ex]
$\leq$&$\displaystyle{2^{w-l-1}\,\biggl\{\,
E\Bigl(\,\Big|\,J(\,U_{j:n}\,) - J\Bigl(\,\dfrac{j}{n+1}\,\Bigr)\,\Big|^{w-l}\,\Bigr)\, + \,
\Big|\,J\Bigl(\,\dfrac{j}{n+1}\,\Bigr) - E\bigl(\,J(\,U_{j:n}\,)\,\bigr)\,\Big|^{w-l}\,\biggr\}}$\\[2ex]
&\hspace*{46.6ex}((\ref{EWD_0_1_05}) with $\nu = 2$ and 
$p = w-l$)\index{H{\"o}lder's inequality!for finite sequences using length $\nu$}\\[2ex]
$=$&$\displaystyle{2^{w-l-1}\,\biggl\{\,
E\Bigl(\,\Big|\,J(\,U_{j:n}\,) - J\Bigl(\,\dfrac{j}{n+1}\,\Bigr)\,\Big|^{w-l}\,\Bigr)\, + \,
\Big|\,E\Bigl(\,J\Bigl(\,\dfrac{j}{n+1}\,\Bigr) - J(\,U_{j:n}\,)\,\Bigr)\,\Big|^{w-l}\,\biggr\}}$\\[3.5ex]
$\leq$&$\displaystyle{2^{w-l}\,
E\Bigl(\,\Big|\,J(\,U_{j:n}\,) - J\Bigl(\,\dfrac{j}{n+1}\,\Bigr)\,\Big|^{w-l}\,\Bigr)}$
\end{tabular}\\[2ex]
and\\[1.5ex]
\hspace*{11.8ex}\begin{tabular}[t]{@{}l@{\hspace*{1ex}}l@{\hspace*{-6.2ex}}r@{}}
&$\displaystyle{\Big|\,E\bigl(\,J(\,U_{j:n}\,)\,\bigr)\,\Big|^l}$\\[2.5ex]
$\leq$&$\displaystyle{\,E\Bigl(\,\Big|\,J(\,U_{j:n}\,)\,\Big|^l\,\Bigr)}$
&(H{\"o}lder's inequality\index{H{\"o}lder's inequality!for random variables})\\[2.5ex]
$\leq$&$\displaystyle{2^{l-1}\,
\biggl\{\,E\Bigl(\,\Big|\,J(\,U_{j:n}\,) - J\Bigl(\,\dfrac{j}{n+1}\,\Bigr)\,\Big|^{l}\,\Bigr)\, + \,
\Big|\,J\Bigl(\,\dfrac{j}{n+1}\,\Bigr)\,\Big|^{l}\,\biggr\}}$\\[2ex]
&&((\ref{EWD_0_1_05}) with $\nu = 2$ and 
$p = l$)\index{H{\"o}lder's inequality!for finite sequences using length $\nu$}.
\end{tabular}\\[2.5ex]
If we use Lemma \ref{EWD_4_4_29}, \ref{EWD_4_4_29_BWb}) with $g = J^w$, we obtain
from the above considerations\\[2ex]
\refstepcounter{DSBcount}
\label{EWD_4_4_32}
\text{\hspace*{-0.8ex}(\theDSBcount)}
\hspace*{-0.8ex}
\begin{tabular}[t]{@{}l@{\hspace*{0.8ex}}l@{}}
&$\displaystyle{\bigg|\,\sum\limits_{j=1}^n\,\Bigl(\,E\bigl(\,J(\,U_{j:n}\,)\,\bigr)\,\Bigr)^w\, -\,
n \int\limits_{0}^{1} J^w(t)\,dt\,\bigg|}$\\[4.5ex]
$\leq$&$\displaystyle{\sum\limits_{j=1}^n\,\Big|\,E\bigl(\,J^w(\,U_{j:n}\,)\,\bigr)\, 
-\, \Bigl(\,E\bigl(\,J(\,U_{j:n}\,)\,\bigr)\,\Bigr)^w\,\Big|}$\\[4.5ex]
$\leq$&$\displaystyle{\sum\limits_{l = 0}^{w-2}\,\dbinom{w}{l}\,\biggl\{\,
\sum\limits_{j=1}^n\,\Big|\,E\bigl(\,J(\,U_{j:n}\,)\,\bigr)\,\Big|^l\,
E\Bigl(\,\Big|\,J(\,U_{j:n}\,) - E\bigl(\,J(\,U_{j:n}\,)\,\bigr)\,\Big|^{w-l}\,\Bigr)\,
\biggr\}}$\\[4.5ex]
$\leq$&$\displaystyle{2^{w-1}\,\sum\limits_{l = 0}^{w-2}\,\dbinom{w}{l}\,\biggl\{\,
\sum\limits_{j=1}^n\,E\Bigl(\,\Big|\,J(\,U_{j:n}\,) - J\Bigl(\,\dfrac{j}{n+1}\,\Bigr)\,\Big|^{l}\,\Bigr)\,
E\Bigl(\,\Big|\,J(\,U_{j:n}\,) - J\Bigl(\,\dfrac{j}{n+1}\,\Bigr)\,\Big|^{w-l}\,\Bigr)}$\\[4.5ex]
&\hspace*{20ex}$\displaystyle{+\ \sum\limits_{j=1}^n\,\Big|\,J\Bigl(\,\dfrac{j}{n+1}\,\Bigr)\,\Big|^{l}\,
E\Bigl(\,\Big|\,J(\,U_{j:n}\,) - J\Bigl(\,\dfrac{j}{n+1}\,\Bigr)\,\Big|^{w-l}\,\Bigr)\,\biggr\}}$.
\end{tabular}\\[3ex]
Now let $h$ be the function that is given by\\[2ex] 
\hspace*{11.8ex}$h'(t) = \Gamma(J,\alpha)\,\bigl(\,t\,(1-t)\,\bigr)^{-1-\alpha}$
\hspace*{2ex}for all $t \in (0,\,1)$
\hspace*{3ex}and\hspace*{3ex}
$h\bigl(\,\dfrac{1}{2}\,\bigr) = 0$\\[2ex]
(for the definition of $\Gamma(J,\alpha)$\index{condition $V_{\alpha}$!$\Gamma(g,\alpha)$} 
see (\ref{EWD_4_4_14})). Then, for all $1 \leq j \leq n$ we have\\[2ex]
\hspace*{11.8ex}\begin{tabular}[t]{@{}l@{\hspace*{1ex}}l@{\hspace*{14.7ex}}r@{}}
&$\Big|\,J(\,U_{j:n}\,) - J\Bigl(\,\dfrac{j}{n+1}\,\Bigr)\,\Big|$\\[3ex]
$=$&$\dfrac{|\,J'(\zeta)\,|}{|\,h'(\zeta)\,|}\, 
\Big|\,h(\,U_{j:n}\,) - h\Bigl(\,\dfrac{j}{n+1}\,\Bigr)\,\Big|$
&(generalized mean value theorem\index{mean value theorem}\index{mean value theorem!generalized})\\[3ex]
$\leq$&$\Big|\,h(\,U_{j:n}\,) - h\Bigl(\,\dfrac{j}{n+1}\,\Bigr)\,\Big|$
&(cf. (\ref{EWD_4_4_14}): $|\,J'(\zeta)\,| \leq |\,h'(\zeta)\,|$)
\end{tabular}\\[2.5ex]
for a $\zeta$ such that
$\min\Bigl\{\,U_{j:n},\,\dfrac{j}{n+1}\,\Bigr\} \leq \zeta \leq \max\Bigl\{\,U_{j:n},\,\dfrac{j}{n+1}\,\Bigr\}$.\\[2ex]
$h$ also fulfils the \textit{Condition} $R_{w}$ from Albers, Bickel \& van Zwet 
\cite{10.1214/aos/1176343350} (see page 150 above) due to \mbox{\rule[-2.5ex]{0ex}{6ex}$\alpha < \dfrac{1}{w}$}, 
so that Lemma A2.3 (with $m = 1$, $k_{1} = l$, $r_{1} = w$ and $N = n$) from this paper
yields for $1 \leq l \leq w$ and $1 \leq j \leq n$:\\[2.5ex]
\refstepcounter{DSBcount}
\label{EWD_4_4_33}
\text{\hspace*{-0.8ex}(\theDSBcount)}
\hspace*{2.3ex}
\begin{tabular}[t]{@{}l@{\hspace*{1ex}}l@{}}
&$\displaystyle{E\Bigl(\,\Big|\,J(\,U_{j:n}\,) - J\Bigl(\,\dfrac{j}{n+1}\,\Bigr)\,\Big|^{l}\,\Bigr)}$\\[4ex]
$\leq$&$\displaystyle{E\Bigl(\,\Big|\,h(\,U_{j:n}\,) - h\Bigl(\,\dfrac{j}{n+1}\,\Bigr)\,\Big|^{l}\,\Bigr)}$
\end{tabular}\\[4ex]
\hspace*{12.1ex}\begin{tabular}[t]{@{\hspace*{-0.3ex}}l@{\hspace*{0.8ex}}l@{}}
$\leq$&$\displaystyle{c_{1}\,\Biggl\{\,\dfrac{1}{n^{l/2}}\,\biggl(\,\dfrac{j}{n+1}\,\Bigl(\,1 - \dfrac{j}{n+1}\,\Bigr)
\,\biggr)^{-\,(l/2)\,-\,l\alpha}}$\\[3ex]
&\hspace*{5ex}$\displaystyle{+\ 
\dfrac{1}{n^{(l+1)/2}}\,\biggl(\,\dfrac{j}{n+1}\,\Bigl(\,1 - \dfrac{j}{n+1}\,\Bigr)
\,\biggr)^{-\,((l+1)/2)\,-\,l\alpha}\,\Biggr\}}$.
\end{tabular}\\[3ex]
Here, $c_{1}$ only depends on $\Gamma(J,\alpha)$\index{condition $V_{\alpha}$!$\Gamma(g,\alpha)$}, 
$\alpha$, $l$ and therefore $J$, $\alpha$, $l$, but \textbf{not} on $j$.\\[2.8ex] 
Because of (\ref{EWD_4_4_10}) we also get\\[2ex]
\refstepcounter{DSBcount}
\label{EWD_4_4_34}
\text{\hspace*{-0.8ex}(\theDSBcount)}
\hspace*{2.8ex}
$\Big|\,J\Bigl(\,\dfrac{j}{n+1}\,\Bigr)\,\Big|^{l}\,
\leq c_{2}\,\biggl(\,\dfrac{j}{n+1}\,\Bigl(\,1 - \dfrac{j}{n+1}\,\Bigr)\,\biggr)^{-l\alpha}$
\hspace*{2ex}for all $t \in (0,\,1)$,\\[2.5ex]
where $c_{2}$ (argumentation similar to (\ref{EWD_4_4_14})) also depends only on $J$, $\alpha$ and $l$.\\[2.8ex]
If we now apply the estimates (\ref{EWD_4_4_33}) and (\ref{EWD_4_4_34}) to the 
terms in (\ref{EWD_4_4_32}) and then use the inequality (\ref{EWD_4_4_07}) with the following five parameters\\[2ex]
\hspace*{12.1ex}\begin{tabular}[t]{@{}l@{}}
$\beta = \dfrac{w}{2} + w\alpha > 1$,\hspace*{1.4ex} 
$\beta = \dfrac{w+1}{2} + w\alpha > 1$,\hspace*{1.4ex} 
$\beta = \dfrac{w+2}{2} + w\alpha > 1$,\\[2.5ex] 
$\beta = \dfrac{w-l}{2} + w\alpha > 1$
\hspace*{1.4ex}and\hspace*{1.4ex} 
$\beta = \dfrac{w-l+1}{2} + w\alpha > 1$
\end{tabular}\\[2.5ex] 
(note $w \geq 2$ and $w - l \geq 2\,\,$!), then follows\\[2ex]
\hspace*{12.1ex}$\displaystyle{\bigg|\,\sum\limits_{j=1}^n\,\Bigl(\,E\bigl(\,J(\,U_{j:n}\,)\,\bigr)\,\Bigr)^w\, -\,
n \int\limits_{0}^{1} J^w(t)\,dt\,\bigg|
\, \leq\, 2^{w-1}\,\sum\limits_{l = 0}^{w-2}\,\dbinom{w}{l}\,c_{3}\,n^{w\alpha}}$,\\[2.5ex]
where $c_{3}$ depends only on $J$, $\alpha$ and $w$ (and all $l \leq w$).\hspace*{1ex}\hfill$\Box$\\[4ex]
The results for the case of the exact scores are now summarized in the following theorem.\\[4ex]
\refstepcounter{DSBcount}
\label{EWD_4_4_35}
\textbf{\hspace*{-0.8ex}\theDSBcount\ Theorem}\index{Theorem!for simple linear rank statistics!with exact scores}\\[0.8ex]
Let $J: (0,\,1) \rightarrow \mathbb{R}$ be a nonconstant, continuously differentiable and 
integrable function and $d_{j} = E\bigl(\,J(\,U_{j:n}\,)\,\bigr)$ for $j = 1,\ldots,n$ and\\[2ex]
\hspace*{12.1ex}$\displaystyle{n_{2}(J) = \sup \bigg\{\,n \in \mathbb{N}\,:\,\sum\limits_{j=1}^n\,
\bigl(\,d_{j} - \bm\bar{d}\,\bigr)^{2}\, =\, 0\,\biggr\}}$.\\[2.5ex]
In addition, let $A = (a_{ij})$ be an $n{\times}n-$matrix such that\\[2ex] 
\hspace*{12.1ex}$n > n_{2}(J)$
\hspace*{3ex}and\hspace*{3ex}
$a_{ij} = e_{i}d_{j}$.\\[2.5ex]
Furthermore, suppose that the following two conditions hold:\\[2.2ex]
(\ref{EWD_4_3_02})
\hspace*{4ex}
\begin{tabular}{@{}l@{}}
$\displaystyle{\sum\limits_{i=1}^{n}\,\big|\,e_{i} - \bm\bar{e}\,\big|^{r}\, \geq\, e n}$,\hspace*{4ex}
$\displaystyle{\sum\limits_{i=1}^{n}\,\big|\,e_{i} - \bm\bar{e}\,\big|^{k}\, \leq\, E n}$\\[2.5ex]
for some $k > 2$, $0 < r < k$ and $e > 0$, $E > 0$.
\end{tabular}\\[3ex]
\refstepcounter{DSBcount}
\label{EWD_4_4_36}
\text{\hspace*{-0.8ex}(\theDSBcount)}
\hspace*{2.8ex}
\begin{tabular}{@{}l@{}}
There exists an $s > 2$ such that $\displaystyle{\ \int\limits_{0}^{1} |\,J(t)\,|^s\,dt < \infty\ }$.
\end{tabular}\vspace*{0.5ex}
\begin{enumerate}
\item\label{EWD_4_4_35_BWa}
Then there exists a constant $\mathcal{E}_{9} > 0$ depending only on $e$, $E$, $r$, $k$ and $J$, $s$
such
\linebreak 
that\\[2ex]
\refstepcounter{DSBcount}
\label{EWD_4_4_37}
\text{\hspace*{-0.8ex}(\theDSBcount)}
\hspace*{2.8ex}
\begin{tabular}[t]{@{}l@{\hspace*{1ex}}c@{\hspace*{1ex}}l@{}}
$||\,\mathscr{F}_{\!A} - e_{1,A}\,||$&$\leq$&
$\displaystyle{\mathcal{E}_{9}\,(\log n)^2\,n^{-\,1\, +\, \bigl((4/k)\, -\, 1\bigr)^{+}\, 
+\, \bigl((4/s)\, -\, 1\bigr)^{+}}}$.
\end{tabular}\vspace*{0.5ex}
\item\label{EWD_4_4_35_BWb}
Let $\epsilon > 0$. Then there exists a constant $\mathcal{E}_{10} > 0$ depending only on 
$e$, $E$, $r$, $k$, $J$, $s$ and $\epsilon$ such that\\[2ex]
\refstepcounter{DSBcount}
\label{EWD_4_4_38}
\text{\hspace*{-0.8ex}(\theDSBcount)}
\hspace*{2.8ex}
\begin{tabular}[t]{@{}l@{\hspace*{1ex}}c@{\hspace*{1ex}}l@{}}
$||\,\mathscr{F}_{\!A} - e_{2,A}\,||$&$\leq$&
$\displaystyle{\mathcal{E}_{10}\,n^{-\,(3/2)\, +\, \epsilon\, +\, \bigl((5/k)\, -\, 1\bigr)^{+}\, 
+\, \bigl((5/s)\, -\, 1\bigr)^{+}}}$.
\end{tabular}\vspace*{0.5ex}
\item\label{EWD_4_4_35_BWc}
Let $k \geq 3$ and suppose that
$J$ satisfies the \textbf{condition} {\boldmath $V_{\alpha}$} (cf. (\ref{EWD_4_4_09})) 
for some $0 < \alpha < \dfrac{1}{3}$.\\[1.5ex]
Then there exists a constant $\mathcal{E}_{11} > 0$ depending only on $e$, $E$, $r$, $k$ and $J$, $\alpha$, $s$
such that\\[2ex]
\refstepcounter{DSBcount}
\label{EWD_4_4_39}
\text{\hspace*{-0.8ex}(\theDSBcount)}
\hfill
\begin{tabular}[t]{@{}l@{\hspace*{1ex}}c@{\hspace*{1ex}}l@{}}
$||\,\mathscr{F}_{\!A} - e_{1,A}^{I}\,||$&$\leq$&
$\displaystyle{\mathcal{E}_{11}\,\max \bigg\{\,n^{-\,(3/2)\, +\, 3\,\alpha}\, ,\,
(\log n)^2\,n^{-\,1\, +\, \bigl((4/k)\, -\, 1\bigr)^{+}\, 
+\, \bigl((4/s)\, -\, 1\bigr)^{+}}\,\bigg\}}$.
\end{tabular}\vspace*{0.5ex}
\item\label{EWD_4_4_35_BWd}
Let $\epsilon > 0$, $k \geq 4$ and suppose that 
$J$ satisfies the \textbf{condition} {\boldmath $V_{\alpha}$} (cf. (\ref{EWD_4_4_09})) 
for some \mbox{\rule[0ex]{0ex}{3.8ex}$0 < \alpha < \dfrac{1}{4}$}.\\[1.5ex]
Then there exists a constant $\mathcal{E}_{12} > 0$ depending only on 
$e$, $E$, $r$, $k$, $J$, $\alpha$, $s$ and $\epsilon$ such that\\[2ex]
\refstepcounter{DSBcount}
\label{EWD_4_4_40}
\text{\hspace*{-0.8ex}(\theDSBcount)}
\hfill
\begin{tabular}[t]{@{}l@{\hspace*{1ex}}c@{\hspace*{1ex}}l@{}}
$||\,\mathscr{F}_{\!A} - e_{2,A}^{I}\,||$&$\leq$&
$\displaystyle{\mathcal{E}_{12}\,\max \bigg\{\,n^{-\,(3/2)\, +\, 3\,\alpha}\, ,\,
n^{-\,(3/2)\, +\, \epsilon\, +\, \bigl((5/k)\, -\, 1\bigr)^{+}\, 
+\, \bigl((5/s)\, -\, 1\bigr)^{+}}\,\bigg\}}$.
\end{tabular}
\end{enumerate}
\vspace*{3.5ex}
\textbf{Proof:}\\[0.8ex]
\textbf{For a) and b):}\\[0.8ex]
To show \ref{EWD_4_4_35_BWa}) and \ref{EWD_4_4_35_BWb}), we use the Theorem \ref{EWD_4_3_13}.
We must therefore prove the existence of positive constants $d$, $D$ and $\delta$ depending 
only on $J$ and $s$ such that\\[2ex]
\refstepcounter{DSBcount}
\label{EWD_4_4_41}
\text{\hspace*{-0.8ex}(\theDSBcount)}
\hspace*{2.8ex}
$\displaystyle{\sum\limits_{j=1}^{n}\,\big|\,d_{j} - \bm\bar{d}\,\big|^{s}\, \leq\, D n}$
\hspace*{3ex}for all $n \in \mathbb{N}$,\\[3ex]
\refstepcounter{DSBcount}
\label{EWD_4_4_42}
\text{\hspace*{-0.8ex}(\theDSBcount)}
\hspace*{2.8ex}
$\displaystyle{\sum\limits_{j=1}^{n}\,\bigl(\,d_{j} - \bm\bar{d}\,\bigr)^{2}\, \geq\, d n}$
\hspace*{3ex}for all $n > n_{2}(J)$,\\[3ex]
\refstepcounter{DSBcount}
\label{EWD_4_4_43}
\text{\hspace*{-0.8ex}(\theDSBcount)}
\hspace*{2.8ex}
\begin{tabular}{@{}l@{}}
$\displaystyle{\lambda\Bigl(\,\bigcup_{j = 1}^{n}\,\Bigl\{\,x \in \mathbb{R}\,:\,\big|\,x - d_{j}\,\big| < \zeta\,
\Bigr\}\,\Bigr)\, \geq\, \delta n \zeta}$\\[3ex]
for some $\zeta \geq n^{-3/2}\,\log\ n$ and all $n \in \mathbb{N}$.
\end{tabular}\\[3.5ex]
In the following we assume $\bm\bar{d} = 0$ without loss of generality. 
This is possible since we can pass from $J$ to 
\mbox{\rule[0ex]{0ex}{5.8ex}$\displaystyle{\widetilde{J} = J - \int\limits_{0}^{1} J(t)\,dt}$ because of
$\displaystyle{\bm\bar{d} = \int\limits_{0}^{1} J(t)\,dt}$} 
(cf. Lemma \ref{EWD_4_4_29}, \ref{EWD_4_4_29_BWb}) with $g = J$).\\[2.5ex]
Furthermore, we can assume in the following that there are $\gamma,\, \theta > 0$ 
and $\theta < a < b < 1 - \theta$ 
\linebreak
with\\[2ex]
\refstepcounter{DSBcount}
\label{EWD_4_4_44}
\text{\hspace*{-0.8ex}(\theDSBcount)}
\hspace*{2.8ex}
$J(x)\, \geq\, \gamma$\hspace*{2ex}and\hspace*{2ex}$J'(x)\, \geq\, \gamma$
\hspace*{3ex}for $x \in T = [\,a - \theta,\,b +\theta\,]$.\\[2.5ex]
The other three possible cases $J(x)\, \leq\, -\, \gamma$, $J'(x)\, \leq\, -\, \gamma\,$ 
and $\,J(x)\, \geq\, \gamma$, $J'(x)\, \leq\, -\, \gamma\,$ and \\
$\,J(x)\, \leq\, -\, \gamma$, $J'(x)\, \geq\, \gamma$ are handled completely analogously
(note that $J^2 = (- J)^2$ and consider $d_{j} - d_{j+1}$ for negative $J'$).\\[2.8ex]
Now, for $a\, \leq\, \dfrac{j}{n}\, \leq\, b$ we obtain\\[2ex]
\refstepcounter{DSBcount}
\label{EWD_4_4_45}
\text{\hspace*{-0.8ex}(\theDSBcount)}
\hspace*{2.8ex}
$P\bigl(\,U_{j:n}\,\leq\, a - \theta\,\bigr)\, \leq\, e^{-\,2n\theta^2}$
\hspace*{2ex}and\hspace*{2ex}
$P\bigl(\,U_{j:n}\,\geq\, b + \theta\,\bigr)\, \leq\, e^{-\,2n\theta^2}$.\\[2.5ex]
For a justification of (\ref{EWD_4_4_45}) see e.g. Okamoto \cite{Okamoto1959}, Theorem 1 on page 33 and\\[2ex]
\hspace*{12.1ex}$P\bigl(\,U_{j:n}\,\leq\, a - \theta\,\bigr)\, =\,
P\Bigl(\,b(n,\,a - \theta)\, \geq\, j\,\Bigr)\, \leq\, P\Bigl(\,\dfrac{1}{n}\,b(n,\,a - \theta)
\, \geq\, (\,a - \theta\,) + \theta\,\Bigr)$,\\[2.5ex]
where $b(n,\,x)$ is a binomially distributed random variable with the parameters $n$ and 
$x$.\index{random variable!binomially distributed}\index{random variable!binomially distributed!$b(n,\,x)$}\\[2.8ex]
After these preliminary considerations, we will now show the three statements above:\\[2.8ex]
\textbf{Proof of (\ref{EWD_4_4_41}):}\\[2ex]
\hspace*{12.1ex}$\displaystyle{\sum\limits_{j=1}^{n}\,|\,d_{j}\,|^s
\, \leq \, \sum\limits_{j=1}^{n}\,E\bigl(\,|\,J(\,U_{j:n}\,)\,|^s\,\bigr)
= n \int\limits_{0}^{1} |\,J(t)\,|^s\,dt}$.\\[2.5ex]
For {''$=$''}, Lemma \ref{EWD_4_4_29}, \ref{EWD_4_4_29_BWb}) with $g = |J|^s$, $s > 2$, was used.
$|J|^s$ is integrable according to condition (\ref{EWD_4_4_36}).\\[2.8ex]
\textbf{Proof of (\ref{EWD_4_4_42}):}\\[0.8ex]
There exists an $n' \in \mathbb{N}$ such that for all $n \geq n'$ and
$j$ with $a\, \leq\, \dfrac{j}{n}\, \leq\, b$ holds:\\[2ex]
\hspace*{12.1ex}\begin{tabular}[t]{@{}l@{\hspace*{1ex}}c@{\hspace*{1ex}}l@{\hspace*{32ex}}r@{}}
$\displaystyle{E\Bigl(\,(J 1_{T})(\,U_{j:n}\,)\,\Bigr)}$&$\geq$&$\gamma\,P\bigl(\,U_{j:n} \in T\,\bigr)$
&(cf. (\ref{EWD_4_4_44}))\\[2.3ex]
&$\geq$&$\displaystyle{\gamma\,\Bigl(\,1 - 2\,e^{- 2 n \theta^{2}}\,\Bigr)}$
&(cf. (\ref{EWD_4_4_45}))\\[2.3ex] 
&$\geq$&$\dfrac{\gamma}{2}$
\end{tabular}\\[2.5ex]
and\\[2.5ex]
\refstepcounter{DSBcount}
\label{EWD_4_4_46}
\text{\hspace*{-0.8ex}(\theDSBcount)}
\hspace*{2.8ex}
\begin{tabular}[t]{@{}l@{\hspace*{1ex}}l@{\hspace*{4.4ex}}r@{}}
&$\Big|\,E\Bigl(\,(J 1_{T^{c}})(\,U_{j:n}\,)\,\Bigr)\,\Big|$\\[2.3ex]
$\leq$&$\Bigl(\,E\bigl(\,J^2(\,U_{j:n}\,)\,\bigr)\,\Bigr)^{1/2}\,
\Bigl(\,P\bigl(\,U_{j:n} \not\in T\,\bigr)\,\Bigr)^{1/2}$
&(H{\"o}lder's inequality with $p = q = 2$\index{H{\"o}lder's inequality!for random variables})\\[2.3ex]
$\leq$&$\displaystyle{\biggl(\,n \int\limits_{0}^{1} J^2(t)\,dt\,\biggr)^{1/2}\,
\Bigl(\,2\,e^{-\,2n\theta^2}\,\Bigr)^{1/2}}$
&(cf. Lemma \ref{EWD_4_4_29}, \ref{EWD_4_4_29_BWb}) with $g = J^2$)\\[4ex]
$\leq$&$\dfrac{\gamma}{4}$.
\end{tabular}\\[2.5ex]
It follows for all $n \geq \max\Bigl\{\,n',\,\dfrac{2}{b-a}\,\Bigr\}$ that\\[2ex]
\hspace*{12.1ex}\begin{tabular}[t]{@{}l@{\hspace*{1ex}}c@{\hspace*{1ex}}l@{}}
$\displaystyle{\sum\limits_{j=1}^{n}\,d_{j}^{\,2}}$&$\geq$&
$\displaystyle{\sum\limits_{a\, \leq \textstyle{\frac{j}{n}}\, \scriptstyle{\leq\, b}}\,
\biggl\{\,E\Bigl(\,(J 1_{T})(\,U_{j:n}\,)\,\Bigr) + E\Bigl(\,(J 1_{T^{c}})(\,U_{j:n}\,)\,\Bigr)\,\biggr\}^2}$\\[4.5ex]
&$\geq$&$\displaystyle{\sum\limits_{na\, \leq\, j\, \leq\, nb}\, 
\Bigl\{\,\dfrac{\gamma}{2} - \dfrac{\gamma}{4}\,\Bigr\}}^2$\\[4ex]
&$\geq$&$\Bigl(\,\dfrac{\gamma}{4}\,\Bigr)^2\,\bigl(\,(b-a)\,n -1\,\bigr)$\\[3ex]
&$\geq$&$\Bigl(\,\dfrac{\gamma}{4}\,\Bigr)^2\,\dfrac{b-a}{2}\,n$.
\end{tabular}\\[2.5ex]
From this we get (\ref{EWD_4_4_42}).\\[2.8ex]
\textbf{Proof of (\ref{EWD_4_4_43}):}\\[0.8ex]
Let $V_{j,n} = \Bigl\{\,U_{j:n} \in T,\,U_{j+1:n} \in T\,\Bigr\}$.
We then obtain, as above, an $n'' \in \mathbb{N}$ such that for all $n \geq n''$ and $j$ 
with $a\, \leq\, \dfrac{j}{n}\, \leq\, b - \dfrac{1}{n}$ and $1\, \leq\, l\, \leq\, n$ holds:\\[2ex]
\hspace*{12.1ex}\begin{tabular}[t]{@{}l@{\hspace*{1ex}}c@{\hspace*{1ex}}l@{\hspace*{26.8ex}}r@{}}
$P\bigl(\,V_{j,n}^{c}\,\bigr)$&$\leq$&$P\bigl(\,U_{j:n} \not\in T\,\bigr)
+ P\bigl(\,U_{j+1:n} \not\in T\,\bigr)$\\[2.3ex]
&$\leq$&$4\,e^{-\,2n\theta^2}$
&(cf. (\ref{EWD_4_4_45}))\\[2.3ex]
&$\leq$&$\dfrac{1}{9n}$
\end{tabular}\\[2.5ex]
and\\[2ex]
\hspace*{12.1ex}\begin{tabular}[t]{@{}l@{\hspace*{1ex}}l@{\hspace*{9.2ex}}r@{}}
&$\Big|\,E\Bigl(\,J(\,U_{l:n}\,)\,1_{V_{j,n}^{c}}\,\Bigr)\,\Big|$\\[2.3ex]
$\leq$&$\Bigl(\,E\bigl(\,J^2(\,U_{l:n}\,)\,\bigr)\,\Bigr)^{1/2}\,
\Bigl(\,P\bigl(\,V_{j,n}^{c}\,\bigr)\,\Bigr)^{1/2}$
&(H{\"o}lder's inequality with $p = q = 2$\index{H{\"o}lder's inequality!for random variables})\\[2.3ex]
$\leq$&$\displaystyle{\biggl(\,n \int\limits_{0}^{1} J^2(t)\,dt\,\biggr)^{1/2}\,
\Bigl(\,4\,e^{-\,2n\theta^2}\,\Bigr)^{1/2}}$
&(cf. Lemma \ref{EWD_4_4_29}, \ref{EWD_4_4_29_BWb}) with $g = J^2$)\\[4ex]
$\leq$&$\dfrac{\gamma}{9n}$.
\end{tabular}\\[2.5ex]
It follows for all $n \geq n''$ and $j$ with $a\, \leq\, \dfrac{j}{n}\, \leq\, b - \dfrac{1}{n}$ that\\[2ex]
\hspace*{12.1ex}\begin{tabular}[t]{@{}l@{\hspace*{1ex}}l@{\hspace*{-22.7ex}}r@{}}
&$d_{j+1} - d_{j}$\\[2.5ex]
$=$&$E\Bigl(\,\Bigl(\,J(\,U_{j+1:n}\,) - J(\,U_{j:n}\,)\,\Bigr)\,1_{V_{j,n}}\,\Bigr)\, +\,
E\Bigl(\,\Bigl(\,J(\,U_{j+1:n}\,) - J(\,U_{j:n}\,)\,\Bigr)\,1_{V_{j,n}^{c}}\,\Bigr)$\\[2.5ex]
$\geq$&$\gamma\,E\Bigl(\,\bigl(\,U_{j+1:n} - U_{j:n}\,\bigr)\,1_{V_{j,n}}\,\Bigr)\, 
-\, \dfrac{2\gamma}{9n}$&(mean value theorem\index{mean value theorem})\\[2.5ex] 
$=$&$\gamma\,E\Bigl(\,U_{j+1:n} - U_{j:n}\,\Bigr)\, 
-\, \gamma\,E\Bigl(\,\bigl(\,U_{j+1:n} - U_{j:n}\,\bigr)\,1_{V_{j,n}^{c}}\,\Bigr)\, 
-\, \dfrac{2\gamma}{9n}$\\[2.5ex]
$\geq$&$\gamma\,E\Bigl(\,U_{j+1:n} - U_{j:n}\,\Bigr)\, 
-\, \gamma\,P\bigl(\,V_{j,n}^{c}\,\bigr)\, 
-\, \dfrac{2\gamma}{9n}$&(since $0 \leq U_{j+1:n} - U_{j:n} \leq 1$)\\[2.5ex]
$\geq$&$\dfrac{\gamma}{n+1}\, -\, \dfrac{\gamma}{9n}\, -\, \dfrac{2\gamma}{9n}$
&(due to Lemma \ref{EWD_4_4_29}, \ref{EWD_4_4_29_BWa}))\\[3ex] 
$\geq$&$\dfrac{\gamma}{2n}\, -\, \dfrac{\gamma}{3n} = \dfrac{\gamma}{6n}$.
\end{tabular}\\[3ex]
From this we get (\ref{EWD_4_4_43}) first for sufficiently large $n \geq n_{0}$ 
(with e.g. $\delta = b - a$ and $\zeta = \dfrac{\gamma}{12n}$)
and then for all $n \in \mathbb{N}$.\\[4ex]
\textbf{For c) and d):}\\[0.8ex]
Since the proofs of \ref{EWD_4_4_35_BWc}) and \ref{EWD_4_4_35_BWd}) are completely analogous, 
we only show \ref{EWD_4_4_35_BWd}).\\[2.8ex]
Because of \ref{EWD_4_4_35_BWb}) it is sufficient to prove the following estimate:\\[2ex]
\refstepcounter{DSBcount}
\label{EWD_4_4_47}
\text{\hspace*{-0.8ex}(\theDSBcount)}
\hspace*{2.8ex}
$\displaystyle{||\,e_{2,A} - e_{2,A}^{I}\,||\, \leq\, c_{1}\,n^{-\,(3/2)\, +\, 3\,\alpha}}$
\hspace*{4ex}for $n > n_{2}(J)$.\\[2.5ex]
Here, $c_{1}$ depends only on $e$, $E$, $r$, $k$ and $J$, $\alpha$.\\[2.8ex]
But now the proof of (\ref{EWD_4_4_47}) is analogous to that of (\ref{EWD_4_4_27}) 
if we use Lemma \ref{EWD_4_4_30} for $w = 2, 3, 4$ and Lemma \ref{EWD_4_4_29}, \ref{EWD_4_4_29_BWb}) 
instead of Corollary \ref{EWD_4_4_15} for $w = 2, 3, 4$.\hspace*{1ex}\hfill$\Box$\\[4ex] 
\refstepcounter{DSBcount}
\label{EWD_4_4_48}
\textbf{\hspace*{-0.8ex}\theDSBcount\ Remarks}
\begin{enumerate}
\item\label{EWD_4_4_48_BWa}
Because of Lemma \ref{EWD_4_4_08}, \ref{EWD_4_4_08_BWb}) the condition (\ref{EWD_4_4_36}) is automatically 
satisfied with $s = 3$ in part \ref{EWD_4_4_35_BWc}) and with $s = 4$ in part \ref{EWD_4_4_35_BWd}).
\item\label{EWD_4_4_48_BWb}
The Remark \ref{EWD_4_4_28}, \ref{EWD_4_4_28_BWa}) applies analogously to the Theorem \ref{EWD_4_4_35}.
\item\label{EWD_4_4_48_BWc}
If $J$ is nonconstant, continuous and 
$\displaystyle{\ \int\limits_{0}^{1} |\,J(t)\,|^p\,dt < \infty\ }$ for a $p > 1$, then 
$n_{2}(J) < \infty$.\\[1.5ex]
This follows as in the {''Proof of (\ref{EWD_4_4_42})''}, if 
H{\"o}lder's inequality\index{H{\"o}lder's inequality!for random variables}
in (\ref{EWD_4_4_46}) is applied with this $p > 1$ instead of $p = 2$.
\item\label{EWD_4_4_48_BWd}
In the Remark \ref{EWD_4_4_28}, \ref{EWD_4_4_28_BWd}) it has already been mentioned that the functions
$J(t) = t$ (Wil\-co\-xon-Scores\index{scores!exact!Wilcoxon}, 
since \mbox{\rule[-3ex]{0ex}{6.8ex}$E\bigl(\,U_{j:n}\,\bigr) = \dfrac{j}{n+1}$}) and
$J(t) = \Phi^{-1}(t)$ (Fisher-Yates-Terry-Hoeff\-ding-Scores\index{scores!exact!Fisher-Yates-Terry-Hoeffding})
satisfy the \textbf{condition} {\boldmath $V_{\alpha}$} (cf. (\ref{EWD_4_4_09})) for any $\alpha > 0$.
\end{enumerate}
                               
\chapter*{Comparison of different numberings}

\addcontentsline{toc}{chapter}{\text{Comparison of different numberings}}

In the following list, citations from \cite{10.1214/aos/1176347258}, which refer to the 
original version of 1987, are compared with the changed citations from this updated version of 2025.\\[2ex]
\hspace*{12.1ex}\begin{tabular}[t]{|l@{\hspace*{5ex}}|l@{\hspace*{5ex}}|} \hline
Original version of 1987&Updated version of 2025\\ \hline \hline
Lemma 2.1.5&Lemma \ref{EWD_2_1_13}\\ \hline
Proposition 4.3.9&Proposition \ref{EWD_4_3_09}\\
pages 148 - 151&pages \pageref{EWD_4_3_09} - \pageref{EWD_4_3_11}\\ \hline
Theorem 4.4.13(a), (b), (d)&Theorem \ref{EWD_4_4_21}, \ref{EWD_4_4_21_BWa}), \ref{EWD_4_4_21_BWb}), 
\ref{EWD_4_4_21_BWd})\\
Theorem 4.4.23(a), (b)&Theorem \ref{EWD_4_4_35}, \ref{EWD_4_4_35_BWa}), \ref{EWD_4_4_35_BWb})\\ \hline
\end{tabular}\\[4ex]
The order of the chapters and sections is the same in the Ph. D. thesis of 1987 and 
the updated version of 2025.\\[1.5ex] 
Only the chapter numbers have been added to the numbering of the sections in this updated version.
For example, section 3 from chapter 1 is now section 1.3.

\renewcommand{\indexname}{Index of keywords and symbols}
\rehead{Index of keywords and symbols} 
\printindex
\rehead{References}    
\printbibliography[heading=bibintoc,title={References}]
	
\end{document}